%% file: ButtenschoenHillenAdhStSt.tex
\documentclass{memo-l}

\usepackage[utf8]{inputenc}
\usepackage[T1]{fontenc}
\makeindex

\usepackage[english]{babel}

\usepackage{csquotes}
\usepackage{enumitem}
\usepackage{pgf}
\usepackage{pgfplots}
\usetikzlibrary{matrix}
\usepgfplotslibrary{groupplots}
\pgfplotsset{compat=newest}
\usepackage{array}
\usepackage{ragged2e}
\usepackage{tabularx}
\usepackage{booktabs}
\usepackage{todonotes}
\usepackage{breqn}
\usepackage{caption}
\usepackage{color}
\definecolor{hillencolor}{rgb}{0.30,0.03,0.66}
\definecolor{hintergrundfarbe}{rgb}{0.8, 0.8,  0.982}
\definecolor{notefarbe}{rgb}{0.9, 0.9,  0.982}

\graphicspath{{./include/}}

\usepackage{xparse}
\usepackage{amsrefs}
\usepackage{amssymb}

\NewDocumentCommand\parencites{mgggggggg}{%
    \IfNoValueTF{#9}{%
    \IfNoValueTF{#8}{%
    \IfNoValueTF{#7}{%
    \IfNoValueTF{#6}{%
    \IfNoValueTF{#5}{%
    \IfNoValueTF{#4}{%
    \IfNoValueTF{#3}{%
    \IfNoValueTF{#2}{%
        \cite{#1}\xspace%
    }{\cites{#1,#2}\xspace}%
    }{\cites{#1,#2,#3}\xspace}%
    }{\cites{#1,#2,#3,#4}\xspace}%
    }{\cites{#1,#2,#3,#4,#5}\xspace}%
    }{\cites{#1,#2,#3,#4,#5,#6}\xspace}%
    }{\cites{#1,#2,#3,#4,#5,#6,#7}\xspace}%
    }{\cites{#1,#2,#3,#4,#5,#6,#7,#8}\xspace}%
    }{\cites{#1,#2,#3,#4,#5,#6,#7,#8,#9}\xspace}%
}

\newcommand{\parencite}[1]{\cite{#1}}

\usepackage{xspace}
\usepackage{mathtools}
\DeclarePairedDelimiter{\abs}{\lvert}{\rvert}
\DeclarePairedDelimiterXPP\seq[1]{}{(}{)}{}{#1}
\DeclarePairedDelimiter{\norm}{\lVert}{\rVert}

\DeclareMathOperator{\codim}{codim}
\DeclareMathOperator{\ind}{ind}
\DeclareMathOperator{\vspan}{span}

\DeclareMathOperator{\coker}{coker}
\DeclareMathOperator{\interior}{int}
\DeclareMathOperator{\dist}{dist}
\DeclareMathOperator{\Bd}{\mathcal{B}}
\DeclareMathOperator{\M}{\mathcal{M}}
\DeclareMathOperator{\K}{\mathcal{K}}
\DeclareMathOperator{\Lap}{\mathcal{L}}
\DeclareMathOperator{\D}{\mathcal{D}}
\DeclareMathOperator{\F}{\mathcal{F}}
\DeclareMathOperator{\Top}{\mathcal{T}}

\DeclareMathOperator{\I}{\mathcal{I}}
\DeclareMathOperator{\Avg}{\mathcal{A}}
\DeclareMathOperator{\Null}{\mathrm{N}}
\DeclareMathOperator{\Range}{\mathrm{R}}
\DeclareMathOperator{\Cont}{\mathfrak{C}}
\DeclareMathOperator{\Soln}{\mathfrak{S}}
\DeclareMathOperator{\Fred}{Fred}
\DeclareMathOperator{\Iso}{Iso}
\DeclareMathOperator{\Lin}{\mathfrak{L}}
\DeclareMathOperator{\Kom}{\mathfrak{K}}
\DeclareMathOperator{\C}{\mathcal{C}}
\DeclareMathOperator{\Pos}{\mathcal{P}}
\DeclareMathOperator{\Nodal}{\mathcal{S}}

\DeclareMathOperator{\O2}{\mathbf{O}(2)}
\DeclareMathOperator{\SO2}{\mathbf{SO}(2)}
\DeclareMathOperator{\Dn}{\mathbf{D}_{n}}
\newcommand*\mean[1]{\bar{#1}}

\makeatletter
\newcommand*{\rom}[1]{\expandafter\@slowromancap\romannumeral #1@}
\makeatother

\newcommand{\ROM}[1]{%
    \textup{\uppercase\expandafter{\romannumeral#1}}%
}

\newcommand{\dd}{\ensuremath{\, \mathrm{d}}}

\newcommand{\Sn}{\mathbb{S}^{n-1}}
\newcommand{\eg}{e.g.\xspace}

\newcommand{\ie}{i.e.\xspace}
\newcommand{\wrt}{w.r.t.\xspace}
\newcommand{\etal}{\emph{et al.}}

\newcommand{\symdif}{\ensuremath{\mathrel{\triangle}}}

\newcommand{\BigOh}[1]{\mathcal{O}(#1)}
\newcommand{\ifrac}[2]{\ensuremath{{}^{#1}\!/_{#2}}}
\newcommand{\N}{\mathbb{N}}
\newcommand{\Z}{\mathbb{Z}}
\newcommand{\R}{\mathbb{R}}
\newcommand{\OO}{\mathcal{O}}

\newcommand{\lb}{\left(}
\newcommand{\rb}{\right)}
\newcommand{\lcb}{\left\{\, }
\newcommand{\rcb}{\,\right\} }
\newcommand{\lsb}{\left[\, }
\newcommand{\rsb}{\,\right] }
\newcommand{\CpInf}{\ensuremath{\C_{\mathrm{per}}^{\infty}}}
\newcommand{\CInf}{\ensuremath{\C^{\infty}}}
\newcommand{\HperOne}{\ensuremath{H_{\mathrm{per}}^{1}}}
\newcommand{\HOne}{\ensuremath{H^{1}}\xspace}

\newcommand{\WperTwo}{\ensuremath{W^{2,2}_{\mathrm{per}}}}
\newcommand{\ssubset}{\ensuremath{\subset\joinrel\subset}}
\newcommand{\sgn}{\operatorname{sgn}}

\newcommand{\nonlocalgradinternal}{%
\ensuremath{\vbox{\offinterlineskip\ialign{%
    \hfil##\hfil\cr
    $\scriptsize\circ$\cr
    \noalign{\kern0ex}
    $\nabla$\cr
}}}}

\newcommand{\nonlocalgrad}[1]{\ensuremath{\mathbin{\nonlocalgradinternal_{\scriptsize#1}}}}

\usepackage{amsmath}
\expandafter\let\csname ver@amsthm.sty\endcsname\relax
\let\theoremstyle\relax
\let\newtheoremstyle\relax
\let\pushQED\relax
\let\popQED\relax
\let\qedhere\relax
\let\mathqed\relax
\let\openbox\relax
\let\proof\relax
\usepackage{amsthm}

\makeatletter
\renewcommand{\thm@space@setup}{
  \thm@preskip=.5\baselineskip\@plus.2\baselineskip \@minus.2\baselineskip
  \thm@postskip=\thm@preskip
}
\newtheoremstyle{amsplain}
{\thm@preskip}
{\thm@postskip}
{\itshape}
{\parindent}
{\scshape}
{.}
{ }
{}
\newtheoremstyle{amsdefinition}
{\thm@preskip}
{\thm@postskip}
{\normalfont}
{\parindent}
{\scshape}
{.}
{ }
{}
\newtheoremstyle{amsremark}
{\thm@preskip}
{\thm@postskip}
{\normalfont}
{\parindent}
{\scshape}
{.}
{ }
{}
\makeatother

\usepackage[capitalise]{cleveref}

\theoremstyle{amsplain}
\newtheorem{theorem}{Theorem}[chapter]
\newtheorem{lemma}[theorem]{Lemma}
\newtheorem{proposition}[theorem]{Proposition}
\newtheorem{corollary}[theorem]{Corollary}

\theoremstyle{amsdefinition}
\newtheorem{definition}[theorem]{Definition}
\newtheorem{example}[theorem]{Example}

\theoremstyle{amsremark}
\newtheorem{remark}[theorem]{Remark}

\numberwithin{section}{chapter}
\numberwithin{equation}{chapter}

\crefname{lemma}{Lemma}{Lemmas}
\Crefname{lemma}{Lemma}{Lemmas}
\crefname{proposition}{Proposition}{Propositions}
\Crefname{proposition}{Proposition}{Propositions}
\crefname{theorem}{Theorem}{Theorems}
\Crefname{theorem}{Theorem}{Theorems}
\crefname{definition}{Definition}{Definitions}
\Crefname{definition}{Definition}{Definitions}
\crefname{remark}{Remark}{Remarks}
\Crefname{remark}{Remark}{Remarks}
\crefname{example}{Example}{Examples}
\Crefname{example}{Example}{Examples}
\crefname{corollary}{Corollary}{Corollaries}
\Crefname{corollary}{Corollary}{Corollaries}

\makeindex

\makeatletter
\def\input@path{{./include/}}
\makeatother

\newcolumntype{b}{X}
\newcolumntype{s}{>{\hsize=.5\hsize}X}
\newcolumntype{Y}{>{\RaggedRight\arraybackslash}X}
\begin{document}

\frontmatter

\title{Non-Local Cell Adhesion Models:\\ Steady States and Bifurcations}


\author{Andreas Buttensch\"{o}n}
\address{University of British Columbia}
\email{abuttens@math.ubc.ca}
\thanks{}

\author{Thomas Hillen}
\address{University of Alberta}
\email{thillen@ualberta.ca}
\thanks{}

\date{}




\begin{abstract}
    Whenever cells form tissues, organs, or organisms, cells interact with each
    other through cellular adhesions. Cell-cell adhesions give the skin its
    stability, it keeps cells together to form organs, they allow immune cells
    to move through the body, and they keep blood inside the vessels. Cell-cell
    adhesions also facilitate diseases such as cancer and cancer metastasis, for
    example. A good understanding of this basic cell mechanism is of upmost
    importance. Here we use the toolbox of mathematical modelling to help to
    gain insight into cell-cell adhesions.

    Mathematical modelling of cellular adhesions has focussed on individual
    based models, where each cell is described as its own entitiy, and
    interaction rules are defined. Early attempts to find continuum models for
    cellular adhesion failed, as they led to problems with backward diffusion.
    The continuum description of cell adhesion remained a challenge until 2006,
    when {\it Armstrong, Painter and Sherratt}\ proposed the use of an
    integro-partial differential equation (iPDE) model for cell adhesion. The
    initial success of the model was the replication of the cell-sorting
    experiments of Steinberg. Since then the non-local adhesion model of
    Armstrong et.\ al.\ has proven popular in applications to embryogenesis, wound
    healing, and cancer cell invasions.

    The mathematical analysis of the non-local adhesion model is challenging. In
    this monograph, we contribute to the analysis of steady states and their
    bifurcation structure.  The importance of steady-states is that these are
    the patterns observed in nature and tissues (\eg\ cell-sorting experiments).
    In the case of periodic boundary conditions, we combine global bifurcation
    results pioneered by Rabinowitz, equivariant bifurcation theory, and the
    mathematical properties of the non-local term to obtain a global bifurcation
    result for the branches of non-trivial solutions.

    We further extend the non-local cell adhesion model to a bounded domain with
    no-flux boundary conditions. Using the derivation of the non-local adhesion
    model we propose several different non-local terms incorporating
    biologically realistic boundary effects. We show that these newly
    constructed non-local operators are weakly differentiable, using the theory
    of distributions. We find that these boundary conditions can include
    boundary adhesion or boundary repulsion, such that boundary-wetting effects
    can be obtained. Finally, we use an asymptotic expansion to study the steady
    states of the non-local cell adhesion model incorporating no-flux boundary
    conditions and we find a similar bifurcation structure as in the periodic
    case.

\end{abstract}

\maketitle

\centerline{\includegraphics[width=8cm]{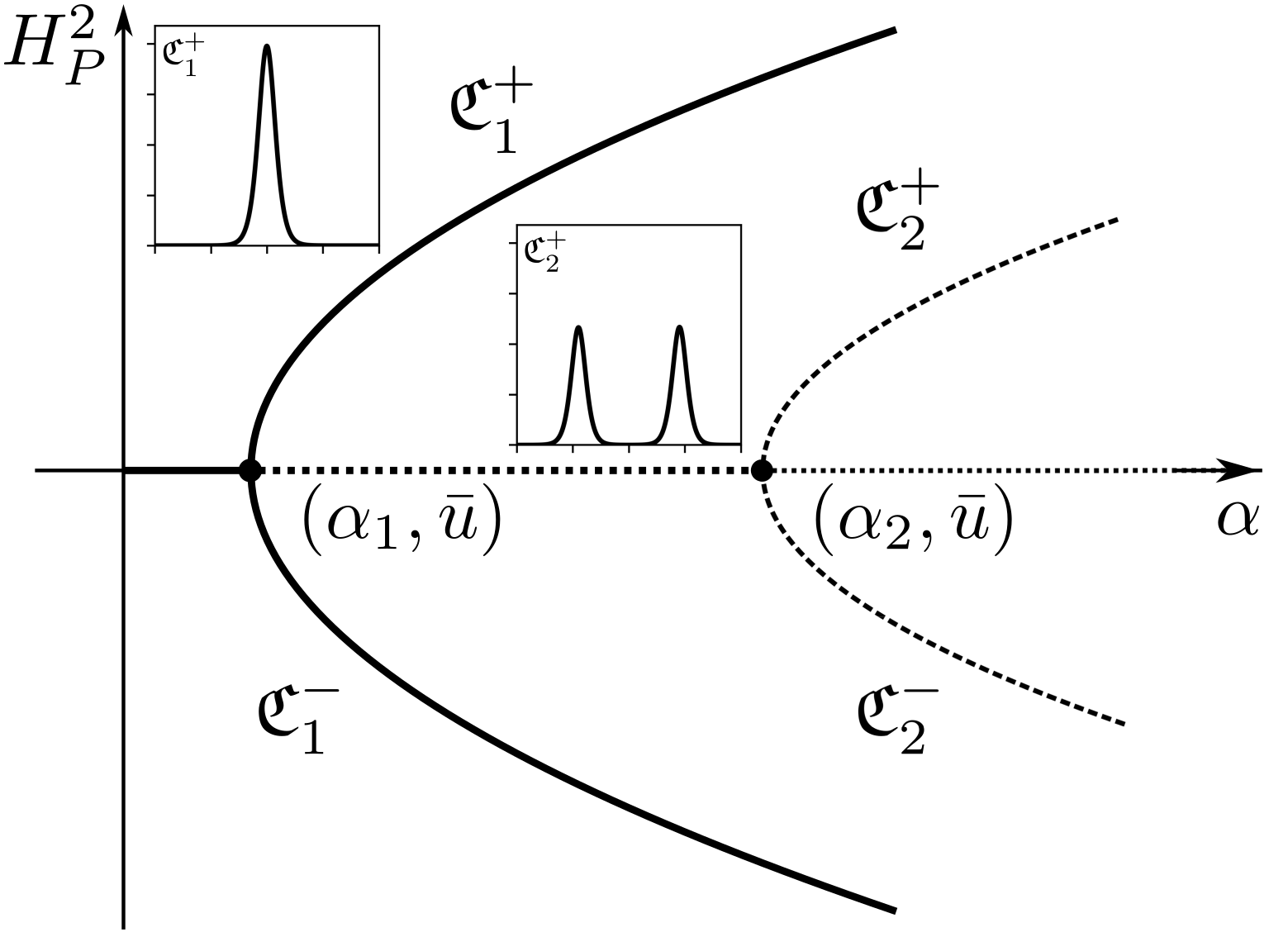}}
\tableofcontents


\mainmatter
\part{Introduction}\label{part:intro}
\include{chapterIntroduction}

\include{chapterPreliminaries}

\part{The Periodic Problem}\label{part:periodic}
\include{chapterNonLocalOperator}
\include{chapterLocalBif}
\include{chapterGlobalBif}

\part{Non-local Equations with Boundary Conditions}\label{part:noflux}
\include{chapterNoFluxBcShort}


\include{chapterDiscussion}

\appendix
\include{AppendixAmorim}

\backmatter
\bibliography{bibliography}
\printindex

\end{document}


%% file: chapterIntroduction.tex
\chapter{Introduction}

Cellular adhesion\index{adhesion} is one of the most important interaction
forces\index{interaction forces} in tissues.  Cells adhere to each other, to
other cells and to the extracellular matrix (ECM). Cell adhesion is responsible
for the formation of tissues, membranes, vasculature, muscle tissue as well as
cell movement and cancer\index{cancer} spread.  At the molecular level, cellular
adhesion is facilitated by a wide range of different cell-membrane
proteins\index{cell-membrane proteins} with integrins\index{integrins} and
cadherins\index{cadherins} being the most prominent adhesion
molecules~\cite{Halbleib2006a,Friedl2017}.  We recognize that cell adhesions are
fundamental for the normal functions of organs\index{organs}, for embryonic
development\index{embryonic development}, wound healing\index{wound healing}, as
well as pathological issues such as cancer metastasis\index{metastasis}
\parencites{McMillen2015}{Graves2016}{Puliafito2012}.

A good understanding of adhesion and its dynamic properties is essential and
mathematical modelling is one powerful tool to gain such an understanding. There
have been several modelling attempts for adhesion and it turns out that one of
the more successful models is the model of Armstrong, Painter and Sherratt from
2006~\cite{Armstrong2006}. It has the form of a non-local partial differential
equation, where the particle flux is an integral-term that arises as balance of
all the adhesion forces acting on a cell.

The present monograph focuses on mathematical properties of the non-local
Armstrong model. The non-local nature of the particle flux term is a challenge
and sophisticated new methods need to be derived. Here we show the existence of
non-trivial steady states, analyse their stability and their bifurcation
structure. The results are largely based on the abstract bifurcation
theory\index{bifurcation theory} of Crandall and Rabinowitz\index{Crandall and
Rabinowitz}~\cite{Rabinowitz1971}.  We show that the non-local term acts like
a non-local derivative\index{non-local derivative}, which allows us to define
non-local gradients\index{non-local gradient} and non-local
curvature\index{non-local curvature}. Furthermore, we discuss the development of
appropriate, and biologically realistic, no-flux boundary
conditions\index{no-flux boundary conditions}, and we show the existence of
non-trivial steady states for this case. As the no-flux boundary conditions are
non-unique, they open the door for further studies of boundary behavior of cells
on tissue boundaries.

\section{The Effect of Cellular Adhesions in Tissues}

Early in the last century the first biological experimenters had begun to
uncover the role of cell adhesion in tissues.  One of the earliest observations
was that if a sponge\index{sponge} is squeezed through a fine mesh (Wilson, 1907
\parencite{Steinberg1963}), it will reform into a functional sponge after
transition. A few years later, Hoftreter observed that different tissues have
different associative preferences \parencite{Steinberg1963}. To describe his
observations, he introduced the concept of
\textquote{tissue affinities\index{tissue affinities}}.
Further, he repeated the earlier observations that previously dissociated
tissues have the ability to regain their form and function after deformation.
Today, this phenomenon is referred to as cell-sorting\index{cell-sorting}, and
we recognize its critical importance in the formation of functional tissues
during organism development.

\begin{figure}[!ht]\centering
    \includegraphics[width=8cm]{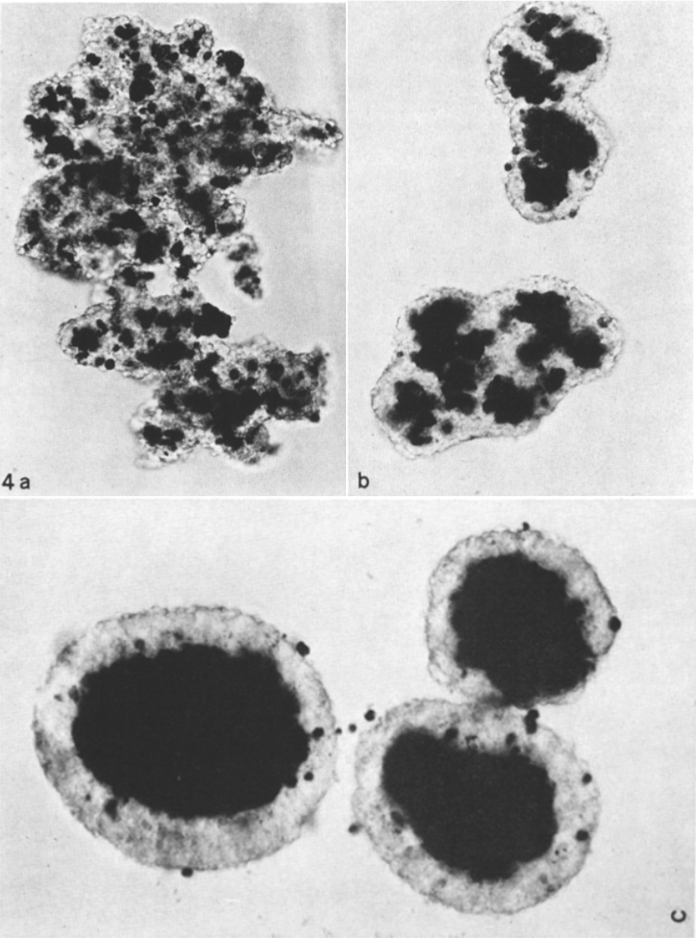}
    \caption[Biological cell-sorting experiment]
    {Two cell populations, black and white, with adhesion molecules of different
    strength on their surfaces. Initially (a) the cells are mixed, and over time
    they slowly (b) re-sort themselves to a sorted final configuration (c) (For
    details on the experimental setup and the figure see
    \parencite{Armstrong1971}).}\label{Fig:ArmstrongExperiment}
\end{figure}

In 1963 Steinberg, proposed the first theory of cell-sorting\index{cell-sorting}
that argued that cell-level properties, namely a cell's adhesion molecules,
drive cell-sorting
\parencites{Steinberg1963}{Steinberg1970}. His theory, capable of explaining the
different cell-sorting patterns is known as the \textit{Differential Adhesion
Hypothesis\index{Differential Adhesion Hypothesis}} (DAH). Steinberg observed
that clusters of cells of the same type behave as if a surface
tension\index{surface tension} holds them together, quite similar to
fluid droplets. In other words, cells rearrange to maximize their
intra-cellular attraction and minimize surface tension. That is, the DAH asserts
that cell-sorting is solely driven by the quantitative differences in the
adhesion potential between cell types (\eg\ cells with the highest potential
of adhesion are found at the centre of aggregates). Interesting is that
Steinberg referred to adhesion being a \textquote{merely close range attraction}
\parencite{Steinberg1963}. An overview of the experimental verifications can be
found in \parencite{Steinberg2007}.

\begin{figure}
    \includegraphics[width=8cm]{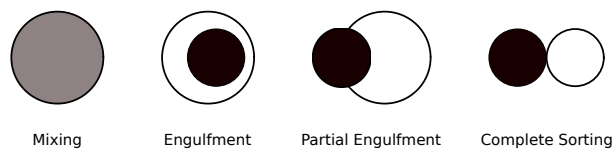}
    \caption[Schematic of outcomes of cell-sorting with two cell populations]
    {The four possible outcomes of cell-sorting with two cell
    populations. The more cohesive cell population is black.
    Mixing occurs with preferential cross-adhesion, engulfment with intermediate
    cross-adhesion, partial engulfment with weak cross-adhesion, and
    cell-sorting with no cross-adhesion \parencites{Foty2004}{Armstrong2006}.
    }\label{Fig:CellSortingOutcomes}
\end{figure}

Further analysis of cell sorting patterns \parencites{Foty2004}{Armstrong2006}
has shown several types of cell clusters as shown in
\cref{Fig:CellSortingOutcomes}. The outcome depends on the relative adhesion
strengths of cells of the same type to each other and to cells of a different
type.

Harris formulated a first critique of the DAH \parencite{Harris1976}. The main
points of his critique were: (1) cells are living objects, and thus open
thermodynamic systems (not closed as assumed by the DAH), (2) cell size and cell
membrane protrusions\index{membrane protrusions} are much larger than individual
adhesion bonds\index{adhesion bonds}, thus making
cellular adhesions a non-local process, and (3) the work of adhesion and
de-adhesion may be different, as cells can stabilize adhesion bonds after their
formation \parencite{Harris1976}. To resolve these issues, Harris proposed the
\textit{Differential Surface Contraction Hypothesis\index{Differential Surface
Contraction Hypothesis}}, arguing the contractile strength of a membrane
completely describes its surface tension. A model similar
to this idea was later implemented in a successful vertex model of cell-sorting
in epithelial tissues \parencites{Brodland2000}{Chen2000}{Brodland2002}.
\medskip

Cellular adhesion is facilitated by a wide range of cell surface molecules and
cell-cell junctions~\cite{Friedl2017}. {\it Adherens junctions\index{adherens
junctions}\/} connect the actomyosin\index{actomyosin} networks of different
cells with each other, they form strong bonds and initiate cell
polarization\index{cell polarization}; {\it tight junctions\index{tight
junctions}\/} are the strongest cell-cell connections and they are used to create
impermeable physical barriers; {\it gap junctions\index{gap junctions}\/} are
intracellular ion channels\index{ion channels} that allow cell-cell
communication\index{cell-cell communication}; {\it integrins\index{integrins}\/}
connect the cell cytoskeleton with the ECM;\ plus a selection of further
adhesion molecules such as {\it cadherins\index{cadherins},
igCAMS\index{igCAMS}, Slit/Robo\index{Slit/Robo}, Ephrin/Eph\index{Ephrin/Eph}},
etc.,~\cite{Friedl2017}.

The transition between tightly packed cells and free moving cells is fluent and
it depends to a large extend on the adhesion processes that are involved. For
example myoblasts\index{myoblasts} and myofibres\index{myofibres} are so tightly
connected that they can exert physical forces (muscles\index{muscles});
epithelial sheets\index{epithelial sheets} form the lining of many organs and
vessels, where it is important to separate an ``inside'' from an ``outside''; in
angiogenesis\index{angiogenesis}, endothelial\index{endothelial} cells move and
sprout new vessels, which requires loose adhesion as compared to mature vessels;
the movement of immune cells\index{immune cells} occurs in small cell clusters
or as individual cells, and it requires highly variable adhesion properties of
the immune cells. Cancer cells\index{cancer cells} often loose their adhesive
properties, which results in local invasion and metastasis\index{metastasis}.
The form of cancer invasion is highly variable and many different types have
been identified, including cluster invasion\index{cluster invasion}, small cell
groups, cancer-immune cell clusters\index{cancer-immune cell clusters}
individual cells and network-type invasions\index{network-type
invasions}~\cite{Friedl2017}.

An important process in tissues is the {\it Epithelial Mesenchymal
Transition\index{Epithelial Mesenchymal Transition}\/} (EMT), where stationary
epithelial cells\index{epithelial cells} loose their adhesive properties and
become invasive, mesenchymal-like cells\index{mesenchymal-like cells}
\parencites{Nieto2011}{Graves2016}. The EMT requires a combination of mechanisms
such as changed cytoskeletal dynamics and changed adhesive properties. The EMT
is a hallmark of cancer\index{hallmark of cancer} metastasis
\parencites{Hanahan2000}{Hanahan2011} and it is very important to understand the
influence of adhesion on the EMT.\@

A detailed mathematical analysis of adhesion models can contribute to our
understanding of this important process, and explain, complement, and enrich
biological and experimental observations.

\section{Prior Modelling of Cellular Adhesions}
Cell adhesions are forces that act on the cell membrane. For a mathematical or
computational model, these forces need to be computed and balanced. Hence an
individual cell modelling approach seems to be a natural way to start.  Indeed,
the first model for cell sorting through differential adhesion was a model of
Graner and Glazier from 1991. They used a Cellular Potts model\index{Cellular
Potts model} approach, where
individual cells are represented as collection of lattice sites in a two
dimensional lattice
\parencites{Graner1992}{Graner1993a}{Graner1993}{Glazier1993}.  Cellular
adhesions are implemented as interface energies of cells that are touching at a
common interface.  Since in this model a single cell contains many lattices
sites the adhesion is a non-local interaction.  The deformation and movement of
the cells is described as an energy minimization\index{energy minimization} approach. Interface energies
are balanced with cell volume and cell shape energies plus a random component
due to noise. At each step random changes in the lattice configuration are
proposed and accepted by a Boltzman like function. Over the years, these
Cellular Potts models have become very successful modelling tools and they have
been widely used in applications \parencite{Scianna2013}. For example Turner
\etal\ \parencite{Turner2002} used a Cellular Potts model to study the effect of
adhesion at the invasion front of a tumour. They observed the formation of
clusters of invasion, and the formation of \textquote{fingering\index{fingering}} invasion
fronts.  In \parencite{Turner2004} they attempted to scale the Cellular Potts
model to a partial differential equation. However, the obtained macroscopic
equations are notoriously difficult to analyse.

In 1996, Byrne \etal\ \parencite{Byrne1996} studied the growth of avascular
tumour spheroids\index{tumour spheroids} in the presence of an external nutrient. The tumour growth is
determined by the balance between proliferative pressure and cell-cell adhesion,
which keep the spheroid compact. The DAH of Steinberg of surface tension\index{surface tension} on cell
clusters is implemented by the Gibbs-Thompson relation\index{Gibbs-Thompson
relation}, which relates the tumour
spheroid's curvature to the external nutrient concentration. It is assumed that
cell-cell adhesion is the force that maintains this curvature\index{curvature}
\parencite{Byrne1996}. Later, this model was modified such that the cell's
proliferation rate depended on the total pressure acting on the cell (due to
adhesion and repulsive forces) \parencite{Byrne2009}. This model was then
successfully compared to a cell-based model of tumour spheroid growth
\parencite{Byrne2009}. In a similar model, Perumpanani \etal\
\parencite{Perumpanani1996}  introduced a density depended diffusion term in a
tumour spheroid model, the idea was that cells in high density areas are slowed
down by the presence of adhesion bonds to neighbours. Since then, this approach
has been used in more complicated models of tumour growth (see
\parencites{Macklin2009}{Lowengrub2010}).  The adhesive mechanism in these
models are purely local. Further, none of these models was able to reproduce
cellular aggregations nor cell sorting commonly linked to adhesive interactions.

Differently to the cellular Potts model, Palsson \etal\ \parencite{Palsson2000}
used a lattice-free model\index{lattice-free model}, resolving the individual
physical forces between the cells using the theory of elasticity. Cells are
represented as deformable ellipsoids\index{deformable ellipsoids} with long and
short range interactions with other cells. This model is a non-local individual
based model and Palsson \etal\ used it to describe chemotaxis\index{chemotaxis}
and slug formation\index{slug formation} in {\it Dictyostelium
discoideum\index{Dictyostelium discoideum}\/}
\parencite{Palsson2000}. In a similar approach cells are modelled as elastic
isotropic spheres, where adhesive and
repulsion forces between adhering elastic spheres are resolved using a modified
Hertz model\index{Hertz model} \parencites{Schluter2012}{Schlu2015} or the
Johnson, Kendall, and Roberts model\index{Johnson, Kendall, and Roberts model}
\parencites{Johnson1971}{Drasdo2005}{Hoehme2010a}. Since these
interactions act over a wide range of cell separations, they are non-local
models.

Brodland \etal\ used a vertex model\index{vertex model} (an individual-based
model) to model cell-sorting in epithelial tissues
\parencites{Brodland2000}{Chen2000}.  Similar
to Steinberg's assumptions,  the model considers surface tension at cell-cell
interfaces. The surface tension in their model depended on the forces of
adhesion, membrane contraction, and circumferential microfilament bundles
\parencite{Brodland2000}. They summarized the findings of their numerical
studies by formulating the \textit{Differential Interfacial Tension
Hypothesis\index{Differential Interfacial Tension Hypothesis}}
of cell-sorting \parencite{Brodland2002}. Once again this was a non-local
description of cell-adhesion.

Since up to this point all cellular adhesion models capable of explaining
aggregations and cell sorting were based on non-local models in cell-based
approaches, Anderson
proposed to combine the continuum and cell-based approaches in a hybrid
model\index{hybrid model}
\parencite{Anderson2005}. The significance of this hybrid approach is that cells
are individually represented (adhesion effects can be taken into account) and
environmental factors such as diffusing proteins and chemokines can be modelled using
well-established reaction diffusion equations. This approach has been popular in
studying the dynamics of tumour spheroids
\parencites{Ramis-Conde2008}{Ramis-Conde2008a}.

\begin{figure}[!ht]\centering
    \includegraphics[width=8cm]{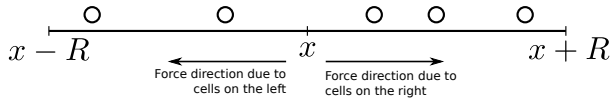}
    \caption[Intuitive explanation of the non-local term of cell adhesion]
    {Intuitive description of the non-local adhesion term. Two cells are pulling
    to the left and three cells are pulling to the right, hence the net force is
    to the right (assuming that all cells generate the same force).
    }\label{Fig:NonLocalTermExplanation}
\end{figure}

In 2006, Armstrong \etal\ proposed the first continuum model of cellular
adhesions capable of explaining adhesion driven cell aggregations and cell sorting
\parencite{Armstrong2006}.  Since this model is the focus of our monograph, we
represent it in its basic and one-dimensional form here. Let $u(x, t)$ denote
the density of a cell population at spatial location $x$ and time $t$. Then its
evolution subject to random motility and cell-cell adhesion is given by the
following non-local integro-partial differential equation\index{non-local
integro-partial differential equation}
\begin{equation}\label{Eqn:ArmstrongModelIntro}
    \frac{\partial}{\partial t} u(x, t) = \underbrace{D \frac{\partial^2}{\partial x^2}
        u(x, t)}_{\text{random motility}}
                - \underbrace{\alpha \frac{\partial}{\partial x}  \lb u(x, t)
                \int_{-R}^{R} h(u(x + r, t)) \Omega(r) \dd r
                \rb}_{\text{cell-cell adhesion}},
\end{equation}
where $D$ is the diffusion coefficient\index{diffusion coefficient}, $\alpha$
the strength of the homotypic adhesion\index{homotypic adhesion}, $h(u)$ is a
possibly nonlinear function describing the nature of the adhesive force,
$\Omega(r)$ an odd function, and $R$ the sensing radius of the
cell. We give a detailed description of this model and the biological meaning of
the terms in \cref{c:randomwalk}, see also \parencite{Buttenschon2017} for a
derivation from a stochastic process.  An intuitive explanation of the non-local
cell-cell adhesion term in equation~\eqref{Eqn:ArmstrongModelIntro} is given in
Figure~\ref{Fig:NonLocalTermExplanation}. The non-local term represents a
tug-of-war\index{tug-of-war} of the cells on the right and the cells on the
left, with the cell at $x$ moving in the direction of largest force. The effect
is that cells move up
\textquote{non-local} gradients of cell population and thus arises the
possibility for formation of cell aggregates. The two-population version of
equation~\eqref{Eqn:ArmstrongModelIntro} was the first continuum model that
correctly replicated cell-sorting experiments \parencite{Armstrong2006}.

As discussed above, cellular adhesions feature prominently in organism
development\index{organism development}, wound healing\index{wound healing} and
cancer invasion\index{cancer invasion} (metastasis). Therefore, it is
unsurprising that model~\eqref{Eqn:ArmstrongModelIntro} has found extensive use
in modelling cancer cell invasion
\parencites{Gerisch2008}{Sherratt2008}{Gerisch2010}{Painter2010}{Chaplain2011}{Andasari2012a},
and developmental processes \parencite{Armstrong2009}. More recently,
spatio-temporal variations of the adhesion strengths \parencite{Domschke2014},
and adhesion strength variations due to signalling proteins
\parencite{Bitsouni2017a} were considered.

The non-local model~\eqref{Eqn:ArmstrongModelIntro} has also been criticized for
oversimplification, namely for its use of a simple diffusion term
\parencite{Murakawa2015}. Supported by experimental data, Murakawa \etal\
\parencite{Murakawa2015}, noticed that under certain conditions
equation~\eqref{Eqn:ArmstrongModelIntro} gave unrealistic solutions. To address
this shortcoming, Murakawa \etal\ modified the modelling assumption
\textquote{cells move randomly} to \textquote{cells move from high pressure to
low pressure regions}. For this reason, they introduced a density-dependent
diffusion\index{density-dependent diffusion} term of porous medium type
\parencite{Murakawa2015}.

Prerequisite for the extensive numerical exploration of the solutions of the
non-local equation~\eqref{Eqn:ArmstrongModelIntro} was the development of
efficient numerical methods to evaluate its integral term. An efficient method
based on the fast Fourier transform\index{fast Fourier transform} was developed
in \parencite{Gerisch2010a}.  Using this efficient algorithm, numerical
solutions of equation~\eqref{Eqn:ArmstrongModelIntro} are implemented using a
spatial finite-volume discretization\index{spatial finite-volume discretization}
and the method of lines\index{method of lines} for the temporal advancement
\parencite{Gerisch2001}.

Existence results for the solutions of the non-local
equation~\eqref{Eqn:ArmstrongModelIntro} and generalizations in any space
dimension were developed in
\parencites{Sherratt2008}{Andasari2012a}{Hillen2017}{Bitsouni2017a}. Most
significant is the general work \parencite{Hillen2017}, who showed local and
global existence of classical solutions. For the case of small adhesion
strength, travelling wave solutions of equation~\eqref{Eqn:ArmstrongModelIntro}
were found and studied in \parencite{Chunhua2013}.

\section{Non-local Partial Differential Equation Models}
Having introduced the non-local adhesion model~\eqref{Eqn:ArmstrongModelIntro},
we take a look at developments of non-local partial differential equation models
in general.  Non-local models basically arise in two different situations:
Firstly, one assumes a priori the existence of non-local interactions, or,
secondly, a non-local term arises after solving a partial differential equation
or taking some asymptotic limit of it \parencites{Mogilner1995}{Iron2000}{Knutsdottir2014}.
For example, in a system of four partial differential equations (two cell
populations and two diffusing signalling molecules), Knutsdottir \etal\
\parencite{Knutsdottir2014} applied a quasi-steady state
assumption\index{quasi-steady state assumption} to the
signalling molecules, which diffuse much faster than the cellular populations.
The resulting elliptic equations were solved using Green's
functions\index{Green's function}, and thus
non-local terms were introduced in the remaining cell density equations. Similar
methods are standard in the analysis of chemotaxis\index{chemotaxis} (see
Horstmann~\cite{Horstmann2003}).

Non-local models arise in the study of Levi walks\index{Levi walks}. Levi walks
have been used to describe anomalous diffusion\index{anomalous diffusion}, for
example in situations where particles show distinct searching and foraging
behaviour. A continuum description of these random walks leads to fractional
Laplacian\index{fractional Laplacian} operators, which are non-local
operators~\cite{Taylor-King2016, Klages}.

In discrete time, integro-difference equations\index{integro-difference
equations} are a common tool in ecological modelling, in particular in
situations where distinct generations can be identified, and only certain parts
of the population spread spatially \parencite{Neubert2000,Kot}. One example is
seed dispersal\index{seed dispersal} of seeds released annually during a certain
season.

Nonlinear integro-partial differential equations have also been derived for
birth-jump processes\index{birth-jump process}~\parencite{Hillen2014}. These are
continuous in time and space and new-born particles are allowed to spread
non-locally. Birth-jump processes have been applied to cancer
spread~\cite{HillenEnderling,Borsi2015,Delgado2017}, to forest fire spread
\parencite{Martin2016}, and to the modelling of evolution and
selection~\cite{Chisholm2015}.

The non-local cell-cell adhesion model~\eqref{Eqn:ArmstrongModelIntro} has the
non-local term in the advection term. The first non-local equations with the
non-local term contained in the advection term were proposed in a series of
papers by Nagai~\etal\ in 1983 \parencites{Nagai1983}{Nagai1983b}{Nagai1983c}.
Their introduction of the non-local term was driven by a desire to model
aggregation processes in ecological systems. For comparison with the non-local
adhesion model~\eqref{Eqn:ArmstrongModelIntro}, the equations of Nagai \etal\
looked like this
\begin{equation}\label{Eqn:NagaiNonLocalEquation}
    u_t = \lb u^m \rb_{xx} - \lsb u \lcb
        \int_{-\infty}^{\infty} K(x - y)u(y,t) \dd y \rcb \rsb_{x},\
        x \in \R,\ t > 0,
\end{equation}
where $m > 1$.  Shortly after, Alt studied generalizations of
equation~\eqref{Eqn:NagaiNonLocalEquation} in \parencites{Alt1985}{Alt1985b}. A
version of equation~\eqref{Eqn:NagaiNonLocalEquation} with finite integration
limits and with the special choice of $K(x - y) = \sgn(x - y)$ was studied by
Ikeda in \parencites{Ikeda1985}{Ikeda1987}. Ikeda established the existence of
solutions of equation~\eqref{Eqn:NagaiNonLocalEquation} on an unbounded domain
and developed spectral results in the special case of $m = 2$, that was used to
give a classification of the steady-states of
equation~\eqref{Eqn:NagaiNonLocalEquation}.  In 1999, Mogilner \etal\
\parencite{Mogilner1999} used such non-local models to develop evolution
equations describing swarms\index{swarms}. More recently, such models were used to
describe the aggregation of plankton \parencite{Adioui2005}, and to model animal
populations featuring long-ranged social attractions and short-ranged dispersal
\parencite{Topaz2006}. Eftimie \etal\ \parencites{Eftimie2007a}{Eftimie2007}
used a Lagrangian formulation to obtain non-local hyperbolic models of
communicating individuals. Since then, equivariant bifurcation theory was used
to study the possible steady states of such communication models
\parencite{Buono2015}. Very recently they discussed the use of Lyapunov-Schmidt
and centre-manifold reduction to study the long-time dynamics of such equations.
The non-local adhesion model falls also in this category of equation
\parencite{Armstrong2006} and more recently it was generalized to include both
aggregations and repulsive behaviour \parencite{Painter2015}.

Another famous non-local model for biological applications is the
{\it aggregation equation\index{aggregation equation}\/}
\begin{equation}\label{aggr}
    u_t = \nabla \cdot (u \nabla (u + W * u )),
\end{equation}
where $W(x)$ is a species interaction kernel that acts non-locally through
convolution~$*$.  The interaction kernel $W$ can be used to describe species
aggregation, repulsion and alignment\index{alignment}~\cite{Topaz2006}. The model has been
successfully employed to describe animal swarming as well as non-biological
applications such as granular matter\index{granular matter},
astrophysics\index{astrophysics},  semiconductors\index{semiconductors}, and opinion
formation\index{opinion formation} (see references in~\cite{Burger2018}). The aggregation
equation~\eqref{aggr} can be seen as gradient flow of an appropriate energy functional,
which contains the non-local interaction potential $W$. This makes the vast
resource of variational methods available to the analysis of~\eqref{aggr}. Many
results are known about aggregations, pattern formation, finite time blow-up,
and the relation of~\eqref{aggr} to chemotaxis
models~\cite{Ducrot2014,Burger2018,Bernoff2013}. As we show in \cref{appendix:adh_pot},
the adhesion model~\eqref{Eqn:ArmstrongModelIntro} studied here can be
formulated as an aggregation
equation with linear diffusion term, and the rich theory of the aggregation
equations becomes available. However, this does not help us in the current study
of steady states and bifurcations, which is the focus of this monograph.

A non-local model for chemotaxis was proposed in~\cite{Othmer2002,Hillen2007} to
replace the chemical gradient in a chemotaxis model by the following non-local
gradient\index{non-local gradient}
\begin{equation}\label{Eqn:NonLocalGradientChemotaxisIntro}
    \nonlocalgrad{R} v(x) \coloneqq \frac{n}{\omega_{D}(x) R}
        \int_{\mathbb{S}^{n-1}_{D}} \sigma v(x + R \sigma) \dd \sigma,
\end{equation}
$\mathbb{S}^{n-1}_{D}(x) = \{ \sigma \in \mathbb{S}^{n-1} : x + \sigma R
\in D \}$, and $\omega_{D}(x) = \abs{\mathbb{S}^{n-1}_{D}(x)}$. It was shown in
\parencite{Hillen2007} that this non-local gradient leads to globally bounded
solutions of the non-local chemotaxis\index{chemotaxis} equation in cases where
the local gradient gives rise to blow-up solutions. Extending the work on steady
states of the
local chemotaxis equation \parencites{Schaaf1985}{Wang2013}, Xiang
\etal~\cite{Xiang2013} used bifurcation techniques to analyse the steady states
of the non-local chemotaxis equation in one dimension.

Non-local terms have also been considered in reaction terms, \ie, in equations of
the form
\[
    u_t = u_{xx} + f(x, u, \bar{u}),
\]
where
\[
    \bar{u} = \int g(x, u) \dd x.
\]
These types of equations were studied in
\cites{Davidson2007a,Davidson2006c,Davidson2006a,Freitas2000a,Freitas1994,
Brandolini2011,Freitas1994a}. Oftentimes such systems occur as the asymptotic
limit (\textquote{shadow system}\index{shadow system}) of a system of
reaction-diffusion equations.  An application of such a model is for example to
Ohmic heating\index{Ohmic heating} in \parencite{Lopez-Gomez1998}.  A different
application considers the growth of phytoplankton\index{phytoplankton} in the
presence of light and nutrients \parencite{Zhang2011}. Yet another application
of such a model is a study on the effect of crop raiding\index{crop raiding} of
large bodied mammals \parencite{Mei2014}.

Many of the above mentioned examples are formulated on infinite domains, to
avoid subtleties in dealing with boundary conditions, or subtleties in ensuring
the non-local term is well-defined. Indeed, even the analysis of local equations
such as the viscous Burger equation on bounded domains have remained unaddressed
until recently \parencite{Watanabe2016}.

\section{Outline of the Main Results}

In this monograph, we consider non-local models of cell-cell adhesion in the
form of an integro-partial differential equations (see
equation~\eqref{Eqn:ArmstrongModelIntro}). The central problems which we would
like to address in this monograph are:

\begin{enumerate}
  \item What are the non-trivial steady states of the non-local cell adhesion
    model  (\ie\ equation~\eqref{Eqn:ArmstrongModelIntro}) in a periodic domain?
  \item What is their bifurcation structure and their stabilities?
  \item How to include biologically realistic boundary conditions that can
    describe no-flux boundaries, boundary adhesion, or boundary repulsion?
  \item What are the non-trivial steady states if boundary conditions are included?
\end{enumerate}

In \cref{Chapter:Preliminaries}, we give a brief summary of the derivative of
non-local cell-cell adhesion models from an underlying stochastic random walk,
and summarize the key bifurcation results we employ.

In \cref{chapter:BasicProperties}, we define the non-local term
in~\eqref{Eqn:ArmstrongModelIntro} as integral operator $\K[u]$ and explore its
mathematical properties. In particular, we establish its continuity properties,
$L^p$-estimates, compactness, and spectral results. These results are key to the
subsequent bifurcation analysis. In fact we show that $\K[u]$ is a
generalization of the classical first-order derivative. Finally, we establish
elementary results on the steady states of the non-local adhesion
model~\eqref{Eqn:ArmstrongModelIntro}, including an a-priori estimate.

In Chapters~\ref{Chapter:LocalBifPeriodic} to~\ref{Chapter:GlobalBifPeriodic},
we investigate the steady-states of a single population non-local model of
cellular adhesions on a periodic domain. We combine global bifurcation
techniques pioneered by Rabinowitz, equivariant bifurcation theory (the equation
is $\O2$-equivariant), and the mathematical properties of the non-local adhesion
term, to obtain the existence of unbounded global bifurcation branches of
non-trivial solutions. In words, our main theorem is:

\begin{theorem}
    For each $\alpha > \alpha_n$ ($\alpha_n$ eigenvalues of the
    \textbf{linearised} problem) the periodic non-local linear adhesion steady state
    problem~\eqref{Eqn:ArmstrongModelIntro} (i.e.\ $h(u) = u$), has at least
    two $n$-spiked solutions (one in each of $\mathcal{S}_{n}^{\pm}$).
\end{theorem}

The solution branches are classified by the derivative's number of zeros (\ie\
the number of extrema remains fixed along a bifurcation branch), or in other
words by their number of spikes. The significance of this result is that it
parallels the seminal classification of solutions of nonlinear Sturm-Liouville
problems (Crandall \& Rabinowitz, 1970) and the classification for equivariant
nonlinear elliptic equations (Healey \& Kielh\"{o}fer, 1991).

In Chapter~\ref{Chapter:AdhNoFlux} we consider the construction of
no-flux boundary conditions for the non-local cell adhesion model, and finally
we explore what we can say about its steady states.  In the past, boundary
conditions for non-local equations were avoided, because their construction
is subtle and requires biological insight. Using the insights from
\parencite{Buttenschon2017}, we construct a non-local operator, which takes
various boundary effects into account, such as no-flux, or boundary adhesion, or
boundary repulsion. These boundary conditions destroy some of the symmetry
properties that were available in the periodic case, and consequently, the
Rabinowitz-bifurcation theory does no longer apply. We find similar steady
states as in the periodic case, but the full bifurcation structure in this case
requires new methods. Numerically, we find that steady states show boundary
behavior, which is well known from fluids that are wetting the boundary
(adhesive), or are repulsed from the boundary.

%% file: chapterPreliminaries.tex
%
\chapter{Preliminaries}\label{Chapter:Preliminaries}
In this section we present some basic results that are needed later. We give a
summary of the derivation of the non-local adhesion model from biological
principles as presented in~\cite{Buttenschon2017}, we introduce some notations
and methods from abstract bifurcation theory as it was developed
in~\cite{Lopez-Gomez2007,Lopez-Gomez2016}, we introduce an averaging operator on
periodic domains, and we cite a global existence results
for~\eqref{Eqn:ArmstrongModelIntro} in $\R^n$ from~\cite{Hillen2017}.
\section{Biological derivation\index{biological derivation} of the non-local adhesion model}\label{c:randomwalk}
The goal of this section is to summarize the derivation of the non-local
adhesion model given in equation~\eqref{Eqn:ArmstrongModelIntro} from a
stochastic position jump process\index{position-jump process}.
The derivation of population level models from
underlying mesoscopic movement models\index{mesoscopic movement model} has a rich
history~\cite{Calvo2015,Hillen2009,okubo2001,Stevens1997,Othmer1988}.

In order to focus on the modelling of the cell-cell interactions in the absence
of boundary effects, we carry out the derivation on an unbounded domain \ie\
$\Omega = \R^n$. Note that in \cref{Chapter:AdhNoFlux}, we discuss
extensions to include boundary effects. Key to the derivation is that
cells send out many membrane protrusions\index{membrane protrusions} to sample their environment.
We assume that these membrane protrusions are more frequent then
translocations of the cell body; the cell body includes the cell nucleus and
most of the cell's mass. In other words, we are most interested in
translocations of the cell body, and not the frequent, but temporary, shifts due
to membrane protrusions. For this reasons, we define the population density
function $u(x, t)$ as follows
\[
    u(x, t) \equiv\mbox{Density of cells with their cell body centred at $x$ at
    time $t$}.
\]
We make two assumptions about the movement of the cells:
\begin{description}
    \item[\bfseries Modelling Assumption 1]
        We assume that in the absence of spatial or temporal heterogeneity
        the movement of individual cells can be described by Brownian
        motion\index{Brownian motion}. It has been shown that
        this is a reasonable assumption for many cell types
        \parencites{Schienbein1994}{Mombach1995}{Beysens2000}.
    \item[\bfseries Modelling Assumption 2]
        The cells' polarization\index{polarization} may be influenced by spatial or temporal
        heterogeneity. We denote the polarization vector\index{polarization vector} by $\vec{p}(x)$.
\end{description}
We describe the evolution of the cell density $u(x, t)$ using a position-jump
process\index{position-jump process}, using stochastically independent jumpers
(\ie\ a continuous time random walk, for which we assume the independence of the
waiting time distribution\index{waiting time distribution} and spatial redistribution).
Here, the waiting time distribution is taken to be the exponential
distribution, with constant mean waiting time. Let $T(x, y)$ denote the rate for
a jump from $y\to x$ with $x, y \in \Omega$. The evolution of $u(x, t)$ is given by
the continuous-time master equation\index{master equation}.
\begin{equation}\label{Eqn:MasterEquation}
    \frac{\partial u}{\partial t}(x,t) = \lambda \int_{D} \lsb T(x, y) u(y, t)
        - T(y, x) u(x, t)\rsb \dd \mu(y),
\end{equation}
where $(D, \mu)$ is a measure space of a physical space (\eg\ a domain or a lattice),
$\lambda$ jump rate, and $T(x,y)$ probability density for jump from $y$ to $x$.
For more details on the derivation of the master equation
see~\cite{Othmer1988,hughes1995random,vanKampen2011}.

For notational convenience we associate to a jump from $y \in \Omega$ to
$x \in \Omega$ the {\it heading\index{heading}\/} $z\coloneqq x - y$. Using the
heading we define $T_y(z) \coloneqq T(y+z,y) = T(x,y)$. Let $D^y$ denote the set
of all possible
headings from $y$. We assume that the set $D^y$ is symmetric (\ie\ if $z \in
D^y$, then so is $-z \in D^y$). We further assume that for every $y \in \Omega$,
the function $T_y$ is non-negative as it represents a rate.

Given $y \in \Omega$ we denote the redistribution kernel at this location
by $T_y(z)$; we assume that $T_y \in L^1(D^y)$ and that $\norm{T_y}_1 = 1$ holds.
This turns $T_y$ into a probability density function (pdf) on $D^y$.

Any function which is defined for both $z$ and $-z$ can be decomposed into even
and odd components\index{even-odd decomposition}, which are denoted by $S_y$ and
$A_y$ respectively.
\begin{lemma}[Lemma~1 in~\cite{Buttenschon2017}]\label{lemma:eins}
    %
    %
    Consider $y \in \Omega$, given $T_y \in L^1(D^y)$, then there exists
    a decomposition as
    \begin{equation}
        T_y(z) =
        \begin{cases}
            S_y(z) + A_y(z) \cdot \frac{z}{\abs{z}} &\mbox{if } z \neq 0\\
            S_y(z)                                  &\mbox{if } z = 0
        \end{cases}
    \end{equation}
    with $S_y \in L^1(D^y)$ and $A_y \in \lb L^1(D^y) \rb^n$.
    The even and odd parts are symmetric such that
    \begin{equation}
        S_y(z) = S_y(-z) \quad \mbox{and} \quad A_y(z) = A_y(-z).
    \end{equation}
\end{lemma}

Using this decomposition we define two properties which are analogous to
Modelling Assumptions 1 and 2 above. First, we define the motility.

\begin{definition}[Motility]\label{defn:motility}
    We define the motility\index{motility} at $y \in \Omega$ as
    \begin{equation}
        M(y) \coloneqq \int_{D^y \setminus \{ 0 \}} T_y(z) \dd z
             = \int_{D^y \setminus \{ 0 \}} S_y(z) \dd z,
    \end{equation}
    where the integration is \wrt\ the measure on $D^y$.
\end{definition}

The motility is the probability of leaving $y$.
This probability is 1 if $0 \notin D^y$; it is also 1 if $0 \in D^y$
and $T_y$ is a continuous pdf, and it may be smaller than 1 if
$0 \in D^y$ and $T_y$ is a discrete pdf. Here we find that the motility depends
solely on the even component $S_y$, in other words solely on modelling
assumption one.

Secondly, we define the polarization vector in a \textit{space-jump} process.

\begin{definition}[Polarization Vector]\label{defn:polarization}
    The polarization vector\index{polarization vector} at $y \in \Omega$ is defined as,
    \begin{equation}
        E(y) \coloneqq \frac{\int_{D^y} z \, T_y(z) \dd z}
                         {\norm{\int_{D^y} z \, T_y(z) \dd z}}
             = \frac{\int_{D^y} z \, A_y(z) \cdot \frac{z}{\abs{z}} \dd z}
                    {\norm{\int_{D^y} z \, A_y(z) \cdot \frac{z}{\abs{z}} \dd z}}.
    \end{equation}
    where the integration is \wrt\ measure on $D^y$.
\end{definition}

The first moment of the pdf $T_y$ can be intuitively understood as the expected
heading\index{expected heading} of a jump originating at $y$. This is in direct
correspondence with a polarized cell which, following polarization, moves in the
direction of the
polarization vector.  The expected heading is solely determined by $A_y$, which
therefore plays the role of the polarization vector $\vec{p}(y)$ in a
space-jump process.  This correspondence motivates us to set $A_y =
\vec{p}(y)$ in the subsequent derivations.\\

We can now derive the macroscopic limit of the master
equation~\eqref{Eqn:MasterEquation}. In addition, to the prior assumptions we
assume that we have a {\it myopic\/} random walk\index{myopic random walk} such
that the jump-probability only depends on the jump location $y$ but not on the
target location $y+z$. We
only consider small jumps of length $h\ll 1$ and we
expand~\eqref{Eqn:MasterEquation} in $h$. For the mathematical details we refer
the reader to~\cite{Buttenschon2017}. At the end of this process we obtain.
\begin{equation}\label{spliteq3}
    u_t(x,t) \approx \frac{\lambda h^{n+1}}{2n}\abs{\Sn}
        \Delta (S_x u(x,t)) - \frac{\lambda h^n}{n}\abs{\Sn}
        \nabla\cdot \lb A_x u(x,t)\rb.
\end{equation}
We then assume that $\vec{A}_x \sim \BigOh{h}$, and let
\[
    D(x) =
        \lim_{\tiny \substack{h \to 0 \\ \lambda \to \infty}} \lambda h^2 S_x,
        \qquad
    \vec{\alpha}(x) =
        \lim_{\tiny \substack{h \to 0 \\ \lambda \to \infty}} 2\lambda
        h^2 \vec{A}_x.
\]

The advection-diffusion limit is:
\begin{equation}\label{Eqn:AdvectionDiffusionLimit}
    \frac{\partial u}{\partial t}(x, t) + \nabla \cdot \lb
        \vec{\alpha}(x) u(x,t) \rb = \Delta \lb D(x) u(x, t) \rb.
\end{equation}
We notice that the diffusive part is defined through the symmetric component
$S_x$ of $T_x$, while the advective part $\vec{\alpha}$ is defined through the
anti-symmetric component $\vec{A}_x$. We identify the anti-symmetric part with
the polarization vector $\vec{A}_x\sim \vec{p}$ and it remains to find a good
model for the polarization vector $\vec{p}$. Here we assume that
the cell's polarization vector is determined by the interactions of adhesion
molecules on cell's protrusions and adhesion molecules present in the
surrounding environment (\eg\ located on the surfaces of other cells).

\begin{figure}
    \begin{minipage}[t]{0.45\textwidth}
        \includegraphics[width=0.9\textwidth]{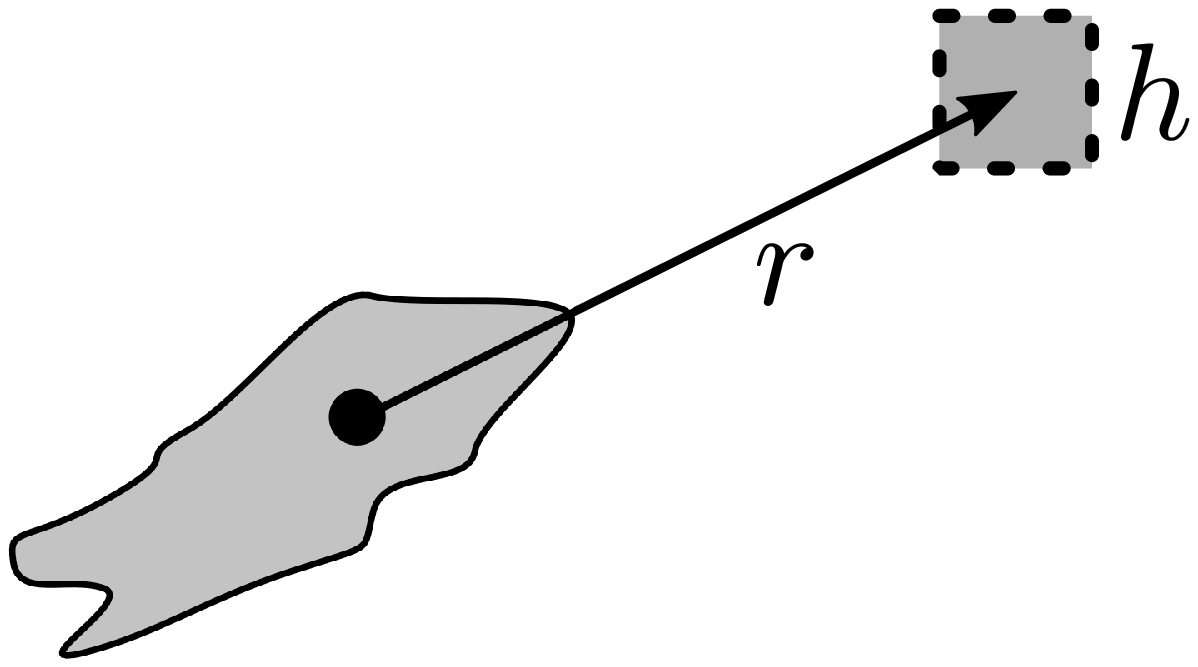}
    \end{minipage}\hfill
    \begin{minipage}[t]{0.45\textwidth}
        \includegraphics[width=0.9\textwidth]{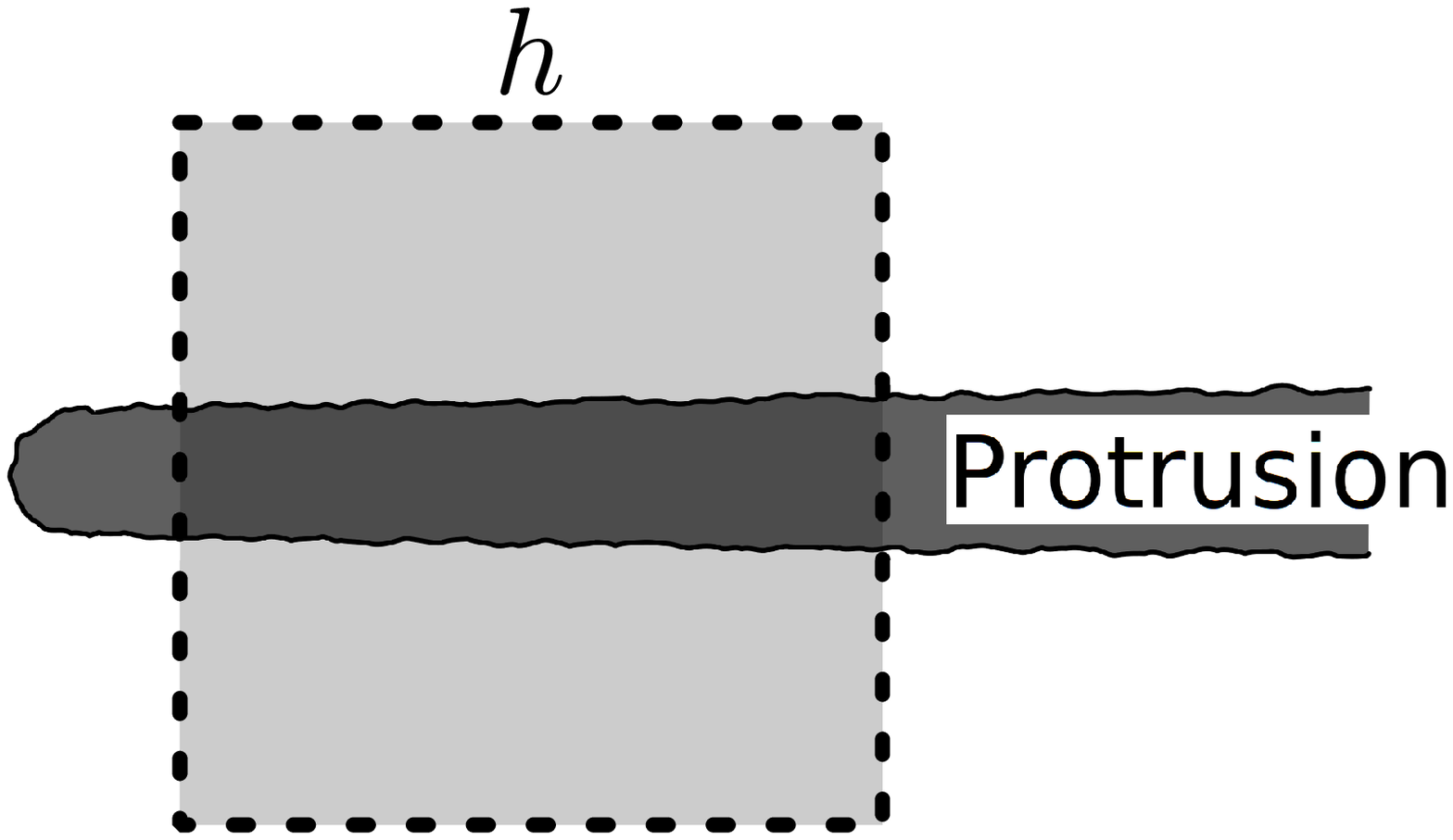}
    \end{minipage}
    \caption{Left: A cell of spatial extend, and a small test volume $V_h$ at a
    distance $r$ from the cell body. Right: A cell protrusion overlapping a
    small test volume $V_h$.}\label{Fig:SmallTestVolume}
\end{figure}

To describe the microscopic interactions between a cell's protrusion and its
surroundings, we consider a small test volume $V_h$ (see
\cref{Fig:SmallTestVolume}). We assume that the factors determining the size of
contribution of the interactions in $V_h$ to the cell's overall polarization
are:

\begin{enumerate}
    \item The distance from $V_h$ to the cell body is $\abs{r}$.
    \item Direction of generated force is $r/\abs{r}$.
    \item The free space\index{free space} $f(r)$.
    \item The part of protrusion in $V_h$ is called $\omega(r)$.
    \item The density of formed adhesion bonds is $N_b$.
    \item The adhesive strength per bond is $\gamma$.
\end{enumerate}
The component of the cell polarization that is generated through forces in the test volume $V_h$ is then
\[
    \vec{p}_h(x + r) = \gamma
    \underbrace{h^n N_b(x + r)}_{\text{\#of adhesion bonds}} \qquad
    \underbrace{f(x + r)}_{\text{free space}} \qquad
    \underbrace{h \omega(r)}_{\text{amt.\ of cell in $V_h$}}\qquad \underbrace{\frac{r}{|r|}}_{\text{direction}}
\]
Summing over all test volumina,
we obtain the cell's polarization:
\begin{equation}\label{pnet}
    \vec{p}_{\text{net}}(x) = \int_{\mathbb{B}_R(0)} N_b(x + r) f(x + r)
    \vec{\Omega}(r) \dd r,
\end{equation}
where $\vec{\Omega}(r) = \omega(r)\frac{r}{|r|}$. Finally we let
$\vec{\alpha}(x) = \vec{p}_{\text{net}}(x)$, and
equation~\eqref{Eqn:AdvectionDiffusionLimit} becomes
\begin{equation}\label{Eqn:FinalAdhesion}
    \frac{\partial u(x,t)}{\partial t} = \frac{\partial}{\partial x} \lsb
    \frac{\partial}{\partial x} \lb D(x) u(x,t) \rb
        - \alpha(x) u(x,t) \vec{p}_{\text{net}}(x) \rsb,
\end{equation}
with the non-local term $\vec{p}_{\text{net}}$ from (\ref{pnet}).

We are now left with choosing the functions $f(\cdot)$, and the reaction kinetic
yielding the density of bound adhesion molecules $N_b(\cdot)$. One particular
choice leading to the original model by Armstrong \etal\ is:

\begin{enumerate}
    \item Let $\omega(r)$ be the uniform distribution \ie\ $\Omega(r) = \sgn(r) / 2R$.
    \item Assume there is always free space $f(x) \equiv 1$.
    \item Assume mass action kinetics for adhesion bonds, i.e. $h(u) = u$.
    \item Assume that the background adhesion bond density is proportional to the population density:  $N_b(x) \propto u(x)$.
\end{enumerate}
With those assumptions equation~\eqref{Eqn:FinalAdhesion} becomes
\[
    \frac{\partial u(x,t)}{\partial t} = \frac{\partial}{\partial x}
    \lsb \frac{\partial}{\partial x} \lb D(x) u(x,t) \rb - \alpha(x) u(x,t)
    \int_{-R}^{R} u(x+r, t) \; \Omega(r) \dd r \rsb.
\]
The general framework in~\eqref{Eqn:FinalAdhesion} can be used to incorporate
more detailed biological features such as receptor binding kinetics, competition
for space, spatial extent of background cells etc (see~\cite{Buttenschon2017}).
Also, extensions to higher space dimensions and to several cell types are
straight forward (see for example~\cite{Armstrong2006,Sherratt2008,
Andasari2012a,Hillen2017}).
\section{Introduction to Nonlinear Analysis}
This section is based on the abstract bifurcation theory\index{bifurcation
theory} as developed by J.~López-Gómez
in~\cite{Lopez-Gomez2007,Lopez-Gomez2016}, and is meant as a quick introduction
to the abstract framework which will be employed in this work.

Banach spaces\index{Banach space} and their subspaces will be denoted using
capital letters that is $X, Y, U, V$ and so on. Operators between function
spaces will be denoted using the calligraphic font for example $\Lap, \F, \K$.
The argument of an operator will be enclosed in square brackets. For instance
\[
    \Lap : X \mapsto Y, \quad \Lap[x] = y.
\]
This is to distinguish the action of the operator from a family of
operators\index{family of operators}. For
example a family of operators may be indexed using a real number $\lambda$.
Then, we have a map from $\R \to \Lin(X, Y)$
\[
    \lambda \mapsto \Lap(\lambda) : X \mapsto Y,
\]
and for each fixed $\lambda$ we may study $\Lap(\lambda)[x] = y$.
The kernel\index{kernel} and range\index{range} of an operator is denoted $\Null[\Lap]$ and $\Range[\Lap]$.

Spaces of operators are denoted using the fraktur font. The most important is
the space of continuous linear operators\index{continuous linear operators}
denoted $\Lin$ and the space of compact operators\index{compact operators}
$\Kom$.  The space of Fredholm operators\index{Fredholm operators} is denoted
$\Fred_{i}$ where $i$ denotes the index\index{index}.

Special subspaces such as a continuum of solutions or the solution set of an
operator equation are also denoted using the fraktur font. For example, $\Cont,
\Soln$.

The following sections are based on \parencites{Lopez-Gomez2007}{Lopez-Gomez2016}.
Let $U, V$ be two real Banach spaces. We denote the space of bounded linear
operators from $U$ to $V$ by $\Lin(U, V)$, and by $\Fred_0(U, V)$ the subset of
$\Lin(U, V)$ containing all Fredholm operators with index $0$. The set of all
isomorphisms between $U$ and $V$ is denoted $\Iso(U, V)$.
The operator $\Lap$ is said to be Fredholm whenever
\[
    \dim \Null[\Lap] < \infty, \quad \codim \Range[\Lap] < \infty,
\]
and recall that
\[
    \codim \Range[\Lap] = \dim V / \Range[\Lap].
\]
The index of a Fredholm operator is defined by
\[
    \ind [\Lap] = \dim \Null[\Lap] - \codim \Range[\Lap].
\]
Therefore if $\Lap \in \Fred_0$, then
\[
    \dim \Null[\Lap] = \codim \Range[\Lap] < \infty.
\]
The most important example of a Fredholm operator with index zero is the
following.
\begin{theorem}[Theorem~4.2 \parencite{gohberg1990}]
    Let $\K \in \Kom(U)$ (compact operators), then $\F = \I - \K$ is a Fredholm
    operator of index zero.
\end{theorem}
\begin{theorem}[Theorem~4.1 \parencite{gohberg1990}]\label{Thm:FredholmCompactPerturbation}
    Let $\Lap \in \Lin(U, V)$ be a Fredholm operator and let $\K \in \Kom(U, V)$
    (compact). Then $\F = \Lap + \K$ is a Fredholm operator with
    \[
	\ind(\Lap + \K) = \ind(\Lap).
    \]
\end{theorem}
When the dynamical system depends on parameters, then it is natural to talk
about operator families\index{operator families}. In \parencite{Lopez-Gomez2007,Lopez-Gomez2016} this
approach is used to define bifurcation values as values in a generalized
spectrum\index{generalized spectrum} of an operator family as defined below.
\begin{definition}
    Let $U, V$ be two Banach spaces over the field $\mathbb{K}$
    and $r \in \mathbb{N}$.
\begin{enumerate}
\item   An \textit{operator family} $\Lap(\Omega)$ of class $\C^{r}$ in $\Omega
    \subset \mathbb{K}$ from $U$ to $V$ is a map
    \[
        \Lap \in \C^{r}(\Omega, \Lin(U, V)).
    \]
    In our application we have $\Omega \subset \R$.
    \item 
    Let $\Lap \in \C(\Omega, \Lin(U, V))$ be an operator family, then
        the point $\lambda_0 \in \Omega$ is a \textit{singular
        value}\index{singular value} of $\Lap$ if
    \[
        \Lap_{0} \coloneqq \Lap(\lambda_0) \notin \Iso(U, V),
    \]
    and it is a \textit{generalized eigenvalue}\index{generalized eigenvalue} of $\Lap$ if
    \[
        \dim \Null[\Lap_{0}] \geq 1.
    \]
    A generalized eigenvalue $\lambda_0$ is simple whenever
    \[
        \dim \Null[\Lap_{0}] = 1.
    \]
\item 
    The set of all singular values of the operator family $\Lap$ is called the
    \textit{spectrum}\index{spectrum} and is defined by
    \[
        \Sigma = \Sigma(\Lap) =
            \lcb \lambda \in \Omega : \Lap(\lambda) \notin \Iso(U, V) \rcb.
    \]
    Similarly, the set of all generalized eigenvalues of the operator
    family $\Lap$ is defined by
    \[
        \textrm{Eig}(\Lap) =
            \lcb \lambda \in \Omega : \dim \Null[\Lap(\lambda)] \geq 1 \rcb.
    \]
    It is obvious that $\textrm{Eig}(\Lap) \subset \Sigma(\Lap)$.
\item
    The \textit{resolvent set}\index{resolvent set} of $\Lap$ is defined by
    \[
        \rho(\Lap) \coloneqq \Omega \setminus \Sigma.
    \]
    \end{enumerate}
\end{definition}
\begin{remark}
    Since $\Lap \in \C(\R, \Lin(U, V))$ and $\Iso(U, V)$ is an open subset of
    $\Lin(U, V)$, we have that $\rho(\Lap)$ is open and possible empty. Thus
    $\Sigma(\Lap)$ is closed.
\end{remark}
\begin{lemma}
    If $\Lap_0 \in \Fred_0(U, V)$ then $\lambda_0 \in \Sigma(\Lap)$ if and only
    if $\lambda_0 \in \mathrm{Eig}(\Lap)$.
\end{lemma}
\begin{proof}
    We only have to prove one direction. Let $\lambda_0 \in \Sigma(\Lap)$, then
    $\Lap_0 : U / N[\Lap_0] \mapsto R[\Lap_0]$ is an isomorphism by the
    open mapping theorem. As $\Lap_0$ is Fredholm with index zero we have that
    $\dim N[\Lap_0] < \infty$, hence $\lambda_0 \in \mathrm{Eig}(\Lap)$.
\end{proof}
Hence, if $\Lap(\Omega) \subset \Fred_0(U, V)$, then $\Sigma(\Lap) =
\mathrm{Eig}(\Lap)$.
\begin{remark}
    Note that the concept of a \textit{generalized eigenvalue} of an operator
    family $\Lap(\Omega)$, should not be confused with the classical notions of
    an eigenvalue and spectrum (denote $\sigma(\Lap(\lambda))$), which is only
    defined for fixed values of $\lambda \in \Omega$. The classical
    spectrum is defined as (see for instance~\cite[Chapter 7]{fabian2013})
    \[
        \sigma(T) \coloneqq \{ \lambda \in \mathbb{K} :
            \lambda \I - T \ \mbox{is not invertible} \}.
    \]
    Note that $\sigma(T)$ can be decomposed into the point, continuum and
    residual spectrum depending on the precise way the operator fails to be
    invertible. It is however, possible to
    recover these classical notions from this more general definition. Indeed,
    suppose that $T \in \Lin(U, V)$, and consider the operator family
    \[
        \Lap^T(\lambda) \coloneqq \lambda I_U - T,
    \]
    then
    \[
        \sigma(T) = \Sigma(\Lap^T).
    \]
\end{remark}
\section{Abstract Bifurcation Theory}\label{s:Bifurcation}
Let $U, V$ be two real Banach spaces, and suppose that we want to analyze the
structure of the solution set\index{solution set} of the nonlinear operator given by
\[
    \F[\lambda, u] = 0, \quad (\lambda, u) \in \R \times U,
\]
where
\[
    \F : \R \times U \mapsto V,
\]
is a continuous map satisfying the following requirements:
\begin{enumerate}[label=\textbf{(F\arabic*)},ref=(F\arabic*),leftmargin=*,labelindent=\parindent]
    \item\label{Assumption:1} For each $\lambda \in \R$, the map
        $\F[\lambda, \cdot]$ is of class
        $\C^1(U, V)$ and
        \[
            \D_u \F(\lambda, u) \in \Fred_0(U, V) \quad\mbox{for all } u \in U.
        \]
    \item\label{Assumption:2} $\D_u \F : \R \times U \mapsto \Lin(U, V)$ is continuous.
    \item\label{Assumption:3} $\F[\lambda, 0] = 0$ for all $\lambda \in \R$.
\end{enumerate}
\begin{definition}
    Assume $\F$ satisfies~\ref{Assumption:1}-\ref{Assumption:3}.
\begin{enumerate}
\item A \textit{component} $\Cont$ is a closed and connected subset of the set
    \[
        \Soln = \lcb (\lambda, u) \in \R \times U : \F[\lambda, u] = 0 \rcb,
    \]
    that is maximal with respect to inclusion.
\item
    As $(\lambda, 0)$ is a known zero, it is referred to as the \textit{trivial
        state}\index{trivial state}.
\item
    Given $\lambda_0 \in \R$ it is said that $(\lambda_0, 0)$ is a
    \textit{bifurcation point}\index{bifurcation point} of $\F = 0$ if there exists a sequence
    $(\lambda_n, u_n) \in \F^{-1}[0]$, with $u_n \neq 0$ for all $n \geq 1$,
    such that
    \[
        \lim_{n \to \infty} (\lambda_n, u_n) = (\lambda_0, 0).
    \]
    \end{enumerate}
\end{definition}
For every map $\F$ satisfying~\ref{Assumption:1}-\ref{Assumption:3}
we denote
\[
    \Lap(\lambda) = \D_u \F(\lambda, 0), \qquad \lambda \in \R.
\]
By property~\ref{Assumption:2} we have that $\Lap \in \C(\R, \Lin(U, V))$
and by~\ref{Assumption:1} we have $\Lap(\lambda) \in \Fred_0(U, V)$, thus
\[
    \Lap(\lambda) \in \Iso(U, V) \quad\iff\quad
        \dim \Null[\Lap(\lambda)] = 0.
\]
\begin{lemma}[\parencite{Lopez-Gomez2016}]\label{Lemma:BifPointsAreEigenvalues}
    Suppose $(\lambda_0, 0)$ is a bifurcation point of $\F = 0$. Then,
    $\lambda_0 \in \Sigma(\Lap)$ ($\Lap = \D_u \F(\lambda, 0)$).
\end{lemma}
\begin{theorem}[Local Bifurcation \parencite{Crandall1971}]\label{Thm:CrandallMain}
    Let $U, V$ be Banach spaces, $W$ a neighbourhood of $0$ in $U$
    and\index{local bifurcation}
    \[
        \F : (-1, 1) \times W \mapsto V,
    \]
    have the properties:
    \begin{enumerate}
        \item $\F[\lambda, 0] = 0$ for $\abs{\lambda} < 1$,
        \item The partial derivatives $\D_\lambda \F, \D_u \F, \D_{\lambda u} \F$ exist and are
            continuous,
        \item $\Null[\D_u \F(0,0)]$ and $V / \Range[\D_u \F(0,0)]$ are one-dimensional, and we write
            \[
                \Null[\D_u \F(0,0)] = \vspan \seq{u_0}.
            \]
        \item $\D_{\lambda u} \F(0,0) [u_0] \notin \Range[\D_u \F(0,0)]$, where

    \end{enumerate}
    If $Z$ is any complement of $\Null[\D_u \F(0,0)]$ in $U$, then there is a
    neighbourhood $N$ of $(0,0)$ in $\R \times U$, an interval $(-a,a)$, and
    continuous functions $\phi : (-a,a) \to \R$, $\psi : (-a,a)\to Z$
    such that $\phi(0)=0, \psi(0)=0$ and
    \[
        \F^{-1}[0] \cap N = \lcb \lb \phi(s), \alpha u_0 + \alpha \psi(s) \rb :
        \abs{s} < a \rcb \cup \lcb (\lambda,0) : (\lambda,0)\in N \rcb.
    \]
    If $\D_{uu}\F$ is continuous then the functions $\phi$ and $\psi$ are once
    continuously differentiable.
\end{theorem}
\begin{theorem}[\parencite{Crandall1971}]\label{Thm:CrandallRegularity}
    In addition to the assumptions of \cref{Thm:CrandallMain}, let $\F$ be twice
    differentiable. If $\phi, \psi$ are the functions of
    \cref{Thm:CrandallMain}, then there is $\delta > 0$ such that
    $\phi^\prime(s) \neq 0$ and $0 < \abs{s} < \delta$ implies that
    $\D_u \F(\phi(s), \alpha u_0 + \alpha \psi(s))$ is an isomorphism of
    $U$ onto $V$.
\end{theorem}
\begin{theorem}[\parencite{Crandall1971}]\label{Thm:CrandallContinuity}
    In addition to the assumptions of \cref{Thm:CrandallMain}, suppose $\F$ has
    $n$ continuous derivatives with respect to $(\lambda, u)$ and $n + 1$ continuous
    derivatives with respect to $u$. Then the functions $(\phi, \psi)$ have $n$
    continuous derivatives with respect to $s$. If
    \[
        \D_u^{(j)}\F(0,0){[u_0]}^{j} = 0 \qquad 1 \leq j \leq n \quad\mbox{then }
            \phi^{(j)}(0) = 0,
    \]
    and
    \[
        \psi^{(j)}(0) = 0 \qquad\mbox{for } 1 \leq j \leq n - 1,
    \]
    and
    \[
        \frac{1}{n+1} \D_u^{(n+1)}\F(0,0){[u_0]}^{n+1} + \D_u\F(0,0)[\psi^{(n)}(0)]
         + \phi^{(n)}(0) \D_{\lambda u} \F(0,0)[u_0] = 0.
    \]
\end{theorem}

\begin{remark}
    $\D_u^{(j)}\F (0,0){[u_0]}^{j}$ means the value of the $j$-th Fr\'echet
    derivative\index{Fr\'echet derivative} of the map $x \to \F(0,x)$ at $(0,0)$
    evaluated at the $j$-tuple each of whose entries is $u_0$.
\end{remark}
The global version of \cref{Thm:CrandallMain} reads.
\begin{theorem}[\parencite{Lopez-Gomez2016}]\label{Thm:RabinowitzMainThm}
    Suppose $\Lap \in \C^1(\R, \Fred_0(U, V))$ and $\lambda_0 \in \R$ is a simple
    eigenvalue of $\Lap$, that is
    \[
        \Null[\Lap(\lambda_0)] = \vspan[\phi_0],
    \]
    and satisfies the following transversality condition
    \begin{equation}\label{Eqn:GlobalTransversalityCondition}
        \Lap^{\prime}(\lambda_0) \phi_0 \notin \Range[\Lap(\lambda_0)].
    \end{equation}
    Then, for every continuous function $\F : \R \times U \to V$
    satisfying~\ref{Assumption:1},~\ref{Assumption:2}, and~\ref{Assumption:3}
    and $\D_u \F(\cdot, 0) = \Lap(\cdot)$, $(\lambda_0, 0)$ is a
    bifurcation point to a continuum $\Cont$ of non-trivial solutions of $\F = 0$.
    For any of these $\F$'s,
    let $\Cont$ be the component of the set of non-trivial solutions of $\F = 0$
    with $(\lambda_0, 0) \in \Cont$. Then either,
    \begin{enumerate}
        \item $\Cont$ is not compact; or
        \item there is another $\Sigma \ni \lambda_1 \neq \lambda_0$ with
            $(\lambda_1, 0) \in \Cont$.
            %
    \end{enumerate}
\end{theorem}

\begin{figure}[!ht]\centering
    \includegraphics[width=.79\textwidth]{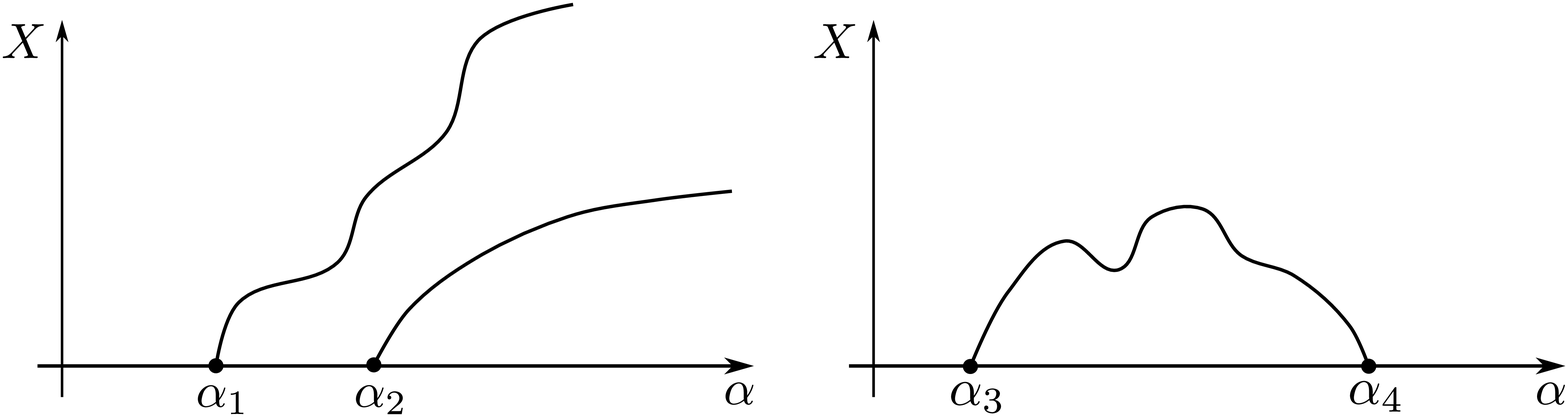}
    \caption[The Rabinowitz alternative]
    {Phasespace plot of the two possible alternatives of
    \cref{Thm:RabinowitzMainThm}. On the left, the bifurcation branches are
    unbounded, while on the right, the non-trivial solution branch connects two
    bifurcation points.}\label{Fig:PossibleBifs}
\end{figure}

And finally we have the following unilateral result.

\begin{theorem}[\parencite{Lopez-Gomez2016}]\label{Thm:UnilateralAbstractBifurcation}
    Suppose the injection $U \hookrightarrow V$ is compact,
    $\F$ satisfies~\ref{Assumption:1}-\ref{Assumption:3}, the map
    \[
      N[\lambda, u] = \F[\lambda, u] - \D_u \F(\lambda, 0)[u], \quad (\lambda,u)
        \in \R \times U
    \]
    admits a continuous extension\index{continuous extension} to $\R \times V$,
    the transversality condition\index{transversality condition}~\eqref{Eqn:GlobalTransversalityCondition}
    holds, and consider a closed subspace $Y \subset U$ such that
    \[
        U = \Null[\Lap(\lambda_0)] \oplus Y.
    \]
    Let $\Cont$ be the component given by \cref{Thm:RabinowitzMainThm} and denote by
    $\Cont^{+}$ and $\Cont^{-}$ the sub-components of $\Cont$ in the directions
    $\phi_0$ and $-\phi_0$ respectively. Then for each
    $\nu \in \{ +, - \}$, $\Cont^\nu$ satisfies some of the following alternatives:
    \begin{enumerate}
        \item $\Cont^\nu$ is not compact in $\R \times U$,
        \item There exists $\lambda_1 \neq \lambda_0$ such that $(\lambda_1, 0)
            \in \Cont^\nu$,
        \item There exists $(\lambda, y) \in \Cont^\nu$ with
            $y \in Y \setminus \{0\}$.
    \end{enumerate}
\end{theorem}
%
%
%
In \cref{Chapter:LocalBifPeriodic}, we apply \cref{Thm:CrandallMain} to find
local bifurcation branches originating from the trivial steady state solution of
the non-local equation~\eqref{Eqn:ArmstrongModelIntro}. In
\cref{Chapter:GlobalBifPeriodic}, we apply \cref{Thm:RabinowitzMainThm} and
\cref{Thm:UnilateralAbstractBifurcation} to obtain a global bifurcation result
for the steady states of equation~\eqref{Eqn:ArmstrongModelIntro}.
\section{The Averaging Operator\index{averaging operator} in Periodic Domains}
As we are dealing with functions on a circle of length $L$, $S^1_L$, we
implicitly defined periodic boundary conditions. Sometimes it is useful to
explicitly use these boundary conditions. For this reason, we define the
boundary operator,\index{boundary operator}
\begin{equation}\label{Eqn:DefnPeriodicBoundaryConditions}
    \Bd[u, u^{\prime}] \coloneqq (u(0) - u(L), u^{\prime}(0) - u^{\prime}(L)),
\end{equation}
which of course has to be equal to zero if we impose periodic boundary
conditions.\index{periodic boundary conditions} The abstract formulation in
terms of an operator equation will be facilitated by the following operators
that we will now define.
\begin{definition} Let $L>0$.
\begin{enumerate}
\item
    We define the set of test functions\index{test functions} to be
    \[
        \CpInf(0, L) \coloneqq \lcb f \in \CInf(0, L) :  f \; \mbox{is}\ L \ \mathrm{periodic} \rcb.
    \]
    Then, for example, $\HperOne(0, L)$ is defined as the completion of
    $\CpInf(0, L)$ with respect to the \HOne\ norm.
\item
    The averaging operator
    \begin{equation}\label{Defn:AvgOperator}
        \Avg : L^{p}(S^1_L) \mapsto \R, \quad
            \Avg[u] \mapsto \frac{1}{L} \int_{0}^{L} u(x) \dd x,
    \end{equation}
     is continuous and compact.
\item
    We define the following sub-manifold of $L^2(S^1_L)$
    \begin{equation}\label{Eqn:AvgZeroFun}
        L^2_{0} 
            \coloneqq \lcb u \in L^2(S^1_L) : \Avg[u] = 0 \rcb.
    \end{equation}
    \end{enumerate}
\end{definition}
\begin{lemma}
    $L^2_0(S^1_L)$ is closed and hence a Banach space.
\end{lemma}
\begin{proof}
    Note that $L^2_0(S^1_L) = \Avg^{-1}(0)$, hence it is closed as $\Avg$ is
    continuous. Finally, a closed subspace of a Banach space is again a Banach
    space.
\end{proof}
\begin{lemma}\label{Defn:Laplacian}
    We consider the Laplace operator\index{Laplace operator} $\Delta $ on $S^1_L$.
\begin{enumerate}
\item   %
    The Laplacian operator
    \[
        \Delta : \WperTwo(S^1_L) \mapsto L^2(S^1_L), \quad
            \Delta[v] \coloneqq v^{\prime\prime}
    \]
    is continuous.
\item
    $\Delta$ is a Fredholm operator. In particular, we have that
    \[
        \Null[\Delta] = \Range[\Avg] = \lcb f \in H^2_p(S^1_L) :
        f(x) \equiv c \in \R \rcb, \quad \Range[\Delta] = \Null[\Avg] = L^2_0.
    \]
    Further, we have that
    \[
        \dim \Null[\Delta] = \dim L^2(S^1_L) / \Range[\Delta] = 1,
    \]
    and thus $\ind \Delta = 0$.
\item
    The restriction operator\index{restriction operator} $\Delta_{\Avg}$
    \[
        \Delta_{\Avg} \coloneqq \Delta \Big|_{\Null[\Avg]} : \Null[\Avg]
	    \mapsto \Range[\Delta]
    \]
    is an isomorphism.
    \end{enumerate}
\end{lemma}
\begin{remark}
    Lemma~\ref{Defn:Laplacian} shows that $\Null[\Delta]$ and
    $\coker[\Delta] = L^2(S^1_L) / \Range[\Delta]$ are finite dimensional. We
    note that $\Avg$ can be used as a projection onto both $\Null[\Delta]$ and
    $\coker[\Delta]$.
\end{remark}
\section{Local and Global Existence}
The local and global existence of solutions to the adhesion model (\ref{Eqn:ArmstrongModelIntro}) in any space dimension
was presented in~\cite{Hillen2017}. We reformulate the multidimensional adhesion
model on $\R^n$ here:
\begin{equation}\label{ArmstrongNd}
u_t = D \Delta u - \alpha \nabla \cdot \lb u \int_{B_R(x)} h (u (x+r) )
    \Omega(r) \dd r \rb,
\end{equation}
where $B_R(x)$ denotes the ball of radius $R>0$ around $x$.
The function $\Omega(r)$ can be written as
\[
    \Omega(r) = \frac{r}{|r|} \omega(|r|).
\]
We assume

\begin{enumerate}[label=\textbf{(A\arabic*)},ref=(A\arabic*),leftmargin=*,labelindent=\parindent]
    \item\label{Existence:Assumption:1} $h\in \C^2(\R^n)$ and there exists a value $b>0$ such that $h(u)=0$ for all $u\geq b$.
    \item\label{Existence:Assumption:2} $\omega\in L^1(\R^n)$.
    \item\label{Existence:Assumption:3} For $p\geq 1$ let $u_0\in X_p\coloneqq \C^0(\R^n) \cap L^\infty(\R^n)\cap L^p(\R^n)$ be non-negative.
\end{enumerate}

\begin{theorem}[Corollary~2.4 in~\cite{Hillen2017}]
    Assume~\ref{Existence:Assumption:1}-\ref{Existence:Assumption:3}. Then there
    exists a unique, global solution
    \[
	u\in \C^0([0, \infty); X_p) \cap \C^{2,1}(\R^n\times (0,\infty))
    \]
   of~\eqref{ArmstrongNd} in the classical sense, with $u(0,x) = u_0(x)$, $x\in \R^n$.
\end{theorem}

%% file: chapterNonLocalOperator.tex
%
\chapter{Basic Properties}\label{chapter:BasicProperties}
In this chapter we define the non-local operator $\K[u]$, and we collect some
basic properties of $\K[u]$ in one spatial dimension. We prove results on
integrability, continuity, regularity, positivity, a-priori estimates and we
show that $\K[u]$ is a compact operator.  We analyse the corresponding spectrum
of $\K$ and we use these properties to derive properties of steady state
solutions such as symmetries, regularities and a-priori estimates. We find that
the non-local term $\K[u]$ acts like a non-local derivative, and the term
$\K[u]'$ acts like a non-local curvature, in a sense made precise later.
\section{Nondimensionalization\index{nondimensionalization} and Mass
Conservation\index{mass conservation}}
We use $S^1_L$ to denote the circle of length $L$, i.e. $S^1_L =\{\R \mod L\}$.
The adhesion model in one dimension on the unit circle $S^1_L$ is given by
\[\label{Eqn:ArmstrongModel}
    u_t(x, t) = D u_{xx}(x, t) - \alpha \lb u(x, t) \K[u(x, t)](x, t) \rb_x
\]
where the operator $\K[u]$ is given by,
\[
    \K[u(x, t)](x, t) = \int_{-R}^{R} h(u(x + r, t)) \Omega(r) \dd r.
\]
To non-dimensionalize the equation, we introduce the following non-dimensional
variables
\[
    x^* = \frac{x}{R}, \qquad t^* = t \frac{D}{R^2}, \qquad
    u^* = \frac{u}{\hat{u}}, \qquad \alpha^* = \frac{\alpha}{\hat{\alpha}},
\]
where $\hat{u}$ depends on the precise choice of the function $h(u)$, and
$\hat{\alpha}$ is given by
\[
    \hat{\alpha} = \frac{D}{R \hat{u}}.
\]
Finally let $\tilde{L} = \ifrac{L}{R}$, and
$\tilde \Omega(\tilde{r}) = \Omega(\tilde{r}R)$ for
$\tilde{r} \in \lsb -1 , 1 \rsb$. The non-dimensionalized model is then given by
\begin{equation}\label{Eqn:3:AdhModel}
    u_t(x, t) = u_{xx}(x, t) - \alpha \lb u(x, t) \int_{-1}^{1}
        \tilde{h}(u(x + \tilde{r}, t)) \tilde{\Omega}(\tilde r) \dd \tilde{r} \rb_{x},
\end{equation}
for $x\in S^1_L$ and $L>2$.
In the following, to make our notation simpler we will drop all tildes from the
previous equation. Typical solutions of equation~\eqref{Eqn:3:AdhModel}
are shown in \cref{Fig:TypicalSinglePopulationSolutions}.

As we do not consider any population dynamics (cell production or cell death)
in equation~\eqref{Eqn:3:AdhModel} it is easy to see that mass in the system
is conserved.
\begin{lemma}\label{Lem:ConsMass}
    Let $u \in \C^2(S^1_L)$ be a classical solution of
    \eqref{Eqn:3:AdhModel}. The total mass of the population $u(x, t)$
    \begin{equation}\label{Eqn:PopMean}
        \mean{u}(t) \coloneqq \frac{1}{L} \int_{0}^{L} u(x, t) \dd x,
    \end{equation}
    is conserved over time.
\end{lemma}
\begin{proof}
    We proceed by computing,
    \[
    \begin{split}
        L \frac{\dd \mean{u}}{\dd t} &= \int_{0}^{L} u_t(x, t) \dd x
            = \int_{0}^{L} \lb u_x - \alpha u \K[u] \rb_x \dd x \\
            &= u_x(L) - u_x(0) - \alpha u(0) \lb \K[u](L) - \K[u](0) \rb = 0.
    \end{split}
    \]
\end{proof}
\begin{figure}\centering
    \includegraphics[width=\textwidth]{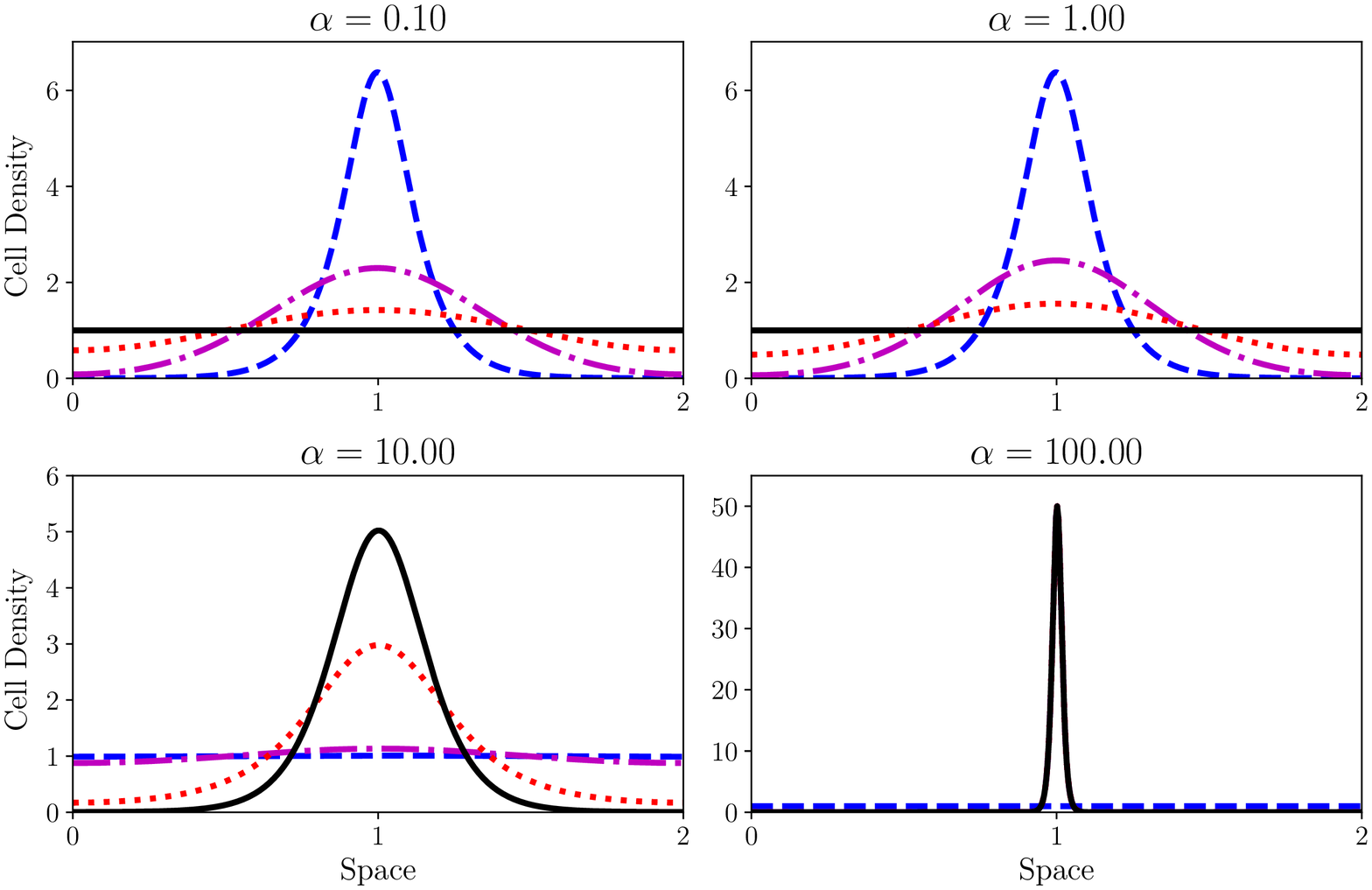}
    \caption[Bifurcation values]{Typical solution behaviour of
    equation~\eqref{Eqn:3:AdhModel} for varying values
    of $\alpha$. The initial condition is shown in blue (dashed), while the final
    steady state solution is shown in black (solid). The remaining curves show
    intermediate times.
    }\label{Fig:TypicalSinglePopulationSolutions}
\end{figure}
\section{The Non-Local Operator in 1-D}\label{section:NonLocalOperator}
The key part of the adhesion model is the non-local adhesion
operator\index{non-local adhesion operator} as defined as follows.
\begin{definition}\label{Defn:NonLocalOperator}
    Let $X, Y$ be Banach spaces of functions, then we define the operator
    $\K : X \to Y$ by
    \begin{equation}\label{Eqn:NonLocalDefn}
        \K[u(x)](x) = \int_{-1}^{1} h(u(x + r)) \Omega(r) \dd r.
    \end{equation}
    The domain of integration $[x-1, x+1]$ is called the \textit{sensing
    domain}\index{sensing domain}.
    \end{definition}
    Notice that we normalized the sensing radius to $R=1$ and we assume that the
    domain is larger than one sensing diameter, \ie\ $L > 2$.
    The directionality function $\Omega$ is assumed to satisfy the following
    conditions:
    \begin{enumerate}[label=\textbf{(K\arabic*)},ref=K\arabic*, leftmargin=*,
            labelindent=\parindent]
        \item\label{OmegaAssumption:1}
            $\Omega(r) = \frac{r}{\abs{r}} \omega(r)$, \mbox{where $\omega(r)$ is an even function.}
        \item\label{OmegaAssumption:2} $\omega(r) \geq 0$,
        \item\label{OmegaAssumption:3} $\omega \in L^{1}(0, 1) \cap L^{\infty}(0, 1)$,
        \item\label{OmegaAssumption:4} $\norm{\omega}_{L^1(0, 1)} = \ifrac{1}{2}$,
        \item\label{OmegaAssumption:5} $M_n(\omega) > 0$ for infinitely many
            integers $n$, where
            \begin{equation}\label{Defn:Mn}
                M_n(\omega) = \int_{0}^{1} \sin\lb\frac{2\pi n r}{L}\rb\omega(r)\dd r.
            \end{equation}
    \end{enumerate}
    Condition (\ref{OmegaAssumption:1}) and condition (\ref{OmegaAssumption:3}) imply that
    \[
        \int_{-1}^1 \Omega (r) \dd r = 0 \qquad\mbox{ and } \qquad \int_{-1}^1
        \omega(\abs{r}) \dd r =1.
    \]
    The Fourier-sine coefficients\index{Fourier-sine coefficients} $M_n(\omega)$ of $\omega$ in
    (\ref{OmegaAssumption:5}) are related to the eigenvalues of $\K$, as we will
    show later.  The function $h(\cdot)$ within the integral describes the nature of
    the adhesive force and is assumed to satisfy:
    \begin{enumerate}[label=\textbf{(H\arabic*)},ref=H\arabic*, leftmargin=*,
            labelindent=\parindent]
        \item\label{NatureForceAssumption:1}
            $h \in \C^2(\R)$,
        \item\label{NatureForceAssumption:2}
            $h(u) \geq 0$ for $u \geq 0$,
        \item\label{NatureForceAssumption:3}
            $h(u) \leq C \lb 1 + u \rb$ for all $u \geq 0$, for some
                $\R \ni C \geq 0$,
        \item\label{NatureForceAssumption:4}
            There exist positive real numbers $\bar u$ such that $h^{\prime}(\mean{u}) \neq 0$.
    \end{enumerate}

\begin{remark}
    Using the assumptions~(\ref{OmegaAssumption:1}) to~(\ref{OmegaAssumption:3}) we
    can rewrite the non-local function defined in
    equation~\eqref{Eqn:NonLocalDefn} as
    \[
        \K[u] : x \mapsto \int_{0}^{1} \lsb h(u(x + r)) - h(u(x - r)) \rsb
        \omega(r) \dd r.
    \]
    This equivalent formulation will be frequently used in the following. Since
    $\omega(r)$ is an even function, we see already here that $\K[u]$ appears as
    a non-local derivative of $h(u(x))$.
\end{remark}

For the even function $\omega(r)$ with $r\in[0,1]$ there are three commonly used forms
(see for instance \parencite{Painter2015}), which we illustrate in \cref{Fig:OmegaExamples}.
\begin{enumerate}[label=\textbf{(O\arabic*)},ref=O\arabic*, leftmargin=*,
        labelindent=\parindent]
    \item\label{OmegaChoice:1} Uniform distribution\index{uniform distribution}
        \[
            \omega(r) = \frac{1}{2}, \qquad 0\leq r\leq 1.
        \]
    \item\label{OmegaChoice:2} Exponential distribution\index{exponential distribution}
        \[
            \omega(r) = \omega_0 \exp\lb -\frac{r}{\xi}\rb, \qquad 0\leq r\leq 1,
        \]
        where $\xi$ is a parameter controlling how quickly $\omega(\cdot)$ goes
        to zero, and $\omega_0$ is a normalization constant.
    \item\label{OmegaChoice:3} Peak signalling\index{peak signalling} a distance $\xi$ away from the
        cell centre.
        \[
            \omega(r) = \omega_0 \frac{r}{\xi}
                \exp\lb -\frac{1}{2} \lb\frac{r}{\xi}\rb^2\rb, \qquad 0\leq r\leq 1,
        \]
        where $\xi$ is a parameter controlling how quickly $\omega(\cdot)$ goes
        to zero, and $\omega_0$ is a normalization constant.
\end{enumerate}

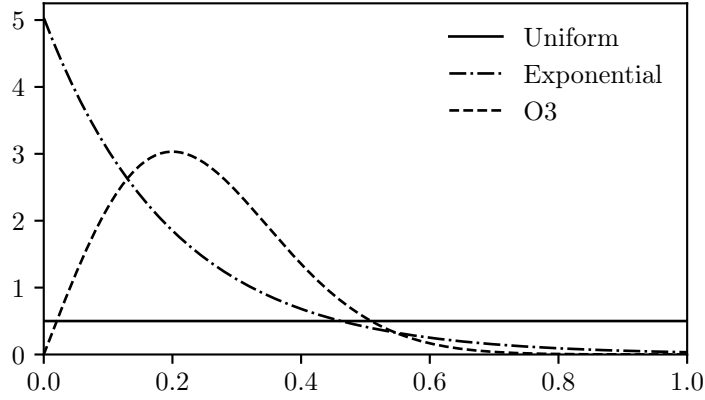
\begin{figure}[!ht]\centering
    \input{DifferentOmegaDistributions.pgf}
    \caption{The different distributions for $\omega(\cdot)$ with $\xi =
    \ifrac{1}{4}$.}\label{Fig:OmegaExamples}
\end{figure}

The next lemma provides an example of a situation in which there are integers
$n$ such that $M_n(\omega) = 0$.
Suppose that $\omega \equiv \ifrac{1}{2}$ and that $L = 2 k R$ where $k \in \N$.
Then we have the following result, indicating a degenerate situation.
\begin{lemma}\label{Lemma:DegenerateCase}
    Let $L = 2k$, for $\N \ni k > 0$, and $\omega(r) \equiv \ifrac{1}{2}$, then
    for $n \in \N$ such that $\frac{n}{k}$ is an even integer, we have that
    \[
        M_n(\omega) = 0.
    \]
\end{lemma}
\begin{proof}
    The term $M_n(\omega)$ is defined in equation~\eqref{Defn:Mn}, as
    \[
        M_n(\omega) = \int_{0}^{1} \sin\lb\frac{2 \pi n x}{L}\rb \Omega(r) \dd r
            = \frac{1}{2} \int_{0}^{1} \sin\lb\frac{\pi n x}{k}\rb \dd r,
    \]
    then this is zero whenever $\frac{n}{k}$ is an even integer.
\end{proof}

The significance of $M_n(\omega)$ being zero is that under those circumstances
we do not have a bifurcation point (see \cref{Chapter:LocalBifPeriodic}). As we
will show later, if $M_n(\omega)$ is zero it is impossible to obtain a
steady-state solution having $n$-peaks. Further note, that both the properties
of $\omega$ and the domain length $L$ determine whether $M_n(\omega)$ is zero.

In the next lemma we establish continuity and $L^p$ properties of the operator
$\K[u]$.
\begin{lemma}\label{Lem:KProp}
    Assume (\ref{OmegaAssumption:1}-\ref{OmegaAssumption:4}) and (\ref{NatureForceAssumption:1}-\ref{NatureForceAssumption:4}).
    We assume $p\geq 1$. The function
    \begin{equation}\label{Eqn:NonLocalFunction}
        \K[u] : x \mapsto \K[u](x),
    \end{equation}
    has the following properties:
    \begin{enumerate}[label=\textbf{(P\arabic*)},ref=(P\arabic*), leftmargin=*,
        labelindent=\parindent]
        \item\label{NonLocalProp:1} Let $\N \ni s \leq 2$. If $u \in \C^s(S^1_L)$ then
            $\K[u]\in \C^s(S^1_L)$.
        \item\label{NonLocalProp:2} If $u \in L^p(S^1_L)$ then $\K[u]\in \C^0(S^1_L)$.
        \item\label{NonLocalProp:3} If $u \in L^p(S^1_L)$ then $\K[u] \in L^p(S^1_L)$,
            and there exists a constant $\tilde C>0$ such that
            \[
                \abs{\K[u]}_p \leq \abs{h(u)}_{p} \leq \tilde C \lb \abs{u}_{p} + L \rb.
            \]
        \item\label{NonLocalProp:4} Assume $k = 1,2$. If $u \in W^{k,p}$,
            then $\K[u] \in W^{k, p}$.
        \item\label{NonLocalProp:5} The operator $\K : L^p \to L^p$ is
            continuously Fr\'echet differentiable in $u$.
        \item\label{NonLocalProp:6}  If $u \in L^p(S^1_L)$ such that $u(x) \geq 0$
            (\ie\ $\abs{u}_1 = L \Avg[u] < \infty$). Then
            \[
                \abs{\K[u](x)} \leq C \abs{\omega}_{\infty} L \lb \Avg[u] + 1 \rb,
                \quad\mbox{with}\quad \Avg[u] = \frac{1}{L}\int_0^L u(x) \dd x.
            \]
    \end{enumerate}
\end{lemma}
\begin{proof}
\begin{enumerate}[leftmargin=*, labelindent=0pt]
    \item If $u\in \C^0(S^1_L)$, then
    \[ K[u](x) = \int_{-1}^1 h(u(x+r)) \Omega(r) \dd r \]
    is continuous in $x$, since $h$ is continuous.
    If $u\in \C^1(S^1_L),$ then
    \[ K[u](x)' = \int_{-1}^1 h'(u(x+r)) u'(x+r)  \Omega(r) \dd r \]
    is continuous in $x$, since $h$ is differentiable.
    If $u\in \C^2(S^1_L)$, then
    \[ K[u](x)'' = \int_{-1}^1 \left(h''(u(x+r))(u'(x+r))^2 + h'(u(x+r)) u''(x+r)\right) \Omega(r) \dd r \]
    is continuous in $x$, since $h$ is $\C^2$.

    \item Let $u \in L^p(S^1_L)$.
        Let $\epsilon > 0$ and let $x_1, x_2 \in S^1_L$ such that $\abs{x_1 - x_2} <
        \delta$. By the density of $\C^0$ in $L^p$ there is a sequence $\seq{u_n}
        \subset \C^0$ such that $u_n \to u$ in $L^p$.
        This means that $\exists N : \forall n \geq N$
        we have that $\abs{u_n - u}_{p} < \ifrac{\epsilon}{3 \abs{\omega}_{\infty}}$.
        Then, we compute
        \[
        \begin{split}
            \abs{\K[u](x_1) - \K[u](x_2)} \leq
                                  &\abs{\K[u](x_2) - \K[u_n](x_2)} +
                                  \abs{\K[u_n](x_1) - \K[u](x_1)} + \\
                                  &\abs{\K[u_n](x_2) - \K[u_n](x_1)}.
        \end{split}
        \]
        The first two terms are treated equivalently. Let $n > N$, and $i = 1,2$,
        then
        \[
        \begin{split}
            \abs{\K[u](x_i) - \K[u_n](x_i)} &\leq
                \int_{0}^{1} \abs{h(u(x_i + r)) - h(u_n(x_i + r))} \omega(r) \dd r \\
            &+   \int_{0}^{1} \abs{h(u(x_i - r)) - h(u_n(x_i - r))} \omega(r) \dd r \\
            &\leq C \int_{0}^{1} \abs{u(x_i + r) - u_n(x_i + r)} \omega(r) \dd r \\
            &+   C  \int_{0}^{1} \abs{u(x_i - r) - u_n(x_i - r)} \omega(r) \dd r \\
            &\leq \lb \int_{0}^{L} \abs{u(x) - u_n(x)}^{p} \dd x \rb^{1/p}
                \lb \int_{0}^{1} \omega^{q}(r) \dd r \rb^{1/q} \\
            &\leq \abs{u - u_n}_{p} \abs{\omega}_{\infty} < \ifrac{\epsilon}{3}.
        \end{split}
        \]
        The last term can be estimated by continuity of $\K[u_n](x)$ as
        \[
            \abs{\K[u_n](x_2) - \K[u_n](x_1)} < \ifrac{\epsilon}{3}.
        \]
        Putting everything together, we obtain
        \[
            \abs{\K[u](x_1) - \K[u](x_2)} < \epsilon.
        \]
        This shows that $\K[u](x) \in \C^0(S^1_L)$.

    \item Let $1 \leq p$, and let $u \in L^p(S^1_L)$.
        By applying the Minkowski Integral Inequality\index{Minkowski Integral
        Inequality} to equation~\eqref{Eqn:NonLocalDefn}, we obtain
        \[
        \begin{split}
            \abs{\K[u](x)}_{p} &=
                \abs*{\int_{-1}^{1} h(u(x + r, t)) \Omega(r) \dd r }_{p} \leq
                \int_{-1}^{1} \abs{h(u(x + r, t))\Omega(r)}_{p} \dd r \\
                &= \int_{-1}^{1} \lcb \int_{0}^{L} \abs{h(u(x+r, t))}^p \dd x
                    \rcb^{1/p} \abs{\Omega(r)} \dd r
                    \leq \abs{h(u)}_{p} \\
                &\leq \tilde C \lb \abs{u}_{p} + L \rb,
        \end{split}
        \]
        using assumption~(\ref{NatureForceAssumption:3}). Note that due to
        assumption~(\ref{OmegaAssumption:3}) we have that ${\abs{\Omega}_{1} = 1}$.
    \item Let $u \in W^{1,p}(S^1_L)$.
        The first derivative with respect to $x$ of the function defined in
        equation~\eqref{Eqn:NonLocalFunction} is given by
        \[
            \lb \K[u](x) \rb^{\prime} = \int_{-1}^{1} h^{\prime}(u(x + r))
            u^{\prime}(x + r) \Omega(r) \dd r.
        \]
        Its $L^p$ norm can be estimated using~\ref{NonLocalProp:3}, we then compute
        \[
        \begin{split}
            \norm{\K[u]}_{1, p} &=
                \lb \abs{\K[u]}_{p}^{p} + \abs{\lb\K[u]\rb^{\prime}}_p^p \rb^{1/p}
                \leq \lb \tilde C^p \lb \abs{u}_{p} + L \rb^p
                + \abs{h^{\prime}(u) u^\prime}_{p}^{p} \rb^{1/p} \\
                &\leq \lb 2^p \tilde  C^p \lb L^p + \abs{u}_{p}^{p} \rb +
                \abs{h^{\prime}}_{\C^0}^{p} \abs{u^{\prime}}_{p}^{p} \rb^{1/p} \\
                &\leq 2 \tilde C \abs{h^{\prime}}_{\C^0} \lb L^p + \abs{u}_{p}^{p}
                    + \abs{u^{\prime}}^{p}_{p} \rb^{1/p}.
        \end{split}
        \]
        Then using that $u \in W^{1, p}(S^1_L)$ and
        assumption~(\ref{NatureForceAssumption:1}), all the terms on the right hand
        side are bounded.

        Now let $u \in W^{2, p}(S^1_L)$.
        The second derivative with respect to $x$ of the function defined in
        equation~\eqref{Eqn:NonLocalFunction} is given by
        \[
            \lb \K[u](x) \rb^{\prime\prime} =
                \int_{-1}^{1} \lb h^{\prime}(u(x + r)) u^{\prime\prime}(x + r) +
                    h^{\prime\prime}(u(x + r)) \lb u^{\prime}(x + r) \rb^2
                    \rb \Omega(r) \dd r.
        \]
        Its $L^p$ norm can be estimated using~\ref{NonLocalProp:3},
        we then compute
        \[
        \begin{split}
            \abs{{\K[u]}^{\prime\prime}}_{p} &\leq
                \abs{h^{\prime}(u) u^{\prime\prime}}_{p} +
                \abs{h^{\prime\prime}(u) \lb u^{\prime} \rb^{2}}_{p} \\
                &\leq \abs{h^{\prime}}_{\C^0} \abs{u^{\prime\prime}}_{p} +
                    \abs{h^{\prime\prime}}_{\C^0} \abs{u^{\prime}}_{p}^{2}.
        \end{split}
        \]
        Then using that $u \in W^{2, p}(S^1_L)$ and
        assumption~(\ref{NatureForceAssumption:1}), all the terms on the right hand
        side are bounded. Combining this result with the estimates of~\ref{NonLocalProp:4},
        we obtain the required result.

    \item Consider the map
        \[
            \K :  L^p(S^1_L) \mapsto L^p(S^1_L).
        \]
        The Fr\'echet derivative of $\K$ is computed by
        \[
        \begin{split}
            \D_{u} \lb \K[u] \rb [w](x) &= \frac{\dd}{\dd \epsilon} \Big|_{\epsilon = 0}
                \int_{-1}^{1} h(u(x+r) + \epsilon w(x+r)) \Omega(r) \dd r \\
                &= \int_{-1}^{1} h^{\prime}(u(x + r)) w(x+r) \Omega(r) \dd r.
        \end{split}
        \]
        We show that the map
        \[
        \begin{split}
            \D_{u} \K : L^p(S^1_L) &\mapsto \Lin(L^p(S^1_L), L^p(S^1_L)) \\
               u &\mapsto \D_{u} \lb \K[u] \rb[w],
        \end{split}
        \]
        is continuous. Meaning, that we have to show that the operator norm of
        $\D_u \K$ is bounded. From~\ref{NonLocalProp:3}, we have
        \[
            \abs{\D_{u}\lb \K[u] \rb[w]}_{p} \leq \tilde C \lb \abs{h^\prime}_{C^0}
            \abs{w}_{p} + L \rb.
        \]
        Hence
        \[
            \norm{\D_{u} \lb\K[u]\rb[w]}_{\text{op}} = \sup_{\abs{w}_{p} = 1}
                \abs{\D_{u}\lb\K[u]\rb[w]}_{p} \leq \tilde C \lb \abs{h^\prime}_{C^0} + L
                \rb.
        \]

    \item Let $1 \leq p$ such that $\Avg[u] < \infty$ and $u(x) \geq 0$. Then,
        we compute
        \[
        \begin{split}
            \K[u](x) &= \int_{-1}^{1} h(u(x + r)) \Omega(r) \dd r \\
                     &= \int_{0}^{1}  h(u(x + r)) \omega(r) \dd r
                     - \int_{-1}^{0} h(u(x + r)) \omega(r) \dd r.
        \end{split}
        \]
        It is easy to see that both integrals on the right hand side are
        non-negative as $h(u(x)) \geq 0$ whenever $u(x) \geq 0$
        (assumption~(\ref{NatureForceAssumption:2}))
        and $\omega(r) \geq 0$. Using assumption~(\ref{NatureForceAssumption:3}),
        it follows that
        \[
            \K[u](x) \leq \int_{0}^{1} h(u(x + r)) \omega(r) \dd r
                    \leq C \abs{\omega}_{\infty} L\lb \Avg[u] + 1 \rb.
        \]
        In the same spirit, we find for
        \[
            \K[u](x) \geq - \int_{-1}^{0} h(u(x + r)) \omega(r) \dd r
                    \geq - C \abs{\omega}_{\infty} L \lb \Avg[u] + 1 \rb.
        \]
\end{enumerate}
\end{proof}
\begin{lemma}
    The operator $\K : L^2(S^1_L) \mapsto L^2(S^1_L)$ is compact.
\end{lemma}
\begin{proof}
    Since $L^2(S^1_L)$ is a Hilbert space, we can use an orthonormal
    system\index{orthonormal system}
    $\{\phi_i\}_{i\in \N} $ to  define a finite truncation\index{finite truncation} of the integral
    operator. We denote the inner product in $L^2$ with $(\cdot, \cdot)$. Given $u\in L^2(S^1_L)$ we consider the expression $h(u(x+r))$ to be a function of two variables $x$ and $r$ and we expand it as
\[ h(u(x+r)) = \sum_{i,j=1}^\infty h_{ij}(u) \phi_i(x) \phi_j(r) \]
with \[ h_{ij}(u) = \Bigl(h(u(x+r)), \phi_i(x)\phi_j(r)\Bigr). \]
We define a finite approximation of $h(u)$ as
\[ \Bigl(h(u(x+r))\Bigr)_n := \sum_{i.j = 1}^n h_{ij}(u) \phi_i(x) \phi_j(r) \] and the corresponding truncated integral operator
\[ \K_n[u](x): = \int_{-1}^1 \Bigl(h(u(x+r))\Bigr)_n \Omega(r) \dd r \]
For each $n>0$ the operator $\K_n$ has finite rank, hence it is compact
    (\cite[Lemma~3.12]{Robinson2001}).  If we can show that $\K_n \to \K$ as
    $n\to\infty$ in the operator norm, then by
    \cite[Theorem~3.10]{Robinson2001} also $\K$ is compact. Indeed
    \begin{eqnarray*}
    \abs{\K-\K_n}_2^2 &\leq & \int_{0}^L\int_{-1}^1 \left| h(u(x+r))-\Bigl(h(u(x+r))\Bigr)_n\right|^2 \Omega^2(r) \dd r \dd x \\
    &\leq& \int_{0}^L \int_{0}^L \left| \sum_{i,j=n+1}^\infty h_{ij}(u) \phi_i(x)\phi_j(r) \right|^2 \Omega^2(r) \dd r \dd x\\
    &\leq& \sum_{i,j=n+1}^\infty |h_{ij}(u)|^2 \|\Omega\|^2_\infty\
    \to 0 \quad\mbox{as}\quad n\to\infty.
    \end{eqnarray*}
\end{proof}
The compactness of $\K$  has far reaching consequences, namely it means that
$\K[u]$ can never be invertible, as otherwise $\I = \K^{-1} \K$ would be
compact. Further a compact operator can never be surjective, and the closed
subspaces of its range must be finite dimensional. Also, the spectral theorems
for compact operators apply to $\K$, which we will employ next.
\section{Spectral Properties}
In this section, we consider the spectral properties of linear non-local
operator $\K[u]$ (\ie\ $h(u) = u$) and of the non-local curvature $\K[u]'$. That
is
\begin{equation}\label{Eqn:LinearK}
    \K[u] = \int_{-1}^{1} u(x + r) \Omega(r) \dd r.
\end{equation}
The results of this section will be used in \cref{Chapter:LocalBifPeriodic} to
study the properties of the linearization of the steady state equation. We start by a result how $\K[u]$ acts
on basis functions of $L^2$.
\begin{lemma}\label{Lemma:IntegralIdentities}
    Let $\Omega$ satisfy (\ref{OmegaAssumption:1}), (\ref{OmegaAssumption:2})
    and (\ref{OmegaAssumption:3}), then
    \[
        \int_{-1}^{1} \cos\lb\frac{2\pi n r}{L}\rb \Omega(r) \dd r = 0,
    \]
    and
    \[
        \int_{-1}^{1} \sin\lb\frac{2\pi n r}{L}\rb \Omega(r) \dd r
            = 2 \int_{0}^{1} \sin\lb\frac{2\pi n r}{L}\rb \omega(r) \dd r.
    \]
\end{lemma}
\begin{proof} Notice that $\Omega(r) = \frac{r}{|r|} \omega(r)$, where $\omega(r)$ is an even function, hence
    both identities follow by integration and symmetry.
\end{proof}
\begin{lemma}\label{Lem:EfuncK}
    Consider the operator $\K[u]$ with $h(u)=u$ as defined in
    equation~\eqref{Eqn:LinearK} as an operator
    \[
        \K : L^2(S^1_L) \mapsto L^2(S^1_L).
    \]
    $\K$ has the following properties:
    \begin{enumerate}
        \item $\K$ is bounded and skew-adjoint\index{skew-adjoint} \ie\ $\K^{*} = - \K$.
        \item $\K$ maps the canonical basis functions of $L^2(S^1_L)$ as follows:
            \begin{align}
                \K[1](x) &= 0,\nonumber \\
                \K\lsb\sin\lb\frac{2\pi n x}{L}\rb\rsb(x) &=
                    2 M_n(\omega)\cos\lb\frac{2\pi n x}{L}\rb,\nonumber \\
                \K\lsb\cos\lb\frac{2\pi n x}{L}\rb\rsb(x) &=
                    -2 M_n(\omega) \sin\lb\frac{2\pi n x}{L}\rb,\nonumber
            \end{align}
            where $M_n(\omega)$ is defined in equation~\eqref{Defn:Mn}.
        \item If instead of $L^2(S^1_L)$ we consider its canonical
            complexification.\index{complexification} Then, we find complex eigenvalues:
            \begin{align*}
                \K[1](x) &= 0,\nonumber \\
                \K\lsb\exp\lb\frac{2\pi n i x}{L}\rb\rsb(x) &=
                    2 i M_n(\omega)\exp\lb\frac{2\pi n i x}{L}\rb.
            \end{align*}
    \end{enumerate}
\end{lemma}
\begin{proof}
$\K$ is bounded by \cref{Lem:KProp}. To see the skew-adjointness,
    we consider for  $u,v\in L^2(S^1_L)$.
\begin{eqnarray*}
(\K[u], v) &=& \int_0^L\int_{-1}^1 u(x+r) \Omega(r) \dd r\, v(x) \dd x\\
&=& \int_{-1}^1 \int_r^{L+r} u(y) \Omega(r) v(y-r) \dd r \dd y \\
&=& \int_0^L \int_{-1}^1 u(y) \Omega(r) v(y-r) \dd r \dd y \\
&=& - \int_0^L \int_{-1}^1 v(y+r) \Omega(r) \dd r\, u(y) \dd y \\
&=& - ( u, \K[v]),
\end{eqnarray*}
where we used periodicity in the second step, and the property that $\Omega$ is odd in the third step.
Properties (2) and (3) follow by direct computation, using the symmetries of $\sin, \cos$ and $\Omega$.
\end{proof}
\begin{remark}
\begin{enumerate}
\item
    Note that in general, skew-adjoint operators have purely imaginary
    eigenvalues.
\item     Note that since any skew-adjoint operator is normal, we have that $\K$ is a
    normal operator. Thus we have that $\K$ is a compact and normal
        operator\index{normal operator} on
    $L^2(S^1_L)$. This means that it is also a compact, and normal operator on
    the canonical complexification of $L^2(S^1_L)$ ($H = L^2 + i L^2$) and hence
    we can apply a spectral theorem \cite[Theorem 7.53]{fabian2013} to obtain an
    orthonormal basis  of $H$ over which the operator $\K$ is
        diagonalizable.\index{diagonalizable}
\item
    Note that it is easy to see from \cref{Lem:EfuncK} that
    the non-local operator $\K[u]$ removes mass, in the sense that $\Avg[\K[u]] = 0$
    (the average operator $\Avg$ was defined in~\eqref{Defn:AvgOperator}).
    Hence, for example if $\K : L^2(S^1_L) \to L^2(S^1_L)$, then we conclude
    that $\Range[\K] = L^2_0$ (where $L^2_0$ contains $L^2$ functions with
    mass zero, as defined in~\eqref{Eqn:AvgZeroFun}).
\item
    For the sake of comparison, the eigenvalues corresponding to the
    eigenfunctions $\exp\lb\frac{2\pi n i x}{L}\rb$ of the derivative operator
    are
    \[
        \lambda_n = \frac{2 \pi i n}{L}.
    \]
    \end{enumerate}
\end{remark}
Let us consider the spectrum of the derivative of $\K[u]$. From
\cref{Lem:KProp}~\ref{NonLocalProp:4}, we know that
${\K[u]}^{\prime}$ maps $H^1(S^1_L)$ into $H^1(S^1_L)$. Here we study the
operator, which we will refer to as the linear {\it non-local
curvature}\index{non-local curvature}
\begin{equation}\label{Eqn:KpFunc}
    {\K[u]}^{\prime} : x \mapsto \lb \int_{-1}^{1} u(x + r) \Omega(r) \dd r
    \rb^{\prime},
\end{equation}
where $(\cdot)^{\prime}$ denotes the spatial derivative with respect to $x$.
Note that using the properties of $\Omega$, this function can be rewritten as
\[
    \lb \K[u] \rb^{\prime} : x \mapsto
        \int_{0}^{1} \lb u^{\prime}(x + r) + u^{\prime}(x - r) \rb
            \omega(r) \dd r.
\]
\begin{lemma}
    The operator $\lb \K \rb^{\prime}$ given above is self-adjoint.\index{self-adjoint}
\end{lemma}
\begin{proof}
    Let $y, z \in L^2(S^1_L)$, we then compute using integration by parts and
    using \cref{Lem:EfuncK} to obtain
    \[
        \lb \lb \K[y] \rb^{\prime}, z \rb = - \lb \K[y], z^{\prime} \rb
            = \lb y, \K[z^{\prime}] \rb
            = \lb y, \lb \K[z] \rb^{\prime} \rb.
    \]
\end{proof}
\begin{lemma}\label{Lem:EfuncKp}
    Let $v_n$ be an eigenfunction of
    \[
    \left\{
        \begin{array}{@{}ll@{}}
            -v_n^{\prime\prime} = \lambda_n v_n, \quad &\mbox{in } \lsb 0, L \rsb \\
            \Bd[v_n, v_n^\prime] = 0,
    \end{array}
    \right.
    \]
    that is
    \[
    v_n\in \left\{ 1, \sin\left(\frac{2 n \pi x}{L}\right), \cos\left(\frac{2 n \pi x}{L}\right),\ n=1, 2, \dots  \right\}
    \]
    Then the
    operator $\lb \K \rb^{\prime}$ from equation~\eqref{Eqn:KpFunc}
    has the same set of eigenfunctions satisfying
    \[
        {\K[v_n]}^{\prime} = \mu_n v_n,
    \]
    where
    \[
        \mu_n = - \frac{4 \pi n}{L} M_n(\omega),
    \]
    where $M_n(\omega)$ is defined in equation~\eqref{Defn:Mn}.
\end{lemma}
\begin{proof}
    For $n = 0$ we have that $v_0 = 1$, then trivially $\K[1]' = 0$ and hence
    $\mu_0 = 0$. Next consider $v_n = \sin\lb\frac{2\pi n x}{L}\rb$, then
    \[
        {\K[v_n]}^{\prime} = \int_{-1}^{1} \lb\sin\lb\frac{2\pi n (x + r)}{L}\rb\rb^\prime
            \Omega(r) \dd r.
    \]
    Using the symmetries of $\sin, \cos$ and $\Omega$, we obtain
    \[
        {\K[v_n]}^{\prime}
            = - \frac{4 \pi n}{L} \sin\lb\frac{2\pi n x}{L}\rb M_n(\omega).
    \]
    Finally, consider $v_n = \cos\lb\frac{2\pi n x}{L}\rb$, then
    \[
        {\K[v_n]}^{\prime} =
            \int_{-1}^{1} \lb\cos\lb\frac{2\pi n (x + r)}{L}\rb\rb^\prime
            \Omega(r) \dd r.
    \]
    Once again using symmetries, we obtain
    \[
        {\K[v_n]}^{\prime} =
            - \frac{4 \pi n}{L} \cos\lb\frac{2\pi n x}{L}\rb M_n(\omega).
    \]
\end{proof}
From the definition of $M_n(\omega)$ in equation~\eqref{Defn:Mn}, we
easily see from (K4) that $\abs{M_n(\omega)} < \ifrac{1}{2}$. Here, we want to understand
in more detail how $M_n(\omega)$ behaves as $n \to \infty$. For this reason, we
introduce the following common definition from Fourier analysis (see for
instance~\cite[Chapter II]{zygmund2002}).
\begin{definition}
    The integral \textit{modulus of continuity}\index{modulus of continuity} is
    defined for periodic ${f \in L^p(S^1_L)}$, $p \geq 1$ by
    \[
        m_{p}(\delta) = \sup_{0 \leq h \leq \delta} \lcb
            \frac{1}{L} \int_{0}^{L} \abs{f(x + h) - f(x)}^p \dd x \rcb^{1/p}, \qquad \delta>0.
    \]
\end{definition}
It is obvious, that as $\delta \to 0$ we have that $m_p(\delta) \to 0$.
\begin{lemma}\label{Lemma:BoundMn}
    Let $\omega(r)$ satisfy (\ref{OmegaAssumption:1}), (\ref{OmegaAssumption:2}),
    and (\ref{OmegaAssumption:3}) and $M_n(\omega)$ be given by equation~\eqref{Defn:Mn}),
    then we have that
    \[
        \abs{M_n(\omega)} \leq \frac{L}{2} m_{1}\lb\frac{L}{2n}\rb.
    \]
\end{lemma}
\begin{proof}
    First, extend $\omega(r)$ to the whole of $[0, L]$ by defining the
    $L$-periodic function
    \[
        \tilde{\omega}(r) =
        \begin{cases} \omega(r) &\mbox{if } r \leq 1 \\
            0 &\mbox{otherwise}
        \end{cases}.
    \]
    Then, we use an common technique from Fourier theory (see for
    instance~\cite[Chapter II]{zygmund2002}).
    \begin{eqnarray*}
        M_n(\omega) & =&  \int_{0}^{L} \sin\lb\frac{2\pi n r}{L}\rb
            \tilde{\omega}(r) \dd r\\
            &=& -\int_0^L \sin\left(\frac{2\pi n}{L} \left(r-\frac{L}{2n}\right)\right) \tilde \omega(r) \dd r \\
            &=& - \int_{0}^{L} \sin\lb\frac{2\pi n r}{L}\rb
                \tilde{\omega}\left(r + \frac{L}{2n}\right) \dd r.
    \end{eqnarray*}
    Taking the average of both integrals, we obtain
    \[
        M_n(\omega) = \frac{1}{2} \int_{0}^{L} \left( \tilde{\omega}(r) -
        \tilde{\omega}\left(r + \frac{L}{2n}\right) \right) \sin\lb\frac{2\pi n r}{L}\rb \dd r.
    \]
    Hence,
    \[
    M_n(\omega) \leq \frac{L}{2} \sup_{0<h\leq \frac{L}{2n}} \frac{1}{L}\int_0^L |\tilde \omega(r)-\tilde \omega(r+h)| \dd r  = \frac{L}{2} m_1\left(\frac{L}{2n}\right).
    \]
\end{proof}
\begin{example}\label{Example:UniformOmega1}
    Suppose that $\omega(r) = \frac{1}{2}$ is chosen to be the uniform
    function, example  \ie~(\ref{OmegaChoice:1}).
    Then, we can compute $M_n(\omega)$ and find that
    \[
        M_n(\omega) = \frac{L}{2 \pi n} \sin^2\lb\frac{\pi n}{L} \rb.
    \]
    Hence, $M_n(\omega) \to 0$ as $n \to \infty$, and thus so do the eigenvalues
    of $\K$ from equation~\eqref{Eqn:LinearK} see
    \cref{Lem:EfuncK} (since $\K$ is compact this is
    expected as the only possible accumulation point of the eigenvalues is
    zero). But the eigenvalues of the non-local
    curvature~\eqref{Eqn:KpFunc} are given by (see
    \cref{Lem:EfuncKp})
    \[
        \mu_n = - \frac{4 \pi n}{L} M_n(\omega)
            = - 2 \sin^2 \lb\frac{\pi n}{L}\rb.
    \]
    Thus, the eigenvalues of the non-local curvature keep oscillating in
    $(-2, 0)$.
\end{example}
%
%
\section{The Behaviour of $\K_{R}$ for $R\to 0$}
In this subsection, we show that the linear non-local operator $\K[u]$ defined
in \cref{Defn:NonLocalOperator} is related to the classical derivative.  In
particular, we show that for smooth functions $h$ and $u$, the non-local
operator approximates the gradient as the sensing radius $R$ converges to zero,
hence we call it a non-local gradient\index{non-local gradient}. For this
section we give up the normalization of a sensing radius of $R=1$ and we
consider the non-local operator for general $R>0$
\[
    \K_{R}[u] \coloneqq \int_{-R}^{R} h(u(x + r)) \Omega(r) \dd r,
\]
where $\Omega(r)$ is defined by
\[
    \Omega(r) = \frac{r}{\abs{r}} \frac{\omega(r / R)}{R},
\]
with $\omega(\cdot)$ being the function introduced in
\cref{Defn:NonLocalOperator}. In this case
\[
    \int_{0}^R \Omega(r) \dd r = \int_{0}^R \frac{r}{|r|}\frac{\omega(r/R)}{R}
    \dd r = \int_0^1 \omega(\sigma) \dd\sigma = \frac{1}{2}.
\]
Our aim is to develop an asymptotic expansion
$\K_{R}[u]$ as $R \to 0$. In preparation, we compute the moments of the
distribution $\tilde{\Omega}$.
\[
    \mu_n = \int_{-R}^{R} r^n \Omega(r) \dd r
          = \begin{cases}
                0 &\mbox{if $n$ even} \\
                2 R^n \int_{0}^{1} \sigma^n \omega(\sigma) \dd \sigma &\mbox{if $n$ odd}
          \end{cases}.
\]
In case where $R$ is small, we can consider a Taylor expansion of the integrand.
\begin{eqnarray*}
    \K_R[u] &=& \int_{-R}^R h(u(x+r)) \Omega(r) \dd r \\
    &\approx&
    \int_{-R}^R \lb h(u(x)) + r \frac{\partial}{\partial x} h(u(x))+
        \frac{r^2}{2} \frac{\partial^2}{\partial x^2} h(u(x)) + \frac{r^3}{6}
        \frac{\partial^3}{\partial x^3} h(u(x)) + \ldots\rb \Omega(r) \dd r \\
    &=& \sum_{j=0}^\infty \frac{\mu_j r^j}{j!}\frac{\partial ^j}{\partial x^j} h(u(x)) \\
    &=& \mu_1 \frac{\partial}{\partial x} h(u(x)) + \frac{\mu_3}{6} \frac{\partial^3}{\partial x^3} h(u(x)) + \ldots.
\end{eqnarray*}
Hence $\K_R[u]$ for small $R$ acts like a derivative. It should be noted that
the moments $\mu_j$ scale as $R^j$ for $j$ odd, hence, for differentiable
functions $h$ and $u$ we find
\[
    \frac{\K_R[u]}{2R}\to c_1 \frac{\partial}{\partial x}h(u(x)), \quad \mbox{for}\quad R\to 0,
\]
with
\[
    c_1 = \int_0^1\sigma \omega(\sigma) \dd\sigma.
\]
We illustrate this relationship for a given test function in
\cref{Fig:SensingRadiusComparision}, where we plot $c_1 u'(x)$ and $K_R[u](x)$
for several values of $R$. We can obtain a similar asymptotic expansion for the
non-local curvature\index{non-local curvature} (\ie, the derivative of $\K_{R}[u]$ as $R \to 0$)
\[
    \frac{{\K_{R}[u](x)}^{\prime}}{4 R^2}  \approx \frac{\partial^2}{\partial x^2} h(u(x)) + \BigOh{R^2}.
\]
\begin{example}
    Suppose that $\omega(\sigma) = \ifrac{1}{2}$ for $\sigma \in [0, 1]$. Then
    $\Omega(r) =\frac{1}{2R} \frac{r}{|r|}$ and
    \[
        \K_{R}[u](x) = \frac{1}{2R} \int_{-R}^{R} h(u(x+r))\frac{r}{|r|} \dd r
            \approx \frac{R}{2} \frac{\partial}{\partial x} h(u(x)),
    \]
    for small $R$.
\end{example}
\begin{example}
    Consider a singular measure $\omega$ of the form
    \[
        \omega(r) =\delta_R(r), \quad \mbox{then}\quad
            \Omega(r) = \delta_R(r) - \delta_{-R}(r).
    \]
    then
    \[
        \K_R[u] = \frac{1}{2}\lsb h(u(x+R)) - h(u(x-R))\rsb,
    \]
    where we used the convention that a $\delta$-distribution on the domain
    boundary carries weight $\ifrac{1}{2}$.
\end{example}

\begin{figure}\centering
    \input{NonLocalTermExamples.pgf}
    \caption[Comparision of non-local term to first derivative]
    {Comparison of the non-local term to the first derivative of a function
    for several values for the sensing radius
    $R$.}\label{Fig:SensingRadiusComparision}
\end{figure}
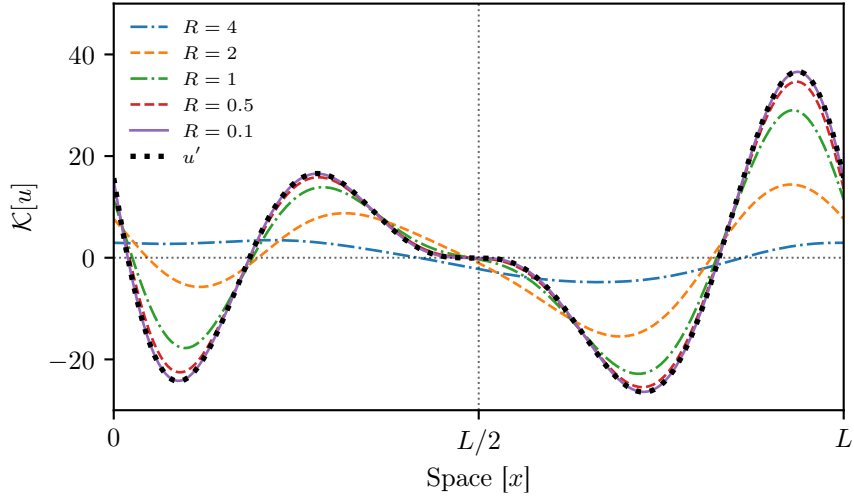

\section{Properties of Steady State Solutions}\label{Chapter:PropertiesOfSolutions}
In this section, we prove several properties of solutions of the steady state
equation
\begin{equation}\label{Eqn:SteadyStates}
    u''(x) = \alpha\Bigl(u(x) \K[u](x)\Bigr)', \quad x\in S^1_L.
\end{equation}
These will be useful later when we carry out the bifurcation analysis.
Throughout this section we assume
(\ref{OmegaAssumption:1}-\ref{OmegaAssumption:5}) and
(\ref{NatureForceAssumption:1}-\ref{NatureForceAssumption:4}).
\begin{lemma}
    \label{Lemma:RightHandSideAvgZero}
    Let $u \in L^{2}(S^1_L)$ be $L$-periodic and the nonlinearity within the
    non-local term is linear, \ie, ($h(u) = u$), then $\Avg[u(x) \K[u](x)] = 0$,
    where $\Avg[\cdot]$ is the averaging operator defined in
    \cref{Defn:AvgOperator}.
\end{lemma}
\begin{proof}
    We proceed by simple calculation, and an application of Fubini's
    theorem.\index{Fubini's theorem}
    \[
    \begin{split}
        \int_{0}^{L} u(x) \K[u](x) \dd x &=
            \int_{0}^{L} u(x) \int_{-1}^{1} u(x + r) \Omega(r) \dd r \dd x \\
            &= \int_{0}^{L} u(x) \int_{0}^{1} \lb u(x + r) - u(x - r) \rb \omega(r) \dd r \dd x \\
            &= \int_{0}^{1} \int_{0}^{L} u(x) \lb u(x + r) - u(x - r) \rb
                \dd x \, \omega(r) \dd r.
    \end{split}
    \]
    Since $u$ is periodic we have
    \[
        \int_{0}^{L} u(x) u(x + r) \dd x = \int_{0}^{L} u(x) u(x - r) \dd x,
    \]
    for all $r \in [0, 1]$ and the result
    follows.
\end{proof}
Note that in the proof of \cref{Lemma:RightHandSideAvgZero} we required that
the nonlinearity within the non-local term is the simple linear function
($h(u) = u$). The previous proof does not work for non-linear $h(u)$. We can
however recover the result by imposing an additional symmetry assumption.
\begin{lemma}\label{Lemma:RightHandSideAvgZeroGeneral}
    Let $u \in L^{2}(S^1_L)$ be $L$-periodic and $u(x) = u(L - x)$
    then the average ${\Avg[u(x) \K[u](x)] = 0}$,
    where $\Avg[\cdot]$ is the averaging operator defined in
    \cref{Defn:AvgOperator} and $h(\cdot)$
    satisfies~\ref{NatureForceAssumption:1}--\ref{NatureForceAssumption:4}.
\end{lemma}
\begin{proof}
    Due to the fact that $u(x) = u(L - x)$ we obtain that
    \[
        u(x) \K[u](x) = - u(L - x) \K[u](L - x).
    \]
    Hence, upon integration we obtain
    \[
    \begin{split}
        \int_{0}^{L} u(x) \K[u](x) \dd x &= \int_{0}^{L/2} u(x) \K[u](x) \dd x
            + \int_{L/2}^{L} u(x) \K[u](x) \dd x \\
            &= \int_{0}^{L/2} u(x) \K[u](x) \dd x
            - \int_{0}^{L/2} u(x) \K[u](x) \dd x = 0.
    \end{split}
    \]
\end{proof}
\begin{remark}\label{Remark:FluxHasZero}
    Note that if $u \in H^2$, then \cref{Lemma:RightHandSideAvgZero} or
    \cref{Lemma:RightHandSideAvgZeroGeneral} imply that the flux defined by
    \[
        J(x) = u^{\prime}(x) - \alpha u(x) \K[u(x)](x),
    \]
    satisfies $\Avg[J] = 0$, and since $J$ is continuous we have that
    $\exists \hat{x} \in S^1_L : J(\hat{x}) = 0$. Note that the continuity
    of $J$ follows as $u \in H^2 \ssubset \C^1$.
\end{remark}
We derive an a priori estimate of positive solutions of
equation~\eqref{Eqn:SteadyStates}. Prior to being able to prove the estimate we
require a technical result.
\begin{lemma}\label{Prop:Reg}
    Let $u \in H^2_{\Bd}$ be a solution of equation~\eqref{Eqn:SteadyStates}.
    Then
    \begin{equation}\label{Eqn:IntegratedMainEquation}
        u^{\prime}(x) = \alpha u(x) \K[u](x).
    \end{equation}
    and
    $u \in \C^3(S^1_L)$.
\end{lemma}
\begin{proof}
    By Sobolev's theorem we have that $u \in \C^{1, 1/2}(S^1_L)$. Integrating
    equation~\eqref{Eqn:SteadyStates} from $\hat{x}$ (the point at which the
    flux is zero, whose existence is guaranteed by \cref{Remark:FluxHasZero}), we
    observe that
    \begin{equation}\label{Eqn:1stInt}
        u^{\prime}(x) = \alpha u(x) \K[u](x).
    \end{equation}
    \cref{Lem:KProp} implies that whenever $u \in H^2$ we
    also have that $\K[u] \in H^2$, hence applying the Banach algebra property
    of $H^2$ and Sobolev's theorem we have that $u \K[u]
    \in H^2 \ssubset \C^{1, 1/2} \subset \C^{1}$. Then by
    equation~\eqref{Eqn:IntegratedMainEquation} we have that $u^{\prime} \in
    \C^{1}(S^1_L)$ and thus $u \in \C^2(S^1_L)$. Further
    from equation~\eqref{Eqn:SteadyStates}, we have that
    \[
        u^{\prime\prime}(x) = \alpha \lb u(x) \K[u](x) \rb^{\prime}.
    \]
    As $u \in \C^2(S^1_L)$ we apply
    \cref{Lem:KProp}~\ref{NonLocalProp:1} to find that
    ${\K[u] \in \C^2(S^1_L)}$, hence ${u \K[u] \in \C^2(S^1_L)}$. This means,
    that $u^{\prime\prime} \in \C^1(S^1_L)$ and finally we have that $u \in
    \C^3(S^1_L)$.
\end{proof}
\begin{lemma}\label{Lemma:AprioriEstimate}
    Let $u \in \C^1(S^1_L)$ be a non-negative solution of
    equation~\eqref{Eqn:SteadyStates}, subject to the integral constraint
    \[
        \Avg[u] = \mean{u},
    \]
    where $\mean{u} > 0$. Then, we have
    \[
        \mean{u} e^{- \alpha \mu L} \leq u(x) \leq \mean{u} e^{\alpha \mu L},
    \]
    where $\mu = C (\mean{u} + L) \abs{\omega}_{\infty}$, and $C$
    from~\ref{NonLocalProp:6}. Further we have that,
    \[
        \abs{u^{\prime}(x)} \leq \alpha \mu \mean{u} e^{\alpha \mu L}.
    \]
    Then,
    \[
        \norm{u}_{\C^1} \leq \lb 1 + \alpha \mu \rb \mean{u} e^{\alpha \mu L}
        \coloneqq \kappa(\alpha, L, \mean{u}, \Omega).
    \]
\end{lemma}
\begin{proof}
    Integrating equation~\eqref{Eqn:SteadyStates} from $\hat{x}$ (given in
    \cref{Remark:FluxHasZero}) to $x$ we obtain equation (\ref{Eqn:IntegratedMainEquation}).
    Thus using \cref{Lem:KProp}~\ref{NonLocalProp:6}
    we obtain the following differential inequality,
    \[
        u^{\prime}(x) \leq \alpha C \abs{\omega}_{\infty} u(x)
            \lb \mean{u} + L \rb.
    \]
    Let us denote $\mu = C (\mean{u} + L) \abs{\omega}_{\infty}$. Next we note that
    if $u(x)$ has mass $\mean{u}$ then there is $\tilde{x} \in [0, L]$ such that
    $u(\tilde{x}) = \mean{u}$. Then integrating from $x$ to $\tilde{x}$, we
    obtain
    \[
        \ln u(\hat{x}) - \ln u(x) \leq \alpha \mu L.
    \]
    Similarly integrating from $\tilde{x}$ to $x$, we obtain
    \[
        \ln u(x) - \ln u(\tilde{x}) \leq \alpha \mu L.
    \]
    Combining both inequalities, we find
    \[
        - \alpha \mu L \leq \ln u(x) - \ln u(\tilde{x}) \leq \alpha \mu L,
    \]
    which yields the a priori estimate.
\end{proof}

Next we show that for steady states $u(x)$ of equation~\eqref{Eqn:SteadyStates},
we find that $u(x)$ has maxima or minima whenever $\K[u](x) = 0$.

\begin{figure}[!ht]\centering
    \input{ApplicationOfNonLocalTerm.pgf}
    \caption[The properties of positive solutions of the non-local equation]
    {Top Left: a typical non-trivial steady state solution. Top Right:
    Examples of the
    non-local term $\K[u]$, and ${\K[u]}^{\prime}$ applied to the
    solution on the left. We observe the properties proven
    in \cref{Lem:ZerosDer} and \cref{Lemma:NonLocalPropertiesMaxMin}.
    We observe that $u'$ and $\K[u]$ have the same positive and negative
    regions. The dashed black line denotes the locations of the zeros of
    $\K[u](x)$ and $u^{\prime}(x)$.
    Bottom Left: The product of $\K[u](x) $ and $u'(x)$ to show that they have the same sign and the same zeroes.
    The dashed black line denotes the locations of the zeros of $\K[u](x)$ and
    $u^{\prime}(x)$.
    Bottom Right: A comparison of the second derivative and ${\K[u]}^{\prime}$ and $u''$. If
    $u$ is convex it is also non-locally convex and if $u$ is non-locally
    concave then it is concave.  The dashed black line denotes the locations of
    the zeros of $\K[u]$ and $u^{\prime}$.
    }\label{Fig:ExamplesOfNonLocalTermApplied}
\end{figure}

\begin{lemma}\label{Lem:ZerosDer}
    Suppose $u(x)$ is a solution of equation~\eqref{Eqn:SteadyStates} and assume
    $u(\hat x) >0$. Then $u^{\prime}(\hat{x}) = 0$ if and only if
    $\K[u](\hat{x}) = 0$ (see \cref{Fig:ExamplesOfNonLocalTermApplied}).
\end{lemma}

\begin{proof}
    Suppose that at $\hat{x} \in S^1_L$ we have $u^{\prime}(\hat{x}) = 0$, then
    equation~\eqref{Eqn:IntegratedMainEquation}, implies that
    \[
        0 = \alpha u(\hat{x}) \K[u](\hat{x}).
    \]
    But both $\alpha \neq 0$ and $u(\hat{x}) \neq 0$, thus $\K[u](\hat{x}) = 0$.
    The other direction follows immediately.
\end{proof}

\begin{lemma}\label{Lemma:NonLocalPropertiesMaxMin}
    Suppose $u(x)$ is a solution of equation~\eqref{Eqn:SteadyStates} and $u(\hat x)>0$,
    then it achieves a non-zero maximum (minimum) at $\hat{x}$ if and only if
    \begin{enumerate}
        \item $\K[u](\hat{x}) = 0$, and
        \item $\lb \K[u](\hat{x}) \rb^{\prime} < (>)\, 0$.
    \end{enumerate}
    See \cref{Fig:ExamplesOfNonLocalTermApplied} for an example.
\end{lemma}
\begin{proof}
    {\bfseries (1)}
    from \cref{Lem:ZerosDer}. \\
    {\bfseries (2)}
    Without loss of generality suppose $u(x)$ achieves a non-zero maximum
    at $\hat{x}$. This means that $u^{\prime\prime}(\hat{x}) < 0$. Thus from
    equation~\eqref{Eqn:SteadyStates} we get that,
    \[
        0 > u^{\prime\prime}(\hat{x}) =
            \alpha u(\hat{x}) \lb \K[u](\hat{x}) \rb^{\prime}.
    \]
    But both $\alpha > 0$, and $u(\hat{x}) > 0$ and thus
    $\lb \K[u] \rb^{\prime} < 0$.
\end{proof}
%
%
\begin{lemma}
    Let $u(x)$ be a positive solution of equation~\eqref{Eqn:SteadyStates} then
    \[
        u^{\prime}(x) \K[u](x) \geq 0.
    \]
    For an example see \cref{Fig:ExamplesOfNonLocalTermApplied}.
\end{lemma}
\begin{proof}
    Substituting equation~\eqref{Eqn:IntegratedMainEquation} into
    equation~\eqref{Eqn:SteadyStates} we obtain,
    \[
        u^{\prime\prime}(x) = \alpha^2 u(x) \lb \K[u](x) \rb^{2} + \alpha u(x)
        \lb \K[u] \rb^{\prime}.
    \]
    Note that the first term on the right hand side is positive and thus we have
    that
    \begin{equation}\label{Eqn:LemmaPositivityProductProof}
        u^{\prime\prime} \geq \alpha u(x) \lb \K[u] \rb^{\prime}.
    \end{equation}
    Then using equation~\eqref{Eqn:SteadyStates} and
    result~\eqref{Eqn:LemmaPositivityProductProof} we obtain that,
    \[
        \alpha u^{\prime}(x) \K[u](x) =
            u^{\prime\prime}(x) - \alpha u(x) \lb \K[u](x) \rb^{\prime} \geq 0.
    \]
\end{proof}

Finally we can relate the local and non-local curvatures for steady state solutions:
\begin{lemma}\label{Lemma:NonLocalConvexity}
    Let $u(x)$ be a solution of equation~\eqref{Eqn:SteadyStates} then
    \begin{enumerate}
        \item If $u^{\prime\prime}(x) \leq 0$, then $\lb \K[u](x) \rb^{\prime}
            \leq 0$,
        \item If $\lb \K[u](x) \rb^{\prime} \geq 0$, then $u^{\prime\prime}(x)
            \geq 0$.
    \end{enumerate}
\end{lemma}
\begin{proof}
    Proof is by equation~\eqref{Eqn:LemmaPositivityProductProof}.
\end{proof}
\section{Summary}
In this section, we studied basic mathematical properties of the non-local
operator $\K[u]$, such as continuity, regularity, growth estimates, and spectral
properties (linear non-local operator). These results are relevant for
subsequent bifurcation analysis.

We saw that the linear non-local operator shares many properties with the first
derivative operator. It is skew-adjoint, it has the complex exponentials as
eigenfunctions with purely complex eigenvalues and it maps sine basis functions
to cosine basis functions and cosine basis functions to negative sine.  The
range of $\K[u]$ is contained in the subspace of zero average functions, and the
null-spaces are all the constant functions. There are however some difference as
well. Most notably $\K[u]$ is a compact operator, thus its eigenvalues
accumulate at zero (see \cref{Example:UniformOmega1}),
while the eigenvalues of ${(\cdot)}^{\prime}$ diverge.

Similarly, the non-local curvature ${\K[u]}^{\prime}$ shares many properties
with a second derivative $(\cdot)''$. Both are  self-adjoint, and both share the
same eigenfunctions (provided we consider $(\cdot)''$ with periodic boundary
conditions as well).  Differences are that ${\K[u]}^{\prime}$ is always a
bounded operator, with bounded eigenvalues.  The eigenvalues of the second
derivative, $-\frac{2 \pi n}{L}$ are modified by a specific factor of $2
M_n(\omega)$, with
\[
    M_n(\omega) = \int_0^1 \sin\left(\frac{2 n \pi r}{L}\right)\omega(r) \dd r,
\]
which are the Fourier-sine coefficients\index{Fourier-sine coefficients} of a
Fourier-sine expansion\index{Fourier-sine expansion} of $\omega(r)$.

In this section, we have seen that the non-local operator $\K[u]$ can be viewed
as a generalization of the local derivative (in the sense that as $R \to 0$,
$\K_{R}[u] \to c_1 u^{\prime}$.

In view of this analogy the estimate ($\abs{\K[u]}_{p} \leq C (\abs{u}_{p} +
L)$) obtained in \cref{Lem:KProp}~\ref{NonLocalProp:3}
can be viewed as a type of reverse Poincar\'e inequality\index{Poincar\'e
inequality}. Intuitively this
estimate \textquote{earns} us an order of differentiability. Similar results for
other non-local operators have been obtained in
\parencites{Ikeda1985}{Ikeda1987}{Hillen2007}.

%% file: DifferentOmegaDistributions.pgf
\begingroup%
\makeatletter%
\begin{pgfpicture}%
\pgfpathrectangle{\pgfpointorigin}{\pgfqpoint{3.922305in}{2.424118in}}%
\pgfusepath{use as bounding box, clip}%
\begin{pgfscope}%
\pgfsetbuttcap%
\pgfsetmiterjoin%
\definecolor{currentfill}{rgb}{1.000000,1.000000,1.000000}%
\pgfsetfillcolor{currentfill}%
\pgfsetlinewidth{0.000000pt}%
\definecolor{currentstroke}{rgb}{1.000000,1.000000,1.000000}%
\pgfsetstrokecolor{currentstroke}%
\pgfsetdash{}{0pt}%
\pgfpathmoveto{\pgfqpoint{0.000000in}{0.000000in}}%
\pgfpathlineto{\pgfqpoint{3.922305in}{0.000000in}}%
\pgfpathlineto{\pgfqpoint{3.922305in}{2.424118in}}%
\pgfpathlineto{\pgfqpoint{0.000000in}{2.424118in}}%
\pgfpathclose%
\pgfusepath{fill}%
\end{pgfscope}%
\begin{pgfscope}%
\pgfsetbuttcap%
\pgfsetmiterjoin%
\definecolor{currentfill}{rgb}{1.000000,1.000000,1.000000}%
\pgfsetfillcolor{currentfill}%
\pgfsetlinewidth{0.000000pt}%
\definecolor{currentstroke}{rgb}{0.000000,0.000000,0.000000}%
\pgfsetstrokecolor{currentstroke}%
\pgfsetstrokeopacity{0.000000}%
\pgfsetdash{}{0pt}%
\pgfpathmoveto{\pgfqpoint{0.316407in}{0.386111in}}%
\pgfpathlineto{\pgfqpoint{3.683902in}{0.386111in}}%
\pgfpathlineto{\pgfqpoint{3.683902in}{2.225507in}}%
\pgfpathlineto{\pgfqpoint{0.316407in}{2.225507in}}%
\pgfpathclose%
\pgfusepath{fill}%
\end{pgfscope}%
\begin{pgfscope}%
\pgfsetbuttcap%
\pgfsetroundjoin%
\definecolor{currentfill}{rgb}{0.000000,0.000000,0.000000}%
\pgfsetfillcolor{currentfill}%
\pgfsetlinewidth{0.803000pt}%
\definecolor{currentstroke}{rgb}{0.000000,0.000000,0.000000}%
\pgfsetstrokecolor{currentstroke}%
\pgfsetdash{}{0pt}%
\pgfsys@defobject{currentmarker}{\pgfqpoint{0.000000in}{-0.048611in}}{\pgfqpoint{0.000000in}{0.000000in}}{%
\pgfpathmoveto{\pgfqpoint{0.000000in}{0.000000in}}%
\pgfpathlineto{\pgfqpoint{0.000000in}{-0.048611in}}%
\pgfusepath{stroke,fill}%
}%
\begin{pgfscope}%
\pgfsys@transformshift{0.316407in}{0.386111in}%
\pgfsys@useobject{currentmarker}{}%
\end{pgfscope}%
\end{pgfscope}%
\begin{pgfscope}%
\pgftext[x=0.316407in,y=0.288889in,,top]{\rmfamily\fontsize{10.000000}{12.000000}\selectfont \(\displaystyle 0.0\)}%
\end{pgfscope}%
\begin{pgfscope}%
\pgfsetbuttcap%
\pgfsetroundjoin%
\definecolor{currentfill}{rgb}{0.000000,0.000000,0.000000}%
\pgfsetfillcolor{currentfill}%
\pgfsetlinewidth{0.803000pt}%
\definecolor{currentstroke}{rgb}{0.000000,0.000000,0.000000}%
\pgfsetstrokecolor{currentstroke}%
\pgfsetdash{}{0pt}%
\pgfsys@defobject{currentmarker}{\pgfqpoint{0.000000in}{-0.048611in}}{\pgfqpoint{0.000000in}{0.000000in}}{%
\pgfpathmoveto{\pgfqpoint{0.000000in}{0.000000in}}%
\pgfpathlineto{\pgfqpoint{0.000000in}{-0.048611in}}%
\pgfusepath{stroke,fill}%
}%
\begin{pgfscope}%
\pgfsys@transformshift{0.989906in}{0.386111in}%
\pgfsys@useobject{currentmarker}{}%
\end{pgfscope}%
\end{pgfscope}%
\begin{pgfscope}%
\pgftext[x=0.989906in,y=0.288889in,,top]{\rmfamily\fontsize{10.000000}{12.000000}\selectfont \(\displaystyle 0.2\)}%
\end{pgfscope}%
\begin{pgfscope}%
\pgfsetbuttcap%
\pgfsetroundjoin%
\definecolor{currentfill}{rgb}{0.000000,0.000000,0.000000}%
\pgfsetfillcolor{currentfill}%
\pgfsetlinewidth{0.803000pt}%
\definecolor{currentstroke}{rgb}{0.000000,0.000000,0.000000}%
\pgfsetstrokecolor{currentstroke}%
\pgfsetdash{}{0pt}%
\pgfsys@defobject{currentmarker}{\pgfqpoint{0.000000in}{-0.048611in}}{\pgfqpoint{0.000000in}{0.000000in}}{%
\pgfpathmoveto{\pgfqpoint{0.000000in}{0.000000in}}%
\pgfpathlineto{\pgfqpoint{0.000000in}{-0.048611in}}%
\pgfusepath{stroke,fill}%
}%
\begin{pgfscope}%
\pgfsys@transformshift{1.663405in}{0.386111in}%
\pgfsys@useobject{currentmarker}{}%
\end{pgfscope}%
\end{pgfscope}%
\begin{pgfscope}%
\pgftext[x=1.663405in,y=0.288889in,,top]{\rmfamily\fontsize{10.000000}{12.000000}\selectfont \(\displaystyle 0.4\)}%
\end{pgfscope}%
\begin{pgfscope}%
\pgfsetbuttcap%
\pgfsetroundjoin%
\definecolor{currentfill}{rgb}{0.000000,0.000000,0.000000}%
\pgfsetfillcolor{currentfill}%
\pgfsetlinewidth{0.803000pt}%
\definecolor{currentstroke}{rgb}{0.000000,0.000000,0.000000}%
\pgfsetstrokecolor{currentstroke}%
\pgfsetdash{}{0pt}%
\pgfsys@defobject{currentmarker}{\pgfqpoint{0.000000in}{-0.048611in}}{\pgfqpoint{0.000000in}{0.000000in}}{%
\pgfpathmoveto{\pgfqpoint{0.000000in}{0.000000in}}%
\pgfpathlineto{\pgfqpoint{0.000000in}{-0.048611in}}%
\pgfusepath{stroke,fill}%
}%
\begin{pgfscope}%
\pgfsys@transformshift{2.336904in}{0.386111in}%
\pgfsys@useobject{currentmarker}{}%
\end{pgfscope}%
\end{pgfscope}%
\begin{pgfscope}%
\pgftext[x=2.336904in,y=0.288889in,,top]{\rmfamily\fontsize{10.000000}{12.000000}\selectfont \(\displaystyle 0.6\)}%
\end{pgfscope}%
\begin{pgfscope}%
\pgfsetbuttcap%
\pgfsetroundjoin%
\definecolor{currentfill}{rgb}{0.000000,0.000000,0.000000}%
\pgfsetfillcolor{currentfill}%
\pgfsetlinewidth{0.803000pt}%
\definecolor{currentstroke}{rgb}{0.000000,0.000000,0.000000}%
\pgfsetstrokecolor{currentstroke}%
\pgfsetdash{}{0pt}%
\pgfsys@defobject{currentmarker}{\pgfqpoint{0.000000in}{-0.048611in}}{\pgfqpoint{0.000000in}{0.000000in}}{%
\pgfpathmoveto{\pgfqpoint{0.000000in}{0.000000in}}%
\pgfpathlineto{\pgfqpoint{0.000000in}{-0.048611in}}%
\pgfusepath{stroke,fill}%
}%
\begin{pgfscope}%
\pgfsys@transformshift{3.010403in}{0.386111in}%
\pgfsys@useobject{currentmarker}{}%
\end{pgfscope}%
\end{pgfscope}%
\begin{pgfscope}%
\pgftext[x=3.010403in,y=0.288889in,,top]{\rmfamily\fontsize{10.000000}{12.000000}\selectfont \(\displaystyle 0.8\)}%
\end{pgfscope}%
\begin{pgfscope}%
\pgfsetbuttcap%
\pgfsetroundjoin%
\definecolor{currentfill}{rgb}{0.000000,0.000000,0.000000}%
\pgfsetfillcolor{currentfill}%
\pgfsetlinewidth{0.803000pt}%
\definecolor{currentstroke}{rgb}{0.000000,0.000000,0.000000}%
\pgfsetstrokecolor{currentstroke}%
\pgfsetdash{}{0pt}%
\pgfsys@defobject{currentmarker}{\pgfqpoint{0.000000in}{-0.048611in}}{\pgfqpoint{0.000000in}{0.000000in}}{%
\pgfpathmoveto{\pgfqpoint{0.000000in}{0.000000in}}%
\pgfpathlineto{\pgfqpoint{0.000000in}{-0.048611in}}%
\pgfusepath{stroke,fill}%
}%
\begin{pgfscope}%
\pgfsys@transformshift{3.683902in}{0.386111in}%
\pgfsys@useobject{currentmarker}{}%
\end{pgfscope}%
\end{pgfscope}%
\begin{pgfscope}%
\pgftext[x=3.683902in,y=0.288889in,,top]{\rmfamily\fontsize{10.000000}{12.000000}\selectfont \(\displaystyle 1.0\)}%
\end{pgfscope}%
\begin{pgfscope}%
\pgfsetbuttcap%
\pgfsetroundjoin%
\definecolor{currentfill}{rgb}{0.000000,0.000000,0.000000}%
\pgfsetfillcolor{currentfill}%
\pgfsetlinewidth{0.803000pt}%
\definecolor{currentstroke}{rgb}{0.000000,0.000000,0.000000}%
\pgfsetstrokecolor{currentstroke}%
\pgfsetdash{}{0pt}%
\pgfsys@defobject{currentmarker}{\pgfqpoint{-0.048611in}{0.000000in}}{\pgfqpoint{0.000000in}{0.000000in}}{%
\pgfpathmoveto{\pgfqpoint{0.000000in}{0.000000in}}%
\pgfpathlineto{\pgfqpoint{-0.048611in}{0.000000in}}%
\pgfusepath{stroke,fill}%
}%
\begin{pgfscope}%
\pgfsys@transformshift{0.316407in}{0.386111in}%
\pgfsys@useobject{currentmarker}{}%
\end{pgfscope}%
\end{pgfscope}%
\begin{pgfscope}%
\pgftext[x=0.149740in,y=0.338283in,left,base]{\rmfamily\fontsize{10.000000}{12.000000}\selectfont \(\displaystyle 0\)}%
\end{pgfscope}%
\begin{pgfscope}%
\pgfsetbuttcap%
\pgfsetroundjoin%
\definecolor{currentfill}{rgb}{0.000000,0.000000,0.000000}%
\pgfsetfillcolor{currentfill}%
\pgfsetlinewidth{0.803000pt}%
\definecolor{currentstroke}{rgb}{0.000000,0.000000,0.000000}%
\pgfsetstrokecolor{currentstroke}%
\pgfsetdash{}{0pt}%
\pgfsys@defobject{currentmarker}{\pgfqpoint{-0.048611in}{0.000000in}}{\pgfqpoint{0.000000in}{0.000000in}}{%
\pgfpathmoveto{\pgfqpoint{0.000000in}{0.000000in}}%
\pgfpathlineto{\pgfqpoint{-0.048611in}{0.000000in}}%
\pgfusepath{stroke,fill}%
}%
\begin{pgfscope}%
\pgfsys@transformshift{0.316407in}{0.736472in}%
\pgfsys@useobject{currentmarker}{}%
\end{pgfscope}%
\end{pgfscope}%
\begin{pgfscope}%
\pgftext[x=0.149740in,y=0.688644in,left,base]{\rmfamily\fontsize{10.000000}{12.000000}\selectfont \(\displaystyle 1\)}%
\end{pgfscope}%
\begin{pgfscope}%
\pgfsetbuttcap%
\pgfsetroundjoin%
\definecolor{currentfill}{rgb}{0.000000,0.000000,0.000000}%
\pgfsetfillcolor{currentfill}%
\pgfsetlinewidth{0.803000pt}%
\definecolor{currentstroke}{rgb}{0.000000,0.000000,0.000000}%
\pgfsetstrokecolor{currentstroke}%
\pgfsetdash{}{0pt}%
\pgfsys@defobject{currentmarker}{\pgfqpoint{-0.048611in}{0.000000in}}{\pgfqpoint{0.000000in}{0.000000in}}{%
\pgfpathmoveto{\pgfqpoint{0.000000in}{0.000000in}}%
\pgfpathlineto{\pgfqpoint{-0.048611in}{0.000000in}}%
\pgfusepath{stroke,fill}%
}%
\begin{pgfscope}%
\pgfsys@transformshift{0.316407in}{1.086833in}%
\pgfsys@useobject{currentmarker}{}%
\end{pgfscope}%
\end{pgfscope}%
\begin{pgfscope}%
\pgftext[x=0.149740in,y=1.039006in,left,base]{\rmfamily\fontsize{10.000000}{12.000000}\selectfont \(\displaystyle 2\)}%
\end{pgfscope}%
\begin{pgfscope}%
\pgfsetbuttcap%
\pgfsetroundjoin%
\definecolor{currentfill}{rgb}{0.000000,0.000000,0.000000}%
\pgfsetfillcolor{currentfill}%
\pgfsetlinewidth{0.803000pt}%
\definecolor{currentstroke}{rgb}{0.000000,0.000000,0.000000}%
\pgfsetstrokecolor{currentstroke}%
\pgfsetdash{}{0pt}%
\pgfsys@defobject{currentmarker}{\pgfqpoint{-0.048611in}{0.000000in}}{\pgfqpoint{0.000000in}{0.000000in}}{%
\pgfpathmoveto{\pgfqpoint{0.000000in}{0.000000in}}%
\pgfpathlineto{\pgfqpoint{-0.048611in}{0.000000in}}%
\pgfusepath{stroke,fill}%
}%
\begin{pgfscope}%
\pgfsys@transformshift{0.316407in}{1.437194in}%
\pgfsys@useobject{currentmarker}{}%
\end{pgfscope}%
\end{pgfscope}%
\begin{pgfscope}%
\pgftext[x=0.149740in,y=1.389367in,left,base]{\rmfamily\fontsize{10.000000}{12.000000}\selectfont \(\displaystyle 3\)}%
\end{pgfscope}%
\begin{pgfscope}%
\pgfsetbuttcap%
\pgfsetroundjoin%
\definecolor{currentfill}{rgb}{0.000000,0.000000,0.000000}%
\pgfsetfillcolor{currentfill}%
\pgfsetlinewidth{0.803000pt}%
\definecolor{currentstroke}{rgb}{0.000000,0.000000,0.000000}%
\pgfsetstrokecolor{currentstroke}%
\pgfsetdash{}{0pt}%
\pgfsys@defobject{currentmarker}{\pgfqpoint{-0.048611in}{0.000000in}}{\pgfqpoint{0.000000in}{0.000000in}}{%
\pgfpathmoveto{\pgfqpoint{0.000000in}{0.000000in}}%
\pgfpathlineto{\pgfqpoint{-0.048611in}{0.000000in}}%
\pgfusepath{stroke,fill}%
}%
\begin{pgfscope}%
\pgfsys@transformshift{0.316407in}{1.787555in}%
\pgfsys@useobject{currentmarker}{}%
\end{pgfscope}%
\end{pgfscope}%
\begin{pgfscope}%
\pgftext[x=0.149740in,y=1.739728in,left,base]{\rmfamily\fontsize{10.000000}{12.000000}\selectfont \(\displaystyle 4\)}%
\end{pgfscope}%
\begin{pgfscope}%
\pgfsetbuttcap%
\pgfsetroundjoin%
\definecolor{currentfill}{rgb}{0.000000,0.000000,0.000000}%
\pgfsetfillcolor{currentfill}%
\pgfsetlinewidth{0.803000pt}%
\definecolor{currentstroke}{rgb}{0.000000,0.000000,0.000000}%
\pgfsetstrokecolor{currentstroke}%
\pgfsetdash{}{0pt}%
\pgfsys@defobject{currentmarker}{\pgfqpoint{-0.048611in}{0.000000in}}{\pgfqpoint{0.000000in}{0.000000in}}{%
\pgfpathmoveto{\pgfqpoint{0.000000in}{0.000000in}}%
\pgfpathlineto{\pgfqpoint{-0.048611in}{0.000000in}}%
\pgfusepath{stroke,fill}%
}%
\begin{pgfscope}%
\pgfsys@transformshift{0.316407in}{2.137917in}%
\pgfsys@useobject{currentmarker}{}%
\end{pgfscope}%
\end{pgfscope}%
\begin{pgfscope}%
\pgftext[x=0.149740in,y=2.090089in,left,base]{\rmfamily\fontsize{10.000000}{12.000000}\selectfont \(\displaystyle 5\)}%
\end{pgfscope}%
\begin{pgfscope}%
\pgfpathrectangle{\pgfqpoint{0.316407in}{0.386111in}}{\pgfqpoint{3.367494in}{1.839396in}}%
\pgfusepath{clip}%
\pgfsetrectcap%
\pgfsetroundjoin%
\pgfsetlinewidth{1.003750pt}%
\definecolor{currentstroke}{rgb}{0.000000,0.000000,0.000000}%
\pgfsetstrokecolor{currentstroke}%
\pgfsetdash{}{0pt}%
\pgfpathmoveto{\pgfqpoint{0.314723in}{0.561292in}}%
\pgfpathlineto{\pgfqpoint{3.683902in}{0.561292in}}%
\pgfpathlineto{\pgfqpoint{3.683902in}{0.561292in}}%
\pgfusepath{stroke}%
\end{pgfscope}%
\begin{pgfscope}%
\pgfpathrectangle{\pgfqpoint{0.316407in}{0.386111in}}{\pgfqpoint{3.367494in}{1.839396in}}%
\pgfusepath{clip}%
\pgfsetbuttcap%
\pgfsetroundjoin%
\pgfsetlinewidth{1.003750pt}%
\definecolor{currentstroke}{rgb}{0.000000,0.000000,0.000000}%
\pgfsetstrokecolor{currentstroke}%
\pgfsetdash{{6.400000pt}{1.600000pt}{1.000000pt}{1.600000pt}}{0.000000pt}%
\pgfpathmoveto{\pgfqpoint{0.314723in}{2.145396in}}%
\pgfpathlineto{\pgfqpoint{0.318092in}{2.145396in}}%
\pgfpathlineto{\pgfqpoint{0.365261in}{2.026400in}}%
\pgfpathlineto{\pgfqpoint{0.412429in}{1.915453in}}%
\pgfpathlineto{\pgfqpoint{0.459598in}{1.812010in}}%
\pgfpathlineto{\pgfqpoint{0.506766in}{1.715564in}}%
\pgfpathlineto{\pgfqpoint{0.557304in}{1.619456in}}%
\pgfpathlineto{\pgfqpoint{0.607841in}{1.530296in}}%
\pgfpathlineto{\pgfqpoint{0.658379in}{1.447581in}}%
\pgfpathlineto{\pgfqpoint{0.708917in}{1.370846in}}%
\pgfpathlineto{\pgfqpoint{0.759454in}{1.299658in}}%
\pgfpathlineto{\pgfqpoint{0.809992in}{1.233617in}}%
\pgfpathlineto{\pgfqpoint{0.863899in}{1.168426in}}%
\pgfpathlineto{\pgfqpoint{0.917806in}{1.108250in}}%
\pgfpathlineto{\pgfqpoint{0.971713in}{1.052703in}}%
\pgfpathlineto{\pgfqpoint{1.025620in}{1.001428in}}%
\pgfpathlineto{\pgfqpoint{1.082896in}{0.951263in}}%
\pgfpathlineto{\pgfqpoint{1.140172in}{0.905188in}}%
\pgfpathlineto{\pgfqpoint{1.200817in}{0.860491in}}%
\pgfpathlineto{\pgfqpoint{1.261462in}{0.819642in}}%
\pgfpathlineto{\pgfqpoint{1.325477in}{0.780333in}}%
\pgfpathlineto{\pgfqpoint{1.389491in}{0.744589in}}%
\pgfpathlineto{\pgfqpoint{1.456875in}{0.710459in}}%
\pgfpathlineto{\pgfqpoint{1.527627in}{0.678114in}}%
\pgfpathlineto{\pgfqpoint{1.601749in}{0.647683in}}%
\pgfpathlineto{\pgfqpoint{1.679240in}{0.619254in}}%
\pgfpathlineto{\pgfqpoint{1.760101in}{0.592878in}}%
\pgfpathlineto{\pgfqpoint{1.847699in}{0.567660in}}%
\pgfpathlineto{\pgfqpoint{1.938667in}{0.544723in}}%
\pgfpathlineto{\pgfqpoint{2.036373in}{0.523304in}}%
\pgfpathlineto{\pgfqpoint{2.140818in}{0.503596in}}%
\pgfpathlineto{\pgfqpoint{2.255370in}{0.485221in}}%
\pgfpathlineto{\pgfqpoint{2.380030in}{0.468474in}}%
\pgfpathlineto{\pgfqpoint{2.514797in}{0.453537in}}%
\pgfpathlineto{\pgfqpoint{2.666410in}{0.439946in}}%
\pgfpathlineto{\pgfqpoint{2.834869in}{0.428032in}}%
\pgfpathlineto{\pgfqpoint{3.026912in}{0.417632in}}%
\pgfpathlineto{\pgfqpoint{3.249278in}{0.408769in}}%
\pgfpathlineto{\pgfqpoint{3.515443in}{0.401372in}}%
\pgfpathlineto{\pgfqpoint{3.683902in}{0.397995in}}%
\pgfpathlineto{\pgfqpoint{3.683902in}{0.397995in}}%
\pgfusepath{stroke}%
\end{pgfscope}%
\begin{pgfscope}%
\pgfpathrectangle{\pgfqpoint{0.316407in}{0.386111in}}{\pgfqpoint{3.367494in}{1.839396in}}%
\pgfusepath{clip}%
\pgfsetbuttcap%
\pgfsetroundjoin%
\pgfsetlinewidth{1.003750pt}%
\definecolor{currentstroke}{rgb}{0.000000,0.000000,0.000000}%
\pgfsetstrokecolor{currentstroke}%
\pgfsetdash{{3.700000pt}{1.600000pt}}{0.000000pt}%
\pgfpathmoveto{\pgfqpoint{0.314723in}{0.381729in}}%
\pgfpathlineto{\pgfqpoint{0.429275in}{0.675592in}}%
\pgfpathlineto{\pgfqpoint{0.489920in}{0.822696in}}%
\pgfpathlineto{\pgfqpoint{0.540458in}{0.937509in}}%
\pgfpathlineto{\pgfqpoint{0.584257in}{1.029830in}}%
\pgfpathlineto{\pgfqpoint{0.624687in}{1.108215in}}%
\pgfpathlineto{\pgfqpoint{0.661748in}{1.173712in}}%
\pgfpathlineto{\pgfqpoint{0.695440in}{1.227606in}}%
\pgfpathlineto{\pgfqpoint{0.729132in}{1.275854in}}%
\pgfpathlineto{\pgfqpoint{0.759454in}{1.314291in}}%
\pgfpathlineto{\pgfqpoint{0.789777in}{1.347900in}}%
\pgfpathlineto{\pgfqpoint{0.820100in}{1.376627in}}%
\pgfpathlineto{\pgfqpoint{0.847053in}{1.398051in}}%
\pgfpathlineto{\pgfqpoint{0.874007in}{1.415621in}}%
\pgfpathlineto{\pgfqpoint{0.900960in}{1.429376in}}%
\pgfpathlineto{\pgfqpoint{0.927913in}{1.439379in}}%
\pgfpathlineto{\pgfqpoint{0.954867in}{1.445714in}}%
\pgfpathlineto{\pgfqpoint{0.981820in}{1.448484in}}%
\pgfpathlineto{\pgfqpoint{1.008774in}{1.447812in}}%
\pgfpathlineto{\pgfqpoint{1.035727in}{1.443837in}}%
\pgfpathlineto{\pgfqpoint{1.062681in}{1.436713in}}%
\pgfpathlineto{\pgfqpoint{1.089634in}{1.426608in}}%
\pgfpathlineto{\pgfqpoint{1.119957in}{1.411900in}}%
\pgfpathlineto{\pgfqpoint{1.153648in}{1.391726in}}%
\pgfpathlineto{\pgfqpoint{1.187340in}{1.367893in}}%
\pgfpathlineto{\pgfqpoint{1.224401in}{1.337925in}}%
\pgfpathlineto{\pgfqpoint{1.264831in}{1.301364in}}%
\pgfpathlineto{\pgfqpoint{1.312000in}{1.254519in}}%
\pgfpathlineto{\pgfqpoint{1.369276in}{1.193066in}}%
\pgfpathlineto{\pgfqpoint{1.450136in}{1.101171in}}%
\pgfpathlineto{\pgfqpoint{1.618595in}{0.908584in}}%
\pgfpathlineto{\pgfqpoint{1.689348in}{0.833241in}}%
\pgfpathlineto{\pgfqpoint{1.749993in}{0.773125in}}%
\pgfpathlineto{\pgfqpoint{1.807269in}{0.720733in}}%
\pgfpathlineto{\pgfqpoint{1.861176in}{0.675607in}}%
\pgfpathlineto{\pgfqpoint{1.915083in}{0.634664in}}%
\pgfpathlineto{\pgfqpoint{1.965621in}{0.600083in}}%
\pgfpathlineto{\pgfqpoint{2.016158in}{0.569108in}}%
\pgfpathlineto{\pgfqpoint{2.066696in}{0.541600in}}%
\pgfpathlineto{\pgfqpoint{2.120603in}{0.515873in}}%
\pgfpathlineto{\pgfqpoint{2.174510in}{0.493614in}}%
\pgfpathlineto{\pgfqpoint{2.231786in}{0.473437in}}%
\pgfpathlineto{\pgfqpoint{2.289062in}{0.456472in}}%
\pgfpathlineto{\pgfqpoint{2.353076in}{0.440855in}}%
\pgfpathlineto{\pgfqpoint{2.420460in}{0.427693in}}%
\pgfpathlineto{\pgfqpoint{2.494582in}{0.416449in}}%
\pgfpathlineto{\pgfqpoint{2.578811in}{0.406976in}}%
\pgfpathlineto{\pgfqpoint{2.673148in}{0.399556in}}%
\pgfpathlineto{\pgfqpoint{2.787700in}{0.393774in}}%
\pgfpathlineto{\pgfqpoint{2.932575in}{0.389711in}}%
\pgfpathlineto{\pgfqpoint{3.141464in}{0.387222in}}%
\pgfpathlineto{\pgfqpoint{3.535658in}{0.386203in}}%
\pgfpathlineto{\pgfqpoint{3.683902in}{0.386144in}}%
\pgfpathlineto{\pgfqpoint{3.683902in}{0.386144in}}%
\pgfusepath{stroke}%
\end{pgfscope}%
\begin{pgfscope}%
\pgfsetrectcap%
\pgfsetmiterjoin%
\pgfsetlinewidth{0.803000pt}%
\definecolor{currentstroke}{rgb}{0.000000,0.000000,0.000000}%
\pgfsetstrokecolor{currentstroke}%
\pgfsetdash{}{0pt}%
\pgfpathmoveto{\pgfqpoint{0.316407in}{0.386111in}}%
\pgfpathlineto{\pgfqpoint{0.316407in}{2.225507in}}%
\pgfusepath{stroke}%
\end{pgfscope}%
\begin{pgfscope}%
\pgfsetrectcap%
\pgfsetmiterjoin%
\pgfsetlinewidth{0.803000pt}%
\definecolor{currentstroke}{rgb}{0.000000,0.000000,0.000000}%
\pgfsetstrokecolor{currentstroke}%
\pgfsetdash{}{0pt}%
\pgfpathmoveto{\pgfqpoint{3.683902in}{0.386111in}}%
\pgfpathlineto{\pgfqpoint{3.683902in}{2.225507in}}%
\pgfusepath{stroke}%
\end{pgfscope}%
\begin{pgfscope}%
\pgfsetrectcap%
\pgfsetmiterjoin%
\pgfsetlinewidth{0.803000pt}%
\definecolor{currentstroke}{rgb}{0.000000,0.000000,0.000000}%
\pgfsetstrokecolor{currentstroke}%
\pgfsetdash{}{0pt}%
\pgfpathmoveto{\pgfqpoint{0.316407in}{0.386111in}}%
\pgfpathlineto{\pgfqpoint{3.683902in}{0.386111in}}%
\pgfusepath{stroke}%
\end{pgfscope}%
\begin{pgfscope}%
\pgfsetrectcap%
\pgfsetmiterjoin%
\pgfsetlinewidth{0.803000pt}%
\definecolor{currentstroke}{rgb}{0.000000,0.000000,0.000000}%
\pgfsetstrokecolor{currentstroke}%
\pgfsetdash{}{0pt}%
\pgfpathmoveto{\pgfqpoint{0.316407in}{2.225507in}}%
\pgfpathlineto{\pgfqpoint{3.683902in}{2.225507in}}%
\pgfusepath{stroke}%
\end{pgfscope}%
\begin{pgfscope}%
\pgfsetrectcap%
\pgfsetroundjoin%
\pgfsetlinewidth{1.003750pt}%
\definecolor{currentstroke}{rgb}{0.000000,0.000000,0.000000}%
\pgfsetstrokecolor{currentstroke}%
\pgfsetdash{}{0pt}%
\pgfpathmoveto{\pgfqpoint{2.439096in}{2.051896in}}%
\pgfpathlineto{\pgfqpoint{2.716874in}{2.051896in}}%
\pgfusepath{stroke}%
\end{pgfscope}%
\begin{pgfscope}%
\pgftext[x=2.827985in,y=2.003285in,left,base]{\rmfamily\fontsize{10.000000}{12.000000}\selectfont Uniform}%
\end{pgfscope}%
\begin{pgfscope}%
\pgfsetbuttcap%
\pgfsetroundjoin%
\pgfsetlinewidth{1.003750pt}%
\definecolor{currentstroke}{rgb}{0.000000,0.000000,0.000000}%
\pgfsetstrokecolor{currentstroke}%
\pgfsetdash{{6.400000pt}{1.600000pt}{1.000000pt}{1.600000pt}}{0.000000pt}%
\pgfpathmoveto{\pgfqpoint{2.439096in}{1.858229in}}%
\pgfpathlineto{\pgfqpoint{2.716874in}{1.858229in}}%
\pgfusepath{stroke}%
\end{pgfscope}%
\begin{pgfscope}%
\pgftext[x=2.827985in,y=1.809618in,left,base]{\rmfamily\fontsize{10.000000}{12.000000}\selectfont Exponential}%
\end{pgfscope}%
\begin{pgfscope}%
\pgfsetbuttcap%
\pgfsetroundjoin%
\pgfsetlinewidth{1.003750pt}%
\definecolor{currentstroke}{rgb}{0.000000,0.000000,0.000000}%
\pgfsetstrokecolor{currentstroke}%
\pgfsetdash{{3.700000pt}{1.600000pt}}{0.000000pt}%
\pgfpathmoveto{\pgfqpoint{2.439096in}{1.664563in}}%
\pgfpathlineto{\pgfqpoint{2.716874in}{1.664563in}}%
\pgfusepath{stroke}%
\end{pgfscope}%
\begin{pgfscope}%
\pgftext[x=2.827985in,y=1.615952in,left,base]{\rmfamily\fontsize{10.000000}{12.000000}\selectfont O3}%
\end{pgfscope}%
\end{pgfpicture}%
\makeatother%
\endgroup%

%% file: NonLocalTermExamples.pgf
\begingroup%
\makeatletter%
\begin{pgfpicture}%
\pgfpathrectangle{\pgfpointorigin}{\pgfqpoint{4.706766in}{2.908942in}}%
\pgfusepath{use as bounding box, clip}%
\begin{pgfscope}%
\pgfsetbuttcap%
\pgfsetmiterjoin%
\definecolor{currentfill}{rgb}{1.000000,1.000000,1.000000}%
\pgfsetfillcolor{currentfill}%
\pgfsetlinewidth{0.000000pt}%
\definecolor{currentstroke}{rgb}{1.000000,1.000000,1.000000}%
\pgfsetstrokecolor{currentstroke}%
\pgfsetdash{}{0pt}%
\pgfpathmoveto{\pgfqpoint{0.000000in}{0.000000in}}%
\pgfpathlineto{\pgfqpoint{4.706766in}{0.000000in}}%
\pgfpathlineto{\pgfqpoint{4.706766in}{2.908942in}}%
\pgfpathlineto{\pgfqpoint{0.000000in}{2.908942in}}%
\pgfpathclose%
\pgfusepath{fill}%
\end{pgfscope}%
\begin{pgfscope}%
\pgfsetbuttcap%
\pgfsetmiterjoin%
\definecolor{currentfill}{rgb}{1.000000,1.000000,1.000000}%
\pgfsetfillcolor{currentfill}%
\pgfsetlinewidth{0.000000pt}%
\definecolor{currentstroke}{rgb}{0.000000,0.000000,0.000000}%
\pgfsetstrokecolor{currentstroke}%
\pgfsetstrokeopacity{0.000000}%
\pgfsetdash{}{0pt}%
\pgfpathmoveto{\pgfqpoint{0.687659in}{0.580556in}}%
\pgfpathlineto{\pgfqpoint{4.508155in}{0.580556in}}%
\pgfpathlineto{\pgfqpoint{4.508155in}{2.710330in}}%
\pgfpathlineto{\pgfqpoint{0.687659in}{2.710330in}}%
\pgfpathclose%
\pgfusepath{fill}%
\end{pgfscope}%
\begin{pgfscope}%
\pgfsetbuttcap%
\pgfsetroundjoin%
\definecolor{currentfill}{rgb}{0.000000,0.000000,0.000000}%
\pgfsetfillcolor{currentfill}%
\pgfsetlinewidth{0.803000pt}%
\definecolor{currentstroke}{rgb}{0.000000,0.000000,0.000000}%
\pgfsetstrokecolor{currentstroke}%
\pgfsetdash{}{0pt}%
\pgfsys@defobject{currentmarker}{\pgfqpoint{0.000000in}{-0.048611in}}{\pgfqpoint{0.000000in}{0.000000in}}{%
\pgfpathmoveto{\pgfqpoint{0.000000in}{0.000000in}}%
\pgfpathlineto{\pgfqpoint{0.000000in}{-0.048611in}}%
\pgfusepath{stroke,fill}%
}%
\begin{pgfscope}%
\pgfsys@transformshift{0.687659in}{0.580556in}%
\pgfsys@useobject{currentmarker}{}%
\end{pgfscope}%
\end{pgfscope}%
\begin{pgfscope}%
\pgftext[x=0.687659in,y=0.483333in,,top]{\rmfamily\fontsize{10.000000}{12.000000}\selectfont \(\displaystyle 0\)}%
\end{pgfscope}%
\begin{pgfscope}%
\pgfsetbuttcap%
\pgfsetroundjoin%
\definecolor{currentfill}{rgb}{0.000000,0.000000,0.000000}%
\pgfsetfillcolor{currentfill}%
\pgfsetlinewidth{0.803000pt}%
\definecolor{currentstroke}{rgb}{0.000000,0.000000,0.000000}%
\pgfsetstrokecolor{currentstroke}%
\pgfsetdash{}{0pt}%
\pgfsys@defobject{currentmarker}{\pgfqpoint{0.000000in}{-0.048611in}}{\pgfqpoint{0.000000in}{0.000000in}}{%
\pgfpathmoveto{\pgfqpoint{0.000000in}{0.000000in}}%
\pgfpathlineto{\pgfqpoint{0.000000in}{-0.048611in}}%
\pgfusepath{stroke,fill}%
}%
\begin{pgfscope}%
\pgfsys@transformshift{2.597907in}{0.580556in}%
\pgfsys@useobject{currentmarker}{}%
\end{pgfscope}%
\end{pgfscope}%
\begin{pgfscope}%
\pgftext[x=2.597907in,y=0.483333in,,top]{\rmfamily\fontsize{10.000000}{12.000000}\selectfont \(\displaystyle L/2\)}%
\end{pgfscope}%
\begin{pgfscope}%
\pgfsetbuttcap%
\pgfsetroundjoin%
\definecolor{currentfill}{rgb}{0.000000,0.000000,0.000000}%
\pgfsetfillcolor{currentfill}%
\pgfsetlinewidth{0.803000pt}%
\definecolor{currentstroke}{rgb}{0.000000,0.000000,0.000000}%
\pgfsetstrokecolor{currentstroke}%
\pgfsetdash{}{0pt}%
\pgfsys@defobject{currentmarker}{\pgfqpoint{0.000000in}{-0.048611in}}{\pgfqpoint{0.000000in}{0.000000in}}{%
\pgfpathmoveto{\pgfqpoint{0.000000in}{0.000000in}}%
\pgfpathlineto{\pgfqpoint{0.000000in}{-0.048611in}}%
\pgfusepath{stroke,fill}%
}%
\begin{pgfscope}%
\pgfsys@transformshift{4.508155in}{0.580556in}%
\pgfsys@useobject{currentmarker}{}%
\end{pgfscope}%
\end{pgfscope}%
\begin{pgfscope}%
\pgftext[x=4.508155in,y=0.483333in,,top]{\rmfamily\fontsize{10.000000}{12.000000}\selectfont \(\displaystyle L\)}%
\end{pgfscope}%
\begin{pgfscope}%
\pgftext[x=2.597907in,y=0.288889in,,top]{\rmfamily\fontsize{10.000000}{12.000000}\selectfont Space [\(\displaystyle x\)]}%
\end{pgfscope}%
\begin{pgfscope}%
\pgfsetbuttcap%
\pgfsetroundjoin%
\definecolor{currentfill}{rgb}{0.000000,0.000000,0.000000}%
\pgfsetfillcolor{currentfill}%
\pgfsetlinewidth{0.803000pt}%
\definecolor{currentstroke}{rgb}{0.000000,0.000000,0.000000}%
\pgfsetstrokecolor{currentstroke}%
\pgfsetdash{}{0pt}%
\pgfsys@defobject{currentmarker}{\pgfqpoint{-0.048611in}{0.000000in}}{\pgfqpoint{0.000000in}{0.000000in}}{%
\pgfpathmoveto{\pgfqpoint{0.000000in}{0.000000in}}%
\pgfpathlineto{\pgfqpoint{-0.048611in}{0.000000in}}%
\pgfusepath{stroke,fill}%
}%
\begin{pgfscope}%
\pgfsys@transformshift{0.687659in}{0.846777in}%
\pgfsys@useobject{currentmarker}{}%
\end{pgfscope}%
\end{pgfscope}%
\begin{pgfscope}%
\pgftext[x=0.343522in,y=0.798950in,left,base]{\rmfamily\fontsize{10.000000}{12.000000}\selectfont \(\displaystyle -20\)}%
\end{pgfscope}%
\begin{pgfscope}%
\pgfsetbuttcap%
\pgfsetroundjoin%
\definecolor{currentfill}{rgb}{0.000000,0.000000,0.000000}%
\pgfsetfillcolor{currentfill}%
\pgfsetlinewidth{0.803000pt}%
\definecolor{currentstroke}{rgb}{0.000000,0.000000,0.000000}%
\pgfsetstrokecolor{currentstroke}%
\pgfsetdash{}{0pt}%
\pgfsys@defobject{currentmarker}{\pgfqpoint{-0.048611in}{0.000000in}}{\pgfqpoint{0.000000in}{0.000000in}}{%
\pgfpathmoveto{\pgfqpoint{0.000000in}{0.000000in}}%
\pgfpathlineto{\pgfqpoint{-0.048611in}{0.000000in}}%
\pgfusepath{stroke,fill}%
}%
\begin{pgfscope}%
\pgfsys@transformshift{0.687659in}{1.379221in}%
\pgfsys@useobject{currentmarker}{}%
\end{pgfscope}%
\end{pgfscope}%
\begin{pgfscope}%
\pgftext[x=0.520992in,y=1.331393in,left,base]{\rmfamily\fontsize{10.000000}{12.000000}\selectfont \(\displaystyle 0\)}%
\end{pgfscope}%
\begin{pgfscope}%
\pgfsetbuttcap%
\pgfsetroundjoin%
\definecolor{currentfill}{rgb}{0.000000,0.000000,0.000000}%
\pgfsetfillcolor{currentfill}%
\pgfsetlinewidth{0.803000pt}%
\definecolor{currentstroke}{rgb}{0.000000,0.000000,0.000000}%
\pgfsetstrokecolor{currentstroke}%
\pgfsetdash{}{0pt}%
\pgfsys@defobject{currentmarker}{\pgfqpoint{-0.048611in}{0.000000in}}{\pgfqpoint{0.000000in}{0.000000in}}{%
\pgfpathmoveto{\pgfqpoint{0.000000in}{0.000000in}}%
\pgfpathlineto{\pgfqpoint{-0.048611in}{0.000000in}}%
\pgfusepath{stroke,fill}%
}%
\begin{pgfscope}%
\pgfsys@transformshift{0.687659in}{1.911665in}%
\pgfsys@useobject{currentmarker}{}%
\end{pgfscope}%
\end{pgfscope}%
\begin{pgfscope}%
\pgftext[x=0.451547in,y=1.863837in,left,base]{\rmfamily\fontsize{10.000000}{12.000000}\selectfont \(\displaystyle 20\)}%
\end{pgfscope}%
\begin{pgfscope}%
\pgfsetbuttcap%
\pgfsetroundjoin%
\definecolor{currentfill}{rgb}{0.000000,0.000000,0.000000}%
\pgfsetfillcolor{currentfill}%
\pgfsetlinewidth{0.803000pt}%
\definecolor{currentstroke}{rgb}{0.000000,0.000000,0.000000}%
\pgfsetstrokecolor{currentstroke}%
\pgfsetdash{}{0pt}%
\pgfsys@defobject{currentmarker}{\pgfqpoint{-0.048611in}{0.000000in}}{\pgfqpoint{0.000000in}{0.000000in}}{%
\pgfpathmoveto{\pgfqpoint{0.000000in}{0.000000in}}%
\pgfpathlineto{\pgfqpoint{-0.048611in}{0.000000in}}%
\pgfusepath{stroke,fill}%
}%
\begin{pgfscope}%
\pgfsys@transformshift{0.687659in}{2.444109in}%
\pgfsys@useobject{currentmarker}{}%
\end{pgfscope}%
\end{pgfscope}%
\begin{pgfscope}%
\pgftext[x=0.451547in,y=2.396281in,left,base]{\rmfamily\fontsize{10.000000}{12.000000}\selectfont \(\displaystyle 40\)}%
\end{pgfscope}%
\begin{pgfscope}%
\pgftext[x=0.287966in,y=1.645443in,,bottom,rotate=90.000000]{\rmfamily\fontsize{10.000000}{12.000000}\selectfont \(\displaystyle \mathcal{K}[u]\)}%
\end{pgfscope}%
\begin{pgfscope}%
\pgfpathrectangle{\pgfqpoint{0.687659in}{0.580556in}}{\pgfqpoint{3.820497in}{2.129775in}}%
\pgfusepath{clip}%
\pgfsetbuttcap%
\pgfsetroundjoin%
\pgfsetlinewidth{1.003750pt}%
\definecolor{currentstroke}{rgb}{0.121569,0.466667,0.705882}%
\pgfsetstrokecolor{currentstroke}%
\pgfsetdash{{6.400000pt}{1.600000pt}{1.000000pt}{1.600000pt}}{0.000000pt}%
\pgfpathmoveto{\pgfqpoint{0.687659in}{1.457460in}}%
\pgfpathlineto{\pgfqpoint{0.802667in}{1.454518in}}%
\pgfpathlineto{\pgfqpoint{0.947478in}{1.451606in}}%
\pgfpathlineto{\pgfqpoint{1.051024in}{1.452777in}}%
\pgfpathlineto{\pgfqpoint{1.165650in}{1.457312in}}%
\pgfpathlineto{\pgfqpoint{1.469792in}{1.471131in}}%
\pgfpathlineto{\pgfqpoint{1.561493in}{1.470994in}}%
\pgfpathlineto{\pgfqpoint{1.647845in}{1.467780in}}%
\pgfpathlineto{\pgfqpoint{1.734961in}{1.461419in}}%
\pgfpathlineto{\pgfqpoint{1.828191in}{1.451460in}}%
\pgfpathlineto{\pgfqpoint{1.936704in}{1.436650in}}%
\pgfpathlineto{\pgfqpoint{2.086100in}{1.412869in}}%
\pgfpathlineto{\pgfqpoint{2.333693in}{1.370089in}}%
\pgfpathlineto{\pgfqpoint{2.522445in}{1.334202in}}%
\pgfpathlineto{\pgfqpoint{2.730300in}{1.295187in}}%
\pgfpathlineto{\pgfqpoint{2.837285in}{1.278395in}}%
\pgfpathlineto{\pgfqpoint{2.930514in}{1.266856in}}%
\pgfpathlineto{\pgfqpoint{3.020687in}{1.258820in}}%
\pgfpathlineto{\pgfqpoint{3.115445in}{1.253524in}}%
\pgfpathlineto{\pgfqpoint{3.210203in}{1.251296in}}%
\pgfpathlineto{\pgfqpoint{3.299229in}{1.252252in}}%
\pgfpathlineto{\pgfqpoint{3.381760in}{1.256186in}}%
\pgfpathlineto{\pgfqpoint{3.459706in}{1.262965in}}%
\pgfpathlineto{\pgfqpoint{3.535359in}{1.272634in}}%
\pgfpathlineto{\pgfqpoint{3.611013in}{1.285418in}}%
\pgfpathlineto{\pgfqpoint{3.690105in}{1.301954in}}%
\pgfpathlineto{\pgfqpoint{3.778367in}{1.323667in}}%
\pgfpathlineto{\pgfqpoint{3.896432in}{1.356196in}}%
\pgfpathlineto{\pgfqpoint{4.072957in}{1.404708in}}%
\pgfpathlineto{\pgfqpoint{4.155870in}{1.424043in}}%
\pgfpathlineto{\pgfqpoint{4.227703in}{1.437660in}}%
\pgfpathlineto{\pgfqpoint{4.294950in}{1.447272in}}%
\pgfpathlineto{\pgfqpoint{4.361051in}{1.453586in}}%
\pgfpathlineto{\pgfqpoint{4.429063in}{1.456950in}}%
\pgfpathlineto{\pgfqpoint{4.504716in}{1.457504in}}%
\pgfpathlineto{\pgfqpoint{4.508155in}{1.457465in}}%
\pgfpathlineto{\pgfqpoint{4.508155in}{1.457465in}}%
\pgfusepath{stroke}%
\end{pgfscope}%
\begin{pgfscope}%
\pgfpathrectangle{\pgfqpoint{0.687659in}{0.580556in}}{\pgfqpoint{3.820497in}{2.129775in}}%
\pgfusepath{clip}%
\pgfsetbuttcap%
\pgfsetroundjoin%
\pgfsetlinewidth{1.003750pt}%
\definecolor{currentstroke}{rgb}{1.000000,0.498039,0.054902}%
\pgfsetstrokecolor{currentstroke}%
\pgfsetdash{{3.700000pt}{1.600000pt}}{0.000000pt}%
\pgfpathmoveto{\pgfqpoint{0.687659in}{1.582879in}}%
\pgfpathlineto{\pgfqpoint{0.750703in}{1.508979in}}%
\pgfpathlineto{\pgfqpoint{0.835144in}{1.410815in}}%
\pgfpathlineto{\pgfqpoint{0.880613in}{1.362761in}}%
\pgfpathlineto{\pgfqpoint{0.919204in}{1.326329in}}%
\pgfpathlineto{\pgfqpoint{0.953974in}{1.297667in}}%
\pgfpathlineto{\pgfqpoint{0.986069in}{1.275146in}}%
\pgfpathlineto{\pgfqpoint{1.016254in}{1.257678in}}%
\pgfpathlineto{\pgfqpoint{1.044911in}{1.244579in}}%
\pgfpathlineto{\pgfqpoint{1.072803in}{1.235174in}}%
\pgfpathlineto{\pgfqpoint{1.100313in}{1.229161in}}%
\pgfpathlineto{\pgfqpoint{1.127442in}{1.226391in}}%
\pgfpathlineto{\pgfqpoint{1.154570in}{1.226708in}}%
\pgfpathlineto{\pgfqpoint{1.182462in}{1.230158in}}%
\pgfpathlineto{\pgfqpoint{1.211119in}{1.236863in}}%
\pgfpathlineto{\pgfqpoint{1.240922in}{1.247035in}}%
\pgfpathlineto{\pgfqpoint{1.272635in}{1.261164in}}%
\pgfpathlineto{\pgfqpoint{1.307023in}{1.279936in}}%
\pgfpathlineto{\pgfqpoint{1.345232in}{1.304410in}}%
\pgfpathlineto{\pgfqpoint{1.389936in}{1.336882in}}%
\pgfpathlineto{\pgfqpoint{1.448395in}{1.383469in}}%
\pgfpathlineto{\pgfqpoint{1.577541in}{1.487560in}}%
\pgfpathlineto{\pgfqpoint{1.625302in}{1.521328in}}%
\pgfpathlineto{\pgfqpoint{1.667332in}{1.547331in}}%
\pgfpathlineto{\pgfqpoint{1.705923in}{1.567644in}}%
\pgfpathlineto{\pgfqpoint{1.742603in}{1.583523in}}%
\pgfpathlineto{\pgfqpoint{1.778137in}{1.595595in}}%
\pgfpathlineto{\pgfqpoint{1.813289in}{1.604301in}}%
\pgfpathlineto{\pgfqpoint{1.848441in}{1.609834in}}%
\pgfpathlineto{\pgfqpoint{1.883976in}{1.612319in}}%
\pgfpathlineto{\pgfqpoint{1.920656in}{1.611791in}}%
\pgfpathlineto{\pgfqpoint{1.959247in}{1.608109in}}%
\pgfpathlineto{\pgfqpoint{2.000130in}{1.601075in}}%
\pgfpathlineto{\pgfqpoint{2.044835in}{1.590194in}}%
\pgfpathlineto{\pgfqpoint{2.094888in}{1.574756in}}%
\pgfpathlineto{\pgfqpoint{2.153730in}{1.553255in}}%
\pgfpathlineto{\pgfqpoint{2.228619in}{1.522399in}}%
\pgfpathlineto{\pgfqpoint{2.342863in}{1.471594in}}%
\pgfpathlineto{\pgfqpoint{2.519006in}{1.389868in}}%
\pgfpathlineto{\pgfqpoint{2.622169in}{1.338839in}}%
\pgfpathlineto{\pgfqpoint{2.704318in}{1.294789in}}%
\pgfpathlineto{\pgfqpoint{2.781118in}{1.250095in}}%
\pgfpathlineto{\pgfqpoint{2.862885in}{1.198832in}}%
\pgfpathlineto{\pgfqpoint{3.001583in}{1.107307in}}%
\pgfpathlineto{\pgfqpoint{3.080293in}{1.057587in}}%
\pgfpathlineto{\pgfqpoint{3.134167in}{1.027089in}}%
\pgfpathlineto{\pgfqpoint{3.178871in}{1.005222in}}%
\pgfpathlineto{\pgfqpoint{3.218227in}{0.989329in}}%
\pgfpathlineto{\pgfqpoint{3.253761in}{0.978217in}}%
\pgfpathlineto{\pgfqpoint{3.287002in}{0.971010in}}%
\pgfpathlineto{\pgfqpoint{3.318334in}{0.967348in}}%
\pgfpathlineto{\pgfqpoint{3.348136in}{0.966923in}}%
\pgfpathlineto{\pgfqpoint{3.377175in}{0.969567in}}%
\pgfpathlineto{\pgfqpoint{3.405832in}{0.975284in}}%
\pgfpathlineto{\pgfqpoint{3.434106in}{0.984061in}}%
\pgfpathlineto{\pgfqpoint{3.462763in}{0.996215in}}%
\pgfpathlineto{\pgfqpoint{3.491801in}{1.011917in}}%
\pgfpathlineto{\pgfqpoint{3.521604in}{1.031575in}}%
\pgfpathlineto{\pgfqpoint{3.552553in}{1.055730in}}%
\pgfpathlineto{\pgfqpoint{3.585031in}{1.085042in}}%
\pgfpathlineto{\pgfqpoint{3.619801in}{1.120677in}}%
\pgfpathlineto{\pgfqpoint{3.657628in}{1.164009in}}%
\pgfpathlineto{\pgfqpoint{3.699657in}{1.216985in}}%
\pgfpathlineto{\pgfqpoint{3.750093in}{1.285785in}}%
\pgfpathlineto{\pgfqpoint{3.831095in}{1.402450in}}%
\pgfpathlineto{\pgfqpoint{3.909041in}{1.512729in}}%
\pgfpathlineto{\pgfqpoint{3.957184in}{1.575576in}}%
\pgfpathlineto{\pgfqpoint{3.997304in}{1.623036in}}%
\pgfpathlineto{\pgfqpoint{4.032838in}{1.660414in}}%
\pgfpathlineto{\pgfqpoint{4.064933in}{1.689838in}}%
\pgfpathlineto{\pgfqpoint{4.094736in}{1.713092in}}%
\pgfpathlineto{\pgfqpoint{4.122629in}{1.731045in}}%
\pgfpathlineto{\pgfqpoint{4.148993in}{1.744453in}}%
\pgfpathlineto{\pgfqpoint{4.174592in}{1.754039in}}%
\pgfpathlineto{\pgfqpoint{4.199046in}{1.759963in}}%
\pgfpathlineto{\pgfqpoint{4.223118in}{1.762670in}}%
\pgfpathlineto{\pgfqpoint{4.247189in}{1.762267in}}%
\pgfpathlineto{\pgfqpoint{4.271261in}{1.758773in}}%
\pgfpathlineto{\pgfqpoint{4.295714in}{1.752109in}}%
\pgfpathlineto{\pgfqpoint{4.320932in}{1.742034in}}%
\pgfpathlineto{\pgfqpoint{4.347296in}{1.728161in}}%
\pgfpathlineto{\pgfqpoint{4.375189in}{1.709967in}}%
\pgfpathlineto{\pgfqpoint{4.405374in}{1.686513in}}%
\pgfpathlineto{\pgfqpoint{4.438615in}{1.656637in}}%
\pgfpathlineto{\pgfqpoint{4.476060in}{1.618688in}}%
\pgfpathlineto{\pgfqpoint{4.508155in}{1.583313in}}%
\pgfpathlineto{\pgfqpoint{4.508155in}{1.583313in}}%
\pgfusepath{stroke}%
\end{pgfscope}%
\begin{pgfscope}%
\pgfpathrectangle{\pgfqpoint{0.687659in}{0.580556in}}{\pgfqpoint{3.820497in}{2.129775in}}%
\pgfusepath{clip}%
\pgfsetbuttcap%
\pgfsetroundjoin%
\pgfsetlinewidth{1.003750pt}%
\definecolor{currentstroke}{rgb}{0.172549,0.627451,0.172549}%
\pgfsetstrokecolor{currentstroke}%
\pgfsetdash{{6.400000pt}{1.600000pt}{1.000000pt}{1.600000pt}}{0.000000pt}%
\pgfpathmoveto{\pgfqpoint{0.687659in}{1.683225in}}%
\pgfpathlineto{\pgfqpoint{0.757199in}{1.463983in}}%
\pgfpathlineto{\pgfqpoint{0.808780in}{1.306728in}}%
\pgfpathlineto{\pgfqpoint{0.844697in}{1.206716in}}%
\pgfpathlineto{\pgfqpoint{0.875264in}{1.130410in}}%
\pgfpathlineto{\pgfqpoint{0.902010in}{1.071569in}}%
\pgfpathlineto{\pgfqpoint{0.926463in}{1.024974in}}%
\pgfpathlineto{\pgfqpoint{0.948625in}{0.989102in}}%
\pgfpathlineto{\pgfqpoint{0.969257in}{0.961368in}}%
\pgfpathlineto{\pgfqpoint{0.988362in}{0.940674in}}%
\pgfpathlineto{\pgfqpoint{1.006320in}{0.925641in}}%
\pgfpathlineto{\pgfqpoint{1.023132in}{0.915449in}}%
\pgfpathlineto{\pgfqpoint{1.039179in}{0.909192in}}%
\pgfpathlineto{\pgfqpoint{1.054845in}{0.906304in}}%
\pgfpathlineto{\pgfqpoint{1.070511in}{0.906526in}}%
\pgfpathlineto{\pgfqpoint{1.086558in}{0.909858in}}%
\pgfpathlineto{\pgfqpoint{1.103370in}{0.916571in}}%
\pgfpathlineto{\pgfqpoint{1.121328in}{0.927167in}}%
\pgfpathlineto{\pgfqpoint{1.140815in}{0.942391in}}%
\pgfpathlineto{\pgfqpoint{1.162212in}{0.963204in}}%
\pgfpathlineto{\pgfqpoint{1.185901in}{0.990738in}}%
\pgfpathlineto{\pgfqpoint{1.212647in}{1.026771in}}%
\pgfpathlineto{\pgfqpoint{1.243596in}{1.073919in}}%
\pgfpathlineto{\pgfqpoint{1.280659in}{1.136284in}}%
\pgfpathlineto{\pgfqpoint{1.331477in}{1.228376in}}%
\pgfpathlineto{\pgfqpoint{1.449924in}{1.444955in}}%
\pgfpathlineto{\pgfqpoint{1.492336in}{1.514903in}}%
\pgfpathlineto{\pgfqpoint{1.528634in}{1.569029in}}%
\pgfpathlineto{\pgfqpoint{1.561493in}{1.612625in}}%
\pgfpathlineto{\pgfqpoint{1.591678in}{1.647675in}}%
\pgfpathlineto{\pgfqpoint{1.619953in}{1.675897in}}%
\pgfpathlineto{\pgfqpoint{1.646699in}{1.698346in}}%
\pgfpathlineto{\pgfqpoint{1.672299in}{1.715909in}}%
\pgfpathlineto{\pgfqpoint{1.697135in}{1.729279in}}%
\pgfpathlineto{\pgfqpoint{1.721588in}{1.738960in}}%
\pgfpathlineto{\pgfqpoint{1.745660in}{1.745194in}}%
\pgfpathlineto{\pgfqpoint{1.769731in}{1.748273in}}%
\pgfpathlineto{\pgfqpoint{1.794185in}{1.748320in}}%
\pgfpathlineto{\pgfqpoint{1.819785in}{1.745238in}}%
\pgfpathlineto{\pgfqpoint{1.846531in}{1.738853in}}%
\pgfpathlineto{\pgfqpoint{1.875188in}{1.728760in}}%
\pgfpathlineto{\pgfqpoint{1.906519in}{1.714341in}}%
\pgfpathlineto{\pgfqpoint{1.941671in}{1.694626in}}%
\pgfpathlineto{\pgfqpoint{1.982936in}{1.667773in}}%
\pgfpathlineto{\pgfqpoint{2.037193in}{1.628488in}}%
\pgfpathlineto{\pgfqpoint{2.189264in}{1.516347in}}%
\pgfpathlineto{\pgfqpoint{2.238553in}{1.485056in}}%
\pgfpathlineto{\pgfqpoint{2.283257in}{1.460299in}}%
\pgfpathlineto{\pgfqpoint{2.326051in}{1.440102in}}%
\pgfpathlineto{\pgfqpoint{2.368845in}{1.423309in}}%
\pgfpathlineto{\pgfqpoint{2.413549in}{1.409111in}}%
\pgfpathlineto{\pgfqpoint{2.463221in}{1.396657in}}%
\pgfpathlineto{\pgfqpoint{2.527794in}{1.383908in}}%
\pgfpathlineto{\pgfqpoint{2.645477in}{1.361091in}}%
\pgfpathlineto{\pgfqpoint{2.693238in}{1.348219in}}%
\pgfpathlineto{\pgfqpoint{2.734503in}{1.333944in}}%
\pgfpathlineto{\pgfqpoint{2.773094in}{1.317320in}}%
\pgfpathlineto{\pgfqpoint{2.810157in}{1.297971in}}%
\pgfpathlineto{\pgfqpoint{2.846837in}{1.275307in}}%
\pgfpathlineto{\pgfqpoint{2.883900in}{1.248744in}}%
\pgfpathlineto{\pgfqpoint{2.922491in}{1.217214in}}%
\pgfpathlineto{\pgfqpoint{2.963374in}{1.179750in}}%
\pgfpathlineto{\pgfqpoint{3.008842in}{1.133784in}}%
\pgfpathlineto{\pgfqpoint{3.064245in}{1.073132in}}%
\pgfpathlineto{\pgfqpoint{3.229307in}{0.889266in}}%
\pgfpathlineto{\pgfqpoint{3.268280in}{0.852364in}}%
\pgfpathlineto{\pgfqpoint{3.301140in}{0.825340in}}%
\pgfpathlineto{\pgfqpoint{3.330560in}{0.805040in}}%
\pgfpathlineto{\pgfqpoint{3.357307in}{0.790279in}}%
\pgfpathlineto{\pgfqpoint{3.381760in}{0.780221in}}%
\pgfpathlineto{\pgfqpoint{3.404685in}{0.774041in}}%
\pgfpathlineto{\pgfqpoint{3.426464in}{0.771287in}}%
\pgfpathlineto{\pgfqpoint{3.447479in}{0.771673in}}%
\pgfpathlineto{\pgfqpoint{3.468112in}{0.775092in}}%
\pgfpathlineto{\pgfqpoint{3.488745in}{0.781634in}}%
\pgfpathlineto{\pgfqpoint{3.509760in}{0.791603in}}%
\pgfpathlineto{\pgfqpoint{3.531156in}{0.805257in}}%
\pgfpathlineto{\pgfqpoint{3.553318in}{0.823182in}}%
\pgfpathlineto{\pgfqpoint{3.576243in}{0.845806in}}%
\pgfpathlineto{\pgfqpoint{3.600314in}{0.874024in}}%
\pgfpathlineto{\pgfqpoint{3.625532in}{0.908433in}}%
\pgfpathlineto{\pgfqpoint{3.652278in}{0.950221in}}%
\pgfpathlineto{\pgfqpoint{3.680935in}{1.000806in}}%
\pgfpathlineto{\pgfqpoint{3.712266in}{1.062581in}}%
\pgfpathlineto{\pgfqpoint{3.746654in}{1.137463in}}%
\pgfpathlineto{\pgfqpoint{3.785245in}{1.229128in}}%
\pgfpathlineto{\pgfqpoint{3.831478in}{1.347267in}}%
\pgfpathlineto{\pgfqpoint{3.900253in}{1.532636in}}%
\pgfpathlineto{\pgfqpoint{3.983166in}{1.754246in}}%
\pgfpathlineto{\pgfqpoint{4.026342in}{1.860765in}}%
\pgfpathlineto{\pgfqpoint{4.061112in}{1.938654in}}%
\pgfpathlineto{\pgfqpoint{4.090915in}{1.998164in}}%
\pgfpathlineto{\pgfqpoint{4.117279in}{2.044181in}}%
\pgfpathlineto{\pgfqpoint{4.140969in}{2.079497in}}%
\pgfpathlineto{\pgfqpoint{4.161984in}{2.105518in}}%
\pgfpathlineto{\pgfqpoint{4.181088in}{2.124469in}}%
\pgfpathlineto{\pgfqpoint{4.198664in}{2.137703in}}%
\pgfpathlineto{\pgfqpoint{4.215094in}{2.146269in}}%
\pgfpathlineto{\pgfqpoint{4.230377in}{2.150833in}}%
\pgfpathlineto{\pgfqpoint{4.245279in}{2.152051in}}%
\pgfpathlineto{\pgfqpoint{4.259798in}{2.150123in}}%
\pgfpathlineto{\pgfqpoint{4.274317in}{2.145094in}}%
\pgfpathlineto{\pgfqpoint{4.289601in}{2.136443in}}%
\pgfpathlineto{\pgfqpoint{4.305649in}{2.123664in}}%
\pgfpathlineto{\pgfqpoint{4.322843in}{2.105816in}}%
\pgfpathlineto{\pgfqpoint{4.341183in}{2.082129in}}%
\pgfpathlineto{\pgfqpoint{4.361051in}{2.051217in}}%
\pgfpathlineto{\pgfqpoint{4.382830in}{2.011336in}}%
\pgfpathlineto{\pgfqpoint{4.406520in}{1.961286in}}%
\pgfpathlineto{\pgfqpoint{4.432884in}{1.898158in}}%
\pgfpathlineto{\pgfqpoint{4.463069in}{1.817552in}}%
\pgfpathlineto{\pgfqpoint{4.499367in}{1.711293in}}%
\pgfpathlineto{\pgfqpoint{4.508155in}{1.684401in}}%
\pgfpathlineto{\pgfqpoint{4.508155in}{1.684401in}}%
\pgfusepath{stroke}%
\end{pgfscope}%
\begin{pgfscope}%
\pgfpathrectangle{\pgfqpoint{0.687659in}{0.580556in}}{\pgfqpoint{3.820497in}{2.129775in}}%
\pgfusepath{clip}%
\pgfsetbuttcap%
\pgfsetroundjoin%
\pgfsetlinewidth{1.003750pt}%
\definecolor{currentstroke}{rgb}{0.839216,0.152941,0.156863}%
\pgfsetstrokecolor{currentstroke}%
\pgfsetdash{{3.700000pt}{1.600000pt}}{0.000000pt}%
\pgfpathmoveto{\pgfqpoint{0.687659in}{1.738522in}}%
\pgfpathlineto{\pgfqpoint{0.739240in}{1.509629in}}%
\pgfpathlineto{\pgfqpoint{0.789676in}{1.291083in}}%
\pgfpathlineto{\pgfqpoint{0.822153in}{1.163787in}}%
\pgfpathlineto{\pgfqpoint{0.849282in}{1.069421in}}%
\pgfpathlineto{\pgfqpoint{0.873353in}{0.996383in}}%
\pgfpathlineto{\pgfqpoint{0.895514in}{0.938680in}}%
\pgfpathlineto{\pgfqpoint{0.916147in}{0.893300in}}%
\pgfpathlineto{\pgfqpoint{0.935634in}{0.857714in}}%
\pgfpathlineto{\pgfqpoint{0.953974in}{0.830529in}}%
\pgfpathlineto{\pgfqpoint{0.971168in}{0.810454in}}%
\pgfpathlineto{\pgfqpoint{0.987215in}{0.796308in}}%
\pgfpathlineto{\pgfqpoint{1.002499in}{0.786829in}}%
\pgfpathlineto{\pgfqpoint{1.017018in}{0.781316in}}%
\pgfpathlineto{\pgfqpoint{1.031538in}{0.779087in}}%
\pgfpathlineto{\pgfqpoint{1.046057in}{0.780017in}}%
\pgfpathlineto{\pgfqpoint{1.060958in}{0.784122in}}%
\pgfpathlineto{\pgfqpoint{1.076624in}{0.791718in}}%
\pgfpathlineto{\pgfqpoint{1.093436in}{0.803412in}}%
\pgfpathlineto{\pgfqpoint{1.111776in}{0.820090in}}%
\pgfpathlineto{\pgfqpoint{1.132027in}{0.842903in}}%
\pgfpathlineto{\pgfqpoint{1.154188in}{0.872674in}}%
\pgfpathlineto{\pgfqpoint{1.179406in}{0.911964in}}%
\pgfpathlineto{\pgfqpoint{1.208062in}{0.962581in}}%
\pgfpathlineto{\pgfqpoint{1.242068in}{1.029242in}}%
\pgfpathlineto{\pgfqpoint{1.285244in}{1.121138in}}%
\pgfpathlineto{\pgfqpoint{1.378855in}{1.330297in}}%
\pgfpathlineto{\pgfqpoint{1.433112in}{1.446403in}}%
\pgfpathlineto{\pgfqpoint{1.473995in}{1.526815in}}%
\pgfpathlineto{\pgfqpoint{1.509147in}{1.589366in}}%
\pgfpathlineto{\pgfqpoint{1.540861in}{1.639677in}}%
\pgfpathlineto{\pgfqpoint{1.569899in}{1.680145in}}%
\pgfpathlineto{\pgfqpoint{1.597028in}{1.712822in}}%
\pgfpathlineto{\pgfqpoint{1.622628in}{1.738955in}}%
\pgfpathlineto{\pgfqpoint{1.647081in}{1.759585in}}%
\pgfpathlineto{\pgfqpoint{1.670389in}{1.775295in}}%
\pgfpathlineto{\pgfqpoint{1.692932in}{1.786854in}}%
\pgfpathlineto{\pgfqpoint{1.715093in}{1.794795in}}%
\pgfpathlineto{\pgfqpoint{1.737254in}{1.799441in}}%
\pgfpathlineto{\pgfqpoint{1.759415in}{1.800912in}}%
\pgfpathlineto{\pgfqpoint{1.781958in}{1.799306in}}%
\pgfpathlineto{\pgfqpoint{1.805266in}{1.794547in}}%
\pgfpathlineto{\pgfqpoint{1.830101in}{1.786260in}}%
\pgfpathlineto{\pgfqpoint{1.856465in}{1.774161in}}%
\pgfpathlineto{\pgfqpoint{1.885504in}{1.757353in}}%
\pgfpathlineto{\pgfqpoint{1.917981in}{1.734883in}}%
\pgfpathlineto{\pgfqpoint{1.956190in}{1.704527in}}%
\pgfpathlineto{\pgfqpoint{2.005480in}{1.661127in}}%
\pgfpathlineto{\pgfqpoint{2.163664in}{1.518857in}}%
\pgfpathlineto{\pgfqpoint{2.207222in}{1.485506in}}%
\pgfpathlineto{\pgfqpoint{2.246577in}{1.459186in}}%
\pgfpathlineto{\pgfqpoint{2.283640in}{1.438069in}}%
\pgfpathlineto{\pgfqpoint{2.319556in}{1.421126in}}%
\pgfpathlineto{\pgfqpoint{2.355472in}{1.407607in}}%
\pgfpathlineto{\pgfqpoint{2.392153in}{1.397127in}}%
\pgfpathlineto{\pgfqpoint{2.430743in}{1.389334in}}%
\pgfpathlineto{\pgfqpoint{2.473919in}{1.383829in}}%
\pgfpathlineto{\pgfqpoint{2.528940in}{1.380134in}}%
\pgfpathlineto{\pgfqpoint{2.651208in}{1.372928in}}%
\pgfpathlineto{\pgfqpoint{2.693238in}{1.366651in}}%
\pgfpathlineto{\pgfqpoint{2.729918in}{1.358093in}}%
\pgfpathlineto{\pgfqpoint{2.763924in}{1.347004in}}%
\pgfpathlineto{\pgfqpoint{2.796402in}{1.333177in}}%
\pgfpathlineto{\pgfqpoint{2.828497in}{1.316123in}}%
\pgfpathlineto{\pgfqpoint{2.860592in}{1.295514in}}%
\pgfpathlineto{\pgfqpoint{2.893070in}{1.270954in}}%
\pgfpathlineto{\pgfqpoint{2.926693in}{1.241623in}}%
\pgfpathlineto{\pgfqpoint{2.962228in}{1.206479in}}%
\pgfpathlineto{\pgfqpoint{3.000819in}{1.163884in}}%
\pgfpathlineto{\pgfqpoint{3.044759in}{1.110603in}}%
\pgfpathlineto{\pgfqpoint{3.099779in}{1.038698in}}%
\pgfpathlineto{\pgfqpoint{3.240006in}{0.853067in}}%
\pgfpathlineto{\pgfqpoint{3.279361in}{0.807524in}}%
\pgfpathlineto{\pgfqpoint{3.312220in}{0.774059in}}%
\pgfpathlineto{\pgfqpoint{3.340877in}{0.749096in}}%
\pgfpathlineto{\pgfqpoint{3.366859in}{0.730426in}}%
\pgfpathlineto{\pgfqpoint{3.390548in}{0.717085in}}%
\pgfpathlineto{\pgfqpoint{3.412709in}{0.708079in}}%
\pgfpathlineto{\pgfqpoint{3.433342in}{0.702933in}}%
\pgfpathlineto{\pgfqpoint{3.453211in}{0.701103in}}%
\pgfpathlineto{\pgfqpoint{3.472697in}{0.702427in}}%
\pgfpathlineto{\pgfqpoint{3.491801in}{0.706842in}}%
\pgfpathlineto{\pgfqpoint{3.510906in}{0.714439in}}%
\pgfpathlineto{\pgfqpoint{3.530392in}{0.725549in}}%
\pgfpathlineto{\pgfqpoint{3.550643in}{0.740760in}}%
\pgfpathlineto{\pgfqpoint{3.571658in}{0.760547in}}%
\pgfpathlineto{\pgfqpoint{3.593437in}{0.785378in}}%
\pgfpathlineto{\pgfqpoint{3.616362in}{0.816265in}}%
\pgfpathlineto{\pgfqpoint{3.640816in}{0.854517in}}%
\pgfpathlineto{\pgfqpoint{3.666798in}{0.901028in}}%
\pgfpathlineto{\pgfqpoint{3.694690in}{0.957433in}}%
\pgfpathlineto{\pgfqpoint{3.724875in}{1.025597in}}%
\pgfpathlineto{\pgfqpoint{3.758117in}{1.108505in}}%
\pgfpathlineto{\pgfqpoint{3.795561in}{1.210494in}}%
\pgfpathlineto{\pgfqpoint{3.839883in}{1.340636in}}%
\pgfpathlineto{\pgfqpoint{3.901018in}{1.530717in}}%
\pgfpathlineto{\pgfqpoint{4.001507in}{1.843289in}}%
\pgfpathlineto{\pgfqpoint{4.043918in}{1.964385in}}%
\pgfpathlineto{\pgfqpoint{4.078306in}{2.053543in}}%
\pgfpathlineto{\pgfqpoint{4.107727in}{2.121574in}}%
\pgfpathlineto{\pgfqpoint{4.133709in}{2.174173in}}%
\pgfpathlineto{\pgfqpoint{4.157016in}{2.214611in}}%
\pgfpathlineto{\pgfqpoint{4.178031in}{2.245065in}}%
\pgfpathlineto{\pgfqpoint{4.196754in}{2.267036in}}%
\pgfpathlineto{\pgfqpoint{4.213948in}{2.282674in}}%
\pgfpathlineto{\pgfqpoint{4.229613in}{2.292956in}}%
\pgfpathlineto{\pgfqpoint{4.244132in}{2.298977in}}%
\pgfpathlineto{\pgfqpoint{4.257888in}{2.301469in}}%
\pgfpathlineto{\pgfqpoint{4.271261in}{2.300812in}}%
\pgfpathlineto{\pgfqpoint{4.284634in}{2.297047in}}%
\pgfpathlineto{\pgfqpoint{4.298389in}{2.289866in}}%
\pgfpathlineto{\pgfqpoint{4.312526in}{2.278923in}}%
\pgfpathlineto{\pgfqpoint{4.327428in}{2.263368in}}%
\pgfpathlineto{\pgfqpoint{4.343093in}{2.242448in}}%
\pgfpathlineto{\pgfqpoint{4.359905in}{2.214629in}}%
\pgfpathlineto{\pgfqpoint{4.377863in}{2.178651in}}%
\pgfpathlineto{\pgfqpoint{4.396968in}{2.133259in}}%
\pgfpathlineto{\pgfqpoint{4.417982in}{2.075017in}}%
\pgfpathlineto{\pgfqpoint{4.440908in}{2.002020in}}%
\pgfpathlineto{\pgfqpoint{4.466508in}{1.909851in}}%
\pgfpathlineto{\pgfqpoint{4.496310in}{1.790511in}}%
\pgfpathlineto{\pgfqpoint{4.508155in}{1.740166in}}%
\pgfpathlineto{\pgfqpoint{4.508155in}{1.740166in}}%
\pgfusepath{stroke}%
\end{pgfscope}%
\begin{pgfscope}%
\pgfpathrectangle{\pgfqpoint{0.687659in}{0.580556in}}{\pgfqpoint{3.820497in}{2.129775in}}%
\pgfusepath{clip}%
\pgfsetrectcap%
\pgfsetroundjoin%
\pgfsetlinewidth{1.003750pt}%
\definecolor{currentstroke}{rgb}{0.580392,0.403922,0.741176}%
\pgfsetstrokecolor{currentstroke}%
\pgfsetdash{}{0pt}%
\pgfpathmoveto{\pgfqpoint{0.687659in}{1.780227in}}%
\pgfpathlineto{\pgfqpoint{0.710202in}{1.656406in}}%
\pgfpathlineto{\pgfqpoint{0.756052in}{1.404130in}}%
\pgfpathlineto{\pgfqpoint{0.784709in}{1.265152in}}%
\pgfpathlineto{\pgfqpoint{0.811455in}{1.150081in}}%
\pgfpathlineto{\pgfqpoint{0.836673in}{1.054739in}}%
\pgfpathlineto{\pgfqpoint{0.860362in}{0.976862in}}%
\pgfpathlineto{\pgfqpoint{0.882523in}{0.914217in}}%
\pgfpathlineto{\pgfqpoint{0.903538in}{0.863831in}}%
\pgfpathlineto{\pgfqpoint{0.923025in}{0.824825in}}%
\pgfpathlineto{\pgfqpoint{0.941365in}{0.794758in}}%
\pgfpathlineto{\pgfqpoint{0.958559in}{0.772280in}}%
\pgfpathlineto{\pgfqpoint{0.974607in}{0.756152in}}%
\pgfpathlineto{\pgfqpoint{0.989890in}{0.745022in}}%
\pgfpathlineto{\pgfqpoint{1.004409in}{0.738153in}}%
\pgfpathlineto{\pgfqpoint{1.018547in}{0.734818in}}%
\pgfpathlineto{\pgfqpoint{1.032684in}{0.734675in}}%
\pgfpathlineto{\pgfqpoint{1.046821in}{0.737603in}}%
\pgfpathlineto{\pgfqpoint{1.061723in}{0.743871in}}%
\pgfpathlineto{\pgfqpoint{1.077770in}{0.754099in}}%
\pgfpathlineto{\pgfqpoint{1.094964in}{0.768835in}}%
\pgfpathlineto{\pgfqpoint{1.113687in}{0.789028in}}%
\pgfpathlineto{\pgfqpoint{1.134701in}{0.816432in}}%
\pgfpathlineto{\pgfqpoint{1.158009in}{0.852105in}}%
\pgfpathlineto{\pgfqpoint{1.184373in}{0.898309in}}%
\pgfpathlineto{\pgfqpoint{1.214940in}{0.958373in}}%
\pgfpathlineto{\pgfqpoint{1.252002in}{1.038392in}}%
\pgfpathlineto{\pgfqpoint{1.302820in}{1.156165in}}%
\pgfpathlineto{\pgfqpoint{1.418593in}{1.426686in}}%
\pgfpathlineto{\pgfqpoint{1.460622in}{1.515698in}}%
\pgfpathlineto{\pgfqpoint{1.496539in}{1.584731in}}%
\pgfpathlineto{\pgfqpoint{1.528634in}{1.639897in}}%
\pgfpathlineto{\pgfqpoint{1.558055in}{1.684463in}}%
\pgfpathlineto{\pgfqpoint{1.585183in}{1.720131in}}%
\pgfpathlineto{\pgfqpoint{1.610783in}{1.748829in}}%
\pgfpathlineto{\pgfqpoint{1.634854in}{1.771331in}}%
\pgfpathlineto{\pgfqpoint{1.658162in}{1.788954in}}%
\pgfpathlineto{\pgfqpoint{1.680323in}{1.801934in}}%
\pgfpathlineto{\pgfqpoint{1.702102in}{1.811157in}}%
\pgfpathlineto{\pgfqpoint{1.723499in}{1.816892in}}%
\pgfpathlineto{\pgfqpoint{1.744896in}{1.819439in}}%
\pgfpathlineto{\pgfqpoint{1.766675in}{1.818894in}}%
\pgfpathlineto{\pgfqpoint{1.789218in}{1.815170in}}%
\pgfpathlineto{\pgfqpoint{1.812525in}{1.808158in}}%
\pgfpathlineto{\pgfqpoint{1.837361in}{1.797431in}}%
\pgfpathlineto{\pgfqpoint{1.864489in}{1.782259in}}%
\pgfpathlineto{\pgfqpoint{1.894292in}{1.761969in}}%
\pgfpathlineto{\pgfqpoint{1.928298in}{1.734984in}}%
\pgfpathlineto{\pgfqpoint{1.969563in}{1.698109in}}%
\pgfpathlineto{\pgfqpoint{2.029169in}{1.640220in}}%
\pgfpathlineto{\pgfqpoint{2.127748in}{1.544512in}}%
\pgfpathlineto{\pgfqpoint{2.174362in}{1.503834in}}%
\pgfpathlineto{\pgfqpoint{2.214482in}{1.472793in}}%
\pgfpathlineto{\pgfqpoint{2.251544in}{1.447961in}}%
\pgfpathlineto{\pgfqpoint{2.286696in}{1.428098in}}%
\pgfpathlineto{\pgfqpoint{2.320702in}{1.412387in}}%
\pgfpathlineto{\pgfqpoint{2.354708in}{1.400062in}}%
\pgfpathlineto{\pgfqpoint{2.389478in}{1.390765in}}%
\pgfpathlineto{\pgfqpoint{2.426158in}{1.384210in}}%
\pgfpathlineto{\pgfqpoint{2.466278in}{1.380237in}}%
\pgfpathlineto{\pgfqpoint{2.515185in}{1.378666in}}%
\pgfpathlineto{\pgfqpoint{2.669166in}{1.376145in}}%
\pgfpathlineto{\pgfqpoint{2.707375in}{1.370905in}}%
\pgfpathlineto{\pgfqpoint{2.741381in}{1.363177in}}%
\pgfpathlineto{\pgfqpoint{2.773476in}{1.352740in}}%
\pgfpathlineto{\pgfqpoint{2.804425in}{1.339434in}}%
\pgfpathlineto{\pgfqpoint{2.834992in}{1.322917in}}%
\pgfpathlineto{\pgfqpoint{2.865559in}{1.302886in}}%
\pgfpathlineto{\pgfqpoint{2.896891in}{1.278630in}}%
\pgfpathlineto{\pgfqpoint{2.929368in}{1.249531in}}%
\pgfpathlineto{\pgfqpoint{2.963756in}{1.214493in}}%
\pgfpathlineto{\pgfqpoint{3.000819in}{1.172229in}}%
\pgfpathlineto{\pgfqpoint{3.042466in}{1.119911in}}%
\pgfpathlineto{\pgfqpoint{3.093284in}{1.050817in}}%
\pgfpathlineto{\pgfqpoint{3.262167in}{0.816795in}}%
\pgfpathlineto{\pgfqpoint{3.298083in}{0.774802in}}%
\pgfpathlineto{\pgfqpoint{3.329032in}{0.743224in}}%
\pgfpathlineto{\pgfqpoint{3.356160in}{0.719780in}}%
\pgfpathlineto{\pgfqpoint{3.380996in}{0.702290in}}%
\pgfpathlineto{\pgfqpoint{3.403539in}{0.690056in}}%
\pgfpathlineto{\pgfqpoint{3.424554in}{0.682034in}}%
\pgfpathlineto{\pgfqpoint{3.444423in}{0.677648in}}%
\pgfpathlineto{\pgfqpoint{3.463527in}{0.676520in}}%
\pgfpathlineto{\pgfqpoint{3.482249in}{0.678481in}}%
\pgfpathlineto{\pgfqpoint{3.500972in}{0.683585in}}%
\pgfpathlineto{\pgfqpoint{3.519694in}{0.691925in}}%
\pgfpathlineto{\pgfqpoint{3.538798in}{0.703850in}}%
\pgfpathlineto{\pgfqpoint{3.558667in}{0.719973in}}%
\pgfpathlineto{\pgfqpoint{3.579300in}{0.740778in}}%
\pgfpathlineto{\pgfqpoint{3.601079in}{0.767252in}}%
\pgfpathlineto{\pgfqpoint{3.624004in}{0.800109in}}%
\pgfpathlineto{\pgfqpoint{3.648457in}{0.840730in}}%
\pgfpathlineto{\pgfqpoint{3.674439in}{0.890039in}}%
\pgfpathlineto{\pgfqpoint{3.702332in}{0.949750in}}%
\pgfpathlineto{\pgfqpoint{3.732517in}{1.021807in}}%
\pgfpathlineto{\pgfqpoint{3.765758in}{1.109328in}}%
\pgfpathlineto{\pgfqpoint{3.803585in}{1.217974in}}%
\pgfpathlineto{\pgfqpoint{3.848289in}{1.356228in}}%
\pgfpathlineto{\pgfqpoint{3.911716in}{1.563550in}}%
\pgfpathlineto{\pgfqpoint{4.004563in}{1.866249in}}%
\pgfpathlineto{\pgfqpoint{4.046975in}{1.993673in}}%
\pgfpathlineto{\pgfqpoint{4.081363in}{2.087673in}}%
\pgfpathlineto{\pgfqpoint{4.110784in}{2.159554in}}%
\pgfpathlineto{\pgfqpoint{4.136766in}{2.215274in}}%
\pgfpathlineto{\pgfqpoint{4.160073in}{2.258251in}}%
\pgfpathlineto{\pgfqpoint{4.181088in}{2.290756in}}%
\pgfpathlineto{\pgfqpoint{4.200192in}{2.314774in}}%
\pgfpathlineto{\pgfqpoint{4.217386in}{2.331611in}}%
\pgfpathlineto{\pgfqpoint{4.233052in}{2.342819in}}%
\pgfpathlineto{\pgfqpoint{4.247571in}{2.349551in}}%
\pgfpathlineto{\pgfqpoint{4.261326in}{2.352578in}}%
\pgfpathlineto{\pgfqpoint{4.274317in}{2.352361in}}%
\pgfpathlineto{\pgfqpoint{4.287308in}{2.349089in}}%
\pgfpathlineto{\pgfqpoint{4.300681in}{2.342462in}}%
\pgfpathlineto{\pgfqpoint{4.314437in}{2.332135in}}%
\pgfpathlineto{\pgfqpoint{4.328956in}{2.317309in}}%
\pgfpathlineto{\pgfqpoint{4.344622in}{2.296730in}}%
\pgfpathlineto{\pgfqpoint{4.361433in}{2.269293in}}%
\pgfpathlineto{\pgfqpoint{4.379392in}{2.233814in}}%
\pgfpathlineto{\pgfqpoint{4.398496in}{2.189041in}}%
\pgfpathlineto{\pgfqpoint{4.419129in}{2.132551in}}%
\pgfpathlineto{\pgfqpoint{4.441290in}{2.062531in}}%
\pgfpathlineto{\pgfqpoint{4.464979in}{1.977132in}}%
\pgfpathlineto{\pgfqpoint{4.485230in}{1.893516in}}%
\pgfpathlineto{\pgfqpoint{4.504716in}{1.800268in}}%
\pgfpathlineto{\pgfqpoint{4.508155in}{1.782196in}}%
\pgfpathlineto{\pgfqpoint{4.508155in}{1.782196in}}%
\pgfusepath{stroke}%
\end{pgfscope}%
\begin{pgfscope}%
\pgfpathrectangle{\pgfqpoint{0.687659in}{0.580556in}}{\pgfqpoint{3.820497in}{2.129775in}}%
\pgfusepath{clip}%
\pgfsetbuttcap%
\pgfsetroundjoin%
\pgfsetlinewidth{2.007500pt}%
\definecolor{currentstroke}{rgb}{0.000000,0.000000,0.000000}%
\pgfsetstrokecolor{currentstroke}%
\pgfsetdash{{2.000000pt}{3.300000pt}}{0.000000pt}%
\pgfpathmoveto{\pgfqpoint{0.687659in}{1.796183in}}%
\pgfpathlineto{\pgfqpoint{0.689569in}{1.785249in}}%
\pgfpathlineto{\pgfqpoint{0.723193in}{1.581735in}}%
\pgfpathlineto{\pgfqpoint{0.754142in}{1.412616in}}%
\pgfpathlineto{\pgfqpoint{0.783180in}{1.270630in}}%
\pgfpathlineto{\pgfqpoint{0.810309in}{1.153049in}}%
\pgfpathlineto{\pgfqpoint{0.835527in}{1.057016in}}%
\pgfpathlineto{\pgfqpoint{0.859216in}{0.978508in}}%
\pgfpathlineto{\pgfqpoint{0.881377in}{0.915297in}}%
\pgfpathlineto{\pgfqpoint{0.902392in}{0.864395in}}%
\pgfpathlineto{\pgfqpoint{0.921878in}{0.824933in}}%
\pgfpathlineto{\pgfqpoint{0.940219in}{0.794456in}}%
\pgfpathlineto{\pgfqpoint{0.957413in}{0.771614in}}%
\pgfpathlineto{\pgfqpoint{0.973460in}{0.755165in}}%
\pgfpathlineto{\pgfqpoint{0.988744in}{0.743746in}}%
\pgfpathlineto{\pgfqpoint{1.003263in}{0.736617in}}%
\pgfpathlineto{\pgfqpoint{1.017400in}{0.733044in}}%
\pgfpathlineto{\pgfqpoint{1.031538in}{0.732678in}}%
\pgfpathlineto{\pgfqpoint{1.045675in}{0.735398in}}%
\pgfpathlineto{\pgfqpoint{1.060576in}{0.741465in}}%
\pgfpathlineto{\pgfqpoint{1.076242in}{0.751216in}}%
\pgfpathlineto{\pgfqpoint{1.093436in}{0.765677in}}%
\pgfpathlineto{\pgfqpoint{1.112158in}{0.785605in}}%
\pgfpathlineto{\pgfqpoint{1.132791in}{0.812216in}}%
\pgfpathlineto{\pgfqpoint{1.155716in}{0.846948in}}%
\pgfpathlineto{\pgfqpoint{1.181698in}{0.892079in}}%
\pgfpathlineto{\pgfqpoint{1.211883in}{0.950970in}}%
\pgfpathlineto{\pgfqpoint{1.248181in}{1.028933in}}%
\pgfpathlineto{\pgfqpoint{1.297089in}{1.141967in}}%
\pgfpathlineto{\pgfqpoint{1.423942in}{1.438499in}}%
\pgfpathlineto{\pgfqpoint{1.464825in}{1.524394in}}%
\pgfpathlineto{\pgfqpoint{1.499977in}{1.591350in}}%
\pgfpathlineto{\pgfqpoint{1.531691in}{1.645315in}}%
\pgfpathlineto{\pgfqpoint{1.560729in}{1.688817in}}%
\pgfpathlineto{\pgfqpoint{1.587858in}{1.724028in}}%
\pgfpathlineto{\pgfqpoint{1.613457in}{1.752268in}}%
\pgfpathlineto{\pgfqpoint{1.637529in}{1.774323in}}%
\pgfpathlineto{\pgfqpoint{1.660454in}{1.791255in}}%
\pgfpathlineto{\pgfqpoint{1.682615in}{1.803873in}}%
\pgfpathlineto{\pgfqpoint{1.704394in}{1.812740in}}%
\pgfpathlineto{\pgfqpoint{1.725791in}{1.818127in}}%
\pgfpathlineto{\pgfqpoint{1.747188in}{1.820331in}}%
\pgfpathlineto{\pgfqpoint{1.768967in}{1.819446in}}%
\pgfpathlineto{\pgfqpoint{1.791510in}{1.815382in}}%
\pgfpathlineto{\pgfqpoint{1.815200in}{1.807890in}}%
\pgfpathlineto{\pgfqpoint{1.840418in}{1.796583in}}%
\pgfpathlineto{\pgfqpoint{1.867546in}{1.780978in}}%
\pgfpathlineto{\pgfqpoint{1.897731in}{1.759979in}}%
\pgfpathlineto{\pgfqpoint{1.932501in}{1.731892in}}%
\pgfpathlineto{\pgfqpoint{1.974913in}{1.693451in}}%
\pgfpathlineto{\pgfqpoint{2.039103in}{1.630544in}}%
\pgfpathlineto{\pgfqpoint{2.125837in}{1.546242in}}%
\pgfpathlineto{\pgfqpoint{2.172834in}{1.504972in}}%
\pgfpathlineto{\pgfqpoint{2.213335in}{1.473431in}}%
\pgfpathlineto{\pgfqpoint{2.250398in}{1.448433in}}%
\pgfpathlineto{\pgfqpoint{2.285550in}{1.428422in}}%
\pgfpathlineto{\pgfqpoint{2.319556in}{1.412581in}}%
\pgfpathlineto{\pgfqpoint{2.353562in}{1.400146in}}%
\pgfpathlineto{\pgfqpoint{2.387950in}{1.390844in}}%
\pgfpathlineto{\pgfqpoint{2.424248in}{1.384245in}}%
\pgfpathlineto{\pgfqpoint{2.463985in}{1.380207in}}%
\pgfpathlineto{\pgfqpoint{2.511746in}{1.378593in}}%
\pgfpathlineto{\pgfqpoint{2.674515in}{1.375867in}}%
\pgfpathlineto{\pgfqpoint{2.711960in}{1.370397in}}%
\pgfpathlineto{\pgfqpoint{2.745584in}{1.362417in}}%
\pgfpathlineto{\pgfqpoint{2.777297in}{1.351761in}}%
\pgfpathlineto{\pgfqpoint{2.808246in}{1.338093in}}%
\pgfpathlineto{\pgfqpoint{2.838813in}{1.321186in}}%
\pgfpathlineto{\pgfqpoint{2.869380in}{1.300739in}}%
\pgfpathlineto{\pgfqpoint{2.900711in}{1.276042in}}%
\pgfpathlineto{\pgfqpoint{2.933571in}{1.246110in}}%
\pgfpathlineto{\pgfqpoint{2.968341in}{1.210128in}}%
\pgfpathlineto{\pgfqpoint{3.005786in}{1.166826in}}%
\pgfpathlineto{\pgfqpoint{3.048197in}{1.112892in}}%
\pgfpathlineto{\pgfqpoint{3.100543in}{1.041017in}}%
\pgfpathlineto{\pgfqpoint{3.254143in}{0.826685in}}%
\pgfpathlineto{\pgfqpoint{3.291587in}{0.781697in}}%
\pgfpathlineto{\pgfqpoint{3.323301in}{0.748254in}}%
\pgfpathlineto{\pgfqpoint{3.351193in}{0.723170in}}%
\pgfpathlineto{\pgfqpoint{3.376411in}{0.704520in}}%
\pgfpathlineto{\pgfqpoint{3.399336in}{0.691266in}}%
\pgfpathlineto{\pgfqpoint{3.420733in}{0.682353in}}%
\pgfpathlineto{\pgfqpoint{3.440984in}{0.677207in}}%
\pgfpathlineto{\pgfqpoint{3.460470in}{0.675445in}}%
\pgfpathlineto{\pgfqpoint{3.479193in}{0.676832in}}%
\pgfpathlineto{\pgfqpoint{3.497915in}{0.681350in}}%
\pgfpathlineto{\pgfqpoint{3.516637in}{0.689096in}}%
\pgfpathlineto{\pgfqpoint{3.535742in}{0.700408in}}%
\pgfpathlineto{\pgfqpoint{3.555610in}{0.715893in}}%
\pgfpathlineto{\pgfqpoint{3.576243in}{0.736038in}}%
\pgfpathlineto{\pgfqpoint{3.597640in}{0.761329in}}%
\pgfpathlineto{\pgfqpoint{3.620183in}{0.792816in}}%
\pgfpathlineto{\pgfqpoint{3.644254in}{0.831861in}}%
\pgfpathlineto{\pgfqpoint{3.669854in}{0.879399in}}%
\pgfpathlineto{\pgfqpoint{3.697365in}{0.937145in}}%
\pgfpathlineto{\pgfqpoint{3.727168in}{1.007063in}}%
\pgfpathlineto{\pgfqpoint{3.760027in}{1.092301in}}%
\pgfpathlineto{\pgfqpoint{3.796708in}{1.196321in}}%
\pgfpathlineto{\pgfqpoint{3.839883in}{1.328483in}}%
\pgfpathlineto{\pgfqpoint{3.897579in}{1.515944in}}%
\pgfpathlineto{\pgfqpoint{4.011059in}{1.886458in}}%
\pgfpathlineto{\pgfqpoint{4.052324in}{2.009288in}}%
\pgfpathlineto{\pgfqpoint{4.085948in}{2.100099in}}%
\pgfpathlineto{\pgfqpoint{4.114987in}{2.170007in}}%
\pgfpathlineto{\pgfqpoint{4.140587in}{2.223936in}}%
\pgfpathlineto{\pgfqpoint{4.163512in}{2.265331in}}%
\pgfpathlineto{\pgfqpoint{4.184145in}{2.296481in}}%
\pgfpathlineto{\pgfqpoint{4.202867in}{2.319373in}}%
\pgfpathlineto{\pgfqpoint{4.219679in}{2.335315in}}%
\pgfpathlineto{\pgfqpoint{4.235344in}{2.346058in}}%
\pgfpathlineto{\pgfqpoint{4.249864in}{2.352334in}}%
\pgfpathlineto{\pgfqpoint{4.263619in}{2.354911in}}%
\pgfpathlineto{\pgfqpoint{4.276610in}{2.354254in}}%
\pgfpathlineto{\pgfqpoint{4.289601in}{2.350524in}}%
\pgfpathlineto{\pgfqpoint{4.302974in}{2.343413in}}%
\pgfpathlineto{\pgfqpoint{4.317111in}{2.332219in}}%
\pgfpathlineto{\pgfqpoint{4.332013in}{2.316267in}}%
\pgfpathlineto{\pgfqpoint{4.347678in}{2.294838in}}%
\pgfpathlineto{\pgfqpoint{4.364490in}{2.266470in}}%
\pgfpathlineto{\pgfqpoint{4.382448in}{2.229979in}}%
\pgfpathlineto{\pgfqpoint{4.401935in}{2.183125in}}%
\pgfpathlineto{\pgfqpoint{4.422950in}{2.124140in}}%
\pgfpathlineto{\pgfqpoint{4.445493in}{2.051189in}}%
\pgfpathlineto{\pgfqpoint{4.469564in}{1.962411in}}%
\pgfpathlineto{\pgfqpoint{4.495546in}{1.854312in}}%
\pgfpathlineto{\pgfqpoint{4.508155in}{1.798285in}}%
\pgfpathlineto{\pgfqpoint{4.508155in}{1.798285in}}%
\pgfusepath{stroke}%
\end{pgfscope}%
\begin{pgfscope}%
\pgfpathrectangle{\pgfqpoint{0.687659in}{0.580556in}}{\pgfqpoint{3.820497in}{2.129775in}}%
\pgfusepath{clip}%
\pgfsetbuttcap%
\pgfsetroundjoin%
\pgfsetlinewidth{0.752812pt}%
\definecolor{currentstroke}{rgb}{0.000000,0.000000,0.000000}%
\pgfsetstrokecolor{currentstroke}%
\pgfsetstrokeopacity{0.600000}%
\pgfsetdash{{0.750000pt}{1.237500pt}}{0.000000pt}%
\pgfpathmoveto{\pgfqpoint{0.687659in}{1.379221in}}%
\pgfpathlineto{\pgfqpoint{4.508155in}{1.379221in}}%
\pgfusepath{stroke}%
\end{pgfscope}%
\begin{pgfscope}%
\pgfpathrectangle{\pgfqpoint{0.687659in}{0.580556in}}{\pgfqpoint{3.820497in}{2.129775in}}%
\pgfusepath{clip}%
\pgfsetbuttcap%
\pgfsetroundjoin%
\pgfsetlinewidth{0.752812pt}%
\definecolor{currentstroke}{rgb}{0.000000,0.000000,0.000000}%
\pgfsetstrokecolor{currentstroke}%
\pgfsetstrokeopacity{0.600000}%
\pgfsetdash{{0.750000pt}{1.237500pt}}{0.000000pt}%
\pgfpathmoveto{\pgfqpoint{2.597907in}{0.580556in}}%
\pgfpathlineto{\pgfqpoint{2.597907in}{2.710330in}}%
\pgfusepath{stroke}%
\end{pgfscope}%
\begin{pgfscope}%
\pgfsetrectcap%
\pgfsetmiterjoin%
\pgfsetlinewidth{0.803000pt}%
\definecolor{currentstroke}{rgb}{0.000000,0.000000,0.000000}%
\pgfsetstrokecolor{currentstroke}%
\pgfsetdash{}{0pt}%
\pgfpathmoveto{\pgfqpoint{0.687659in}{0.580556in}}%
\pgfpathlineto{\pgfqpoint{0.687659in}{2.710330in}}%
\pgfusepath{stroke}%
\end{pgfscope}%
\begin{pgfscope}%
\pgfsetrectcap%
\pgfsetmiterjoin%
\pgfsetlinewidth{0.803000pt}%
\definecolor{currentstroke}{rgb}{0.000000,0.000000,0.000000}%
\pgfsetstrokecolor{currentstroke}%
\pgfsetdash{}{0pt}%
\pgfpathmoveto{\pgfqpoint{4.508155in}{0.580556in}}%
\pgfpathlineto{\pgfqpoint{4.508155in}{2.710330in}}%
\pgfusepath{stroke}%
\end{pgfscope}%
\begin{pgfscope}%
\pgfsetrectcap%
\pgfsetmiterjoin%
\pgfsetlinewidth{0.803000pt}%
\definecolor{currentstroke}{rgb}{0.000000,0.000000,0.000000}%
\pgfsetstrokecolor{currentstroke}%
\pgfsetdash{}{0pt}%
\pgfpathmoveto{\pgfqpoint{0.687659in}{0.580556in}}%
\pgfpathlineto{\pgfqpoint{4.508155in}{0.580556in}}%
\pgfusepath{stroke}%
\end{pgfscope}%
\begin{pgfscope}%
\pgfsetrectcap%
\pgfsetmiterjoin%
\pgfsetlinewidth{0.803000pt}%
\definecolor{currentstroke}{rgb}{0.000000,0.000000,0.000000}%
\pgfsetstrokecolor{currentstroke}%
\pgfsetdash{}{0pt}%
\pgfpathmoveto{\pgfqpoint{0.687659in}{2.710330in}}%
\pgfpathlineto{\pgfqpoint{4.508155in}{2.710330in}}%
\pgfusepath{stroke}%
\end{pgfscope}%
\begin{pgfscope}%
\pgfsetbuttcap%
\pgfsetroundjoin%
\pgfsetlinewidth{1.003750pt}%
\definecolor{currentstroke}{rgb}{0.121569,0.466667,0.705882}%
\pgfsetstrokecolor{currentstroke}%
\pgfsetdash{{6.400000pt}{1.600000pt}{1.000000pt}{1.600000pt}}{0.000000pt}%
\pgfpathmoveto{\pgfqpoint{0.775159in}{2.588803in}}%
\pgfpathlineto{\pgfqpoint{0.969603in}{2.588803in}}%
\pgfusepath{stroke}%
\end{pgfscope}%
\begin{pgfscope}%
\pgftext[x=1.047381in,y=2.554775in,left,base]{\rmfamily\fontsize{7.000000}{8.400000}\selectfont \(\displaystyle R = 4\)}%
\end{pgfscope}%
\begin{pgfscope}%
\pgfsetbuttcap%
\pgfsetroundjoin%
\pgfsetlinewidth{1.003750pt}%
\definecolor{currentstroke}{rgb}{1.000000,0.498039,0.054902}%
\pgfsetstrokecolor{currentstroke}%
\pgfsetdash{{3.700000pt}{1.600000pt}}{0.000000pt}%
\pgfpathmoveto{\pgfqpoint{0.775159in}{2.453236in}}%
\pgfpathlineto{\pgfqpoint{0.969603in}{2.453236in}}%
\pgfusepath{stroke}%
\end{pgfscope}%
\begin{pgfscope}%
\pgftext[x=1.047381in,y=2.419208in,left,base]{\rmfamily\fontsize{7.000000}{8.400000}\selectfont \(\displaystyle R = 2\)}%
\end{pgfscope}%
\begin{pgfscope}%
\pgfsetbuttcap%
\pgfsetroundjoin%
\pgfsetlinewidth{1.003750pt}%
\definecolor{currentstroke}{rgb}{0.172549,0.627451,0.172549}%
\pgfsetstrokecolor{currentstroke}%
\pgfsetdash{{6.400000pt}{1.600000pt}{1.000000pt}{1.600000pt}}{0.000000pt}%
\pgfpathmoveto{\pgfqpoint{0.775159in}{2.317670in}}%
\pgfpathlineto{\pgfqpoint{0.969603in}{2.317670in}}%
\pgfusepath{stroke}%
\end{pgfscope}%
\begin{pgfscope}%
\pgftext[x=1.047381in,y=2.283642in,left,base]{\rmfamily\fontsize{7.000000}{8.400000}\selectfont \(\displaystyle R = 1\)}%
\end{pgfscope}%
\begin{pgfscope}%
\pgfsetbuttcap%
\pgfsetroundjoin%
\pgfsetlinewidth{1.003750pt}%
\definecolor{currentstroke}{rgb}{0.839216,0.152941,0.156863}%
\pgfsetstrokecolor{currentstroke}%
\pgfsetdash{{3.700000pt}{1.600000pt}}{0.000000pt}%
\pgfpathmoveto{\pgfqpoint{0.775159in}{2.182103in}}%
\pgfpathlineto{\pgfqpoint{0.969603in}{2.182103in}}%
\pgfusepath{stroke}%
\end{pgfscope}%
\begin{pgfscope}%
\pgftext[x=1.047381in,y=2.148076in,left,base]{\rmfamily\fontsize{7.000000}{8.400000}\selectfont \(\displaystyle R = 0.5\)}%
\end{pgfscope}%
\begin{pgfscope}%
\pgfsetrectcap%
\pgfsetroundjoin%
\pgfsetlinewidth{1.003750pt}%
\definecolor{currentstroke}{rgb}{0.580392,0.403922,0.741176}%
\pgfsetstrokecolor{currentstroke}%
\pgfsetdash{}{0pt}%
\pgfpathmoveto{\pgfqpoint{0.775159in}{2.046537in}}%
\pgfpathlineto{\pgfqpoint{0.969603in}{2.046537in}}%
\pgfusepath{stroke}%
\end{pgfscope}%
\begin{pgfscope}%
\pgftext[x=1.047381in,y=2.012509in,left,base]{\rmfamily\fontsize{7.000000}{8.400000}\selectfont \(\displaystyle R = 0.1\)}%
\end{pgfscope}%
\begin{pgfscope}%
\pgfsetbuttcap%
\pgfsetroundjoin%
\pgfsetlinewidth{2.007500pt}%
\definecolor{currentstroke}{rgb}{0.000000,0.000000,0.000000}%
\pgfsetstrokecolor{currentstroke}%
\pgfsetdash{{2.000000pt}{3.300000pt}}{0.000000pt}%
\pgfpathmoveto{\pgfqpoint{0.775159in}{1.910487in}}%
\pgfpathlineto{\pgfqpoint{0.969603in}{1.910487in}}%
\pgfusepath{stroke}%
\end{pgfscope}%
\begin{pgfscope}%
\pgftext[x=1.047381in,y=1.876460in,left,base]{\rmfamily\fontsize{7.000000}{8.400000}\selectfont \(\displaystyle u^{\prime}\)}%
\end{pgfscope}%
\end{pgfpicture}%
\makeatother%
\endgroup%

%% file: chapterLocalBif.tex
\chapter{Local Bifurcation}\label{Chapter:LocalBifPeriodic}
The success of the Armstrong-Painter-Sherratt adhesion model~\eqref{ArmstrongNd}
is that it can replicate the complicated patterns observed in cell-sorting
experiments\index{cell-sorting experiments}
\parencite{Armstrong2006}. In mathematical terms, these patterns are steady
states of equation~\eqref{ArmstrongNd}. Thus, understanding the
conditions, under which these steady states form and become stable is important.
Furthermore, understanding the set of steady-states is one of the first steps
toward understanding the equation's global attractor.

Up to this point, the steady states of equation~\eqref{ArmstrongNd} have
only been studied numerically and using linear stability analysis
\parencites{Armstrong2006}. Closely related to
equation~\eqref{ArmstrongNd} are the local and non-local chemotaxis
equations\index{chemotaxis equations}
\parencites{Othmer2002}{Hillen2007}{Hillen2009}. For both the local and
non-local chemotaxis equations a global bifurcation analysis, to understand
their steady-states, was carried out \parencites{Wang2013}{Xiang2013}. Inspired
by their results, we present here an exploration of the set of
non-homogeneous steady-state solutions of equation~\eqref{ArmstrongNd} in one
dimension, i.e.\ \eqref{Eqn:ArmstrongModel}.

Central to our analysis is the abstract bifurcation
theory\index{bifurcation theory} of Crandall and Rabinowitz\index{Crandall and
Rabinowitz} \parencite{Crandall1971,Rabinowitz1971}, in particular the local
bifurcation Theorem~\ref{Thm:CrandallMain} in \parencite{Crandall1971}
and the global bifurcation Theorem~\ref{Thm:UnilateralAbstractBifurcation}
in \parencite{Rabinowitz1971}.
Preceding the formulation of these general theorems, Crandall \etal\ studied
the set of solutions of non-linear Sturm-Liouville problems\index{Sturm-Liouville problem}
\parencites{Crandall1970}{Rabinowitz1970a}.  For linear Sturm-Liouville
eigenvalue problems it is well known that the eigenfunctions can be classified
by their number of zeros \parencite{Coddington1955}. Crandall \etal\
\cite{Crandall1971} showed, that under rather weak assumptions the same
classification holds for non-linear eigenvalue problems. In fact, each global
solution branch inherits the number of zeros of the eigenfunction spanning the
nullspace at the bifurcation point.
Furthermore, they showed that each global solution branch is unbounded and that
branches do not meet (since it is impossible for solutions of Sturm-Liouville
problems to have degenerate zeros). In \parencite{Rabinowitz1971}, Rabinowitz
formulates a general global bifurcation theorem, providing two alternatives for
the global structure of the solution branches.
They are either bounded, connecting two bifurcation
points, or are unbounded (see Theorem~1.3 in \parencite{Rabinowitz1971} and
Theorem~\ref{Thm:UnilateralAbstractBifurcation} in Chapter~2). This is now
known as the {\it Rabinowitz-alternative}\index{Rabinowitz-alternative}.  An
extension of the global bifurcation theorem to study so called unilateral
branches\index{unilateral branches} (sub-branches in only the positive or
negative direction of the eigenfunction at the bifurcation point) was
originally reported in \parencite{Rabinowitz1971}. The original proof contained
holes that were filled in by \parencites{Lopez2001}{Shi2009a}{Lopez-Gomez2016}.
Since the original formulation of these bifurcation theorems, similar theorems,
which apply in more general settings, have been developed. Here, we will use
bifurcation theorems that are applicable to Fredholm operators
\parencites{Lopez2001}{Shi2009a}{Lopez-Gomez2016} (see
\cref{Chapter:GlobalBifPeriodic}).

While, both the local and non-local chemotaxis model in \parencite{Wang2013} and
\parencite{Xiang2013} were formulated with no-flux boundary conditions, the
formulation of no-flux boundary conditions for
equation~\eqref{Eqn:ArmstrongModel} is more challenging due to the
complicated non-local structure of $\K[u]$. Hence we first study and analyze the
formation of non-homogeneous solutions in isolation of boundary considerations, and
we formulate equation~\eqref{Eqn:ArmstrongModel} on a circle (equivalent to
periodic boundary conditions). A more detailed discussion on the challenges of
construction of no-flux boundary conditions can be found in \cref{part:noflux}.

At the same time however, formulating our model on a circle gives rise to some
additional challenges. The most critical being that eigenvalues of the Laplacian
on $S^1_L$ need not be simple. This is a challenge, because bifurcations require
eigenvalues of odd multiplicity (the theorems in
\parencites{Shi2009a}{Lopez-Gomez2016} require simple eigenvalues).  This
challenge was previously observed by Matano, who studied nonlinear reaction
diffusion equations\index{reaction diffusion equations} on the circle
\parencite{matano1988}.  The solution of Matano \parencite{matano1988} was to
impose symmetry requirements on the
non-linear term, so that the equation would be $\O2$ equivariant i.e.\ invariant
under translations and reflections (see also \cref{subsec:Symmetries}). This
will also be our approach here.

Around the same time as Matano \parencite{matano1988}, Healey~\cite{Healey1988}
extended the Rabinowitz alternative to so called $G$-reduced
problems\index{$G$-reduced problems} ($G$ being a symmetry group), which is a
nonlinear bifurcation problem formulated on a Banach space whose elements are
fix-points of the symmetries defined by $G$ \parencite{Healey1988}.  Thus,
showing that certain symmetries\index{symmetries} persist along global
bifurcation branches\index{global bifurcation branches} (similar results are
often referred to as the equivariant branching lemma\index{equivariant branching
lemma} \parencite{golubitsky2003}).  Subsequently, these ideas were used in a
series of papers, which studied the condition under which solutions of
non-linear elliptic equations are classifiable by their number of zeros. The key
to these results, was to impose sufficient conditions on the nonlinear terms
such that the resulting equation was equivariant\index{equivariant} under
actions of $\O2 \times\ \Z_2$ ($\Z_2$ describes the action of the negative
identity, \ie, a reflection through the $x$-axis).
\parencites{Healey1991}{Healey1993}.  Intuitively, this symmetry requirement
ensured that the zeros of solutions were \textquote{frozen} (\ie, fixed
location), and thus the number of zeros is preserved along the global
bifurcation branch.  Furthermore, this result easily shows that global solution
branches do not meet.  Much more recently, Buono \etal\ used $\O2$ equivariance
to compare the accessible bifurcations in a non-local hyperbolic model of
swarming\index{swarming} and the equation's formal parabolic limit
\parencite{Buono2015}.

In a similar spirit, we will show that the steady-state
equation~\eqref{Eqn:StSt} of equation~\eqref{Eqn:ArmstrongModel} is indeed
equivariant under actions of $\O2$.  Using the properties of the non-local term
$\K[u]$, we will then show that this leads to \textquote{frozen} maxima and
minima (equivalently frozen zeros of the derivative $u^{\prime}$). Since we
consider equation~\eqref{Eqn:StSt} on a periodic domain, we can prescribe the
location of one maxima (or minima) without restricting possible solutions. This
ensures that only simple eigenvalues occur in our subsequent analysis. In
Chapter~\ref{Chapter:GlobalBifPeriodic} we prove a global bifurcation result for
the steady-state solutions of the non-local cell-cell adhesion equation for the
first bifurcation branch. Finally we discuss how this proof could be extended to
all bifurcation branches.
\section{The Abstract Bifurcation Problem}
The steady states of equation~\eqref{Eqn:3:AdhModel} are solutions of
the following non-local equation.

\begin{subequations}\label{Eqn:StSt}
\begin{equation}\label{Eqn:StSt2}
    u^{\prime\prime}(x) = \alpha \Bigl( u(x) \K[u(x)](x) \Bigr)^{\prime} \quad
    \mbox{in } S^1_L.
\end{equation}
Some typical steady state solutions of this equation are shown in
\cref{Fig:TypicalSinglePopulationSolutions}.
Due to the mass conserving property of equation~\eqref{Eqn:3:AdhModel}
we find families of steady states that are parametrized by their average
\begin{equation}\label{Eqn:IntConst}
    \mean{u} = \Avg[u] = \frac{1}{L} \int_{0}^{L} u(x) \dd x.
\end{equation}
\end{subequations}
We understand this equation as constraint, which ensures that bifurcations arise
in a  given family of solutions which have the same mass (note that the mass is
$L \bar u$). We define the following function space
\[
    H^2_{\Bd}(S^1_L) \coloneqq \lcb u \in H^2(S^1_L) : \Bd[u, u^{\prime}] = 0 \rcb,
\]
where the boundary operator $\Bd[\cdot, \cdot]$ was defined in
equation~\eqref{Eqn:DefnPeriodicBoundaryConditions} for periodic boundary
conditions. For given $\bar u >0$ we define the following operator
\begin{subequations}\label{Eqn:F}
\begin{equation}
    \F\ \colon\ \R \times H^2_{\Bd}(S^1_L) \mapsto L^2(S^1_L) \times \R,
\end{equation}
\begin{equation}\label{Eqn:OpF}
    \F[\alpha, u] =
        \begin{bmatrix}
            \lb - u^{\prime} + \alpha u \K[u] \rb^{\prime} \\
            \Avg[u] - \mean{u}
        \end{bmatrix}.
\end{equation}
\end{subequations}
We will need the Fr\'echet derivative\index{Fr\'echet derivative} of $\F$.
\begin{lemma}\label{Lem:FrechF}
    The Fr\'echet derivative
    $\D_u \F\ \colon\ \R \times H^2_{\Bd} \mapsto \Lin(H^2_{\Bd}, L^2 \times \R)$
    of the operator $\F$ is given by
    \begin{equation}\label{Eqn:DF}
        \D_u \F(\alpha, u)[w] =
            \begin{pmatrix}
                \lsb -w^{\prime} + \alpha \lb
                    u \K_{h}[w] +
                    w \K[u] \rb \rsb^{\prime} \\
                \Avg[w]
            \end{pmatrix}.
    \end{equation}
    where
    \begin{equation}\label{Eqn:Kh}
        \K_{h}[w] = \int_{-1}^{1} h^{\prime}(u(x + r)) w(x + r) \Omega(r) \dd r.
    \end{equation}
\end{lemma}
\begin{proof}
    Let $u, v, w \in H^2_{\Bd}$ and we compute
    \begin{equation}\label{Fdiff}
    \F[\alpha, u + w] - \F[\alpha, u] - \D_{u} \F(\alpha, u)[w]
    \end{equation}
    The second component is
    \[ \Avg[u+w] -\bar u -(\Avg[u]-\bar u) -\Avg[w]=0.\]
    In the first component of \eqref{Fdiff} the local derivative terms cancel directly:
    \[ (-(u+w)' -(-u') - (-w'))' = 0. \]
    It remains to check the non-local terms in the first component of \eqref{Fdiff}:
    \[
    \begin{split}
        \Bigl(\F[\alpha, u + w] &- \F[\alpha, u] - \D_{u} \F(\alpha, u)[w]\Bigr)_1  \\
        & = \alpha\Bigl(u \K[u+w] + w\K[u+w] - u \K[u]-u\K_h[w]+w\K[u]\Bigr)'\\
        &=         \alpha \left[ \lb
                \underbrace{%
                    u \int_{-1}^{1} \lb h(u + w) - h(w) -
                    h^{\prime}(u) w \rb \Omega(r) \dd r}_{\ROM{1}}
                    \rb^{\prime} \right. \\
                    &\quad + \left. \lb
                    \underbrace{%
                    w \int_{-1}^{1} \lb h(u + w) - h(u) \rb \Omega(r) \dd r}_{\ROM{2}} \rb^{\prime}
                    \right].
    \end{split}
    \]

    We require an $L^2$ estimate of
    the previous terms. For this we consider the two terms separately (we denote
    them by $\ROM{1}$ and $\ROM{2}$ respectively).

    \[
        \abs{{(\ROM{1})}^{\prime}}_{2} \leq \norm{\ROM{1}}_{H^1} \leq
            \norm{u}_{H^1} \norm{w}_{H^1}
            \norm*{\frac{h(u + w) - h(u) - h^{\prime}(u) w}{w}}_{H^1},
    \]
    where we use the Banach algebra property of
    $H^1$. Now, by assumption~\ref{NatureForceAssumption:1}\ $h(\cdot) \in \C^2$,
    hence we have that
    \[
        \lim_{\norm{w} \to 0} \norm*{%
            \frac{h(u + w) - h(u) - h^{\prime}(u) w}{w}
            }_{H^1} = 0.
   \]
    For the second term, we proceed similarly
    \[
        \abs{{(\ROM{2})}^{\prime}}_{2} \leq \norm{\ROM{2}}_{H^1} \leq
            \norm{w}_{H^1} \norm{h(u + w) - h(u)}_{H^1},
    \]
    which goes to 0 as $\|w\|\to 0$, since $h(u)$ is continuous.
    Together we find
    \begin{eqnarray*}
        &&\frac{%
            \abs{\F[\alpha, u + w] - \F[\alpha, u] - \D_{u} \F(\alpha, u)[w]}_{2}
        }{\norm{w}_{H^2}} \\
        &\leq & \frac{\|w\|_{H^1}}{\|w\|_{H^2}} \alpha \|u\|_{H^1} \left(\left\|\frac{h(u+w) - h(u) - h'(u+w)}{w}\right\|_{H^1} + \|h(u+w) - h(w) \|_{H^1} \right)\\
        &\to& 0 \quad\mbox{ as } \quad \|w\|_{H^2} \to 0.
    \end{eqnarray*}
    Hence $\D_u \F(\alpha, u)$ is the Fr\'echet derivative of $\F$.
\end{proof}

We then prove a series of properties of $\F$ that allow us to apply the
bifurcation theorems from \cref{Chapter:Preliminaries}.

\begin{lemma}\label{Lem:PropF}
    For each $\R \ni \mean{u} > 0$, the operator $\F$ be defined as in
    equation~\eqref{Eqn:F} has the following properties:
    \begin{enumerate}
        \item $\F[\alpha, \mean{u}] = 0$ for all $\alpha \in \R$.
        \item The first component of $\F$ maps into $L^2_0(S^1_L)$
            (defined in \cref{Eqn:AvgZeroFun}).
        \item $\D_u \F(\alpha, u)$ is Fredholm\index{Fredholm operator} with index 0, for
            each $\alpha \in \R$.
        \item $\F[\alpha, u]$ is $\C^1$ smooth in $u$.
        \item $\D_{\alpha u} \F(\alpha, u)$ exists and is continuous in $u$.
    \end{enumerate}
\end{lemma}
\begin{proof}
    \begin{enumerate}[leftmargin=\parindent, labelindent=*]
        \item We note that $\K[\mean{u}] = 0$, hence the conclusion follows.
        \item This is easily seen by integrating the equation, and using the
            periodic boundary conditions $\Bd$.
        \item The Fr\'echet derivative of $\F$ with respect to $u$ was found in
            \cref{Lem:FrechF} to be
            \[
                \D_u \F(\alpha, u)[w] =
                    \begin{pmatrix}
                        \lb -w^{\prime} + \alpha \lb u \K_{h}[w] + w \K[u] \rb \rb^{\prime} \\
                        \Avg[w]
                    \end{pmatrix}.
            \]
            We split this operator as follows
            \[
                \D_u \F(\alpha, u)[w] = \Top_1(\alpha, u)[w] + \Top_2(\alpha,
                u)[w],
            \]
            where
            \[
                \Top_1(\alpha, u)[w] = \begin{pmatrix}
                                        - w^{\prime\prime} \\ 0
                                    \end{pmatrix},
            \]
            and
            \[
                \Top_2(\alpha, u)[w] = \begin{pmatrix}
                    \alpha \lb u \K_{h}[w] + w \K[u] \rb^{\prime} \\
                                        \Avg[w]
                                    \end{pmatrix}.
            \]

            $\Top_1(\alpha, u)$ is Fredholm with index 0, by
            \cref{Defn:Laplacian}.

            The operator $\Top_2(\alpha, u)$ is compact. First, consider the first
            component of $\Top_2$.
            Now ${\D\ \colon\ H^1 \mapsto L^2}$, ${w \mapsto w^\prime}$ is continuous.
            Then let $\Top(u)[w] = u \K_h[w] + w \K[u]$, where
            $u, w \in H^2_{\Bd}$. Since $h^{\prime}(u) \in \C^1$ we have that
            $\K_{h}[w] \in H^2_{\Bd}$. Since $H^2_{\Bd}$ is a Banach algebra we
            have that $u \K_{h}[w] \in H^2_{\Bd}$ and $w \K[u] \in H^2_{\Bd}$.
            Since, $H^2 \ssubset \C^1$, we can conclude that the first component
            is given by the composition $\D \circ \Top(u)$ and hence is compact.

            For the second component we let $\seq{w_n} \subset H^2$ be bounded,
            then $\Avg[w_n] = \frac{1}{L} \int_{0}^{L} w_n(x) \dd x$,
            but as $\seq{w_n} \subset L^{\infty}$, we have that
            $\Avg[w_n]$ is a bounded sequence in $\R$ and
            thus has a convergent subsequence, making $\Avg$ compact.

            Finally, we recall the well known results that the compact
            perturbation of a Fredholm operator is Fredholm with the same index
            (see \cref{Thm:FredholmCompactPerturbation}).

            Hence, $\D_u \F(\alpha, u)$ is Fredholm with index 0.
        \item  For this we have to check that the mapping
            \[
                 u \mapsto \D_u \F(\alpha, u)
            \]
            is continuous.
            The $u$-dependence in the linearization $\D_u \F(\alpha, u)$ arises in the
            integral terms. For each $w\in H^2$ the map
            \[ u\mapsto \D_u \F(\alpha, u)[w] =\left(\begin{array}{c}
                -w'' +\alpha(u\K_h[w] + w \K[u])' \\ \Avg[w]\end{array} \right)
            \]
            is continuous in $u$, since $\K_h[w]$ is bounded and $\K[u]$ is continuous by
            \cref{Lem:KProp}.

        \item We simply compute to find
            \[
                \D_{\alpha u} \F(\alpha, u) =
                \begin{pmatrix}
                    \lb u \K_{h}[w] + w \K[u] \rb^{\prime} \\
                        0
                \end{pmatrix}.
            \]
            Its continuity follows from item (4).
    \end{enumerate}
\end{proof}
\begin{remark}
    Note that \cref{Lem:PropF} implies that $\F$
    satisfies properties~\ref{Assumption:1},~\ref{Assumption:2}
    and~\ref{Assumption:3} from Section~\ref{s:Bifurcation}.
\end{remark}
\section{Symmetries and equivariant flows\index{equivariant flows}}\label{subsec:Symmetries}
For the following bifurcation analysis we will require a detailed understanding
of symmetries of our solutions. We use group theory to describe these
symmetries. We obtain the following result.
\begin{lemma}\label{Lem:Fequiv}
    The operator $\F$ defined in equation~\eqref{Eqn:OpF} is
    equivariant under the actions of $\O2$, \ie,
    \[
        \F[\alpha, \gamma u] = \gamma \F[\alpha, u], \quad
            \forall \gamma \in \O2.
    \]
\end{lemma}
\begin{remark}
    The group $\O2$ is generated by $\SO2$ and a reflection. In more
    detail $\SO2$ can be represented in $\R^2$ by rotations
    \[
        \sigma_{\theta} = \begin{pmatrix} \cos \theta & - \sin \theta \\
                                      \sin \theta & \cos \theta
        \end{pmatrix}, \quad\theta \in {\lsb0, 2\pi\rb}
    \]
    and a reflection given by
    \[
        \rho = \begin{pmatrix} 1 & 0 \\ 0 & -1 \end{pmatrix}.
    \]
    It is easy to see that this group is compact and hence proper. Here we
    represent the group by its action on functions on $S^1_L$ as
    \begin{equation}\label{Eqn:ActionO2}
    \begin{split}
        \sigma_{a} u(x) &= u(x - a) \quad a\in [0,L] \quad\mbox{translations\index{translations}}  \\
        \rho u(x)       &= u(L - x) \quad\mbox{reflection\index{reflection}}.
    \end{split}
    \end{equation}
    In the following, we will denote reflection about $a \in S^1_L$ by
    \[
        \rho_{a} u(x) = u(2a - x),
    \]
    Since $\rho_a = \sigma_{L-2a}\rho$ this operation is in $\O2$.
\end{remark}

\begin{figure}[!ht]\centering
    \input{SampleFunctions.pgf}
    \caption[Examples of the symmetries of $\O2$]
    {Examples of the actions of $\sigma_a$ (left) and $\rho_a$ and
    (right). Here $a = \ifrac{L}{4}$. In both subplots, the original function is
    shown in dashed black, while the shifted and respectively reflected
    functions are solid.
    }\label{Fig:SymExmp}
\end{figure}
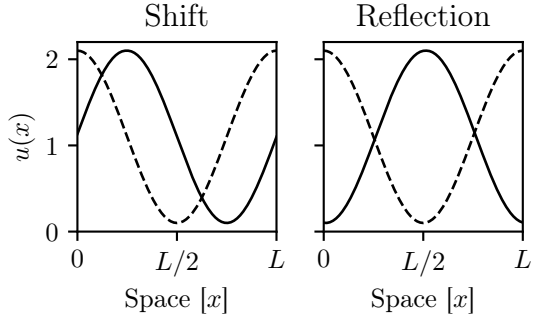

For the proof of this lemma we require the following lemma describing how the
non-local operators behaves under actions of $\O2$.
\begin{lemma}\label{Lem:NLSym}
    Let the non-local operator $\K[u]$ be defined as in
    \cref{Defn:NonLocalOperator}, then
    \[
        \K[\sigma_{\theta} u] = \sigma_{\theta} \K[u], \qquad
        \K[\rho_a u] = - \rho_a \K[u],
    \]
    and for the non-local curvature operator, we have that
    \[
        \lb\K[\sigma_{\theta} u]\rb^{\prime} = \sigma_{\theta} \lb \K[u]
        \rb^{\prime}, \qquad
        \lb\K[\rho_a u]\rb^{\prime} = \rho_a \lb \K[u] \rb^{\prime}.
    \]
\end{lemma}
\begin{proof}
    The results for $\sigma_{\theta}$ are trivial. For $\rho_a$ we first deal
    with $\K[u](x)$. Then
    \[
    \begin{split}
        \K[\rho_a u](x) &=
            \int_{0}^{1} \lsb h(u(2a - x - r)) - h(u(2a - x + r)) \rsb \omega(r) \dd r
            \\ &= - \K[u](2a - x) = - \rho_a \lb \K[u] \rb.
    \end{split}
    \]
    Second, we show the same for ${\K[u]}^{\prime}$ using a simple change of
    variables
    \[
    \begin{split}
        {\K[\rho_a u]}^{\prime} &= \lb \int_{-1}^{1} - h(u(2a - x - r)) \Omega(r) \dd r
        \rb^{\prime}  \\
        &= \lb \int_{-1}^{1} h(u(2a - x + y)) \Omega(y) \dd y
        \rb^{\prime} = {\K[u]}^{\prime}(2a - x) = \rho_a {\K[u]}^{\prime}.
    \end{split}
   \]
\end{proof}
Now we can complete the proof of \cref{Lem:Fequiv}.
\begin{proof}[Proof of \cref{Lem:Fequiv}]
    The elements of $\SO2$ are given by the translations. As
    \[
        \frac{\partial}{\partial x} u(x - a) =
            \frac{\partial}{\partial (x - a)} u(x - a),
    \]
    it is trivial to see that $\F$ is equivariant under actions of elements in
    $\SO2$. Note that to obtain all elements in $\O2$ we only need the
    reflection through $a = \ifrac{L}{2}$ as defined in
    equation~\eqref{Eqn:ActionO2}, we find that
    \[
        \lb \rho v(x) \rb^{\prime} = \lb v(2a - x) \rb^{\prime}
            =  - v^{\prime}(2a - x) = - (\rho v^{\prime})(x),
    \]
    and hence ${(\rho v(x))}^{\prime\prime} = (\rho v^{\prime\prime})(x)$. For
    the non-local term we apply \cref{Lem:NLSym}
    with $a = \ifrac{L}{2}$. Then substitute these into
    equation~\eqref{Eqn:OpF} and we obtain the result.
\end{proof}
The linearization of the steady state equation~\eqref{Eqn:StSt2} is given by
$\D_u \F(\alpha, u)[w] =0$.  We saw in \cref{Lem:EfuncK}  that the
non-local operator $\K$ maps $\sin(\frac{n\pi x}{L})$ to $\cos(\frac{n\pi
x}{L})$ and vice versa. Hence, as we will show below, the eigenfunctions of the
linearization $\D_u \F$ are sine and cosine functions. As we are working on the
circle of length $L$, a sine function can easily be shifted to a cosine function
with an appropriate phase shift. For the bifurcation analysis we need to remove
this symmetry, which we achieve by stipulating reflection symmetry $u(x) =
u(L-x)$. This makes the corresponding eigenspace one dimensional and we can
apply the abstract bifurcation theory outlined before.  Once these bifurcation
branches are identified, we can shift the solutions around the circle to obtain
other solutions, that are not reflection symmetric through the domain's center.

For this reason, we define new function spaces
\begin{equation}\label{Eqn:SymH2}
    H^2_{P} \coloneqq \lcb u \in H^2_{\Bd}(S^1_L)\ \colon\ u(x) = u(L - x) \rcb,
\end{equation}
and $L^2_P$ accordingly. It is then easy to see that the operator $\F$ now maps
$\R \times H^2_{P} \mapsto L^2_{P} \times \R$, since if $u \in H^{2}_{P}$ we find
that
\[
    \F[\alpha, u] = \F[\alpha, \rho u] = \rho \F[\alpha, u].
\]
\section{Generalized Eigenvalues\index{generalized eigenvalue} of $\F$}\label{subsec:NonLinearEigenvalues}
Due to \cref{Lemma:BifPointsAreEigenvalues} we have that possible bifurcation
points\index{bifurcation point} are those that are generalized eigenvalues of
the linearization of $\F$ evaluated at the trivial solution $(\alpha,
\mean{u})$. The main
bifurcation result of Crandall and Rabinowitz \cref{Thm:CrandallMain} is
formulated for a trivial steady state at $0$. Hence we need to shift our
solutions by $-\bar u$ to obtain a trivial steady state at $0$.
\begin{equation}\label{Eqn:ChangeOfVariables}
    v(x) \coloneqq u(x) - \mean{u}.
\end{equation}
Under this change of variable the operator $\F$ on $H^2_P$ becomes
\begin{equation}\label{Eqn:FChange}
    \F[\alpha, v] =
        \begin{bmatrix}
            \lb - v^{\prime} + \alpha (v + \mean{u}) \K[v] \rb^{\prime} \\
            \Avg[v]
        \end{bmatrix},
\end{equation}
and its linearization becomes
\[
    \D_v \F(\alpha, v)[w] = \begin{pmatrix}
        \lb -w^{\prime} + \alpha \lb (v + \mean{u}) \K_{h}[w] +
            w \K[v] \rb \rb^{\prime} \\
        \Avg[w]
    \end{pmatrix}.
  \]
Note that we treat the linearizations as an operator family\index{operator family} indexed by $\alpha \in \R$. In
the next two lemmas we characterize the values of $\alpha$ which lie in
the generalized spectrum $\Sigma\lb \D_v \F(\alpha, 0) \rb$.
\begin{lemma}\label{Lem:OpEvals}
    Let the operator $\F$ on $H^2_P$ be as defined in~\eqref{Eqn:FChange},
    and let $\alpha > 0$.  Its Fr\'echet derivative has been shown to exist in
    \cref{Lem:PropF}. Assume that $M_n(\omega) > 0$ (\ie\ the
    Fourier sine coefficients of $\omega$ as defined in
    equation~\eqref{Defn:Mn}), and define
    \begin{equation}\label{Eqn:NonLinEvals}
        \alpha_{n} \coloneqq
            \frac{n \pi}{\mean{u} L M_n(\omega) h^{\prime}(\mean{u})}
            \quad\mbox{for}\ n \in \N.
    \end{equation}
    Then, we have that
    \[
        \dim \Null[\D_v \F(\alpha_n, 0)] = 1.
    \]
    Thus, the generalized spectrum\index{generalized spectrum} of the linearization is given by the $\alpha_n$:
    \[
        \Sigma\lb \D_v \F(\alpha, 0) \rb
            = \lcb \alpha_n\ \colon\ n \in \N \setminus \lcb 0 \rcb \rcb.
    \]
\end{lemma}
\begin{remark}
    Note that $\alpha = 0$ is not an eigenvalue of $\D_v \F(\alpha, 0)$, since
    in this case, the only solution of equation~\eqref{Eqn:PerPrb} is
    the zero solution.
\end{remark}
\begin{remark}
    Note that since we assume that $L \geq 2$ we have that $M_1(\omega) > 0$,
    since then we have that $\sin\lb\frac{2\pi x}{L}\rb > 0$ on $(0,1)$ and
    $\omega(r) \geq 0$ by assumption. Thus, there is always one such bifurcation
    point.
\end{remark}
\begin{proof}
    The nullspace of $\D_v \F (\alpha, 0)$ is given by the solution of the
    following equation
    \begin{equation}\label{Eqn:PerPrb}
    \left\{
        \begin{array}{@{}ll@{}}
            -w^{\prime\prime} + \alpha \mean{u} h^{\prime}(\mean{u})
                \lb \int_{-1}^{1} w(x+r) \Omega(r) \dd r \rb^{\prime} = 0
            \quad &\mbox{in } \lsb 0, L \rsb  \\
            \Bd[w, w^\prime] = 0, \quad \Avg[w] = 0.
    \end{array}
    \right.
    \end{equation}
    We solve this system using an eigenfunction ansatz.
    \[
        w(x) = a_0 + \sum_{n=1}^{\infty} a_n \cos\lb\frac{2 n \pi x}{L}\rb
            + \sum_{n=1}^{\infty} b_n \sin\lb\frac{2 n \pi x}{L}\rb,
    \]
    then because $\Avg[w] = 0$ we have that $a_0 = 0$ and as $w \in H^2_{P}$ has reflection symmetry,
    we have that $b_n = 0, \forall n \in \N$. Hence
    \[
        w(x) = \sum_{n=1}^{\infty} a_n \cos\lb\frac{2 n \pi x}{L}\rb,
    \]
    which is then substituted into equation~\eqref{Eqn:PerPrb} to
    obtain
    \[
        \sum_{n=1}^{\infty} a_n \frac{2\pi n}{L} \cos\lb\frac{2 n \pi x}{L}\rb
            \lcb 2 \alpha \mean{u} h^{\prime}(\mean{u})
            \int_{0}^{1} \sin\lb\frac{2 n \pi x}{L}\rb \omega(r)\dd r
            - \frac{2 n \pi}{L} \rcb = 0.
    \]
    To find a non trivial solution we require that
    \[
        \lcb \alpha \mean{u} h^{\prime}(\mean{u})
        \int_{0}^{1} \sin\lb\frac{2 n \pi x}{L}\rb \omega(r)\dd r
            - \frac{n \pi}{L} \rcb  = 0.
    \]
    Hence, $\D_v \F(\alpha, 0)$ is not an isomorphism whenever $\alpha$ equals
    one of the following
    \[
        \alpha_{n} = \frac{n \pi}{L \mean{u} M_n(\omega) h^{\prime}(\mean{u})},
    \]
    where $M_n(\omega)$ is defined in equation~\eqref{Defn:Mn}.
    The corresponding eigenfunctions are
    \[ e_n (x) = \cos \left(\frac{2\pi n x}{L}\right),\ n=1,2,\dots.\]
\end{proof}
\begin{remark}
    Note that the values $\alpha_n$ found in \cref{Lem:OpEvals}
    are exactly the values at which an eigenvalue $\lambda$ of the linear
    operator
    \begin{equation}
        v^{\prime\prime} - \alpha \mean{u} \lb \K[u] \rb^{\prime} = \lambda v,
    \end{equation}
    crosses through $0$.
\end{remark}
\begin{lemma}\label{Lemma:Transversality}
    Let $\alpha_n$ be a generalized eigenvalue of $\D_v \F(\alpha, 0)$ as found
    in \cref{Lem:OpEvals}. Then, we have that
    \[
        \D_{\alpha v} \F (\alpha_n, 0)[e_n] \notin \Range[\D_v \F(\alpha_n, 0)],
    \]
    where $e_n(x)$ is the eigenfunction corresponding to $\alpha_n$
    from~\eqref{Eqn:NonLinEvals}.
\end{lemma}

\begin{proof}
    We proceed by contradiction. We assume that
    \begin{equation}\label{Eqn:Transversality}
        \D_{\alpha v} \F (\alpha_n, 0)[e_n] \in \Range[\D_v \F(\alpha_n, 0)].
    \end{equation}
    That means that there is $z \in H^2_{\Bd}$ such that
    \begin{equation}\label{Eqn:TransHelp}
        \D_v \F(\alpha_n, 0)[z] =
        \D_{\alpha v} \F(\alpha_n, 0)[e_n] = \begin{bmatrix}
        \mean{u} {\K_{h}[e_n]}^{\prime} \\ 0 \end{bmatrix}.
    \end{equation}
    Then equation~\eqref{Eqn:TransHelp} is equivalent to
    \begin{equation}\label{Eqn:TransPrb}
        \left\{
            \begin{array}{@{}ll@{}}
                - z^{\prime\prime} + \alpha_n \mean{u} h^{\prime}(\mean{u})
                    \lb \int_{-1}^{1} w(x + r) \Omega(r) \dd r \rb^{\prime} =
                   \mean{u} {\K_{h}[e_n]}^{\prime}
                \quad &\mbox{in } \lsb 0, L \rsb  \\
                \Bd[z, z^\prime] = 0, \quad \Avg[z] = 0.
        \end{array}
        \right.
    \end{equation}
    We note that (refer to \cref{Lem:EfuncKp}),
    \[
        {\K_{h}[e_n]}^{\prime} = -  h^{\prime}(\mean{u})
        \lb \frac{4 \pi n}{L} \rb \cos\lb \frac{2 \pi n x}{L} \rb M_n(\omega).
    \]
    Like in the previous lemma, we once again use the ansatz
    \[
        z(x) = \sum_{j=1}^{\infty} t_j \cos\lb \frac{2 \pi n x}{L}\rb.
    \]
    Upon substitution into equation~\eqref{Eqn:TransPrb}, we obtain
    \[
    \begin{split}
        \sum_{j=1}^{\infty} t_j \lb \frac{4 \pi j}{L}\rb \cos\lb \frac{2\pi j x}{L} \rb
            &\lcb M_j(\omega) \alpha_n\bar u h'(\bar u)  - \frac{\pi j}{L} \rcb \\
            &= - \mean{u} h^{\prime}(\mean{u})
            \lb \frac{4 \pi n}{L} \rb \cos\lb\frac{2 \pi n x}{L}\rb M_n(\omega).
    \end{split}
    \]
    Then this equation has a solution if and only if $t_j=0, \forall j \neq n$ and
    \[
        t_n \lcb \alpha_n M_n(\omega) \mean{u} h^{\prime}(\mean{u}) - \frac{\pi n}{L} \rcb
            = - \mean{u} h^{\prime}(\mean{u}) M_n(\omega)
    \]
    is satisfied. But
    $\alpha_n = \frac{\pi n}{\mean{u} h^{\prime}(\mean{u}) M_n(\omega) L}$,
    hence the term on the left hand side is zero, and we have found a
    contradiction.
\end{proof}

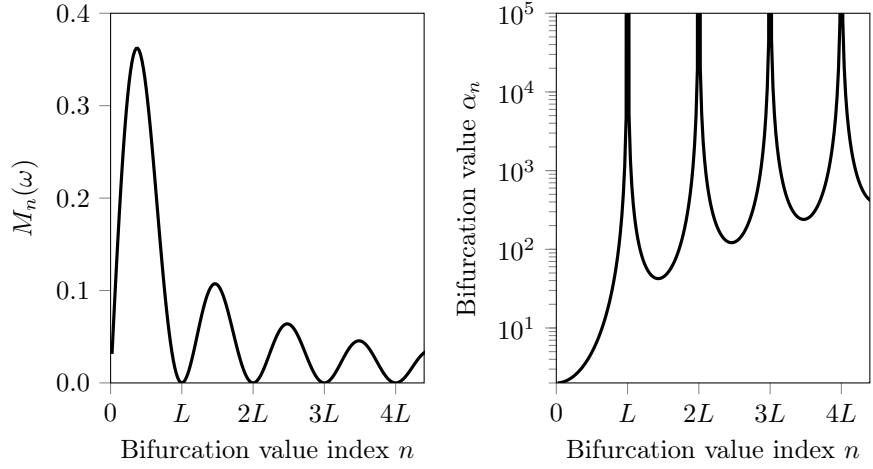
\begin{figure}\centering
    \input{BifurcationValuesExample.tikz}
    \caption[Bifurcation values]{Left: Shows the function $M_n(\omega)$ defined
    in equation~\eqref{Defn:Mn} for uniform $\omega$
    (see~\ref{OmegaChoice:1}). When $L = n$ the term $M_n(\omega) = 0$.
    Right: The bifurcation values $\alpha_n$ as given in
    equation~\eqref{Eqn:NonLinEvals}. The bifurcation values blow up
    whenever $L = n$. Note that even if $L \neq n$, we still observe very large
    bifurcation values for $n \sim k L$ ($k \in \N$).
    }\label{Fig:TypBifValues}
\end{figure}

\begin{example}
    Suppose that $h(u) = u$, and that $\omega$ is uniform. Using the precise
    form for $M_n(\omega)$ from \cref{Example:UniformOmega1} the bifurcation
    values are
    \[
        \alpha_n = \frac{2}{\bar{u}} \lb \frac{L}{n\pi}
        \sin\lb\frac{n\pi}{L}\rb\rb^{-2}.
    \]
    This means that $\alpha_n$ is a function of $L, \bar{u}, n$. We investigate
    now the behaviour of $\alpha_n$ as a function of $L$ and $\bar{u}$.
    The domain of $L$ is $[2, \infty)$. When $L = 2$
    \[
        \alpha_n = \begin{cases}
                        \dfrac{n^2 \pi^2}{2 \bar{u}} &\mbox{if } n\ \mbox{even}\\
                        0 &\mbox{else}
                   \end{cases},
    \]
    and when $L \to \infty$, we have that $\alpha_n$ decreases towards
    $\ifrac{2}{\bar{u}}$. When $\bar{u} \to \infty$ we find that $\alpha_n \to 0$.
    That means that increasing the domain size (equivalently decreasing the sensing
    radius) or increasing the total mass in the system, decreases the threshold
    adhesion strength required for aggregation.
    Finally, we note that the map $n \mapsto \alpha_n$ is not monotone, because
    $M_n(\omega)$ approaches zero when $n$ is close to a multiple of $L$ (see
    the left side of \cref{Fig:TypBifValues}).
\end{example}
\section{Local Bifurcation Result}\label{subsec:LocalBifurcation}
Based on the previous lemmas we can now formulate the local bifurcation\index{local bifurcation} result.
\begin{theorem}\label{Thm:LocBif}
    Let $\F$ be given on $H^2_P$  as in~\eqref{Eqn:F}. Then by
    \cref{Lem:PropF}, \cref{Lem:OpEvals} and
    \cref{Lemma:Transversality} all requirements of \cref{Thm:CrandallMain} are
    satisfied. Then there are continuous functions
    \begin{equation}\label{Eqn:LocBif}
        (\alpha_k(s), u_k(s)) : (-\delta_k, \delta_k) \to \R \times
        H^2_{P},
    \end{equation}
    with $\alpha_k(0) = \alpha_k$ such that
    \begin{equation}\label{Eqn:LocBifFun}
        u_k(x, s) = \mean{u} + s \alpha_k \cos\lb\frac{2 \pi k x}{L}\rb + o(s),
    \end{equation}
    where $(\alpha_k(s), u_k(s))$ is a solution of
    the steady state equation~\eqref{Eqn:StSt} and all
    non-trivial solutions near the bifurcation point $(\alpha_k, \mean{u})$
    lie on the curve $\Gamma_k = (\alpha_k(s), u_k(s))$.
\end{theorem}
\begin{proof}
    To be able to apply \cref{Thm:CrandallMain} we set $U = H^2_{P}$ and
    $V = L^2_{0, P} \times \R$ ($L^2_{0, P}$ is the sub-space of functions
    in $L^2$ with the symmetry as in equation~\eqref{Eqn:SymH2} and
    average zero as in equation~\eqref{Eqn:AvgZeroFun}) and $W = H^2_{P}$.
    Then the operator defined in
    equation~\eqref{Eqn:FChange} satisfies all the properties
    required by \cref{Thm:CrandallMain}. These properties are proved by
    \cref{Lem:PropF}, \cref{Lem:OpEvals} and
    \cref{Lemma:Transversality}. Finally, we revert the change of variables
    given in equation~\eqref{Eqn:ChangeOfVariables}.
\end{proof}

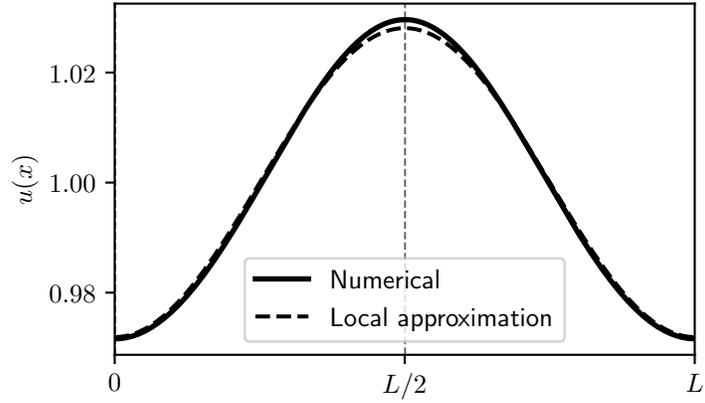
\begin{figure}\centering
 \input{local_single_peak.pgf}
 \caption[Agreement local bifurcation result to numerical solutions]
 {Comparison of numerical solution to the solution approximation
 given in equation~\eqref{Eqn:LocBifFun} near the first
 bifurcation point i.e.\ $(\alpha_1, \mean{u})$. The correct value of $s$
 for equation~\eqref{Eqn:LocBifFun} was
 estimated using the asymptotic expansion introduced in
 \cref{sec:BifurcationType}. Here $L = 5$.}\label{Fig:TypBifSoln}
\end{figure}

\section{Summary}\label{sec:LocalBifurcationSummary}

In this chapter, we identified the points at which the trivial steady state
$\bar u$ of the non-local cell-cell adhesion model~\eqref{Eqn:StSt},
bifurcates to non-homogeneous solutions. To apply the abstract bifurcation
results outlined in \cref{Chapter:Preliminaries}, we cast
equation~\eqref{Eqn:StSt} as an abstract operator equation
$\F[\alpha, u] = 0$ (see equation~\eqref{Eqn:F}), whose
linearization is shown to be a Fredholm operator with index zero. To ensure that
the eigenspaces of the linearization $\D_u \F(\alpha, u)$ at possible
bifurcation points are one-dimensional, we imposed that solutions be reflection
symmetric through the domain's centre. Since the domain is periodic, this does
not reduce the size of $\F^{-1}(0)$. (Indeed, it's equivalent to the well known
decomposition of any periodic function $u \in H^1$ into $u = \mean{u} + v$ where
$v \in H^1_0$). Interestingly, the existence of bifurcation points depends on
the integration kernel of the non-local operator $\K[u]$, since bifurcation
points only exist when the quantity $M_n(\omega)$ is positive. Using the
previously established properties of the non-local operator $\K$ we are able to
apply the abstract local bifurcation theorem of Crandall \etal\
(\cref{Thm:CrandallMain}).

%% file: SampleFunctions.pgf
\begingroup%
\makeatletter%
\begin{pgfpicture}%
\pgfpathrectangle{\pgfpointorigin}{\pgfqpoint{3.148194in}{1.945691in}}%
\pgfusepath{use as bounding box, clip}%
\begin{pgfscope}%
\pgfsetbuttcap%
\pgfsetmiterjoin%
\definecolor{currentfill}{rgb}{1.000000,1.000000,1.000000}%
\pgfsetfillcolor{currentfill}%
\pgfsetlinewidth{0.000000pt}%
\definecolor{currentstroke}{rgb}{1.000000,1.000000,1.000000}%
\pgfsetstrokecolor{currentstroke}%
\pgfsetdash{}{0pt}%
\pgfpathmoveto{\pgfqpoint{0.000000in}{0.000000in}}%
\pgfpathlineto{\pgfqpoint{3.148194in}{0.000000in}}%
\pgfpathlineto{\pgfqpoint{3.148194in}{1.945691in}}%
\pgfpathlineto{\pgfqpoint{0.000000in}{1.945691in}}%
\pgfpathclose%
\pgfusepath{fill}%
\end{pgfscope}%
\begin{pgfscope}%
\pgfsetbuttcap%
\pgfsetmiterjoin%
\definecolor{currentfill}{rgb}{1.000000,1.000000,1.000000}%
\pgfsetfillcolor{currentfill}%
\pgfsetlinewidth{0.000000pt}%
\definecolor{currentstroke}{rgb}{0.000000,0.000000,0.000000}%
\pgfsetstrokecolor{currentstroke}%
\pgfsetstrokeopacity{0.000000}%
\pgfsetdash{}{0pt}%
\pgfpathmoveto{\pgfqpoint{0.618473in}{0.580556in}}%
\pgfpathlineto{\pgfqpoint{1.660417in}{0.580556in}}%
\pgfpathlineto{\pgfqpoint{1.660417in}{1.573469in}}%
\pgfpathlineto{\pgfqpoint{0.618473in}{1.573469in}}%
\pgfpathclose%
\pgfusepath{fill}%
\end{pgfscope}%
\begin{pgfscope}%
\pgfsetbuttcap%
\pgfsetroundjoin%
\definecolor{currentfill}{rgb}{0.000000,0.000000,0.000000}%
\pgfsetfillcolor{currentfill}%
\pgfsetlinewidth{0.803000pt}%
\definecolor{currentstroke}{rgb}{0.000000,0.000000,0.000000}%
\pgfsetstrokecolor{currentstroke}%
\pgfsetdash{}{0pt}%
\pgfsys@defobject{currentmarker}{\pgfqpoint{0.000000in}{-0.048611in}}{\pgfqpoint{0.000000in}{0.000000in}}{%
\pgfpathmoveto{\pgfqpoint{0.000000in}{0.000000in}}%
\pgfpathlineto{\pgfqpoint{0.000000in}{-0.048611in}}%
\pgfusepath{stroke,fill}%
}%
\begin{pgfscope}%
\pgfsys@transformshift{0.618473in}{0.580556in}%
\pgfsys@useobject{currentmarker}{}%
\end{pgfscope}%
\end{pgfscope}%
\begin{pgfscope}%
\pgftext[x=0.618473in,y=0.483333in,,top]{\rmfamily\fontsize{10.000000}{12.000000}\selectfont \(\displaystyle 0\)}%
\end{pgfscope}%
\begin{pgfscope}%
\pgfsetbuttcap%
\pgfsetroundjoin%
\definecolor{currentfill}{rgb}{0.000000,0.000000,0.000000}%
\pgfsetfillcolor{currentfill}%
\pgfsetlinewidth{0.803000pt}%
\definecolor{currentstroke}{rgb}{0.000000,0.000000,0.000000}%
\pgfsetstrokecolor{currentstroke}%
\pgfsetdash{}{0pt}%
\pgfsys@defobject{currentmarker}{\pgfqpoint{0.000000in}{-0.048611in}}{\pgfqpoint{0.000000in}{0.000000in}}{%
\pgfpathmoveto{\pgfqpoint{0.000000in}{0.000000in}}%
\pgfpathlineto{\pgfqpoint{0.000000in}{-0.048611in}}%
\pgfusepath{stroke,fill}%
}%
\begin{pgfscope}%
\pgfsys@transformshift{1.139445in}{0.580556in}%
\pgfsys@useobject{currentmarker}{}%
\end{pgfscope}%
\end{pgfscope}%
\begin{pgfscope}%
\pgftext[x=1.139445in,y=0.483333in,,top]{\rmfamily\fontsize{10.000000}{12.000000}\selectfont \(\displaystyle L/2\)}%
\end{pgfscope}%
\begin{pgfscope}%
\pgfsetbuttcap%
\pgfsetroundjoin%
\definecolor{currentfill}{rgb}{0.000000,0.000000,0.000000}%
\pgfsetfillcolor{currentfill}%
\pgfsetlinewidth{0.803000pt}%
\definecolor{currentstroke}{rgb}{0.000000,0.000000,0.000000}%
\pgfsetstrokecolor{currentstroke}%
\pgfsetdash{}{0pt}%
\pgfsys@defobject{currentmarker}{\pgfqpoint{0.000000in}{-0.048611in}}{\pgfqpoint{0.000000in}{0.000000in}}{%
\pgfpathmoveto{\pgfqpoint{0.000000in}{0.000000in}}%
\pgfpathlineto{\pgfqpoint{0.000000in}{-0.048611in}}%
\pgfusepath{stroke,fill}%
}%
\begin{pgfscope}%
\pgfsys@transformshift{1.660417in}{0.580556in}%
\pgfsys@useobject{currentmarker}{}%
\end{pgfscope}%
\end{pgfscope}%
\begin{pgfscope}%
\pgftext[x=1.660417in,y=0.483333in,,top]{\rmfamily\fontsize{10.000000}{12.000000}\selectfont \(\displaystyle L\)}%
\end{pgfscope}%
\begin{pgfscope}%
\pgftext[x=1.139445in,y=0.288889in,,top]{\rmfamily\fontsize{10.000000}{12.000000}\selectfont Space [\(\displaystyle x\)]}%
\end{pgfscope}%
\begin{pgfscope}%
\pgfsetbuttcap%
\pgfsetroundjoin%
\definecolor{currentfill}{rgb}{0.000000,0.000000,0.000000}%
\pgfsetfillcolor{currentfill}%
\pgfsetlinewidth{0.803000pt}%
\definecolor{currentstroke}{rgb}{0.000000,0.000000,0.000000}%
\pgfsetstrokecolor{currentstroke}%
\pgfsetdash{}{0pt}%
\pgfsys@defobject{currentmarker}{\pgfqpoint{-0.048611in}{0.000000in}}{\pgfqpoint{0.000000in}{0.000000in}}{%
\pgfpathmoveto{\pgfqpoint{0.000000in}{0.000000in}}%
\pgfpathlineto{\pgfqpoint{-0.048611in}{0.000000in}}%
\pgfusepath{stroke,fill}%
}%
\begin{pgfscope}%
\pgfsys@transformshift{0.618473in}{0.580556in}%
\pgfsys@useobject{currentmarker}{}%
\end{pgfscope}%
\end{pgfscope}%
\begin{pgfscope}%
\pgftext[x=0.451806in,y=0.532728in,left,base]{\rmfamily\fontsize{10.000000}{12.000000}\selectfont \(\displaystyle 0\)}%
\end{pgfscope}%
\begin{pgfscope}%
\pgfsetbuttcap%
\pgfsetroundjoin%
\definecolor{currentfill}{rgb}{0.000000,0.000000,0.000000}%
\pgfsetfillcolor{currentfill}%
\pgfsetlinewidth{0.803000pt}%
\definecolor{currentstroke}{rgb}{0.000000,0.000000,0.000000}%
\pgfsetstrokecolor{currentstroke}%
\pgfsetdash{}{0pt}%
\pgfsys@defobject{currentmarker}{\pgfqpoint{-0.048611in}{0.000000in}}{\pgfqpoint{0.000000in}{0.000000in}}{%
\pgfpathmoveto{\pgfqpoint{0.000000in}{0.000000in}}%
\pgfpathlineto{\pgfqpoint{-0.048611in}{0.000000in}}%
\pgfusepath{stroke,fill}%
}%
\begin{pgfscope}%
\pgfsys@transformshift{0.618473in}{1.031880in}%
\pgfsys@useobject{currentmarker}{}%
\end{pgfscope}%
\end{pgfscope}%
\begin{pgfscope}%
\pgftext[x=0.451806in,y=0.984052in,left,base]{\rmfamily\fontsize{10.000000}{12.000000}\selectfont \(\displaystyle 1\)}%
\end{pgfscope}%
\begin{pgfscope}%
\pgfsetbuttcap%
\pgfsetroundjoin%
\definecolor{currentfill}{rgb}{0.000000,0.000000,0.000000}%
\pgfsetfillcolor{currentfill}%
\pgfsetlinewidth{0.803000pt}%
\definecolor{currentstroke}{rgb}{0.000000,0.000000,0.000000}%
\pgfsetstrokecolor{currentstroke}%
\pgfsetdash{}{0pt}%
\pgfsys@defobject{currentmarker}{\pgfqpoint{-0.048611in}{0.000000in}}{\pgfqpoint{0.000000in}{0.000000in}}{%
\pgfpathmoveto{\pgfqpoint{0.000000in}{0.000000in}}%
\pgfpathlineto{\pgfqpoint{-0.048611in}{0.000000in}}%
\pgfusepath{stroke,fill}%
}%
\begin{pgfscope}%
\pgfsys@transformshift{0.618473in}{1.483204in}%
\pgfsys@useobject{currentmarker}{}%
\end{pgfscope}%
\end{pgfscope}%
\begin{pgfscope}%
\pgftext[x=0.451806in,y=1.435376in,left,base]{\rmfamily\fontsize{10.000000}{12.000000}\selectfont \(\displaystyle 2\)}%
\end{pgfscope}%
\begin{pgfscope}%
\pgftext[x=0.396251in,y=1.077012in,,bottom,rotate=90.000000]{\rmfamily\fontsize{10.000000}{12.000000}\selectfont \(\displaystyle u(x)\)}%
\end{pgfscope}%
\begin{pgfscope}%
\pgfpathrectangle{\pgfqpoint{0.618473in}{0.580556in}}{\pgfqpoint{1.041944in}{0.992913in}}%
\pgfusepath{clip}%
\pgfsetbuttcap%
\pgfsetroundjoin%
\pgfsetlinewidth{1.003750pt}%
\definecolor{currentstroke}{rgb}{0.000000,0.000000,0.000000}%
\pgfsetstrokecolor{currentstroke}%
\pgfsetdash{{3.700000pt}{1.600000pt}}{0.000000pt}%
\pgfpathmoveto{\pgfqpoint{0.618473in}{1.528336in}}%
\pgfpathlineto{\pgfqpoint{0.634754in}{1.526163in}}%
\pgfpathlineto{\pgfqpoint{0.651034in}{1.519664in}}%
\pgfpathlineto{\pgfqpoint{0.667314in}{1.508903in}}%
\pgfpathlineto{\pgfqpoint{0.683595in}{1.493981in}}%
\pgfpathlineto{\pgfqpoint{0.699875in}{1.475045in}}%
\pgfpathlineto{\pgfqpoint{0.721582in}{1.443870in}}%
\pgfpathlineto{\pgfqpoint{0.743290in}{1.406418in}}%
\pgfpathlineto{\pgfqpoint{0.770423in}{1.351761in}}%
\pgfpathlineto{\pgfqpoint{0.797557in}{1.289765in}}%
\pgfpathlineto{\pgfqpoint{0.835545in}{1.193824in}}%
\pgfpathlineto{\pgfqpoint{0.911520in}{0.988963in}}%
\pgfpathlineto{\pgfqpoint{0.954934in}{0.877397in}}%
\pgfpathlineto{\pgfqpoint{0.987495in}{0.802263in}}%
\pgfpathlineto{\pgfqpoint{1.014629in}{0.747607in}}%
\pgfpathlineto{\pgfqpoint{1.036336in}{0.710155in}}%
\pgfpathlineto{\pgfqpoint{1.058043in}{0.678980in}}%
\pgfpathlineto{\pgfqpoint{1.074324in}{0.660043in}}%
\pgfpathlineto{\pgfqpoint{1.090604in}{0.645122in}}%
\pgfpathlineto{\pgfqpoint{1.106884in}{0.634360in}}%
\pgfpathlineto{\pgfqpoint{1.123165in}{0.627861in}}%
\pgfpathlineto{\pgfqpoint{1.139445in}{0.625688in}}%
\pgfpathlineto{\pgfqpoint{1.155726in}{0.627861in}}%
\pgfpathlineto{\pgfqpoint{1.172006in}{0.634360in}}%
\pgfpathlineto{\pgfqpoint{1.188286in}{0.645122in}}%
\pgfpathlineto{\pgfqpoint{1.204567in}{0.660043in}}%
\pgfpathlineto{\pgfqpoint{1.220847in}{0.678980in}}%
\pgfpathlineto{\pgfqpoint{1.242554in}{0.710155in}}%
\pgfpathlineto{\pgfqpoint{1.264261in}{0.747607in}}%
\pgfpathlineto{\pgfqpoint{1.291395in}{0.802263in}}%
\pgfpathlineto{\pgfqpoint{1.318529in}{0.864259in}}%
\pgfpathlineto{\pgfqpoint{1.356517in}{0.960201in}}%
\pgfpathlineto{\pgfqpoint{1.432492in}{1.165061in}}%
\pgfpathlineto{\pgfqpoint{1.475906in}{1.276628in}}%
\pgfpathlineto{\pgfqpoint{1.508467in}{1.351761in}}%
\pgfpathlineto{\pgfqpoint{1.535601in}{1.406418in}}%
\pgfpathlineto{\pgfqpoint{1.557308in}{1.443870in}}%
\pgfpathlineto{\pgfqpoint{1.579015in}{1.475045in}}%
\pgfpathlineto{\pgfqpoint{1.595296in}{1.493981in}}%
\pgfpathlineto{\pgfqpoint{1.611576in}{1.508903in}}%
\pgfpathlineto{\pgfqpoint{1.627856in}{1.519664in}}%
\pgfpathlineto{\pgfqpoint{1.644137in}{1.526163in}}%
\pgfpathlineto{\pgfqpoint{1.660417in}{1.528336in}}%
\pgfpathlineto{\pgfqpoint{1.660417in}{1.528336in}}%
\pgfusepath{stroke}%
\end{pgfscope}%
\begin{pgfscope}%
\pgfpathrectangle{\pgfqpoint{0.618473in}{0.580556in}}{\pgfqpoint{1.041944in}{0.992913in}}%
\pgfusepath{clip}%
\pgfsetrectcap%
\pgfsetroundjoin%
\pgfsetlinewidth{1.003750pt}%
\definecolor{currentstroke}{rgb}{0.000000,0.000000,0.000000}%
\pgfsetstrokecolor{currentstroke}%
\pgfsetdash{}{0pt}%
\pgfpathmoveto{\pgfqpoint{0.618473in}{1.091779in}}%
\pgfpathlineto{\pgfqpoint{0.678168in}{1.249726in}}%
\pgfpathlineto{\pgfqpoint{0.710729in}{1.327755in}}%
\pgfpathlineto{\pgfqpoint{0.737863in}{1.385534in}}%
\pgfpathlineto{\pgfqpoint{0.759570in}{1.425891in}}%
\pgfpathlineto{\pgfqpoint{0.781277in}{1.460278in}}%
\pgfpathlineto{\pgfqpoint{0.802984in}{1.488107in}}%
\pgfpathlineto{\pgfqpoint{0.819265in}{1.504385in}}%
\pgfpathlineto{\pgfqpoint{0.835545in}{1.516546in}}%
\pgfpathlineto{\pgfqpoint{0.851825in}{1.524475in}}%
\pgfpathlineto{\pgfqpoint{0.868106in}{1.528095in}}%
\pgfpathlineto{\pgfqpoint{0.884386in}{1.528095in}}%
\pgfpathlineto{\pgfqpoint{0.900666in}{1.524475in}}%
\pgfpathlineto{\pgfqpoint{0.916947in}{1.516546in}}%
\pgfpathlineto{\pgfqpoint{0.933227in}{1.504385in}}%
\pgfpathlineto{\pgfqpoint{0.949508in}{1.488107in}}%
\pgfpathlineto{\pgfqpoint{0.965788in}{1.467870in}}%
\pgfpathlineto{\pgfqpoint{0.987495in}{1.435072in}}%
\pgfpathlineto{\pgfqpoint{1.009202in}{1.396147in}}%
\pgfpathlineto{\pgfqpoint{1.036336in}{1.339898in}}%
\pgfpathlineto{\pgfqpoint{1.068897in}{1.263277in}}%
\pgfpathlineto{\pgfqpoint{1.106884in}{1.165061in}}%
\pgfpathlineto{\pgfqpoint{1.220847in}{0.864259in}}%
\pgfpathlineto{\pgfqpoint{1.253408in}{0.790695in}}%
\pgfpathlineto{\pgfqpoint{1.280542in}{0.737689in}}%
\pgfpathlineto{\pgfqpoint{1.302249in}{0.701750in}}%
\pgfpathlineto{\pgfqpoint{1.323956in}{0.672232in}}%
\pgfpathlineto{\pgfqpoint{1.340236in}{0.654615in}}%
\pgfpathlineto{\pgfqpoint{1.356517in}{0.641066in}}%
\pgfpathlineto{\pgfqpoint{1.372797in}{0.631716in}}%
\pgfpathlineto{\pgfqpoint{1.389078in}{0.626654in}}%
\pgfpathlineto{\pgfqpoint{1.405358in}{0.625930in}}%
\pgfpathlineto{\pgfqpoint{1.421638in}{0.629549in}}%
\pgfpathlineto{\pgfqpoint{1.437919in}{0.637478in}}%
\pgfpathlineto{\pgfqpoint{1.454199in}{0.649640in}}%
\pgfpathlineto{\pgfqpoint{1.470479in}{0.665917in}}%
\pgfpathlineto{\pgfqpoint{1.486760in}{0.686154in}}%
\pgfpathlineto{\pgfqpoint{1.508467in}{0.718953in}}%
\pgfpathlineto{\pgfqpoint{1.530174in}{0.757878in}}%
\pgfpathlineto{\pgfqpoint{1.557308in}{0.814126in}}%
\pgfpathlineto{\pgfqpoint{1.589869in}{0.890748in}}%
\pgfpathlineto{\pgfqpoint{1.627856in}{0.988963in}}%
\pgfpathlineto{\pgfqpoint{1.660417in}{1.077012in}}%
\pgfpathlineto{\pgfqpoint{1.660417in}{1.077012in}}%
\pgfusepath{stroke}%
\end{pgfscope}%
\begin{pgfscope}%
\pgfsetrectcap%
\pgfsetmiterjoin%
\pgfsetlinewidth{0.803000pt}%
\definecolor{currentstroke}{rgb}{0.000000,0.000000,0.000000}%
\pgfsetstrokecolor{currentstroke}%
\pgfsetdash{}{0pt}%
\pgfpathmoveto{\pgfqpoint{0.618473in}{0.580556in}}%
\pgfpathlineto{\pgfqpoint{0.618473in}{1.573469in}}%
\pgfusepath{stroke}%
\end{pgfscope}%
\begin{pgfscope}%
\pgfsetrectcap%
\pgfsetmiterjoin%
\pgfsetlinewidth{0.803000pt}%
\definecolor{currentstroke}{rgb}{0.000000,0.000000,0.000000}%
\pgfsetstrokecolor{currentstroke}%
\pgfsetdash{}{0pt}%
\pgfpathmoveto{\pgfqpoint{1.660417in}{0.580556in}}%
\pgfpathlineto{\pgfqpoint{1.660417in}{1.573469in}}%
\pgfusepath{stroke}%
\end{pgfscope}%
\begin{pgfscope}%
\pgfsetrectcap%
\pgfsetmiterjoin%
\pgfsetlinewidth{0.803000pt}%
\definecolor{currentstroke}{rgb}{0.000000,0.000000,0.000000}%
\pgfsetstrokecolor{currentstroke}%
\pgfsetdash{}{0pt}%
\pgfpathmoveto{\pgfqpoint{0.618473in}{0.580556in}}%
\pgfpathlineto{\pgfqpoint{1.660417in}{0.580556in}}%
\pgfusepath{stroke}%
\end{pgfscope}%
\begin{pgfscope}%
\pgfsetrectcap%
\pgfsetmiterjoin%
\pgfsetlinewidth{0.803000pt}%
\definecolor{currentstroke}{rgb}{0.000000,0.000000,0.000000}%
\pgfsetstrokecolor{currentstroke}%
\pgfsetdash{}{0pt}%
\pgfpathmoveto{\pgfqpoint{0.618473in}{1.573469in}}%
\pgfpathlineto{\pgfqpoint{1.660417in}{1.573469in}}%
\pgfusepath{stroke}%
\end{pgfscope}%
\begin{pgfscope}%
\pgftext[x=1.139445in,y=1.656802in,,base]{\rmfamily\fontsize{12.000000}{14.400000}\selectfont Shift}%
\end{pgfscope}%
\begin{pgfscope}%
\pgfsetbuttcap%
\pgfsetmiterjoin%
\definecolor{currentfill}{rgb}{1.000000,1.000000,1.000000}%
\pgfsetfillcolor{currentfill}%
\pgfsetlinewidth{0.000000pt}%
\definecolor{currentstroke}{rgb}{0.000000,0.000000,0.000000}%
\pgfsetstrokecolor{currentstroke}%
\pgfsetstrokeopacity{0.000000}%
\pgfsetdash{}{0pt}%
\pgfpathmoveto{\pgfqpoint{1.907639in}{0.580556in}}%
\pgfpathlineto{\pgfqpoint{2.949583in}{0.580556in}}%
\pgfpathlineto{\pgfqpoint{2.949583in}{1.573469in}}%
\pgfpathlineto{\pgfqpoint{1.907639in}{1.573469in}}%
\pgfpathclose%
\pgfusepath{fill}%
\end{pgfscope}%
\begin{pgfscope}%
\pgfsetbuttcap%
\pgfsetroundjoin%
\definecolor{currentfill}{rgb}{0.000000,0.000000,0.000000}%
\pgfsetfillcolor{currentfill}%
\pgfsetlinewidth{0.803000pt}%
\definecolor{currentstroke}{rgb}{0.000000,0.000000,0.000000}%
\pgfsetstrokecolor{currentstroke}%
\pgfsetdash{}{0pt}%
\pgfsys@defobject{currentmarker}{\pgfqpoint{0.000000in}{-0.048611in}}{\pgfqpoint{0.000000in}{0.000000in}}{%
\pgfpathmoveto{\pgfqpoint{0.000000in}{0.000000in}}%
\pgfpathlineto{\pgfqpoint{0.000000in}{-0.048611in}}%
\pgfusepath{stroke,fill}%
}%
\begin{pgfscope}%
\pgfsys@transformshift{1.907639in}{0.580556in}%
\pgfsys@useobject{currentmarker}{}%
\end{pgfscope}%
\end{pgfscope}%
\begin{pgfscope}%
\pgftext[x=1.907639in,y=0.483333in,,top]{\rmfamily\fontsize{10.000000}{12.000000}\selectfont \(\displaystyle 0\)}%
\end{pgfscope}%
\begin{pgfscope}%
\pgfsetbuttcap%
\pgfsetroundjoin%
\definecolor{currentfill}{rgb}{0.000000,0.000000,0.000000}%
\pgfsetfillcolor{currentfill}%
\pgfsetlinewidth{0.803000pt}%
\definecolor{currentstroke}{rgb}{0.000000,0.000000,0.000000}%
\pgfsetstrokecolor{currentstroke}%
\pgfsetdash{}{0pt}%
\pgfsys@defobject{currentmarker}{\pgfqpoint{0.000000in}{-0.048611in}}{\pgfqpoint{0.000000in}{0.000000in}}{%
\pgfpathmoveto{\pgfqpoint{0.000000in}{0.000000in}}%
\pgfpathlineto{\pgfqpoint{0.000000in}{-0.048611in}}%
\pgfusepath{stroke,fill}%
}%
\begin{pgfscope}%
\pgfsys@transformshift{2.428611in}{0.580556in}%
\pgfsys@useobject{currentmarker}{}%
\end{pgfscope}%
\end{pgfscope}%
\begin{pgfscope}%
\pgftext[x=2.428611in,y=0.483333in,,top]{\rmfamily\fontsize{10.000000}{12.000000}\selectfont \(\displaystyle L/2\)}%
\end{pgfscope}%
\begin{pgfscope}%
\pgfsetbuttcap%
\pgfsetroundjoin%
\definecolor{currentfill}{rgb}{0.000000,0.000000,0.000000}%
\pgfsetfillcolor{currentfill}%
\pgfsetlinewidth{0.803000pt}%
\definecolor{currentstroke}{rgb}{0.000000,0.000000,0.000000}%
\pgfsetstrokecolor{currentstroke}%
\pgfsetdash{}{0pt}%
\pgfsys@defobject{currentmarker}{\pgfqpoint{0.000000in}{-0.048611in}}{\pgfqpoint{0.000000in}{0.000000in}}{%
\pgfpathmoveto{\pgfqpoint{0.000000in}{0.000000in}}%
\pgfpathlineto{\pgfqpoint{0.000000in}{-0.048611in}}%
\pgfusepath{stroke,fill}%
}%
\begin{pgfscope}%
\pgfsys@transformshift{2.949583in}{0.580556in}%
\pgfsys@useobject{currentmarker}{}%
\end{pgfscope}%
\end{pgfscope}%
\begin{pgfscope}%
\pgftext[x=2.949583in,y=0.483333in,,top]{\rmfamily\fontsize{10.000000}{12.000000}\selectfont \(\displaystyle L\)}%
\end{pgfscope}%
\begin{pgfscope}%
\pgftext[x=2.428611in,y=0.288889in,,top]{\rmfamily\fontsize{10.000000}{12.000000}\selectfont Space [\(\displaystyle x\)]}%
\end{pgfscope}%
\begin{pgfscope}%
\pgfsetbuttcap%
\pgfsetroundjoin%
\definecolor{currentfill}{rgb}{0.000000,0.000000,0.000000}%
\pgfsetfillcolor{currentfill}%
\pgfsetlinewidth{0.803000pt}%
\definecolor{currentstroke}{rgb}{0.000000,0.000000,0.000000}%
\pgfsetstrokecolor{currentstroke}%
\pgfsetdash{}{0pt}%
\pgfsys@defobject{currentmarker}{\pgfqpoint{-0.048611in}{0.000000in}}{\pgfqpoint{0.000000in}{0.000000in}}{%
\pgfpathmoveto{\pgfqpoint{0.000000in}{0.000000in}}%
\pgfpathlineto{\pgfqpoint{-0.048611in}{0.000000in}}%
\pgfusepath{stroke,fill}%
}%
\begin{pgfscope}%
\pgfsys@transformshift{1.907639in}{0.580556in}%
\pgfsys@useobject{currentmarker}{}%
\end{pgfscope}%
\end{pgfscope}%
\begin{pgfscope}%
\pgfsetbuttcap%
\pgfsetroundjoin%
\definecolor{currentfill}{rgb}{0.000000,0.000000,0.000000}%
\pgfsetfillcolor{currentfill}%
\pgfsetlinewidth{0.803000pt}%
\definecolor{currentstroke}{rgb}{0.000000,0.000000,0.000000}%
\pgfsetstrokecolor{currentstroke}%
\pgfsetdash{}{0pt}%
\pgfsys@defobject{currentmarker}{\pgfqpoint{-0.048611in}{0.000000in}}{\pgfqpoint{0.000000in}{0.000000in}}{%
\pgfpathmoveto{\pgfqpoint{0.000000in}{0.000000in}}%
\pgfpathlineto{\pgfqpoint{-0.048611in}{0.000000in}}%
\pgfusepath{stroke,fill}%
}%
\begin{pgfscope}%
\pgfsys@transformshift{1.907639in}{1.031880in}%
\pgfsys@useobject{currentmarker}{}%
\end{pgfscope}%
\end{pgfscope}%
\begin{pgfscope}%
\pgfsetbuttcap%
\pgfsetroundjoin%
\definecolor{currentfill}{rgb}{0.000000,0.000000,0.000000}%
\pgfsetfillcolor{currentfill}%
\pgfsetlinewidth{0.803000pt}%
\definecolor{currentstroke}{rgb}{0.000000,0.000000,0.000000}%
\pgfsetstrokecolor{currentstroke}%
\pgfsetdash{}{0pt}%
\pgfsys@defobject{currentmarker}{\pgfqpoint{-0.048611in}{0.000000in}}{\pgfqpoint{0.000000in}{0.000000in}}{%
\pgfpathmoveto{\pgfqpoint{0.000000in}{0.000000in}}%
\pgfpathlineto{\pgfqpoint{-0.048611in}{0.000000in}}%
\pgfusepath{stroke,fill}%
}%
\begin{pgfscope}%
\pgfsys@transformshift{1.907639in}{1.483204in}%
\pgfsys@useobject{currentmarker}{}%
\end{pgfscope}%
\end{pgfscope}%
\begin{pgfscope}%
\pgfpathrectangle{\pgfqpoint{1.907639in}{0.580556in}}{\pgfqpoint{1.041944in}{0.992913in}}%
\pgfusepath{clip}%
\pgfsetbuttcap%
\pgfsetroundjoin%
\pgfsetlinewidth{1.003750pt}%
\definecolor{currentstroke}{rgb}{0.000000,0.000000,0.000000}%
\pgfsetstrokecolor{currentstroke}%
\pgfsetdash{{3.700000pt}{1.600000pt}}{0.000000pt}%
\pgfpathmoveto{\pgfqpoint{1.907639in}{1.528336in}}%
\pgfpathlineto{\pgfqpoint{1.923920in}{1.526163in}}%
\pgfpathlineto{\pgfqpoint{1.940200in}{1.519664in}}%
\pgfpathlineto{\pgfqpoint{1.956480in}{1.508903in}}%
\pgfpathlineto{\pgfqpoint{1.972761in}{1.493981in}}%
\pgfpathlineto{\pgfqpoint{1.989041in}{1.475045in}}%
\pgfpathlineto{\pgfqpoint{2.010748in}{1.443870in}}%
\pgfpathlineto{\pgfqpoint{2.032456in}{1.406418in}}%
\pgfpathlineto{\pgfqpoint{2.059589in}{1.351761in}}%
\pgfpathlineto{\pgfqpoint{2.086723in}{1.289765in}}%
\pgfpathlineto{\pgfqpoint{2.124711in}{1.193824in}}%
\pgfpathlineto{\pgfqpoint{2.200686in}{0.988963in}}%
\pgfpathlineto{\pgfqpoint{2.244100in}{0.877397in}}%
\pgfpathlineto{\pgfqpoint{2.276661in}{0.802263in}}%
\pgfpathlineto{\pgfqpoint{2.303795in}{0.747607in}}%
\pgfpathlineto{\pgfqpoint{2.325502in}{0.710155in}}%
\pgfpathlineto{\pgfqpoint{2.347209in}{0.678980in}}%
\pgfpathlineto{\pgfqpoint{2.363490in}{0.660043in}}%
\pgfpathlineto{\pgfqpoint{2.379770in}{0.645122in}}%
\pgfpathlineto{\pgfqpoint{2.396051in}{0.634360in}}%
\pgfpathlineto{\pgfqpoint{2.412331in}{0.627861in}}%
\pgfpathlineto{\pgfqpoint{2.428611in}{0.625688in}}%
\pgfpathlineto{\pgfqpoint{2.444892in}{0.627861in}}%
\pgfpathlineto{\pgfqpoint{2.461172in}{0.634360in}}%
\pgfpathlineto{\pgfqpoint{2.477452in}{0.645122in}}%
\pgfpathlineto{\pgfqpoint{2.493733in}{0.660043in}}%
\pgfpathlineto{\pgfqpoint{2.510013in}{0.678980in}}%
\pgfpathlineto{\pgfqpoint{2.531720in}{0.710155in}}%
\pgfpathlineto{\pgfqpoint{2.553427in}{0.747607in}}%
\pgfpathlineto{\pgfqpoint{2.580561in}{0.802263in}}%
\pgfpathlineto{\pgfqpoint{2.607695in}{0.864259in}}%
\pgfpathlineto{\pgfqpoint{2.645683in}{0.960201in}}%
\pgfpathlineto{\pgfqpoint{2.721658in}{1.165061in}}%
\pgfpathlineto{\pgfqpoint{2.765072in}{1.276628in}}%
\pgfpathlineto{\pgfqpoint{2.797633in}{1.351761in}}%
\pgfpathlineto{\pgfqpoint{2.824767in}{1.406418in}}%
\pgfpathlineto{\pgfqpoint{2.846474in}{1.443870in}}%
\pgfpathlineto{\pgfqpoint{2.868181in}{1.475045in}}%
\pgfpathlineto{\pgfqpoint{2.884462in}{1.493981in}}%
\pgfpathlineto{\pgfqpoint{2.900742in}{1.508903in}}%
\pgfpathlineto{\pgfqpoint{2.917022in}{1.519664in}}%
\pgfpathlineto{\pgfqpoint{2.933303in}{1.526163in}}%
\pgfpathlineto{\pgfqpoint{2.949583in}{1.528336in}}%
\pgfpathlineto{\pgfqpoint{2.949583in}{1.528336in}}%
\pgfusepath{stroke}%
\end{pgfscope}%
\begin{pgfscope}%
\pgfpathrectangle{\pgfqpoint{1.907639in}{0.580556in}}{\pgfqpoint{1.041944in}{0.992913in}}%
\pgfusepath{clip}%
\pgfsetrectcap%
\pgfsetroundjoin%
\pgfsetlinewidth{1.003750pt}%
\definecolor{currentstroke}{rgb}{0.000000,0.000000,0.000000}%
\pgfsetstrokecolor{currentstroke}%
\pgfsetdash{}{0pt}%
\pgfpathmoveto{\pgfqpoint{1.907639in}{0.626654in}}%
\pgfpathlineto{\pgfqpoint{1.923920in}{0.625930in}}%
\pgfpathlineto{\pgfqpoint{1.940200in}{0.629549in}}%
\pgfpathlineto{\pgfqpoint{1.956480in}{0.637478in}}%
\pgfpathlineto{\pgfqpoint{1.972761in}{0.649640in}}%
\pgfpathlineto{\pgfqpoint{1.989041in}{0.665917in}}%
\pgfpathlineto{\pgfqpoint{2.005322in}{0.686154in}}%
\pgfpathlineto{\pgfqpoint{2.027029in}{0.718953in}}%
\pgfpathlineto{\pgfqpoint{2.048736in}{0.757878in}}%
\pgfpathlineto{\pgfqpoint{2.075870in}{0.814126in}}%
\pgfpathlineto{\pgfqpoint{2.108431in}{0.890748in}}%
\pgfpathlineto{\pgfqpoint{2.146418in}{0.988963in}}%
\pgfpathlineto{\pgfqpoint{2.260381in}{1.289765in}}%
\pgfpathlineto{\pgfqpoint{2.292941in}{1.363329in}}%
\pgfpathlineto{\pgfqpoint{2.320075in}{1.416336in}}%
\pgfpathlineto{\pgfqpoint{2.341783in}{1.452275in}}%
\pgfpathlineto{\pgfqpoint{2.363490in}{1.481793in}}%
\pgfpathlineto{\pgfqpoint{2.379770in}{1.499409in}}%
\pgfpathlineto{\pgfqpoint{2.396051in}{1.512958in}}%
\pgfpathlineto{\pgfqpoint{2.412331in}{1.522308in}}%
\pgfpathlineto{\pgfqpoint{2.428611in}{1.527370in}}%
\pgfpathlineto{\pgfqpoint{2.444892in}{1.528336in}}%
\pgfpathlineto{\pgfqpoint{2.461172in}{1.526163in}}%
\pgfpathlineto{\pgfqpoint{2.477452in}{1.519664in}}%
\pgfpathlineto{\pgfqpoint{2.493733in}{1.508903in}}%
\pgfpathlineto{\pgfqpoint{2.510013in}{1.493981in}}%
\pgfpathlineto{\pgfqpoint{2.526293in}{1.475045in}}%
\pgfpathlineto{\pgfqpoint{2.548001in}{1.443870in}}%
\pgfpathlineto{\pgfqpoint{2.569708in}{1.406418in}}%
\pgfpathlineto{\pgfqpoint{2.596842in}{1.351761in}}%
\pgfpathlineto{\pgfqpoint{2.623976in}{1.289765in}}%
\pgfpathlineto{\pgfqpoint{2.661963in}{1.193824in}}%
\pgfpathlineto{\pgfqpoint{2.737938in}{0.988963in}}%
\pgfpathlineto{\pgfqpoint{2.781353in}{0.877397in}}%
\pgfpathlineto{\pgfqpoint{2.813913in}{0.802263in}}%
\pgfpathlineto{\pgfqpoint{2.841047in}{0.747607in}}%
\pgfpathlineto{\pgfqpoint{2.862755in}{0.710155in}}%
\pgfpathlineto{\pgfqpoint{2.884462in}{0.678980in}}%
\pgfpathlineto{\pgfqpoint{2.900742in}{0.660043in}}%
\pgfpathlineto{\pgfqpoint{2.917022in}{0.645122in}}%
\pgfpathlineto{\pgfqpoint{2.933303in}{0.634360in}}%
\pgfpathlineto{\pgfqpoint{2.949583in}{0.627861in}}%
\pgfpathlineto{\pgfqpoint{2.949583in}{0.627861in}}%
\pgfusepath{stroke}%
\end{pgfscope}%
\begin{pgfscope}%
\pgfsetrectcap%
\pgfsetmiterjoin%
\pgfsetlinewidth{0.803000pt}%
\definecolor{currentstroke}{rgb}{0.000000,0.000000,0.000000}%
\pgfsetstrokecolor{currentstroke}%
\pgfsetdash{}{0pt}%
\pgfpathmoveto{\pgfqpoint{1.907639in}{0.580556in}}%
\pgfpathlineto{\pgfqpoint{1.907639in}{1.573469in}}%
\pgfusepath{stroke}%
\end{pgfscope}%
\begin{pgfscope}%
\pgfsetrectcap%
\pgfsetmiterjoin%
\pgfsetlinewidth{0.803000pt}%
\definecolor{currentstroke}{rgb}{0.000000,0.000000,0.000000}%
\pgfsetstrokecolor{currentstroke}%
\pgfsetdash{}{0pt}%
\pgfpathmoveto{\pgfqpoint{2.949583in}{0.580556in}}%
\pgfpathlineto{\pgfqpoint{2.949583in}{1.573469in}}%
\pgfusepath{stroke}%
\end{pgfscope}%
\begin{pgfscope}%
\pgfsetrectcap%
\pgfsetmiterjoin%
\pgfsetlinewidth{0.803000pt}%
\definecolor{currentstroke}{rgb}{0.000000,0.000000,0.000000}%
\pgfsetstrokecolor{currentstroke}%
\pgfsetdash{}{0pt}%
\pgfpathmoveto{\pgfqpoint{1.907639in}{0.580556in}}%
\pgfpathlineto{\pgfqpoint{2.949583in}{0.580556in}}%
\pgfusepath{stroke}%
\end{pgfscope}%
\begin{pgfscope}%
\pgfsetrectcap%
\pgfsetmiterjoin%
\pgfsetlinewidth{0.803000pt}%
\definecolor{currentstroke}{rgb}{0.000000,0.000000,0.000000}%
\pgfsetstrokecolor{currentstroke}%
\pgfsetdash{}{0pt}%
\pgfpathmoveto{\pgfqpoint{1.907639in}{1.573469in}}%
\pgfpathlineto{\pgfqpoint{2.949583in}{1.573469in}}%
\pgfusepath{stroke}%
\end{pgfscope}%
\begin{pgfscope}%
\pgftext[x=2.428611in,y=1.656802in,,base]{\rmfamily\fontsize{12.000000}{14.400000}\selectfont Reflection}%
\end{pgfscope}%
\end{pgfpicture}%
\makeatother%
\endgroup%

%% file: BifurcationValuesExample.tikz
\begin{tikzpicture}

\definecolor{color0}{rgb}{0.,0.,0.}

\begin{groupplot}[group style={group size=2 by 1, horizontal sep = 50pt}]
\nextgroupplot[
xlabel={Bifurcation value index $n$},
ylabel={$M_n(\omega)$},
xmin=0, xmax=220,
ymin=0, ymax=0.4,
width=5.75cm,
height=6.5cm,
xtick={0,50,100,150,200},
xticklabels={$0$,$L$,$2L$,$3L$,$4L$},
ytick={0,0.1,0.2,0.3,0.4},
yticklabels={$0.0$,$0.1$,$0.2$,$0.3$,$0.4$},
tick align=outside,
tick pos=left,
x grid style={lightgray!92.026143790849673!black},
y grid style={lightgray!92.026143790849673!black}
]
\addplot [very thick, color0, forget plot]
table {%
1 0.031374606588922
2 0.0625018150337495
3 0.0931368282680688
4 0.123040020615416
5 0.151979446958872
6 0.179733261161931
7 0.206092015278692
8 0.230860812593439
9 0.253861289369858
10 0.274933402344305
11 0.293937001437545
12 0.310753169853119
13 0.325285316642873
14 0.337460009912984
15 0.347227542076781
16 0.354562221891628
17 0.35946239140278
18 0.361950169313551
19 0.362070925664378
20 0.35989249599002
21 0.355504146292053
22 0.349015303172298
23 0.340554066283794
24 0.330265522833577
25 0.318309886183791
26 0.304860482615609
27 0.290101612019528
28 0.274226309635377
29 0.257434036970107
30 0.239928330660014
31 0.221914438310425
32 0.203596970238873
33 0.185177595571129
34 0.166852810301943
35 0.148811803747192
36 0.131234448299494
37 0.114289435577225
38 0.0981325799536164
39 0.0829053080977691
40 0.0687333505860763
41 0.0557256488860664
42 0.043973488113036
43 0.0335498629523452
44 0.0245090810675361
45 0.0168866052176524
46 0.0106991322274275
47 0.00594490393200439
48 0.00260424229307291
49 0.000640298093651473
50 2.38694182861061e-33
51 0.000615188364488665
52 0.00240391596283652
53 0.00527189593970202
54 0.00911407560114188
55 0.0138163133598974
56 0.0192571351244926
57 0.0253095457359797
58 0.0318428707025433
59 0.0387246034631986
60 0.0458222337240508
61 0.0530050330461147
62 0.060145774810281
63 0.067122366926307
64 0.0738193771684652
65 0.0801294327869495
66 0.0859544780343341
67 0.091206875430556
68 0.0958103389359401
69 0.0997006896756983
70 0.10282642742572
71 0.105149113692016
72 0.106643564858202
73 0.107297856500374
74 0.10711314254062
75 0.106103295394597
76 0.104294375631656
77 0.101723941876978
78 0.0984402137152636
79 0.0945011021789001
80 0.0899731239975051
81 0.0849302171311504
82 0.0794524761907796
83 0.0736248271547862
84 0.0675356613126911
85 0.0612754486017849
86 0.0549353504509509
87 0.0486058519121535
88 0.042375432252698
89 0.0363292923125055
90 0.0305481558160339
91 0.0251071604871288
92 0.0200748532689948
93 0.0155123022252779
94 0.0114723358188467
95 0.00799891826099327
96 0.00512666752564234
97 0.00288052046189904
98 0.00127554724558673
99 0.000316915218069917
100 4.77388365722122e-33
101 0.000310639669197241
102 0.00122552578497548
103 0.00271272315343889
104 0.00473230848520829
105 0.00723711652185103
106 0.0101735808204867
107 0.0134826552051481
108 0.0171008009328473
109 0.020961023892924
110 0.0249939456676641
111 0.0291288920343513
112 0.0332949824842627
113 0.037422204569534
114 0.0414424573577348
115 0.0452905489665367
116 0.0489051340540177
117 0.0522295782380107
118 0.0552127376918976
119 0.0578096435934721
120 0.0599820826650034
121 0.0616990667118439
122 0.0629371858179554
123 0.0636808416628233
124 0.0639223592581117
125 0.0636619772367581
126 0.062907718634967
127 0.0616751458624194
128 0.0599870052327387
129 0.0578727680010318
130 0.055368076306157
131 0.0525141037223144
132 0.0493568412700297
133 0.0459463207056184
134 0.0423357876885526
135 0.0385808380085313
136 0.0347385304322188
137 0.0308664899004185
138 0.0270220147698364
139 0.0232612015526115
140 0.0196381001674504
141 0.0162039120874377
142 0.0130062429630106
143 0.0100884203283276
144 0.00748888588174714
145 0.0052406705847887
146 0.00337095946891551
147 0.00190075159730751
148 0.000844619122077691
149 0.000210567829455852
150 7.16082548583184e-33
151 0.00020777885158226
152 0.000822392303075637
153 0.00182621231898175
154 0.0031958446913095
155 0.00490256280512489
156 0.00691281773699733
157 0.00918881596783973
158 0.0116891550680222
159 0.0143695069454636
160 0.017183337646519
161 0.0200826522721304
162 0.0230187533224532
163 0.0259430007138487
164 0.02880756182184
165 0.0315661401887982
166 0.0341746719895546
167 0.0365919799631572
168 0.0387803752835948
169 0.0407061987433324
170 0.0423402936458847
171 0.0436584039306029
172 0.0446414922662242
173 0.0452759741302154
174 0.0455538652184244
175 0.0454728408833987
176 0.0450362076591241
177 0.0442527882741654
178 0.0431367228639919
179 0.0417071903471123
180 0.0399880551100023
181 0.0380074452354872
182 0.0357972694925491
183 0.0333926811685643
184 0.0308314975557938
185 0.028153584492712
186 0.025400215799902
187 0.0226134177345312
188 0.0198353087140288
189 0.0171074445281111
190 0.0144701790707529
191 0.0119620502844436
192 0.00961920052472666
193 0.00747483993238784
194 0.00555876065449273
195 0.00389690889638134
196 0.00251102082888605
197 0.00141832733403152
198 0.000631331464987366
199 0.000157661339642825
200 9.54776731444245e-33
201 0.000156092570094135
202 0.00061882985181929
203 0.00137640632908476
204 0.00241254942383168
205 0.00370681577948467
206 0.00523494935423099
207 0.00696929520266105
208 0.00887926202282458
209 0.0109318258580321
210 0.0130920667783002
211 0.015323729932763
212 0.0175898020671577
213 0.0198530944429922
214 0.0220768230784195
215 0.024225177354194
216 0.0262638682882687
217 0.0281606481744113
218 0.0298857937965317
219 0.0314125460622063
220 0.0327174996354564
221 0.0337809369779779
222 0.0345871021161737
223 0.0351244104238891
224 0.035385591732169
};
\nextgroupplot[
xlabel={Bifurcation value index $n$},
ylabel={Bifurcation value $\alpha_n$},
xmin=0, xmax=220,
ymin=2.00263397387048, ymax=100000,
ymode=log,
width=5.75cm,
height=6.5cm,
xtick={0,50,100,150,200},
xticklabels={$0$,$L$,$2L$,$3L$,$4L$},
ytick={0.1,1,10,100,1000,10000,100000,1000000},
yticklabels={,,$10^{1}$,$10^{2}$,$10^{3}$,$10^{4}$,$10^{5}$,},
tick align=outside,
tick pos=left,
x grid style={lightgray!92.026143790849673!black},
y grid style={lightgray!92.026143790849673!black}
]
\addplot [very thick, color0, forget plot]
table {%
1 2.00263397387048
2 2.01056091052294
3 2.02385632751907
4 2.04264767699246
5 2.0671167822054
6 2.09750335577077
7 2.134109712634
8 2.17730683231881
9 2.22754197400412
10 2.28534810743404
11 2.35135549593816
12 2.42630585946373
13 2.51106966143854
14 2.60666721141852
15 2.71429446650442
16 2.83535466295675
17 2.97149723522445
18 3.12466591033028
19 3.29715843981
20 3.49170120365823
21 3.71154296868242
22 3.96057351931388
23 4.24347486559456
24 4.56591551181373
25 4.93480220054468
26 5.35860917705917
27 5.84781318907076
28 6.41547446103735
29 7.07802185184107
30 7.856327708231
31 8.77720016802606
32 9.87548732153766
33 11.1970951181987
34 12.8033983999141
35 14.7778254287463
36 17.2359219693796
37 20.3411500976885
38 24.3304559796224
39 29.5571215646441
40 36.5655697189446
41 46.2283710901348
42 60.0120195657872
43 80.5299767073522
44 112.799069354783
45 167.436459358637
46 270.140155282252
47 496.744291943291
48 1158.08308446121
49 4808.324171263
50 1.31615802946419e+33
51 5208.85096603714
52 1359.13917551351
53 631.668047111227
54 372.272539131877
55 250.121130646854
56 182.715847880476
57 141.504539925467
58 114.444690373757
59 95.7293038457819
60 82.2725318676255
61 72.3090396726107
62 64.7688869706847
63 58.9729910428968
64 54.4740249896441
65 50.9684183154719
66 48.2453316868734
67 46.1558861208393
68 44.5939974363185
69 43.4841311133944
70 42.7733397448133
71 42.4260501250069
72 42.4206882729817
73 42.7475946290235
74 43.4079051088403
75 44.4132198049021
76 45.7859861045766
77 47.5606095994521
78 49.7853910981521
79 52.5254867744829
80 55.8672192585316
81 59.9242562980426
82 64.8464616699406
83 70.8326797697621
84 78.1494628977455
85 87.1590112021975
86 98.3617892635272
87 112.463232351647
88 130.481337331161
89 153.926337878729
90 185.113196702157
91 227.731791194174
92 287.948828574162
93 376.692140909621
94 514.820545877515
95 746.229158376148
96 1176.56506194765
97 2115.82935395854
98 4827.35674616554
99 19627.8155779047
100 1.31615802946419e+33
101 20428.8691674532
102 5229.46892827017
103 2385.67686429443
104 1380.8298296469
105 911.598500952583
106 654.654102928898
107 498.64126734585
108 396.814170190104
109 326.733656705469
110 276.527120999519
111 239.430174094662
112 211.358199313285
113 189.726914242063
114 172.837995304057
115 159.540196975649
116 149.033329471664
117 140.750261775043
118 134.283481899502
119 129.338118188779
120 125.701243331696
121 123.221543320815
122 121.795818722041
123 121.360172479357
124 121.884577968139
125 123.370055013617
126 125.848046294366
127 129.381864096738
128 134.070323430658
129 140.053937736608
130 147.52437585456
131 156.73832682986
132 168.037588955539
133 181.878251189915
134 198.873548156451
135 219.857851786858
136 245.984269093862
137 278.877316423314
138 320.878949914079
139 375.459004438192
140 447.928229057071
141 546.73780228612
142 685.987733857447
143 890.620602319443
144 1208.16193292396
145 1738.44521383474
146 2721.31736767319
147 4859.27904237331
148 11009.8315460237
149 44460.4768538986
150 1.31615802946419e+33
151 45662.057238222
152 11612.9998191807
153 5264.05030787707
154 3027.71451922208
155 1986.49922769939
156 1417.91226850105
157 1073.54429197379
158 849.285746281355
159 695.240600553063
160 585.049115503115
161 503.714758762093
162 442.194242887266
163 394.772839258939
164 357.698578154656
165 328.429630447036
166 305.199348017329
167 286.754624197836
168 272.193119299882
169 260.859119665965
170 252.275411964307
171 246.098022556106
172 242.085964866491
173 240.081208416597
174 239.995934089709
175 241.805307826689
176 245.544789746427
177 251.31157668093
178 259.270271458556
179 269.663374738216
180 282.827797496316
181 299.219411763477
182 319.448869178333
183 344.333809378656
184 374.975654176029
185 412.874347182728
186 460.103361460382
187 519.583402312689
188 595.523293728377
189 694.155121243033
190 825.010666783687
191 1003.24640436598
192 1254.12873541563
193 1622.3153609368
194 2192.82323049413
195 3144.0846257344
196 4904.39707245887
197 8727.09335718738
198 19705.5074840354
199 79306.3079992447
200 1.31615802946419e+33
201 80908.4151783437
202 20509.7318482449
203 9266.78837785922
204 5312.9266078574
205 3474.82331088729
206 2472.49034459691
207 1866.21361380928
208 1471.85941864751
209 1201.25013538857
210 1007.83851537841
211 865.162793674901
212 757.277017692624
213 674.110766093521
214 609.055773541441
215 557.636718729631
216 516.743387323868
217 484.168973389219
218 458.322909637458
219 438.047135544951
220 422.495845642645
221 411.055487830879
222 403.291126706333
223 398.910702440741
224 397.741973473552
};
\end{groupplot}

\end{tikzpicture}

%% file: local_single_peak.pgf
\begingroup%
\makeatletter%
\begin{pgfpicture}%
\pgfpathrectangle{\pgfpointorigin}{\pgfqpoint{3.922305in}{2.424118in}}%
\pgfusepath{use as bounding box, clip}%
\begin{pgfscope}%
\pgfsetbuttcap%
\pgfsetmiterjoin%
\definecolor{currentfill}{rgb}{1.000000,1.000000,1.000000}%
\pgfsetfillcolor{currentfill}%
\pgfsetlinewidth{0.000000pt}%
\definecolor{currentstroke}{rgb}{1.000000,1.000000,1.000000}%
\pgfsetstrokecolor{currentstroke}%
\pgfsetdash{}{0pt}%
\pgfpathmoveto{\pgfqpoint{0.000000in}{0.000000in}}%
\pgfpathlineto{\pgfqpoint{3.922305in}{0.000000in}}%
\pgfpathlineto{\pgfqpoint{3.922305in}{2.424118in}}%
\pgfpathlineto{\pgfqpoint{0.000000in}{2.424118in}}%
\pgfpathclose%
\pgfusepath{fill}%
\end{pgfscope}%
\begin{pgfscope}%
\pgfsetbuttcap%
\pgfsetmiterjoin%
\definecolor{currentfill}{rgb}{1.000000,1.000000,1.000000}%
\pgfsetfillcolor{currentfill}%
\pgfsetlinewidth{0.000000pt}%
\definecolor{currentstroke}{rgb}{0.000000,0.000000,0.000000}%
\pgfsetstrokecolor{currentstroke}%
\pgfsetstrokeopacity{0.000000}%
\pgfsetdash{}{0pt}%
\pgfpathmoveto{\pgfqpoint{0.687659in}{0.386111in}}%
\pgfpathlineto{\pgfqpoint{3.723694in}{0.386111in}}%
\pgfpathlineto{\pgfqpoint{3.723694in}{2.225507in}}%
\pgfpathlineto{\pgfqpoint{0.687659in}{2.225507in}}%
\pgfpathclose%
\pgfusepath{fill}%
\end{pgfscope}%
\begin{pgfscope}%
\pgfsetbuttcap%
\pgfsetroundjoin%
\definecolor{currentfill}{rgb}{0.000000,0.000000,0.000000}%
\pgfsetfillcolor{currentfill}%
\pgfsetlinewidth{0.803000pt}%
\definecolor{currentstroke}{rgb}{0.000000,0.000000,0.000000}%
\pgfsetstrokecolor{currentstroke}%
\pgfsetdash{}{0pt}%
\pgfsys@defobject{currentmarker}{\pgfqpoint{0.000000in}{-0.048611in}}{\pgfqpoint{0.000000in}{0.000000in}}{%
\pgfpathmoveto{\pgfqpoint{0.000000in}{0.000000in}}%
\pgfpathlineto{\pgfqpoint{0.000000in}{-0.048611in}}%
\pgfusepath{stroke,fill}%
}%
\begin{pgfscope}%
\pgfsys@transformshift{0.687659in}{0.386111in}%
\pgfsys@useobject{currentmarker}{}%
\end{pgfscope}%
\end{pgfscope}%
\begin{pgfscope}%
\pgftext[x=0.687659in,y=0.288889in,,top]{\sffamily\fontsize{10.000000}{12.000000}\selectfont \(\displaystyle 0\)}%
\end{pgfscope}%
\begin{pgfscope}%
\pgfsetbuttcap%
\pgfsetroundjoin%
\definecolor{currentfill}{rgb}{0.000000,0.000000,0.000000}%
\pgfsetfillcolor{currentfill}%
\pgfsetlinewidth{0.803000pt}%
\definecolor{currentstroke}{rgb}{0.000000,0.000000,0.000000}%
\pgfsetstrokecolor{currentstroke}%
\pgfsetdash{}{0pt}%
\pgfsys@defobject{currentmarker}{\pgfqpoint{0.000000in}{-0.048611in}}{\pgfqpoint{0.000000in}{0.000000in}}{%
\pgfpathmoveto{\pgfqpoint{0.000000in}{0.000000in}}%
\pgfpathlineto{\pgfqpoint{0.000000in}{-0.048611in}}%
\pgfusepath{stroke,fill}%
}%
\begin{pgfscope}%
\pgfsys@transformshift{2.205676in}{0.386111in}%
\pgfsys@useobject{currentmarker}{}%
\end{pgfscope}%
\end{pgfscope}%
\begin{pgfscope}%
\pgftext[x=2.205676in,y=0.288889in,,top]{\sffamily\fontsize{10.000000}{12.000000}\selectfont \(\displaystyle L/2\)}%
\end{pgfscope}%
\begin{pgfscope}%
\pgfsetbuttcap%
\pgfsetroundjoin%
\definecolor{currentfill}{rgb}{0.000000,0.000000,0.000000}%
\pgfsetfillcolor{currentfill}%
\pgfsetlinewidth{0.803000pt}%
\definecolor{currentstroke}{rgb}{0.000000,0.000000,0.000000}%
\pgfsetstrokecolor{currentstroke}%
\pgfsetdash{}{0pt}%
\pgfsys@defobject{currentmarker}{\pgfqpoint{0.000000in}{-0.048611in}}{\pgfqpoint{0.000000in}{0.000000in}}{%
\pgfpathmoveto{\pgfqpoint{0.000000in}{0.000000in}}%
\pgfpathlineto{\pgfqpoint{0.000000in}{-0.048611in}}%
\pgfusepath{stroke,fill}%
}%
\begin{pgfscope}%
\pgfsys@transformshift{3.723694in}{0.386111in}%
\pgfsys@useobject{currentmarker}{}%
\end{pgfscope}%
\end{pgfscope}%
\begin{pgfscope}%
\pgftext[x=3.723694in,y=0.288889in,,top]{\sffamily\fontsize{10.000000}{12.000000}\selectfont \(\displaystyle L\)}%
\end{pgfscope}%
\begin{pgfscope}%
\pgfsetbuttcap%
\pgfsetroundjoin%
\definecolor{currentfill}{rgb}{0.000000,0.000000,0.000000}%
\pgfsetfillcolor{currentfill}%
\pgfsetlinewidth{0.803000pt}%
\definecolor{currentstroke}{rgb}{0.000000,0.000000,0.000000}%
\pgfsetstrokecolor{currentstroke}%
\pgfsetdash{}{0pt}%
\pgfsys@defobject{currentmarker}{\pgfqpoint{-0.048611in}{0.000000in}}{\pgfqpoint{0.000000in}{0.000000in}}{%
\pgfpathmoveto{\pgfqpoint{0.000000in}{0.000000in}}%
\pgfpathlineto{\pgfqpoint{-0.048611in}{0.000000in}}%
\pgfusepath{stroke,fill}%
}%
\begin{pgfscope}%
\pgfsys@transformshift{0.687659in}{0.711473in}%
\pgfsys@useobject{currentmarker}{}%
\end{pgfscope}%
\end{pgfscope}%
\begin{pgfscope}%
\pgftext[x=0.343522in,y=0.661331in,left,base]{\sffamily\fontsize{10.000000}{12.000000}\selectfont \(\displaystyle 0.98\)}%
\end{pgfscope}%
\begin{pgfscope}%
\pgfsetbuttcap%
\pgfsetroundjoin%
\definecolor{currentfill}{rgb}{0.000000,0.000000,0.000000}%
\pgfsetfillcolor{currentfill}%
\pgfsetlinewidth{0.803000pt}%
\definecolor{currentstroke}{rgb}{0.000000,0.000000,0.000000}%
\pgfsetstrokecolor{currentstroke}%
\pgfsetdash{}{0pt}%
\pgfsys@defobject{currentmarker}{\pgfqpoint{-0.048611in}{0.000000in}}{\pgfqpoint{0.000000in}{0.000000in}}{%
\pgfpathmoveto{\pgfqpoint{0.000000in}{0.000000in}}%
\pgfpathlineto{\pgfqpoint{-0.048611in}{0.000000in}}%
\pgfusepath{stroke,fill}%
}%
\begin{pgfscope}%
\pgfsys@transformshift{0.687659in}{1.287983in}%
\pgfsys@useobject{currentmarker}{}%
\end{pgfscope}%
\end{pgfscope}%
\begin{pgfscope}%
\pgftext[x=0.343522in,y=1.237841in,left,base]{\sffamily\fontsize{10.000000}{12.000000}\selectfont \(\displaystyle 1.00\)}%
\end{pgfscope}%
\begin{pgfscope}%
\pgfsetbuttcap%
\pgfsetroundjoin%
\definecolor{currentfill}{rgb}{0.000000,0.000000,0.000000}%
\pgfsetfillcolor{currentfill}%
\pgfsetlinewidth{0.803000pt}%
\definecolor{currentstroke}{rgb}{0.000000,0.000000,0.000000}%
\pgfsetstrokecolor{currentstroke}%
\pgfsetdash{}{0pt}%
\pgfsys@defobject{currentmarker}{\pgfqpoint{-0.048611in}{0.000000in}}{\pgfqpoint{0.000000in}{0.000000in}}{%
\pgfpathmoveto{\pgfqpoint{0.000000in}{0.000000in}}%
\pgfpathlineto{\pgfqpoint{-0.048611in}{0.000000in}}%
\pgfusepath{stroke,fill}%
}%
\begin{pgfscope}%
\pgfsys@transformshift{0.687659in}{1.864493in}%
\pgfsys@useobject{currentmarker}{}%
\end{pgfscope}%
\end{pgfscope}%
\begin{pgfscope}%
\pgftext[x=0.343522in,y=1.814351in,left,base]{\sffamily\fontsize{10.000000}{12.000000}\selectfont \(\displaystyle 1.02\)}%
\end{pgfscope}%
\begin{pgfscope}%
\pgftext[x=0.287966in,y=1.305809in,,bottom,rotate=90.000000]{\sffamily\fontsize{10.000000}{12.000000}\selectfont \(\displaystyle u(x)\)}%
\end{pgfscope}%
\begin{pgfscope}%
\pgfpathrectangle{\pgfqpoint{0.687659in}{0.386111in}}{\pgfqpoint{3.036036in}{1.839396in}}%
\pgfusepath{clip}%
\pgfsetrectcap%
\pgfsetroundjoin%
\pgfsetlinewidth{2.007500pt}%
\definecolor{currentstroke}{rgb}{0.000000,0.000000,0.000000}%
\pgfsetstrokecolor{currentstroke}%
\pgfsetdash{}{0pt}%
\pgfpathmoveto{\pgfqpoint{0.687659in}{0.469722in}}%
\pgfpathlineto{\pgfqpoint{0.718505in}{0.471166in}}%
\pgfpathlineto{\pgfqpoint{0.749352in}{0.475732in}}%
\pgfpathlineto{\pgfqpoint{0.780199in}{0.483406in}}%
\pgfpathlineto{\pgfqpoint{0.812232in}{0.494641in}}%
\pgfpathlineto{\pgfqpoint{0.844265in}{0.509164in}}%
\pgfpathlineto{\pgfqpoint{0.877485in}{0.527648in}}%
\pgfpathlineto{\pgfqpoint{0.911891in}{0.550395in}}%
\pgfpathlineto{\pgfqpoint{0.947483in}{0.577691in}}%
\pgfpathlineto{\pgfqpoint{0.984262in}{0.609797in}}%
\pgfpathlineto{\pgfqpoint{1.023414in}{0.648167in}}%
\pgfpathlineto{\pgfqpoint{1.064938in}{0.693366in}}%
\pgfpathlineto{\pgfqpoint{1.108836in}{0.745894in}}%
\pgfpathlineto{\pgfqpoint{1.155106in}{0.806154in}}%
\pgfpathlineto{\pgfqpoint{1.206122in}{0.877865in}}%
\pgfpathlineto{\pgfqpoint{1.261883in}{0.961776in}}%
\pgfpathlineto{\pgfqpoint{1.325950in}{1.064061in}}%
\pgfpathlineto{\pgfqpoint{1.404253in}{1.195293in}}%
\pgfpathlineto{\pgfqpoint{1.684247in}{1.670361in}}%
\pgfpathlineto{\pgfqpoint{1.743568in}{1.761594in}}%
\pgfpathlineto{\pgfqpoint{1.795770in}{1.836289in}}%
\pgfpathlineto{\pgfqpoint{1.842040in}{1.897273in}}%
\pgfpathlineto{\pgfqpoint{1.884751in}{1.948618in}}%
\pgfpathlineto{\pgfqpoint{1.925089in}{1.992307in}}%
\pgfpathlineto{\pgfqpoint{1.961868in}{2.027757in}}%
\pgfpathlineto{\pgfqpoint{1.997461in}{2.057840in}}%
\pgfpathlineto{\pgfqpoint{2.030680in}{2.081989in}}%
\pgfpathlineto{\pgfqpoint{2.062713in}{2.101549in}}%
\pgfpathlineto{\pgfqpoint{2.093560in}{2.116826in}}%
\pgfpathlineto{\pgfqpoint{2.123220in}{2.128148in}}%
\pgfpathlineto{\pgfqpoint{2.152881in}{2.136115in}}%
\pgfpathlineto{\pgfqpoint{2.181355in}{2.140574in}}%
\pgfpathlineto{\pgfqpoint{2.209829in}{2.141887in}}%
\pgfpathlineto{\pgfqpoint{2.238303in}{2.140049in}}%
\pgfpathlineto{\pgfqpoint{2.266777in}{2.135068in}}%
\pgfpathlineto{\pgfqpoint{2.295251in}{2.126966in}}%
\pgfpathlineto{\pgfqpoint{2.324911in}{2.115244in}}%
\pgfpathlineto{\pgfqpoint{2.354571in}{2.100228in}}%
\pgfpathlineto{\pgfqpoint{2.385418in}{2.081193in}}%
\pgfpathlineto{\pgfqpoint{2.417451in}{2.057840in}}%
\pgfpathlineto{\pgfqpoint{2.451857in}{2.028828in}}%
\pgfpathlineto{\pgfqpoint{2.487450in}{1.994724in}}%
\pgfpathlineto{\pgfqpoint{2.525415in}{1.954009in}}%
\pgfpathlineto{\pgfqpoint{2.565753in}{1.906175in}}%
\pgfpathlineto{\pgfqpoint{2.609651in}{1.849223in}}%
\pgfpathlineto{\pgfqpoint{2.657107in}{1.782532in}}%
\pgfpathlineto{\pgfqpoint{2.710496in}{1.702050in}}%
\pgfpathlineto{\pgfqpoint{2.772189in}{1.603267in}}%
\pgfpathlineto{\pgfqpoint{2.850493in}{1.471639in}}%
\pgfpathlineto{\pgfqpoint{3.115063in}{1.021724in}}%
\pgfpathlineto{\pgfqpoint{3.177943in}{0.923627in}}%
\pgfpathlineto{\pgfqpoint{3.233705in}{0.842192in}}%
\pgfpathlineto{\pgfqpoint{3.284721in}{0.773110in}}%
\pgfpathlineto{\pgfqpoint{3.330991in}{0.715501in}}%
\pgfpathlineto{\pgfqpoint{3.374888in}{0.665704in}}%
\pgfpathlineto{\pgfqpoint{3.416413in}{0.623258in}}%
\pgfpathlineto{\pgfqpoint{3.455564in}{0.587621in}}%
\pgfpathlineto{\pgfqpoint{3.492343in}{0.558187in}}%
\pgfpathlineto{\pgfqpoint{3.527936in}{0.533560in}}%
\pgfpathlineto{\pgfqpoint{3.562342in}{0.513462in}}%
\pgfpathlineto{\pgfqpoint{3.595561in}{0.497586in}}%
\pgfpathlineto{\pgfqpoint{3.627595in}{0.485616in}}%
\pgfpathlineto{\pgfqpoint{3.659628in}{0.476967in}}%
\pgfpathlineto{\pgfqpoint{3.690475in}{0.471802in}}%
\pgfpathlineto{\pgfqpoint{3.721321in}{0.469757in}}%
\pgfpathlineto{\pgfqpoint{3.723694in}{0.469729in}}%
\pgfpathlineto{\pgfqpoint{3.723694in}{0.469729in}}%
\pgfusepath{stroke}%
\end{pgfscope}%
\begin{pgfscope}%
\pgfpathrectangle{\pgfqpoint{0.687659in}{0.386111in}}{\pgfqpoint{3.036036in}{1.839396in}}%
\pgfusepath{clip}%
\pgfsetbuttcap%
\pgfsetroundjoin%
\pgfsetlinewidth{0.752812pt}%
\definecolor{currentstroke}{rgb}{0.000000,0.000000,0.000000}%
\pgfsetstrokecolor{currentstroke}%
\pgfsetstrokeopacity{0.600000}%
\pgfsetdash{{2.775000pt}{1.200000pt}}{0.000000pt}%
\pgfpathmoveto{\pgfqpoint{0.688843in}{0.386111in}}%
\pgfpathlineto{\pgfqpoint{0.688843in}{2.225507in}}%
\pgfusepath{stroke}%
\end{pgfscope}%
\begin{pgfscope}%
\pgfpathrectangle{\pgfqpoint{0.687659in}{0.386111in}}{\pgfqpoint{3.036036in}{1.839396in}}%
\pgfusepath{clip}%
\pgfsetbuttcap%
\pgfsetroundjoin%
\pgfsetlinewidth{0.752812pt}%
\definecolor{currentstroke}{rgb}{0.000000,0.000000,0.000000}%
\pgfsetstrokecolor{currentstroke}%
\pgfsetstrokeopacity{0.600000}%
\pgfsetdash{{2.775000pt}{1.200000pt}}{0.000000pt}%
\pgfpathmoveto{\pgfqpoint{2.207463in}{0.386111in}}%
\pgfpathlineto{\pgfqpoint{2.207463in}{2.225507in}}%
\pgfusepath{stroke}%
\end{pgfscope}%
\begin{pgfscope}%
\pgfpathrectangle{\pgfqpoint{0.687659in}{0.386111in}}{\pgfqpoint{3.036036in}{1.839396in}}%
\pgfusepath{clip}%
\pgfsetbuttcap%
\pgfsetroundjoin%
\pgfsetlinewidth{0.752812pt}%
\definecolor{currentstroke}{rgb}{0.000000,0.000000,0.000000}%
\pgfsetstrokecolor{currentstroke}%
\pgfsetstrokeopacity{0.600000}%
\pgfsetdash{{2.775000pt}{1.200000pt}}{0.000000pt}%
\pgfpathmoveto{\pgfqpoint{0.687659in}{0.386111in}}%
\pgfpathlineto{\pgfqpoint{0.687659in}{2.225507in}}%
\pgfusepath{stroke}%
\end{pgfscope}%
\begin{pgfscope}%
\pgfpathrectangle{\pgfqpoint{0.687659in}{0.386111in}}{\pgfqpoint{3.036036in}{1.839396in}}%
\pgfusepath{clip}%
\pgfsetbuttcap%
\pgfsetroundjoin%
\pgfsetlinewidth{1.505625pt}%
\definecolor{currentstroke}{rgb}{0.000000,0.000000,0.000000}%
\pgfsetstrokecolor{currentstroke}%
\pgfsetdash{{5.550000pt}{2.400000pt}}{0.000000pt}%
\pgfpathmoveto{\pgfqpoint{0.687659in}{0.478604in}}%
\pgfpathlineto{\pgfqpoint{0.717319in}{0.480127in}}%
\pgfpathlineto{\pgfqpoint{0.746979in}{0.484690in}}%
\pgfpathlineto{\pgfqpoint{0.777826in}{0.492644in}}%
\pgfpathlineto{\pgfqpoint{0.808673in}{0.503834in}}%
\pgfpathlineto{\pgfqpoint{0.840706in}{0.518834in}}%
\pgfpathlineto{\pgfqpoint{0.873926in}{0.537954in}}%
\pgfpathlineto{\pgfqpoint{0.908332in}{0.561488in}}%
\pgfpathlineto{\pgfqpoint{0.943924in}{0.589704in}}%
\pgfpathlineto{\pgfqpoint{0.980703in}{0.622833in}}%
\pgfpathlineto{\pgfqpoint{1.019855in}{0.662324in}}%
\pgfpathlineto{\pgfqpoint{1.061379in}{0.708685in}}%
\pgfpathlineto{\pgfqpoint{1.106463in}{0.763845in}}%
\pgfpathlineto{\pgfqpoint{1.155106in}{0.828455in}}%
\pgfpathlineto{\pgfqpoint{1.208495in}{0.904693in}}%
\pgfpathlineto{\pgfqpoint{1.269002in}{0.996692in}}%
\pgfpathlineto{\pgfqpoint{1.341373in}{1.112586in}}%
\pgfpathlineto{\pgfqpoint{1.451710in}{1.295929in}}%
\pgfpathlineto{\pgfqpoint{1.583402in}{1.513411in}}%
\pgfpathlineto{\pgfqpoint{1.655773in}{1.626838in}}%
\pgfpathlineto{\pgfqpoint{1.715094in}{1.714265in}}%
\pgfpathlineto{\pgfqpoint{1.767296in}{1.785939in}}%
\pgfpathlineto{\pgfqpoint{1.814753in}{1.846083in}}%
\pgfpathlineto{\pgfqpoint{1.858650in}{1.896944in}}%
\pgfpathlineto{\pgfqpoint{1.900175in}{1.940442in}}%
\pgfpathlineto{\pgfqpoint{1.939326in}{1.977053in}}%
\pgfpathlineto{\pgfqpoint{1.976105in}{2.007334in}}%
\pgfpathlineto{\pgfqpoint{2.011698in}{2.032676in}}%
\pgfpathlineto{\pgfqpoint{2.046104in}{2.053339in}}%
\pgfpathlineto{\pgfqpoint{2.078137in}{2.069094in}}%
\pgfpathlineto{\pgfqpoint{2.110170in}{2.081421in}}%
\pgfpathlineto{\pgfqpoint{2.141017in}{2.090000in}}%
\pgfpathlineto{\pgfqpoint{2.171864in}{2.095314in}}%
\pgfpathlineto{\pgfqpoint{2.201524in}{2.097324in}}%
\pgfpathlineto{\pgfqpoint{2.231184in}{2.096288in}}%
\pgfpathlineto{\pgfqpoint{2.260845in}{2.092210in}}%
\pgfpathlineto{\pgfqpoint{2.291691in}{2.084758in}}%
\pgfpathlineto{\pgfqpoint{2.322538in}{2.074063in}}%
\pgfpathlineto{\pgfqpoint{2.354571in}{2.059570in}}%
\pgfpathlineto{\pgfqpoint{2.386605in}{2.041690in}}%
\pgfpathlineto{\pgfqpoint{2.419824in}{2.019654in}}%
\pgfpathlineto{\pgfqpoint{2.455417in}{1.992214in}}%
\pgfpathlineto{\pgfqpoint{2.492196in}{1.959852in}}%
\pgfpathlineto{\pgfqpoint{2.531347in}{1.921136in}}%
\pgfpathlineto{\pgfqpoint{2.572872in}{1.875542in}}%
\pgfpathlineto{\pgfqpoint{2.616769in}{1.822639in}}%
\pgfpathlineto{\pgfqpoint{2.664226in}{1.760505in}}%
\pgfpathlineto{\pgfqpoint{2.716428in}{1.686925in}}%
\pgfpathlineto{\pgfqpoint{2.775749in}{1.597719in}}%
\pgfpathlineto{\pgfqpoint{2.845747in}{1.486573in}}%
\pgfpathlineto{\pgfqpoint{2.943033in}{1.325713in}}%
\pgfpathlineto{\pgfqpoint{3.105572in}{1.056838in}}%
\pgfpathlineto{\pgfqpoint{3.176757in}{0.945525in}}%
\pgfpathlineto{\pgfqpoint{3.236078in}{0.858329in}}%
\pgfpathlineto{\pgfqpoint{3.288280in}{0.786901in}}%
\pgfpathlineto{\pgfqpoint{3.335736in}{0.727012in}}%
\pgfpathlineto{\pgfqpoint{3.379634in}{0.676413in}}%
\pgfpathlineto{\pgfqpoint{3.421158in}{0.633181in}}%
\pgfpathlineto{\pgfqpoint{3.460310in}{0.596838in}}%
\pgfpathlineto{\pgfqpoint{3.497089in}{0.566821in}}%
\pgfpathlineto{\pgfqpoint{3.532681in}{0.541743in}}%
\pgfpathlineto{\pgfqpoint{3.565901in}{0.521984in}}%
\pgfpathlineto{\pgfqpoint{3.597934in}{0.506354in}}%
\pgfpathlineto{\pgfqpoint{3.629967in}{0.494156in}}%
\pgfpathlineto{\pgfqpoint{3.660814in}{0.485702in}}%
\pgfpathlineto{\pgfqpoint{3.691661in}{0.480514in}}%
\pgfpathlineto{\pgfqpoint{3.721321in}{0.478625in}}%
\pgfpathlineto{\pgfqpoint{3.723694in}{0.478606in}}%
\pgfpathlineto{\pgfqpoint{3.723694in}{0.478606in}}%
\pgfusepath{stroke}%
\end{pgfscope}%
\begin{pgfscope}%
\pgfsetrectcap%
\pgfsetmiterjoin%
\pgfsetlinewidth{0.803000pt}%
\definecolor{currentstroke}{rgb}{0.000000,0.000000,0.000000}%
\pgfsetstrokecolor{currentstroke}%
\pgfsetdash{}{0pt}%
\pgfpathmoveto{\pgfqpoint{0.687659in}{0.386111in}}%
\pgfpathlineto{\pgfqpoint{0.687659in}{2.225507in}}%
\pgfusepath{stroke}%
\end{pgfscope}%
\begin{pgfscope}%
\pgfsetrectcap%
\pgfsetmiterjoin%
\pgfsetlinewidth{0.803000pt}%
\definecolor{currentstroke}{rgb}{0.000000,0.000000,0.000000}%
\pgfsetstrokecolor{currentstroke}%
\pgfsetdash{}{0pt}%
\pgfpathmoveto{\pgfqpoint{3.723694in}{0.386111in}}%
\pgfpathlineto{\pgfqpoint{3.723694in}{2.225507in}}%
\pgfusepath{stroke}%
\end{pgfscope}%
\begin{pgfscope}%
\pgfsetrectcap%
\pgfsetmiterjoin%
\pgfsetlinewidth{0.803000pt}%
\definecolor{currentstroke}{rgb}{0.000000,0.000000,0.000000}%
\pgfsetstrokecolor{currentstroke}%
\pgfsetdash{}{0pt}%
\pgfpathmoveto{\pgfqpoint{0.687659in}{0.386111in}}%
\pgfpathlineto{\pgfqpoint{3.723694in}{0.386111in}}%
\pgfusepath{stroke}%
\end{pgfscope}%
\begin{pgfscope}%
\pgfsetrectcap%
\pgfsetmiterjoin%
\pgfsetlinewidth{0.803000pt}%
\definecolor{currentstroke}{rgb}{0.000000,0.000000,0.000000}%
\pgfsetstrokecolor{currentstroke}%
\pgfsetdash{}{0pt}%
\pgfpathmoveto{\pgfqpoint{0.687659in}{2.225507in}}%
\pgfpathlineto{\pgfqpoint{3.723694in}{2.225507in}}%
\pgfusepath{stroke}%
\end{pgfscope}%
\begin{pgfscope}%
\pgfsetbuttcap%
\pgfsetmiterjoin%
\definecolor{currentfill}{rgb}{1.000000,1.000000,1.000000}%
\pgfsetfillcolor{currentfill}%
\pgfsetfillopacity{0.800000}%
\pgfsetlinewidth{1.003750pt}%
\definecolor{currentstroke}{rgb}{0.800000,0.800000,0.800000}%
\pgfsetstrokecolor{currentstroke}%
\pgfsetstrokeopacity{0.800000}%
\pgfsetdash{}{0pt}%
\pgfpathmoveto{\pgfqpoint{1.396213in}{0.455556in}}%
\pgfpathlineto{\pgfqpoint{3.015139in}{0.455556in}}%
\pgfpathquadraticcurveto{\pgfqpoint{3.042917in}{0.455556in}}{\pgfqpoint{3.042917in}{0.483333in}}%
\pgfpathlineto{\pgfqpoint{3.042917in}{0.862901in}}%
\pgfpathquadraticcurveto{\pgfqpoint{3.042917in}{0.890679in}}{\pgfqpoint{3.015139in}{0.890679in}}%
\pgfpathlineto{\pgfqpoint{1.396213in}{0.890679in}}%
\pgfpathquadraticcurveto{\pgfqpoint{1.368435in}{0.890679in}}{\pgfqpoint{1.368435in}{0.862901in}}%
\pgfpathlineto{\pgfqpoint{1.368435in}{0.483333in}}%
\pgfpathquadraticcurveto{\pgfqpoint{1.368435in}{0.455556in}}{\pgfqpoint{1.396213in}{0.455556in}}%
\pgfpathclose%
\pgfusepath{stroke,fill}%
\end{pgfscope}%
\begin{pgfscope}%
\pgfsetrectcap%
\pgfsetroundjoin%
\pgfsetlinewidth{2.007500pt}%
\definecolor{currentstroke}{rgb}{0.000000,0.000000,0.000000}%
\pgfsetstrokecolor{currentstroke}%
\pgfsetdash{}{0pt}%
\pgfpathmoveto{\pgfqpoint{1.423991in}{0.783450in}}%
\pgfpathlineto{\pgfqpoint{1.701769in}{0.783450in}}%
\pgfusepath{stroke}%
\end{pgfscope}%
\begin{pgfscope}%
\pgftext[x=1.812880in,y=0.734839in,left,base]{\sffamily\fontsize{10.000000}{12.000000}\selectfont Numerical}%
\end{pgfscope}%
\begin{pgfscope}%
\pgfsetbuttcap%
\pgfsetroundjoin%
\pgfsetlinewidth{1.505625pt}%
\definecolor{currentstroke}{rgb}{0.000000,0.000000,0.000000}%
\pgfsetstrokecolor{currentstroke}%
\pgfsetdash{{5.550000pt}{2.400000pt}}{0.000000pt}%
\pgfpathmoveto{\pgfqpoint{1.423991in}{0.586722in}}%
\pgfpathlineto{\pgfqpoint{1.701769in}{0.586722in}}%
\pgfusepath{stroke}%
\end{pgfscope}%
\begin{pgfscope}%
\pgftext[x=1.812880in,y=0.538111in,left,base]{\sffamily\fontsize{10.000000}{12.000000}\selectfont Local approximation}%
\end{pgfscope}%
\end{pgfpicture}%
\makeatother%
\endgroup%

%% file: chapterGlobalBif.tex
\chapter{Global Bifurcation\index{global bifurcation}}\label{Chapter:GlobalBifPeriodic}
For each $\mean{u}>0$ we found local bifurcations at $(\alpha_n, \mean{u})$ with
non-trivial eigenfunctions $e_n$ of $\D_u \F(\alpha_n, \mean{u})$ on $H^2_P$ be given by
\begin{equation}\label{Eqn:GenEfunc}
    \alpha_n = \frac{n \pi}{\mean{u} h'(\mean{u}) L M_n(\omega)} , \qquad
    e_n(x) = \cos\lb\frac{2 \pi n x}{L}\rb,
\end{equation}
where $M_n(\omega)$ are the Fourier-sine coefficients of $\omega$ (see~\ref{Defn:Mn}).
For many examples of PDEs
\parencites{Healey1988}{Healey1991}{Healey1993}{golubitsky2003} the
symmetries of the unstable modes $e_n$ are conserved along the bifurcating
solution branch.  We will show that this is the case here as well. For this we
define the so called {\it isotropy subgroup\/} associated with $e_n$. The isotropy
subgroup\index{isotropy subgroup} contains all the $\O2$-group actions that leave $e_n$ invariant.
\[
    \Sigma_n \coloneqq \lcb \gamma \in \O2 : \gamma e_n = e_n \rcb.
\]
It is easy to see that for the eigenfunction $e_n$ the isotropy subgroup is given by,
\[
    \Sigma_n \coloneqq \lcb \sigma_{\frac{mL}{n}},\ \rho \sigma_{\frac{mL}{n}} :
    m \in \N,\ 0 \leq m \leq n-1 \rcb = \Dn,
\]
where the shift $\sigma_a$ and the reflection $\rho$ were defined in
\eqref{Eqn:ActionO2}. $\Dn$ is the dihedral group\index{dihedral group} of order
$2n$, which is the group of symmetries of regular polygons with $n$ sides.
We can define the group action directly on elements of the circle $S^1_L$. For
given $x\in S_L^1$ we write
\[
  \sigma_{\frac{mL}{n}} x = \lb x-\frac{mL}{n}\rb \mod L, \qquad \rho x = (L-x) \mod L
\]
For convenience, we introduce the square bracket to denote mod $L$:
\[ [x] = x \mod L.\]
We explicitly classify the orbits\index{orbit} of $\Sigma_n$ on $S_L^1$.
\begin{lemma}\label{lemma:orbits}
    Given $x\in S_L^1$ then the orbit of $x$ under $\Sigma_n$ is
    \[ \OO_n (x) = \left\{ \left[\pm \left(x-\frac{mL}{n}\right)\right],\ m=0,\dots, n-1\right\}.\]
    \begin{enumerate}
        \item If $x=\kappa\frac{L}{n}$ for some $\kappa\in \{0, \dots, n-1\}$, then
        \[ \OO_n (x) =\left\{ 0, \frac{L}{n}, \frac{2L}{n}, \dots , \frac{(n-1)L}{n}\right\}\]
        and $|\OO_n(x)| = n$.
        \item If $x=\kappa \frac{L}{2n}$ for an odd $\kappa$, then
        \[ \OO_n (x) =\left\{ \frac{L}{2n}, \frac{3L}{2n}, \dots , \frac{(2n-1)L}{2n}\right\}\]
        and $|\OO_n(x)| = n$.
        \item If $x\neq \kappa\frac{L}{2n}$ for any $\kappa\in\{0, \dots, 2n-1\}$, then
        \begin{eqnarray*}
        \OO_n (x) &=&\left\{ x, \left[x-\frac{L}{n}\right], \left[x-\frac{2L}{n}\right], \dots, \left[x-\frac{(n-1)L}{n}\right],\right. \\
        && \hspace*{0.5cm} \left.[-x], \left[\frac{L}{n}-x\right], \dots, \left[\frac{(n-1)L}{n}-x\right]\right\} \end{eqnarray*}
        and $|\OO_n(x)| = 2n$.
    \end{enumerate}
\end{lemma}
\begin{proof}
We check when two elements of the orbit are identical. Given $m,k\in\{0,\dots,n-1\}$ then
\[ \left[x-\frac{mL}{n}\right] = \left[ x-\frac{kL}{n}\right]\]
if and only if $m=k$. Moreover,
\[
\left[x-\frac{mL}{n}\right] = \left[\frac{kL}{n}-x\right]
\]
if there is a $\theta\in \Z$ such that
\[ x-\frac{mL}{n} = \frac{kL}{n} -x + \theta L,\]
which implies
\[ x =(k+m + \theta n) \frac{L}{2n},\]
i.e.
\[ [x] = \kappa \frac{L}{2n}, \qquad \kappa \in \{ 0, \dots, 2n-1\}.\]
This shows that multiples of $\frac{L}{n}$ and multiples of $\frac{L}{2n}$ form
their own classes of orbits of lengths $n$. All other orbits have length $2n$.
\end{proof}

Using the isotropy subgroup we define the
fixed-point subspaces\index{fixed-point subspace}
(containing all the functions invariant under actions of the isotropy subgroup).
\begin{alignat*}{2}
    H^2_{\Sigma_n} &= \{ u \in H^2 : \sigma u &&= u,\ \forall \sigma \in
        \Sigma_n \}, \\
    L^2_{\Sigma n} &= \{ u \in L^2 : \sigma u &&= u,\ \forall \sigma \in
        \Sigma_n \}.
\end{alignat*}
Note that both of the above spaces are again Banach spaces, since both are closed
subspaces, which follows from the fact that $\Sigma_n$ is a topological group
and hence the action generated by $\sigma \in \Sigma_n$ is continuous.
Using these symmetries we make the observation.

\begin{lemma}\label{l:sameuofx}
Given $u\in H^2_{\Sigma_n}$ and $x\in S_L^1$ then
\[ u(x) = u(\tilde x) \quad \mbox{for all}\ \quad \tilde x\in \OO_n(x).\]
\end{lemma}

For each $n \in \N$ we obtain a $\Sigma_n$ reduced problem of $\F$ such
that
\[
    \F : \R \times H^2_{\Sigma_n} \mapsto L^2_{\Sigma_n} \times \R
\]
since, whenever $u \in H^2_{\Sigma_n}$ we have that
\[
    \F[\alpha, u] = \F[\alpha, \sigma u] = \sigma \F[\alpha, u].
\]
For each $\Sigma_n$ we define the symmetry-preserving steady state problem as
\begin{equation}\label{Eqn:RedPrb}
    \F[\alpha, u] = 0,\ u\in H^2_{\Sigma_n}.
\end{equation}
Then this problem (\ref{Eqn:RedPrb}) has bifurcation points as
multiples of $n$, \ie\ at ${\alpha_{kn}, 1 \leq k}$.
Based on the symmetry we can identify zeroes of the local and non-local derivatives:
\begin{lemma}\label{Lem:FixMaxMin}
    Suppose that $u \in H^2_{\Sigma_n}$, then

    \[
        \K[u]\lb\frac{mL}{2n}\rb = 0,
        \quad u^{\prime}\lb\frac{mL}{2n}\rb = 0,
        \quad \K[u]'' \lb\frac{mL}{2n}\rb = 0,
        \quad 0\leq m\leq 2n-1.
    \]
\end{lemma}
\begin{proof}
 From \cref{Lem:NLSym} we have that for each $a>0$ and $\rho_a u(x) = u(2a-x)$ that
 \[
     \K[\rho_a u] = - \rho_a \K[u].
 \]
 Since $u \in H^2_{\Sigma_n}$, we have that
 $\rho_{\frac{mL}{2n}} u = u$, and so we find that
 \[
     \K[u](x) = \K[\rho_{\frac{mL}{2n}} u] = -\rho_{\frac{mL}{2n}} \K[u](x)
          = - \K[u]\lb\frac{mL}{n} - x\rb.
 \]
 Letting $x = \frac{mL}{2n}$, we obtain
 that
 \[
     \K[u]\lb\frac{mL}{2n}\rb = - \K[u]\lb\frac{mL}{2n}\rb.
 \]
 Then from Lemma~\ref{Prop:Reg} it follows that also
 \[ u'\left(\frac{mL}{2n}\right) =0.\]
 Similarly  for
 the second derivative of the non-local term ${\K[u]}^{\prime\prime}$ we find%
 \[
     {\K[u]}^{\prime\prime}\lb\frac{mL}{2n}\rb =
         -{\K[u]}^{\prime\prime}\lb\frac{mL}{2n}\rb,
 \]
 and thus we must have that ${\K[u]}^{\prime\prime}\lb\frac{mL}{2n}\rb = 0$.
 \end{proof}
The ordering of minima and maxima imposed on $\Sigma_n$ by the dihedral group
motivates the following definition of a \textit{tiling}\index{tiling} of the
domain $S^1_L$. Intuitively, the tiling segregates the domain into pieces on
which the function $u(x)$ is increasing and decreasing.
\begin{definition}[Domain tiling]\label{Defn:Tiling}\index{domain tiling}
    For $n \in \N$ we define, a \textit{tile} by
    \[
        T_i \coloneqq \lb \frac{iL}{2n}, \frac{(i + 1)L}{2n} \rb,\
            i = 0, \dots, 2n - 1.
    \]
    The collection $\mathbf{T}_n \coloneqq \{ T_i \}_{i=0}^{2n-1}$ is called the
    \textit{tiling} of $S^1_L$, and we have that
    \[
        S^1_L = \mathrm{cl} \lb \bigcup_{i = 0}^{2n - 1} T_i \rb.
    \]
    The collection of tile boundaries is denoted by
    \[
        \partial \mathbf{T}^n = S^1_L \setminus \bigcup_{i = 0}^{2n - 1} T_i.
    \]
\end{definition}

A key motivation for this tiling is the fact that each tile contains exactly one point of a canonical orbit.
\begin{lemma}\label{orbits}
Consider $x\in S_L^1$.
\begin{enumerate}
\item  If $x\neq \kappa\frac{L}{2n}$ for any $\kappa\in\{0, 2n-1\}$, then each tile $T_i$ contains exactly one element of the orbit $\OO_n(x)$.
\item If $x=\kappa \frac{L}{2n}$ for some $\kappa\in\{0, 2n-1\}$ then
\[ \OO_n(x) \subset \partial\mathbf{T}^n. \]
\end{enumerate}
\end{lemma}

\begin{proof}
The second case is obvious, hence we consider  $x\neq \kappa\frac{L}{2n}$. Note
 that we have $2n$ tiles and $|\OO_n(x)|=2n$, thus if we can show that $T_i\cap
 \OO_n(x)\neq\emptyset$ for all $i=0,\dots 2n-1$, then we are done.  Assume
 there is an index $j$ such that

\begin{equation}\label{emptytile}
T_j\cap \OO_n(x)=\emptyset.
\end{equation}
If $j$ is even, then, by symmetry, this implies that the orbit $\OO_n(x)$ does
not intersect any even tile. If $j$ is odd, then this means that $\OO_n(x)$ does
not intersect any odd tile. In either case, we know that for $x\in T_0$ we have
$\rho x = L-x\in T_{2n-1}$. While $T_0$ is an even tile, $T_{2n-1}$ is an odd
tile. Consequently, the orbit does not have an intersection with any tile, which
means that it has to contain only boundary points. But this is case 2. Hence the
assumption \eqref{emptytile} is false.
\end{proof}

\begin{figure}[!ht]\centering
    \includegraphics[width=10cm]{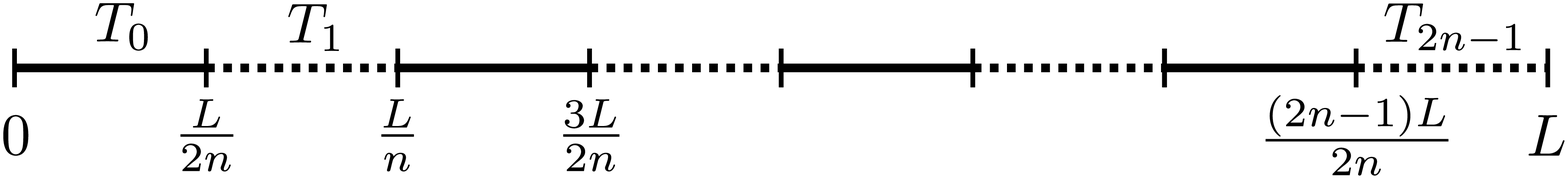}
    \caption[Domain tiling]
    {A tiling for $S^1_L$ as defined in \cref{Defn:Tiling} and where
    $n \in \N$. The solid components form $T^1$ while the dashed components
    form $T^2$ (see \cref{Defn:Spiky}).
    $\partial \mathbf{T}$ is denoted by the vertical lines.}
\end{figure}

Motivated by \cref{Lem:FixMaxMin}, we use the tiling to define function
spaces of functions which have alternating regions on which they are increasing
and decreasing.
\begin{definition}[Spaces of spiky functions\index{space of spiky functions}]\label{Defn:Spiky}
    Let $n \in \N$, and let $\mathbf{T}_n$ be the tiling of $S^1_L$ as defined
    in \cref{Defn:Tiling}. We define two additional collections of tiles
    \[
        T^{1} \coloneqq \bigcup_{i = 0}^{n - 1} T_{2i}, \qquad
        T^{2} \coloneqq \bigcup_{i = 0}^{n - 1} T_{2i + 1}.
    \]
    For each $n \in \N$ we define two spaces of spiky functions to
    be, the space of functions whose derivative has $2n$ simple zeros located on
    $\partial \mathbf{T}_n$
    \[
        \Nodal^{+}_{n} = \lcb u \in \C^2 :
                            u^{\prime} > 0 \ \mbox{in } T^{1},\
                            u^{\prime} < 0 \ \mbox{in } T^{2},\
                            u^{\prime\prime}(x) \neq 0, \  x \in \partial \mathbf{T}_n
                            \rcb,
    \]
    and
    \[
        \Nodal^{-}_{n} = \lcb u \in \C^2 : -u \in \Nodal^{+}_{n} \rcb.
    \]
\end{definition}
\begin{remark}
    It is easy to see that both $\Nodal^{\pm}_{n}$ are nonempty since
    $\cos\lb\frac{2\pi n x}{L}\rb \in \Nodal^{-}_{n}$ while the negative of this
    function is an element of $\Nodal^{+}_{n}$ (see \cref{Fig:CosSpike}).

    \begin{figure}[!ht]\centering
        \input{nodal_test_function.pgf}
        \caption{The solid black line is an example of a function in
        $\Nodal_{2}^{-}$, while the dashed line is a function in
        $\Nodal_{2}^{+}$.}\label{Fig:CosSpike}
    \end{figure}
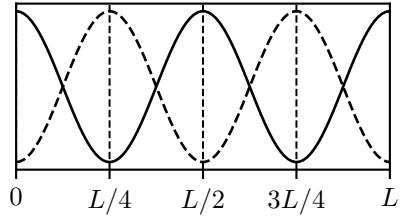
\end{remark}
\begin{lemma}\label{Lem:SolnMass}
    Let $u \in H^2_{\Sigma_n}$ be a positive solution of
    equation~\eqref{Eqn:RedPrb} and $m\in\{0,\dots,n\}$ a natural number. Then we have
     \[ 
        \frac{1}{L}\int_0^{\frac{mL}{n}} u(x) \dd x  
            =  \bar u \frac{m}{n},\quad \mbox{and} \quad          
            \frac{1}{L} \int_{0}^{\frac{L}{2n}} u(x) \dd x = \frac{\bar{u}}{2n}.
     \]
\end{lemma}
\begin{proof}
    Let $u \in H^2_{\Sigma_n}$ be a positive solution of
    equation~\eqref{Eqn:RedPrb}, then
    \[
    \int_0^{\frac{m L}{n}} u(y) \dd y = 
        \sum_{i=0}^{m-1} \int_{i\frac{L}{n}}^{(i+1)\frac{L}{n}} u(x) \dd x =  \bar u \frac{mL}{n},
    \]
    since $u$ is $\frac{L}{n}$-periodic. Furthermore
    \[
        L \bar{u} = \int_{0}^{\frac{L}{2n}} u(x) \dd x + \sum_{i=2}^{2n} \int_{\frac{i-1}{2n}L}^{\frac{i}{2n} L} u(x) \dd x.
    \]
    For the individual integrals we write
    \[
        \int_{\frac{i-1}{2n}L}^{\frac{i}{2n} L} u(y) \dd y = \int_{0}^{\frac{L}{2n}} u\left(x +\frac{i-1}{2n} L \right) \dd x
    \]
    If $i$ is an odd integer, then $i-1$ is even and for $m$ large enough we find
    \begin{eqnarray*}
    u\left(x+\frac{i-1}{2n} L\right) &=& \sigma_{\frac{mL}{n}} u\left(x+\frac{i-1}{2n} L\right)\\
    &=& u\left( x-\frac{2m-i+1}{2n} L\right)\\
    &=& u\left(x-\frac{M}{n}L\right)\\
    &=& \sigma_{\frac{ML}{n}} u(x) \\
    &=& u(x)
    \end{eqnarray*}
    for $M=\frac{1}{2}(2m-i+1)$. Hence
     \[
        \int_{\frac{i-1}{2n}L}^{\frac{i}{2n} L} u(y) \dd y = \int_{0}^{\frac{L}{2n}} u(x) \dd x.
    \]
    If $i$ is even, we have
    \begin{eqnarray*}
    u\left(x+\frac{i-1}{2n}L\right) &=& u\left(x-\frac{1}{2n}L + \frac{i}{2n} L\right) \\
    &=& \sigma_{\frac{i/2}{n}} u\left(x-\frac{L}{2n}\right) \\
    &=& u\left(x-\frac{L}{2n}\right)
    \end{eqnarray*}
    Again
    \[
        \int_{\frac{i-1}{2n}L}^{\frac{i}{2n} L} u(y) \dd y =
        \int_{0}^{\frac{L}{2n}} u\left(x-\frac{L}{2n}\right) \dd x =
        \int_0^{\frac{L}{2n}} u(x) \dd x.
    \]
    We obtain that
    \[
        L \bar{u} = 2n \int_{0}^{\frac{L}{2n}} u(x) \dd x.
    \]
\end{proof}
\section{An Area Function}
In this section we explore a type of non-local convexity\index{non-local
convexity} that is exhibited by an area function of solutions $u(x)$ of
equation~\eqref{Eqn:StSt2}.  The area function will be defined below in
Definition~\ref{Defn:AreaF} and we will see that it is intimately tied to the
non-local operator $\K[u]$. In fact, in \cref{Lem:KnotZero} we use this area
function to prove a non-local maximum principle.

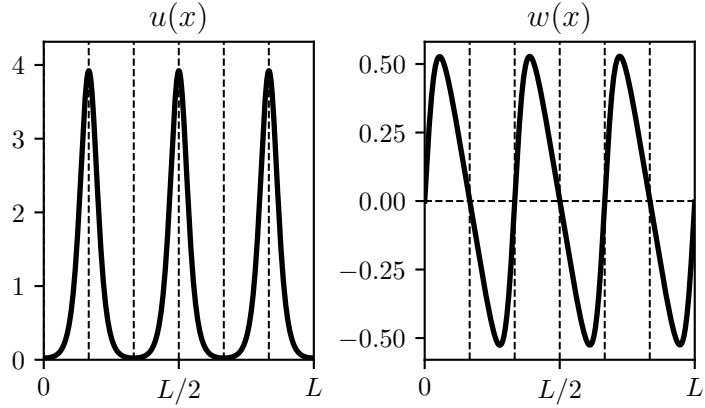
\begin{figure}\centering
    \input{typical_solution_with_w.pgf}
    \caption{Left: A typical three peaked solution of
    equation~\eqref{Eqn:F}. Right: The function $w(x)$
    as defined in \cref{Defn:AreaF} corresponding to the
    solution $u(x)$ on the left.
    }\label{Fig:SolnW}
\end{figure}

\begin{definition}\label{Defn:AreaF}
    Let $u \in H^2_{\Sigma_n}$, be a solution of
    equation~\eqref{Eqn:StSt}, then its
    \textit{modified area function}\index{modified area function} is defined to be
    \[
        w(x) \coloneqq \int_{0}^{x} u(y) \dd y - \bar{u} x.
    \]
    where $\bar u$ is defined in equation~\eqref{Eqn:PopMean}. An
    example of $u(x)$ and the corresponding $w(x)$ is shown in \cref{Fig:SolnW}.
\end{definition}

\begin{lemma}\label{Lem:wProp}
    Let $u(x) \in H^2_{\Sigma_n}$ be a positive solution of
    equation~\eqref{Eqn:StSt}, and let $\mathbf{T}_n$ be a
    tiling of $S^1_L$. Then the function $w(x)$ defined in
    \cref{Defn:AreaF} has the properties:
    \begin{enumerate}
        \item\label{Lem:ContW} $w(x)$ is periodic with period $\frac{L}{n}$  and $w \in \C^2$.

        \item\label{Lem:wSym} $w(x)$ has the following symmetry
            property
            \[
                w(x) = - w\lb \frac{mL}{n} -x \rb,\ m\in\{0,\dots, n-1\}.
            \]
            From this it follows that $w(x) = 0$ for
            $x \in \partial \mathbf{T}_n$.

        \item\label{Lem:wNonLocal} If $h(u)=u$ and $\omega(r) \equiv \ifrac{1}{2}$ then
            \[
                \K[u](x) = \Delta_1 w(x) \coloneqq \frac{1}{2} \lsb w(x + 1) + w(x - 1) - 2 w(x) \rsb.
            \]
    \end{enumerate}
\end{lemma}
\begin{proof}
    \begin{enumerate}
        \item Since $u \in H^2 \subset \C^1$ we have that $w \in \C^2$. Since $u(x) \in H^2_{\Sigma_n}$ we have
            \begin{eqnarray*}
             w\left(x+\frac{L}{n}\right) &=& \int_0^{x+\frac{L}{n}} u(y) \dd y - \bar u \left(x+\frac{L}{n}\right)\\
             &=& \int_0^x u(y) \dd y - \bar u x +\left(\int_0^{\frac{L}{n}} u(y)
                \dd y  -\bar u \frac{L}{n}\right)\\
             &=& w(x),
             \end{eqnarray*}
            where we use \cref{Lem:SolnMass}.
            we observe that $w(0) = w(L) = 0$ and $w'(0) = w'(L)$, so
            that $w(x)$ is periodic.

        \item 
            We note that if $u \in H^2_{\Sigma_n}$ we have from Lemma~\ref{Lem:SolnMass} that
            \[
                \int_{0}^{x} u(y) \dd y = L \bar{u} \frac{m}{n}
                    - \int_{0}^{\frac{mL}{n} - x} u(y) \dd y.
            \]
            Using this identity it's easy to verify that $w(x) +
            w\lb\frac{mL}{n}-x\rb=0$.

        \item 
            In this case $h(u)=u$ is linear and  the linear factors cancel out such that only the integral
            terms of $w(x)$ remain.
    \end{enumerate}
\end{proof}

\begin{remark}\label{Remark:LimitationLinear}
\cref{Lem:wProp}~(3) makes a strong connection between the area
function $w$ and the non-local operator $\K$ in the linear case. Such a
connection is not so easily obtained if $h$ is more general. In that case we
would define \[ w(x) :=\int_0^x h(u(y)) \dd y  - \overline{h(u)} x \]
and the relations to $\K$ are less clear.
Hence, in the following, we restrict to the linear case of $h(u)=u$.

\end{remark}
\subsection{A non-local maximum principle\index{non-local maximum principle}}
We assume $h(u)=u$.\\

\begin{lemma}\label{Lem:PosW}
    Let $n \in \N$, $\mathbf{T}_n$ be a tiling of $S^1_L$, and
    \[ u \in \lcb v \in \C^1 : v' \geq 0\ \mbox{in}\ T^1,\ v' \leq 0\ \mbox{in}\
    T^2 \rcb, \]
    where the $T^i$ are defined in \cref{Defn:Tiling}. Let 
    $T \in \mathbf{T}_n$ be any tile, then
    \begin{enumerate}
        \item If $u'(x) \equiv 0$ in $T$, then $w(x) \equiv 0$ on $T$.
        \item If $u'(\tilde x) \neq 0$ for some $\tilde x \in T$, then
        $u'(\tilde x)  \geq 0$ implies $w(x) < 0$ in $T$ and  $u'(\tilde x)  \leq 0$ implies $w(x) > 0$ in $T$.
    \end{enumerate}
\end{lemma}
\begin{proof}
    Choose any tile $T\subset \mathbf{T}_n$. We assume the tile has been chosen
    such that $u'(x)\leq 0$. The alternative case follows a similar argument.
    Note that on $\partial T$ we have due to \cref{Lem:wProp} that $w = 0$.
    Then there are two cases to consider:
    \begin{enumerate}
        \item If $u'(x) \equiv 0$, we find that $w''(x) = u'(x) \equiv 0$,
            and $w$ satisfies the bounadry value problem
            \begin{equation*} \left\{
            \begin{alignedat}{2}
                w''(x) &= 0\quad &\mbox{in}\ &T \\
                w &= 0 &\mbox{on}\ &\partial T,
            \end{alignedat}
            \right.
            \end{equation*}
            hence, by direct integration, we find that $w \equiv 0$.
        \item If $u'(x) \leq 0$ in $T$ and $u'(\tilde x) <0$ for some $\tilde
         x\in T$, then $w''(\tilde x) = u'(\tilde x) <0$ and  $w$ satisfies the
         boundary value problem
            \begin{equation*} \left\{
            \begin{alignedat}{2}
                w''(x) &\leq 0\quad &\mbox{in}\ &T \\
                w &= 0 &\mbox{on}\ &\partial T.
            \end{alignedat}
            \right.
            \end{equation*}
            Then by the elliptic minimum principle (Theorem~1.4 from~\cite{LopezGomez2013})
            $w$ attains its minimum on the boundaries, and $w(x)>0$ in $T$. Since $w$
            is differentiable, we also have that
            \[
                \frac{\partial w}{\partial \vec{n}} < 0,
            \]
            where $\vec{n}$ is the unit outward normal at $\partial T$.
\end{enumerate}
  If $u'(x) \geq 0$ on $T$ and $u'(\tilde x)>0$ for some $\tilde x\in T$, then
  the elliptic maximum principle shows that $w<0$ on $T$ and
  \[
   \frac{\partial w}{\partial \vec{n}} > 0,
  \]
  where $\vec{n}$ is the unit outward normal at $\partial T$.
\end{proof}

On a given tile $T$, we now combine the concavity of $w$ and the non-locality
of $\K[u]$ to obtain a type of non-local maximum principle for the solutions of
the steady state equations~\eqref{Eqn:StSt}. Key in establishing this result,
are the symmetries of $H^2_{\Sigma_n}$, which allow us to restrict the non-local
operator $\K[u]$ to the canonical tile $T=(0,\frac{L}{2n})$. In more detail, any time
the sensing domain (refer to \cref{Defn:NonLocalOperator}) of $\K[u]$ reaches a
tile boundary, it is reflected. The precise final locations of the sensing
domain endpoints, are important in what follows. This is formalized by
introducing the reflection (subsequent definition), and its properties
(subsequent Lemma).

Let $T_0\in\mathbf{T}_n$ be a given tile and consider $x\in T_0$. The sensing
domain of the non-local operator is
\[
  E(x) = [x-1, x+1] \subset S^1_L,
\]
understood as a closed interval on the circle. Since we are in a symmetric
domain, and since the domain length $L>2$, we can assume without loss of
generality that $0<x-1<L$ and $0<x+1<L$, such that we do not need to worry about
the domain wrapping around the boundary of the parameterization on $[0,L]$.

The end points of integration $x-1$ and $x+1$ might reach into neighboring
tiles, hence we consider the orbits of these points
\[ \OO_{n}(x-1)\qquad \mbox{and}\qquad \OO_n(x+1).\]
We had shown in \cref{orbits} that there are three types of orbits. Those inside
the tiles, which have length $2n$, those on even boundary points
$\theta\frac{L}{n}$ and those residing on odd boundary points
$(2l-1)\frac{L}{2n}$. If $x-1$ and $x+1$ are not in $\partial \mathbf T_n$, then
each orbit has exactly one representative in our starting tile $T_0$, and we can
compute this representative explicitly.

\begin{figure}\centering
 \input{orbits_small.pgf}
 \caption{The orbits of $x$ as defined in \cref{defn:orbits}. The orbits of the
 form $\left[x-\theta\frac{L}{n}\right]$ denoted by open circles ($\circ$), and
 $\left[ \theta\frac{L}{n}-x\right]$ denoted by filled circles ($\bullet$).}
\end{figure}

\begin{definition}\label{defn:orbits}
Given a tile $T_0$. For $x\notin \partial \mathbf{T}_n $ we denote the
representative of $x$ in $T_0$ as $R(x)$.
\end{definition}
The points on the orbits can be written in one of two ways,
\[
 \mbox{either}\qquad \left[x-\theta\frac{L}{n}\right] \qquad \mbox{or} \qquad \left[ \theta\frac{L}{n}-x\right],\quad \theta \in \{0,\dots, n-1\}
\]
where we again use the square brackets to indicate mod $L$. We can compute the
representative as follows.
\begin{lemma}
Consider a given tile $T_0$ and some $x\notin \partial \mathbf{T}_n $ then there
 exists a unique index $l\in\{0,\dots, 2n-1\} $ such that
\begin{eqnarray*}
  \mbox{either} & R(x)& = \left[x-l\frac{L}{2n}\right] \qquad\mbox{and} \qquad l\mbox{ is even,}\\
  \mbox{or} & R(x) & = \left[(l-1)\frac{L}{2n}-x \right] \qquad\mbox{and} \qquad l\mbox{ is odd.}
\end{eqnarray*}
\end{lemma}

\begin{proof}
If $R(x)$ has the form of $[x-\theta\frac{L}{n}]$, then we chose an even $l$ and write
\[ R(x) = \left[x-\theta\frac{L}{n} \right] = \left[ x- l\frac{L}{2n}\right].\]
If $R(x) $ is of the form $[\theta\frac{L}{n}-x]$, then we chose an odd $l$ and write \[ R(x) = \left[\theta\frac{L}{n}-x \right] = \left[ (l-1)\frac{L}{2n}-x\right].\]
As the representative is unique, only one of these cases can arise.
\end{proof}

Then we get explicit equations for the representation of the domain end points
\begin{align}
R(x-1) &= \begin{cases}
[x-1-l\frac{L}{2n}] & \mbox{if $l$ is even}\\
[(l-1)\frac{L}{2n}-(x-1)] & \mbox{if $l$ is odd.}
\end{cases}\label{Rxminusone}\\
R(x+1) &= \begin{cases}
[x+1-k\frac{L}{2n}] & \mbox{if $k$ is even}\\
[(k-1)\frac{L}{2n}-(x+1)] & \mbox{if $k$ is odd.}
\end{cases}\label{Rxplusone}
\end{align}
The corresponding numbers of reflections $l$ and $k$ are denoted as $l_x$ and
$k_x$, where $x$ is the center of the sensing domain.
See \cref{fig:representative} for an illustrative example.

\begin{figure}\centering
 \input{representatives_small.pgf}
 \caption{The two different orbits of $R(x-1)$ denoted by open and closed
 diamonds ($\diamond$) and $R(x+1)$ denoted by open and closed triangles
 ($\triangle$). The location of $x$ is denoted by a filled square ($\square$).
 Key is that each tile $T_i$ contains a representative of each of $R(x-1)$ and
 $R(x+1)$.}\label{fig:representative}
\end{figure}

\begin{lemma}\label{Lem:IntRefl}
 Let $x\in T$, $x\notin \partial \mathbf{T}_n$ and let $R(x-1),\ R(x+1)$
 denote the representatives of the left and right integral boundary points in
 $T$. Let $l_x$ denote the corresponding index for $R(x-1)$ and $r_x$ the index
 for $R(x+1)$, as defined in (\ref{Rxminusone},~\ref{Rxplusone}).
\begin{enumerate} \item Then   $\abs{x - R(x+1)} \leq 1$ and $\abs{x -
    R(x-1)} \leq 1$, where in both cases equality is achieved if and only if $k_x,l_x=0$.
    \item Moreover $|l_x-k_x|\leq 1. $
    \item $\displaystyle w(x-1) = (-1)^{l_x} w(R(x-1)),\ w(x+1) = (-1)^{k_x} w(R(x+1))$
    \item If $h(u) =u $ then
\begin{equation}\label{Kureflected}
    \K[u] = \frac{1}{2} \left( (-1)^{l_x} w(R(x-1)) + (-1)^{k_x} w(R(x+1)) - 2 w(x) \right)
\end{equation}
\end{enumerate}
\end{lemma}

\begin{proof}
\begin{enumerate}
\item Since $x-1$ and $x+1$ are the integration end-points of $\K[u]$ they are
 separated from $x$ by at most one sensing radius, \ie\ ${\abs{R(x+1) - x}}$,
  ${\abs{R(x-1) - x}} \leq 1$. Finally, note that equality is only achieved when
  $k_x,l_x=0$, since if $k_x\neq0$ or $l_x\neq0$ we shift the point $x-1$, or
 $x+1$ closer to $x$ by units of $\frac{L}{n}$, hence reducing the distance.
\item The number of tiles covered to the right or left cannot vary by more than
 one, since the sensing domain is symmetric and the tiles have uniform length.
\item As we have seen earlier $w(x+\frac{L}{2n}) = - w(x)$, hence a shift by
 increments of $\frac{L}{2n}$ leads to a sign change in $w$. Hence
  \[w(x-1) = (-1)^{l_x} w(R(x-1)), \quad w(x+1) = (-1)^{k_x} w(R(x+1)).\]
\item Equation~\eqref{Kureflected} follows directly from the previous item and
    \cref{Lem:wProp}-\ref{Lem:wNonLocal}.
\end{enumerate}
\end{proof}

\begin{proposition}\label{Lem:KnotZero}
    Let $L > 2$, and $\omega(r) \equiv \ifrac{1}{2}$, and $n \in \N$ be such
    that $M_n(\omega) > 0$. Let $\mathbf{T}_n$ be a tiling of $S^1_L$.
    Let $u \in H^2_{\Sigma_n}$ be a positive solution of
    equation~\eqref{Eqn:RedPrb}.
    Further suppose that $h(u) = u$, then $\K[u](x) \neq 0$ for all
    $x \in \interior T_{i}$ for any $T_i \in \mathbf{T}_n$.
\end{proposition}

For the proof of \cref{Lem:KnotZero} we require the following
additional lemmas, which we state and prove now.

\begin{lemma}\label{Lemma:LinearW}
    Let $u(x) \in H^2$ be a positive solution of
    equation~\eqref{Eqn:StSt}, and let
    $w(x)$ be its modified area function defined \cref{Defn:AreaF}. Then
    if $w(x)$ is linear on $[a,b]$, then $\Delta_1 w \equiv 0$ on $[a,b]$.
\end{lemma}

\begin{proof}
    Since $w(x)$ is linear we have that $w'' = 0$. Since $w(x)$ is the modified
    area function of a solution of equation~\eqref{Eqn:StSt}
    we have that $u' = 0$. By \cref{Lem:ZerosDer} we must have that
    $\Delta_1 w = 0$ on $[a,b]$ as well.
\end{proof}

The next lemma recalls a useful property of concave functions.

\begin{lemma}\label{Lem:Concave}
    Let $f \in \C^2[a,b]$ be concave, and let $a \leq x_1 < x_2 \leq x_3 < x_4 \leq b$.
    If the function $f(x)$ satisfies
    \[
        \frac{f(x_2) - f(x_1)}{x_2 - x_1} = \frac{f(x_4) - f(x_3)}{x_4 - x_3} = M \in \R,
    \]
    then $f(x)$ must be linear on $[x_1, x_4]$. If in addition $f(x_1) = f(x_4)$,
    then $f(x)$ must be constant on $[x_1, x_4]$.
\end{lemma}

\begin{proof}
    Let $y \in (x_1, x_2)$ then by the concavity of $f(x)$ we have
    \[
        M = \frac{f(x_2) - f(x_1)}{x_2 - x_1} \geq \frac{f(x_2) - f(y)}{x_2 - y}
            \geq \frac{f(x_4) - f(x_3)}{x_4-x_3} = M.
    \]
    Rearranging we find that
    \[
        f(y) = f(x_2) - M(x_2 - y) = f(x_1) + M(y - x_1),
        \quad\mbox{for}\ y \in (x_1, x_2).
    \]
    We repeat the same argument for any $y \in (x_2, x_3)$,
    and $y \in (x_3, x_4)$, to find that
    \[
        f(y) = f(x_1) + M(y - x_1), \quad\mbox{for}\ y \in (x_1, x_4).
    \]

    If $f(x_1) = f(x_4)$ then we have
    \[
        f(x_4) = f(x_1) + M(x_4 - x_1) = f(x_1),
    \]
    thus $M = 0$, and $f(y) = f(x_1) = f(x_4)$.
\end{proof}

\begin{lemma}\label{Lemma:Helper}
    Let $L > 2$, and $\omega(r) \equiv \ifrac{1}{2}$, and $n \in \N$ be such
    that $M_n(\omega) > 0$, and $h(u) = u$. Let $\mathbf{T}_n$ be a tiling of
    $S^1_L$, and denote the canonical tile by $T$. Let $u \in H^2_{\Sigma_n}$ be
    a positive solution of equation~\eqref{Eqn:RedPrb}, then the set
    \[
        N \coloneqq \{ x \in T : \Delta_1 w(x) = 0 \}
    \]
    is both relatively open and closed in $T$.
\end{lemma}

\begin{proof}
    The canonical tile is given by $T = (0, \ifrac{L}{2n})$. Note that when
    $w \equiv 0$ the result follows immediately. Thus we consider $w \neq 0$.
    Without loss of generality assume that on $T$ we have that $u' \leq 0$, and
    thus $w > 0$ by \cref{Lem:PosW}. \\

    {\noindent\bf $N$ is closed:} We note that $N$ can be written as
    \[
        N = \lb \Delta_1 w \rb^{-1}(0).
    \]
    Since $w \in \C^2$ we have that $\Delta_1 w$ is continuous, implying that
    $N$ must be closed. \\

    {\noindent\bf $N$ is open:} Let $\tilde{x} \in N$, since $\Delta_1 w(\tilde{x}) = 0$ we have
    \[
     \K[u](x) = \frac{1}{2}\lsb (-1)^{k_x} w(R(x+1)(k_x)) + (-1)^{l_x} w(R(x-1)(l_x)) - 2w(\tilde{x}) \rsb.
    \]
    where $R(x-1)(l_x)$ and $R(x+1)(k_x)$ are the endpoints of the sensing domain of $\K$,
    reflected by the appropriate $(l_x,k_x)$-reflection. For a graphical
    representation of the setup of this proof see \cref{Fig:SecantLines}.
    We consider two cases:\\

    {\noindent\bf $l_x, k_x$ have equal parity:} Since $|k_x-l_x|\leq 1$ we must
    have that $l_x=k_x$.

    \begin{itemize}[leftmargin=\parindent, labelindent=*]
        \item If $k_x$ odd, we have that
            \[
                \Delta_1 w(\tilde{x}) = \frac{1}{2}\lsb -w(R(x+1)) - w(R(x-1)) - 2w(\tilde{x}) \rsb = 0.
            \]
            Since either $w > 0$ or $w \equiv 0$ on $T$, we must have that
            $w \equiv 0$ on $T$. Then we have
            $B(\tilde{x}, \ifrac{\varepsilon}{2}) \subset N$ for
            $\varepsilon = \dist(\tilde{x}, \partial T)$.\\

        \item If $k_x$ even, we have that
            \[
                \Delta_1 w(\tilde{x}) = \frac{1}{2}\lsb w(R(x+1)) + w(R(x-1)) - 2w(\tilde{x}) \rsb
            \]
            If $k_x \geq 2$, the separation of $R(x+1)$ and $R(x-1)$ is
            \[
                R(x-1) - R(x+1) = \frac{k_x L - 2n}{n}.
            \]
            The case in which the separation is zero coincides with
            $R(x-1) = R(x+1) = \tilde{x}$, but this is impossible since
            $L \neq \ifrac{2n}{k_x}$ (since $M_n(\omega) > 0$). It's now
            straightforward to show that $\tilde{x} - R(x+1) = -(\tilde{x} - R(x-1))$,
            and therefore we either have $R(x-1) < \tilde{x} < R(x+1)$ or
            $R(x+1) < \tilde{x} < R(x-1)$. If $\Delta_1 w(\tilde{x}) = 0$ we have
            \[
                \R \ni M \coloneqq \frac{w(R(x+1)) - w(\tilde{x})}{R(x+1) - \tilde{x}} =
                    \frac{w(\tilde{x}) - w(R(x-1))}{\tilde{x} - R(x-1)}
            \]
            Applying \cref{Lem:Concave} we find that $w$ must be linear
            on $(R(x-1), R(x+1))$ or $(R(x+1), R(x-1))$, a subsequent application of
            \cref{Lemma:LinearW} implies that $\Delta_1 w = 0$ on the same
            intervals. This means that
            $B(\tilde{x}, \ifrac{\varepsilon}{2}) \subset N$, for
            $\varepsilon = |\tilde{x} - R(x+1)|$.
    \end{itemize}

    {\noindent\bf $l_x, k_x$ have different parity:}\index{parity} Since $|l_x-k_x|\leq1$, we
    have from equation~\eqref{Rxminusone} and~\eqref{Rxplusone} that
    \[
        R(x-1) - R(x+1) = \begin{cases}
         2(-1)^{k_x} \lb L/2n - x\rb, &\mbox{if}\ k_x = l_x+1\\
         2(-1)^{k_x+1} x, &\mbox{if}\ k_x = l_x-1
                   \end{cases}
    \]
    Note that since $x \notin \partial T_n$ we have $R(x-1) - R(x+1) \neq 0$.
    Without loss of generality consider the case in which $l_x$ is odd and $k_x$
    is even (the inverse case follows by the same argument with a sign flip).
    In this case,
    \[
        \Delta_1 w(\tilde{x}) = 0 \quad\iff\quad w(R(x-1)) - w(R(x+1)) = -2w(\tilde{x}).
    \]

\begin{figure}\centering
    \input{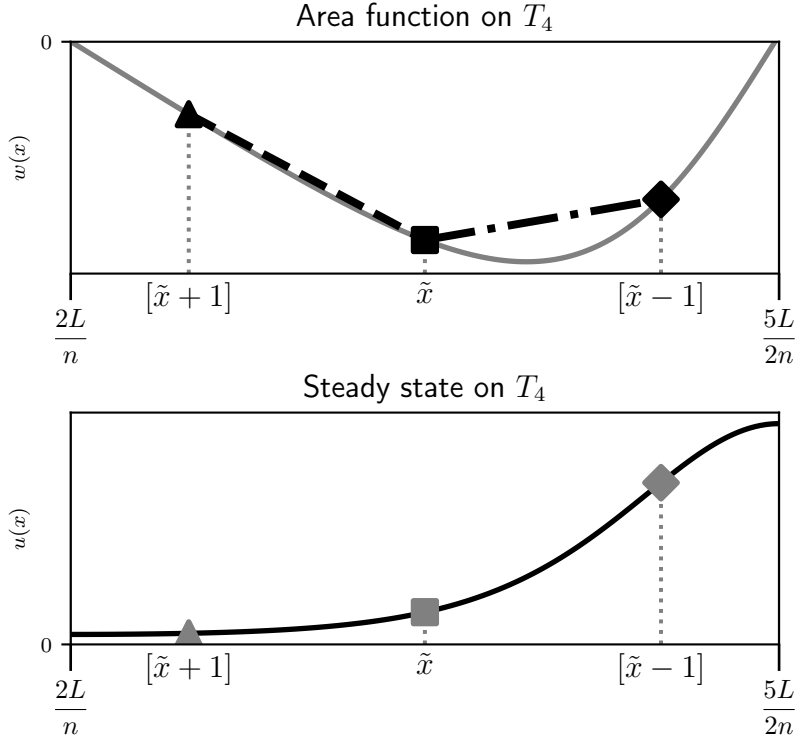}
    \caption[Secant lines of area function]{Top: Shows the area function $w(x)$
    on a selected tile. Bottom: The corresponding spiky function $u(x)$.
    In both the location of $x$ is indicated by a square, the location of the
    representative of $[x-1]$ is indicated by a diamond, and the representative
    of $[x+1]$ is indicated by a triangle. The two secant lines of the area
    function connecting $([x+1], x)$ and $(x, [x-1])$ are shown by a dashed and
    dashed dotted line respectively.
    }\label{Fig:SecantLines}
\end{figure}

    \begin{itemize}[leftmargin=\parindent, labelindent=*]
        \item If $k_x = l_x+1$, then $k_x \geq 2$ and $R(x-1) > R(x+1)$, and since
            $R(x-1) \leq L_n$ and $R(x-1) - R(x+1) = 2(L/2n - \tilde{x})$ we have
            $\tilde{x} - R(x+1) \geq L/2n - \tilde{x}$, and $R(x+1) < \tilde{x}$.
            We must consider two configurations
            (1) $R(x+1) < R(x-1) < \tilde{x} < L/2n$,
            and (2) $R(x+1) < \tilde{x} < R(x-1) < L/2n$.
            For configuration (1), $\Delta_1 w(\tilde{x}) = 0$ is equivalent to
            \[
                \frac{w(R(x-1)) - w(R(x+1))}{2(L/2n-\tilde{x})} =
                \frac{w(R(x-1)) - w(R(x+1))}{R(x-1) - R(x+1)} =
                  \frac{-2w(\tilde{x})}{2(L/2n-\tilde{x})}.
            \]
            By \cref{Lem:Concave} we have that $w$ is linear on
            $(R(x+1), L/2n)$. It follows by \cref{Lemma:LinearW} that $\Delta_1 w = 0$
            on the same interval. Then
            $B(\tilde{x}, \ifrac{\varepsilon}{2}) \subset N$, for
            $\varepsilon = L/2n - \tilde{x}$.

            For configuration (2), the properties of $R(x-1), R(x+1)$, and
            the concavity of $w$ imply that
            \[
                \frac{R(x-1) - \tilde{x}}{\tilde{x} - R(x+1)} \leq 1,\qquad
                \frac{w(\tilde{x}) - w(R(x+1))}{\tilde{x} - R(x+1)} \geq
                    \frac{w(R(x-1)) - w(\tilde{x})}{R(x-1) - \tilde{x}}.
            \]
            For $R(x-1) < L/2n$ we have $w(R(x-1)) > 0$, resulting in
            \[
                w(\tilde{x}) - w(R(x+1)) > -w(R(x-1)) - w(\tilde{x})
                    \quad\iff\quad \Delta_1 w(\tilde{x}) < 0.
            \]
            If $\Delta_1 w(\tilde{x}) = 0$, we must have that $w \equiv 0$ on
            $[R(x+1), R(x-1)]$. Then $B(\tilde{x}, \ifrac{\varepsilon}{2}) \subset N$
            for $\varepsilon = R(x-1) - \tilde{x}$.

            If $R(x-1) = L/2n$, then $\Delta_1 w(\tilde{x}) = 0$
            is equivalent to
            \[
                \frac{w(R(x-1)) - w(\tilde{x})}{L/2n-\tilde{x}}
                    = \frac{w(\tilde{x}) - w(R(x+1))}{L/2n-\tilde{x}}.
            \]
            Once again apply \cref{Lem:Concave}, and
            \cref{Lemma:LinearW} to find that $\Delta_1 w = 0$ on $(R(x+1), L/2n)$.
            Thus $B(\tilde{x}, \ifrac{\varepsilon}{2}) \subset N$, for
            $\varepsilon = L/2n - \tilde{x}$.

        \item If $k_x = l_x - 1$, then $R(x-1) < R(x+1)$, and $2x \leq R(x+1)$. We once again
            have to consider two configurations: (1) $x < R(x-1) < R(x+1)$ and (2)
            $R(x-1) < x < R(x+1)$. For configuration (1), $\Delta_1 w(\tilde{x}) = 0$ is
            equivalent to:
            \[
                \frac{w(R(x+1)) - w(R(x-1))}{R(x+1) - R(x-1)} = \frac{w(R(x+1)) - w(R(x-1))}{2x}
                    = \frac{2w(x) - 2w(0)}{2x - 0}.
            \]
            By \cref{Lem:Concave} we have that $w$ is linear on $(0,
            R(x+1))$, and $\Delta_1 w = 0$ on the same interval by
            \cref{Lem:PosW}. Then
            $B(\tilde{x}, \ifrac{\varepsilon}{2}) \subset N$, for $\varepsilon =
            \tilde{x}$.

            For configuration (2), the properties of $R(x-1), R(x+1)$ and the concavity
            of $w$ imply that
            \[
                \frac{\tilde{x} - R(x-1)}{\tilde{x} - R(x+1)} \geq -1, \qquad
                \frac{w(x) - w(R(x-1))}{x - R(x-1)} \geq \frac{w(R(x+1)) - w(x)}{R(x+1) - x}.
            \]
            For $R(x-1) > 0$, we have that $w(R(x-1)) > 0$, hence it follows that
            \[
                w(\tilde{x}) + w(R(x-1)) > w(\tilde{x}) - w(R(x-1))
                    \geq w(R(x+1)) - w(\tilde{x}).
            \]
            Thus we find that
            \[
                w(R(x+1)) - w(R(x-1)) - 2w(\tilde{x}) < 0 \quad\iff\quad \Delta_1
                w(\tilde{x})< 0.
            \]
            If $\Delta_1 w(\tilde{x}) = 0$, we must have that $w \equiv 0$ on
            $(\tilde{x}, R(x+1))$. Then $B(\tilde{x}, \ifrac{\varepsilon}{2}) \subset N$
            for $\varepsilon = \tilde{x} - R(x-1)$.

            If $R(x-1) = 0$, $\Delta_1 w(\tilde{x}) = 0$ is equivalent to
            \[
                \frac{w(R(x+1)) - w(\tilde{x})}{\tilde{x}}
                    = \frac{w(\tilde{x}) - w(0)}{\tilde{x}}
            \]
            then we can once again use \cref{Lem:Concave}, and
            \cref{Lemma:LinearW} to find that $\Delta_1 w = 0$ on $(0, R(x+1))$.
            Thus $B(\tilde{x}, \ifrac{\varepsilon}{2}) \subset N$, for
            $\varepsilon = \tilde{x}$.
    \end{itemize}

\end{proof}

\begin{proof}[Proof of \cref{Lem:KnotZero}]
    Without loss of generality, let $x \in \lb 0, \ifrac{L}{2n}\rb \eqqcolon T$
    be the canonical tile on which $w'' \leq 0$. Since $u$ is non-constant,
    we have that $w > 0$ on $T$. Suppose that the set
    \[
        N = \{ x \in \interior{T} : u'(x) = 0 \} \neq \emptyset.
    \]
    Note that by \cref{Lem:ZerosDer} we have that
    $u'(x) = 0$ if and only if $\Delta_1 w(x) = 0$, thus the set $N$ here and
    the set $N$ in \cref{Lemma:Helper} coincide. By \cref{Lemma:Helper} the set
    $N$ is both relatively open and closed in $T$. But $T$ is connected, so that
    $N$ must equal $T$. But this means that $u' \equiv 0$ on $T$, a
    contradiction. This means that $N$ must be empty.
\end{proof}

\section{Global bifurcation\index{global bifurcation} branches for linear adhesion function}

\begin{theorem}\label{Thm:GlobalBif}
    Let $\F$ be given as in equation~\eqref{Eqn:F} with $h(u) =
    u$, fix $\N \ni n > 0$, assume the assumptions from
    \cref{Thm:LocBif} hold, and let $\bar u>0$ be given.  Further let $\Gamma_n = (\alpha_n(s),
    u_n(x, s))$, $s\in(-\delta_n, \delta_n)$ denote the local bifurcation branch
    from Theorem~\ref{Thm:LocBif}. Then the
    set of solutions of equation~\eqref{Eqn:RedPrb} contains a
    closed connected set $\Cont \subset \R \times H^2_{\Sigma_n}$ such that
    \begin{enumerate}[label= (\arabic*)]
        \item\label{GlbBif:Res1}
            $\Cont$ contains $(\alpha_n(s), u_n(x, s))$ for
            $s \in (-\delta_n, \delta_n)$, where $\alpha_n$ and $u_n$ were
            defined in \eqref{Eqn:LocBif} and
            \eqref{Eqn:LocBifFun}.
        \item\label{GlbBif:Res2} For any $(\alpha, u) \in \Cont$,
            we have that $\alpha > 0$, $u > 0$.
        \item\label{GlbBif:Res3}
            $\Cont = \Cont^{+} \cup \Cont^{-}$ can be split into the positive and negative direction of the $n$-th
            eigenfunction given in~\eqref{Eqn:GenEfunc}.
            $\Cont^{\pm}$ are closed and connected subsets of $\Cont$ such that
            $\Cont^{+} \cap \Cont^{-} = \{ (\alpha_n, \mean{u}) \}$. We denote
            $\Cont_{n}^{\pm} \coloneqq \Cont^{\pm} \setminus \{ (\alpha_n,
            \mean{u}) \}$, and we have $\Cont_{n}^{+} \subset \R \times
            \Nodal_{n}^{-}$, and that $\Cont_{n}^{-} \subset \R \times
            \Nodal_{n}^{+}$, where $\Nodal_{n}^{\pm}$ are defined in
            \cref{Defn:Spiky}.
        \item\label{GlbBif:Res4} The unilateral branches\index{unilateral branches}
            $\Cont^{\pm}_{n}$ are unbounded, that is for any
            $\alpha \geq \alpha_n$ there exists $(\alpha, u) \in \Cont^{\pm}_{n}$.
    \end{enumerate}
\end{theorem}
\begin{proof}
\begin{enumerate}
    \item From \cref{Thm:RabinowitzMainThm} it follows that there exists the
        component $\Cont$, which is the maximal, connected and closed subset of
        the closure of the set of non-trivial solutions
        \[
            \Soln = \lcb (\alpha, u) \in \R \times H^2_{\Sigma_n} :
                \F[\alpha, u] = 0,\ u \neq \mean{u} \rcb.
        \]
        containing $(\alpha_n, \mean{u})$, this of course is also a consequence
        of \cref{Thm:LocBif}. This proves~\ref{GlbBif:Res1}.

    \item
        Since $(\alpha_n, \mean{u}) \in \Cont$ we have $\alpha_n > 0$.
        Indeed, if $\alpha \leq 0$,
        then by the connectedness of $\Cont$ there would have to be a point in
        $\Cont$ at which $\alpha = 0$. Thus suppose that $(0, u) \in \Cont$, but
        then from equation~\eqref{Eqn:StSt2} we have that
        $u \equiv \mean{u}$. Thus $(0, \mean{u})$ is a bifurcation point. But
        this contradicts \cref{Lem:OpEvals} which showed that all
        bifurcation values are non-zero.

        We show the positivity of $u$ by considering
        \[
            \Pos \coloneqq \lcb u \in H^2 : u > 0 \ \mbox{in } S^1_L
            \rcb.
        \]
        Then we want to show that $\Cont \subset \R \times \Pos$. We first note
        that $\mean{u} \in \Pos$ and by \cref{Thm:LocBif} we have that
        the solution component around the bifurcation point is also in
        $\Pos$. Since, $\Cont$ is connected and $\Pos$ is open, we have that if
        $\Cont \not\subset \R \times \Pos$, then there exists a
        $(\alpha, u) \in \R \times \partial \Pos$ such that $0 \leq u$.
        First, suppose that there exists $\hat{x} \in S^1_L$ such that
        $u(\hat{x}) = 0$. Then note that
        equation~\eqref{Eqn:StSt2} can be written as,
        \[
            \left\{
                    \begin{array}{@{}ll@{}}
                        -u^{\prime\prime} + a(x) u^{\prime}(x) + b(x) u(x) = 0
                            \quad &\mbox{in } \lsb 0, L \rsb  \\
                        \Bd[u, u^\prime] = 0.
                \end{array}
            \right.
        \]
        where
        \[
            a(x) = \alpha \K[u](x) \qquad b(x) = \alpha \K[u^{\prime}](x).
        \]
        Due to \cref{Prop:Reg} we have that both
        $a(x), b(x) \in \C^2(S^1_L)$. Then, the
        maximum principle for non-negative functions
        implies that $u \equiv 0$. However, this contradicts the integral
        constraint in equation~\eqref{Eqn:StSt},
        which must hold on $\Cont$.
        Thus we have that, $\Cont \subset \R \times \Pos$.

    \item\label{BifThmP:3}
        Then consider the decomposition of $\Cont$ into subcontinua such that
        $\Cont^{+} \cap \Cont^{-} = \{(\alpha_n, \mean{u})\}$. Since the cosine is decreasing in the first tile, we claim
        that, $\Cont^{+}_{n} \subset \R \times \Nodal^{-}_{n}$.
        Since, $\Cont^{+}_{n}$ is a connected topological space it suffices
        to show that $\Cont^{+}_{n} \cap (\R \times \Nodal^{-}_{n})$ is
        nonempty, relatively open and relatively closed\index{relatively closed} in $\Cont^{+}_{n}$.
        Then we conclude that this space is all of $\Cont^{+}_{n}$. We split the
        proof of this claim into three pieces. \\

        {\bfseries (1)} To show that $\Cont^{+}_{n} \cap (\R \times \Nodal^{-}_{n})$ is
        non-empty, we note that due to the local bifurcation result
        \cref{Thm:LocBif} for small $s$ the solution along the branch
        is given by equation~\eqref{Eqn:LocBifFun} with $n$,
        thus the solution branch is in
        $\Cont^{+}_{n} \cap (\R \times \Nodal^{-}_{n})$.  \\

        {\bfseries (2)} To show that $\Cont^{+}_{n} \cap (\R \times \Nodal^{-}_{n})$
        is relatively open\index{relatively open} in $\Cont^{+}_{n}$, we use sequential openness. Let
        $(\hat{\alpha}, \hat{u}) \in \Cont^{+}_{n} \cap (\R \times \Nodal^{-}_{n})$,
        and consider any sequence $(\alpha_k, u_k) \subset \Cont^{+}_{n}$
        convergent to $(\hat{\alpha}, \hat{u})$. Then we are left with showing
        that the tail of this sequence is contained in
        $\Cont^{+}_{n} \cap (\R \times \Nodal^{-}_{n})$. Due to the reflections,
        we must only consider the canonical tile $T_0 \coloneqq (0, L/2n)$.

        Since $\hat{u} \in \Nodal^{-}_{n}$ we have that $\hat{u}^{\prime} < 0$
        on $T_0$, and that
        $u^{\prime\prime} \neq 0$ on $\partial T_0$. Further, we must have
        that $\hat{u}^{\prime}$ is decreasing at the left boundary
        ($R(x-1) = 0$), and increasing at the right boundary ($R(x+1) = L / 2n$).
        Hence, we have that

        \begin{equation}\label{Eqn:GlbOpenContra}
            \hat{u}^{\prime\prime}(R(x-1)) < 0 < \hat{u}^{\prime\prime}(R(x+1)).
        \end{equation}

        \begin{figure}[!ht]\centering
            \includegraphics[width=8cm]{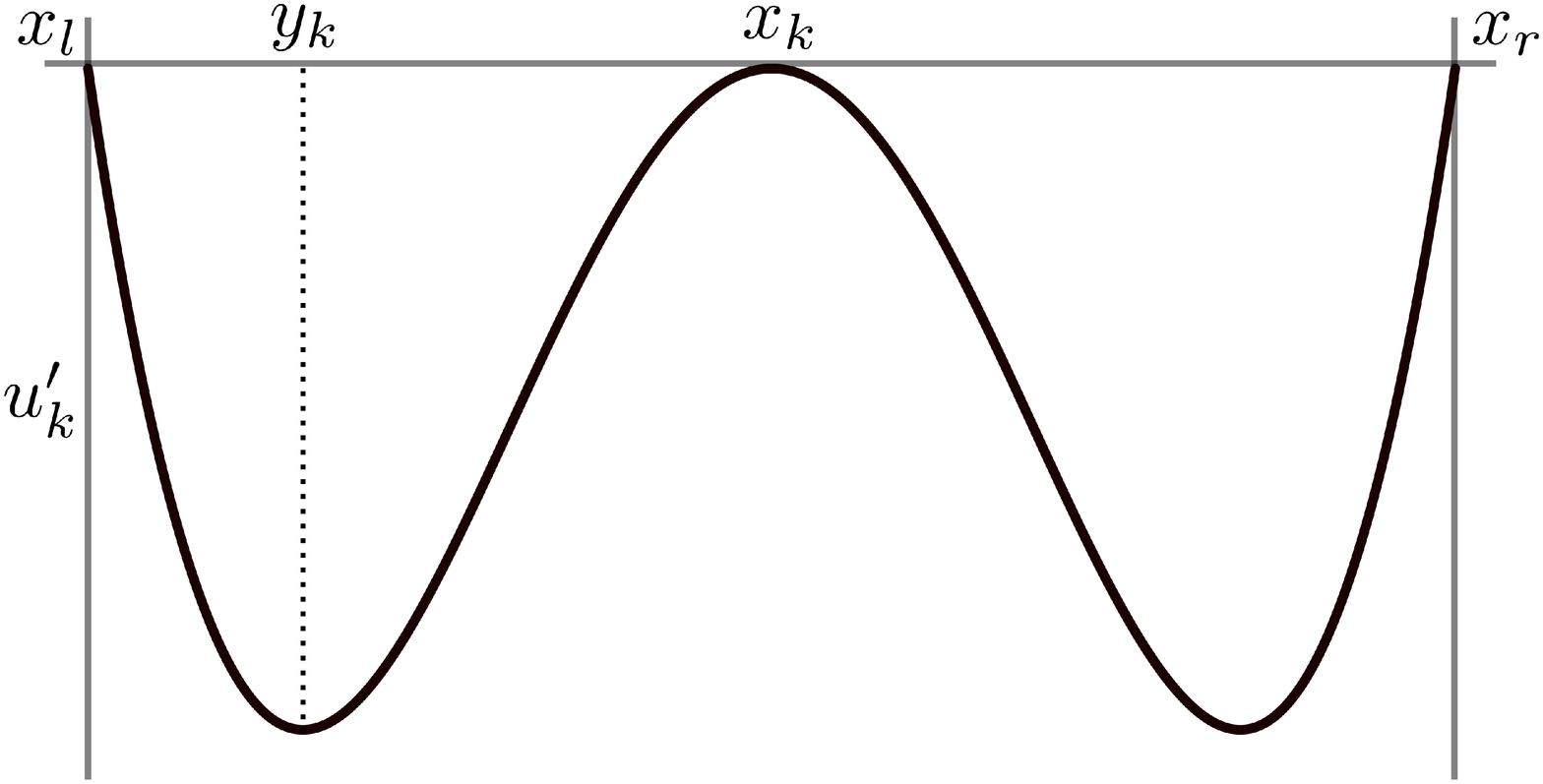}
            \caption[Proof setup of global first bifurcation branch]
            {Plot of the derivative $u^{\prime}$. The point $x_k$
            is such that $u^{\prime}(x_k) = 0$, this implies the existence of
            point $y_k$ at which ${u^{\prime\prime}(y_k) = 0}$.}\label{Fig:OpennessProof}
        \end{figure}

        Since $u_{k}$ converges to $\hat{u}$, we have that eventually
        $u_{k}^{\prime} \leq 0$ on $T_0$. Suppose that there exist a
        $x_{k_n}$ ($k_n$ a subsequence of $k \to \infty$) such that
        $u^{\prime}(x_{k_n}) = 0$. This implies that
        there exists $y_{k_n}$ such that $u_{k}^{\prime\prime}(y_{k_n}) = 0$
        (where either $y_{k_n} \in (R(x-1), x_{k_n})$ or
        $y_{k_n} \in (x_{k_n}, R(x+1))$).

        Without loss of generality suppose that $y_{k_n} \in (R(x-1), x_{k_n})$
        (see \cref{Fig:OpennessProof}). Since, $\hat{u} \in \Nodal_{n}^{-}$
        we must have that $x_{k_n} \to R(x-1)$ as $k \to \infty$, which
        implies that $y_{k_n} \to R(x-1)$. Hence, we have that
        $u^{\prime\prime}(y_{k_n}) = 0$ for $y_{k_n} \to R(x-1)$
        contradicting~\eqref{Eqn:GlbOpenContra}. \\

        {\bfseries (3)} To show that $\Cont^{+}_{n} \cap (\R \times \Nodal^{-}_{n})$
        is relatively closed,\index{relatively closed} we consider the sequence
        $(\alpha_k, u_k) \subset \Cont^{+}_{n} \cap (\R \times \Nodal^{-}_{n})$
        convergent to $(\hat{\alpha}, \hat{u}) \in \Cont^{+}_{n}$. Then again
        consider $T_0$ where $u_{k}^{\prime} < 0, \forall k$. Then we
        must have that, $\hat{u}^{\prime} \leq 0$. Suppose that there exists
        $\tilde{x} \in T_0$ such that $\hat{u}^{\prime}(\tilde{x}) = 0$.
        This means, since $\hat{u}^{\prime} \leq 0$, that
        $\hat{u}^{\prime\prime}(\tilde{x}) = 0$. Evaluating
        equation~\eqref{Eqn:1stInt} at $\tilde{x}$ we
        obtain that $\K[\hat{u}](\tilde{x}) = 0$.

        If $\tilde{x} \in \interior T_0$ and if $\hat{u}$ is non-constant,
        then \cref{Lem:KnotZero} implies that
        $\K[\hat{u}](\tilde{x}) \neq 0$. Thus we have a contradiction.

        If $\tilde{x} \in \partial T_0$ and if $\hat{u}$ is non-constant, then
        $u_k'(\tilde{x}) = 0, \forall k$, but $u_k''(\tilde{x}) \neq 0, \forall k$.
        If $\hat{u}''(\tilde{x}) = 0$, equation~\eqref{Eqn:StSt2}
        implies that\ $\K[\hat{u}]'(\tilde{x}) = 0$. Writing $\K[u]$ in terms of
        $w(x)$ and differentiating we get that
        \[
            \K[u]'(x) = u(x+1) + u(x-1) - 2u(x).
        \]
        Now since $\hat{u}'(\tilde{x}) = 0$ and $\hat{u}$ is non-constant, thus
        $\hat{u}$ must either have a global extrema at $\tilde{x}$. Since $L
        \neq \frac{2n}{k}$, it's impossible for an extrema to be located at
        $\tilde{x}$ and $\tilde{x} \pm 1$, thus $\K[\hat{u}]'(\tilde{x})$ cannot
        be zero.

        In the previous we excluded the possibility that $\hat{u}$
        could be constant. Indeed if $\hat{u}$ is constant then $\hat{u} \equiv
        \mean{u}$, from the integral constraint in
        equation~\eqref{Eqn:StSt}. This means that
        $(\hat{\alpha}, \mean{u})$ is a bifurcation point. Therefore we must
        have that $\hat{\alpha} = \alpha_{jn}$ for some $\N \ni j \geq 1$.  The
        case $j = 1$ cannot occur since $\Cont^{+}_{n}$ does not contain
        $(\alpha_n, \mean{u})$. For $j \geq 2$ we have from the local
        bifurcation result (see \cref{Thm:LocBif}) that in a small
        neighborhood of $(\alpha_{jn}, \mean{u})$ the solution branch must be
        given by~\eqref{Eqn:LocBif}. But this means that
        $u_{k}^{\prime} > 0$ and $u_{k}^{\prime} < 0$ on $T_0$, which is a
        contradiction. Thus, $\Cont^{+}_{n} \cap (\R \times \Nodal_{n}^{-})$ is
        closed.

    Thus, we have proven that $\Cont^{+}_{n} \subset \R \times \Nodal^{-}_{n}$.
    For $\Cont^{-}_{n}$ we proceed analogously.

\item \cref{Thm:UnilateralAbstractBifurcation} implies that each $\Cont_{n}^{\pm}$
        satisfies one of the following alternatives:
        \begin{enumerate}[label= (\roman*)]
            \item\label{Item:Alt1} it is not compact in
                $\R \times H^2_{\Sigma_n}$,
            \item\label{Item:Alt2} it contains a point
                $(\hat{\alpha}, \mean{u})$ where $\hat{\alpha} \neq \alpha_n$,
            \item\label{Item:Alt3} it contains a point
                $(\alpha, \mean{u} + \tilde{u})$ where
                $\tilde{u} \in Y_n \setminus \{ 0 \}$,
        \end{enumerate}
        where
        \[
            Y_n = \lcb u \in H^2_{\Sigma_n} : \int_{0}^{L} u(x)
            \cos\lb\frac{2 \pi n x}{L}\rb \dd x = 0 \rcb.
        \]
        If alternative~\ref{Item:Alt2} holds, then $\hat{\alpha}$ is
        a bifurcation point, which is impossible after the proof
        of~\ref{GlbBif:Res3}.

        If alternative~\ref{Item:Alt3} holds, then there is a
        $(\alpha, \mean{u} + \tilde{u}) \in \Cont_{n}^{\pm}$ with $\tilde{u} \in
        Y_n \setminus \{ 0 \}$. This means the following holds, where we integrate
        by parts
        \begin{equation}\label{Eqn:GlbBifAlt3Contra}
            0 = \int_{0}^{L} \tilde{u}(x) \cos\lb\frac{2 \pi n x}{L}\rb \dd x
              = - \frac{L}{2\pi n} \int_{0}^{L} \tilde{u}^{\prime}(x)
              \sin\lb\frac{2\pi n x}{L}\rb \dd x,
        \end{equation}
        then we note that since $\tilde{u} \in \Nodal_{n}^{\pm}$ we have that
        $\sin\lb\frac{2\pi n x}{L}\rb$ and $\tilde{u}^{\prime}$ have the same
        zeros. Further, since $\tilde{u}^{\prime\prime} \neq 0$ on $\partial T$
        we have that $\tilde{u}^{\prime}$ must change sign at those points.
        Hence we have two cases, either
        $\sin\lb\frac{2\pi n x}{L}\rb \tilde{u}^{\prime} \leq 0$, or
        $\sin\lb\frac{2\pi n x}{L}\rb \tilde{u}^{\prime} \geq 0$.
        In both cases we then find that
        \[
            \int_{0}^{L} \tilde{u}^{\prime}(x) \sin\lb\frac{2\pi n x}{L}\rb \dd x
            \neq 0,
        \]
        contradicting equation~\eqref{Eqn:GlbBifAlt3Contra}.
        Therefore only alternative~\ref{Item:Alt1} holds, and thus
        $\Cont^{\pm}$ are non compact.

        Finally, to complete~\ref{GlbBif:Res4} we note that since
        $\Cont^{\pm}$ are connected its projection onto the $\alpha$ coordinate
        are intervals containing $\alpha_n$. From the a-priori estimate for
        positive solution derived in \cref{Lemma:AprioriEstimate} we find that
        for bounded $\alpha$ the solution $u$ is bounded. In particular,
        uniformly bounded in $\C^1$. From equation,
        \[
            u^{\prime\prime}
                = \alpha u^{\prime} \K[u] + \alpha u {\K[u]}^{\prime}
        \]
        we then also have that $u^{\prime\prime}$ is bounded, since $u$ and
        $u^{\prime}$ are bounded. Iterating this process one more time for
        $u^{\prime\prime\prime}$ we consider equation,
        \[
            u^{\prime\prime\prime} = \alpha u^{\prime\prime} \K[u]
                    + 2 \alpha u^{\prime} {\K[u]}^{\prime}
                    + \alpha u {\K[u]}^{\prime\prime}.
        \]
        The first two terms on the right hand side are bounded, and the
        second-derivative of $\K[u]$ is also bounded since $u \in \C^3$ and
        $h(\cdot) \in \C^2$. Hence $u$ is bounded in the norm of
        $\C^3$ and hence $H^3$ for bounded $\alpha$. But $H^3 \ssubset H^2$, and if $\alpha$ were contained in a bounded interval, then
        $\Cont^{\pm}$ would be
        compact in $\R \times H^2$, which contradicts our earlier observation. Thus the $\alpha$ coordinate
        must be unbounded.
\end{enumerate}
\end{proof}

\section{Bifurcation type for linear adhesion function}\label{sec:BifurcationType}
We study the stability and the type of the first bifurcating branch. Note that
this work is limited to the case $h(u) = u$ as we make heavy use of
\cref{Lem:EfuncK} and \cref{Lem:EfuncKp}.

For each $n \in \N$, we let $e_n \coloneqq \cos\lb\frac{2\pi n x}{L}\rb$, which
is the function spanning the nullspace of $\D_u \F(\alpha_n, \mean{u})$. Now
we decompose the function space $H^2$ into two pieces
\[
    H^2 = \vspan\lcb \cos\lb\frac{2\pi n x}{L}\rb \rcb \oplus Y_n,
\]
where
\[
    Y_n \coloneqq \lcb u \in H^2 : \int_{0}^{L} u(x) e_n(x) \dd x = 0 \rcb.
\]

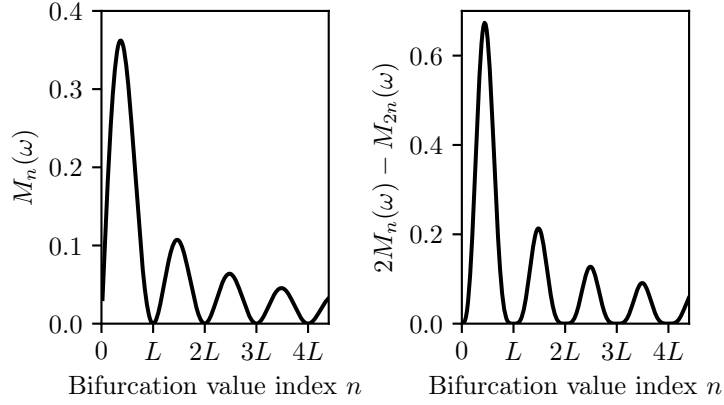
\begin{figure}\centering
    \input{BifurcationTypeDetermination.pgf}
    \caption[Bifurcation type]{Left: Shows the function $M_n(\omega)$ defined
    in equation~\eqref{Defn:Mn} for uniform $\omega$
    (see~\ref{OmegaChoice:1}). When $L = n$ the term $M_n(\omega) = 0$.
    Right: Show the function $2 M_n(\omega) - M_{2n}(\omega)$ whose sign
    determines whether we have a sub or supercritical pitchfork
    bifurcation\index{subcritical bifurcation}\index{supercritical
    bifurcation}\index{pitchfork bifurcation}.
    Here the function is always non-negative, hence we always have
    supercritical pitchfork bifurcations.
    }\label{Fig:TypBifs}
\end{figure}

\begin{figure}\centering
 \includegraphics[width=0.9\textwidth]{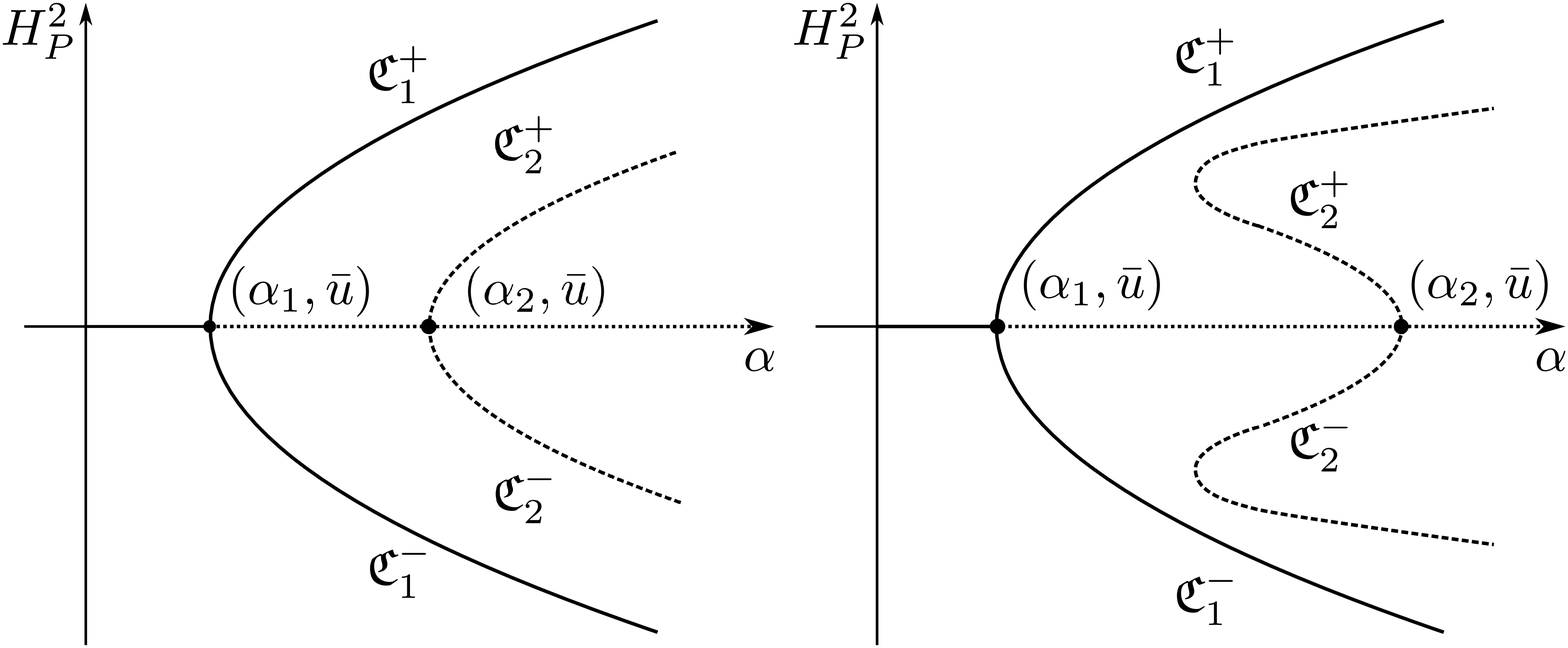}
    \caption[Bifurcation diagrams supercritical and
    subcritical\index{subcritical bifurcation} pitchfork bifurcations]
    {Left: Sample bifurcation diagram of the supercritical\index{supercritical
    bifurcation} pitchfork
    bifurcation,\index{pitchfork bifurcation} \ie,
    we have that $\alpha_{3,1}, \alpha_{3,2} > 0$ and have a switch of stability at the
    bifurcation point.
    Right: Sample bifurcation diagram of the subcritical pitchfork bifurcation, \ie,
    we have that $\alpha_{3,1}>0,\ \alpha_{3,2} < 0$.}\label{Fig:SuperPitchfork}
\end{figure}

\begin{theorem}\label{Thm:BifType}
    Suppose that all the assumptions of \cref{Thm:GlobalBif} are
    satisfied. Then the bifurcation branch\index{bifurcation branch} $\Gamma_n$ in a neighborhood of
    $(\alpha_n, \mean{u})$ is parameterized by
    \[
        \begin{bmatrix}
            u_n(x, s) \\
            \alpha_n(s)
        \end{bmatrix} =
        \begin{bmatrix}
            \mean{u} \\
            \alpha_n
        \end{bmatrix}
        + s \begin{bmatrix}
            \alpha_n \\
            0
        \end{bmatrix}
        \cos\lb\frac{2 \pi n x}{L}\rb
        + s^2 \begin{bmatrix}
            p_1(x) \\
            \alpha_{n, 3}
        \end{bmatrix}
        + o(s^2)
    \]
    for $s \in (-\delta, \delta)$, where $p_1(x)$ is determined below. Since
    ${\alpha_n}^{\prime}(0) = 0$ this is a pitchfork bifurcation, whose
    direction is determined by the sign of $\alpha_{3,n}$, its value is given by
    \[
        \alpha_{3, n} = \frac{1}{4 \mean{u}^5} \lb\frac{\pi n}{L}\rb^3
            \lsb M_n^2(\omega) \lb 2M_n(\omega) - M_{2n}(\omega)\rb \rsb^{-1}.
    \]
    where $M_n(\omega)$ is defined in \cref{Lem:EfuncK}.
    The sign of $\alpha_{3, n}$ is thus determined by the sign of the difference
    $\Delta M_n \coloneqq 2M_n(\omega) - M_{2n}(\omega)$.
\end{theorem}
\begin{proof}
From \cref{Thm:CrandallContinuity} we have that the bifurcating branches
are of class $\C^3$ in particular $(\alpha_n(s), u_n(x, s)) \in \C^3$ with
respect to $s$. Then we can write an asymptotic expansion of
$u_n(x,s), \alpha_n(s)$ for the $n$-th bifurcation branch
(see \cref{Thm:LocBif}), that is,
\[
\begin{split}
    u_n(s, x) &= \mean{u} + s \alpha_n \cos\lb\frac{2\pi n x}{L}\rb + s^2 p_1(x)
        + s^3 p_2(x) + o(s^3) \\
    \alpha_n(s) &= \alpha_n + s \alpha_{n, 2} + s^2 \alpha_{n, 3} + o(s^3),
\end{split}
\]
where $p_i \in Y_n$. Since $\Avg[u_n(s,x)] = \mean{u}$, we have that $\Avg[p_i] = 0$.
For the following discussion, we will also need Fourier expressions for the
functions $p_i(x)$. That is, we express both by
\begin{equation}\label{Eqn:pFourExp}
    p_i(x) = \sum_{\substack{k = 1 \\ k\neq n}}^{\infty} b_k^i \cos\lb\frac{2\pi k x}{L}\rb.
\end{equation}
Then we substitute the asymptotic expansion into equation~\eqref{Eqn:StSt2}, and
we group the result by associating terms of equal powers of $s$.

The terms of order $\BigOh{s}$ give
\[
    \alpha_n e_{n}^{\prime\prime} - \alpha_{n}^{2} \mean{u} {\K[e_n]}^{\prime} = 0.
\]
It is straightforward to verify that this equation is satisfied precisely when
$\alpha_n$ is a bifurcation point (see \cref{Lem:OpEvals}).

The terms of order $\BigOh{s^2}$ give:
\begin{equation}\label{Eqn:2ndExp}
    \underbrace{p_1^{\prime\prime}(x)}_{\ROM{1}} -
    \alpha_n \underbrace{\mean{u} {\K[p_1]}^{\prime}}_{\ROM{2}} -
        \alpha_n^3 \underbrace{\lb e_n \K[e_n] \rb^{\prime}}_{\ROM{3}} -
        \alpha_n^2 \alpha_{2, n} \underbrace{\mean{u} {\K[e_n]}^{\prime}}_{\ROM{4}}
        = 0.
\end{equation}
We will project this equation onto the space spanned by $e_n$, and we use roman
numerals (\ie\ $\ROM{1}$) to refer to those projected terms.  First we show
that $\alpha_{2, n} = 0$. Indeed, by projecting equation~\eqref{Eqn:2ndExp}
onto the nullspace of $\D_u \F(\alpha_n, \mean{u})$, we will obtain the result.
We will show the work term by term in order as they appear in
equation~\eqref{Eqn:2ndExp}. First, since $p_1 \in Y_n$
\[
    \ROM{1} = \int_{0}^{L}
        p_1^{\prime\prime}(x) \cos\lb\frac{2\pi n x}{L}\rb \dd x = 0,
\]
by integration by parts. For the second term, we apply the results of
\cref{Lem:EfuncKp} to obtain
\[
    {\K[p_1]}^{\prime} = - \frac{4 \pi}{L}
        \sum_{\substack{k = 1 \\k\neq n}}^{\infty}
        k b_k^1 M_k(\omega) \cos\lb\frac{2\pi k x}{L}\rb.
\]
Then projecting those terms on the nullspace easily shows that
\[
    \ROM{2} = 0.
\]
We proceed similarly for the third term. First, using
\cref{Lem:EfuncK}, we obtain
\begin{equation}\label{Eqn:ComSimpl}
    \lb e_n \K[e_n] \rb^{\prime} = \lb\frac{-4\pi n}{L}\rb M_n(\omega) \cos\lb\frac{4\pi n x}{L}\rb.
\end{equation}
Then carrying out the projection onto $e_n$, we obtain
\[
\ROM{3} = \int_{0}^{L} \lb e_n \K[e_n] \rb^\prime \cos\lb\frac{2\pi n x}{L}\rb \dd x
        = 0.
\]
Finally, for the fourth term,  using \cref{Lem:EfuncKp}, we obtain
\[
    \ROM{4} = \frac{- 4\pi n}{L} M_n(\omega)
        \int_{0}^{L} \cos^2\lb\frac{2\pi n x}{L}\rb \dd x = -2 n \pi M_n(\omega).
\]
Substituting all these results into equation~\eqref{Eqn:2ndExp},
we obtain that
\[
    \alpha_{2, n} = 0.
\]

The terms of order $\BigOh{s^3}$ are
\begin{equation}\label{Eqn:3rdExp}
    \underbrace{p_2^{\prime\prime}(x)}_{\ROM{1}}
    - \alpha_n \underbrace{\mean{u} {\K[p_2]}^{\prime}}_{\ROM{2}}
        - \alpha_n^2 \underbrace{\lb e_n \K[p_1] \rb^{\prime}}_{\ROM{3}}
        - \alpha_n^2 \underbrace{\lb p_1 \K[e_n] \rb^{\prime}}_{\ROM{4}}
        - \alpha_n \alpha_{3, n} \mean{u} \underbrace{{\K[e_n]}^{\prime}}_{\ROM{5}} = 0.
\end{equation}
Again we project this equation onto the nullspace of $\D_u \F(\alpha_n, \mean{u})$
term by term. The first term is
\[
    \ROM{1} = \int_{0}^{L} p_2^{\prime\prime}(x)
        \cos\lb\frac{2\pi n x}{L}\rb \dd x = 0.
\]
The second term, projected onto the nullspace is given by
\[
    \ROM{2} = - \frac{4 \pi}{L} \sum_{\substack{k = 1\\k \neq n}}^{\infty}
        \int_{0}^{L} k b_k^2 M_k(\omega) \cos\lb\frac{2\pi k x}{L}\rb
        \cos\lb\frac{2\pi n x}{L}\rb \dd x = 0.
\]
The third term is the same as the third term above, hence $\ROM{3} = 0$.

The fourth term is more interesting. Following integration by parts, we obtain
\[
    \ROM{4} = \frac{2\pi n}{L} \int_{0}^{L} p_1(x) \K[e_n]
        \sin\lb\frac{2\pi n x}{L}\rb \dd x.
\]
In this case, after substituting in the Fourier expansion for $p_1$ and using
\cref{Lem:EfuncKp} to rewrite the non-local term we
encounter integral terms of the form
\[
    \int_{0}^{L} \cos\lb\frac{2\pi k x}{L}\rb\sin^{2}\lb\frac{2\pi n x}{L}\rb \dd x
         = \begin{cases}
             0 &\mbox{if } k \neq 2n \\
             - L / 4 &\mbox{if } k = 2n.
            \end{cases}
\]
This means only the term $k = 2n$ remains from the Fourier expansion, we
obtain
\[
    \ROM{4} = n \pi M_n(\omega) b_{2n}^1.
\]
Finally, we deal with the last term
\[
\begin{split}
    \ROM{5} &= \int_{0}^{L} {\K[e_n]}^{\prime} \cos\lb\frac{2\pi n x}{L}\rb \dd x \\
        &= -\frac{4\pi n}{L} M_n(\omega)
        \int_{0}^{L} \cos^2\lb\frac{2\pi n x}{L}\rb \dd x
        = - 2 M_n(\omega) \pi n.
\end{split}
\]
Substituting each of the projections into
equation~\eqref{Eqn:3rdExp}, we can solve for $\alpha_{3, n}$.
\[
    \alpha_{3, n} = \frac{\alpha_n}{2 \mean{u}}
        \underbrace{\int_{0}^{L} p_1(x) \cos\lb\frac{4 \pi n x}{L}\rb \dd x}_{b_{2n}^1}.
\]
Thus next we will have to find $b_{2n}^1$. To find it, we solve the equation of
order $\BigOh{s^2}$ for $p_1(x)$, and recalling the result from
equation~\eqref{Eqn:ComSimpl}, we then obtain
\[
    p_1^{\prime\prime}(x) - \alpha_n \mean{u} {\K[p_1]}^{\prime}
    = \alpha_n^3 \lb e_n \K[e_n] \rb^\prime
    = - \alpha_n^3 \frac{4 \pi n}{L} M_n(\omega) \cos\lb\frac{4\pi n x}{L}\rb.
\]
Using the Fourier expansion~\eqref{Eqn:pFourExp} of $p_1(x)$
and substituting it into the previous equation and matching modes, we obtain
\[
    b_{2n}^1 = \frac{1}{2 \mean{u}^3} \lb\frac{\pi n}{L}\rb^2
        \lsb 2 M_n^2(\omega) - M_n(\omega) M_{2n}(\omega) \rsb^{-1}.
\]
which we can substitute into our expression for $\alpha_{3,n}$. We find that
\[
    \alpha_{3, n} = \frac{1}{4 \mean{u}^5} \lb\frac{\pi n}{L}\rb^3
        \lsb M_n^2(\omega) \lb 2M_n(\omega) - M_{2n}(\omega)\rb \rsb^{-1}.
\]
\end{proof}
\begin{example}\label{Exmp:UnifOmega}
    Suppose that $\omega$ is chosen to be the uniform function\index{uniform function}
    \ie~\ref{OmegaChoice:1}. Using the explicit form of $M_n(\omega)$
    computed in \cref{Example:UniformOmega1}, we can compute
    $\alpha_{3, n}$. In this case, and find that
    \[
        \alpha_{3, n} = \frac{1}{\mean{u}^5} \lb\frac{\pi n}{L}\rb^6
        \csc^8\lb\frac{\pi n}{L}\rb,
    \]
    which is always positive. For a plot $M_n(\omega)$ and $2 M_n - M_{2n}$
    refer to \cref{Fig:TypBifs}.
    If in addition, $L = 2, \mean{u} = 1$, then we
    only have bifurcations for $n$ odd (since $M_n(\omega) = 0$ for even $n$).
    We find that
    \[
        \alpha_n        = \frac{{(n \pi)}^2}{2}, \qquad
        \alpha_{3, n}   = \lb\frac{n \pi}{2}\rb^6.
    \]
    A bifurcation diagram for this situation is shown in \cref{Fig:SuperPitchfork}.
\end{example}

\begin{corollary}\label{Cor:DelMPos}
    Let $\omega \equiv \ifrac{1}{2}$, and $n$ be such that $M_n(\omega) > 0$,
    then $\Delta M_n > 0$. In other words, when $\omega$ is the uniform
    integration kernel only super-critical pitchfork bifurcations are possible.
\end{corollary}
\begin{proof}
    Since $M_n(\omega) > 0$ we must have that $L \neq n$. Using the expression
    for $M_n(\omega)$ found in \cref{Example:UniformOmega1} we find
    \[
        \frac{M_{2n}(\omega)}{2 M_n(\omega)} = \cos^2\lb\frac{\pi n}{L}\rb < 1.
    \]
\end{proof}

While the proof of \cref{Thm:GlobalBif} is currently restricted to the
situations of uniform $\omega(r)$, we can show that if it holds for general
$\omega(r)$ we can prove the following corollary, which shows that the first
bifurcation point must always be super-critical.

\begin{corollary}
    The $n$-th bifurcation must be a super-critical pitchfork bifurcation i.e.\
    at bifurcation point $(\alpha_n, \mean{u})$, whenever $n < \ifrac{L}{2}$.
    Note that since $L > 2$, this means that the first bifurcation is always
    supercritical\index{supercritical bifurcation}.
\end{corollary}

\begin{proof}
    By \cref{Thm:BifType} it requires to check the sign of
    $\Delta M_n(\omega)$. We have that
    \[
        \Delta M_n(\omega) = \int_{0}^{1} \lsb 2\sin\lb\frac{2\pi n r}{L}\rb
            - \sin\lb\frac{4\pi n r}{L}\rb \rsb \omega(r) \dd r.
    \]
    Using a change of variables we can rewrite this as
    \[
        \Delta M_n(\omega) = \frac{4L}{\pi n}
            \int_{0}^{\pi n / L} \sin^3(y) \cos(y) \omega\lb\frac{L}{\pi n} y\rb
            \dd y.
    \]
    Then since $n < \ifrac{L}{2}$, we have that $\ifrac{\pi n}{L} <
    \ifrac{\pi}{2}$. This means that the integration limits are confined to the
    first quadrant, and hence the integral must be strictly positive.
\end{proof}

\begin{figure}\centering
    \input{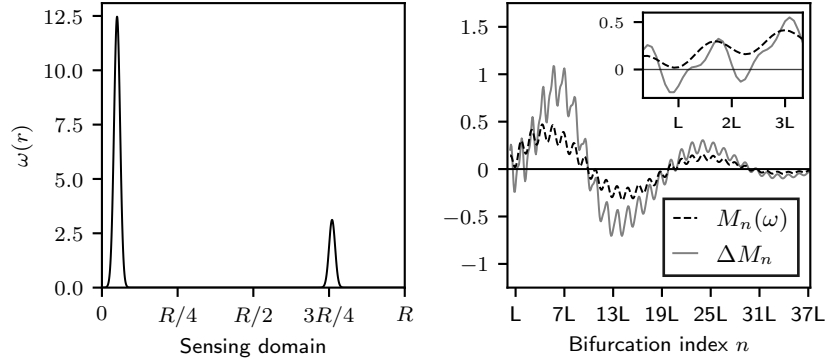}
    \caption[Subcritical Example]{Left: The integral kernel $\omega(r)$ for
    which we have a subcritical bifurcation points. For its construction refer
    to \cref{Ex:SubOmg}. Right: Plot of $M_n(\omega)$
    and $\Delta M_n(\omega)$. It is evident that due to the oscillations in
    $M_n(\omega)$ there are locations at which $\Delta M_n(\omega)$ is negative,
    leading to negative values of $\alpha_{3,n}$.
    In dashed black we have $M_n(\omega)$ and in solid gray we have $\Delta
    M_n(\omega)$.}\label{Fig:SubcriticalExample}
\end{figure}

\begin{example}[Construction of sub-critical bifurcation $\omega(r)$]\label{Ex:SubOmg}
 The goal is to construct a $\omega(r)$ such that $M_n(\omega) > 0$ but $\Delta
 M_n(\omega) < 0$. Recall that
 \[
  M_n(\omega) = \int_{0}^{1} \sin\lb\frac{2\pi nr}{L}\rb \omega(r) \dd r
 \]
and
 \[
  \Delta M_n(\omega) = 4 \int_{0}^{1} \sin\lb\frac{2\pi nr}{L}\rb\omega(r)\sin^2\lb\frac{\pi nr}{L}\rb \dd r.
 \]
 Next we define
 \[
  I^+ = \lcb x \in [0, 1] : \sin\lb\frac{2\pi nr}{L}\rb > 0\rcb, \quad
  I^- = \lcb x \in [0, 1] : \sin\lb\frac{2\pi nr}{L}\rb < 0\rcb.
 \]
 Assuming that $M_n(\omega) > 0$ implies that
 \[
  \int_{I^+} \sin\lb\frac{2\pi nr}{L}\rb\omega(r)\dd r >
  -\int_{I^-}\sin\lb\frac{2\pi nr}{L}\rb\omega(r)\dd r.
 \]

 Next we consider the sum of two point measures
 \[
  \omega(r) = a_1 \delta_{r_1}(r) + a_2 \delta_{r_2}(r)
 \]
 where we require that $a_1 + a_2 = \ifrac{1}{2}$, and $r_1 \in I^+$, $r_2 \in
 I^-$. We want to choose $a_1, a_2$ such that $\Delta M_n < 0$ i.e.\
 \[
  \int_{I^+} \sin\lb\frac{2\pi nr}{L}\rb\omega(r)\sin^2\lb\frac{\pi nr}{L}\rb\dd r
   < -\int_{I^-}\sin\lb\frac{2\pi nr}{L}\rb\omega(r)\sin^2\lb\frac{\pi nr}{L}\rb\dd r.
 \]
 Plugging in $\omega(r)$ this gives two conditions for $a_1$ and $a_2$.
 \[
  s_1 \coloneqq
  - \frac{\sin\lb\dfrac{2\pi n r_2}{L}\rb}{\sin\lb\dfrac{2\pi n r_1}{L}\rb} <
  \frac{a_1}{a_2}, \quad
  s_2 \coloneqq - \frac{\sin\lb\dfrac{2\pi n r_2}{L}\rb\sin^2\lb\dfrac{\pi n
  r_2}{L}\rb}{\sin\lb\dfrac{2\pi n r_1}{L}\rb\sin^2\lb\dfrac{\pi n r_1}{L}\rb}
  > \frac{a_1}{a_2}
 \]
 Then, in the limit of $r_1\to 0^+$ and $r_2\to \frac{L}{2n}^+$ we find
 \[
  \frac{a_1}{a_2} > \lim_{\substack{r_1 \to 0\\ r_2 \to \ifrac{L}{2n}}} s_1 = 1,\quad
  \frac{a_1}{a_2} < \lim_{\substack{r_1 \to 0\\ r_2 \to \ifrac{L}{2n}}} s_2 = \infty.
 \]
 Hence for $r_1$ close to $0$ and $r_2$ close to $\frac{L}{2n}$ we can easily
 find $a_1$ and $a_2$ that satisfy the above condition.  In
 \cref{Fig:SubcriticalExample} we chose $\ifrac{a_1}{a_2} = 8$, $L=3$, $n=2$,
 which gives
 \[ a_1 = \frac{8}{18} \qquad\mbox{and}\qquad a_2 = \frac{1}{18}.\]
 We mollify the point measures with sharp Gaussians, denoting them with
 $\delta^h$ as mollified versions of $\delta$ and we get
 \[ \omega(r) =\frac{8}{18}\delta_{r_1}^h(r) + \frac{1}{18}\delta_{r_2}^h(r), \]
 which is shown in \cref{Fig:SubcriticalExample} on the left. In
 \cref{Fig:SubcriticalExample} on the right we plot $M_n(\omega)$ and $\Delta
 M_n(\omega)$ as functions of $n$. We see that $\Delta M_n<0$ in areas where
 $M_n>0$, leading to possible backward bifurcations\index{backward bifurcation}.
\end{example}

\section{Stability of solutions}
So far we have studied the set of solutions of the equation $\F[\alpha, u] = 0$,
which are the steady states of the evolution equation
\begin{equation}\label{Eqn:DynamicalSys}
    \frac{\dd u}{\dd t} = -\F[\alpha, u].
\end{equation}
In this section, we are interested in the linear stability\index{linear stability} of these steady state
solutions. Then the linear stability of such a solution along a branch, is
determined by the sign of the eigenvalue of $\D_u \F(\alpha(s), u(s))$. An
eigenvalue perturbation result proven in \parencite{Crandall1973}, shows that the
eigenvalue along the trivial solution branch is related to the eigenvalues along
the non-trivial solution branch near a bifurcation point. An application of the
main result of \parencite{Crandall1973} is the goal of this subsection. Again we
limit this to the first solution branch. The eigenvalue problem resulting from
linearizing equation~\eqref{Eqn:DynamicalSys} around solutions in $\Gamma_n$
near the bifurcation point $(\alpha_n, \mean{u})$ is given by
\begin{equation}\label{Eqn:EVP}
    - w^{\prime\prime} + \alpha_n(s) \lb w \K[u(s)] + u(s) \K[w] \rb^{\prime}
        = -\lambda w.
\end{equation}
When $s = 0$, $u(0) = \mean{u}$, we have that
\begin{equation}\label{Eqn:EVP0}
    - w^{\prime\prime} + \alpha_n \mean{u} \lb \K[w] \rb^{\prime} = -\lambda w.
\end{equation}
\begin{lemma}\label{Lem:EvalZero}
    The $k$-th eigenvalue of the eigenvalue problem in
    equation~\eqref{Eqn:EVP0} is given by
    \[
        \lambda_k = \begin{cases}
                        0 &\mbox{if } n = k \\
                        \lb\dfrac{2k \pi}{L}\rb^2 \lb
                        \dfrac{n}{k} \dfrac{M_k(\omega)}{M_n(\omega)} - 1 \rb &\mbox{if } n \neq k
                    \end{cases}.
    \]
    Since $M_k(\omega) \to 0$ as $k \to \infty$, we see that
    $\lambda_k \to - \infty$ as $k \to \infty$.
\end{lemma}
\begin{proof}
    Apply \cref{Lem:EfuncKp} and \cref{Lemma:BoundMn}.
\end{proof}
\begin{example}\label{Exmp:UnifEvals}
    Suppose that $\omega$ is chosen to be the uniform
    function, \ie~\ref{OmegaChoice:1} (see \cref{section:NonLocalOperator}), and
    compute $\lambda_k$. We find that $\lambda_k > 0$ for $k < n$, and
    $\lambda_k < 0$ for $k > n$.
\end{example}
To apply the results of \parencite{Crandall1973} we introduce the following definition.
\begin{definition}[Definition 1.2
    \parencite{Crandall1973}]\label{Def:MSimpEval}
    Let $\mathcal{T}, \mathcal{M} \in \Lin(X, Y)$. Then $\mu \in \R$ is a
    $\mathcal{M}$-simple eigenvalue\index{$\mathcal{M}$-simple eigenvalue}
    of $\mathcal{T}$ if
    \[
        \dim N[\mathcal{T} - \mu \mathcal{M}] =
        \codim R[\mathcal{T} - \mu \mathcal{M}] = 1,
    \]
    and if $N[\mathcal{T} - \mu \mathcal{M}] = \vspan\{x_0\}$ we have that
    \[
        \mathcal{M} x_0 \notin R[\mathcal{T} - \mu \mathcal{M}].
    \]
\end{definition}
For the purpose here, we define the operator
$\M : H^2 \to L^2_0 \times \R$ by
\[
    \M[w] = \begin{pmatrix} w(x) - \Avg[w] \\
                            0
            \end{pmatrix}.
\]
\begin{lemma}
    $\lambda = 0$ is a $\mathcal{M}$-simple eigenvalue of
    $\D_u \F(\alpha_n,\mean{u})$.
\end{lemma}
\begin{proof}
    Recall that $\D_u \F(\alpha_n, \mean{u}) : H^2 \mapsto L^2_0 \times \R$ is
    Fredholm\index{Fredholm operator} with index~0. Thus, the operator satisfies the
    first condition in \cref{Def:MSimpEval}. We establish the second condition
    by contradiction. Suppose that
    \[
        \mathcal{M}[e_n] \in R[\D_u \F(\alpha, \mean{u})] =
            \{ u \in L^2 : \Avg[u] = 0 \}.
    \]
    In other words, there exists $w \in H^2$ such that
    \begin{equation}\label{Eqn:MEvalPrb}
        \left\{
            \begin{array}{@{}ll@{}}
                - w^{\prime\prime} + \alpha_n \mean{u} {\K[w]}^\prime
                    = \alpha_n \cos\lb\frac{2n \pi x}{L}\rb
                \quad &\mbox{in } \lsb 0, L \rsb  \\
                \Avg[w] = 0.
        \end{array}
        \right.
    \end{equation}
    We expand $w(x)$ using a Fourier series
    \[
        w(x) = \sum_{k = 1}^{\infty} w_k \cos\lb\frac{2\pi k x}{L}\rb,
    \]
    and substitute that into equation~\eqref{Eqn:MEvalPrb}. We then
    obtain the equation
    \[
        \lb \frac{2 k \pi}{L} \rb^2 \lsb 1 - \frac{n}{k}\frac{M_n(\omega)}{M_k(\omega)} \rsb
         = \begin{cases} 0 &\mbox{if } n \neq k \\
                        \alpha_n &\mbox{if } n = k,
            \end{cases}
    \]
    which leads to a contradiction when $n = k$.
\end{proof}
\begin{theorem}\label{Thm:EvalSign}
    Let all the assumptions of \cref{Thm:GlobalBif} hold.
    Then for any $\mean{u} > 0$, and for $s \in (-\delta, 0) \cup
    (0, \delta)$ the sign of the smallest magnitude eigenvalue
    of the solution $u(s, x)$ given by equation~\eqref{Eqn:LocBifFun}
    of~\eqref{Eqn:StSt2}, has opposite sign of $\alpha_{3,n}$,
    in the class of functions such that $\Avg[u] = \mean{u}$.
\end{theorem}
\begin{proof}
Now we are ready to apply~\cite[Theorem 1.16]{Crandall1973}. That implies
that there are open intervals $I, J$ with $\alpha_n \in I,\ 0 \in J$, chosen
such that
\[
    \gamma : I \mapsto \R \qquad \mu : J \mapsto \R,
\]
satisfying
\[
    \gamma(\alpha_n) = \mu(0) = 0,
\]
and
\[
    u : I \mapsto H^2 \qquad w : J \mapsto H^2
\]
satisfying
\[
    u(\alpha_n) = e_n = w(0), \qquad u(\alpha) - e_n \in Y_n, \qquad
    w(s) - e_n \in Y_n.
\]
Then we have two eigenvalue problems
\begin{align}\label{Eqn:EvalPrbs}
    \D_u \F(\alpha, \mean{u})[u(\alpha)] &= -\gamma(\alpha) \mathcal{M}
    [u(\alpha)],\quad\mbox{for}\ \alpha \in I, \\
    \D_u \F(\alpha(s), u(s))[w(\alpha)] &= -\mu(s) \mathcal{M} [w(s)],
        \quad\mbox{for}\ s \in J.
\end{align}
Whenever $\mu(s) \neq 0$, we have that
\begin{equation}\label{Eqn:EvalLim}
    \lim_{s \to 0, \mu(s) \neq 0} \frac{-s \alpha^{\prime}(s)
    \gamma^{\prime}(\alpha_n)}{\mu(s)} = 1.
\end{equation}
Thus we are left with computing $\gamma^{\prime}(\alpha_n)$.
 We differentiate
equation~\eqref{Eqn:EvalPrbs} with respect to $\alpha$ and
then set $\alpha = \alpha_n$ to obtain
\[
    \underbrace{- \dot{u}^{\prime\prime}}_{\ROM{1}} +
    \underbrace{\mean{u} {\K[e_n]}^{\prime}}_{\ROM{2}} +
    \alpha_n \mean{u} \underbrace{{\K[\dot{u}]}^{\prime}}_{\ROM{3}}
    = -\underbrace{\gamma^{\prime}(\alpha_n) e_n}_{\ROM{4}}.
\]
We again project to the null space and use roman numerals to denote the projected terms.
We multiply this equation by $e_n$ and integrate by parts. We obtain
the following results term-wise. The first term gives us
\[
    \ROM{1} = \lb\frac{2 n \pi}{L}\rb^2 \int_{0}^{L} \dot{u}
        \cos\lb\frac{2n \pi x}{L}\rb \dd x.
\]
The second term gives us, recalling the result from
\cref{Lem:EfuncKp}
\[
    \ROM{2} = - 2 \mean{u} n \pi M_n(\omega).
\]
The third term gives us, applying \cref{Lem:EfuncK} that
\[
\begin{split}
    \ROM{3} &= \lb \K[\dot{u}]', e_n \rb_{L^2}
        = \lb \dot{u}, \K[e_n]' \rb_{L^2} \\
        &= \frac{- 4n \pi}{L} M_n(\omega) \int_{0}^{L} \dot{u}
        \cos\lb\frac{2n \pi x}{L}\rb \dd x.
\end{split}
\]
Finally, the last term gives
\[
    \ROM{4} = -\gamma^{\prime}(\alpha_n) \frac{L}{2}.
\]
Combining all the terms we get
\[
    \underbrace{\lsb \lb \frac{2n \pi}{L} \rb^2 - \alpha_n \mean{u} \frac{4n\pi}{L}
        M_n(\omega) \rsb}_{\ROM{5}} \int_{0}^{L} \dot{u} \cos\lb\frac{2n\pi x}{L}\rb \dd x
        - 2n\mean{u} \pi M_n(\omega) = -\gamma^{\prime}(\alpha_n) \frac{L}{2}.
\]
Note that when substituting $\alpha_n$ into $\ROM{5}$, we see that this term
is zero, and thus we obtain
\[
    \gamma^{\prime}(\alpha_n) = \lb \frac{2n\pi}{L} \rb^2 \frac{1}{\alpha_n}.
\]
Substituting $\gamma'$ and the local expansion for $\alpha(s)$ from
\cref{Thm:BifType} into equation~\eqref{Eqn:EvalLim}, we
find that
\[
    \lim_{s \to 0, \mu(s) \neq 0} -\frac{2}{\alpha_n} \lb \frac{2 n \pi s}{L} \rb^2
        \frac{\alpha_{3,n}}{\mu(s)} = 1.
\]
We conclude, that
\[
    \sgn \mu(s) = - \sgn \alpha_{3,n},
\]
for $s \in (0, \delta) \cup (-\delta, 0)$.
\end{proof}

\begin{corollary}\label{Cor:StabBranch}
    Let all the assumptions of \cref{Thm:EvalSign} hold, and let
    $\omega \equiv \ifrac{1}{2}$. Then the first bifurcation branch is linearly
    asymptotically stable\index{linearly asymptotically stable}
    near the bifurcation point $(\alpha_1, \mean{u})$.
    All higher higher bifurcation branches are saddles\index{saddles}.
\end{corollary}

\begin{proof}
    Since $n$ is such that $M_n(\omega) > 0$ we have that from
    \cref{Cor:DelMPos} that $\Delta M_n(\omega) > 0$, implying that
    $\alpha_{3,n} > 0$ for all possible $n$. Applying
    \cref{Thm:EvalSign} we find that $\mu(s) < 0$.
    In the case that $n=1$ we have from \cref{Exmp:UnifEvals} that all
    eigenvalues $\lambda_k < 0$ for $k \geq 2$. By standard eigenvalue
    perturbation, we conclude that all eigenvalues have negative sign in a small
    neighborhood of $(\alpha_1, \mean{u})$. Thus, solutions in $\Gamma_1$ near
    the bifurcation point are linearly asymptotically stable.

    In the cases $n \geq 2$, we find eigenvalues $\lambda_k > 0$ for $k < n$ and
    $\lambda_k > 0$ for $k > n$. Thus, in a small neighborhood of the
    bifurcation point the solutions in the branches $\Gamma_n,\ n\geq 2$ are
    saddles.
\end{proof}

\section{Numerical verification}\label{subsec:numer}
\begin{figure}\centering
    \includegraphics[width=\textwidth]{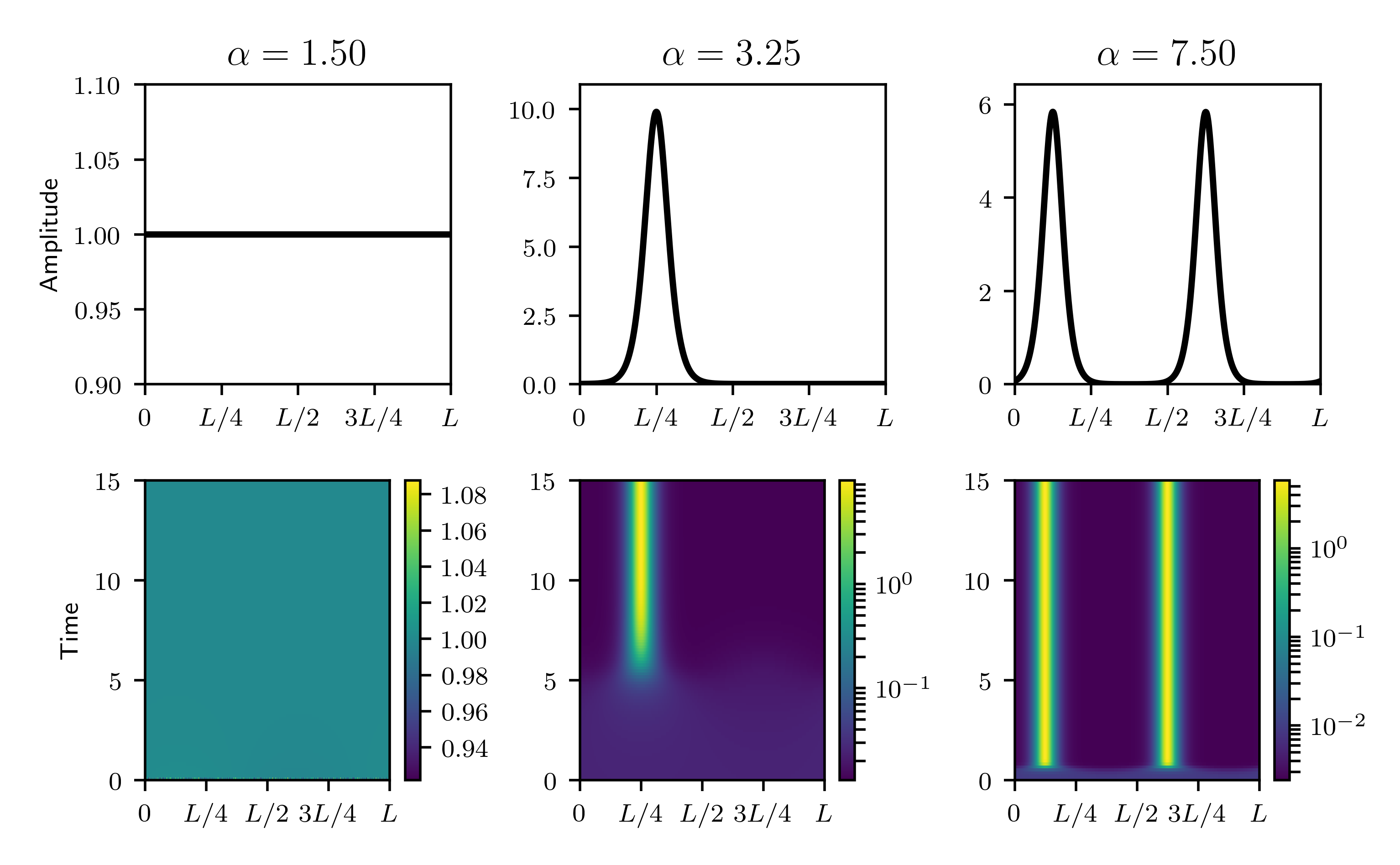}
    \caption{Numerical solutions of equation~\eqref{Eqn:3:AdhModel}
    on the periodic domain $S_L^1$. In the top row we show the final solution
    profiles, while below are the kymographs.
    (Left) $\alpha = 1.5$, (Middle) $\alpha = 3.25$, (Right) $\alpha = 7.5$.
    Here $L = 5$.}\label{Fig:PerNumRes}
\end{figure}

In this section, we verify the predictions of \cref{Thm:GlobalBif} on the steady
states of equation~\eqref{Eqn:3:AdhModel} by solving this equation numerically.
In addition, we also explore numerical solutions in cases not covered by
\cref{Thm:GlobalBif} such as, non-constant $\omega(r)$ or nonlinearities $h(u)$
within the non-local term. For all numerical solutions in this section, we pick
the domain size $L = 5$, and $\mean{u} = 1$.  From \cref{Lem:OpEvals} we know
that the first three bifurcation points are located at
\[
    \alpha_1 = \frac{16\pi^2}{25(5 - \sqrt{5})}, \qquad
    \alpha_2 = \frac{64\pi^2}{25(5 + \sqrt{5})}, \quad
    \alpha_3 = \frac{144\pi^2}{25(5 + \sqrt{5})}.
\]
This roughly means that $\alpha_1 \sim 2.28$, $\alpha_2 \sim 3.49$, $\alpha_3
\sim 7.85$. For all subsequent numerical simulations we pick a value of $\alpha$
from each of the intervals $(0, \alpha_1)$, $(\alpha_1, \alpha_2)$, and
$(\alpha_2, \alpha_3)$. The short-time numerical solutions of
equation~\eqref{Eqn:3:AdhModel} are presented in \cref{Fig:PerNumRes}. The top
row shows the final solution profiles, while the bottom row shows a
kymograph\index{kymograph}
with the spatial information on the $x$-axis, and time on the $y$-axis. These
numerical solutions\index{numerical solutions} have three key properties: (1) when the value of $\alpha <
\alpha_1$ the solution is constant and equals $\mean{u}$, (2) when $\alpha \in
(\alpha_1, \alpha_2)$ the solution has a single peak, while when $\alpha \in
(\alpha_2, \alpha_3)$ the solution has two peaks, (3) the peaks are uniformly
spaced on the domain. In addition, it is straightforward to check that these
solutions are symmetric under the actions of $\O2,\ \mathbf{D}_1$, and
$\mathbf{D}_2$ respectively.  These observations match the predictions of
\cref{Thm:LocBif} and \cref{Thm:GlobalBif}, and (4) the translational
symmetry\index{translational symmetry} of equation~\eqref{Eqn:StSt} is on
display, since solution peaks may form anywhere in the domain.

\begin{figure}\centering
 \includegraphics[width=0.9\textwidth]{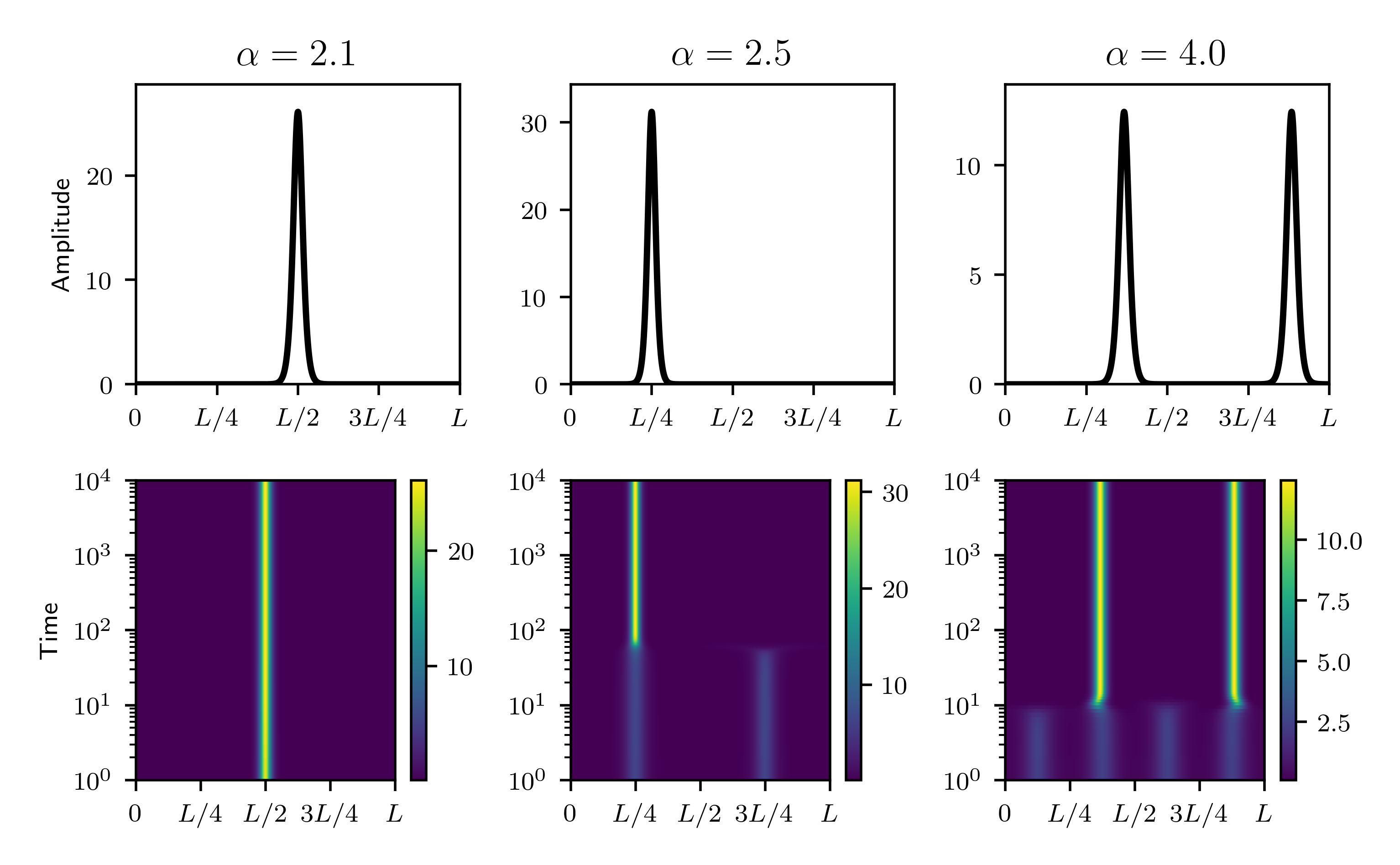}
 \captionof{figure}{Long term numerical solutions of
    equation~\eqref{Eqn:3:AdhModel} on the periodic domain $S_L^1$.
    The final solution profiles remain unchanged by either extending the
    simulation up to $10^{13}$, or reducing the solver's error tolerance by a
    factor of $100$. (Left) for $\alpha = 2.1$ the single
    peak that formed initially is stable, while (Middle) for $\alpha = 2.5$
    the initial pattern with two peaks coarsens into a singled peaked pattern.
    (Right) $\alpha = 4.0$ the initial pattern of four peaks, coarsens once to a
    double peak solution, which does not undergo additional coarsening.
    Here $L = 10$.}\label{Fig:PerNumResLong}
\end{figure}

Next we ask the question whether these patterns persist for long time. In
similar equations, such as models of chemotaxis, it is well known that patterns
with many peaks undergo coarsening\index{coarsening} \cite{Painter2011}. Indeed,
our numerical solutions of equation~\eqref{Eqn:3:AdhModel} exhibit similar
coarsening (see
\cref{Fig:PerNumResLong}). \cref{Cor:StabBranch} states that near the
bifurcation point solutions with two or more spikes are saddles of
equation~\eqref{Eqn:3:AdhModel}. While \cref{Cor:StabBranch} is not valid far
away from bifurcation points, the numerical results suggest this to remain true.
Interestingly however, the right kymograph in \cref{Fig:PerNumResLong} shows a
coarsening from $4$ to $2$ spikes, without a further coarsening to a single
spike (even if the simulation time is increased to $t_f = 10^{13}$).

\begin{figure}\centering
    \includegraphics[width=\textwidth]{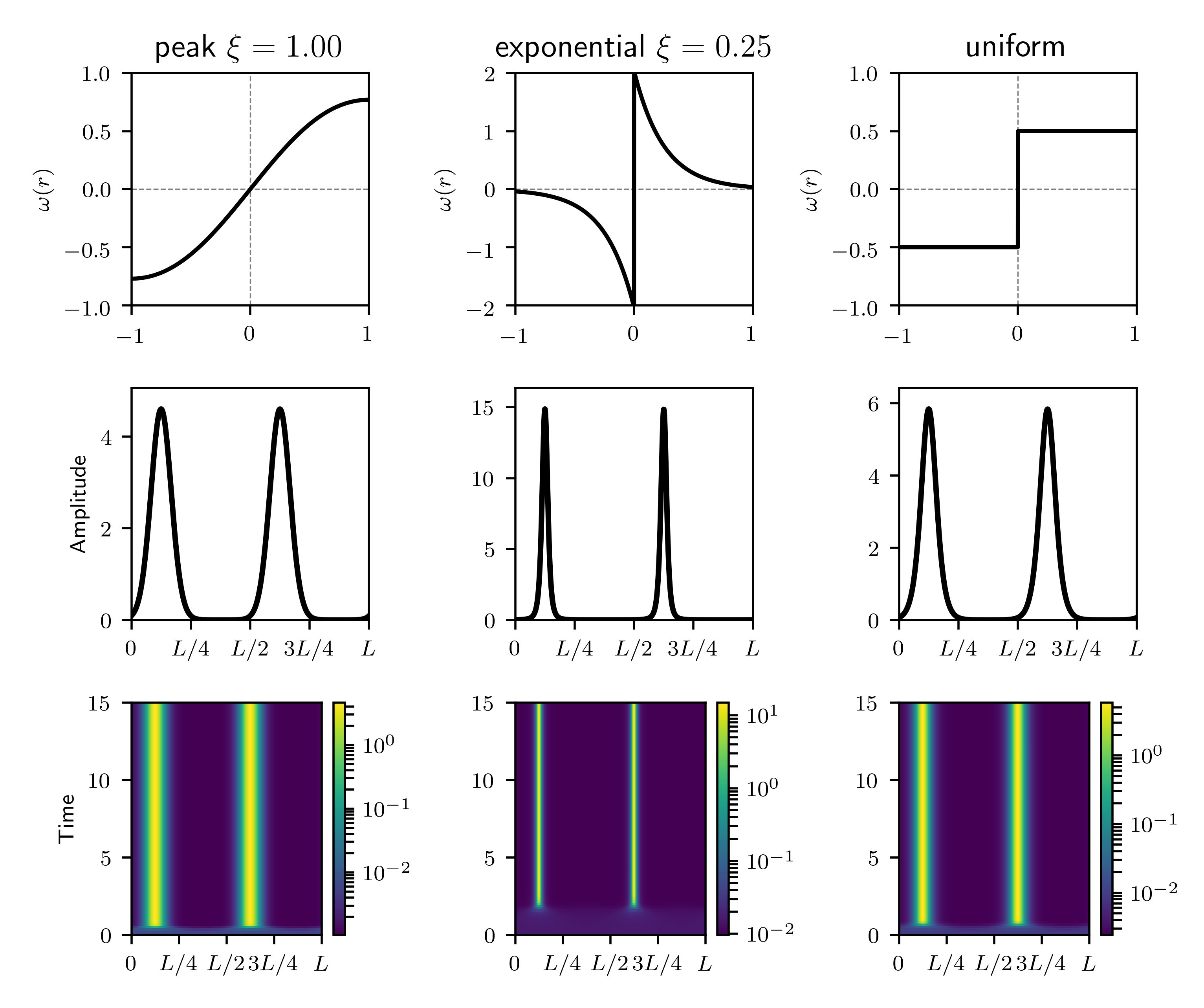}
    \caption{A comparison of solutions for the three different choices
    of the kernel $\omega(r)$~\ref{OmegaChoice:1}-\ref{OmegaChoice:3}
    introduced in \cref{section:NonLocalOperator}. For each simulation $L = 5$,
    $\alpha = 7.5$, which in each case is greater than $\alpha_2$ \ie\ the
    second bifurcation point. It was verified by computing $\Delta M_n(\omega)$
    that in each case a switch of stability occurs at the bifurcation point
    $\alpha_2$. Here $L = 5$.}\label{Fig:KernelComp}
\end{figure}

While \cref{Thm:GlobalBif} is not valid for non-constant $\omega(r)$, we can use
\cref{Lem:OpEvals} to determine the bifurcation points from the trivial
solutions. In three selected cases of non-constant $\omega(r)$, we numerically explore the
solutions of equation~\eqref{Eqn:3:AdhModel}. The final solution profiles, and
the corresponding kymographs for $\alpha \in (\alpha_2, \alpha_3)$ are shown in
\cref{Fig:KernelComp}. In each case it is straight forward to verify that the
solution profiles are $\mathbf{D}_2$ symmetric, and resemble the solutions for constant
$\omega(r)$. Additional, numerical simulations (not shown here), suggest that
once again only the single spiked solution is stable on long time scales.

\begin{figure}\centering
 \includegraphics[width=\textwidth]{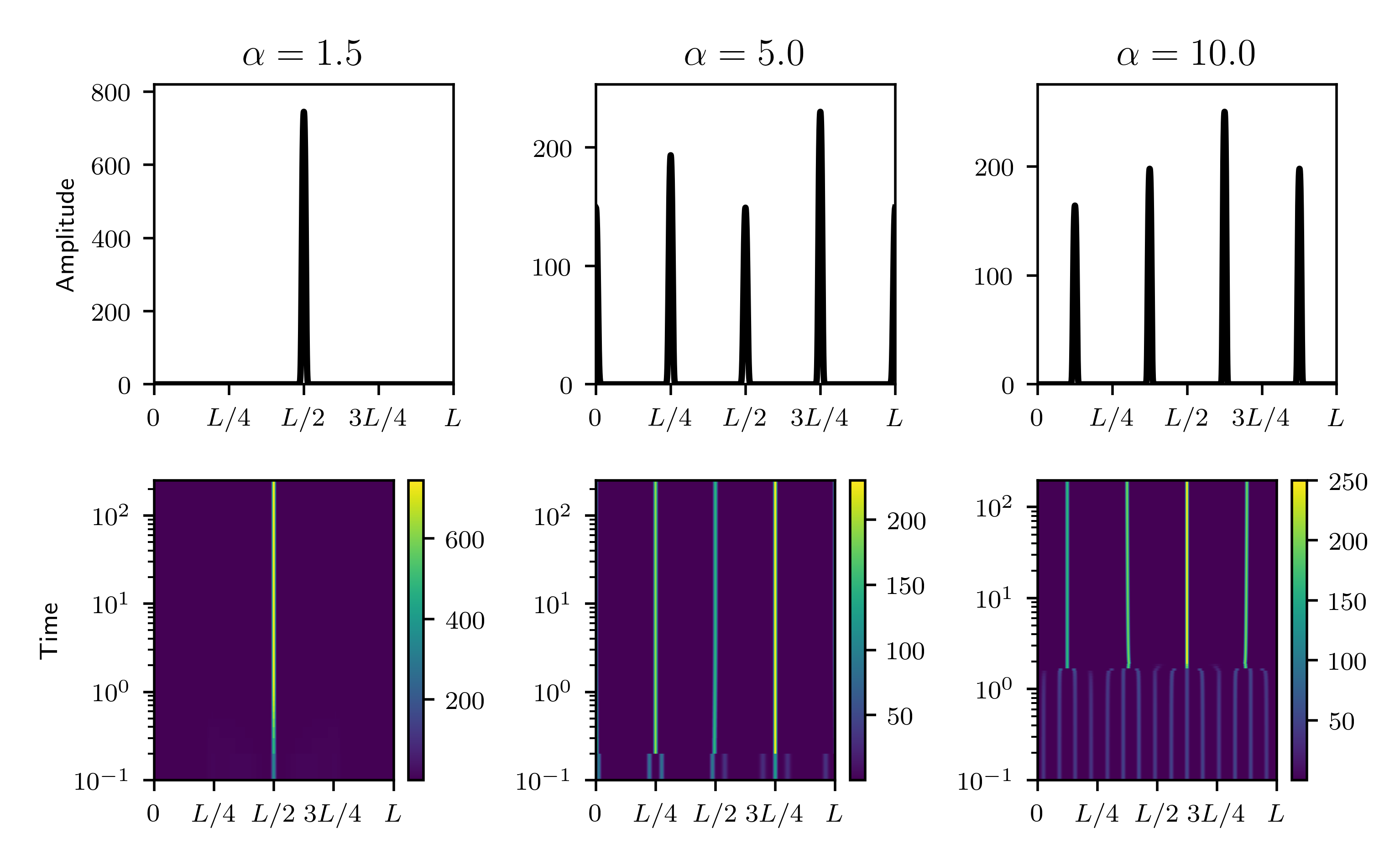}
 \caption{The same $\omega(r)$ as in \cref{Fig:SubcriticalExample}. Here $L=3$.
 The structure of the solutions are much more complicated.}\label{Fig:SubcriticalNumerical}
\end{figure}

Finally to ascertain whether the limitation to constant $\omega(r)$ is a
technical short-coming we compute numerical solutions of
equation~\eqref{Eqn:3:AdhModel} using the integration kernel $\omega(r)$
constructed in \cref{Ex:SubOmg}, for which sub-critical bifurcations are a
possibility. The final solution profiles and kymographs are shown in
\cref{Fig:SubcriticalNumerical}. We notice that the spacing between adjacent
peaks, and peak heights vary. Some symmetric features are identifiable, but are
more complicated. Thus \cref{Thm:GlobalBif} does not apply. In conclusion this
means that \cref{Thm:GlobalBif} cannot blindly be extended to general
$\omega(r)$, but we must in detail understand what properties of $\omega(r)$
yield regular symmetric solutions, and which lead to more complicated solutions
of equation~\eqref{Eqn:StSt2}.

\subsection{Numerical implementation\index{numerical implementation}}
Equation~\eqref{Eqn:3:AdhModel} is solved using a method of lines\index{method
of lines} approach,
where the spatial derivatives are discretized to yield a large system of
time-dependent ODEs (MOL-ODEs\index{MOL-ODEs}). The domain $[0, L]$ is discretized into a
cell-centered grid of uniform length $h = 1/N$, where $N$ is the number of grid
cells per unit length. Here we set $N = 1024$. The discretization\index{discretization} of the
advection term utilities a high-order upwind scheme\index{upwind scheme} augmented with a
flux-limiting scheme\index{flux-limiting scheme} to ensure positivity of solutions. For full details on the
numerical method we refer to~\cite{Gerisch2001}.

The non-local term in equation~\eqref{Eqn:ArmstrongModel} presents challenges to
its efficient and accurate evaluation. Here we employ the scheme based on the
Fast Fourier Transform introduced in \parencite{Gerisch2010}. The resulting
system of ordinary differential equations are integrated using the ROWMAP
integrator\index{ROWMAP integrator} introduced in \parencite{Weiner1997a}. Here
we use the implementation provided by the authors of
\parencite{Weiner1997a}\footnote{\url{http://www.mathematik.uni-halle.de/wissenschaftliches_rechnen/forschung/software/}}.
The integrator (written in Fortran) was wrapped using
\texttt{f2py}\footnote{\url{https://docs.scipy.org/doc/numpy-dev/f2py/}} into a
\texttt{scipy}\footnote{\url{https://www.scipy.org/}} integrate class.  The
spatial discretization (right hand side of ODE) is implemented using
\texttt{NumPy}\footnote{\url{www.numpy.org}}. The integrators error tolerance is
set to $v_{tol} = 10^{-6}$.
\section{Summary}\label{sec:SummaryGlobalBif}
In this chapter, we established the existence of global solution branches of the
non-local equation~\eqref{Eqn:StSt}, where each solution branch originates from
the homogeneous solution (see \cref{Chapter:LocalBifPeriodic}).  While it is
straightforward to apply the abstract global bifurcation theorem
(\cref{Thm:UnilateralAbstractBifurcation}) following the establishment of the
local result, the true challenge lies in discerning which of its three
alternatives hold. For nonlinear Sturm-Liouville problems Rabinowitz \etal\
classified solution branches by the solution's number of zeros. This is based
on the fact that its solutions can only have simple zeros (\ie\ $f(x_0) = 0$,
but $f'(x_0) \neq 0$). This does not help in our case, since firstly, we are
searching for positive solutions of equation~\eqref{Eqn:StSt}. Secondly, for
second-order elliptic PDEs the nodal separation\index{nodal separation} is proven by transforming the
equation into a two-dimensional initial value problem (IVP). But in the case of
non-local equations this is impossible, since we must look both backwards and
forward. A different approach to classifying solution branches is to use
symmetries, see for instance Healey \etal\ (equivariant nonlinear
elliptic equations)\index{equivariant nonlinear equations}
\parencite{Healey1991}. Symmetries have also extensively used in applications
of nonlinear systems of reaction diffusion equations
\cite{Nishiura1982,Fujii1982,Fujii1983}. Here we also use
symmetries\index{symmetries}.

Indeed, we show that equation~\eqref{Eqn:StSt} is $\O2$ equivariant. Employing
ideas of equivariant bifurcation theory\index{equivariant bifurcation theory}
\cite{Healey1991,golubitsky2003,Buono2015}, we first collect at each bifurcation
point all elements of $\O2$ that leave the local nullspace invariant. This leads
to the isotropy subgroup\index{isotropy subgroup}, which here is the dihedral
group\index{dihedral group}.  Next for each
isotropy subgroup we collect all functions left invariant under its action, in a
fixed-point function space.  Combining the symmetry properties of the non-local
operator $\K$ (\cref{Lem:FixMaxMin}), and with the particular properties of
positive solutions of equation~\eqref{Eqn:StSt}
(\cref{Chapter:PropertiesOfSolutions}), we find that the zeros of the solution's
derivative have fixed spatial positions. Using those locations we introduce the
domain tiling\index{domain tiling} $\Omega_i$, on the top of which we construct
the spaces of spiky\index{space of spiky functions}
functions ($\Nodal^{\pm}_{n}$). In other words, we classify solution branches by
the solution's derivative number of zeros. What remains to be shown is that the
branches obtained from the abstract bifurcation theorem
\cref{Thm:UnilateralAbstractBifurcation} are contained in $\Nodal^{\pm}_{n}$.
The challenge is to ascertain that different branches (with different number of
zeros) cannot intersect. In other words, we must show that $u'(x)$ can only have
simple zeros. Finally note that due to \cref{Lem:ZerosDer} this question is
closely tied to the zeros of the non-local operator $\K[u]$.

Using the equation's symmetries, we reduce the problem to having to show that
$u'(x) \neq 0$ in each tile. Classically, for local equations this is
established using a maximum principle (see for instance
\cite{Wang2013,Xiang2013}). In the case of uniform integration kernel $\omega(r)
= \ifrac{1}{2}$, the non-local operator $\K$ and the solution's area function
$w$ are related by $\K[u] = \Delta_1 w$. On each tile the area function $w(x)$
satisfies a Dirichlet problem and the maximum principle implies that $w > 0$, or
$w < 0$. In a novel use of symmetries we confine the non-local operator $\K[u]$
to a single tile. Then the non-locality of the operator $\K[u]$ establishes that
$\K[u] \neq 0$ in any tile, and hence $u'(x) \neq 0$.  Due to its similarity to
classical maximum principles we refer to this result as a type of ``non-local''
maximum principle\index{non-local maximum principle}. This allows us to
establish that the branches obtained from the abstract bifurcation theorem
(\cref{Thm:UnilateralAbstractBifurcation}) are
contained in $\Nodal^{\pm}_{n}$, and establish \cref{Thm:GlobalBif}.

Using an asymptotic expansion near each bifurcation point, we find that
each bifurcation is of pitchfork type\index{pitchfork bifurcation}. Both
supercritical\index{supercritical bifurcation} and
subcritical\index{subcritical bifurcation} bifurcations are possible depending
on the sign of $\alpha_{3,n}$ (see
\cref{Thm:BifType}).  Using an eigenvalue perturbation result by Crandall \etal\
\parencite{Crandall1973} we show that the sign of the smallest magnitude
eigenvalue of the linearization about the non-trivial solution is opposite to the sign of
 $\alpha_{3,n}$. Interestingly, the sign of $\alpha_{3,n}$ is solely
determined by the properties of the integration kernel $\omega(r)$ (through
$M_n(\omega)$ and $M_{2n}(\omega)$). In all cases, we have that the bifurcation
of the first mode is always supercritical (\ie\ $\alpha_{3,1} > 0$ always).  For
uniform integration kernel (\ie\ $\omega(r) = \ifrac{1}{2}$) we have that all
bifurcations are supercritical (\ie\ $\alpha_{3,n}>0$ always). Some possible
bifurcation diagrams are shown in \cref{Fig:SuperPitchfork}).

The global bifurcation result \cref{Thm:GlobalBif} is limited to the linear
adhesion model (\ie\ $h(u) = u$) and to the uniform integration kernel (\ie\
$\omega(r) = \ifrac{1}{2}$).  Improving both limitations is not easy. Allowing
for more general $h(u)$ requires understanding how the area function\index{area function} $w(x)$ is
modified (see also \cref{Remark:LimitationLinear}), while allowing for different
$\omega(r)$ requires generalizing the ``non-local'' maximum principle (\ie\
result \cref{Lem:KnotZero}). Numerical simulations suggest that for commonly
used $\omega(r)$ a similar result should hold (see \cref{Fig:KernelComp}).
However \cref{Thm:BifType} and \cref{Lem:OpEvals} suggest that $\omega(r)$ may
strongly modify the equation's solutions. In fact it is possible to generate
integral kernels $\omega(r)$ for which \cref{Thm:BifType} predicts sub-critical
bifurcations (see \cref{Fig:SubcriticalExample}, and \cref{Ex:SubOmg}).
Numerical solutions such as in \cref{Fig:SubcriticalExample} suggest that in such cases
the solution's symmetries are more complicated, and \cref{Thm:GlobalBif} does no
longer hold. It is a highly desirable future goal to understand for which
integration kernels $\omega(r)$ \cref{Thm:GlobalBif} continues to hold and for
which ones this result has to be modified.

It has to be noted, that this result does not exclude the possibility of further
secondary bifurcations along the solution branches, which may break further
symmetries. The analysis of secondary bifurcations\index{secondary bifurcations} has remained a challenge.  A
possible mechanism for identifying further bifurcations is to monitor the
Leray-Schauder degree\index{Leray-Schauder degree} for further sign changes along solution branches. For this
task, one has to rely on a combination of numerical exploration and mathematical
ingenuity in using the equation's structure to identify secondary bifurcations.

Spiky-type patterns have been studied extensively in two component reaction
diffusion equations with no advection terms. An analytical theory to determine
the stability of spike patterns has been developed for the
Gierer-Meinhardt\index{Gierer-Meinhardt model}
\cite{Wei2008} and for Gray-Scott\index{Gray-Scott model} models. Similar, developments have been
carried out for reaction-diffusion models\index{reaction-diffusion model} of species segregation and
cross-diffusion\index{cross-diffusion} \cite{Kolokolnikov2011}. Common to these studies, is that the
analysis of the spectrum involves various classes of non-local eigenvalue
problems. Recently, similar techniques have been extended to include
reaction-diffusion equations with chemotactic drift, such as models of
urban crime\index{models of urban crime}
\cite{Kolokolnikov2014,Tse2016,Buttenschoen2018}, others have focused on
the well-known Keller-Segel model\index{Keller-Segel model}
\cite{Hillen2004,Kang2005,Potapov2005,Sleeman2005}. A highly desirable goal
would be to extend these methods to include non-local adhesion
models~\eqref{Eqn:3:AdhModel}, and improve upon \cref{Cor:StabBranch}.

This abstract bifurcation analysis gives rise to several interesting modelling
observations. Most noticeable in our analysis were the mathematical properties
of the integral kernel $\Omega(r)$ of the non-local operator $\K[u]$. Following
the modelling work in \parencite{Buttenschon2017} this term determines how
likely it is that a cell protrusion reaches a particular target. In our analysis
the properties of $\Omega$ enter through the quantity $M_n(\omega)$ at several
key moments: (1) the sign of $M_n(\omega)$ determines whether or not we have a
bifurcation and (2) it determines whether we immediately observe a switch of
stability. The minimum adhesion strength ($\alpha_1$ in \cref{Lem:OpEvals})
which allows the formation of cell aggregates (non-trivial solutions) is reduced
with increasing domain size $L$, magnitude of $M_n(\omega)$, and size of
$h^{\prime}(\mean{u})$, while no parameter increases $\alpha_1$. Note that these statements only hold for the
non-dimensionalized equation.

Extensions to multiple population systems\index{multiple population systems} are highly desirable to study the
possible cell-sorting patterns\index{cell-sorting patterns}. In addition, such extensions promise more
intricate dynamics. A different extension could be to consider the model
variations proposed by Murakawa \etal\ \parencite{Murakawa2015}, and Carrillo
\etal\ \parencite{Carrillo2019}, who added a density dependent diffusion term.
An extension to higher spatial dimension would be worth while, to more
realistically study the formation of tissues. In fact, the notion of a tiling in
higher dimensions has been considered prior (see for example the work by Courant
who defined tilings in higher dimensions \parencite{Courant1923}). Finally,
beyond global existence of the time-dependent solutions of
equation~\eqref{Eqn:ArmstrongModel} little is known about the qualitative
long-time behaviour of its solutions, or the structure of its global attractor.

%% file: nodal_test_function.pgf
\begingroup%
\makeatletter%
\begin{pgfpicture}%
\pgfpathrectangle{\pgfpointorigin}{\pgfqpoint{2.353383in}{1.454471in}}%
\pgfusepath{use as bounding box, clip}%
\begin{pgfscope}%
\pgfsetbuttcap%
\pgfsetmiterjoin%
\definecolor{currentfill}{rgb}{1.000000,1.000000,1.000000}%
\pgfsetfillcolor{currentfill}%
\pgfsetlinewidth{0.000000pt}%
\definecolor{currentstroke}{rgb}{1.000000,1.000000,1.000000}%
\pgfsetstrokecolor{currentstroke}%
\pgfsetdash{}{0pt}%
\pgfpathmoveto{\pgfqpoint{0.000000in}{0.000000in}}%
\pgfpathlineto{\pgfqpoint{2.353383in}{0.000000in}}%
\pgfpathlineto{\pgfqpoint{2.353383in}{1.454471in}}%
\pgfpathlineto{\pgfqpoint{0.000000in}{1.454471in}}%
\pgfpathclose%
\pgfusepath{fill}%
\end{pgfscope}%
\begin{pgfscope}%
\pgfsetbuttcap%
\pgfsetmiterjoin%
\definecolor{currentfill}{rgb}{1.000000,1.000000,1.000000}%
\pgfsetfillcolor{currentfill}%
\pgfsetlinewidth{0.000000pt}%
\definecolor{currentstroke}{rgb}{0.000000,0.000000,0.000000}%
\pgfsetstrokecolor{currentstroke}%
\pgfsetstrokeopacity{0.000000}%
\pgfsetdash{}{0pt}%
\pgfpathmoveto{\pgfqpoint{0.198611in}{0.386111in}}%
\pgfpathlineto{\pgfqpoint{2.154772in}{0.386111in}}%
\pgfpathlineto{\pgfqpoint{2.154772in}{1.255860in}}%
\pgfpathlineto{\pgfqpoint{0.198611in}{1.255860in}}%
\pgfpathclose%
\pgfusepath{fill}%
\end{pgfscope}%
\begin{pgfscope}%
\pgfsetbuttcap%
\pgfsetroundjoin%
\definecolor{currentfill}{rgb}{0.000000,0.000000,0.000000}%
\pgfsetfillcolor{currentfill}%
\pgfsetlinewidth{0.803000pt}%
\definecolor{currentstroke}{rgb}{0.000000,0.000000,0.000000}%
\pgfsetstrokecolor{currentstroke}%
\pgfsetdash{}{0pt}%
\pgfsys@defobject{currentmarker}{\pgfqpoint{0.000000in}{-0.048611in}}{\pgfqpoint{0.000000in}{0.000000in}}{%
\pgfpathmoveto{\pgfqpoint{0.000000in}{0.000000in}}%
\pgfpathlineto{\pgfqpoint{0.000000in}{-0.048611in}}%
\pgfusepath{stroke,fill}%
}%
\begin{pgfscope}%
\pgfsys@transformshift{0.198611in}{0.386111in}%
\pgfsys@useobject{currentmarker}{}%
\end{pgfscope}%
\end{pgfscope}%
\begin{pgfscope}%
\pgftext[x=0.198611in,y=0.288889in,,top]{\rmfamily\fontsize{10.000000}{12.000000}\selectfont \(\displaystyle 0\)}%
\end{pgfscope}%
\begin{pgfscope}%
\pgfsetbuttcap%
\pgfsetroundjoin%
\definecolor{currentfill}{rgb}{0.000000,0.000000,0.000000}%
\pgfsetfillcolor{currentfill}%
\pgfsetlinewidth{0.803000pt}%
\definecolor{currentstroke}{rgb}{0.000000,0.000000,0.000000}%
\pgfsetstrokecolor{currentstroke}%
\pgfsetdash{}{0pt}%
\pgfsys@defobject{currentmarker}{\pgfqpoint{0.000000in}{-0.048611in}}{\pgfqpoint{0.000000in}{0.000000in}}{%
\pgfpathmoveto{\pgfqpoint{0.000000in}{0.000000in}}%
\pgfpathlineto{\pgfqpoint{0.000000in}{-0.048611in}}%
\pgfusepath{stroke,fill}%
}%
\begin{pgfscope}%
\pgfsys@transformshift{0.687651in}{0.386111in}%
\pgfsys@useobject{currentmarker}{}%
\end{pgfscope}%
\end{pgfscope}%
\begin{pgfscope}%
\pgftext[x=0.687651in,y=0.288889in,,top]{\rmfamily\fontsize{10.000000}{12.000000}\selectfont \(\displaystyle L/4\)}%
\end{pgfscope}%
\begin{pgfscope}%
\pgfsetbuttcap%
\pgfsetroundjoin%
\definecolor{currentfill}{rgb}{0.000000,0.000000,0.000000}%
\pgfsetfillcolor{currentfill}%
\pgfsetlinewidth{0.803000pt}%
\definecolor{currentstroke}{rgb}{0.000000,0.000000,0.000000}%
\pgfsetstrokecolor{currentstroke}%
\pgfsetdash{}{0pt}%
\pgfsys@defobject{currentmarker}{\pgfqpoint{0.000000in}{-0.048611in}}{\pgfqpoint{0.000000in}{0.000000in}}{%
\pgfpathmoveto{\pgfqpoint{0.000000in}{0.000000in}}%
\pgfpathlineto{\pgfqpoint{0.000000in}{-0.048611in}}%
\pgfusepath{stroke,fill}%
}%
\begin{pgfscope}%
\pgfsys@transformshift{1.176692in}{0.386111in}%
\pgfsys@useobject{currentmarker}{}%
\end{pgfscope}%
\end{pgfscope}%
\begin{pgfscope}%
\pgftext[x=1.176692in,y=0.288889in,,top]{\rmfamily\fontsize{10.000000}{12.000000}\selectfont \(\displaystyle L/2\)}%
\end{pgfscope}%
\begin{pgfscope}%
\pgfsetbuttcap%
\pgfsetroundjoin%
\definecolor{currentfill}{rgb}{0.000000,0.000000,0.000000}%
\pgfsetfillcolor{currentfill}%
\pgfsetlinewidth{0.803000pt}%
\definecolor{currentstroke}{rgb}{0.000000,0.000000,0.000000}%
\pgfsetstrokecolor{currentstroke}%
\pgfsetdash{}{0pt}%
\pgfsys@defobject{currentmarker}{\pgfqpoint{0.000000in}{-0.048611in}}{\pgfqpoint{0.000000in}{0.000000in}}{%
\pgfpathmoveto{\pgfqpoint{0.000000in}{0.000000in}}%
\pgfpathlineto{\pgfqpoint{0.000000in}{-0.048611in}}%
\pgfusepath{stroke,fill}%
}%
\begin{pgfscope}%
\pgfsys@transformshift{1.665732in}{0.386111in}%
\pgfsys@useobject{currentmarker}{}%
\end{pgfscope}%
\end{pgfscope}%
\begin{pgfscope}%
\pgftext[x=1.665732in,y=0.288889in,,top]{\rmfamily\fontsize{10.000000}{12.000000}\selectfont \(\displaystyle 3L/4\)}%
\end{pgfscope}%
\begin{pgfscope}%
\pgfsetbuttcap%
\pgfsetroundjoin%
\definecolor{currentfill}{rgb}{0.000000,0.000000,0.000000}%
\pgfsetfillcolor{currentfill}%
\pgfsetlinewidth{0.803000pt}%
\definecolor{currentstroke}{rgb}{0.000000,0.000000,0.000000}%
\pgfsetstrokecolor{currentstroke}%
\pgfsetdash{}{0pt}%
\pgfsys@defobject{currentmarker}{\pgfqpoint{0.000000in}{-0.048611in}}{\pgfqpoint{0.000000in}{0.000000in}}{%
\pgfpathmoveto{\pgfqpoint{0.000000in}{0.000000in}}%
\pgfpathlineto{\pgfqpoint{0.000000in}{-0.048611in}}%
\pgfusepath{stroke,fill}%
}%
\begin{pgfscope}%
\pgfsys@transformshift{2.154772in}{0.386111in}%
\pgfsys@useobject{currentmarker}{}%
\end{pgfscope}%
\end{pgfscope}%
\begin{pgfscope}%
\pgftext[x=2.154772in,y=0.288889in,,top]{\rmfamily\fontsize{10.000000}{12.000000}\selectfont \(\displaystyle L\)}%
\end{pgfscope}%
\begin{pgfscope}%
\pgfpathrectangle{\pgfqpoint{0.198611in}{0.386111in}}{\pgfqpoint{1.956161in}{0.869749in}}%
\pgfusepath{clip}%
\pgfsetrectcap%
\pgfsetroundjoin%
\pgfsetlinewidth{1.003750pt}%
\definecolor{currentstroke}{rgb}{0.000000,0.000000,0.000000}%
\pgfsetstrokecolor{currentstroke}%
\pgfsetdash{}{0pt}%
\pgfpathmoveto{\pgfqpoint{0.198611in}{1.216326in}}%
\pgfpathlineto{\pgfqpoint{0.218370in}{1.213145in}}%
\pgfpathlineto{\pgfqpoint{0.238130in}{1.203654in}}%
\pgfpathlineto{\pgfqpoint{0.257889in}{1.188007in}}%
\pgfpathlineto{\pgfqpoint{0.277648in}{1.166453in}}%
\pgfpathlineto{\pgfqpoint{0.297407in}{1.139341in}}%
\pgfpathlineto{\pgfqpoint{0.317166in}{1.107107in}}%
\pgfpathlineto{\pgfqpoint{0.336926in}{1.070268in}}%
\pgfpathlineto{\pgfqpoint{0.356685in}{1.029419in}}%
\pgfpathlineto{\pgfqpoint{0.376444in}{0.985216in}}%
\pgfpathlineto{\pgfqpoint{0.396203in}{0.938370in}}%
\pgfpathlineto{\pgfqpoint{0.415962in}{0.889636in}}%
\pgfpathlineto{\pgfqpoint{0.435722in}{0.839796in}}%
\pgfpathlineto{\pgfqpoint{0.455481in}{0.789655in}}%
\pgfpathlineto{\pgfqpoint{0.475240in}{0.740017in}}%
\pgfpathlineto{\pgfqpoint{0.494999in}{0.691682in}}%
\pgfpathlineto{\pgfqpoint{0.514758in}{0.645428in}}%
\pgfpathlineto{\pgfqpoint{0.534518in}{0.601998in}}%
\pgfpathlineto{\pgfqpoint{0.554277in}{0.562093in}}%
\pgfpathlineto{\pgfqpoint{0.574036in}{0.526352in}}%
\pgfpathlineto{\pgfqpoint{0.593795in}{0.495353in}}%
\pgfpathlineto{\pgfqpoint{0.613554in}{0.469593in}}%
\pgfpathlineto{\pgfqpoint{0.633314in}{0.449487in}}%
\pgfpathlineto{\pgfqpoint{0.653073in}{0.435359in}}%
\pgfpathlineto{\pgfqpoint{0.672832in}{0.427435in}}%
\pgfpathlineto{\pgfqpoint{0.692591in}{0.425844in}}%
\pgfpathlineto{\pgfqpoint{0.712350in}{0.430611in}}%
\pgfpathlineto{\pgfqpoint{0.732110in}{0.441659in}}%
\pgfpathlineto{\pgfqpoint{0.751869in}{0.458811in}}%
\pgfpathlineto{\pgfqpoint{0.771628in}{0.481790in}}%
\pgfpathlineto{\pgfqpoint{0.791387in}{0.510227in}}%
\pgfpathlineto{\pgfqpoint{0.811146in}{0.543664in}}%
\pgfpathlineto{\pgfqpoint{0.830906in}{0.581564in}}%
\pgfpathlineto{\pgfqpoint{0.850665in}{0.623315in}}%
\pgfpathlineto{\pgfqpoint{0.870424in}{0.668248in}}%
\pgfpathlineto{\pgfqpoint{0.890183in}{0.715638in}}%
\pgfpathlineto{\pgfqpoint{0.909942in}{0.764723in}}%
\pgfpathlineto{\pgfqpoint{0.929702in}{0.814713in}}%
\pgfpathlineto{\pgfqpoint{0.949461in}{0.864804in}}%
\pgfpathlineto{\pgfqpoint{0.969220in}{0.914190in}}%
\pgfpathlineto{\pgfqpoint{0.988979in}{0.962077in}}%
\pgfpathlineto{\pgfqpoint{1.008738in}{1.007693in}}%
\pgfpathlineto{\pgfqpoint{1.028498in}{1.050305in}}%
\pgfpathlineto{\pgfqpoint{1.048257in}{1.089227in}}%
\pgfpathlineto{\pgfqpoint{1.068016in}{1.123834in}}%
\pgfpathlineto{\pgfqpoint{1.087775in}{1.153567in}}%
\pgfpathlineto{\pgfqpoint{1.107534in}{1.177949in}}%
\pgfpathlineto{\pgfqpoint{1.127294in}{1.196587in}}%
\pgfpathlineto{\pgfqpoint{1.147053in}{1.209181in}}%
\pgfpathlineto{\pgfqpoint{1.166812in}{1.215530in}}%
\pgfpathlineto{\pgfqpoint{1.186571in}{1.215530in}}%
\pgfpathlineto{\pgfqpoint{1.206330in}{1.209181in}}%
\pgfpathlineto{\pgfqpoint{1.226090in}{1.196587in}}%
\pgfpathlineto{\pgfqpoint{1.245849in}{1.177949in}}%
\pgfpathlineto{\pgfqpoint{1.265608in}{1.153567in}}%
\pgfpathlineto{\pgfqpoint{1.285367in}{1.123834in}}%
\pgfpathlineto{\pgfqpoint{1.305126in}{1.089227in}}%
\pgfpathlineto{\pgfqpoint{1.324886in}{1.050305in}}%
\pgfpathlineto{\pgfqpoint{1.344645in}{1.007693in}}%
\pgfpathlineto{\pgfqpoint{1.364404in}{0.962077in}}%
\pgfpathlineto{\pgfqpoint{1.384163in}{0.914190in}}%
\pgfpathlineto{\pgfqpoint{1.403922in}{0.864804in}}%
\pgfpathlineto{\pgfqpoint{1.423682in}{0.814713in}}%
\pgfpathlineto{\pgfqpoint{1.443441in}{0.764723in}}%
\pgfpathlineto{\pgfqpoint{1.463200in}{0.715638in}}%
\pgfpathlineto{\pgfqpoint{1.482959in}{0.668248in}}%
\pgfpathlineto{\pgfqpoint{1.502718in}{0.623315in}}%
\pgfpathlineto{\pgfqpoint{1.522478in}{0.581564in}}%
\pgfpathlineto{\pgfqpoint{1.542237in}{0.543664in}}%
\pgfpathlineto{\pgfqpoint{1.561996in}{0.510227in}}%
\pgfpathlineto{\pgfqpoint{1.581755in}{0.481790in}}%
\pgfpathlineto{\pgfqpoint{1.601514in}{0.458811in}}%
\pgfpathlineto{\pgfqpoint{1.621274in}{0.441659in}}%
\pgfpathlineto{\pgfqpoint{1.641033in}{0.430611in}}%
\pgfpathlineto{\pgfqpoint{1.660792in}{0.425844in}}%
\pgfpathlineto{\pgfqpoint{1.680551in}{0.427435in}}%
\pgfpathlineto{\pgfqpoint{1.700310in}{0.435359in}}%
\pgfpathlineto{\pgfqpoint{1.720070in}{0.449487in}}%
\pgfpathlineto{\pgfqpoint{1.739829in}{0.469593in}}%
\pgfpathlineto{\pgfqpoint{1.759588in}{0.495353in}}%
\pgfpathlineto{\pgfqpoint{1.779347in}{0.526352in}}%
\pgfpathlineto{\pgfqpoint{1.799106in}{0.562093in}}%
\pgfpathlineto{\pgfqpoint{1.818866in}{0.601998in}}%
\pgfpathlineto{\pgfqpoint{1.838625in}{0.645428in}}%
\pgfpathlineto{\pgfqpoint{1.858384in}{0.691682in}}%
\pgfpathlineto{\pgfqpoint{1.878143in}{0.740017in}}%
\pgfpathlineto{\pgfqpoint{1.897902in}{0.789655in}}%
\pgfpathlineto{\pgfqpoint{1.917662in}{0.839796in}}%
\pgfpathlineto{\pgfqpoint{1.937421in}{0.889636in}}%
\pgfpathlineto{\pgfqpoint{1.957180in}{0.938370in}}%
\pgfpathlineto{\pgfqpoint{1.976939in}{0.985216in}}%
\pgfpathlineto{\pgfqpoint{1.996698in}{1.029419in}}%
\pgfpathlineto{\pgfqpoint{2.016458in}{1.070268in}}%
\pgfpathlineto{\pgfqpoint{2.036217in}{1.107107in}}%
\pgfpathlineto{\pgfqpoint{2.055976in}{1.139341in}}%
\pgfpathlineto{\pgfqpoint{2.075735in}{1.166453in}}%
\pgfpathlineto{\pgfqpoint{2.095494in}{1.188007in}}%
\pgfpathlineto{\pgfqpoint{2.115254in}{1.203654in}}%
\pgfpathlineto{\pgfqpoint{2.135013in}{1.213145in}}%
\pgfpathlineto{\pgfqpoint{2.154772in}{1.216326in}}%
\pgfusepath{stroke}%
\end{pgfscope}%
\begin{pgfscope}%
\pgfpathrectangle{\pgfqpoint{0.198611in}{0.386111in}}{\pgfqpoint{1.956161in}{0.869749in}}%
\pgfusepath{clip}%
\pgfsetbuttcap%
\pgfsetroundjoin%
\pgfsetlinewidth{1.003750pt}%
\definecolor{currentstroke}{rgb}{0.000000,0.000000,0.000000}%
\pgfsetstrokecolor{currentstroke}%
\pgfsetdash{{3.700000pt}{1.600000pt}}{0.000000pt}%
\pgfpathmoveto{\pgfqpoint{0.198611in}{0.425645in}}%
\pgfpathlineto{\pgfqpoint{0.218370in}{0.428826in}}%
\pgfpathlineto{\pgfqpoint{0.238130in}{0.438316in}}%
\pgfpathlineto{\pgfqpoint{0.257889in}{0.453964in}}%
\pgfpathlineto{\pgfqpoint{0.277648in}{0.475518in}}%
\pgfpathlineto{\pgfqpoint{0.297407in}{0.502630in}}%
\pgfpathlineto{\pgfqpoint{0.317166in}{0.534864in}}%
\pgfpathlineto{\pgfqpoint{0.336926in}{0.571703in}}%
\pgfpathlineto{\pgfqpoint{0.356685in}{0.612552in}}%
\pgfpathlineto{\pgfqpoint{0.376444in}{0.656755in}}%
\pgfpathlineto{\pgfqpoint{0.396203in}{0.703601in}}%
\pgfpathlineto{\pgfqpoint{0.415962in}{0.752335in}}%
\pgfpathlineto{\pgfqpoint{0.435722in}{0.802174in}}%
\pgfpathlineto{\pgfqpoint{0.455481in}{0.852316in}}%
\pgfpathlineto{\pgfqpoint{0.475240in}{0.901954in}}%
\pgfpathlineto{\pgfqpoint{0.494999in}{0.950289in}}%
\pgfpathlineto{\pgfqpoint{0.514758in}{0.996543in}}%
\pgfpathlineto{\pgfqpoint{0.534518in}{1.039972in}}%
\pgfpathlineto{\pgfqpoint{0.554277in}{1.079878in}}%
\pgfpathlineto{\pgfqpoint{0.574036in}{1.115618in}}%
\pgfpathlineto{\pgfqpoint{0.593795in}{1.146618in}}%
\pgfpathlineto{\pgfqpoint{0.613554in}{1.172378in}}%
\pgfpathlineto{\pgfqpoint{0.633314in}{1.192484in}}%
\pgfpathlineto{\pgfqpoint{0.653073in}{1.206612in}}%
\pgfpathlineto{\pgfqpoint{0.672832in}{1.214536in}}%
\pgfpathlineto{\pgfqpoint{0.692591in}{1.216127in}}%
\pgfpathlineto{\pgfqpoint{0.712350in}{1.211360in}}%
\pgfpathlineto{\pgfqpoint{0.732110in}{1.200312in}}%
\pgfpathlineto{\pgfqpoint{0.751869in}{1.183160in}}%
\pgfpathlineto{\pgfqpoint{0.771628in}{1.160181in}}%
\pgfpathlineto{\pgfqpoint{0.791387in}{1.131744in}}%
\pgfpathlineto{\pgfqpoint{0.811146in}{1.098307in}}%
\pgfpathlineto{\pgfqpoint{0.830906in}{1.060407in}}%
\pgfpathlineto{\pgfqpoint{0.850665in}{1.018656in}}%
\pgfpathlineto{\pgfqpoint{0.870424in}{0.973723in}}%
\pgfpathlineto{\pgfqpoint{0.890183in}{0.926333in}}%
\pgfpathlineto{\pgfqpoint{0.909942in}{0.877248in}}%
\pgfpathlineto{\pgfqpoint{0.929702in}{0.827258in}}%
\pgfpathlineto{\pgfqpoint{0.949461in}{0.777167in}}%
\pgfpathlineto{\pgfqpoint{0.969220in}{0.727780in}}%
\pgfpathlineto{\pgfqpoint{0.988979in}{0.679894in}}%
\pgfpathlineto{\pgfqpoint{1.008738in}{0.634278in}}%
\pgfpathlineto{\pgfqpoint{1.028498in}{0.591666in}}%
\pgfpathlineto{\pgfqpoint{1.048257in}{0.552743in}}%
\pgfpathlineto{\pgfqpoint{1.068016in}{0.518137in}}%
\pgfpathlineto{\pgfqpoint{1.087775in}{0.488404in}}%
\pgfpathlineto{\pgfqpoint{1.107534in}{0.464022in}}%
\pgfpathlineto{\pgfqpoint{1.127294in}{0.445384in}}%
\pgfpathlineto{\pgfqpoint{1.147053in}{0.432789in}}%
\pgfpathlineto{\pgfqpoint{1.166812in}{0.426441in}}%
\pgfpathlineto{\pgfqpoint{1.186571in}{0.426441in}}%
\pgfpathlineto{\pgfqpoint{1.206330in}{0.432789in}}%
\pgfpathlineto{\pgfqpoint{1.226090in}{0.445384in}}%
\pgfpathlineto{\pgfqpoint{1.245849in}{0.464022in}}%
\pgfpathlineto{\pgfqpoint{1.265608in}{0.488404in}}%
\pgfpathlineto{\pgfqpoint{1.285367in}{0.518137in}}%
\pgfpathlineto{\pgfqpoint{1.305126in}{0.552743in}}%
\pgfpathlineto{\pgfqpoint{1.324886in}{0.591666in}}%
\pgfpathlineto{\pgfqpoint{1.344645in}{0.634278in}}%
\pgfpathlineto{\pgfqpoint{1.364404in}{0.679894in}}%
\pgfpathlineto{\pgfqpoint{1.384163in}{0.727780in}}%
\pgfpathlineto{\pgfqpoint{1.403922in}{0.777167in}}%
\pgfpathlineto{\pgfqpoint{1.423682in}{0.827258in}}%
\pgfpathlineto{\pgfqpoint{1.443441in}{0.877248in}}%
\pgfpathlineto{\pgfqpoint{1.463200in}{0.926333in}}%
\pgfpathlineto{\pgfqpoint{1.482959in}{0.973723in}}%
\pgfpathlineto{\pgfqpoint{1.502718in}{1.018656in}}%
\pgfpathlineto{\pgfqpoint{1.522478in}{1.060407in}}%
\pgfpathlineto{\pgfqpoint{1.542237in}{1.098307in}}%
\pgfpathlineto{\pgfqpoint{1.561996in}{1.131744in}}%
\pgfpathlineto{\pgfqpoint{1.581755in}{1.160181in}}%
\pgfpathlineto{\pgfqpoint{1.601514in}{1.183160in}}%
\pgfpathlineto{\pgfqpoint{1.621274in}{1.200312in}}%
\pgfpathlineto{\pgfqpoint{1.641033in}{1.211360in}}%
\pgfpathlineto{\pgfqpoint{1.660792in}{1.216127in}}%
\pgfpathlineto{\pgfqpoint{1.680551in}{1.214536in}}%
\pgfpathlineto{\pgfqpoint{1.700310in}{1.206612in}}%
\pgfpathlineto{\pgfqpoint{1.720070in}{1.192484in}}%
\pgfpathlineto{\pgfqpoint{1.739829in}{1.172378in}}%
\pgfpathlineto{\pgfqpoint{1.759588in}{1.146618in}}%
\pgfpathlineto{\pgfqpoint{1.779347in}{1.115618in}}%
\pgfpathlineto{\pgfqpoint{1.799106in}{1.079878in}}%
\pgfpathlineto{\pgfqpoint{1.818866in}{1.039972in}}%
\pgfpathlineto{\pgfqpoint{1.838625in}{0.996543in}}%
\pgfpathlineto{\pgfqpoint{1.858384in}{0.950289in}}%
\pgfpathlineto{\pgfqpoint{1.878143in}{0.901954in}}%
\pgfpathlineto{\pgfqpoint{1.897902in}{0.852316in}}%
\pgfpathlineto{\pgfqpoint{1.917662in}{0.802174in}}%
\pgfpathlineto{\pgfqpoint{1.937421in}{0.752335in}}%
\pgfpathlineto{\pgfqpoint{1.957180in}{0.703601in}}%
\pgfpathlineto{\pgfqpoint{1.976939in}{0.656755in}}%
\pgfpathlineto{\pgfqpoint{1.996698in}{0.612552in}}%
\pgfpathlineto{\pgfqpoint{2.016458in}{0.571703in}}%
\pgfpathlineto{\pgfqpoint{2.036217in}{0.534864in}}%
\pgfpathlineto{\pgfqpoint{2.055976in}{0.502630in}}%
\pgfpathlineto{\pgfqpoint{2.075735in}{0.475518in}}%
\pgfpathlineto{\pgfqpoint{2.095494in}{0.453964in}}%
\pgfpathlineto{\pgfqpoint{2.115254in}{0.438316in}}%
\pgfpathlineto{\pgfqpoint{2.135013in}{0.428826in}}%
\pgfpathlineto{\pgfqpoint{2.154772in}{0.425645in}}%
\pgfusepath{stroke}%
\end{pgfscope}%
\begin{pgfscope}%
\pgfpathrectangle{\pgfqpoint{0.198611in}{0.386111in}}{\pgfqpoint{1.956161in}{0.869749in}}%
\pgfusepath{clip}%
\pgfsetbuttcap%
\pgfsetroundjoin%
\pgfsetlinewidth{0.752812pt}%
\definecolor{currentstroke}{rgb}{0.000000,0.000000,0.000000}%
\pgfsetstrokecolor{currentstroke}%
\pgfsetdash{{2.775000pt}{1.200000pt}}{0.000000pt}%
\pgfpathmoveto{\pgfqpoint{0.687651in}{0.386111in}}%
\pgfpathlineto{\pgfqpoint{0.687651in}{1.255860in}}%
\pgfusepath{stroke}%
\end{pgfscope}%
\begin{pgfscope}%
\pgfpathrectangle{\pgfqpoint{0.198611in}{0.386111in}}{\pgfqpoint{1.956161in}{0.869749in}}%
\pgfusepath{clip}%
\pgfsetbuttcap%
\pgfsetroundjoin%
\pgfsetlinewidth{0.752812pt}%
\definecolor{currentstroke}{rgb}{0.000000,0.000000,0.000000}%
\pgfsetstrokecolor{currentstroke}%
\pgfsetdash{{2.775000pt}{1.200000pt}}{0.000000pt}%
\pgfpathmoveto{\pgfqpoint{1.176692in}{0.386111in}}%
\pgfpathlineto{\pgfqpoint{1.176692in}{1.255860in}}%
\pgfusepath{stroke}%
\end{pgfscope}%
\begin{pgfscope}%
\pgfpathrectangle{\pgfqpoint{0.198611in}{0.386111in}}{\pgfqpoint{1.956161in}{0.869749in}}%
\pgfusepath{clip}%
\pgfsetbuttcap%
\pgfsetroundjoin%
\pgfsetlinewidth{0.752812pt}%
\definecolor{currentstroke}{rgb}{0.000000,0.000000,0.000000}%
\pgfsetstrokecolor{currentstroke}%
\pgfsetdash{{2.775000pt}{1.200000pt}}{0.000000pt}%
\pgfpathmoveto{\pgfqpoint{1.665732in}{0.386111in}}%
\pgfpathlineto{\pgfqpoint{1.665732in}{1.255860in}}%
\pgfusepath{stroke}%
\end{pgfscope}%
\begin{pgfscope}%
\pgfpathrectangle{\pgfqpoint{0.198611in}{0.386111in}}{\pgfqpoint{1.956161in}{0.869749in}}%
\pgfusepath{clip}%
\pgfsetbuttcap%
\pgfsetroundjoin%
\pgfsetlinewidth{0.752812pt}%
\definecolor{currentstroke}{rgb}{0.000000,0.000000,0.000000}%
\pgfsetstrokecolor{currentstroke}%
\pgfsetdash{{2.775000pt}{1.200000pt}}{0.000000pt}%
\pgfpathmoveto{\pgfqpoint{0.198611in}{0.386111in}}%
\pgfpathlineto{\pgfqpoint{0.198611in}{1.255860in}}%
\pgfusepath{stroke}%
\end{pgfscope}%
\begin{pgfscope}%
\pgfsetrectcap%
\pgfsetmiterjoin%
\pgfsetlinewidth{0.803000pt}%
\definecolor{currentstroke}{rgb}{0.000000,0.000000,0.000000}%
\pgfsetstrokecolor{currentstroke}%
\pgfsetdash{}{0pt}%
\pgfpathmoveto{\pgfqpoint{0.198611in}{0.386111in}}%
\pgfpathlineto{\pgfqpoint{0.198611in}{1.255860in}}%
\pgfusepath{stroke}%
\end{pgfscope}%
\begin{pgfscope}%
\pgfsetrectcap%
\pgfsetmiterjoin%
\pgfsetlinewidth{0.803000pt}%
\definecolor{currentstroke}{rgb}{0.000000,0.000000,0.000000}%
\pgfsetstrokecolor{currentstroke}%
\pgfsetdash{}{0pt}%
\pgfpathmoveto{\pgfqpoint{2.154772in}{0.386111in}}%
\pgfpathlineto{\pgfqpoint{2.154772in}{1.255860in}}%
\pgfusepath{stroke}%
\end{pgfscope}%
\begin{pgfscope}%
\pgfsetrectcap%
\pgfsetmiterjoin%
\pgfsetlinewidth{0.803000pt}%
\definecolor{currentstroke}{rgb}{0.000000,0.000000,0.000000}%
\pgfsetstrokecolor{currentstroke}%
\pgfsetdash{}{0pt}%
\pgfpathmoveto{\pgfqpoint{0.198611in}{0.386111in}}%
\pgfpathlineto{\pgfqpoint{2.154772in}{0.386111in}}%
\pgfusepath{stroke}%
\end{pgfscope}%
\begin{pgfscope}%
\pgfsetrectcap%
\pgfsetmiterjoin%
\pgfsetlinewidth{0.803000pt}%
\definecolor{currentstroke}{rgb}{0.000000,0.000000,0.000000}%
\pgfsetstrokecolor{currentstroke}%
\pgfsetdash{}{0pt}%
\pgfpathmoveto{\pgfqpoint{0.198611in}{1.255860in}}%
\pgfpathlineto{\pgfqpoint{2.154772in}{1.255860in}}%
\pgfusepath{stroke}%
\end{pgfscope}%
\end{pgfpicture}%
\makeatother%
\endgroup%

%% file: typical_solution_with_w.pgf
\begingroup%
\makeatletter%
\begin{pgfpicture}%
\pgfpathrectangle{\pgfpointorigin}{\pgfqpoint{3.922305in}{2.424118in}}%
\pgfusepath{use as bounding box, clip}%
\begin{pgfscope}%
\pgfsetbuttcap%
\pgfsetmiterjoin%
\definecolor{currentfill}{rgb}{1.000000,1.000000,1.000000}%
\pgfsetfillcolor{currentfill}%
\pgfsetlinewidth{0.000000pt}%
\definecolor{currentstroke}{rgb}{1.000000,1.000000,1.000000}%
\pgfsetstrokecolor{currentstroke}%
\pgfsetdash{}{0pt}%
\pgfpathmoveto{\pgfqpoint{0.000000in}{0.000000in}}%
\pgfpathlineto{\pgfqpoint{3.922305in}{0.000000in}}%
\pgfpathlineto{\pgfqpoint{3.922305in}{2.424118in}}%
\pgfpathlineto{\pgfqpoint{0.000000in}{2.424118in}}%
\pgfpathclose%
\pgfusepath{fill}%
\end{pgfscope}%
\begin{pgfscope}%
\pgfsetbuttcap%
\pgfsetmiterjoin%
\definecolor{currentfill}{rgb}{1.000000,1.000000,1.000000}%
\pgfsetfillcolor{currentfill}%
\pgfsetlinewidth{0.000000pt}%
\definecolor{currentstroke}{rgb}{0.000000,0.000000,0.000000}%
\pgfsetstrokecolor{currentstroke}%
\pgfsetstrokeopacity{0.000000}%
\pgfsetdash{}{0pt}%
\pgfpathmoveto{\pgfqpoint{0.316407in}{0.386111in}}%
\pgfpathlineto{\pgfqpoint{1.729920in}{0.386111in}}%
\pgfpathlineto{\pgfqpoint{1.729920in}{2.051896in}}%
\pgfpathlineto{\pgfqpoint{0.316407in}{2.051896in}}%
\pgfpathclose%
\pgfusepath{fill}%
\end{pgfscope}%
\begin{pgfscope}%
\pgfsetbuttcap%
\pgfsetroundjoin%
\definecolor{currentfill}{rgb}{0.000000,0.000000,0.000000}%
\pgfsetfillcolor{currentfill}%
\pgfsetlinewidth{0.803000pt}%
\definecolor{currentstroke}{rgb}{0.000000,0.000000,0.000000}%
\pgfsetstrokecolor{currentstroke}%
\pgfsetdash{}{0pt}%
\pgfsys@defobject{currentmarker}{\pgfqpoint{0.000000in}{-0.048611in}}{\pgfqpoint{0.000000in}{0.000000in}}{%
\pgfpathmoveto{\pgfqpoint{0.000000in}{0.000000in}}%
\pgfpathlineto{\pgfqpoint{0.000000in}{-0.048611in}}%
\pgfusepath{stroke,fill}%
}%
\begin{pgfscope}%
\pgfsys@transformshift{0.316407in}{0.386111in}%
\pgfsys@useobject{currentmarker}{}%
\end{pgfscope}%
\end{pgfscope}%
\begin{pgfscope}%
\pgftext[x=0.316407in,y=0.288889in,,top]{\rmfamily\fontsize{10.000000}{12.000000}\selectfont \(\displaystyle 0\)}%
\end{pgfscope}%
\begin{pgfscope}%
\pgfsetbuttcap%
\pgfsetroundjoin%
\definecolor{currentfill}{rgb}{0.000000,0.000000,0.000000}%
\pgfsetfillcolor{currentfill}%
\pgfsetlinewidth{0.803000pt}%
\definecolor{currentstroke}{rgb}{0.000000,0.000000,0.000000}%
\pgfsetstrokecolor{currentstroke}%
\pgfsetdash{}{0pt}%
\pgfsys@defobject{currentmarker}{\pgfqpoint{0.000000in}{-0.048611in}}{\pgfqpoint{0.000000in}{0.000000in}}{%
\pgfpathmoveto{\pgfqpoint{0.000000in}{0.000000in}}%
\pgfpathlineto{\pgfqpoint{0.000000in}{-0.048611in}}%
\pgfusepath{stroke,fill}%
}%
\begin{pgfscope}%
\pgfsys@transformshift{1.023164in}{0.386111in}%
\pgfsys@useobject{currentmarker}{}%
\end{pgfscope}%
\end{pgfscope}%
\begin{pgfscope}%
\pgftext[x=1.023164in,y=0.288889in,,top]{\rmfamily\fontsize{10.000000}{12.000000}\selectfont \(\displaystyle L/2\)}%
\end{pgfscope}%
\begin{pgfscope}%
\pgfsetbuttcap%
\pgfsetroundjoin%
\definecolor{currentfill}{rgb}{0.000000,0.000000,0.000000}%
\pgfsetfillcolor{currentfill}%
\pgfsetlinewidth{0.803000pt}%
\definecolor{currentstroke}{rgb}{0.000000,0.000000,0.000000}%
\pgfsetstrokecolor{currentstroke}%
\pgfsetdash{}{0pt}%
\pgfsys@defobject{currentmarker}{\pgfqpoint{0.000000in}{-0.048611in}}{\pgfqpoint{0.000000in}{0.000000in}}{%
\pgfpathmoveto{\pgfqpoint{0.000000in}{0.000000in}}%
\pgfpathlineto{\pgfqpoint{0.000000in}{-0.048611in}}%
\pgfusepath{stroke,fill}%
}%
\begin{pgfscope}%
\pgfsys@transformshift{1.729920in}{0.386111in}%
\pgfsys@useobject{currentmarker}{}%
\end{pgfscope}%
\end{pgfscope}%
\begin{pgfscope}%
\pgftext[x=1.729920in,y=0.288889in,,top]{\rmfamily\fontsize{10.000000}{12.000000}\selectfont \(\displaystyle L\)}%
\end{pgfscope}%
\begin{pgfscope}%
\pgfsetbuttcap%
\pgfsetroundjoin%
\definecolor{currentfill}{rgb}{0.000000,0.000000,0.000000}%
\pgfsetfillcolor{currentfill}%
\pgfsetlinewidth{0.803000pt}%
\definecolor{currentstroke}{rgb}{0.000000,0.000000,0.000000}%
\pgfsetstrokecolor{currentstroke}%
\pgfsetdash{}{0pt}%
\pgfsys@defobject{currentmarker}{\pgfqpoint{-0.048611in}{0.000000in}}{\pgfqpoint{0.000000in}{0.000000in}}{%
\pgfpathmoveto{\pgfqpoint{0.000000in}{0.000000in}}%
\pgfpathlineto{\pgfqpoint{-0.048611in}{0.000000in}}%
\pgfusepath{stroke,fill}%
}%
\begin{pgfscope}%
\pgfsys@transformshift{0.316407in}{0.386111in}%
\pgfsys@useobject{currentmarker}{}%
\end{pgfscope}%
\end{pgfscope}%
\begin{pgfscope}%
\pgftext[x=0.149740in,y=0.338283in,left,base]{\rmfamily\fontsize{10.000000}{12.000000}\selectfont \(\displaystyle 0\)}%
\end{pgfscope}%
\begin{pgfscope}%
\pgfsetbuttcap%
\pgfsetroundjoin%
\definecolor{currentfill}{rgb}{0.000000,0.000000,0.000000}%
\pgfsetfillcolor{currentfill}%
\pgfsetlinewidth{0.803000pt}%
\definecolor{currentstroke}{rgb}{0.000000,0.000000,0.000000}%
\pgfsetstrokecolor{currentstroke}%
\pgfsetdash{}{0pt}%
\pgfsys@defobject{currentmarker}{\pgfqpoint{-0.048611in}{0.000000in}}{\pgfqpoint{0.000000in}{0.000000in}}{%
\pgfpathmoveto{\pgfqpoint{0.000000in}{0.000000in}}%
\pgfpathlineto{\pgfqpoint{-0.048611in}{0.000000in}}%
\pgfusepath{stroke,fill}%
}%
\begin{pgfscope}%
\pgfsys@transformshift{0.316407in}{0.772162in}%
\pgfsys@useobject{currentmarker}{}%
\end{pgfscope}%
\end{pgfscope}%
\begin{pgfscope}%
\pgftext[x=0.149740in,y=0.724334in,left,base]{\rmfamily\fontsize{10.000000}{12.000000}\selectfont \(\displaystyle 1\)}%
\end{pgfscope}%
\begin{pgfscope}%
\pgfsetbuttcap%
\pgfsetroundjoin%
\definecolor{currentfill}{rgb}{0.000000,0.000000,0.000000}%
\pgfsetfillcolor{currentfill}%
\pgfsetlinewidth{0.803000pt}%
\definecolor{currentstroke}{rgb}{0.000000,0.000000,0.000000}%
\pgfsetstrokecolor{currentstroke}%
\pgfsetdash{}{0pt}%
\pgfsys@defobject{currentmarker}{\pgfqpoint{-0.048611in}{0.000000in}}{\pgfqpoint{0.000000in}{0.000000in}}{%
\pgfpathmoveto{\pgfqpoint{0.000000in}{0.000000in}}%
\pgfpathlineto{\pgfqpoint{-0.048611in}{0.000000in}}%
\pgfusepath{stroke,fill}%
}%
\begin{pgfscope}%
\pgfsys@transformshift{0.316407in}{1.158213in}%
\pgfsys@useobject{currentmarker}{}%
\end{pgfscope}%
\end{pgfscope}%
\begin{pgfscope}%
\pgftext[x=0.149740in,y=1.110385in,left,base]{\rmfamily\fontsize{10.000000}{12.000000}\selectfont \(\displaystyle 2\)}%
\end{pgfscope}%
\begin{pgfscope}%
\pgfsetbuttcap%
\pgfsetroundjoin%
\definecolor{currentfill}{rgb}{0.000000,0.000000,0.000000}%
\pgfsetfillcolor{currentfill}%
\pgfsetlinewidth{0.803000pt}%
\definecolor{currentstroke}{rgb}{0.000000,0.000000,0.000000}%
\pgfsetstrokecolor{currentstroke}%
\pgfsetdash{}{0pt}%
\pgfsys@defobject{currentmarker}{\pgfqpoint{-0.048611in}{0.000000in}}{\pgfqpoint{0.000000in}{0.000000in}}{%
\pgfpathmoveto{\pgfqpoint{0.000000in}{0.000000in}}%
\pgfpathlineto{\pgfqpoint{-0.048611in}{0.000000in}}%
\pgfusepath{stroke,fill}%
}%
\begin{pgfscope}%
\pgfsys@transformshift{0.316407in}{1.544264in}%
\pgfsys@useobject{currentmarker}{}%
\end{pgfscope}%
\end{pgfscope}%
\begin{pgfscope}%
\pgftext[x=0.149740in,y=1.496436in,left,base]{\rmfamily\fontsize{10.000000}{12.000000}\selectfont \(\displaystyle 3\)}%
\end{pgfscope}%
\begin{pgfscope}%
\pgfsetbuttcap%
\pgfsetroundjoin%
\definecolor{currentfill}{rgb}{0.000000,0.000000,0.000000}%
\pgfsetfillcolor{currentfill}%
\pgfsetlinewidth{0.803000pt}%
\definecolor{currentstroke}{rgb}{0.000000,0.000000,0.000000}%
\pgfsetstrokecolor{currentstroke}%
\pgfsetdash{}{0pt}%
\pgfsys@defobject{currentmarker}{\pgfqpoint{-0.048611in}{0.000000in}}{\pgfqpoint{0.000000in}{0.000000in}}{%
\pgfpathmoveto{\pgfqpoint{0.000000in}{0.000000in}}%
\pgfpathlineto{\pgfqpoint{-0.048611in}{0.000000in}}%
\pgfusepath{stroke,fill}%
}%
\begin{pgfscope}%
\pgfsys@transformshift{0.316407in}{1.930315in}%
\pgfsys@useobject{currentmarker}{}%
\end{pgfscope}%
\end{pgfscope}%
\begin{pgfscope}%
\pgftext[x=0.149740in,y=1.882487in,left,base]{\rmfamily\fontsize{10.000000}{12.000000}\selectfont \(\displaystyle 4\)}%
\end{pgfscope}%
\begin{pgfscope}%
\pgfpathrectangle{\pgfqpoint{0.316407in}{0.386111in}}{\pgfqpoint{1.413513in}{1.665785in}}%
\pgfusepath{clip}%
\pgfsetrectcap%
\pgfsetroundjoin%
\pgfsetlinewidth{2.007500pt}%
\definecolor{currentstroke}{rgb}{0.000000,0.000000,0.000000}%
\pgfsetstrokecolor{currentstroke}%
\pgfsetdash{}{0pt}%
\pgfpathmoveto{\pgfqpoint{0.316407in}{0.394974in}}%
\pgfpathlineto{\pgfqpoint{0.340802in}{0.396539in}}%
\pgfpathlineto{\pgfqpoint{0.360594in}{0.400782in}}%
\pgfpathlineto{\pgfqpoint{0.376244in}{0.407265in}}%
\pgfpathlineto{\pgfqpoint{0.389131in}{0.415975in}}%
\pgfpathlineto{\pgfqpoint{0.400178in}{0.427107in}}%
\pgfpathlineto{\pgfqpoint{0.410764in}{0.442342in}}%
\pgfpathlineto{\pgfqpoint{0.420890in}{0.462740in}}%
\pgfpathlineto{\pgfqpoint{0.430556in}{0.489448in}}%
\pgfpathlineto{\pgfqpoint{0.440222in}{0.525594in}}%
\pgfpathlineto{\pgfqpoint{0.449888in}{0.574099in}}%
\pgfpathlineto{\pgfqpoint{0.460014in}{0.641960in}}%
\pgfpathlineto{\pgfqpoint{0.470140in}{0.731601in}}%
\pgfpathlineto{\pgfqpoint{0.481187in}{0.859079in}}%
\pgfpathlineto{\pgfqpoint{0.493154in}{1.035518in}}%
\pgfpathlineto{\pgfqpoint{0.508343in}{1.309807in}}%
\pgfpathlineto{\pgfqpoint{0.532738in}{1.752027in}}%
\pgfpathlineto{\pgfqpoint{0.541023in}{1.849722in}}%
\pgfpathlineto{\pgfqpoint{0.547007in}{1.889609in}}%
\pgfpathlineto{\pgfqpoint{0.551149in}{1.900100in}}%
\pgfpathlineto{\pgfqpoint{0.552990in}{1.900100in}}%
\pgfpathlineto{\pgfqpoint{0.555292in}{1.896054in}}%
\pgfpathlineto{\pgfqpoint{0.558513in}{1.882935in}}%
\pgfpathlineto{\pgfqpoint{0.563116in}{1.849722in}}%
\pgfpathlineto{\pgfqpoint{0.569560in}{1.777464in}}%
\pgfpathlineto{\pgfqpoint{0.578766in}{1.633718in}}%
\pgfpathlineto{\pgfqpoint{0.597637in}{1.274445in}}%
\pgfpathlineto{\pgfqpoint{0.615128in}{0.969934in}}%
\pgfpathlineto{\pgfqpoint{0.628015in}{0.796538in}}%
\pgfpathlineto{\pgfqpoint{0.639983in}{0.675746in}}%
\pgfpathlineto{\pgfqpoint{0.651029in}{0.593600in}}%
\pgfpathlineto{\pgfqpoint{0.662076in}{0.533781in}}%
\pgfpathlineto{\pgfqpoint{0.672662in}{0.492442in}}%
\pgfpathlineto{\pgfqpoint{0.683249in}{0.462740in}}%
\pgfpathlineto{\pgfqpoint{0.693835in}{0.441564in}}%
\pgfpathlineto{\pgfqpoint{0.704882in}{0.426016in}}%
\pgfpathlineto{\pgfqpoint{0.716849in}{0.414487in}}%
\pgfpathlineto{\pgfqpoint{0.729737in}{0.406305in}}%
\pgfpathlineto{\pgfqpoint{0.744926in}{0.400366in}}%
\pgfpathlineto{\pgfqpoint{0.763797in}{0.396478in}}%
\pgfpathlineto{\pgfqpoint{0.786811in}{0.394976in}}%
\pgfpathlineto{\pgfqpoint{0.811206in}{0.396417in}}%
\pgfpathlineto{\pgfqpoint{0.830998in}{0.400502in}}%
\pgfpathlineto{\pgfqpoint{0.846647in}{0.406778in}}%
\pgfpathlineto{\pgfqpoint{0.859535in}{0.415220in}}%
\pgfpathlineto{\pgfqpoint{0.871042in}{0.426557in}}%
\pgfpathlineto{\pgfqpoint{0.881628in}{0.441564in}}%
\pgfpathlineto{\pgfqpoint{0.891754in}{0.461661in}}%
\pgfpathlineto{\pgfqpoint{0.901420in}{0.487982in}}%
\pgfpathlineto{\pgfqpoint{0.911086in}{0.523618in}}%
\pgfpathlineto{\pgfqpoint{0.920752in}{0.571459in}}%
\pgfpathlineto{\pgfqpoint{0.930878in}{0.638436in}}%
\pgfpathlineto{\pgfqpoint{0.941004in}{0.726992in}}%
\pgfpathlineto{\pgfqpoint{0.952051in}{0.853097in}}%
\pgfpathlineto{\pgfqpoint{0.964018in}{1.028004in}}%
\pgfpathlineto{\pgfqpoint{0.979207in}{1.300933in}}%
\pgfpathlineto{\pgfqpoint{1.004062in}{1.752027in}}%
\pgfpathlineto{\pgfqpoint{1.012347in}{1.849722in}}%
\pgfpathlineto{\pgfqpoint{1.018331in}{1.889609in}}%
\pgfpathlineto{\pgfqpoint{1.022473in}{1.900100in}}%
\pgfpathlineto{\pgfqpoint{1.024314in}{1.900100in}}%
\pgfpathlineto{\pgfqpoint{1.026616in}{1.896054in}}%
\pgfpathlineto{\pgfqpoint{1.029838in}{1.882935in}}%
\pgfpathlineto{\pgfqpoint{1.034441in}{1.849722in}}%
\pgfpathlineto{\pgfqpoint{1.040884in}{1.777464in}}%
\pgfpathlineto{\pgfqpoint{1.050090in}{1.633718in}}%
\pgfpathlineto{\pgfqpoint{1.068961in}{1.274445in}}%
\pgfpathlineto{\pgfqpoint{1.086452in}{0.969934in}}%
\pgfpathlineto{\pgfqpoint{1.099340in}{0.796538in}}%
\pgfpathlineto{\pgfqpoint{1.111307in}{0.675746in}}%
\pgfpathlineto{\pgfqpoint{1.122354in}{0.593600in}}%
\pgfpathlineto{\pgfqpoint{1.133400in}{0.533781in}}%
\pgfpathlineto{\pgfqpoint{1.143987in}{0.492442in}}%
\pgfpathlineto{\pgfqpoint{1.154573in}{0.462740in}}%
\pgfpathlineto{\pgfqpoint{1.165159in}{0.441564in}}%
\pgfpathlineto{\pgfqpoint{1.176206in}{0.426016in}}%
\pgfpathlineto{\pgfqpoint{1.188173in}{0.414487in}}%
\pgfpathlineto{\pgfqpoint{1.201061in}{0.406305in}}%
\pgfpathlineto{\pgfqpoint{1.216250in}{0.400366in}}%
\pgfpathlineto{\pgfqpoint{1.235122in}{0.396478in}}%
\pgfpathlineto{\pgfqpoint{1.258135in}{0.394976in}}%
\pgfpathlineto{\pgfqpoint{1.282530in}{0.396417in}}%
\pgfpathlineto{\pgfqpoint{1.302322in}{0.400502in}}%
\pgfpathlineto{\pgfqpoint{1.317972in}{0.406778in}}%
\pgfpathlineto{\pgfqpoint{1.330859in}{0.415220in}}%
\pgfpathlineto{\pgfqpoint{1.342366in}{0.426557in}}%
\pgfpathlineto{\pgfqpoint{1.352953in}{0.441564in}}%
\pgfpathlineto{\pgfqpoint{1.363079in}{0.461661in}}%
\pgfpathlineto{\pgfqpoint{1.372745in}{0.487982in}}%
\pgfpathlineto{\pgfqpoint{1.382410in}{0.523618in}}%
\pgfpathlineto{\pgfqpoint{1.392076in}{0.571459in}}%
\pgfpathlineto{\pgfqpoint{1.402202in}{0.638436in}}%
\pgfpathlineto{\pgfqpoint{1.412328in}{0.726992in}}%
\pgfpathlineto{\pgfqpoint{1.423375in}{0.853097in}}%
\pgfpathlineto{\pgfqpoint{1.435342in}{1.028004in}}%
\pgfpathlineto{\pgfqpoint{1.450532in}{1.300933in}}%
\pgfpathlineto{\pgfqpoint{1.475387in}{1.752027in}}%
\pgfpathlineto{\pgfqpoint{1.483672in}{1.849722in}}%
\pgfpathlineto{\pgfqpoint{1.489655in}{1.889609in}}%
\pgfpathlineto{\pgfqpoint{1.493798in}{1.900100in}}%
\pgfpathlineto{\pgfqpoint{1.495639in}{1.900100in}}%
\pgfpathlineto{\pgfqpoint{1.497940in}{1.896054in}}%
\pgfpathlineto{\pgfqpoint{1.501162in}{1.882935in}}%
\pgfpathlineto{\pgfqpoint{1.505765in}{1.849722in}}%
\pgfpathlineto{\pgfqpoint{1.512209in}{1.777464in}}%
\pgfpathlineto{\pgfqpoint{1.521414in}{1.633718in}}%
\pgfpathlineto{\pgfqpoint{1.540286in}{1.274445in}}%
\pgfpathlineto{\pgfqpoint{1.557776in}{0.969934in}}%
\pgfpathlineto{\pgfqpoint{1.570664in}{0.796538in}}%
\pgfpathlineto{\pgfqpoint{1.582631in}{0.675746in}}%
\pgfpathlineto{\pgfqpoint{1.593678in}{0.593600in}}%
\pgfpathlineto{\pgfqpoint{1.604725in}{0.533781in}}%
\pgfpathlineto{\pgfqpoint{1.615311in}{0.492442in}}%
\pgfpathlineto{\pgfqpoint{1.625897in}{0.462740in}}%
\pgfpathlineto{\pgfqpoint{1.636484in}{0.441564in}}%
\pgfpathlineto{\pgfqpoint{1.647530in}{0.426016in}}%
\pgfpathlineto{\pgfqpoint{1.659498in}{0.414487in}}%
\pgfpathlineto{\pgfqpoint{1.672385in}{0.406305in}}%
\pgfpathlineto{\pgfqpoint{1.687575in}{0.400366in}}%
\pgfpathlineto{\pgfqpoint{1.706446in}{0.396478in}}%
\pgfpathlineto{\pgfqpoint{1.729460in}{0.394976in}}%
\pgfpathlineto{\pgfqpoint{1.729920in}{0.394975in}}%
\pgfpathlineto{\pgfqpoint{1.729920in}{0.394975in}}%
\pgfusepath{stroke}%
\end{pgfscope}%
\begin{pgfscope}%
\pgfpathrectangle{\pgfqpoint{0.316407in}{0.386111in}}{\pgfqpoint{1.413513in}{1.665785in}}%
\pgfusepath{clip}%
\pgfsetbuttcap%
\pgfsetroundjoin%
\pgfsetlinewidth{0.752812pt}%
\definecolor{currentstroke}{rgb}{0.000000,0.000000,0.000000}%
\pgfsetstrokecolor{currentstroke}%
\pgfsetdash{{2.775000pt}{1.200000pt}}{0.000000pt}%
\pgfpathmoveto{\pgfqpoint{0.551839in}{0.386111in}}%
\pgfpathlineto{\pgfqpoint{0.551839in}{2.051896in}}%
\pgfusepath{stroke}%
\end{pgfscope}%
\begin{pgfscope}%
\pgfpathrectangle{\pgfqpoint{0.316407in}{0.386111in}}{\pgfqpoint{1.413513in}{1.665785in}}%
\pgfusepath{clip}%
\pgfsetbuttcap%
\pgfsetroundjoin%
\pgfsetlinewidth{0.752812pt}%
\definecolor{currentstroke}{rgb}{0.000000,0.000000,0.000000}%
\pgfsetstrokecolor{currentstroke}%
\pgfsetdash{{2.775000pt}{1.200000pt}}{0.000000pt}%
\pgfpathmoveto{\pgfqpoint{0.787502in}{0.386111in}}%
\pgfpathlineto{\pgfqpoint{0.787502in}{2.051896in}}%
\pgfusepath{stroke}%
\end{pgfscope}%
\begin{pgfscope}%
\pgfpathrectangle{\pgfqpoint{0.316407in}{0.386111in}}{\pgfqpoint{1.413513in}{1.665785in}}%
\pgfusepath{clip}%
\pgfsetbuttcap%
\pgfsetroundjoin%
\pgfsetlinewidth{0.752812pt}%
\definecolor{currentstroke}{rgb}{0.000000,0.000000,0.000000}%
\pgfsetstrokecolor{currentstroke}%
\pgfsetdash{{2.775000pt}{1.200000pt}}{0.000000pt}%
\pgfpathmoveto{\pgfqpoint{1.023164in}{0.386111in}}%
\pgfpathlineto{\pgfqpoint{1.023164in}{2.051896in}}%
\pgfusepath{stroke}%
\end{pgfscope}%
\begin{pgfscope}%
\pgfpathrectangle{\pgfqpoint{0.316407in}{0.386111in}}{\pgfqpoint{1.413513in}{1.665785in}}%
\pgfusepath{clip}%
\pgfsetbuttcap%
\pgfsetroundjoin%
\pgfsetlinewidth{0.752812pt}%
\definecolor{currentstroke}{rgb}{0.000000,0.000000,0.000000}%
\pgfsetstrokecolor{currentstroke}%
\pgfsetdash{{2.775000pt}{1.200000pt}}{0.000000pt}%
\pgfpathmoveto{\pgfqpoint{1.258826in}{0.386111in}}%
\pgfpathlineto{\pgfqpoint{1.258826in}{2.051896in}}%
\pgfusepath{stroke}%
\end{pgfscope}%
\begin{pgfscope}%
\pgfpathrectangle{\pgfqpoint{0.316407in}{0.386111in}}{\pgfqpoint{1.413513in}{1.665785in}}%
\pgfusepath{clip}%
\pgfsetbuttcap%
\pgfsetroundjoin%
\pgfsetlinewidth{0.752812pt}%
\definecolor{currentstroke}{rgb}{0.000000,0.000000,0.000000}%
\pgfsetstrokecolor{currentstroke}%
\pgfsetdash{{2.775000pt}{1.200000pt}}{0.000000pt}%
\pgfpathmoveto{\pgfqpoint{1.494488in}{0.386111in}}%
\pgfpathlineto{\pgfqpoint{1.494488in}{2.051896in}}%
\pgfusepath{stroke}%
\end{pgfscope}%
\begin{pgfscope}%
\pgfpathrectangle{\pgfqpoint{0.316407in}{0.386111in}}{\pgfqpoint{1.413513in}{1.665785in}}%
\pgfusepath{clip}%
\pgfsetbuttcap%
\pgfsetroundjoin%
\pgfsetlinewidth{0.752812pt}%
\definecolor{currentstroke}{rgb}{0.000000,0.000000,0.000000}%
\pgfsetstrokecolor{currentstroke}%
\pgfsetdash{{2.775000pt}{1.200000pt}}{0.000000pt}%
\pgfpathmoveto{\pgfqpoint{0.316407in}{0.386111in}}%
\pgfpathlineto{\pgfqpoint{0.316407in}{2.051896in}}%
\pgfusepath{stroke}%
\end{pgfscope}%
\begin{pgfscope}%
\pgfsetrectcap%
\pgfsetmiterjoin%
\pgfsetlinewidth{0.803000pt}%
\definecolor{currentstroke}{rgb}{0.000000,0.000000,0.000000}%
\pgfsetstrokecolor{currentstroke}%
\pgfsetdash{}{0pt}%
\pgfpathmoveto{\pgfqpoint{0.316407in}{0.386111in}}%
\pgfpathlineto{\pgfqpoint{0.316407in}{2.051896in}}%
\pgfusepath{stroke}%
\end{pgfscope}%
\begin{pgfscope}%
\pgfsetrectcap%
\pgfsetmiterjoin%
\pgfsetlinewidth{0.803000pt}%
\definecolor{currentstroke}{rgb}{0.000000,0.000000,0.000000}%
\pgfsetstrokecolor{currentstroke}%
\pgfsetdash{}{0pt}%
\pgfpathmoveto{\pgfqpoint{1.729920in}{0.386111in}}%
\pgfpathlineto{\pgfqpoint{1.729920in}{2.051896in}}%
\pgfusepath{stroke}%
\end{pgfscope}%
\begin{pgfscope}%
\pgfsetrectcap%
\pgfsetmiterjoin%
\pgfsetlinewidth{0.803000pt}%
\definecolor{currentstroke}{rgb}{0.000000,0.000000,0.000000}%
\pgfsetstrokecolor{currentstroke}%
\pgfsetdash{}{0pt}%
\pgfpathmoveto{\pgfqpoint{0.316407in}{0.386111in}}%
\pgfpathlineto{\pgfqpoint{1.729920in}{0.386111in}}%
\pgfusepath{stroke}%
\end{pgfscope}%
\begin{pgfscope}%
\pgfsetrectcap%
\pgfsetmiterjoin%
\pgfsetlinewidth{0.803000pt}%
\definecolor{currentstroke}{rgb}{0.000000,0.000000,0.000000}%
\pgfsetstrokecolor{currentstroke}%
\pgfsetdash{}{0pt}%
\pgfpathmoveto{\pgfqpoint{0.316407in}{2.051896in}}%
\pgfpathlineto{\pgfqpoint{1.729920in}{2.051896in}}%
\pgfusepath{stroke}%
\end{pgfscope}%
\begin{pgfscope}%
\pgftext[x=1.023164in,y=2.135229in,,base]{\rmfamily\fontsize{12.000000}{14.400000}\selectfont \(\displaystyle u(x)\)}%
\end{pgfscope}%
\begin{pgfscope}%
\pgfsetbuttcap%
\pgfsetmiterjoin%
\definecolor{currentfill}{rgb}{1.000000,1.000000,1.000000}%
\pgfsetfillcolor{currentfill}%
\pgfsetlinewidth{0.000000pt}%
\definecolor{currentstroke}{rgb}{0.000000,0.000000,0.000000}%
\pgfsetstrokecolor{currentstroke}%
\pgfsetstrokeopacity{0.000000}%
\pgfsetdash{}{0pt}%
\pgfpathmoveto{\pgfqpoint{2.310181in}{0.386111in}}%
\pgfpathlineto{\pgfqpoint{3.723694in}{0.386111in}}%
\pgfpathlineto{\pgfqpoint{3.723694in}{2.051896in}}%
\pgfpathlineto{\pgfqpoint{2.310181in}{2.051896in}}%
\pgfpathclose%
\pgfusepath{fill}%
\end{pgfscope}%
\begin{pgfscope}%
\pgfsetbuttcap%
\pgfsetroundjoin%
\definecolor{currentfill}{rgb}{0.000000,0.000000,0.000000}%
\pgfsetfillcolor{currentfill}%
\pgfsetlinewidth{0.803000pt}%
\definecolor{currentstroke}{rgb}{0.000000,0.000000,0.000000}%
\pgfsetstrokecolor{currentstroke}%
\pgfsetdash{}{0pt}%
\pgfsys@defobject{currentmarker}{\pgfqpoint{0.000000in}{-0.048611in}}{\pgfqpoint{0.000000in}{0.000000in}}{%
\pgfpathmoveto{\pgfqpoint{0.000000in}{0.000000in}}%
\pgfpathlineto{\pgfqpoint{0.000000in}{-0.048611in}}%
\pgfusepath{stroke,fill}%
}%
\begin{pgfscope}%
\pgfsys@transformshift{2.310181in}{0.386111in}%
\pgfsys@useobject{currentmarker}{}%
\end{pgfscope}%
\end{pgfscope}%
\begin{pgfscope}%
\pgftext[x=2.310181in,y=0.288889in,,top]{\rmfamily\fontsize{10.000000}{12.000000}\selectfont \(\displaystyle 0\)}%
\end{pgfscope}%
\begin{pgfscope}%
\pgfsetbuttcap%
\pgfsetroundjoin%
\definecolor{currentfill}{rgb}{0.000000,0.000000,0.000000}%
\pgfsetfillcolor{currentfill}%
\pgfsetlinewidth{0.803000pt}%
\definecolor{currentstroke}{rgb}{0.000000,0.000000,0.000000}%
\pgfsetstrokecolor{currentstroke}%
\pgfsetdash{}{0pt}%
\pgfsys@defobject{currentmarker}{\pgfqpoint{0.000000in}{-0.048611in}}{\pgfqpoint{0.000000in}{0.000000in}}{%
\pgfpathmoveto{\pgfqpoint{0.000000in}{0.000000in}}%
\pgfpathlineto{\pgfqpoint{0.000000in}{-0.048611in}}%
\pgfusepath{stroke,fill}%
}%
\begin{pgfscope}%
\pgfsys@transformshift{3.016938in}{0.386111in}%
\pgfsys@useobject{currentmarker}{}%
\end{pgfscope}%
\end{pgfscope}%
\begin{pgfscope}%
\pgftext[x=3.016938in,y=0.288889in,,top]{\rmfamily\fontsize{10.000000}{12.000000}\selectfont \(\displaystyle L/2\)}%
\end{pgfscope}%
\begin{pgfscope}%
\pgfsetbuttcap%
\pgfsetroundjoin%
\definecolor{currentfill}{rgb}{0.000000,0.000000,0.000000}%
\pgfsetfillcolor{currentfill}%
\pgfsetlinewidth{0.803000pt}%
\definecolor{currentstroke}{rgb}{0.000000,0.000000,0.000000}%
\pgfsetstrokecolor{currentstroke}%
\pgfsetdash{}{0pt}%
\pgfsys@defobject{currentmarker}{\pgfqpoint{0.000000in}{-0.048611in}}{\pgfqpoint{0.000000in}{0.000000in}}{%
\pgfpathmoveto{\pgfqpoint{0.000000in}{0.000000in}}%
\pgfpathlineto{\pgfqpoint{0.000000in}{-0.048611in}}%
\pgfusepath{stroke,fill}%
}%
\begin{pgfscope}%
\pgfsys@transformshift{3.723694in}{0.386111in}%
\pgfsys@useobject{currentmarker}{}%
\end{pgfscope}%
\end{pgfscope}%
\begin{pgfscope}%
\pgftext[x=3.723694in,y=0.288889in,,top]{\rmfamily\fontsize{10.000000}{12.000000}\selectfont \(\displaystyle L\)}%
\end{pgfscope}%
\begin{pgfscope}%
\pgfsetbuttcap%
\pgfsetroundjoin%
\definecolor{currentfill}{rgb}{0.000000,0.000000,0.000000}%
\pgfsetfillcolor{currentfill}%
\pgfsetlinewidth{0.803000pt}%
\definecolor{currentstroke}{rgb}{0.000000,0.000000,0.000000}%
\pgfsetstrokecolor{currentstroke}%
\pgfsetdash{}{0pt}%
\pgfsys@defobject{currentmarker}{\pgfqpoint{-0.048611in}{0.000000in}}{\pgfqpoint{0.000000in}{0.000000in}}{%
\pgfpathmoveto{\pgfqpoint{0.000000in}{0.000000in}}%
\pgfpathlineto{\pgfqpoint{-0.048611in}{0.000000in}}%
\pgfusepath{stroke,fill}%
}%
\begin{pgfscope}%
\pgfsys@transformshift{2.310181in}{0.499980in}%
\pgfsys@useobject{currentmarker}{}%
\end{pgfscope}%
\end{pgfscope}%
\begin{pgfscope}%
\pgftext[x=1.858020in,y=0.452152in,left,base]{\rmfamily\fontsize{10.000000}{12.000000}\selectfont \(\displaystyle -0.50\)}%
\end{pgfscope}%
\begin{pgfscope}%
\pgfsetbuttcap%
\pgfsetroundjoin%
\definecolor{currentfill}{rgb}{0.000000,0.000000,0.000000}%
\pgfsetfillcolor{currentfill}%
\pgfsetlinewidth{0.803000pt}%
\definecolor{currentstroke}{rgb}{0.000000,0.000000,0.000000}%
\pgfsetstrokecolor{currentstroke}%
\pgfsetdash{}{0pt}%
\pgfsys@defobject{currentmarker}{\pgfqpoint{-0.048611in}{0.000000in}}{\pgfqpoint{0.000000in}{0.000000in}}{%
\pgfpathmoveto{\pgfqpoint{0.000000in}{0.000000in}}%
\pgfpathlineto{\pgfqpoint{-0.048611in}{0.000000in}}%
\pgfusepath{stroke,fill}%
}%
\begin{pgfscope}%
\pgfsys@transformshift{2.310181in}{0.858790in}%
\pgfsys@useobject{currentmarker}{}%
\end{pgfscope}%
\end{pgfscope}%
\begin{pgfscope}%
\pgftext[x=1.858020in,y=0.810962in,left,base]{\rmfamily\fontsize{10.000000}{12.000000}\selectfont \(\displaystyle -0.25\)}%
\end{pgfscope}%
\begin{pgfscope}%
\pgfsetbuttcap%
\pgfsetroundjoin%
\definecolor{currentfill}{rgb}{0.000000,0.000000,0.000000}%
\pgfsetfillcolor{currentfill}%
\pgfsetlinewidth{0.803000pt}%
\definecolor{currentstroke}{rgb}{0.000000,0.000000,0.000000}%
\pgfsetstrokecolor{currentstroke}%
\pgfsetdash{}{0pt}%
\pgfsys@defobject{currentmarker}{\pgfqpoint{-0.048611in}{0.000000in}}{\pgfqpoint{0.000000in}{0.000000in}}{%
\pgfpathmoveto{\pgfqpoint{0.000000in}{0.000000in}}%
\pgfpathlineto{\pgfqpoint{-0.048611in}{0.000000in}}%
\pgfusepath{stroke,fill}%
}%
\begin{pgfscope}%
\pgfsys@transformshift{2.310181in}{1.217601in}%
\pgfsys@useobject{currentmarker}{}%
\end{pgfscope}%
\end{pgfscope}%
\begin{pgfscope}%
\pgftext[x=1.966045in,y=1.169773in,left,base]{\rmfamily\fontsize{10.000000}{12.000000}\selectfont \(\displaystyle 0.00\)}%
\end{pgfscope}%
\begin{pgfscope}%
\pgfsetbuttcap%
\pgfsetroundjoin%
\definecolor{currentfill}{rgb}{0.000000,0.000000,0.000000}%
\pgfsetfillcolor{currentfill}%
\pgfsetlinewidth{0.803000pt}%
\definecolor{currentstroke}{rgb}{0.000000,0.000000,0.000000}%
\pgfsetstrokecolor{currentstroke}%
\pgfsetdash{}{0pt}%
\pgfsys@defobject{currentmarker}{\pgfqpoint{-0.048611in}{0.000000in}}{\pgfqpoint{0.000000in}{0.000000in}}{%
\pgfpathmoveto{\pgfqpoint{0.000000in}{0.000000in}}%
\pgfpathlineto{\pgfqpoint{-0.048611in}{0.000000in}}%
\pgfusepath{stroke,fill}%
}%
\begin{pgfscope}%
\pgfsys@transformshift{2.310181in}{1.576412in}%
\pgfsys@useobject{currentmarker}{}%
\end{pgfscope}%
\end{pgfscope}%
\begin{pgfscope}%
\pgftext[x=1.966045in,y=1.528584in,left,base]{\rmfamily\fontsize{10.000000}{12.000000}\selectfont \(\displaystyle 0.25\)}%
\end{pgfscope}%
\begin{pgfscope}%
\pgfsetbuttcap%
\pgfsetroundjoin%
\definecolor{currentfill}{rgb}{0.000000,0.000000,0.000000}%
\pgfsetfillcolor{currentfill}%
\pgfsetlinewidth{0.803000pt}%
\definecolor{currentstroke}{rgb}{0.000000,0.000000,0.000000}%
\pgfsetstrokecolor{currentstroke}%
\pgfsetdash{}{0pt}%
\pgfsys@defobject{currentmarker}{\pgfqpoint{-0.048611in}{0.000000in}}{\pgfqpoint{0.000000in}{0.000000in}}{%
\pgfpathmoveto{\pgfqpoint{0.000000in}{0.000000in}}%
\pgfpathlineto{\pgfqpoint{-0.048611in}{0.000000in}}%
\pgfusepath{stroke,fill}%
}%
\begin{pgfscope}%
\pgfsys@transformshift{2.310181in}{1.935222in}%
\pgfsys@useobject{currentmarker}{}%
\end{pgfscope}%
\end{pgfscope}%
\begin{pgfscope}%
\pgftext[x=1.966045in,y=1.887394in,left,base]{\rmfamily\fontsize{10.000000}{12.000000}\selectfont \(\displaystyle 0.50\)}%
\end{pgfscope}%
\begin{pgfscope}%
\pgfpathrectangle{\pgfqpoint{2.310181in}{0.386111in}}{\pgfqpoint{1.413513in}{1.665785in}}%
\pgfusepath{clip}%
\pgfsetrectcap%
\pgfsetroundjoin%
\pgfsetlinewidth{2.007500pt}%
\definecolor{currentstroke}{rgb}{0.000000,0.000000,0.000000}%
\pgfsetstrokecolor{currentstroke}%
\pgfsetdash{}{0pt}%
\pgfpathmoveto{\pgfqpoint{2.310181in}{1.217601in}}%
\pgfpathlineto{\pgfqpoint{2.333656in}{1.608194in}}%
\pgfpathlineto{\pgfqpoint{2.346083in}{1.766245in}}%
\pgfpathlineto{\pgfqpoint{2.356670in}{1.864267in}}%
\pgfpathlineto{\pgfqpoint{2.365875in}{1.922318in}}%
\pgfpathlineto{\pgfqpoint{2.373700in}{1.953521in}}%
\pgfpathlineto{\pgfqpoint{2.380144in}{1.968154in}}%
\pgfpathlineto{\pgfqpoint{2.385207in}{1.973469in}}%
\pgfpathlineto{\pgfqpoint{2.389810in}{1.974107in}}%
\pgfpathlineto{\pgfqpoint{2.394412in}{1.971178in}}%
\pgfpathlineto{\pgfqpoint{2.399936in}{1.963502in}}%
\pgfpathlineto{\pgfqpoint{2.407300in}{1.947217in}}%
\pgfpathlineto{\pgfqpoint{2.416966in}{1.917320in}}%
\pgfpathlineto{\pgfqpoint{2.429854in}{1.866179in}}%
\pgfpathlineto{\pgfqpoint{2.447344in}{1.782971in}}%
\pgfpathlineto{\pgfqpoint{2.473580in}{1.641480in}}%
\pgfpathlineto{\pgfqpoint{2.520528in}{1.368485in}}%
\pgfpathlineto{\pgfqpoint{2.630535in}{0.725905in}}%
\pgfpathlineto{\pgfqpoint{2.655850in}{0.596973in}}%
\pgfpathlineto{\pgfqpoint{2.672420in}{0.527116in}}%
\pgfpathlineto{\pgfqpoint{2.684387in}{0.488919in}}%
\pgfpathlineto{\pgfqpoint{2.693132in}{0.470322in}}%
\pgfpathlineto{\pgfqpoint{2.699576in}{0.463074in}}%
\pgfpathlineto{\pgfqpoint{2.704179in}{0.461878in}}%
\pgfpathlineto{\pgfqpoint{2.708322in}{0.464014in}}%
\pgfpathlineto{\pgfqpoint{2.712924in}{0.470343in}}%
\pgfpathlineto{\pgfqpoint{2.718448in}{0.484026in}}%
\pgfpathlineto{\pgfqpoint{2.724892in}{0.509323in}}%
\pgfpathlineto{\pgfqpoint{2.732256in}{0.551884in}}%
\pgfpathlineto{\pgfqpoint{2.740541in}{0.618727in}}%
\pgfpathlineto{\pgfqpoint{2.750207in}{0.723193in}}%
\pgfpathlineto{\pgfqpoint{2.762174in}{0.889943in}}%
\pgfpathlineto{\pgfqpoint{2.779204in}{1.177599in}}%
\pgfpathlineto{\pgfqpoint{2.804980in}{1.609129in}}%
\pgfpathlineto{\pgfqpoint{2.817407in}{1.767180in}}%
\pgfpathlineto{\pgfqpoint{2.827994in}{1.865202in}}%
\pgfpathlineto{\pgfqpoint{2.837199in}{1.923253in}}%
\pgfpathlineto{\pgfqpoint{2.845024in}{1.954457in}}%
\pgfpathlineto{\pgfqpoint{2.851468in}{1.969089in}}%
\pgfpathlineto{\pgfqpoint{2.856531in}{1.974404in}}%
\pgfpathlineto{\pgfqpoint{2.861134in}{1.975042in}}%
\pgfpathlineto{\pgfqpoint{2.865737in}{1.972113in}}%
\pgfpathlineto{\pgfqpoint{2.871260in}{1.964437in}}%
\pgfpathlineto{\pgfqpoint{2.878624in}{1.948152in}}%
\pgfpathlineto{\pgfqpoint{2.888290in}{1.918255in}}%
\pgfpathlineto{\pgfqpoint{2.901178in}{1.867114in}}%
\pgfpathlineto{\pgfqpoint{2.918669in}{1.783906in}}%
\pgfpathlineto{\pgfqpoint{2.944904in}{1.642415in}}%
\pgfpathlineto{\pgfqpoint{2.991853in}{1.369420in}}%
\pgfpathlineto{\pgfqpoint{3.101859in}{0.726840in}}%
\pgfpathlineto{\pgfqpoint{3.127174in}{0.597908in}}%
\pgfpathlineto{\pgfqpoint{3.143744in}{0.528051in}}%
\pgfpathlineto{\pgfqpoint{3.155712in}{0.489854in}}%
\pgfpathlineto{\pgfqpoint{3.164457in}{0.471257in}}%
\pgfpathlineto{\pgfqpoint{3.170901in}{0.464009in}}%
\pgfpathlineto{\pgfqpoint{3.175503in}{0.462813in}}%
\pgfpathlineto{\pgfqpoint{3.179646in}{0.464949in}}%
\pgfpathlineto{\pgfqpoint{3.184249in}{0.471278in}}%
\pgfpathlineto{\pgfqpoint{3.189772in}{0.484961in}}%
\pgfpathlineto{\pgfqpoint{3.196216in}{0.510258in}}%
\pgfpathlineto{\pgfqpoint{3.203580in}{0.552819in}}%
\pgfpathlineto{\pgfqpoint{3.211865in}{0.619662in}}%
\pgfpathlineto{\pgfqpoint{3.221531in}{0.724128in}}%
\pgfpathlineto{\pgfqpoint{3.233498in}{0.890878in}}%
\pgfpathlineto{\pgfqpoint{3.250529in}{1.178534in}}%
\pgfpathlineto{\pgfqpoint{3.276304in}{1.610064in}}%
\pgfpathlineto{\pgfqpoint{3.288732in}{1.768116in}}%
\pgfpathlineto{\pgfqpoint{3.299318in}{1.866137in}}%
\pgfpathlineto{\pgfqpoint{3.308524in}{1.924188in}}%
\pgfpathlineto{\pgfqpoint{3.316348in}{1.955392in}}%
\pgfpathlineto{\pgfqpoint{3.322792in}{1.970024in}}%
\pgfpathlineto{\pgfqpoint{3.327855in}{1.975339in}}%
\pgfpathlineto{\pgfqpoint{3.332458in}{1.975977in}}%
\pgfpathlineto{\pgfqpoint{3.337061in}{1.973048in}}%
\pgfpathlineto{\pgfqpoint{3.342584in}{1.965372in}}%
\pgfpathlineto{\pgfqpoint{3.349949in}{1.949087in}}%
\pgfpathlineto{\pgfqpoint{3.359615in}{1.919190in}}%
\pgfpathlineto{\pgfqpoint{3.372502in}{1.868049in}}%
\pgfpathlineto{\pgfqpoint{3.389993in}{1.784841in}}%
\pgfpathlineto{\pgfqpoint{3.416229in}{1.643350in}}%
\pgfpathlineto{\pgfqpoint{3.463177in}{1.370355in}}%
\pgfpathlineto{\pgfqpoint{3.573183in}{0.727775in}}%
\pgfpathlineto{\pgfqpoint{3.598499in}{0.598843in}}%
\pgfpathlineto{\pgfqpoint{3.615069in}{0.528986in}}%
\pgfpathlineto{\pgfqpoint{3.627036in}{0.490789in}}%
\pgfpathlineto{\pgfqpoint{3.635781in}{0.472192in}}%
\pgfpathlineto{\pgfqpoint{3.642225in}{0.464944in}}%
\pgfpathlineto{\pgfqpoint{3.646828in}{0.463748in}}%
\pgfpathlineto{\pgfqpoint{3.650970in}{0.465884in}}%
\pgfpathlineto{\pgfqpoint{3.655573in}{0.472213in}}%
\pgfpathlineto{\pgfqpoint{3.661096in}{0.485896in}}%
\pgfpathlineto{\pgfqpoint{3.667540in}{0.511193in}}%
\pgfpathlineto{\pgfqpoint{3.674905in}{0.553754in}}%
\pgfpathlineto{\pgfqpoint{3.683190in}{0.620597in}}%
\pgfpathlineto{\pgfqpoint{3.692856in}{0.725063in}}%
\pgfpathlineto{\pgfqpoint{3.704823in}{0.891813in}}%
\pgfpathlineto{\pgfqpoint{3.721853in}{1.179469in}}%
\pgfpathlineto{\pgfqpoint{3.723694in}{1.212213in}}%
\pgfpathlineto{\pgfqpoint{3.723694in}{1.212213in}}%
\pgfusepath{stroke}%
\end{pgfscope}%
\begin{pgfscope}%
\pgfpathrectangle{\pgfqpoint{2.310181in}{0.386111in}}{\pgfqpoint{1.413513in}{1.665785in}}%
\pgfusepath{clip}%
\pgfsetbuttcap%
\pgfsetroundjoin%
\pgfsetlinewidth{0.752812pt}%
\definecolor{currentstroke}{rgb}{0.000000,0.000000,0.000000}%
\pgfsetstrokecolor{currentstroke}%
\pgfsetdash{{2.775000pt}{1.200000pt}}{0.000000pt}%
\pgfpathmoveto{\pgfqpoint{2.545614in}{0.386111in}}%
\pgfpathlineto{\pgfqpoint{2.545614in}{2.051896in}}%
\pgfusepath{stroke}%
\end{pgfscope}%
\begin{pgfscope}%
\pgfpathrectangle{\pgfqpoint{2.310181in}{0.386111in}}{\pgfqpoint{1.413513in}{1.665785in}}%
\pgfusepath{clip}%
\pgfsetbuttcap%
\pgfsetroundjoin%
\pgfsetlinewidth{0.752812pt}%
\definecolor{currentstroke}{rgb}{0.000000,0.000000,0.000000}%
\pgfsetstrokecolor{currentstroke}%
\pgfsetdash{{2.775000pt}{1.200000pt}}{0.000000pt}%
\pgfpathmoveto{\pgfqpoint{2.781276in}{0.386111in}}%
\pgfpathlineto{\pgfqpoint{2.781276in}{2.051896in}}%
\pgfusepath{stroke}%
\end{pgfscope}%
\begin{pgfscope}%
\pgfpathrectangle{\pgfqpoint{2.310181in}{0.386111in}}{\pgfqpoint{1.413513in}{1.665785in}}%
\pgfusepath{clip}%
\pgfsetbuttcap%
\pgfsetroundjoin%
\pgfsetlinewidth{0.752812pt}%
\definecolor{currentstroke}{rgb}{0.000000,0.000000,0.000000}%
\pgfsetstrokecolor{currentstroke}%
\pgfsetdash{{2.775000pt}{1.200000pt}}{0.000000pt}%
\pgfpathmoveto{\pgfqpoint{3.016938in}{0.386111in}}%
\pgfpathlineto{\pgfqpoint{3.016938in}{2.051896in}}%
\pgfusepath{stroke}%
\end{pgfscope}%
\begin{pgfscope}%
\pgfpathrectangle{\pgfqpoint{2.310181in}{0.386111in}}{\pgfqpoint{1.413513in}{1.665785in}}%
\pgfusepath{clip}%
\pgfsetbuttcap%
\pgfsetroundjoin%
\pgfsetlinewidth{0.752812pt}%
\definecolor{currentstroke}{rgb}{0.000000,0.000000,0.000000}%
\pgfsetstrokecolor{currentstroke}%
\pgfsetdash{{2.775000pt}{1.200000pt}}{0.000000pt}%
\pgfpathmoveto{\pgfqpoint{3.252600in}{0.386111in}}%
\pgfpathlineto{\pgfqpoint{3.252600in}{2.051896in}}%
\pgfusepath{stroke}%
\end{pgfscope}%
\begin{pgfscope}%
\pgfpathrectangle{\pgfqpoint{2.310181in}{0.386111in}}{\pgfqpoint{1.413513in}{1.665785in}}%
\pgfusepath{clip}%
\pgfsetbuttcap%
\pgfsetroundjoin%
\pgfsetlinewidth{0.752812pt}%
\definecolor{currentstroke}{rgb}{0.000000,0.000000,0.000000}%
\pgfsetstrokecolor{currentstroke}%
\pgfsetdash{{2.775000pt}{1.200000pt}}{0.000000pt}%
\pgfpathmoveto{\pgfqpoint{3.488262in}{0.386111in}}%
\pgfpathlineto{\pgfqpoint{3.488262in}{2.051896in}}%
\pgfusepath{stroke}%
\end{pgfscope}%
\begin{pgfscope}%
\pgfpathrectangle{\pgfqpoint{2.310181in}{0.386111in}}{\pgfqpoint{1.413513in}{1.665785in}}%
\pgfusepath{clip}%
\pgfsetbuttcap%
\pgfsetroundjoin%
\pgfsetlinewidth{0.752812pt}%
\definecolor{currentstroke}{rgb}{0.000000,0.000000,0.000000}%
\pgfsetstrokecolor{currentstroke}%
\pgfsetdash{{2.775000pt}{1.200000pt}}{0.000000pt}%
\pgfpathmoveto{\pgfqpoint{2.310181in}{1.217601in}}%
\pgfpathlineto{\pgfqpoint{3.723694in}{1.217601in}}%
\pgfusepath{stroke}%
\end{pgfscope}%
\begin{pgfscope}%
\pgfsetrectcap%
\pgfsetmiterjoin%
\pgfsetlinewidth{0.803000pt}%
\definecolor{currentstroke}{rgb}{0.000000,0.000000,0.000000}%
\pgfsetstrokecolor{currentstroke}%
\pgfsetdash{}{0pt}%
\pgfpathmoveto{\pgfqpoint{2.310181in}{0.386111in}}%
\pgfpathlineto{\pgfqpoint{2.310181in}{2.051896in}}%
\pgfusepath{stroke}%
\end{pgfscope}%
\begin{pgfscope}%
\pgfsetrectcap%
\pgfsetmiterjoin%
\pgfsetlinewidth{0.803000pt}%
\definecolor{currentstroke}{rgb}{0.000000,0.000000,0.000000}%
\pgfsetstrokecolor{currentstroke}%
\pgfsetdash{}{0pt}%
\pgfpathmoveto{\pgfqpoint{3.723694in}{0.386111in}}%
\pgfpathlineto{\pgfqpoint{3.723694in}{2.051896in}}%
\pgfusepath{stroke}%
\end{pgfscope}%
\begin{pgfscope}%
\pgfsetrectcap%
\pgfsetmiterjoin%
\pgfsetlinewidth{0.803000pt}%
\definecolor{currentstroke}{rgb}{0.000000,0.000000,0.000000}%
\pgfsetstrokecolor{currentstroke}%
\pgfsetdash{}{0pt}%
\pgfpathmoveto{\pgfqpoint{2.310181in}{0.386111in}}%
\pgfpathlineto{\pgfqpoint{3.723694in}{0.386111in}}%
\pgfusepath{stroke}%
\end{pgfscope}%
\begin{pgfscope}%
\pgfsetrectcap%
\pgfsetmiterjoin%
\pgfsetlinewidth{0.803000pt}%
\definecolor{currentstroke}{rgb}{0.000000,0.000000,0.000000}%
\pgfsetstrokecolor{currentstroke}%
\pgfsetdash{}{0pt}%
\pgfpathmoveto{\pgfqpoint{2.310181in}{2.051896in}}%
\pgfpathlineto{\pgfqpoint{3.723694in}{2.051896in}}%
\pgfusepath{stroke}%
\end{pgfscope}%
\begin{pgfscope}%
\pgftext[x=3.016938in,y=2.135229in,,base]{\rmfamily\fontsize{12.000000}{14.400000}\selectfont \(\displaystyle w(x)\)}%
\end{pgfscope}%
\end{pgfpicture}%
\makeatother%
\endgroup%

%% file: BifurcationTypeDetermination.pgf
\begingroup%
\makeatletter%
\begin{pgfpicture}%
\pgfpathrectangle{\pgfpointorigin}{\pgfqpoint{3.922305in}{2.424118in}}%
\pgfusepath{use as bounding box, clip}%
\begin{pgfscope}%
\pgfsetbuttcap%
\pgfsetmiterjoin%
\definecolor{currentfill}{rgb}{1.000000,1.000000,1.000000}%
\pgfsetfillcolor{currentfill}%
\pgfsetlinewidth{0.000000pt}%
\definecolor{currentstroke}{rgb}{1.000000,1.000000,1.000000}%
\pgfsetstrokecolor{currentstroke}%
\pgfsetdash{}{0pt}%
\pgfpathmoveto{\pgfqpoint{0.000000in}{0.000000in}}%
\pgfpathlineto{\pgfqpoint{3.922305in}{0.000000in}}%
\pgfpathlineto{\pgfqpoint{3.922305in}{2.424118in}}%
\pgfpathlineto{\pgfqpoint{0.000000in}{2.424118in}}%
\pgfpathclose%
\pgfusepath{fill}%
\end{pgfscope}%
\begin{pgfscope}%
\pgfsetbuttcap%
\pgfsetmiterjoin%
\definecolor{currentfill}{rgb}{1.000000,1.000000,1.000000}%
\pgfsetfillcolor{currentfill}%
\pgfsetlinewidth{0.000000pt}%
\definecolor{currentstroke}{rgb}{0.000000,0.000000,0.000000}%
\pgfsetstrokecolor{currentstroke}%
\pgfsetstrokeopacity{0.000000}%
\pgfsetdash{}{0pt}%
\pgfpathmoveto{\pgfqpoint{0.618473in}{0.580556in}}%
\pgfpathlineto{\pgfqpoint{1.807889in}{0.580556in}}%
\pgfpathlineto{\pgfqpoint{1.807889in}{2.218562in}}%
\pgfpathlineto{\pgfqpoint{0.618473in}{2.218562in}}%
\pgfpathclose%
\pgfusepath{fill}%
\end{pgfscope}%
\begin{pgfscope}%
\pgfsetbuttcap%
\pgfsetroundjoin%
\definecolor{currentfill}{rgb}{0.000000,0.000000,0.000000}%
\pgfsetfillcolor{currentfill}%
\pgfsetlinewidth{0.803000pt}%
\definecolor{currentstroke}{rgb}{0.000000,0.000000,0.000000}%
\pgfsetstrokecolor{currentstroke}%
\pgfsetdash{}{0pt}%
\pgfsys@defobject{currentmarker}{\pgfqpoint{0.000000in}{-0.048611in}}{\pgfqpoint{0.000000in}{0.000000in}}{%
\pgfpathmoveto{\pgfqpoint{0.000000in}{0.000000in}}%
\pgfpathlineto{\pgfqpoint{0.000000in}{-0.048611in}}%
\pgfusepath{stroke,fill}%
}%
\begin{pgfscope}%
\pgfsys@transformshift{0.618473in}{0.580556in}%
\pgfsys@useobject{currentmarker}{}%
\end{pgfscope}%
\end{pgfscope}%
\begin{pgfscope}%
\pgftext[x=0.618473in,y=0.483333in,,top]{\rmfamily\fontsize{10.000000}{12.000000}\selectfont \(\displaystyle 0\)}%
\end{pgfscope}%
\begin{pgfscope}%
\pgfsetbuttcap%
\pgfsetroundjoin%
\definecolor{currentfill}{rgb}{0.000000,0.000000,0.000000}%
\pgfsetfillcolor{currentfill}%
\pgfsetlinewidth{0.803000pt}%
\definecolor{currentstroke}{rgb}{0.000000,0.000000,0.000000}%
\pgfsetstrokecolor{currentstroke}%
\pgfsetdash{}{0pt}%
\pgfsys@defobject{currentmarker}{\pgfqpoint{0.000000in}{-0.048611in}}{\pgfqpoint{0.000000in}{0.000000in}}{%
\pgfpathmoveto{\pgfqpoint{0.000000in}{0.000000in}}%
\pgfpathlineto{\pgfqpoint{0.000000in}{-0.048611in}}%
\pgfusepath{stroke,fill}%
}%
\begin{pgfscope}%
\pgfsys@transformshift{0.888795in}{0.580556in}%
\pgfsys@useobject{currentmarker}{}%
\end{pgfscope}%
\end{pgfscope}%
\begin{pgfscope}%
\pgftext[x=0.888795in,y=0.483333in,,top]{\rmfamily\fontsize{10.000000}{12.000000}\selectfont \(\displaystyle L\)}%
\end{pgfscope}%
\begin{pgfscope}%
\pgfsetbuttcap%
\pgfsetroundjoin%
\definecolor{currentfill}{rgb}{0.000000,0.000000,0.000000}%
\pgfsetfillcolor{currentfill}%
\pgfsetlinewidth{0.803000pt}%
\definecolor{currentstroke}{rgb}{0.000000,0.000000,0.000000}%
\pgfsetstrokecolor{currentstroke}%
\pgfsetdash{}{0pt}%
\pgfsys@defobject{currentmarker}{\pgfqpoint{0.000000in}{-0.048611in}}{\pgfqpoint{0.000000in}{0.000000in}}{%
\pgfpathmoveto{\pgfqpoint{0.000000in}{0.000000in}}%
\pgfpathlineto{\pgfqpoint{0.000000in}{-0.048611in}}%
\pgfusepath{stroke,fill}%
}%
\begin{pgfscope}%
\pgfsys@transformshift{1.159117in}{0.580556in}%
\pgfsys@useobject{currentmarker}{}%
\end{pgfscope}%
\end{pgfscope}%
\begin{pgfscope}%
\pgftext[x=1.159117in,y=0.483333in,,top]{\rmfamily\fontsize{10.000000}{12.000000}\selectfont \(\displaystyle 2L\)}%
\end{pgfscope}%
\begin{pgfscope}%
\pgfsetbuttcap%
\pgfsetroundjoin%
\definecolor{currentfill}{rgb}{0.000000,0.000000,0.000000}%
\pgfsetfillcolor{currentfill}%
\pgfsetlinewidth{0.803000pt}%
\definecolor{currentstroke}{rgb}{0.000000,0.000000,0.000000}%
\pgfsetstrokecolor{currentstroke}%
\pgfsetdash{}{0pt}%
\pgfsys@defobject{currentmarker}{\pgfqpoint{0.000000in}{-0.048611in}}{\pgfqpoint{0.000000in}{0.000000in}}{%
\pgfpathmoveto{\pgfqpoint{0.000000in}{0.000000in}}%
\pgfpathlineto{\pgfqpoint{0.000000in}{-0.048611in}}%
\pgfusepath{stroke,fill}%
}%
\begin{pgfscope}%
\pgfsys@transformshift{1.429438in}{0.580556in}%
\pgfsys@useobject{currentmarker}{}%
\end{pgfscope}%
\end{pgfscope}%
\begin{pgfscope}%
\pgftext[x=1.429438in,y=0.483333in,,top]{\rmfamily\fontsize{10.000000}{12.000000}\selectfont \(\displaystyle 3L\)}%
\end{pgfscope}%
\begin{pgfscope}%
\pgfsetbuttcap%
\pgfsetroundjoin%
\definecolor{currentfill}{rgb}{0.000000,0.000000,0.000000}%
\pgfsetfillcolor{currentfill}%
\pgfsetlinewidth{0.803000pt}%
\definecolor{currentstroke}{rgb}{0.000000,0.000000,0.000000}%
\pgfsetstrokecolor{currentstroke}%
\pgfsetdash{}{0pt}%
\pgfsys@defobject{currentmarker}{\pgfqpoint{0.000000in}{-0.048611in}}{\pgfqpoint{0.000000in}{0.000000in}}{%
\pgfpathmoveto{\pgfqpoint{0.000000in}{0.000000in}}%
\pgfpathlineto{\pgfqpoint{0.000000in}{-0.048611in}}%
\pgfusepath{stroke,fill}%
}%
\begin{pgfscope}%
\pgfsys@transformshift{1.699760in}{0.580556in}%
\pgfsys@useobject{currentmarker}{}%
\end{pgfscope}%
\end{pgfscope}%
\begin{pgfscope}%
\pgftext[x=1.699760in,y=0.483333in,,top]{\rmfamily\fontsize{10.000000}{12.000000}\selectfont \(\displaystyle 4L\)}%
\end{pgfscope}%
\begin{pgfscope}%
\pgftext[x=1.213181in,y=0.305123in,,top]{\rmfamily\fontsize{10.000000}{12.000000}\selectfont Bifurcation value index \(\displaystyle n\)}%
\end{pgfscope}%
\begin{pgfscope}%
\pgfsetbuttcap%
\pgfsetroundjoin%
\definecolor{currentfill}{rgb}{0.000000,0.000000,0.000000}%
\pgfsetfillcolor{currentfill}%
\pgfsetlinewidth{0.803000pt}%
\definecolor{currentstroke}{rgb}{0.000000,0.000000,0.000000}%
\pgfsetstrokecolor{currentstroke}%
\pgfsetdash{}{0pt}%
\pgfsys@defobject{currentmarker}{\pgfqpoint{-0.048611in}{0.000000in}}{\pgfqpoint{0.000000in}{0.000000in}}{%
\pgfpathmoveto{\pgfqpoint{0.000000in}{0.000000in}}%
\pgfpathlineto{\pgfqpoint{-0.048611in}{0.000000in}}%
\pgfusepath{stroke,fill}%
}%
\begin{pgfscope}%
\pgfsys@transformshift{0.618473in}{0.580556in}%
\pgfsys@useobject{currentmarker}{}%
\end{pgfscope}%
\end{pgfscope}%
\begin{pgfscope}%
\pgftext[x=0.343781in,y=0.532728in,left,base]{\rmfamily\fontsize{10.000000}{12.000000}\selectfont \(\displaystyle 0.0\)}%
\end{pgfscope}%
\begin{pgfscope}%
\pgfsetbuttcap%
\pgfsetroundjoin%
\definecolor{currentfill}{rgb}{0.000000,0.000000,0.000000}%
\pgfsetfillcolor{currentfill}%
\pgfsetlinewidth{0.803000pt}%
\definecolor{currentstroke}{rgb}{0.000000,0.000000,0.000000}%
\pgfsetstrokecolor{currentstroke}%
\pgfsetdash{}{0pt}%
\pgfsys@defobject{currentmarker}{\pgfqpoint{-0.048611in}{0.000000in}}{\pgfqpoint{0.000000in}{0.000000in}}{%
\pgfpathmoveto{\pgfqpoint{0.000000in}{0.000000in}}%
\pgfpathlineto{\pgfqpoint{-0.048611in}{0.000000in}}%
\pgfusepath{stroke,fill}%
}%
\begin{pgfscope}%
\pgfsys@transformshift{0.618473in}{0.990057in}%
\pgfsys@useobject{currentmarker}{}%
\end{pgfscope}%
\end{pgfscope}%
\begin{pgfscope}%
\pgftext[x=0.343781in,y=0.942229in,left,base]{\rmfamily\fontsize{10.000000}{12.000000}\selectfont \(\displaystyle 0.1\)}%
\end{pgfscope}%
\begin{pgfscope}%
\pgfsetbuttcap%
\pgfsetroundjoin%
\definecolor{currentfill}{rgb}{0.000000,0.000000,0.000000}%
\pgfsetfillcolor{currentfill}%
\pgfsetlinewidth{0.803000pt}%
\definecolor{currentstroke}{rgb}{0.000000,0.000000,0.000000}%
\pgfsetstrokecolor{currentstroke}%
\pgfsetdash{}{0pt}%
\pgfsys@defobject{currentmarker}{\pgfqpoint{-0.048611in}{0.000000in}}{\pgfqpoint{0.000000in}{0.000000in}}{%
\pgfpathmoveto{\pgfqpoint{0.000000in}{0.000000in}}%
\pgfpathlineto{\pgfqpoint{-0.048611in}{0.000000in}}%
\pgfusepath{stroke,fill}%
}%
\begin{pgfscope}%
\pgfsys@transformshift{0.618473in}{1.399559in}%
\pgfsys@useobject{currentmarker}{}%
\end{pgfscope}%
\end{pgfscope}%
\begin{pgfscope}%
\pgftext[x=0.343781in,y=1.351731in,left,base]{\rmfamily\fontsize{10.000000}{12.000000}\selectfont \(\displaystyle 0.2\)}%
\end{pgfscope}%
\begin{pgfscope}%
\pgfsetbuttcap%
\pgfsetroundjoin%
\definecolor{currentfill}{rgb}{0.000000,0.000000,0.000000}%
\pgfsetfillcolor{currentfill}%
\pgfsetlinewidth{0.803000pt}%
\definecolor{currentstroke}{rgb}{0.000000,0.000000,0.000000}%
\pgfsetstrokecolor{currentstroke}%
\pgfsetdash{}{0pt}%
\pgfsys@defobject{currentmarker}{\pgfqpoint{-0.048611in}{0.000000in}}{\pgfqpoint{0.000000in}{0.000000in}}{%
\pgfpathmoveto{\pgfqpoint{0.000000in}{0.000000in}}%
\pgfpathlineto{\pgfqpoint{-0.048611in}{0.000000in}}%
\pgfusepath{stroke,fill}%
}%
\begin{pgfscope}%
\pgfsys@transformshift{0.618473in}{1.809061in}%
\pgfsys@useobject{currentmarker}{}%
\end{pgfscope}%
\end{pgfscope}%
\begin{pgfscope}%
\pgftext[x=0.343781in,y=1.761233in,left,base]{\rmfamily\fontsize{10.000000}{12.000000}\selectfont \(\displaystyle 0.3\)}%
\end{pgfscope}%
\begin{pgfscope}%
\pgfsetbuttcap%
\pgfsetroundjoin%
\definecolor{currentfill}{rgb}{0.000000,0.000000,0.000000}%
\pgfsetfillcolor{currentfill}%
\pgfsetlinewidth{0.803000pt}%
\definecolor{currentstroke}{rgb}{0.000000,0.000000,0.000000}%
\pgfsetstrokecolor{currentstroke}%
\pgfsetdash{}{0pt}%
\pgfsys@defobject{currentmarker}{\pgfqpoint{-0.048611in}{0.000000in}}{\pgfqpoint{0.000000in}{0.000000in}}{%
\pgfpathmoveto{\pgfqpoint{0.000000in}{0.000000in}}%
\pgfpathlineto{\pgfqpoint{-0.048611in}{0.000000in}}%
\pgfusepath{stroke,fill}%
}%
\begin{pgfscope}%
\pgfsys@transformshift{0.618473in}{2.218562in}%
\pgfsys@useobject{currentmarker}{}%
\end{pgfscope}%
\end{pgfscope}%
\begin{pgfscope}%
\pgftext[x=0.343781in,y=2.170735in,left,base]{\rmfamily\fontsize{10.000000}{12.000000}\selectfont \(\displaystyle 0.4\)}%
\end{pgfscope}%
\begin{pgfscope}%
\pgftext[x=0.288226in,y=1.399559in,,bottom,rotate=90.000000]{\rmfamily\fontsize{10.000000}{12.000000}\selectfont \(\displaystyle M_n(\omega)\)}%
\end{pgfscope}%
\begin{pgfscope}%
\pgfpathrectangle{\pgfqpoint{0.618473in}{0.580556in}}{\pgfqpoint{1.189415in}{1.638007in}}%
\pgfusepath{clip}%
\pgfsetrectcap%
\pgfsetroundjoin%
\pgfsetlinewidth{1.505625pt}%
\definecolor{currentstroke}{rgb}{0.000000,0.000000,0.000000}%
\pgfsetstrokecolor{currentstroke}%
\pgfsetdash{}{0pt}%
\pgfpathmoveto{\pgfqpoint{0.623880in}{0.709035in}}%
\pgfpathlineto{\pgfqpoint{0.650912in}{1.316566in}}%
\pgfpathlineto{\pgfqpoint{0.667131in}{1.620122in}}%
\pgfpathlineto{\pgfqpoint{0.677944in}{1.784233in}}%
\pgfpathlineto{\pgfqpoint{0.688757in}{1.912604in}}%
\pgfpathlineto{\pgfqpoint{0.699570in}{2.002458in}}%
\pgfpathlineto{\pgfqpoint{0.704976in}{2.032494in}}%
\pgfpathlineto{\pgfqpoint{0.710383in}{2.052560in}}%
\pgfpathlineto{\pgfqpoint{0.715789in}{2.062748in}}%
\pgfpathlineto{\pgfqpoint{0.721196in}{2.063242in}}%
\pgfpathlineto{\pgfqpoint{0.726602in}{2.054321in}}%
\pgfpathlineto{\pgfqpoint{0.732008in}{2.036351in}}%
\pgfpathlineto{\pgfqpoint{0.737415in}{2.009779in}}%
\pgfpathlineto{\pgfqpoint{0.748228in}{1.932999in}}%
\pgfpathlineto{\pgfqpoint{0.759041in}{1.828964in}}%
\pgfpathlineto{\pgfqpoint{0.775260in}{1.634752in}}%
\pgfpathlineto{\pgfqpoint{0.834731in}{0.862020in}}%
\pgfpathlineto{\pgfqpoint{0.845544in}{0.760628in}}%
\pgfpathlineto{\pgfqpoint{0.856356in}{0.680921in}}%
\pgfpathlineto{\pgfqpoint{0.867169in}{0.624369in}}%
\pgfpathlineto{\pgfqpoint{0.872576in}{0.604900in}}%
\pgfpathlineto{\pgfqpoint{0.877982in}{0.591220in}}%
\pgfpathlineto{\pgfqpoint{0.883389in}{0.583178in}}%
\pgfpathlineto{\pgfqpoint{0.888795in}{0.580556in}}%
\pgfpathlineto{\pgfqpoint{0.894201in}{0.583075in}}%
\pgfpathlineto{\pgfqpoint{0.899608in}{0.590400in}}%
\pgfpathlineto{\pgfqpoint{0.905014in}{0.602144in}}%
\pgfpathlineto{\pgfqpoint{0.915827in}{0.637134in}}%
\pgfpathlineto{\pgfqpoint{0.926640in}{0.684199in}}%
\pgfpathlineto{\pgfqpoint{0.948266in}{0.797612in}}%
\pgfpathlineto{\pgfqpoint{0.969891in}{0.908687in}}%
\pgfpathlineto{\pgfqpoint{0.980704in}{0.954049in}}%
\pgfpathlineto{\pgfqpoint{0.991517in}{0.988832in}}%
\pgfpathlineto{\pgfqpoint{1.002330in}{1.011143in}}%
\pgfpathlineto{\pgfqpoint{1.007737in}{1.017263in}}%
\pgfpathlineto{\pgfqpoint{1.013143in}{1.019942in}}%
\pgfpathlineto{\pgfqpoint{1.018549in}{1.019186in}}%
\pgfpathlineto{\pgfqpoint{1.023956in}{1.015050in}}%
\pgfpathlineto{\pgfqpoint{1.029362in}{1.007643in}}%
\pgfpathlineto{\pgfqpoint{1.040175in}{0.983670in}}%
\pgfpathlineto{\pgfqpoint{1.050988in}{0.948997in}}%
\pgfpathlineto{\pgfqpoint{1.067207in}{0.882050in}}%
\pgfpathlineto{\pgfqpoint{1.115865in}{0.662762in}}%
\pgfpathlineto{\pgfqpoint{1.126678in}{0.627535in}}%
\pgfpathlineto{\pgfqpoint{1.137491in}{0.601549in}}%
\pgfpathlineto{\pgfqpoint{1.148304in}{0.585779in}}%
\pgfpathlineto{\pgfqpoint{1.153710in}{0.581853in}}%
\pgfpathlineto{\pgfqpoint{1.159117in}{0.580556in}}%
\pgfpathlineto{\pgfqpoint{1.164523in}{0.581828in}}%
\pgfpathlineto{\pgfqpoint{1.169930in}{0.585574in}}%
\pgfpathlineto{\pgfqpoint{1.180742in}{0.599934in}}%
\pgfpathlineto{\pgfqpoint{1.191555in}{0.622217in}}%
\pgfpathlineto{\pgfqpoint{1.207775in}{0.666391in}}%
\pgfpathlineto{\pgfqpoint{1.251026in}{0.794437in}}%
\pgfpathlineto{\pgfqpoint{1.261839in}{0.817287in}}%
\pgfpathlineto{\pgfqpoint{1.272652in}{0.833214in}}%
\pgfpathlineto{\pgfqpoint{1.283465in}{0.841330in}}%
\pgfpathlineto{\pgfqpoint{1.288871in}{0.842319in}}%
\pgfpathlineto{\pgfqpoint{1.294277in}{0.841252in}}%
\pgfpathlineto{\pgfqpoint{1.299684in}{0.838164in}}%
\pgfpathlineto{\pgfqpoint{1.310497in}{0.826203in}}%
\pgfpathlineto{\pgfqpoint{1.321310in}{0.807289in}}%
\pgfpathlineto{\pgfqpoint{1.337529in}{0.768707in}}%
\pgfpathlineto{\pgfqpoint{1.391593in}{0.621868in}}%
\pgfpathlineto{\pgfqpoint{1.402406in}{0.602016in}}%
\pgfpathlineto{\pgfqpoint{1.413219in}{0.588339in}}%
\pgfpathlineto{\pgfqpoint{1.424032in}{0.581418in}}%
\pgfpathlineto{\pgfqpoint{1.429438in}{0.580556in}}%
\pgfpathlineto{\pgfqpoint{1.434845in}{0.581406in}}%
\pgfpathlineto{\pgfqpoint{1.445658in}{0.588034in}}%
\pgfpathlineto{\pgfqpoint{1.456470in}{0.600632in}}%
\pgfpathlineto{\pgfqpoint{1.467283in}{0.618184in}}%
\pgfpathlineto{\pgfqpoint{1.483503in}{0.650922in}}%
\pgfpathlineto{\pgfqpoint{1.521348in}{0.730400in}}%
\pgfpathlineto{\pgfqpoint{1.532161in}{0.747248in}}%
\pgfpathlineto{\pgfqpoint{1.542973in}{0.759337in}}%
\pgfpathlineto{\pgfqpoint{1.553786in}{0.765961in}}%
\pgfpathlineto{\pgfqpoint{1.564599in}{0.766768in}}%
\pgfpathlineto{\pgfqpoint{1.575412in}{0.761771in}}%
\pgfpathlineto{\pgfqpoint{1.586225in}{0.751347in}}%
\pgfpathlineto{\pgfqpoint{1.597038in}{0.736197in}}%
\pgfpathlineto{\pgfqpoint{1.613257in}{0.706811in}}%
\pgfpathlineto{\pgfqpoint{1.661915in}{0.611165in}}%
\pgfpathlineto{\pgfqpoint{1.672728in}{0.596513in}}%
\pgfpathlineto{\pgfqpoint{1.683541in}{0.586364in}}%
\pgfpathlineto{\pgfqpoint{1.694354in}{0.581201in}}%
\pgfpathlineto{\pgfqpoint{1.705166in}{0.581195in}}%
\pgfpathlineto{\pgfqpoint{1.715979in}{0.586192in}}%
\pgfpathlineto{\pgfqpoint{1.726792in}{0.595735in}}%
\pgfpathlineto{\pgfqpoint{1.737605in}{0.609095in}}%
\pgfpathlineto{\pgfqpoint{1.753824in}{0.634168in}}%
\pgfpathlineto{\pgfqpoint{1.791669in}{0.695874in}}%
\pgfpathlineto{\pgfqpoint{1.807889in}{0.714534in}}%
\pgfpathlineto{\pgfqpoint{1.818701in}{0.722190in}}%
\pgfpathlineto{\pgfqpoint{1.821778in}{0.723442in}}%
\pgfpathlineto{\pgfqpoint{1.821778in}{0.723442in}}%
\pgfusepath{stroke}%
\end{pgfscope}%
\begin{pgfscope}%
\pgfsetrectcap%
\pgfsetmiterjoin%
\pgfsetlinewidth{0.803000pt}%
\definecolor{currentstroke}{rgb}{0.000000,0.000000,0.000000}%
\pgfsetstrokecolor{currentstroke}%
\pgfsetdash{}{0pt}%
\pgfpathmoveto{\pgfqpoint{0.618473in}{0.580556in}}%
\pgfpathlineto{\pgfqpoint{0.618473in}{2.218562in}}%
\pgfusepath{stroke}%
\end{pgfscope}%
\begin{pgfscope}%
\pgfsetrectcap%
\pgfsetmiterjoin%
\pgfsetlinewidth{0.803000pt}%
\definecolor{currentstroke}{rgb}{0.000000,0.000000,0.000000}%
\pgfsetstrokecolor{currentstroke}%
\pgfsetdash{}{0pt}%
\pgfpathmoveto{\pgfqpoint{1.807889in}{0.580556in}}%
\pgfpathlineto{\pgfqpoint{1.807889in}{2.218562in}}%
\pgfusepath{stroke}%
\end{pgfscope}%
\begin{pgfscope}%
\pgfsetrectcap%
\pgfsetmiterjoin%
\pgfsetlinewidth{0.803000pt}%
\definecolor{currentstroke}{rgb}{0.000000,0.000000,0.000000}%
\pgfsetstrokecolor{currentstroke}%
\pgfsetdash{}{0pt}%
\pgfpathmoveto{\pgfqpoint{0.618473in}{0.580556in}}%
\pgfpathlineto{\pgfqpoint{1.807889in}{0.580556in}}%
\pgfusepath{stroke}%
\end{pgfscope}%
\begin{pgfscope}%
\pgfsetrectcap%
\pgfsetmiterjoin%
\pgfsetlinewidth{0.803000pt}%
\definecolor{currentstroke}{rgb}{0.000000,0.000000,0.000000}%
\pgfsetstrokecolor{currentstroke}%
\pgfsetdash{}{0pt}%
\pgfpathmoveto{\pgfqpoint{0.618473in}{2.218562in}}%
\pgfpathlineto{\pgfqpoint{1.807889in}{2.218562in}}%
\pgfusepath{stroke}%
\end{pgfscope}%
\begin{pgfscope}%
\pgfsetbuttcap%
\pgfsetmiterjoin%
\definecolor{currentfill}{rgb}{1.000000,1.000000,1.000000}%
\pgfsetfillcolor{currentfill}%
\pgfsetlinewidth{0.000000pt}%
\definecolor{currentstroke}{rgb}{0.000000,0.000000,0.000000}%
\pgfsetstrokecolor{currentstroke}%
\pgfsetstrokeopacity{0.000000}%
\pgfsetdash{}{0pt}%
\pgfpathmoveto{\pgfqpoint{2.504626in}{0.580556in}}%
\pgfpathlineto{\pgfqpoint{3.694041in}{0.580556in}}%
\pgfpathlineto{\pgfqpoint{3.694041in}{2.218562in}}%
\pgfpathlineto{\pgfqpoint{2.504626in}{2.218562in}}%
\pgfpathclose%
\pgfusepath{fill}%
\end{pgfscope}%
\begin{pgfscope}%
\pgfsetbuttcap%
\pgfsetroundjoin%
\definecolor{currentfill}{rgb}{0.000000,0.000000,0.000000}%
\pgfsetfillcolor{currentfill}%
\pgfsetlinewidth{0.803000pt}%
\definecolor{currentstroke}{rgb}{0.000000,0.000000,0.000000}%
\pgfsetstrokecolor{currentstroke}%
\pgfsetdash{}{0pt}%
\pgfsys@defobject{currentmarker}{\pgfqpoint{0.000000in}{-0.048611in}}{\pgfqpoint{0.000000in}{0.000000in}}{%
\pgfpathmoveto{\pgfqpoint{0.000000in}{0.000000in}}%
\pgfpathlineto{\pgfqpoint{0.000000in}{-0.048611in}}%
\pgfusepath{stroke,fill}%
}%
\begin{pgfscope}%
\pgfsys@transformshift{2.504626in}{0.580556in}%
\pgfsys@useobject{currentmarker}{}%
\end{pgfscope}%
\end{pgfscope}%
\begin{pgfscope}%
\pgftext[x=2.504626in,y=0.483333in,,top]{\rmfamily\fontsize{10.000000}{12.000000}\selectfont \(\displaystyle 0\)}%
\end{pgfscope}%
\begin{pgfscope}%
\pgfsetbuttcap%
\pgfsetroundjoin%
\definecolor{currentfill}{rgb}{0.000000,0.000000,0.000000}%
\pgfsetfillcolor{currentfill}%
\pgfsetlinewidth{0.803000pt}%
\definecolor{currentstroke}{rgb}{0.000000,0.000000,0.000000}%
\pgfsetstrokecolor{currentstroke}%
\pgfsetdash{}{0pt}%
\pgfsys@defobject{currentmarker}{\pgfqpoint{0.000000in}{-0.048611in}}{\pgfqpoint{0.000000in}{0.000000in}}{%
\pgfpathmoveto{\pgfqpoint{0.000000in}{0.000000in}}%
\pgfpathlineto{\pgfqpoint{0.000000in}{-0.048611in}}%
\pgfusepath{stroke,fill}%
}%
\begin{pgfscope}%
\pgfsys@transformshift{2.774948in}{0.580556in}%
\pgfsys@useobject{currentmarker}{}%
\end{pgfscope}%
\end{pgfscope}%
\begin{pgfscope}%
\pgftext[x=2.774948in,y=0.483333in,,top]{\rmfamily\fontsize{10.000000}{12.000000}\selectfont \(\displaystyle L\)}%
\end{pgfscope}%
\begin{pgfscope}%
\pgfsetbuttcap%
\pgfsetroundjoin%
\definecolor{currentfill}{rgb}{0.000000,0.000000,0.000000}%
\pgfsetfillcolor{currentfill}%
\pgfsetlinewidth{0.803000pt}%
\definecolor{currentstroke}{rgb}{0.000000,0.000000,0.000000}%
\pgfsetstrokecolor{currentstroke}%
\pgfsetdash{}{0pt}%
\pgfsys@defobject{currentmarker}{\pgfqpoint{0.000000in}{-0.048611in}}{\pgfqpoint{0.000000in}{0.000000in}}{%
\pgfpathmoveto{\pgfqpoint{0.000000in}{0.000000in}}%
\pgfpathlineto{\pgfqpoint{0.000000in}{-0.048611in}}%
\pgfusepath{stroke,fill}%
}%
\begin{pgfscope}%
\pgfsys@transformshift{3.045269in}{0.580556in}%
\pgfsys@useobject{currentmarker}{}%
\end{pgfscope}%
\end{pgfscope}%
\begin{pgfscope}%
\pgftext[x=3.045269in,y=0.483333in,,top]{\rmfamily\fontsize{10.000000}{12.000000}\selectfont \(\displaystyle 2L\)}%
\end{pgfscope}%
\begin{pgfscope}%
\pgfsetbuttcap%
\pgfsetroundjoin%
\definecolor{currentfill}{rgb}{0.000000,0.000000,0.000000}%
\pgfsetfillcolor{currentfill}%
\pgfsetlinewidth{0.803000pt}%
\definecolor{currentstroke}{rgb}{0.000000,0.000000,0.000000}%
\pgfsetstrokecolor{currentstroke}%
\pgfsetdash{}{0pt}%
\pgfsys@defobject{currentmarker}{\pgfqpoint{0.000000in}{-0.048611in}}{\pgfqpoint{0.000000in}{0.000000in}}{%
\pgfpathmoveto{\pgfqpoint{0.000000in}{0.000000in}}%
\pgfpathlineto{\pgfqpoint{0.000000in}{-0.048611in}}%
\pgfusepath{stroke,fill}%
}%
\begin{pgfscope}%
\pgfsys@transformshift{3.315591in}{0.580556in}%
\pgfsys@useobject{currentmarker}{}%
\end{pgfscope}%
\end{pgfscope}%
\begin{pgfscope}%
\pgftext[x=3.315591in,y=0.483333in,,top]{\rmfamily\fontsize{10.000000}{12.000000}\selectfont \(\displaystyle 3L\)}%
\end{pgfscope}%
\begin{pgfscope}%
\pgfsetbuttcap%
\pgfsetroundjoin%
\definecolor{currentfill}{rgb}{0.000000,0.000000,0.000000}%
\pgfsetfillcolor{currentfill}%
\pgfsetlinewidth{0.803000pt}%
\definecolor{currentstroke}{rgb}{0.000000,0.000000,0.000000}%
\pgfsetstrokecolor{currentstroke}%
\pgfsetdash{}{0pt}%
\pgfsys@defobject{currentmarker}{\pgfqpoint{0.000000in}{-0.048611in}}{\pgfqpoint{0.000000in}{0.000000in}}{%
\pgfpathmoveto{\pgfqpoint{0.000000in}{0.000000in}}%
\pgfpathlineto{\pgfqpoint{0.000000in}{-0.048611in}}%
\pgfusepath{stroke,fill}%
}%
\begin{pgfscope}%
\pgfsys@transformshift{3.585913in}{0.580556in}%
\pgfsys@useobject{currentmarker}{}%
\end{pgfscope}%
\end{pgfscope}%
\begin{pgfscope}%
\pgftext[x=3.585913in,y=0.483333in,,top]{\rmfamily\fontsize{10.000000}{12.000000}\selectfont \(\displaystyle 4L\)}%
\end{pgfscope}%
\begin{pgfscope}%
\pgftext[x=3.099334in,y=0.305123in,,top]{\rmfamily\fontsize{10.000000}{12.000000}\selectfont Bifurcation value index \(\displaystyle n\)}%
\end{pgfscope}%
\begin{pgfscope}%
\pgfsetbuttcap%
\pgfsetroundjoin%
\definecolor{currentfill}{rgb}{0.000000,0.000000,0.000000}%
\pgfsetfillcolor{currentfill}%
\pgfsetlinewidth{0.803000pt}%
\definecolor{currentstroke}{rgb}{0.000000,0.000000,0.000000}%
\pgfsetstrokecolor{currentstroke}%
\pgfsetdash{}{0pt}%
\pgfsys@defobject{currentmarker}{\pgfqpoint{-0.048611in}{0.000000in}}{\pgfqpoint{0.000000in}{0.000000in}}{%
\pgfpathmoveto{\pgfqpoint{0.000000in}{0.000000in}}%
\pgfpathlineto{\pgfqpoint{-0.048611in}{0.000000in}}%
\pgfusepath{stroke,fill}%
}%
\begin{pgfscope}%
\pgfsys@transformshift{2.504626in}{0.580556in}%
\pgfsys@useobject{currentmarker}{}%
\end{pgfscope}%
\end{pgfscope}%
\begin{pgfscope}%
\pgftext[x=2.229934in,y=0.532728in,left,base]{\rmfamily\fontsize{10.000000}{12.000000}\selectfont \(\displaystyle 0.0\)}%
\end{pgfscope}%
\begin{pgfscope}%
\pgfsetbuttcap%
\pgfsetroundjoin%
\definecolor{currentfill}{rgb}{0.000000,0.000000,0.000000}%
\pgfsetfillcolor{currentfill}%
\pgfsetlinewidth{0.803000pt}%
\definecolor{currentstroke}{rgb}{0.000000,0.000000,0.000000}%
\pgfsetstrokecolor{currentstroke}%
\pgfsetdash{}{0pt}%
\pgfsys@defobject{currentmarker}{\pgfqpoint{-0.048611in}{0.000000in}}{\pgfqpoint{0.000000in}{0.000000in}}{%
\pgfpathmoveto{\pgfqpoint{0.000000in}{0.000000in}}%
\pgfpathlineto{\pgfqpoint{-0.048611in}{0.000000in}}%
\pgfusepath{stroke,fill}%
}%
\begin{pgfscope}%
\pgfsys@transformshift{2.504626in}{1.048558in}%
\pgfsys@useobject{currentmarker}{}%
\end{pgfscope}%
\end{pgfscope}%
\begin{pgfscope}%
\pgftext[x=2.229934in,y=1.000730in,left,base]{\rmfamily\fontsize{10.000000}{12.000000}\selectfont \(\displaystyle 0.2\)}%
\end{pgfscope}%
\begin{pgfscope}%
\pgfsetbuttcap%
\pgfsetroundjoin%
\definecolor{currentfill}{rgb}{0.000000,0.000000,0.000000}%
\pgfsetfillcolor{currentfill}%
\pgfsetlinewidth{0.803000pt}%
\definecolor{currentstroke}{rgb}{0.000000,0.000000,0.000000}%
\pgfsetstrokecolor{currentstroke}%
\pgfsetdash{}{0pt}%
\pgfsys@defobject{currentmarker}{\pgfqpoint{-0.048611in}{0.000000in}}{\pgfqpoint{0.000000in}{0.000000in}}{%
\pgfpathmoveto{\pgfqpoint{0.000000in}{0.000000in}}%
\pgfpathlineto{\pgfqpoint{-0.048611in}{0.000000in}}%
\pgfusepath{stroke,fill}%
}%
\begin{pgfscope}%
\pgfsys@transformshift{2.504626in}{1.516559in}%
\pgfsys@useobject{currentmarker}{}%
\end{pgfscope}%
\end{pgfscope}%
\begin{pgfscope}%
\pgftext[x=2.229934in,y=1.468732in,left,base]{\rmfamily\fontsize{10.000000}{12.000000}\selectfont \(\displaystyle 0.4\)}%
\end{pgfscope}%
\begin{pgfscope}%
\pgfsetbuttcap%
\pgfsetroundjoin%
\definecolor{currentfill}{rgb}{0.000000,0.000000,0.000000}%
\pgfsetfillcolor{currentfill}%
\pgfsetlinewidth{0.803000pt}%
\definecolor{currentstroke}{rgb}{0.000000,0.000000,0.000000}%
\pgfsetstrokecolor{currentstroke}%
\pgfsetdash{}{0pt}%
\pgfsys@defobject{currentmarker}{\pgfqpoint{-0.048611in}{0.000000in}}{\pgfqpoint{0.000000in}{0.000000in}}{%
\pgfpathmoveto{\pgfqpoint{0.000000in}{0.000000in}}%
\pgfpathlineto{\pgfqpoint{-0.048611in}{0.000000in}}%
\pgfusepath{stroke,fill}%
}%
\begin{pgfscope}%
\pgfsys@transformshift{2.504626in}{1.984561in}%
\pgfsys@useobject{currentmarker}{}%
\end{pgfscope}%
\end{pgfscope}%
\begin{pgfscope}%
\pgftext[x=2.229934in,y=1.936734in,left,base]{\rmfamily\fontsize{10.000000}{12.000000}\selectfont \(\displaystyle 0.6\)}%
\end{pgfscope}%
\begin{pgfscope}%
\pgftext[x=2.174378in,y=1.399559in,,bottom,rotate=90.000000]{\rmfamily\fontsize{10.000000}{12.000000}\selectfont \(\displaystyle 2M_n(\omega) - M_{2n}(\omega)\)}%
\end{pgfscope}%
\begin{pgfscope}%
\pgfpathrectangle{\pgfqpoint{2.504626in}{0.580556in}}{\pgfqpoint{1.189415in}{1.638007in}}%
\pgfusepath{clip}%
\pgfsetrectcap%
\pgfsetroundjoin%
\pgfsetlinewidth{1.505625pt}%
\definecolor{currentstroke}{rgb}{0.000000,0.000000,0.000000}%
\pgfsetstrokecolor{currentstroke}%
\pgfsetdash{}{0pt}%
\pgfpathmoveto{\pgfqpoint{2.510032in}{0.581134in}}%
\pgfpathlineto{\pgfqpoint{2.515439in}{0.585150in}}%
\pgfpathlineto{\pgfqpoint{2.520845in}{0.595860in}}%
\pgfpathlineto{\pgfqpoint{2.526252in}{0.616169in}}%
\pgfpathlineto{\pgfqpoint{2.531658in}{0.648475in}}%
\pgfpathlineto{\pgfqpoint{2.537065in}{0.694545in}}%
\pgfpathlineto{\pgfqpoint{2.547877in}{0.831310in}}%
\pgfpathlineto{\pgfqpoint{2.558690in}{1.025097in}}%
\pgfpathlineto{\pgfqpoint{2.574910in}{1.389521in}}%
\pgfpathlineto{\pgfqpoint{2.596535in}{1.872409in}}%
\pgfpathlineto{\pgfqpoint{2.607348in}{2.045423in}}%
\pgfpathlineto{\pgfqpoint{2.612755in}{2.104023in}}%
\pgfpathlineto{\pgfqpoint{2.618161in}{2.141424in}}%
\pgfpathlineto{\pgfqpoint{2.623567in}{2.156603in}}%
\pgfpathlineto{\pgfqpoint{2.628974in}{2.149319in}}%
\pgfpathlineto{\pgfqpoint{2.634380in}{2.120111in}}%
\pgfpathlineto{\pgfqpoint{2.639787in}{2.070252in}}%
\pgfpathlineto{\pgfqpoint{2.650600in}{1.916910in}}%
\pgfpathlineto{\pgfqpoint{2.666819in}{1.596200in}}%
\pgfpathlineto{\pgfqpoint{2.693851in}{1.036383in}}%
\pgfpathlineto{\pgfqpoint{2.704664in}{0.864787in}}%
\pgfpathlineto{\pgfqpoint{2.715477in}{0.738203in}}%
\pgfpathlineto{\pgfqpoint{2.726290in}{0.655433in}}%
\pgfpathlineto{\pgfqpoint{2.737103in}{0.609020in}}%
\pgfpathlineto{\pgfqpoint{2.742509in}{0.596100in}}%
\pgfpathlineto{\pgfqpoint{2.747915in}{0.588102in}}%
\pgfpathlineto{\pgfqpoint{2.753322in}{0.583652in}}%
\pgfpathlineto{\pgfqpoint{2.764135in}{0.580747in}}%
\pgfpathlineto{\pgfqpoint{2.785760in}{0.580732in}}%
\pgfpathlineto{\pgfqpoint{2.796573in}{0.583194in}}%
\pgfpathlineto{\pgfqpoint{2.801980in}{0.586730in}}%
\pgfpathlineto{\pgfqpoint{2.807386in}{0.592769in}}%
\pgfpathlineto{\pgfqpoint{2.812793in}{0.602029in}}%
\pgfpathlineto{\pgfqpoint{2.818199in}{0.615142in}}%
\pgfpathlineto{\pgfqpoint{2.829012in}{0.654646in}}%
\pgfpathlineto{\pgfqpoint{2.839825in}{0.712460in}}%
\pgfpathlineto{\pgfqpoint{2.856044in}{0.826001in}}%
\pgfpathlineto{\pgfqpoint{2.877670in}{0.983925in}}%
\pgfpathlineto{\pgfqpoint{2.888483in}{1.042221in}}%
\pgfpathlineto{\pgfqpoint{2.893889in}{1.062125in}}%
\pgfpathlineto{\pgfqpoint{2.899296in}{1.074824in}}%
\pgfpathlineto{\pgfqpoint{2.904702in}{1.079871in}}%
\pgfpathlineto{\pgfqpoint{2.910108in}{1.077121in}}%
\pgfpathlineto{\pgfqpoint{2.915515in}{1.066731in}}%
\pgfpathlineto{\pgfqpoint{2.920921in}{1.049147in}}%
\pgfpathlineto{\pgfqpoint{2.931734in}{0.995470in}}%
\pgfpathlineto{\pgfqpoint{2.947953in}{0.884985in}}%
\pgfpathlineto{\pgfqpoint{2.969579in}{0.733193in}}%
\pgfpathlineto{\pgfqpoint{2.980392in}{0.673488in}}%
\pgfpathlineto{\pgfqpoint{2.991205in}{0.629949in}}%
\pgfpathlineto{\pgfqpoint{3.002018in}{0.602360in}}%
\pgfpathlineto{\pgfqpoint{3.012831in}{0.587831in}}%
\pgfpathlineto{\pgfqpoint{3.023644in}{0.582039in}}%
\pgfpathlineto{\pgfqpoint{3.039863in}{0.580561in}}%
\pgfpathlineto{\pgfqpoint{3.066895in}{0.581925in}}%
\pgfpathlineto{\pgfqpoint{3.077708in}{0.587008in}}%
\pgfpathlineto{\pgfqpoint{3.083114in}{0.591995in}}%
\pgfpathlineto{\pgfqpoint{3.093927in}{0.608721in}}%
\pgfpathlineto{\pgfqpoint{3.104740in}{0.635945in}}%
\pgfpathlineto{\pgfqpoint{3.115553in}{0.673622in}}%
\pgfpathlineto{\pgfqpoint{3.137179in}{0.768261in}}%
\pgfpathlineto{\pgfqpoint{3.153398in}{0.834467in}}%
\pgfpathlineto{\pgfqpoint{3.164211in}{0.864761in}}%
\pgfpathlineto{\pgfqpoint{3.169617in}{0.873902in}}%
\pgfpathlineto{\pgfqpoint{3.175024in}{0.878534in}}%
\pgfpathlineto{\pgfqpoint{3.180430in}{0.878495in}}%
\pgfpathlineto{\pgfqpoint{3.185837in}{0.873804in}}%
\pgfpathlineto{\pgfqpoint{3.191243in}{0.864662in}}%
\pgfpathlineto{\pgfqpoint{3.202056in}{0.834650in}}%
\pgfpathlineto{\pgfqpoint{3.218275in}{0.769671in}}%
\pgfpathlineto{\pgfqpoint{3.239901in}{0.677076in}}%
\pgfpathlineto{\pgfqpoint{3.250714in}{0.639817in}}%
\pgfpathlineto{\pgfqpoint{3.261527in}{0.612309in}}%
\pgfpathlineto{\pgfqpoint{3.272339in}{0.594683in}}%
\pgfpathlineto{\pgfqpoint{3.283152in}{0.585305in}}%
\pgfpathlineto{\pgfqpoint{3.293965in}{0.581531in}}%
\pgfpathlineto{\pgfqpoint{3.310184in}{0.580559in}}%
\pgfpathlineto{\pgfqpoint{3.337217in}{0.581481in}}%
\pgfpathlineto{\pgfqpoint{3.348030in}{0.584940in}}%
\pgfpathlineto{\pgfqpoint{3.358842in}{0.593252in}}%
\pgfpathlineto{\pgfqpoint{3.369655in}{0.608339in}}%
\pgfpathlineto{\pgfqpoint{3.380468in}{0.631038in}}%
\pgfpathlineto{\pgfqpoint{3.396687in}{0.677246in}}%
\pgfpathlineto{\pgfqpoint{3.423720in}{0.759787in}}%
\pgfpathlineto{\pgfqpoint{3.434532in}{0.782143in}}%
\pgfpathlineto{\pgfqpoint{3.439939in}{0.789120in}}%
\pgfpathlineto{\pgfqpoint{3.445345in}{0.792908in}}%
\pgfpathlineto{\pgfqpoint{3.450752in}{0.793369in}}%
\pgfpathlineto{\pgfqpoint{3.456158in}{0.790495in}}%
\pgfpathlineto{\pgfqpoint{3.461565in}{0.784406in}}%
\pgfpathlineto{\pgfqpoint{3.472377in}{0.763674in}}%
\pgfpathlineto{\pgfqpoint{3.483190in}{0.734326in}}%
\pgfpathlineto{\pgfqpoint{3.515629in}{0.636794in}}%
\pgfpathlineto{\pgfqpoint{3.526442in}{0.613086in}}%
\pgfpathlineto{\pgfqpoint{3.537255in}{0.596629in}}%
\pgfpathlineto{\pgfqpoint{3.548068in}{0.586897in}}%
\pgfpathlineto{\pgfqpoint{3.558880in}{0.582297in}}%
\pgfpathlineto{\pgfqpoint{3.575100in}{0.580602in}}%
\pgfpathlineto{\pgfqpoint{3.607538in}{0.581254in}}%
\pgfpathlineto{\pgfqpoint{3.618351in}{0.583876in}}%
\pgfpathlineto{\pgfqpoint{3.629164in}{0.590200in}}%
\pgfpathlineto{\pgfqpoint{3.639977in}{0.601724in}}%
\pgfpathlineto{\pgfqpoint{3.650790in}{0.619131in}}%
\pgfpathlineto{\pgfqpoint{3.667009in}{0.654760in}}%
\pgfpathlineto{\pgfqpoint{3.694041in}{0.719053in}}%
\pgfpathlineto{\pgfqpoint{3.704854in}{0.736740in}}%
\pgfpathlineto{\pgfqpoint{3.707930in}{0.739936in}}%
\pgfpathlineto{\pgfqpoint{3.707930in}{0.739936in}}%
\pgfusepath{stroke}%
\end{pgfscope}%
\begin{pgfscope}%
\pgfsetrectcap%
\pgfsetmiterjoin%
\pgfsetlinewidth{0.803000pt}%
\definecolor{currentstroke}{rgb}{0.000000,0.000000,0.000000}%
\pgfsetstrokecolor{currentstroke}%
\pgfsetdash{}{0pt}%
\pgfpathmoveto{\pgfqpoint{2.504626in}{0.580556in}}%
\pgfpathlineto{\pgfqpoint{2.504626in}{2.218562in}}%
\pgfusepath{stroke}%
\end{pgfscope}%
\begin{pgfscope}%
\pgfsetrectcap%
\pgfsetmiterjoin%
\pgfsetlinewidth{0.803000pt}%
\definecolor{currentstroke}{rgb}{0.000000,0.000000,0.000000}%
\pgfsetstrokecolor{currentstroke}%
\pgfsetdash{}{0pt}%
\pgfpathmoveto{\pgfqpoint{3.694041in}{0.580556in}}%
\pgfpathlineto{\pgfqpoint{3.694041in}{2.218562in}}%
\pgfusepath{stroke}%
\end{pgfscope}%
\begin{pgfscope}%
\pgfsetrectcap%
\pgfsetmiterjoin%
\pgfsetlinewidth{0.803000pt}%
\definecolor{currentstroke}{rgb}{0.000000,0.000000,0.000000}%
\pgfsetstrokecolor{currentstroke}%
\pgfsetdash{}{0pt}%
\pgfpathmoveto{\pgfqpoint{2.504626in}{0.580556in}}%
\pgfpathlineto{\pgfqpoint{3.694041in}{0.580556in}}%
\pgfusepath{stroke}%
\end{pgfscope}%
\begin{pgfscope}%
\pgfsetrectcap%
\pgfsetmiterjoin%
\pgfsetlinewidth{0.803000pt}%
\definecolor{currentstroke}{rgb}{0.000000,0.000000,0.000000}%
\pgfsetstrokecolor{currentstroke}%
\pgfsetdash{}{0pt}%
\pgfpathmoveto{\pgfqpoint{2.504626in}{2.218562in}}%
\pgfpathlineto{\pgfqpoint{3.694041in}{2.218562in}}%
\pgfusepath{stroke}%
\end{pgfscope}%
\end{pgfpicture}%
\makeatother%
\endgroup%

%% file: chapterNoFluxBcShort.tex
\chapter[No-flux Boundary Conditions]{No-flux boundary conditions\index{no-flux
boundary conditions} for non-local operators}\label{Chapter:AdhNoFlux}
It is a challenge to define boundary conditions for non-local models on bounded
domains. The periodic case, which we studied in the previous chapters, is an
exception, since we can work with periodic extensions outside of the domain.
However, no-flux conditions\index{no-flux boundary conditions} or
Dirichlet\index{Dirichlet} or Robin boundary conditions\index{Robin boundary
conditions} need special attention.

In the literature there are several ideas to deal with boundary conditions of
non-local operators. For example, the inclusion of a potential that diverges
close to the boundary can keep particles away from
it~\cite{fetecau2017swarm,wu2015nonlocal}. One could also simply cut off the
part of the integral that reaches outside of the domain, as done, for example
in~\cite{Hillen2007}. Another idea, often used in numerical simulations, is to
introduce ghost points\index{ghost points} outside of the domain with certain symmetry
properties~\cite{Allen1998}.

In our case, we are interested in the motion of biological objects such as cells
or bacteria, hence we will formulate the non-local boundary conditions from a
biological point of view. A cell can attach to the boundary, it could slip off
it, or be neutral relative to the boundary, and our boundary conditions need to
be able to describe these cases.

Next we briefly summarize previous ideas in models for biological populations.
Hillen~\etal\ \parencite{Hillen2007} considered the global existence of a
non-local chemotaxis equations. To correctly define the non-local
chemotaxis\index{non-local chemotaxis}
equation on a bounded domain, they limited the set over which the non-local
operator~\eqref{Eqn:NonLocalGradientChemotaxis} integrates
\begin{equation}\label{Eqn:NonLocalGradientChemotaxis}
    \nonlocalgrad{R} v(x) \coloneqq \frac{n}{\omega_{D}(x) R}
        \int_{\mathbb{S}^{n-1}_{D}} \sigma v(x + R \sigma) \dd \sigma,
\end{equation}
where $\mathbb{S}^{n-1}_{D}(x) = \{ \sigma \in \mathbb{S}^{n-1} : x +
\sigma R \in D \}$, and $\omega_{D}(x) = \abs{\mathbb{S}^{n-1}_{D}(x)}$. The
same approach is briefly discussed in \parencite{Dyson2010}. While the non-local
gradient as defined in equation~\eqref{Eqn:NonLocalGradientChemotaxis} ensures
that it is well-defined it does not satisfy no-flux boundary conditions.

Xiang studied global bifurcations of the non-local chemotaxis equation using the
global bifurcation analysis by Rabinowitz and Crandall (same approach as in
\cref{Chapter:GlobalBifPeriodic}), modified the non-local
gradient~\eqref{Eqn:NonLocalGradientChemotaxis} in 1D such that a no-flux
condition is satisfied \parencite{Xiang2013}. The construction assumed that two
domains are in contact on the boundary, then using a reflection
argument\index{reflection argument} through $x = 0, L$ Xiang obtained
\[
    \nonlocalgrad{R} v(x) \coloneqq \frac{1}{2R}
    \begin{cases}
        v(x + R) - v(R - x) &\mbox{if } 0 \leq x \leq R, \\
        v(x + R) - v(x - R) &\mbox{if } R < x < L - R \\
        v(2L-x-R) -v(x - R) &\mbox{if } L - R \leq x \leq L
    \end{cases}.
\]
A similar reflection approach is briefly discussed by Topaz \etal\
\parencite{Topaz2006}. Another class of non-local models to study species
aggregations are the so called \textit{aggregation equations}\index{aggregation
equation} \cite{Mogilner1999,fetecau2017swarm,wu2015nonlocal}. In these equations the
non-local term is the results of interactions between individuals. The
interactions can be described using a potential energy. Recently, such equations
have been studied on bounded domains \cite{fetecau2017swarm,wu2015nonlocal}. The
resulting boundary conditions are very similar to ours.

\section{Non-local no-flux boundary conditions}\label{sec:NonLocalTerm}

In this chapter we present a comprehensive theory for non-local boundary
conditions of the adhesion model in one dimension. The key idea, presented
in~\cite{Hillen2019} is to introduce a spatially dependent sampling domain
called $E(x)$. It ensures that cells measure adhesive forces inside the domain
and also on the domain boundary. Within this framework we can define naive,
no-flux, neutral, adhesive and repellent boundary conditions, each relating to
their own biological reality. We show results on local and global existence of
solutions, steady states and pattern formations.

We develop non-local boundary conditions for a one-dimensional domain $D =[0,L],
L>0$, where $L > 2R$ such that a single cell cannot touch both boundaries at the
same time. Extensions into higher dimensions, are not straight forward, and we
will discuss them later.  Let $u(x,t)$ be the density of a cell population at
location $x \in D$ and time $t$, which diffuses and adheres to each other.
\begin{equation}\label{Eqn:ArmstrongModelNoFlux}
    u_t(x, t) = d u_{xx}(x,t)
        - \alpha \lb u(x,t) \int_{-R}^{R} h(u(x + r,t)) \Omega(r) \dd r \rb_{x},
\end{equation}
where $d$ is the diffusion coefficient, $R$ the cell sensing radius, $\alpha$
the strength of the homotypic adhesions, and $h(\cdot)$ is a possibly nonlinear
function describing the nature of the adhesive force. We denote the
corresponding cell flux\index{cell flux} by
\[
    J(x, t) \coloneqq - d u_{x}(x, t) +
        \alpha u(x, t) \int_{-R}^{R} h(u(x+r, t)) \Omega(r) \dd r,
\]
which has two components the diffusive flux\index{diffusive flux}, and the
adhesive flux\index{adhesive flux}. We use
the cell flux to define boundary conditions. In a closed container, or a petry
dish, it is reasonable to assume that no cells can leave or enter through the
domain boundary, hence $J\cdot \vec{n} = 0$, where $\vec{n}$ denotes the outward
pointing normal on $\partial D$. In our one-dimensional case we have
$\vec{n}=-1$ for $x=0$, and $\vec{n}=+1$ for $x=L$. Hence we require
\begin{equation}\label{Eqn:NoFluxBoundaryConditions}
    J(0,t) = J(L,t) =0.
\end{equation}
As in~\cite{Hillen2019} we consider two cases; (1) the diffusive flux and the
adhesive flux are independent, and (2) the diffusive flux and adhesive flux
depend on each other.

\subsection{Independent fluxes\index{Independent fluxes}}\label{subsec:IndepFlux}
To satisfy the zero flux condition~\eqref{Eqn:NoFluxBoundaryConditions}, we
assume that the diffusive and adhesive fluxes independently satisfy:
\begin{equation}\label{Eqn:IndepBc}
    u_x(0, t)= u_x(L,t) = 0,\qquad
    \K[u](0) =\K[u](L)=0.
\end{equation}
The condition for the diffusive flux is easily included in the mathematical
problem formulation by restricting to the appropriate function space. For
instance, given a function space $X$, we define a boundary operator $\Bd$ and
construct a function space satisfying Neumann boundary conditions\index{Neumann
boundary conditions}.
\[
    \Bd[u] = (u^{\prime}(0), u^{\prime}(L)), \qquad
    X_{\Bd} \coloneqq X \cap \Null[\Bd].
\]
The condition for the non-local term is more challenging, and in subsequent
section we show that this leads to a spatial dependence in the integration
limits of the non-local term.

\subsection{Dependent flux}
In certain situations we can relax conditions~\eqref{Eqn:IndepBc} and allow
non-zero adhesion fluxes on the boundary \ie\
\begin{equation}\label{Eqn:BoundaryCondition3}
\begin{split}
    d u_x(0,t) &= u(0, t) \K[u](0, t) \neq 0, \\
    d u_x(L,t) &= u(L, t) \K[u](L, t) \neq 0.
\end{split}
\end{equation}
This means that the first derivative can no longer be zero on the
domain's boundary.  This condition, however, is less amenable to mathematical
analysis, as it shows the hyperbolic nature of the non-local drift term.
We need to distinguish the incoming particle flux on the boundary from the
outgoing part, i.e.\ at each location along the boundary we would need to
identify the direction of the net flux and stipulate boundary conditions only at
influx points. It is not clear how this can be done. Hence in this chapter we
focus on the independent flux case~\eqref{subsec:IndepFlux}.

In both cases the definition of appropriate boundary conditions for the
non-local operator $\K[u]$ is both a mathematical and modelling challenge.
We need to ensure that the non-local operator $\K[u]$ is well-defined
near the boundary and satisfies the no-flux boundary conditions.
To determine the near boundary behaviour we look at the
biology. When a cell encounters a boundary, it might attach to it, be repelled
by it, or form a neutral attachment.

\section{Naive boundary conditions\index{naive boundary conditions}}\label{subsec:NaiveBc}
\begin{figure}\centering
 \includegraphics[width=0.9\textwidth]{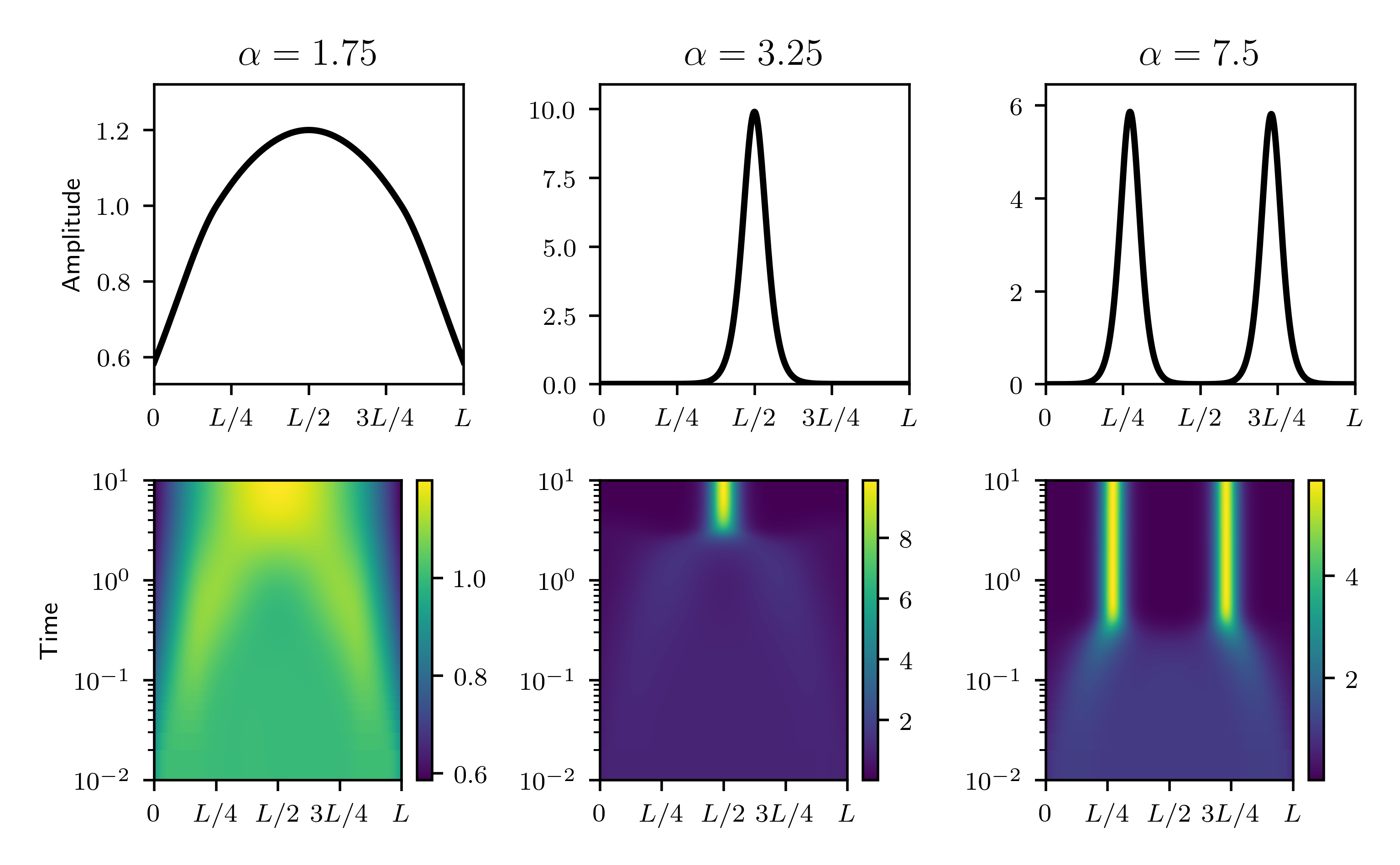}
 \caption[Naive non-local operator numerical solution]{Numerical solutions
 of equation~\eqref{Eqn:ArmstrongModelNoFlux}, with the naive sensing domain
 $E_0(x)$. In the top row we show the final solution profiles, below are the kymographs.
 (Left) $\alpha = 1.75$, (Middle) $\alpha = 3.25$, and (Right) $\alpha = 7.5$.
 The solutions feature several similarities to the periodic solutions in
 \cref{Fig:PerNumRes} (\cref{subsec:numer}). For small $\alpha$ we note
 the effect of the naive boundary conditions, repelling particles into the
 interior of the domain. For all shown $\alpha$ there is no translational
 symmetry.}\label{Fig:NaiveBcNumericalSolution}
\end{figure}

The simplest case to define sensible boundary conditions is a case used in  \parencites{Hillen2007}{Dyson2010}, where we remove any points of the
sampling domain that fall outside the domain $D$.
\[
  E_0(x) = \lcb y \in V : x + y \in D \rcb.
\]
We call this the naive case, as this does not employ any biological reasoning.
We simply restrict the integration domain.
The resulting operator is well defined for all $x\in D$, but it does not
necessarily satisfy the zero-flux condition \eqref{Eqn:IndepBc}. Indeed, this is
easily observed when $\Omega(r) = \frac{1}{2R}\frac{r}{|r|}$, and when
$\K$ is applied to a constant function $c > 0$. We observe that for $x\in[0,R)$
we have
\[
 \K[u] = \frac{1}{2R} \int_0^x - c \;\dd y + \frac{1}{2R} \int_x^{x+R} c\; \dd y
 = c \lb\frac{R-x}{2R}\rb,
\]
which does not go to zero as $x\to 0$. A similar computation can be made on the
right boundary $x=L$.  In this case, $\K$ satisfies the relaxed boundary
conditions~\eqref{Eqn:BoundaryCondition3} and the adhesive flux is pointing into
the domain. Thus cells are repelled from the boundaries, and we would expect
them to accumulate in the domain's interior.  This is nicely demonstrated by the
numerical solutions in this situation see \cref{Fig:NaiveBcNumericalSolution}.
We call boundary conditions with inward flow to be {\it repellent boundary
conditions}\index{repellent boundary conditions}.

\section{No-flux boundary conditions\index{no-flux boundary conditions}}

In this section we show how a more general sensing domain is constructed, based on biological
principles, such that the corresponding non-local operator satisfies the
independent boundary conditions~\eqref{Eqn:IndepBc}.

In~\cref{c:randomwalk} we summarize the derivation of the non-local adhesion
model from biological principles as developed in \cite{Buttenschon2017}. In the
derivation the non-local term resulted from the careful construction of the
adhesive forces between cells. A key
requisite for adhesive interaction between cells is the formation of adhesive
bonds between cells. The magnitude of the adhesive force created at distance $r$
from the cell center, is a function of the density of adhesion bonds $N_b(r)$,
available free space $f(r)$, and the distance weight function $\Omega(r)$.
In an unbounded domain the explicit form of the non-local term is
\[
    \K[u](x) = \int_{-R}^{R} N_b(x+r) f(x+r) \Omega(r) \dd r.
\]
Subsequently we focus on the number of formed adhesion bonds $N_b$, and assume
that $f\equiv 1$. We make an additional assumption that the limits
\begin{equation}\label{Eqn:Omegalimits}
    \Omega^- \coloneqq \lim_{x\to 0-} \Omega(x) \quad\ \mbox{and }\quad
        \Omega^+ \coloneqq \lim_{x\to0+} \Omega(x)
\end{equation}
exists and are non-zero. Since $\Omega$ is an odd function we have that
\[ \Omega^-=-\Omega^+.\]
As discussed in the model derivation \cref{c:randomwalk}, the number of adhesion
bonds $N_b(x+r)$ depends on the local cell population density $u(x+r)$. Hence we
write $N_b(x+r) = h(u(x+r))$ for an appropriate function $h\geq 0$. Often $h(u)$
is a linear function.  For a bounded domain, the number of adhesion bonds is no
longer homogeneous in space, since cells in the interior bind to neighboring
cells, but cells close to the boundary bind to the boundary and to other cells.
Hence the function that describes the net adhesion bonds needs to explicitly
depend on both, the location $x$ of the probing cell and the sampling location
$x+r$. Hence we use $N_b = h(x,u(x+r))$ and obtain
\begin{equation}\label{Eqn:NonLocalTermDerivationDefn}
    \K[u](x) = \int_{E(x)} h(x,u(x+r)) \Omega(r) \dd r.
\end{equation}
To simplify the subsequent discussion, we introduce the following change of
variables under the integral $y \coloneqq x + r$. Then the sampling domain
$E(x)$, in terms of $y$ is
\[
  \tilde E(x) = \{ y\in D : |x-y|\leq R\}
\]
With this the non-local term~\eqref{Eqn:NonLocalTermDerivationDefn} becomes
\begin{equation}\label{Eqn:NonLocalTermChangedVar}
    \K[u](x) = \int_{\tilde E(x)} h(x,u(y)) \Omega(y - x) \dd y.
\end{equation}
The boundary is solid. For this reason, cell protrusions cannot pass through it.
Thus, the number of adhesion bonds formed beyond the wall has to be zero.
Secondly, the solid wall may have repulsive or adhesive properties, which we
describe by two additional functions $h^0$ and $h^L$, one for each boundary.  As
a result the number of adhesion bonds for a cell located at $x\in D$ and binding
at $y\in \tilde E(x)$ is given as
\[
    h(x,u(y)) = \begin{cases}
                    H(u(y)) &\mbox{if } y \in \text{int}(D)  \\
                    h^0(x) \delta(y) &\mbox{if } y = 0 \\
                    h^L(x) \delta(L-y) &\mbox{if } y = L \\
                    0 &\mbox{else}
            \end{cases}.
\]
Notice that the function $H(u(y))$ in the first case does not explicitly depend
on $x$, since in the inner region we are in the homogeneous situation that cells
interact with other cells within the sensing range. However, the functions
$h^i(x)$ can be chosen such that the non-local term satisfies boundary
conditions~\eqref{Eqn:IndepBc}.  First consider the non-local term $\K[u](x)$
defined in equation~\eqref{Eqn:NonLocalTermChangedVar} in the interval $(0, R]$.
There, the non-local term~\eqref{Eqn:NonLocalTermChangedVar} can be decomposed
as follows
\[
    \K[u](x) = \int_{x}^{x + R} H(u(y)) \Omega(y - x) \dd y
                 + \int_{0}^{x} H(u(y)) \Omega(y - x) \dd y
                 + h^0(x) \Omega(-x).
\]
Then, we choose the function $h^0(x)$ so that the boundary
condition~\eqref{Eqn:IndepBc} is satisfied. We take the limit $x\to 0+$,
using \eqref{Eqn:Omegalimits} and obtain
\begin{equation}\label{Eqn:BoundaryConditionDetailedLeft}
    h^0(0) =\frac{-1}{\Omega^-}
        \int_{0}^{R} H(u(y)) \Omega(y - x) \dd y = \frac{1}{\Omega^+}
        \int_{0}^{R} H(u(y)) \Omega(y - x) \dd y.
\end{equation}
Similarly, on the interval $[L-R, L)$, we have the
following decomposition of the non-local term
\[
 \K[u](x) = \int_{x - R}^{x} H(u(y)) \Omega(y - x) \dd y +
            \int_{x}^{L} H(u(y)) \Omega(y - x) \dd y + h^L(x) \Omega(L-x).
\]
To satisfy the no-flux boundary condition~\eqref{Eqn:IndepBc}
at $x = L$ we consider the limit as $x\to L$ and find
\begin{equation}\label{Eqn:BoundaryConditionDetailedRight}
   h^L(L) = \frac{1}{\Omega^-}
        \int_{L-R}^{L} H(u(y)) \Omega(y - x) \dd y.
\end{equation}
It should be noted that if $\Omega^-=0=\Omega^+$ then we would need to require
that the integrals in \eqref{Eqn:BoundaryConditionDetailedLeft} and
\eqref{Eqn:BoundaryConditionDetailedRight} are zero, which will not lead to a
suitable boundary condition. Hence the case of $\Omega^-=0=\Omega^+$ cannot be
studied with the method presented here.

\begin{figure}
    \includegraphics[width=10cm]{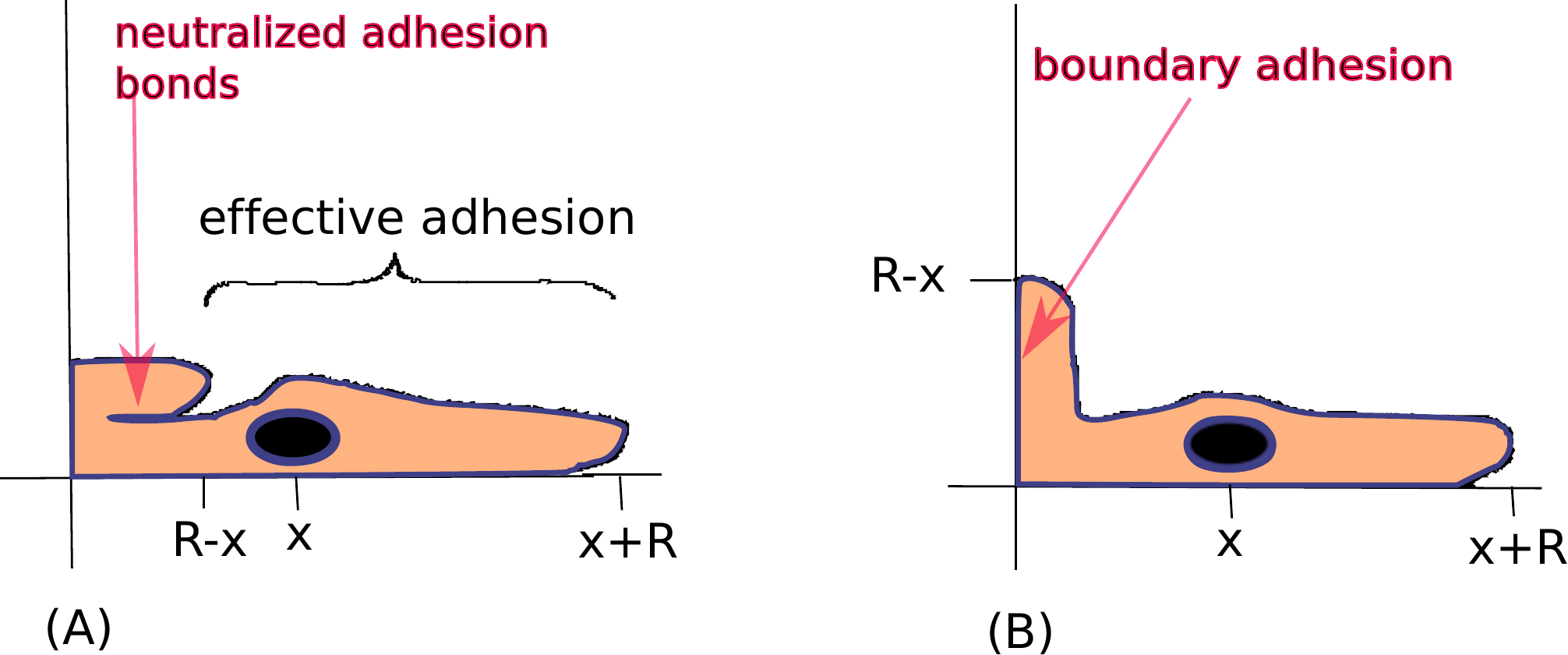}
    \caption{(A): The filopodia\index{filopodia} of cell are reflected or stopped at the boundary.
    As a result the cell starts to form adhesion bonds with itself, which are then
    not contributing to the net adhesion force. Note that only one cell is shown in
    this sketch.  (B) The weak adhesive case\index{weakly adhesive boundary
    conditions}. Cells
    make contact to the boundary in a well balanced way, such that the net flux
    is still zero.\label{fig:bioboundary1}}
\end{figure}

The boundary condition~\eqref{Eqn:IndepBc}, defines the value of the
$h^0(x)$ at $0$ and $h^L(x)$ at $L$. We still have a choice to define these
functions in the interval of size $R$ from the boundary. We will study several
choices.

Depending on the adhesive properties of the boundary a cell may either be
attracted, not-affected (neutral), or repelled. In the following, we propose
various different behaviours of the non-local term $\K[u]$ within one sensing
radius of the boundary.\\

There are many possible choices of the function $h^0(x)$ such that
condition~\eqref{Eqn:BoundaryConditionDetailedLeft} is satisfied. A natural
choice is to assume that the adhesion near the boundary follows the same
mechanism $H(u)$ as in the domain's interior. Therefore, it is natural to pick a
function $h^0(x)$ that is continuous between $x =0, R$. We choose
\begin{equation}\label{Eqn:LeftCorrectionTerm}
    h^0(x) = \frac{1}{\Omega(x)}
        \int_{0}^{R-x} H(u(y)) \Omega(y - x) \dd y,\ x \in [0,R],
\end{equation}
where we also assume that $\Omega(r)\neq 0$ for all $r\in V\backslash \{0\}$.
Thus resulting in a non-local gradient for $x \in (0, R]$ of
\begin{equation}\label{Eqn:ExmpNonLocalLeft}
    \K[u](x) = \int_{R - x}^{x + R} H(u(y)) \Omega(y - x) \dd y.
\end{equation}

We make the same natural choice for the right boundary i.e.\ we set $h^L(x)$
such that it satisfies condition~\eqref{Eqn:BoundaryConditionDetailedRight} at
$x = L$ and $h^L(L-R) = 0$. The resulting non-local operator for $x \in [L - R,
L)$ is given by
\begin{equation}\label{Eqn:ExmpNonLocalRight}
    \K[u](x) = \int_{x - R}^{2L - R - x} H(u(y)) \Omega(y - x) \dd y.
\end{equation}
Combining the boundary non-local terms~\eqref{Eqn:ExmpNonLocalLeft}
and~\eqref{Eqn:ExmpNonLocalRight}, and reverting the change of variables, we
obtain
\begin{equation}\label{NonLocalBoundary}
    \K[u](x) = \begin{cases}
        \begin{split}
            \int_{R-2x}^{R} H(u(x+r)) \Omega(r) \dd r
        \end{split}
            &\mbox{if } x \in (0, R] \\[20pt]
        \begin{split}
            \int_{-R}^{R} H(u(x+r)) \Omega(r) \dd r
        \end{split}
        &\mbox{if } x \in [R, L-R] \\[20pt]
        \begin{split}
            \int_{-R}^{2L - R - 2x}H(u(x+r)) \Omega(r) \dd r
        \end{split}
            &\mbox{if } x \in [L-R, L)\\[10pt]
            0 & \mbox{if } x=0, L
    \end{cases}.
\end{equation}
It is easy to observe that $\K$ is indeed continuous (which we will discuss in
detail in \cref{subsec:cont}). The $\Omega(r)$--terms preceding the integrals in
equation~\eqref{Eqn:LeftCorrectionTerm} cancel out for $r\neq 0$.

\begin{figure}\centering
 \includegraphics[width=0.9\textwidth]{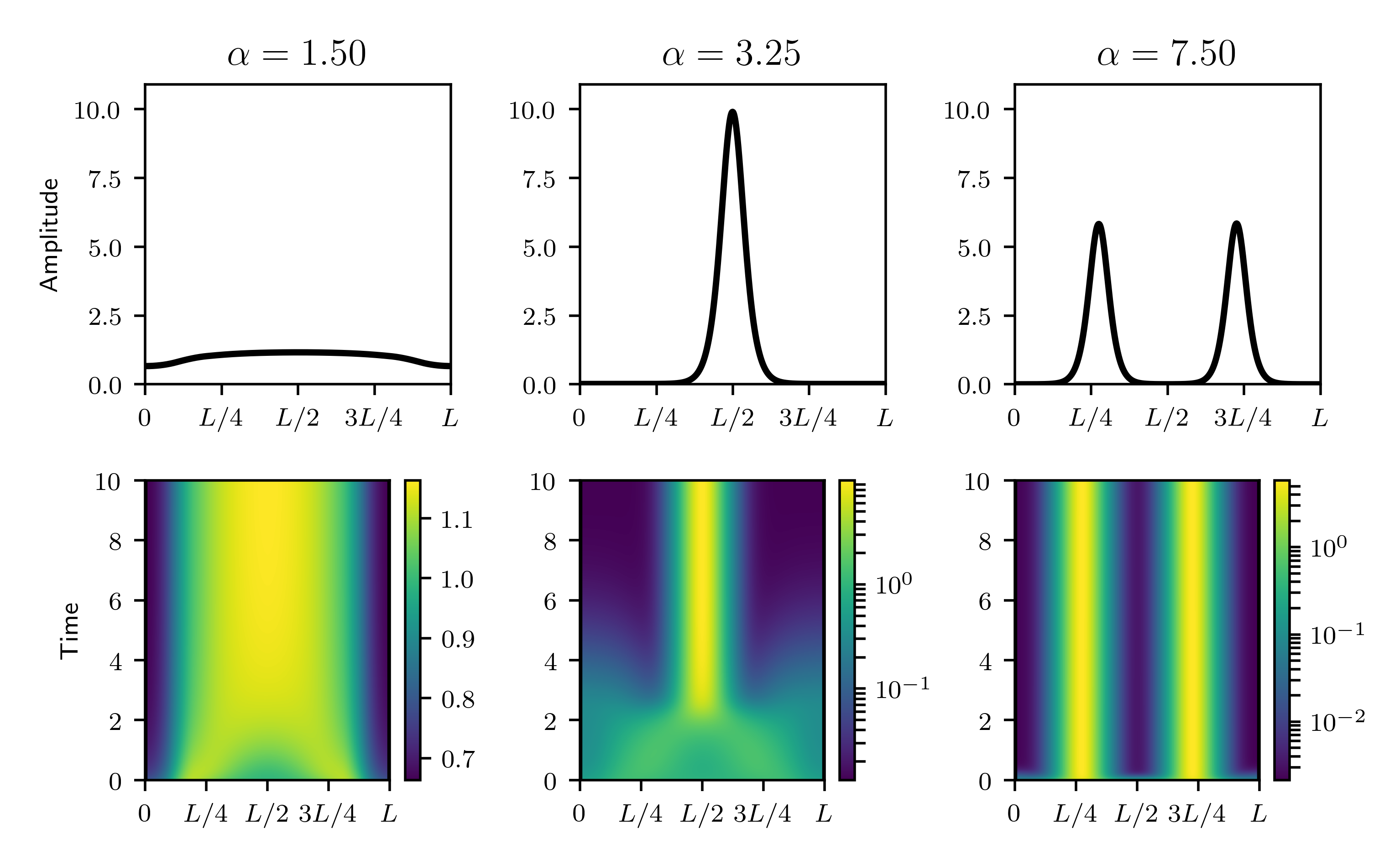}
 \caption{Numerical simulations of
 equation~\eqref{Eqn:ArmstrongModelNoFlux}, with the no-flux sensing domain
 $E_f(x)$. In the top row we show the final solution profiles, below are the
 kymographs. (Left) $\alpha = 1.75$, (Middle) $\alpha = 3.25$, and (Right) $\alpha = 7.5$.
 The solutions feature several similarities to the periodic solutions in
 \cref{Fig:PerNumRes} (\cref{subsec:numer}). For small $\alpha$ we note
 the weakly repulsive effect of the boundaries. Finally, we note that for all
 values of $\alpha$ the solutions have no translation
 symmetry.}\label{fig:noflux}
\end{figure}

This explicit form \eqref{NonLocalBoundary} allows us to write the integral
operator as integral over an effective sensing domain.\index{effective sensing
domain} We define
\begin{eqnarray*}
    f_1(x) &=& \left\{\begin{array}{ll} R-2 x& \quad\mbox{for } x\in[0,R)\\
                                       -R  & \quad\mbox{for } x\in [R,L]
    \end{array} \right., \\
    f_2(x) &=& \left\{\begin{array}{ll} R & \quad\mbox{for } x\in[0,L-R)\\
                              2L-R-2x   & \quad\mbox{for } x\in (L-R,L] \end{array}, \right.
\end{eqnarray*}
and the sampling domain\index{sampling domain}
\[
 E_f(x) \coloneqq [f_1(x), f_2(x)],
\]
and write
\begin{equation}\label{Eqn:Kf}
    \K[u]= \int_{E_f(x)} h(u(x+r)) \Omega(r) \dd r.
\end{equation}
The sensing domain for this case is shown on the right of \cref{Fig:DomainOfIntegration}.
Biologically, we interpret $E_f(x)$ as describing how the
filopodia\index{filopodia} upon hitting the domain's boundary are reflected or
stopped at the boundary. Consequently the
cell starts to form adhesion bonds with itself, which are then not contributing
to the net adhesion force (see also \cref{fig:bioboundary1}).
Numerical solutions of equation~\eqref{Eqn:ArmstrongModelNoFlux} for this
particular choice of sensing domain are shown in \cref{fig:noflux}.

\subsection{Approximate steady states for no-flux non-local term}\label{subsec:SteadyStatesRepellent}
In this section we consider a singular perturbation of the steady state of the
no-flux adhesion model for small adhesive strengths. The steady states of
equation~\eqref{Eqn:ArmstrongModelNoFlux} with repellent non-local term $\K[u]$,
are given by the solutions of the following equation.
\begin{subequations}
\begin{equation}\label{Eqn:SteadyStateEquationNoFlux}
    u_{xx}(x,t) - \alpha \lb
        u(x,t) \int_{E_f(x)} u(x + r) \Omega(r) \dd r \rb_{x} = 0,
\end{equation}
where $E_f(x)$ is no-flux sensing domain.
As equation~\eqref{Eqn:ArmstrongModelNoFlux} exhibits mass conservation we
impose the following mass constraint\index{mass constraint} on the solutions of
equation~\eqref{Eqn:SteadyStateEquationNoFlux}.
\begin{equation}\label{Eqn:IntegralConstraintNoFlux}
    \Avg[u] = \mean{u},
\end{equation}
where $\R \ni \mean{u} > 0$, is the mass per unit length of population
$u(x)$ in $D$.
\end{subequations}
To be able to easily carry out the subsequent asymptotic
expansion\index{asymptotic expansion} we assume that the function $h(\cdot)$
under the integral in
equation~\eqref{Eqn:SteadyStateEquationNoFlux} is linear (\ie\ $h(u) = u$).

When $\alpha = 0$ the boundary conditions are given by the classical Neumann
boundary conditions. It is then easy to see
that equation~\eqref{Eqn:ArmstrongModelNoFlux} admits a constant steady state
solution $u \equiv \mean{u}$. When $\alpha \neq 0$ the situation is much
more complicated.

Here we approximate the ground steady state of
equation~\eqref{Eqn:SteadyStateEquationNoFlux} by using an asymptotic expansion
for small values of $\alpha$. Here we set $R = 1$. In the following, we assume
that $\alpha = \epsilon$. Then we consider the following asymptotic
expansion\index{asymptotic expansion}
\begin{equation}\label{Eqn:AsymptoticExpansion}
    u(x) = u_0(x) + \epsilon u_1(x) + \epsilon^2 u_2(x) + \BigOh{\epsilon^3}.
\end{equation}
Substituting into equation~\eqref{Eqn:SteadyStateEquationNoFlux}, we obtain
\[
\begin{split}
    0 &\approx
    \lb u_0(x) + \epsilon u_1(x) +\epsilon^2 u_2(x) \rb_{xx}  \\
        &- \epsilon
        \lb \lb u_0 + \epsilon u_1 \rb \int_{-1}^{1} \lb u_0 + \epsilon u_1 +
        \epsilon^2 u_2 \rb (x + r) \chi_{E_f(x)}(r) \Omega(r) \dd r \rb_{x}.
\end{split}
\]
Separating the scales of $\epsilon$ we obtain for the zeroth order equation
$\lb u_0 \rb_{xx} = 0$. Due to the Neumann boundary conditions, and the
mass constraint~\eqref{Eqn:IntegralConstraintNoFlux} we have $u_0 \equiv
\mean{u}$. The first order equation is given by
\begin{equation}\label{Eqn:FirstOrder}
    \lb u_1 \rb_{xx} - \lb u_0 \int_{-1}^{1} u_0(x + r) \chi_{E_f(x)}(r)
        \Omega(r) \dd r \rb_{x} = 0.
\end{equation}
Equation~\eqref{Eqn:FirstOrder} requires the solvability condition $\Avg[u_1] =
0$. Using the properties of the function $\Omega(r) = \frac{r}{|r|}\omega(r)$,
we rewrite equation~\eqref{Eqn:FirstOrder} as a piecewise condition
\[
 \lb u_1 \rb_{xx} = \mean{u}^2 \begin{cases}
        2 \omega(1 - 2x)       &\mbox{for } x \in [0, 1/2] \\
        -2 \omega(1 - 2x)      &\mbox{for } x \in [1/2, 1] \\
        0                      &\mbox{for } x \in [1, L - 1] \\
        -2 \omega(2L - 2x - 1) &\mbox{for } x \in [L - 1, L - 1/2] \\
        2 \omega(2L - 2x - 1)  &\mbox{for } x \in [L - 1/2, L]
    \end{cases}.
\]
We solve this differential equation by integrating and after some algebra we
obtain:
\[
    u_1(x) =
    \begin{cases}
    \begin{split}
     u(0) + \mean{u}^2 \int_{0}^{x} \int_{1-2s}^{1} \omega(r) \dd r \dd s
    \end{split}
            &\mbox{for } x \in [0, 1/2] \\[15pt]
    \begin{split}
     A    - \mean{u}^2 \int_{x}^{1} \int_{-1}^{1-2s} \omega(r) \dd r \dd s
    \end{split}
            &\mbox{for } x \in [1/2, 1] \\[15pt]
        A &\mbox{for } x \in [1, L -1] \\[15pt]
    \begin{split}
     A    - \mean{u}^2 \int_{L-1}^{x} \int_{2L-2s-1}^{1} \omega(r) \dd r \dd s
    \end{split}
            &\mbox{for } x \in [L-1, L-1/2] \\[15pt]
    \begin{split}
     u(L) + \mean{u}^2 \int_{x}^{L} \int_{-1}^{2L-2s-1} \omega(r) \dd r \dd s
    \end{split}
            &\mbox{for } x \in [L - 1/2, L]
    \end{cases}
\]
where
\[
 A = \frac{\mean{u}^2}{L} \int_{0}^{1} r \omega(r) \dd r,
\]
and
\[
 u(0) = u(L) = \mean{u}^2 \lsb \frac{1 - L}{L} \rsb \int_{0}^{1} r \omega(r) \dd r.
\]

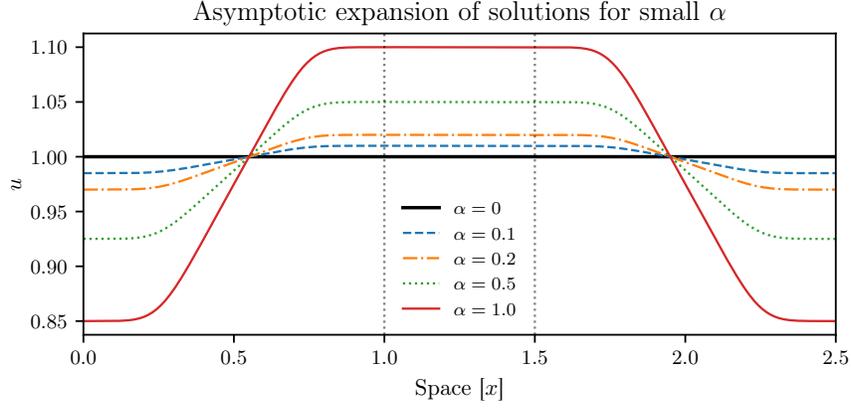
\begin{figure}\centering
 \resizebox{0.95\textwidth}{!}{%
    \input{AsymptoticExpansion.pgf}
 }
 \caption[The asymptotic expansion of the steady states with the repellent
 non-local operator]
 {The solution $u(x)$ given in
 equation~\eqref{Eqn:AsymptoticExpansion} for various values of $\epsilon$ or
 $\alpha$. The vertical black dotted lines denote the boundaries of the areas that are
 within one sensing radius of the domain boundary, at $x = R$ and $x = L-R$.
 }\label{Fig:AsymptoticExpansion}
\end{figure}

For a visual depiction of the asymptotic
expansion~\eqref{Eqn:AsymptoticExpansion} see \cref{Fig:AsymptoticExpansion}.
The domain boundary in this case is weakly repellent. At steady-state, cells
experience adhesive forces from the domain interior and less from the domain
boundary hence leading to local minima at the boundary.

\section{General Sampling Domain}\label{subsec:cont}
We can generalize the above approach and define a general class of suitable
sampling domains. These will include the naive case and the no-flux case, and
further cases for weak adhesion and repulsion as we discuss later. To  define a
general suitable domain of integration for the non-local term we define.

\begin{definition}[Sensing domain\index{sensing domain} \cite{Hillen2019}]\label{def:slice}

 \begin{enumerate}
  \item Two continuous functions $f_{1,2}:D\to \R$ define a suitable sampling
      domain
      \[E(x) = [f_1(x), f_2(x)] \],
      if they satisfy
 \begin{enumerate}
  \item $-R\leq f_1(x) \leq f_2(x)\leq R $ for all $x\in[0,L]$.
  \item $f_1(x) = -R$ for $x\in[R, L]$
  \item $f_2(x) = R$ for $x\in[0,L-R]$.
  \item $f_1(x)$, and $f_2(x)$ are non-increasing and have uniformly bounded one-sided derivatives.
 \end{enumerate}
  \item A suitable sampling domain $E(x)$ satisfies the second
      condition~\eqref{Eqn:IndepBc} if in addition
 \begin{enumerate}[resume]
  \item $f_1(0) = R, \quad \mbox{and} \quad  f_2(L) = -R$. 
 \end{enumerate}
 \end{enumerate}
\end{definition}

\begin{figure}\centering
 \begin{minipage}[t]{0.48\textwidth}
  \resizebox{0.95\textwidth}{!}{%
    \input{SamplingDomainNaive.tikz}
  }
 \end{minipage}\hfill
 \begin{minipage}[t]{0.48\textwidth}
  \resizebox{0.95\textwidth}{!}{%
    \input{SamplingDomainNoFlux.tikz}
  }
 \end{minipage}
 \caption[Varying sensing domains]{Two examples of the spaces $D \times V$, with
    the spatial domain on the $x$-axis and the suitable sensing domain\index{sensing domain} $E(x)$ on
    the $y$-axis.  The shaded region is the set $\{ (x,y): x\in D, y \in E(x)\}$
    and a sample sensing domain $E(x)$ of thickness $\dd x$ is shown in the
    darker grey.  Left: Here the sensing domain is of naive type.  Right: Here
    the sensing domain is of no-flux type.}\label{Fig:DomainOfIntegration}
\end{figure}
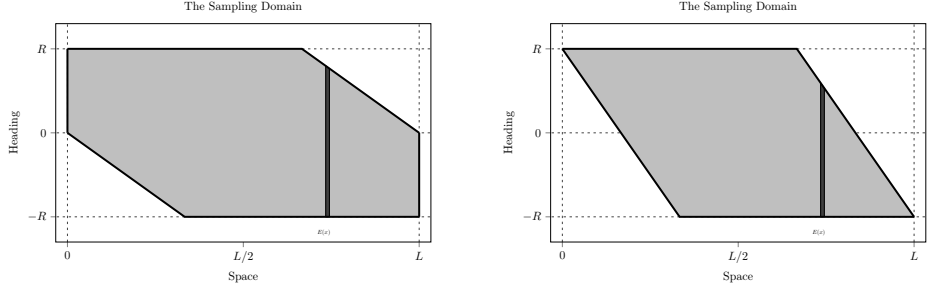

Two examples of suitable sampling domains are shown in
\cref{Fig:DomainOfIntegration}.  We view the non-local operator as function on
$L^p(D)$.
\begin{equation}\label{Eqn:Kchi}
 \begin{split}
  \K &: D \times L^p(D) \mapsto \R, \ p\geq 1. \\
  \K[u(x)](x) &\coloneqq \int_{E(x)} h(u(x+r)) \Omega(r) \dd r,
 \end{split}
\end{equation}
which satisfies the following assumptions:

\begin{enumerate}[label=\textbf{(A\arabic*)},ref=(A\arabic*), leftmargin=*,
        labelindent=\parindent]
    \item\label{OmegaAssumptionNoFlux:1}
        $\Omega(r) = \frac{r}{\abs{r}} \omega(r)$,
        where $\omega(r) = \omega(-r)$,
    \item\label{OmegaAssumptionNoFlux:2} $\omega(r) \geq 0$, $\omega(R)=0$,
    \item\label{OmegaAssumptionNoFlux:3}
        $\omega \in L^{1}(V) \cap L^{\infty}(V)$, and  $\norm{\omega}_{L^1([0,R])} = \ifrac{1}{2}$.
    \item\label{slicecondition} The sampling domain $E(x)$ is suitable according
        to \cref{def:slice}.
    \item\label{OmegaAssumptionNoFlux:4}
        The adhesion function $h(u)$ is linearly bounded\index{linearly bounded} and differentiable and
    \[ |h(u)| \leq k_1 (1+|u|), \quad \mbox{and} \quad |h'(u)|\leq k_2, \quad k_1, k_2 >0.\]
\end{enumerate}

\subsection{Set Convergence\index{set convergence}}\label{subsec:SetConvergence}
The spatial dependence in the integration limits, can equivalently, be viewed as
introducing an indicator function under the integral. To establish the
continuity and differentiability of $\K[u]$ we recall the calculus of indicator
functions in general metric spaces. Let $(V, \mathcal{A}, \mu)$ be a
sigma-finite measure space. Then let $A, B \in \mathcal{A}$, and we define the
symmetric difference\index{symmetric difference} of $A$ and $B$ by
\[ A \symdif B = (A \cup B) \setminus (A \cap B).  \]
The function
\begin{equation}\label{Eqn:PseudoMetric}
  d : \mathcal{A} \times \mathcal{A} \mapsto \R,
\end{equation}
defined by
\[
    d(A, B) = \mu(A \symdif B) \quad \mbox{for } A, B \in \mathcal{A},
\]
is a pseudometric\index{pseudometric} and is known as the Fr\'echet-Nikodym
metric\index{Fr\'echet-Nikodym metric}
\parencite{Bogachev2007}, allowing us to view any sigma-finite measure space as
a pseudometric space.  We can turn this into a metric by introducing the
following equivalence relation
\begin{equation}\label{Eqn:SetEquivalenceRelation}
    X \sim Y \quad\iff\quad \mu(X \symdif Y) = 0.
\end{equation}
It is easy to see, that $A \sim B$ if and only if they differ by a set of
measure zero. If we denote the set of equivalence classes of this relation by
$\tilde{\mathcal{A}} = \mathcal{A} / \mu$ then the
function~\eqref{Eqn:PseudoMetric} is extended by setting $d(\tilde{A},
\tilde{B}) = d(A, B)$ where $\tilde{A}, \tilde{B} \in \mathcal{A} / \mu$. This
turns $\mathcal{A} / \mu$ into a complete metric space
$\tilde{\mathcal{A}}$~\cite[Theorem 1.12.6]{Bogachev2007}.

\begin{lemma}\label{Lemma:SetConvergence}
 Let $(V, \mathcal{A}, \mu)$ be a measure space, with
 pseudometric~\eqref{Eqn:PseudoMetric} and $\tilde{\mathcal{A}}$ the
    corresponding metric space of equivalence classes\index{equivalence class}
 of~\eqref{Eqn:SetEquivalenceRelation}.  Suppose that $E(x)$ is a suitable
 sample slice as defined in Definition~\ref{def:slice}. Then $E(x)$ is continuous
 in $\tilde{\mathcal{A}}$ in the following sense, if $\seq{x_n} \subset \R$ such
 that $x_n \to x$, then $\mu(E(x_n) \symdif E(x)) \to 0$.
\end{lemma}
\begin{proof}
We consider a sequence $x_n \to x \in D$ and define the sets
\[
    E_n \coloneqq [f_1(x_n), f_2(x_n)]
\]
The symmetric difference of $E(x)$ and $E_n$ is
\[
\begin{split}
    E_n \symdif E(x) = &\left\{ y \in V : \min(f_1(x), f_1(x_n)) \leq y \leq
        \max(f_1(x), f_1(x_n)), \right. \\
        &\left. \min(f_2(x), f_2(x_n)) \leq y \leq \max(f_2(x), f_2(x_n)) \right\}
\end{split}
\]
whose measure is computed to be
\[
    \mu(E_n \symdif E(x)) = \abs{f_1(x_n) - f_1(x)} + \abs{f_2(x_n) - f_2(x)}.
\]
Since both $f_{1,2}$ are continuous, we have that
\[
    \mu(E_n \symdif E(x)) \to 0 \quad\mbox{as}\quad x_n \to x.
\]
\end{proof}
Next we collect some properties of the indicator functions of symmetric
differences. Suppose that $A, B \in \mathcal{A}$, then
\begin{eqnarray*}
    \chi_{A \symdif B}(r) &=& \chi_{A}(r) + \chi_{B}(r) - 2 \chi_{A \cap B}(r),\\
    \chi_{A \cap B}(r) &=& \chi_{A}(r) \chi_{B}(r) = \min(\chi_{A}(r), \chi_{B}(r)),\\
    \chi_{A \symdif B}(r) &=& \abs{\chi_{A}(r) - \chi_{B}(r)}.
\end{eqnarray*}
Having established the continuity of the sampling domains $E(x)$ we can consider
the continuity of the integral operator \eqref{Eqn:Kchi}.

\subsection{Continuity near the boundary}
In this section, we explore continuity of the no-flux non-local operator defined
in equation~\eqref{Eqn:Kchi}. To distinguish between the norms over the sensing
domain $V = [-R, R]$ and the spatial domain $D$ we introduce norm notation that
indicates the set over which the norm is taken. For example, in this subsection
the $L^p$ norm over that set $A$ is denoted as
\[
    \abs{u}_{p, A} \coloneqq \lb \int_{A} \abs{u(x)}^{p} \dd x \rb^{1/p}.
\]
\begin{theorem}\label{th:continuity}
    Let $p \geq 1$, $u \in L^p(D)$, and assume
    that~\ref{OmegaAssumptionNoFlux:1}--\ref{OmegaAssumptionNoFlux:4} are
    satisfied.
    \begin{enumerate}
    \item The map $\Bigl(r\mapsto \chi_{E(x)} (r) H(u(x+r)) \Omega(r)\Bigr)\in
        L^p(D) $ and there exists a constant $C>0$ such that
        \[
            |\chi_{E(x)} (r) H(u(x+r))|_{p,V} \leq C(1 + |u|_{p,D}).
        \]
    \item The map $\Bigl(x\mapsto \K[u](x)\Bigr) \in L^p(D)$ and there exists a constant $C>0$ such that
        \[
            |\K[u]|_{p,D} \leq C (1+|u|_{p,D})
        \]
    \item $\K: L^p(D)\mapsto L^p(D)$ is continuous.
    \item If $u\in \C^0(D)$ then the map $\Bigl(x\mapsto \K[u](x)\Bigr) \in \C^0(D)$.
    \item The map $\Bigl(x\mapsto \K[u](x)\Bigr) \in \C^0(D)$.
\end{enumerate}
\end{theorem}
\begin{proof}
    Fix $x \in D$. Then
    \begin{eqnarray*}
        \int_{V} \abs{\chi_{E(x)}(r) H(u(x + r)) \Omega(r) }^p \dd r
            &\leq& |\Omega|_{\infty}^p \int_{D} \abs{H(u(x))}^p \dd x\\
            &\leq& C k_1^p \int_{D} (1+|u(x)| )^p \dd x\\
            &\leq& C (1+|u|_{p,D})^p.
    \end{eqnarray*}
    which proves item (1).
    To prove item (2) we use item (1) and estimate
    \begin{eqnarray*}
    |\K[u]|_{p,D}^p &= & \int_D \left|\int_V \chi_{E(x)(r)} H(u(x+r)) \Omega(r)
        \dd r \right|^p \dd x \\
    &\leq& \int_D \int_V \Bigl( C (1+|u|_{p,D})|\Omega|_{\infty,V}\Bigr)^p \dd r \dd x\\
    &\leq & \Bigl( C (1+|u|_{p,D}) |\Omega|_{\infty,V} \Bigr)^p |V| |D|\\
    &\leq & \Bigl( C (1+|u|_{p,D}) \Bigr)^p.
    \end{eqnarray*}
To proof item (3) we consider a sequence $u_n\to u$ in $L^p(D)$ and use the
Lipschitz continuity of $H$ to estimate
\begin{eqnarray*}
|\K[u_n]-\K[u]|_{p,D} &\leq & \left|\int_D \chi_{E(x)}(r) \Bigl( H(u_n(x+r)) -
    H(u(x+r))\Bigr) \Omega(r) \dd r  \right|_{p,D}\\
&\leq & k_2 \left|\int_D \chi_{E(x)}(r) | u_n(x+r) - u(x+r)| \Omega(r) \dd r  \right|_{p,D}\\
&\leq & C |u_n-u|_{p,D}
\end{eqnarray*}
To show item (4) we consider $u\in \C^0(D)$ and a sequence $\seq{x_n} \subset D$
    such that $x_n \to x \in D$. Then
    \begin{align*}
        \K[u](x) - \K[u](x_n) &= \int_{V} \lcb \chi_{E(x)}(r) H(u(x+r)) - \chi_{E(x_n)}
            H(u(x_n + r)) \rcb \Omega(r) \dd r \\
                &= \int_{V} \lcb \chi_{E(x)}(r) H(u(x+r)) - \chi_{E(x_n)} H(u(x_n + r))
                \right. \\
                &\quad+ \chi_{E(x) \cap E(x_n)}(r) H(u(x_n + r)) - \chi_{E(x) \cap E(x_n)}(r) H(u(x_n + r)) \\
                &\quad+ \left. \chi_{E(x) \cap E(x_n)}(r) H(u(x + r)) - \chi_{E(x) \cap E(x_n)}(r) H(u(x + r))
            \rcb \Omega(r) \dd r  \\
        &= \int_{V} \lcb H(u(x + r)) \chi_{E(x)}(r) \lb \chi_{E(x)}(r) - \chi_{E(x_n)}(r) \rb \right. \\
                &\quad+ H(u(x_n + r)) \chi_{E(x_n)}(r) \lb \chi_{E(x_n)}(r) - \chi_{E(x)}(r) \rb \\
                &\quad+ \left. \chi_{E(x) \cap E(x_n)}(r) \lb H( u(x + r)) - H(u(x_n + r)) \rb
            \rcb \Omega(r) \dd r.
    \end{align*}
    Then, we obtain
    \begin{align*}
        \abs{\K[u](x) - \K[u](x_n)} \leq
        k_1(1+|u|_\infty) \abs{\Omega}_{\infty} &\lsb 2
        \int_{V} \abs*{\chi_{E(x)}(r) - \chi_{E(x_n)}(r)} \dd r \right. \\
            & \left.+ \int_{E(x) \cap E(x_n)} \abs{H(u(x+r)) - H(u(x_n +r))} \dd r \rsb.
    \end{align*}
    Then, the integrals
    \begin{align*}
        \int_{V} \abs{\chi_{E(x)}(r) - \chi_{E(x_n)}(r)} \dd r =
            \int_{E(x) \symdif E(x_n)} \dd r = \mu(E(x) \symdif E(x_n)).
    \end{align*}
    By \cref{Lemma:SetConvergence}, we have that $\mu(E(x) \symdif E(x_n)) \to 0$
    as $n \to \infty$. The second integral converges by the continuity of $H$ and
    $u$.

    Finally, to show item (5) we approximate $u\in L^p(D)$ with functions
    ${u_n \in \C^0(D)}$, $u_n\to u$ in $L^p(D)$. Then
    \begin{align*}
        \K[u](x) - \K[u](x_n) &= \K[u](x)  - \K[u_m](x) \\
                            &\qquad + \K[u_m](x_n) - \K[u](x_n)
                            + \K[u_m](x)   - \K[u_m](x_n) \\
            &\leq 2 \abs{u_m - u}_{p} \abs{\Omega}_{\infty, V}
                  + \abs{\K[u_m](x) - \K[u_m](x_n)}.
    \end{align*}
    The first term vanished due to the density of the smooth functions in
    $L^p(D)$ and item (3). The second terms vanished by item (4).
\end{proof}

In summary, we have shown that whenever $u \in L^p(D)$ and the sensing domain
$E(x)$ is continuous in the sense of \cref{Lemma:SetConvergence} we have that
the no-flux non-local term in equation~\eqref{Eqn:Kchi} is
continuous.

\subsection{Differentiation near the boundary}
The integral limits of the non-local operator $\K$ are spatially dependent,
hence taking the spatial derivative of the non-local term (as done in the model)
requires the differential of indicator functions, i.e.\ we consider a weak
formulation.

\begin{lemma}\label{Lemma:NonLocalDerivative}
In the sense of distributional derivatives we find
\begin{eqnarray}\label{dist-deriv}
    \K[u]'(x) &=& \frac{\dd}{\dd x}\left(\int_{E(x)} H(u(x+r)) \Omega(r) \dd r\right)\nonumber\\
              &=&  \int_{E(x)} \frac{\dd}{\dd x}H(u(x + r)) \Omega(r) \dd r\nonumber\\
              && + f_2^\prime(x) H(u(x + f_2(x)))\Omega(f_2(x))
                  - f_1^\prime(x) H(u(x + f_1(x)))\Omega(f_1(x)).
\end{eqnarray}
\end{lemma}
\begin{proof}
 First, we use interval notation to rewrite the indicator function as
 \[
     \chi_{E(x)}(r) = \chi_{(-R, f_2(x))}(r) \chi_{(f_1(x), R)}(r) =
     \mathcal{H}(r-f_1(x)\mathcal{H}(f_2(x)-r),\ r\in[-R,R].
 \]
 We use the Heaviside function $\mathcal{H}$ and the Dirac-delta distribution $\delta$ to write
 \begin{eqnarray*}
     \frac{\dd}{\dd x}\chi_{E(x)}(r) &=&
     \frac{\dd}{\dd x} \lcb \mathcal{H}(r - f_1(x)) \mathcal{H}(f_2(x) - r) \rcb,\\
        &=& \mathcal{H}(r - f_1(x)) \delta(f_2(x) - r) f_2^\prime(x)
         - \mathcal{H}(f_2(x)-r) \delta(r - f_1(x)) f_1^\prime(x).
 \end{eqnarray*}
 Then
 \begin{eqnarray*}
    &&  \int_{V} \frac{\dd}{\dd x} \chi_{E(x)} H(u(x+r)) \Omega(r) \dd r
        \\
     &=&    \int_{V} \mathcal{H}(r - f_1(x)) \delta(f_2(x) - r) f_2^\prime(x) H(u(x+r))
             \Omega(r) \dd r \\
         &&- \int_{V} \mathcal{H}(f_2(x)-r) \delta(r - f_1(x)) f_1^\prime(x) H(u(x+r))
             \Omega(r) \dd r \\
         &=& \mathcal{H}(f_2(x)-f_1(x))\Bigl(f_2^\prime(x)H(u(x+f_2(x))) \Omega(f_2(x))
             - f_1^\prime(x)H(u(x+f_1(x))) \Omega(f_1(x))\Bigr),\\
         &=& f_2^\prime(x)H(u(x+f_2(x))) \Omega(f_2(x))
             - f_1^\prime(x)H(u(x+f_1(x))) \Omega(f_1(x)).
 \end{eqnarray*}
 Note that $f_2(x)\geq f_1(x)$ implies $\mathcal{H}(f_2(x)-f_1(x)) =1$.
 Together with the product rule, we obtain (\ref{dist-deriv}).
\end{proof}
Next we study the continuity properties of the derivative $\K[u]'$ and we use
short-hand notations for the two parts of the derivative:
\begin{eqnarray*}
 \K[u]'(x)    &=& D_1\K[u](x) + D_2\K[u](x) \\
 D_1 \K[u](x) &=& \int_{E(x)} \frac{\dd}{\dd x}H(u(x + r)) \Omega(r) \dd r\\
 D_2\K[u](x)  &=& f_2^\prime(x) H(u(x + f_2(x)))\Omega(f_2(x))
                - f_1^\prime(x) H(u(x + f_1(x)))\Omega(f_1(x)),
\end{eqnarray*}
and study the two terms separately.
The first term can be treated the same way as $\K[u]$ in \cref{th:continuity}
\begin{corollary}\label{cor:continuity}
Let $p \geq 1$, $u \in W^{1,p}(D)$, and assume that~\ref{OmegaAssumptionNoFlux:1}--\ref{OmegaAssumptionNoFlux:4} are satisfied.
\begin{enumerate}
    \item The map $\Bigl(r\mapsto \chi_{E(x)} (r) \frac{\dd}{\dd x} H(u(x+r))
        \Omega(r)\Bigr)\in L^p(D) $ and there exists a constant $C>0$ such
        that
        \[
            \left|\chi_{E(x)} (r) \frac{\dd}{\dd x} H(u(x+r))\right|_{p,V} \leq C(1 + |u|_{W^{1,p}}).
        \]
    \item The map $\Bigl(x\mapsto D_1\K[u](x)\Bigr) \in L^p(D)$ and there exists
        a constant $C>0$ such that
        \[
            |D_1\K[u]|_{p,D} \leq C (1+|u|_{W^{1,p}})
        \]
    \item $D_1\K: W^{1,p}(D)\to L^p(D)$ is continuous.
    \item The map $\Bigl(x\mapsto \K[u](x)\Bigr) \in C^0(D)$.
\end{enumerate}
\end{corollary}
\begin{proof}
    The proof follows directly from \cref{th:continuity}.
\end{proof}

Let us consider $D_2\K$.
\begin{lemma}\label{Lemma:PropertiesK} Let $p \geq 1$, $u \in W^{1,p}(D)$, and
 assume that~\ref{OmegaAssumptionNoFlux:1}--\ref{OmegaAssumptionNoFlux:4} are satisfied.
The map
\[ x\mapsto D_2\K[u](x) \]
is continuous for $x\in [0,L]$.
\end{lemma}
\begin{proof} Since the sampling slice $E(x) = [f_1(x), f_2(x)]$ is defined by
    suitable functions $f_1(x), f_2(x)$ given by \cref{def:slice}, the only
    possible points of discontinuity of $D_2\K[u](x)$ are $x=R$ and $x=L-R$. Let
    us compute
    \begin{eqnarray*}
        \lim_{x\to R^-} D_2\K[u](x) &=& -\lim_{x\to R^-} f_1'(x) H(u(u-R)) \Omega(-R) \\
        \lim_{x\to R^+} D_2\K[u](x) &=& 0,
    \end{eqnarray*}
    which is continuous, since we assumed $f'_1$ is bounded and $\Omega(\pm R)=0$.
    A limit for $x\to L-R$ leads to the condition $\Omega(R)=0$.
\end{proof}

Hence the weak distributional derivative\index{weak distributional derivative}
$\K[u]'(x)$ is in fact a classical derivative: 
\begin{corollary}\label{cor:K_C1}
Let $p \geq 1$, $u \in W^{1,p}(D)$, and assume
that~\ref{OmegaAssumptionNoFlux:1}--\ref{OmegaAssumptionNoFlux:4} are satisfied, then
$\K[u]'(x)$ is continuous in $x\in [0,L]$. Hence $\K[u]'(x)$ is a classical derivative.
\end{corollary}

\section{Local and global existence\index{local existence}\index{global
existence}}

Having established the regularity of the non-local operator $\K[u]$ in
\cref{cor:K_C1} we are ready to state a result that ensures the global existence
of solutions of equation~\eqref{Eqn:ArmstrongModelNoFlux}.

\begin{theorem}[Theorem from \cite{Hillen2019}]\label{Thm:GlobalExistence}
 Assume~\ref{OmegaAssumptionNoFlux:1}--\ref{OmegaAssumptionNoFlux:4}, then
 equation~\eqref{Eqn:ArmstrongModelNoFlux} has a unique solution in
 \[
    u \in \C^0\lb[0,\infty), H^2(0,L)\rb.
 \]
\end{theorem}
In the case of independent fluxes the proof of this theorem is given
in~\cite{Hillen2019}. The generalization of the existence result to the case of
dependent fluxes is an open problem.

\section{Neutral boundary conditions\index{neutral boundary conditions}}\label{subsec:neutralbc}

A key feature of standard partial differential equations with homogeneous
no-flux boundary conditions is the existence of constant steady states. For the
non-local adhesion model \eqref{Eqn:ArmstrongModelNoFlux} a constant steady
state $\bar u$ must satisfy $\K[\bar u]=0$.  However, looking at the non-local
operators that we defined above: naive case, and weakly repellent case, then
neither of these admits a constant steady state. To allow for constant
steady states we need to modify the non-local operator at the boundary to
compensate for the repellent flux near the boundary.

We start with the non-local operator together with the no-flux sensing domain,
and we add an additional correction function $c(x)$. We require that
$c(0) =c(L)= 0$ such that $\K[u]$ still satisfies boundary
condition~\eqref{Eqn:IndepBc}.
\begin{equation}
    \K[u](x) \coloneqq \int_{E_f(x)}H(u(x + r)) \Omega(r) \dd r - c(x).
\end{equation}
We consider a constant $\bar u>0$ and we want to choose $c(x)$ such that
$\K[\bar u](x) = 0$ for any $x \in D$, where $c$ is only non-zero near the boundary.
Hence we set
\begin{equation}
    c(x) \coloneqq \begin{cases}
                c^0(x)          &\mbox{if } x \in [0, R) \\
                0               &\mbox{if } x \in [R, L - R] \\
                c^L(x)          &\mbox{if } x \in (L - R, L]
           \end{cases},
\end{equation}
To determine $c^0(x)$, we compute $\K[\bar{u}]$ in the boundary regions
and solve for $c^0(x)$. We find that
\[
    c^0(x) \coloneqq \begin{cases}
             \begin{split}
                H(\bar u) \int_{R-2x}^{R} \omega(r) \dd r
             \end{split}
                 &\mbox{if } x \in [0, R/2] \\[20pt]
             \begin{split}
                H(\bar u) \int_{2x-R}^{R} \omega(r) \dd r
             \end{split}
                 &\mbox{if } x \in [R/2,R]
            \end{cases}.
\]
To determine $c^L(x)$ we repeat this process on the right boundary, and find that
\[
    c^L(x) \coloneqq \begin{cases}
             \begin{split}
                - H(\bar u) \int_{2L - 2x - R}^{R} \omega(r) \dd r
             \end{split}
             &\mbox{if } x \in [L - R, L - R/2] \\[20pt]
             \begin{split}
                - H(\bar u) \int_{-R}^{2L - 2x - R} \omega(r) \dd r
             \end{split}
                &\mbox{if } x \in [L - R/2, L].
            \end{cases}
\]

With that definition of the function $c(x)$, we can rewrite the neutral
repellent non-local operator $\K[u]$ as
\begin{equation}\label{Eqn:NonLocalBCNeutral}
    \K_n[u](x) \coloneqq \int_{E_f(x)} \lsb H(u(x + r)) -  H(\bar u)  \rsb
        \Omega(r) \dd r.
\end{equation}
Since this operator respects constant solutions, we call it the {\it neutral
version\/}\index{neutral version} of the repellent no-flux operator.

\begin{figure}\centering
 \includegraphics[width=0.9\textwidth]{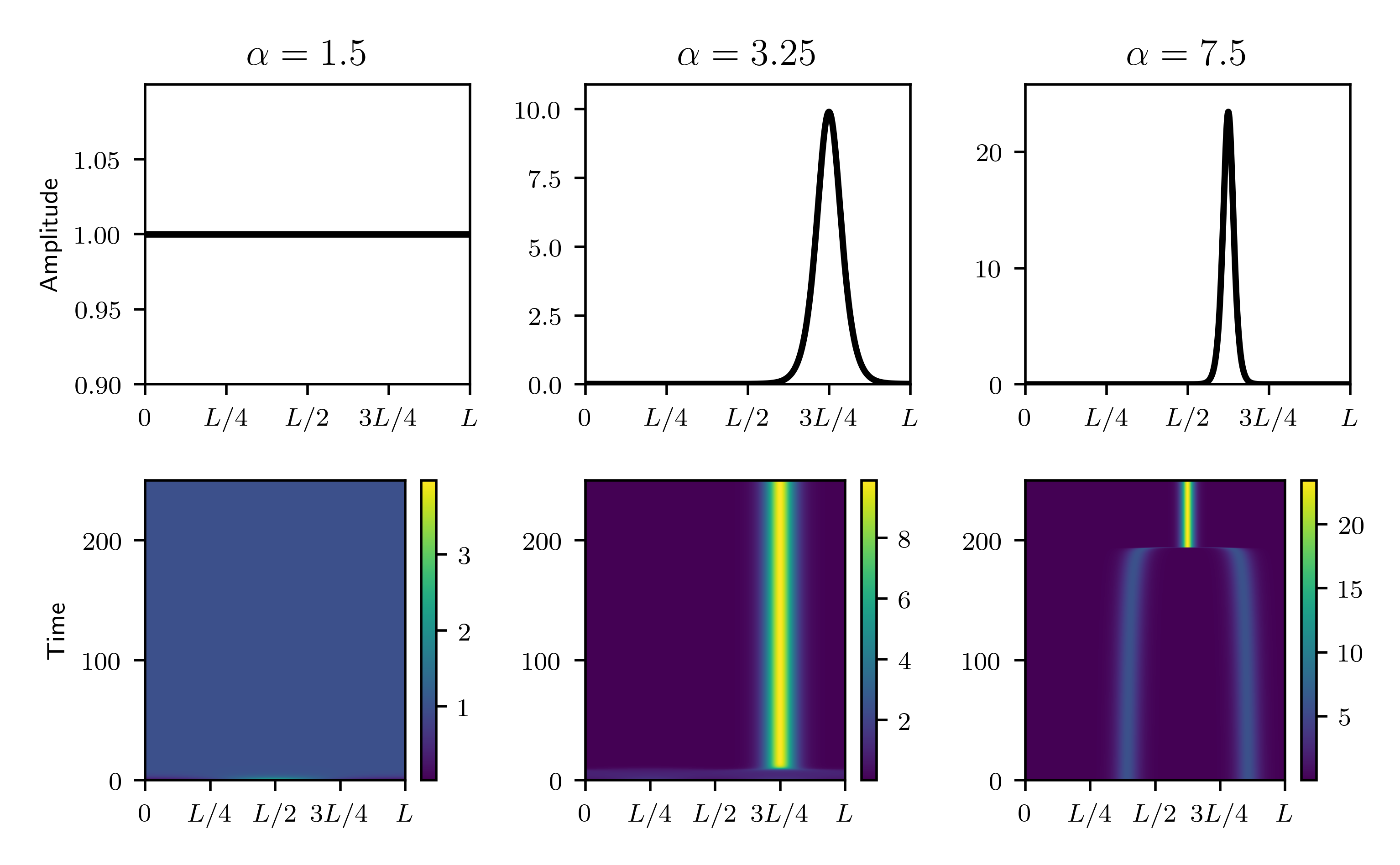}
 \caption[Numerical solution with neutral boundary conditions]{Numerical
 solutions of equation~\eqref{Eqn:ArmstrongModelNoFlux}, with the no-flux
 sensing domain $E_f(x)$, and a neutral non-local operator as constructed in
 equation~\eqref{Eqn:NonLocalBCNeutral}. In the top row we show the final
 solution profiles, below are the kymographs.
 (Left) $\alpha = 1.5$, (Middle) $\alpha = 3.25$, and (Right) $\alpha = 7.5$.
 Compared to the no-flux solutions in \cref{fig:noflux}, the solutions here do
 not feel the boundary, and peaks may form anywhere in the domain.
 }\label{fig:neutralbc}
\end{figure}

Note that this definition does not change the definition of the non-local
term in the interior of the domain (\ie, when $x$ is at least a distance of $R$
from the boundary) since there $\K[\bar u] = 0$. It only changes the
definition of the operator near the boundary. It is interesting to note that the
neutral operator $\K_n$ compared to the original operator $\K$ can be written as
\[
    \K_n[u] = \K[u]-\K[\bar u].
\]
This operation can be easily done for the other cases that we studied and
we get a neutral naive operator\index{neutral naive operator}, and a neutral
weakly adhesive\index{neutral weakly adhesive operator} operator. In the
periodic case we have $\K[\bar u] =0$, hence the non-local operator in the
periodic case is already neutral.

Some numerical solutions for three choices of $\alpha$ are shown in
\cref{fig:neutralbc}. Compared to the numerical solutions with no-flux boundary
condition, the neutral boundary condition solutions do not feel the boundary
(see \cref{fig:neutralbc}).

\section{Weakly adhesive and repulsive boundary conditions\index{weakly adhesive
boundary conditions}\index{weakly repulsive boundary conditions}}

The framework developed here can be used to explicitly model adhesion or
repulsion by the domain boundary. For that we assume that the interaction
force with the boundary is proportional to the extent of cell protrusions
that attach to the boundary, which corresponds to the amount of cell
protrusions that would reach out of the domain if there was no boundary
(see \cref{fig:bioboundary1}~(B)). For
example at $x\in (0,R)$. If the cell extends to $x-R$, then the interval
$[x-R,0)$ is outside of the domain. We assume that instead of leaving the
domain, the protrusion interacts with the boundary, given boundary adhesion
terms of the form
\begin{eqnarray*}
    a^0(x) \coloneqq \beta^0\int_{-R}^{-x} \Omega(r) \dd r, && x\in[0,R)\\
    a^L(x) \coloneqq \beta^L\int_{L-x}^R \Omega(r) \dd r, && x\in(L-R,L]
\end{eqnarray*}
where $\beta^0$ and $\beta^L$ are constants of proportionality.
$\beta^0,\beta^L>0$ describes boundary adhesion, while $\beta^0, \beta^L<0$
describes boundary repulsion.

In this case we define the adhesion operator as linear combination of all
relevant adhesive effects. Using indicator functions $\chi_A(r)$ we can
write
\begin{equation}\label{Fadhesive}
\begin{split}
    \K[u](x,t) &= \int_{E_0(x)} H(u(x+r, t)) \Omega(r) \dd r\\
       & \hspace*{5mm} + \beta^0 \chi_{[0,R]}(x)   \int_{-R}^{-x} \Omega(r) \dd r
        + \beta^L \chi_{[L-R,L]}(x) \int_{L-x}^{R} \Omega(r) \dd r \nonumber\\
    &= \int_{-R}^R\left(\chi_{E_0(x)}H(u(x+r,t))+\beta^0\chi_{[-R,-x]}(x) + \beta^L\chi_{[L-x,R]} \right)\Omega(r) \dd r,
\end{split}
\end{equation}
where we omitted the $r$-dependence in the indicator functions for brevity.
Here $E(x)$ is any suitable sampling domain as defined in
Definition~\ref{def:slice}. Further we note that whenever
\[
    \beta^0 = \frac{1}{2} \int_{E(0)} H(u(r, t))\Omega(r) \dd r,
\]
a similar expression can be found for $\beta^L$, then $\K[u]$ satisfies
condition~\eqref{Eqn:IndepBc}.

\begin{figure}\centering
 \includegraphics[width=0.9\textwidth]{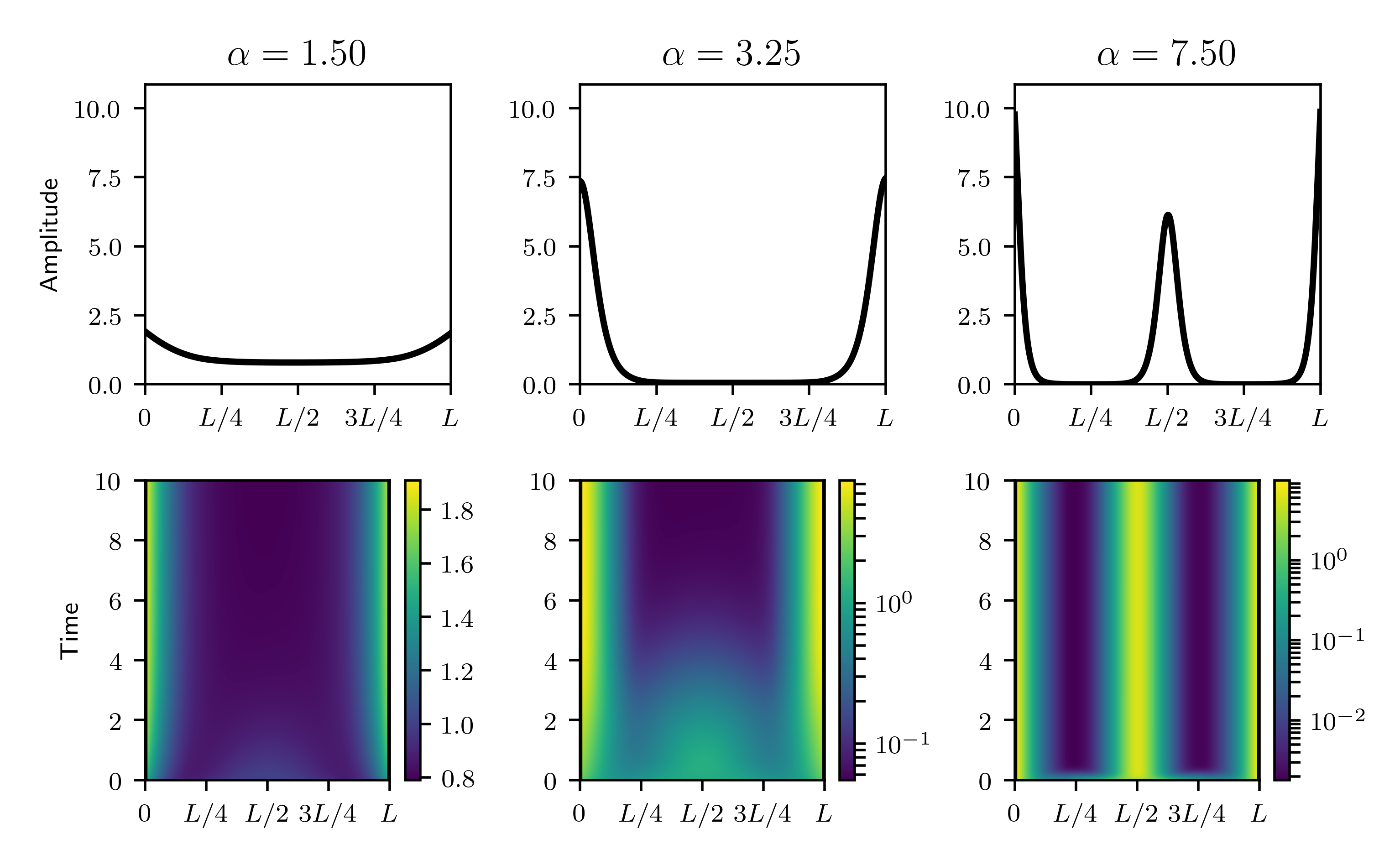}
 \caption{Numerical solutions of equation~\eqref{Eqn:ArmstrongModelNoFlux}, with
 the naive sensing domain\index{naive sensing domain} $E_0(x)$. The non-local
 operator is modified as given in equation~\eqref{Fadhesive}, here $\beta^0,
 \beta^L = 2$. In the top row we
 show the final solution profiles, below are the kymographs. (Left) $\alpha =
 1.5$, (Middle) $\alpha = 3.25$, and (Right) $\alpha = 7.5$. In these
 simulations the weakly adhesive nature of the boundaries attracts cells, and
 leads to aggregations on the boundary.
 }\label{fig:weakadh}
\end{figure}

\begin{figure}\centering
 \includegraphics[width=0.9\textwidth]{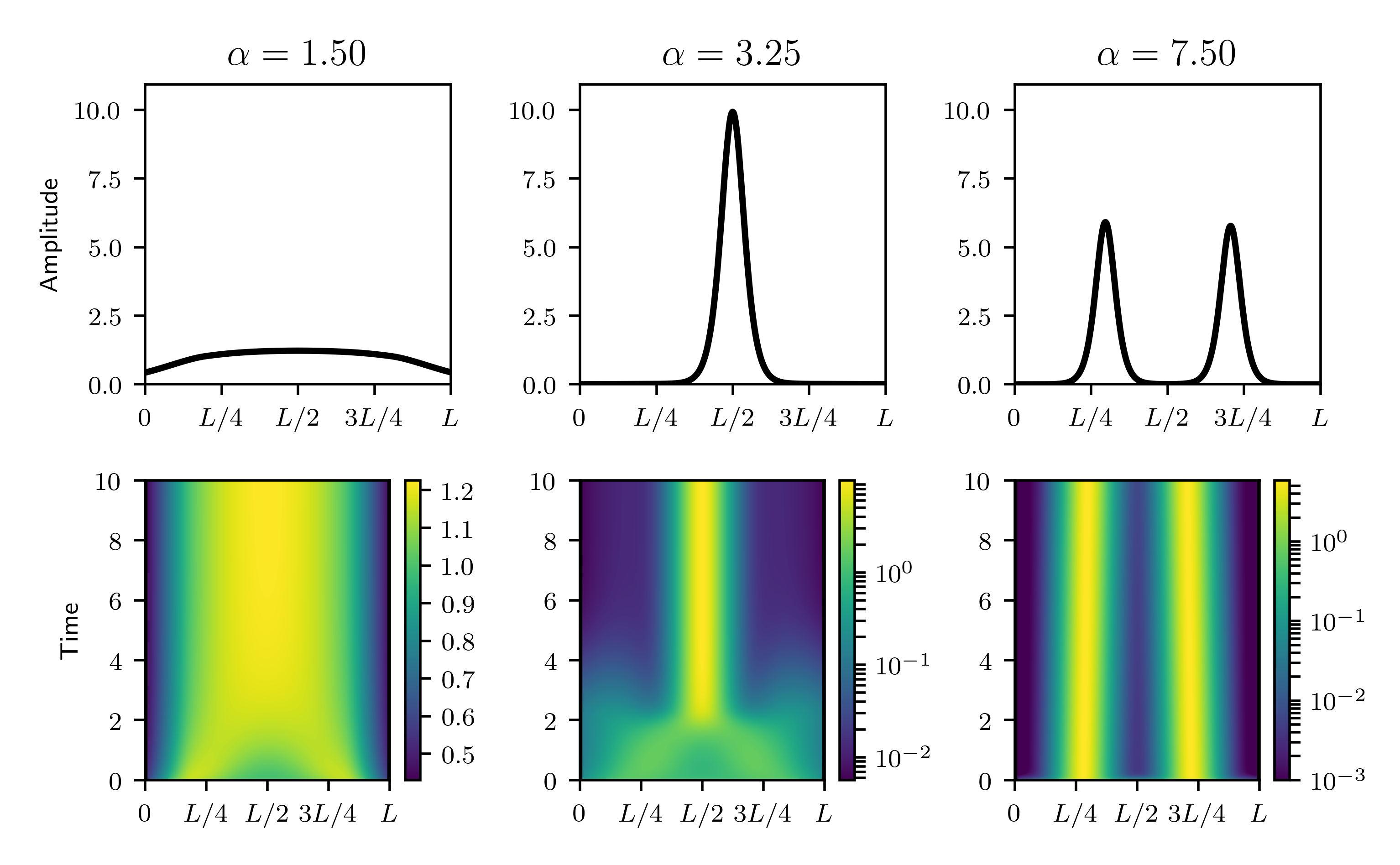}
 \caption{Numerical solutions of equation~\eqref{Eqn:ArmstrongModelNoFlux}, with
 the naive sensing domain $E_0(x)$. The non-local operator is modified as given
 in equation~\eqref{Fadhesive}, here $\beta^0, \beta^L = -1$. In the top row we
 show the final solution profiles, below are the kymographs. (Left) $\alpha = 1.5$, 
 (Middle) $\alpha = 3.25$, and (Right) $\alpha = 7.5$. In these simulations the
 weakly repulsive nature of the boundaries pushes cells away, leading to cell
 free boundaries.
 }\label{fig:weakrep}
\end{figure}

For this choice of a non-local operator we explore the solutions of
equation~\eqref{Eqn:ArmstrongModelNoFlux} numerically. Numerical solutions for
different values of $\alpha$ are shown in \cref{fig:weakadh} for $\beta^{0,L} > 0$,
and \cref{fig:weakrep} for $\beta^{0,L} < 0$.

\section{Conclusion}

In this chapter we considered the non-local cell adhesion
model~\eqref{Eqn:ArmstrongModelNoFlux} on a bounded domain with no-flux boundary
conditions. The presence of the boundary meant that the definition of the
non-local operator in the boundary region (within one sensing radius of the
boundary) had to be revisited. There are many different ways to construct non-local
operators that satisfy the no-flux boundary
condition~\eqref{Eqn:NoFluxBoundaryConditions} and we developed a general
framework through the sampling domain $E(x)$.

We recall while discussing the periodic case in \cref{c:randomwalk}, we argued
that the non-local term can be seen as cell 
polarization\index{cell polarization}. This understanding
also helps to visualize the effect of the boundary conditions discussed here.
The cell polarization points in the direction in which the most new adhesion
bonds were formed.  In the no-flux case, more adhesion bonds are formed in the
domain inside than along the boundary, hence cells at the boundary are facing
inward. This is consistent
with biological observations of cells approaching solid boundaries
\parencite{Paksa2016}.

\begin{table}\centering\small
    \def\arraystretch{1.5}
    \setlength\extrarowheight{2.5pt}
\begin{tabular}{@{\extracolsep{\fill}} |l|c|c|c|}
\hline
\bf Case & $\K[u]$ & $f_1(x)$ & $ f_2(x)$ \\
\hline
\hline
 periodic & $\int_{V} h(u) \Omega \dd r $ & $-R $ & $R$\\
\hline
naive & $ \int_{E_0(x)} h(u) \Omega \dd r $ &
$\left\{\begin{array}{ll} -x & I_1\\ -R &  I_2\end{array}\right. $ & $ \left\{\begin{array}{ll} R, &   I_3\\ L-x, &I_4\end{array}\right.$\\
\hline
non-flux & $ \int_{E_f(x)} H(u) \Omega \dd r $ &
$\left\{\begin{array}{ll} R-2x, &  I_1\\ -R, &  I_2\end{array}\right. $ & $ \left\{\begin{array}{ll} R, &  I_3\\ 2L-R-2x, & I_4\end{array}\right.$\\
\hline
weakly & $\int_{E(x)} H(u) \Omega \dd r + a_0 + a_L$ & $f_1=$naive & $f_2 =$ naive\\
adhesive & $a_0 = -2 h^0(0) \int_{-R}^{-x} \Omega(r) \dd r$ & & \\
& $ a_L = 2 h^L(L)\int_{L-x}^R \Omega(r) \dd r$ & & \\
\hline
neutral & $ \K_n[u]= \K[u]-\K[\bar u]$ & \multicolumn{2}{l|}{any of the  above cases}  \\
\hline
\end{tabular}
\caption{The different cases of suitable boundary conditions on $[0,L]$. The
 sensing domain is defined as $E(x) = [f_1(x), f_2(x)]$. The abbreviations $I_1,
 I_2, I_3, I_4$ stand for $x\in[ 0,R], x\in(R,L], x\in[0,L-R], x\in (L-R,L]$,
 respectively. \label{tab:examples}}
\end{table}

In \cref{tab:examples} we summarize the specific boundary conditions that we
studied here. The periodic case is fully translational symmetric. The naive
case\index{naive case}
and the no-flux case are weakly repellent\index{weakly repellent}. The (neutral)
non-local operator\index{(neutral) non-local operator}
switches between a repulsive and attractive behaviour, depending on whether the
population size is above or below the reference mass $\mean{u}$.  The neutral
non-local operator has the important property that $\K[\mean{u}] = 0$ for all $x
\in D$. For this reason, the bifurcation approach used in
\cref{Chapter:GlobalBifPeriodic} might be applied to the neutral case.  It is
interesting to note, that only recently Watanabe \etal\ also used the comparison
with the mass per unit length to obtain solutions of Burger's equation on a
bounded domain with no-flux boundary conditions \parencite{Watanabe2016}.

In the case of the repellent non-local operator\index{repellent non-local
operator}, we asymptotically approximated
the steady states for small values of adhesion strength $\alpha$. A more
detailed exploration of the steady-states of these equations is hindered by the
fact that there is no constant steady state that exists for all values of
$\alpha$. This means that the approach pioneered by Rabinowitz
\parencite{Rabinowitz1971} (used in \cref{Chapter:GlobalBifPeriodic}) cannot be directly
applied. In addition, even the introduction of the neutral non-local operator
does not improve the situation. This is because due to the spatial dependence of
the non-local operator we are unable to embed the no-flux solutions into a
periodic problem (see for instance \cite{Fujii1982} where this is done for
reaction-diffusion systems).

In the construction of the non-local operators that include the no-flux boundary
conditions, we had freedom to choose their behaviour in the boundary region
(within one sensing radius of the boundary). From biological experiments, it is
known that the cell polarization adapts when cells encounter physical boundaries
\parencite{Paksa2016}. Responsible for this adaptation are the intra-cellular
signalling networks\index{intra-cellular signalling networks}. In
\parencite{Buttenschon2017}, we argued that the non-local operator is a model of
cell polarization. Thus a natural extension would be to explicitly include the
intra-cellular chemical network in a detailed multi-scale model of cell
adhesion.

Mathematically novel is the non-local operator in which the integration limits
are spatially dependent, whose mathematical properties we investigated further.
The spatial dependence of the integration limits posed a particular challenge,
since properties such as continuity require a notion of convergence of sets
(integration domain).  For this reason we make use of the Fr\'echet-Nikodym
metric (see \cref{subsec:SetConvergence}). Using this metric we extended the
estimates of \cref{section:NonLocalOperator}, to the non-local operator with
spatially dependent limits of integration. Differentiation of this non-local
operator was equally made more challenging by the spatially dependent
integration limits.  Using the theory of distributions we computed the non-local
term's weak derivative, which coincides with the classical derivative if the
integration kernel $\Omega(\cdot)$ is zero on the boundary of the sensing domain
($\partial V$).

A heat equation with non-local Robin type boundary conditions 
\index{non-local Robin type boundary conditions} was studied by
Arendt \etal\ in \cite{Arendt2018}. In their case the boundary condition is
linear and can be seen as a perturbation of the heat equation semigroup. In our
case, even for linear $h(u)$, the non-local term $\alpha(u\mathcal{K} [u])_x$ is
nonlinear and Arendt's perturbation approach would not work.

Another challenge is the higher dimensional case\index{higher dimensional case}.
Boundary conditions need to distinguish between transversal\index{transversal}
and tangential effects\index{tangential effects} and combinations of
slip\index{slip}, adhesion, repulsion, and friction\index{friction} are
biologically possible. We hope that the current formalism can form a good basis
to investigate higher dimensional adhesion models on bounded domains.

%% file: AsymptoticExpansion.pgf
\begingroup%
\makeatletter%
\begin{pgfpicture}%
\pgfpathrectangle{\pgfpointorigin}{\pgfqpoint{5.500000in}{2.750000in}}%
\pgfusepath{use as bounding box, clip}%
\begin{pgfscope}%
\pgfsetbuttcap%
\pgfsetmiterjoin%
\definecolor{currentfill}{rgb}{1.000000,1.000000,1.000000}%
\pgfsetfillcolor{currentfill}%
\pgfsetlinewidth{0.000000pt}%
\definecolor{currentstroke}{rgb}{1.000000,1.000000,1.000000}%
\pgfsetstrokecolor{currentstroke}%
\pgfsetdash{}{0pt}%
\pgfpathmoveto{\pgfqpoint{0.000000in}{0.000000in}}%
\pgfpathlineto{\pgfqpoint{5.500000in}{0.000000in}}%
\pgfpathlineto{\pgfqpoint{5.500000in}{2.750000in}}%
\pgfpathlineto{\pgfqpoint{0.000000in}{2.750000in}}%
\pgfpathclose%
\pgfusepath{fill}%
\end{pgfscope}%
\begin{pgfscope}%
\pgfsetbuttcap%
\pgfsetmiterjoin%
\definecolor{currentfill}{rgb}{1.000000,1.000000,1.000000}%
\pgfsetfillcolor{currentfill}%
\pgfsetlinewidth{0.000000pt}%
\definecolor{currentstroke}{rgb}{0.000000,0.000000,0.000000}%
\pgfsetstrokecolor{currentstroke}%
\pgfsetstrokeopacity{0.000000}%
\pgfsetdash{}{0pt}%
\pgfpathmoveto{\pgfqpoint{0.653314in}{0.551844in}}%
\pgfpathlineto{\pgfqpoint{5.268228in}{0.551844in}}%
\pgfpathlineto{\pgfqpoint{5.268228in}{2.401359in}}%
\pgfpathlineto{\pgfqpoint{0.653314in}{2.401359in}}%
\pgfpathclose%
\pgfusepath{fill}%
\end{pgfscope}%
\begin{pgfscope}%
\pgfsetbuttcap%
\pgfsetroundjoin%
\definecolor{currentfill}{rgb}{0.000000,0.000000,0.000000}%
\pgfsetfillcolor{currentfill}%
\pgfsetlinewidth{0.803000pt}%
\definecolor{currentstroke}{rgb}{0.000000,0.000000,0.000000}%
\pgfsetstrokecolor{currentstroke}%
\pgfsetdash{}{0pt}%
\pgfsys@defobject{currentmarker}{\pgfqpoint{0.000000in}{-0.048611in}}{\pgfqpoint{0.000000in}{0.000000in}}{%
\pgfpathmoveto{\pgfqpoint{0.000000in}{0.000000in}}%
\pgfpathlineto{\pgfqpoint{0.000000in}{-0.048611in}}%
\pgfusepath{stroke,fill}%
}%
\begin{pgfscope}%
\pgfsys@transformshift{0.653314in}{0.551844in}%
\pgfsys@useobject{currentmarker}{}%
\end{pgfscope}%
\end{pgfscope}%
\begin{pgfscope}%
\pgftext[x=0.653314in,y=0.454622in,,top]{\rmfamily\fontsize{9.000000}{10.800000}\selectfont \(\displaystyle 0.0\)}%
\end{pgfscope}%
\begin{pgfscope}%
\pgfsetbuttcap%
\pgfsetroundjoin%
\definecolor{currentfill}{rgb}{0.000000,0.000000,0.000000}%
\pgfsetfillcolor{currentfill}%
\pgfsetlinewidth{0.803000pt}%
\definecolor{currentstroke}{rgb}{0.000000,0.000000,0.000000}%
\pgfsetstrokecolor{currentstroke}%
\pgfsetdash{}{0pt}%
\pgfsys@defobject{currentmarker}{\pgfqpoint{0.000000in}{-0.048611in}}{\pgfqpoint{0.000000in}{0.000000in}}{%
\pgfpathmoveto{\pgfqpoint{0.000000in}{0.000000in}}%
\pgfpathlineto{\pgfqpoint{0.000000in}{-0.048611in}}%
\pgfusepath{stroke,fill}%
}%
\begin{pgfscope}%
\pgfsys@transformshift{1.576296in}{0.551844in}%
\pgfsys@useobject{currentmarker}{}%
\end{pgfscope}%
\end{pgfscope}%
\begin{pgfscope}%
\pgftext[x=1.576296in,y=0.454622in,,top]{\rmfamily\fontsize{9.000000}{10.800000}\selectfont \(\displaystyle 0.5\)}%
\end{pgfscope}%
\begin{pgfscope}%
\pgfsetbuttcap%
\pgfsetroundjoin%
\definecolor{currentfill}{rgb}{0.000000,0.000000,0.000000}%
\pgfsetfillcolor{currentfill}%
\pgfsetlinewidth{0.803000pt}%
\definecolor{currentstroke}{rgb}{0.000000,0.000000,0.000000}%
\pgfsetstrokecolor{currentstroke}%
\pgfsetdash{}{0pt}%
\pgfsys@defobject{currentmarker}{\pgfqpoint{0.000000in}{-0.048611in}}{\pgfqpoint{0.000000in}{0.000000in}}{%
\pgfpathmoveto{\pgfqpoint{0.000000in}{0.000000in}}%
\pgfpathlineto{\pgfqpoint{0.000000in}{-0.048611in}}%
\pgfusepath{stroke,fill}%
}%
\begin{pgfscope}%
\pgfsys@transformshift{2.499279in}{0.551844in}%
\pgfsys@useobject{currentmarker}{}%
\end{pgfscope}%
\end{pgfscope}%
\begin{pgfscope}%
\pgftext[x=2.499279in,y=0.454622in,,top]{\rmfamily\fontsize{9.000000}{10.800000}\selectfont \(\displaystyle 1.0\)}%
\end{pgfscope}%
\begin{pgfscope}%
\pgfsetbuttcap%
\pgfsetroundjoin%
\definecolor{currentfill}{rgb}{0.000000,0.000000,0.000000}%
\pgfsetfillcolor{currentfill}%
\pgfsetlinewidth{0.803000pt}%
\definecolor{currentstroke}{rgb}{0.000000,0.000000,0.000000}%
\pgfsetstrokecolor{currentstroke}%
\pgfsetdash{}{0pt}%
\pgfsys@defobject{currentmarker}{\pgfqpoint{0.000000in}{-0.048611in}}{\pgfqpoint{0.000000in}{0.000000in}}{%
\pgfpathmoveto{\pgfqpoint{0.000000in}{0.000000in}}%
\pgfpathlineto{\pgfqpoint{0.000000in}{-0.048611in}}%
\pgfusepath{stroke,fill}%
}%
\begin{pgfscope}%
\pgfsys@transformshift{3.422262in}{0.551844in}%
\pgfsys@useobject{currentmarker}{}%
\end{pgfscope}%
\end{pgfscope}%
\begin{pgfscope}%
\pgftext[x=3.422262in,y=0.454622in,,top]{\rmfamily\fontsize{9.000000}{10.800000}\selectfont \(\displaystyle 1.5\)}%
\end{pgfscope}%
\begin{pgfscope}%
\pgfsetbuttcap%
\pgfsetroundjoin%
\definecolor{currentfill}{rgb}{0.000000,0.000000,0.000000}%
\pgfsetfillcolor{currentfill}%
\pgfsetlinewidth{0.803000pt}%
\definecolor{currentstroke}{rgb}{0.000000,0.000000,0.000000}%
\pgfsetstrokecolor{currentstroke}%
\pgfsetdash{}{0pt}%
\pgfsys@defobject{currentmarker}{\pgfqpoint{0.000000in}{-0.048611in}}{\pgfqpoint{0.000000in}{0.000000in}}{%
\pgfpathmoveto{\pgfqpoint{0.000000in}{0.000000in}}%
\pgfpathlineto{\pgfqpoint{0.000000in}{-0.048611in}}%
\pgfusepath{stroke,fill}%
}%
\begin{pgfscope}%
\pgfsys@transformshift{4.345245in}{0.551844in}%
\pgfsys@useobject{currentmarker}{}%
\end{pgfscope}%
\end{pgfscope}%
\begin{pgfscope}%
\pgftext[x=4.345245in,y=0.454622in,,top]{\rmfamily\fontsize{9.000000}{10.800000}\selectfont \(\displaystyle 2.0\)}%
\end{pgfscope}%
\begin{pgfscope}%
\pgfsetbuttcap%
\pgfsetroundjoin%
\definecolor{currentfill}{rgb}{0.000000,0.000000,0.000000}%
\pgfsetfillcolor{currentfill}%
\pgfsetlinewidth{0.803000pt}%
\definecolor{currentstroke}{rgb}{0.000000,0.000000,0.000000}%
\pgfsetstrokecolor{currentstroke}%
\pgfsetdash{}{0pt}%
\pgfsys@defobject{currentmarker}{\pgfqpoint{0.000000in}{-0.048611in}}{\pgfqpoint{0.000000in}{0.000000in}}{%
\pgfpathmoveto{\pgfqpoint{0.000000in}{0.000000in}}%
\pgfpathlineto{\pgfqpoint{0.000000in}{-0.048611in}}%
\pgfusepath{stroke,fill}%
}%
\begin{pgfscope}%
\pgfsys@transformshift{5.268228in}{0.551844in}%
\pgfsys@useobject{currentmarker}{}%
\end{pgfscope}%
\end{pgfscope}%
\begin{pgfscope}%
\pgftext[x=5.268228in,y=0.454622in,,top]{\rmfamily\fontsize{9.000000}{10.800000}\selectfont \(\displaystyle 2.5\)}%
\end{pgfscope}%
\begin{pgfscope}%
\pgftext[x=2.960771in,y=0.288676in,,top]{\rmfamily\fontsize{10.000000}{12.000000}\selectfont Space [\(\displaystyle x\)]}%
\end{pgfscope}%
\begin{pgfscope}%
\pgfsetbuttcap%
\pgfsetroundjoin%
\definecolor{currentfill}{rgb}{0.000000,0.000000,0.000000}%
\pgfsetfillcolor{currentfill}%
\pgfsetlinewidth{0.803000pt}%
\definecolor{currentstroke}{rgb}{0.000000,0.000000,0.000000}%
\pgfsetstrokecolor{currentstroke}%
\pgfsetdash{}{0pt}%
\pgfsys@defobject{currentmarker}{\pgfqpoint{-0.048611in}{0.000000in}}{\pgfqpoint{0.000000in}{0.000000in}}{%
\pgfpathmoveto{\pgfqpoint{0.000000in}{0.000000in}}%
\pgfpathlineto{\pgfqpoint{-0.048611in}{0.000000in}}%
\pgfusepath{stroke,fill}%
}%
\begin{pgfscope}%
\pgfsys@transformshift{0.653314in}{0.635849in}%
\pgfsys@useobject{currentmarker}{}%
\end{pgfscope}%
\end{pgfscope}%
\begin{pgfscope}%
\pgftext[x=0.327698in,y=0.592804in,left,base]{\rmfamily\fontsize{9.000000}{10.800000}\selectfont \(\displaystyle 0.85\)}%
\end{pgfscope}%
\begin{pgfscope}%
\pgfsetbuttcap%
\pgfsetroundjoin%
\definecolor{currentfill}{rgb}{0.000000,0.000000,0.000000}%
\pgfsetfillcolor{currentfill}%
\pgfsetlinewidth{0.803000pt}%
\definecolor{currentstroke}{rgb}{0.000000,0.000000,0.000000}%
\pgfsetstrokecolor{currentstroke}%
\pgfsetdash{}{0pt}%
\pgfsys@defobject{currentmarker}{\pgfqpoint{-0.048611in}{0.000000in}}{\pgfqpoint{0.000000in}{0.000000in}}{%
\pgfpathmoveto{\pgfqpoint{0.000000in}{0.000000in}}%
\pgfpathlineto{\pgfqpoint{-0.048611in}{0.000000in}}%
\pgfusepath{stroke,fill}%
}%
\begin{pgfscope}%
\pgfsys@transformshift{0.653314in}{0.972145in}%
\pgfsys@useobject{currentmarker}{}%
\end{pgfscope}%
\end{pgfscope}%
\begin{pgfscope}%
\pgftext[x=0.327698in,y=0.929100in,left,base]{\rmfamily\fontsize{9.000000}{10.800000}\selectfont \(\displaystyle 0.90\)}%
\end{pgfscope}%
\begin{pgfscope}%
\pgfsetbuttcap%
\pgfsetroundjoin%
\definecolor{currentfill}{rgb}{0.000000,0.000000,0.000000}%
\pgfsetfillcolor{currentfill}%
\pgfsetlinewidth{0.803000pt}%
\definecolor{currentstroke}{rgb}{0.000000,0.000000,0.000000}%
\pgfsetstrokecolor{currentstroke}%
\pgfsetdash{}{0pt}%
\pgfsys@defobject{currentmarker}{\pgfqpoint{-0.048611in}{0.000000in}}{\pgfqpoint{0.000000in}{0.000000in}}{%
\pgfpathmoveto{\pgfqpoint{0.000000in}{0.000000in}}%
\pgfpathlineto{\pgfqpoint{-0.048611in}{0.000000in}}%
\pgfusepath{stroke,fill}%
}%
\begin{pgfscope}%
\pgfsys@transformshift{0.653314in}{1.308442in}%
\pgfsys@useobject{currentmarker}{}%
\end{pgfscope}%
\end{pgfscope}%
\begin{pgfscope}%
\pgftext[x=0.327698in,y=1.265397in,left,base]{\rmfamily\fontsize{9.000000}{10.800000}\selectfont \(\displaystyle 0.95\)}%
\end{pgfscope}%
\begin{pgfscope}%
\pgfsetbuttcap%
\pgfsetroundjoin%
\definecolor{currentfill}{rgb}{0.000000,0.000000,0.000000}%
\pgfsetfillcolor{currentfill}%
\pgfsetlinewidth{0.803000pt}%
\definecolor{currentstroke}{rgb}{0.000000,0.000000,0.000000}%
\pgfsetstrokecolor{currentstroke}%
\pgfsetdash{}{0pt}%
\pgfsys@defobject{currentmarker}{\pgfqpoint{-0.048611in}{0.000000in}}{\pgfqpoint{0.000000in}{0.000000in}}{%
\pgfpathmoveto{\pgfqpoint{0.000000in}{0.000000in}}%
\pgfpathlineto{\pgfqpoint{-0.048611in}{0.000000in}}%
\pgfusepath{stroke,fill}%
}%
\begin{pgfscope}%
\pgfsys@transformshift{0.653314in}{1.644739in}%
\pgfsys@useobject{currentmarker}{}%
\end{pgfscope}%
\end{pgfscope}%
\begin{pgfscope}%
\pgftext[x=0.327698in,y=1.601694in,left,base]{\rmfamily\fontsize{9.000000}{10.800000}\selectfont \(\displaystyle 1.00\)}%
\end{pgfscope}%
\begin{pgfscope}%
\pgfsetbuttcap%
\pgfsetroundjoin%
\definecolor{currentfill}{rgb}{0.000000,0.000000,0.000000}%
\pgfsetfillcolor{currentfill}%
\pgfsetlinewidth{0.803000pt}%
\definecolor{currentstroke}{rgb}{0.000000,0.000000,0.000000}%
\pgfsetstrokecolor{currentstroke}%
\pgfsetdash{}{0pt}%
\pgfsys@defobject{currentmarker}{\pgfqpoint{-0.048611in}{0.000000in}}{\pgfqpoint{0.000000in}{0.000000in}}{%
\pgfpathmoveto{\pgfqpoint{0.000000in}{0.000000in}}%
\pgfpathlineto{\pgfqpoint{-0.048611in}{0.000000in}}%
\pgfusepath{stroke,fill}%
}%
\begin{pgfscope}%
\pgfsys@transformshift{0.653314in}{1.981036in}%
\pgfsys@useobject{currentmarker}{}%
\end{pgfscope}%
\end{pgfscope}%
\begin{pgfscope}%
\pgftext[x=0.327698in,y=1.937991in,left,base]{\rmfamily\fontsize{9.000000}{10.800000}\selectfont \(\displaystyle 1.05\)}%
\end{pgfscope}%
\begin{pgfscope}%
\pgfsetbuttcap%
\pgfsetroundjoin%
\definecolor{currentfill}{rgb}{0.000000,0.000000,0.000000}%
\pgfsetfillcolor{currentfill}%
\pgfsetlinewidth{0.803000pt}%
\definecolor{currentstroke}{rgb}{0.000000,0.000000,0.000000}%
\pgfsetstrokecolor{currentstroke}%
\pgfsetdash{}{0pt}%
\pgfsys@defobject{currentmarker}{\pgfqpoint{-0.048611in}{0.000000in}}{\pgfqpoint{0.000000in}{0.000000in}}{%
\pgfpathmoveto{\pgfqpoint{0.000000in}{0.000000in}}%
\pgfpathlineto{\pgfqpoint{-0.048611in}{0.000000in}}%
\pgfusepath{stroke,fill}%
}%
\begin{pgfscope}%
\pgfsys@transformshift{0.653314in}{2.317332in}%
\pgfsys@useobject{currentmarker}{}%
\end{pgfscope}%
\end{pgfscope}%
\begin{pgfscope}%
\pgftext[x=0.327698in,y=2.274287in,left,base]{\rmfamily\fontsize{9.000000}{10.800000}\selectfont \(\displaystyle 1.10\)}%
\end{pgfscope}%
\begin{pgfscope}%
\pgftext[x=0.272142in,y=1.476601in,,bottom,rotate=90.000000]{\rmfamily\fontsize{10.000000}{12.000000}\selectfont \(\displaystyle u\)}%
\end{pgfscope}%
\begin{pgfscope}%
\pgfpathrectangle{\pgfqpoint{0.653314in}{0.551844in}}{\pgfqpoint{4.614914in}{1.849515in}} %
\pgfusepath{clip}%
\pgfsetrectcap%
\pgfsetroundjoin%
\pgfsetlinewidth{1.505625pt}%
\definecolor{currentstroke}{rgb}{0.000000,0.000000,0.000000}%
\pgfsetstrokecolor{currentstroke}%
\pgfsetdash{}{0pt}%
\pgfpathmoveto{\pgfqpoint{0.653314in}{1.644739in}}%
\pgfpathlineto{\pgfqpoint{5.268228in}{1.644739in}}%
\pgfpathlineto{\pgfqpoint{5.268228in}{1.644739in}}%
\pgfusepath{stroke}%
\end{pgfscope}%
\begin{pgfscope}%
\pgfpathrectangle{\pgfqpoint{0.653314in}{0.551844in}}{\pgfqpoint{4.614914in}{1.849515in}} %
\pgfusepath{clip}%
\pgfsetbuttcap%
\pgfsetroundjoin%
\pgfsetlinewidth{1.003750pt}%
\definecolor{currentstroke}{rgb}{0.121569,0.466667,0.705882}%
\pgfsetstrokecolor{currentstroke}%
\pgfsetdash{{3.700000pt}{1.600000pt}}{0.000000pt}%
\pgfpathmoveto{\pgfqpoint{0.653314in}{1.543856in}}%
\pgfpathlineto{\pgfqpoint{0.990419in}{1.545360in}}%
\pgfpathlineto{\pgfqpoint{1.073343in}{1.548998in}}%
\pgfpathlineto{\pgfqpoint{1.143648in}{1.555127in}}%
\pgfpathlineto{\pgfqpoint{1.217559in}{1.564720in}}%
\pgfpathlineto{\pgfqpoint{1.316707in}{1.580975in}}%
\pgfpathlineto{\pgfqpoint{1.581704in}{1.628910in}}%
\pgfpathlineto{\pgfqpoint{1.915204in}{1.688113in}}%
\pgfpathlineto{\pgfqpoint{1.998128in}{1.699508in}}%
\pgfpathlineto{\pgfqpoint{2.068434in}{1.706109in}}%
\pgfpathlineto{\pgfqpoint{2.142344in}{1.709909in}}%
\pgfpathlineto{\pgfqpoint{2.241493in}{1.711669in}}%
\pgfpathlineto{\pgfqpoint{2.510095in}{1.711994in}}%
\pgfpathlineto{\pgfqpoint{3.759367in}{1.710490in}}%
\pgfpathlineto{\pgfqpoint{3.842291in}{1.706852in}}%
\pgfpathlineto{\pgfqpoint{3.912597in}{1.700724in}}%
\pgfpathlineto{\pgfqpoint{3.986507in}{1.691130in}}%
\pgfpathlineto{\pgfqpoint{4.085656in}{1.674875in}}%
\pgfpathlineto{\pgfqpoint{4.350653in}{1.626940in}}%
\pgfpathlineto{\pgfqpoint{4.684153in}{1.567738in}}%
\pgfpathlineto{\pgfqpoint{4.767077in}{1.556342in}}%
\pgfpathlineto{\pgfqpoint{4.837382in}{1.549742in}}%
\pgfpathlineto{\pgfqpoint{4.911293in}{1.545941in}}%
\pgfpathlineto{\pgfqpoint{5.010441in}{1.544182in}}%
\pgfpathlineto{\pgfqpoint{5.268228in}{1.543856in}}%
\pgfpathlineto{\pgfqpoint{5.268228in}{1.543856in}}%
\pgfusepath{stroke}%
\end{pgfscope}%
\begin{pgfscope}%
\pgfpathrectangle{\pgfqpoint{0.653314in}{0.551844in}}{\pgfqpoint{4.614914in}{1.849515in}} %
\pgfusepath{clip}%
\pgfsetbuttcap%
\pgfsetroundjoin%
\pgfsetlinewidth{1.003750pt}%
\definecolor{currentstroke}{rgb}{1.000000,0.498039,0.054902}%
\pgfsetstrokecolor{currentstroke}%
\pgfsetdash{{6.400000pt}{1.600000pt}{1.000000pt}{1.600000pt}}{0.000000pt}%
\pgfpathmoveto{\pgfqpoint{0.653314in}{1.442974in}}%
\pgfpathlineto{\pgfqpoint{0.952562in}{1.444490in}}%
\pgfpathlineto{\pgfqpoint{1.022867in}{1.448050in}}%
\pgfpathlineto{\pgfqpoint{1.076948in}{1.453737in}}%
\pgfpathlineto{\pgfqpoint{1.127424in}{1.462138in}}%
\pgfpathlineto{\pgfqpoint{1.177899in}{1.473680in}}%
\pgfpathlineto{\pgfqpoint{1.235586in}{1.490175in}}%
\pgfpathlineto{\pgfqpoint{1.311299in}{1.515323in}}%
\pgfpathlineto{\pgfqpoint{1.460924in}{1.569091in}}%
\pgfpathlineto{\pgfqpoint{1.884558in}{1.721624in}}%
\pgfpathlineto{\pgfqpoint{1.954864in}{1.743223in}}%
\pgfpathlineto{\pgfqpoint{2.008945in}{1.756708in}}%
\pgfpathlineto{\pgfqpoint{2.059420in}{1.766131in}}%
\pgfpathlineto{\pgfqpoint{2.111699in}{1.772620in}}%
\pgfpathlineto{\pgfqpoint{2.171188in}{1.776675in}}%
\pgfpathlineto{\pgfqpoint{2.250506in}{1.778715in}}%
\pgfpathlineto{\pgfqpoint{2.427171in}{1.779246in}}%
\pgfpathlineto{\pgfqpoint{3.716102in}{1.777883in}}%
\pgfpathlineto{\pgfqpoint{3.788210in}{1.774446in}}%
\pgfpathlineto{\pgfqpoint{3.844094in}{1.768727in}}%
\pgfpathlineto{\pgfqpoint{3.894570in}{1.760440in}}%
\pgfpathlineto{\pgfqpoint{3.945045in}{1.749006in}}%
\pgfpathlineto{\pgfqpoint{4.002732in}{1.732606in}}%
\pgfpathlineto{\pgfqpoint{4.078445in}{1.707527in}}%
\pgfpathlineto{\pgfqpoint{4.224464in}{1.655099in}}%
\pgfpathlineto{\pgfqpoint{4.655309in}{1.500004in}}%
\pgfpathlineto{\pgfqpoint{4.723812in}{1.478999in}}%
\pgfpathlineto{\pgfqpoint{4.777893in}{1.465514in}}%
\pgfpathlineto{\pgfqpoint{4.828369in}{1.456092in}}%
\pgfpathlineto{\pgfqpoint{4.880647in}{1.449603in}}%
\pgfpathlineto{\pgfqpoint{4.940136in}{1.445548in}}%
\pgfpathlineto{\pgfqpoint{5.019455in}{1.443508in}}%
\pgfpathlineto{\pgfqpoint{5.196120in}{1.442977in}}%
\pgfpathlineto{\pgfqpoint{5.268228in}{1.442974in}}%
\pgfpathlineto{\pgfqpoint{5.268228in}{1.442974in}}%
\pgfusepath{stroke}%
\end{pgfscope}%
\begin{pgfscope}%
\pgfpathrectangle{\pgfqpoint{0.653314in}{0.551844in}}{\pgfqpoint{4.614914in}{1.849515in}} %
\pgfusepath{clip}%
\pgfsetbuttcap%
\pgfsetroundjoin%
\pgfsetlinewidth{1.003750pt}%
\definecolor{currentstroke}{rgb}{0.172549,0.627451,0.172549}%
\pgfsetstrokecolor{currentstroke}%
\pgfsetdash{{1.000000pt}{1.650000pt}}{0.000000pt}%
\pgfpathmoveto{\pgfqpoint{0.653314in}{1.140326in}}%
\pgfpathlineto{\pgfqpoint{0.907494in}{1.141830in}}%
\pgfpathlineto{\pgfqpoint{0.966984in}{1.145295in}}%
\pgfpathlineto{\pgfqpoint{1.010248in}{1.150750in}}%
\pgfpathlineto{\pgfqpoint{1.046302in}{1.158222in}}%
\pgfpathlineto{\pgfqpoint{1.080554in}{1.168474in}}%
\pgfpathlineto{\pgfqpoint{1.113002in}{1.181423in}}%
\pgfpathlineto{\pgfqpoint{1.145451in}{1.197666in}}%
\pgfpathlineto{\pgfqpoint{1.179702in}{1.218257in}}%
\pgfpathlineto{\pgfqpoint{1.217559in}{1.244645in}}%
\pgfpathlineto{\pgfqpoint{1.262626in}{1.279937in}}%
\pgfpathlineto{\pgfqpoint{1.323918in}{1.332240in}}%
\pgfpathlineto{\pgfqpoint{1.432080in}{1.429402in}}%
\pgfpathlineto{\pgfqpoint{1.864728in}{1.820314in}}%
\pgfpathlineto{\pgfqpoint{1.920612in}{1.865798in}}%
\pgfpathlineto{\pgfqpoint{1.963877in}{1.897121in}}%
\pgfpathlineto{\pgfqpoint{1.999931in}{1.919624in}}%
\pgfpathlineto{\pgfqpoint{2.034182in}{1.937439in}}%
\pgfpathlineto{\pgfqpoint{2.066631in}{1.950934in}}%
\pgfpathlineto{\pgfqpoint{2.100882in}{1.961718in}}%
\pgfpathlineto{\pgfqpoint{2.136936in}{1.969661in}}%
\pgfpathlineto{\pgfqpoint{2.176596in}{1.975166in}}%
\pgfpathlineto{\pgfqpoint{2.225269in}{1.978720in}}%
\pgfpathlineto{\pgfqpoint{2.297377in}{1.980576in}}%
\pgfpathlineto{\pgfqpoint{2.477647in}{1.981014in}}%
\pgfpathlineto{\pgfqpoint{3.676443in}{1.979510in}}%
\pgfpathlineto{\pgfqpoint{3.735932in}{1.976045in}}%
\pgfpathlineto{\pgfqpoint{3.779197in}{1.970590in}}%
\pgfpathlineto{\pgfqpoint{3.815251in}{1.963118in}}%
\pgfpathlineto{\pgfqpoint{3.849502in}{1.952866in}}%
\pgfpathlineto{\pgfqpoint{3.881951in}{1.939917in}}%
\pgfpathlineto{\pgfqpoint{3.914399in}{1.923674in}}%
\pgfpathlineto{\pgfqpoint{3.948651in}{1.903083in}}%
\pgfpathlineto{\pgfqpoint{3.986507in}{1.876695in}}%
\pgfpathlineto{\pgfqpoint{4.031575in}{1.841403in}}%
\pgfpathlineto{\pgfqpoint{4.092867in}{1.789100in}}%
\pgfpathlineto{\pgfqpoint{4.201029in}{1.691938in}}%
\pgfpathlineto{\pgfqpoint{4.633677in}{1.301026in}}%
\pgfpathlineto{\pgfqpoint{4.689561in}{1.255542in}}%
\pgfpathlineto{\pgfqpoint{4.732825in}{1.224219in}}%
\pgfpathlineto{\pgfqpoint{4.768879in}{1.201716in}}%
\pgfpathlineto{\pgfqpoint{4.803131in}{1.183901in}}%
\pgfpathlineto{\pgfqpoint{4.835579in}{1.170406in}}%
\pgfpathlineto{\pgfqpoint{4.869831in}{1.159622in}}%
\pgfpathlineto{\pgfqpoint{4.905885in}{1.151679in}}%
\pgfpathlineto{\pgfqpoint{4.945544in}{1.146174in}}%
\pgfpathlineto{\pgfqpoint{4.994217in}{1.142620in}}%
\pgfpathlineto{\pgfqpoint{5.066325in}{1.140764in}}%
\pgfpathlineto{\pgfqpoint{5.246595in}{1.140326in}}%
\pgfpathlineto{\pgfqpoint{5.268228in}{1.140326in}}%
\pgfpathlineto{\pgfqpoint{5.268228in}{1.140326in}}%
\pgfusepath{stroke}%
\end{pgfscope}%
\begin{pgfscope}%
\pgfpathrectangle{\pgfqpoint{0.653314in}{0.551844in}}{\pgfqpoint{4.614914in}{1.849515in}} %
\pgfusepath{clip}%
\pgfsetrectcap%
\pgfsetroundjoin%
\pgfsetlinewidth{1.003750pt}%
\definecolor{currentstroke}{rgb}{0.839216,0.152941,0.156863}%
\pgfsetstrokecolor{currentstroke}%
\pgfsetdash{}{0pt}%
\pgfpathmoveto{\pgfqpoint{0.653314in}{0.635913in}}%
\pgfpathlineto{\pgfqpoint{0.876848in}{0.637406in}}%
\pgfpathlineto{\pgfqpoint{0.930929in}{0.640853in}}%
\pgfpathlineto{\pgfqpoint{0.968786in}{0.646184in}}%
\pgfpathlineto{\pgfqpoint{0.999432in}{0.653404in}}%
\pgfpathlineto{\pgfqpoint{1.026473in}{0.662724in}}%
\pgfpathlineto{\pgfqpoint{1.051710in}{0.674505in}}%
\pgfpathlineto{\pgfqpoint{1.076948in}{0.689731in}}%
\pgfpathlineto{\pgfqpoint{1.102186in}{0.708746in}}%
\pgfpathlineto{\pgfqpoint{1.127424in}{0.731733in}}%
\pgfpathlineto{\pgfqpoint{1.154464in}{0.760768in}}%
\pgfpathlineto{\pgfqpoint{1.183308in}{0.796494in}}%
\pgfpathlineto{\pgfqpoint{1.215756in}{0.841886in}}%
\pgfpathlineto{\pgfqpoint{1.253613in}{0.900473in}}%
\pgfpathlineto{\pgfqpoint{1.302286in}{0.982030in}}%
\pgfpathlineto{\pgfqpoint{1.374394in}{1.109677in}}%
\pgfpathlineto{\pgfqpoint{1.560072in}{1.447046in}}%
\pgfpathlineto{\pgfqpoint{1.832280in}{1.939791in}}%
\pgfpathlineto{\pgfqpoint{1.891769in}{2.041019in}}%
\pgfpathlineto{\pgfqpoint{1.933231in}{2.105973in}}%
\pgfpathlineto{\pgfqpoint{1.967482in}{2.154323in}}%
\pgfpathlineto{\pgfqpoint{1.998128in}{2.192435in}}%
\pgfpathlineto{\pgfqpoint{2.025169in}{2.221469in}}%
\pgfpathlineto{\pgfqpoint{2.052209in}{2.245945in}}%
\pgfpathlineto{\pgfqpoint{2.077447in}{2.264681in}}%
\pgfpathlineto{\pgfqpoint{2.102685in}{2.279648in}}%
\pgfpathlineto{\pgfqpoint{2.127923in}{2.291200in}}%
\pgfpathlineto{\pgfqpoint{2.154963in}{2.300313in}}%
\pgfpathlineto{\pgfqpoint{2.185609in}{2.307351in}}%
\pgfpathlineto{\pgfqpoint{2.221663in}{2.312349in}}%
\pgfpathlineto{\pgfqpoint{2.268533in}{2.315520in}}%
\pgfpathlineto{\pgfqpoint{2.342444in}{2.317030in}}%
\pgfpathlineto{\pgfqpoint{2.587612in}{2.317290in}}%
\pgfpathlineto{\pgfqpoint{3.645797in}{2.315797in}}%
\pgfpathlineto{\pgfqpoint{3.699878in}{2.312349in}}%
\pgfpathlineto{\pgfqpoint{3.737735in}{2.307018in}}%
\pgfpathlineto{\pgfqpoint{3.768381in}{2.299798in}}%
\pgfpathlineto{\pgfqpoint{3.795421in}{2.290478in}}%
\pgfpathlineto{\pgfqpoint{3.820659in}{2.278698in}}%
\pgfpathlineto{\pgfqpoint{3.845897in}{2.263471in}}%
\pgfpathlineto{\pgfqpoint{3.871134in}{2.244457in}}%
\pgfpathlineto{\pgfqpoint{3.896372in}{2.221469in}}%
\pgfpathlineto{\pgfqpoint{3.923413in}{2.192435in}}%
\pgfpathlineto{\pgfqpoint{3.952256in}{2.156708in}}%
\pgfpathlineto{\pgfqpoint{3.984705in}{2.111316in}}%
\pgfpathlineto{\pgfqpoint{4.022561in}{2.052729in}}%
\pgfpathlineto{\pgfqpoint{4.071234in}{1.971172in}}%
\pgfpathlineto{\pgfqpoint{4.143342in}{1.843525in}}%
\pgfpathlineto{\pgfqpoint{4.329020in}{1.506156in}}%
\pgfpathlineto{\pgfqpoint{4.601228in}{1.013411in}}%
\pgfpathlineto{\pgfqpoint{4.660717in}{0.912183in}}%
\pgfpathlineto{\pgfqpoint{4.702180in}{0.847229in}}%
\pgfpathlineto{\pgfqpoint{4.736431in}{0.798879in}}%
\pgfpathlineto{\pgfqpoint{4.767077in}{0.760768in}}%
\pgfpathlineto{\pgfqpoint{4.794117in}{0.731733in}}%
\pgfpathlineto{\pgfqpoint{4.821158in}{0.707257in}}%
\pgfpathlineto{\pgfqpoint{4.846396in}{0.688521in}}%
\pgfpathlineto{\pgfqpoint{4.871633in}{0.673554in}}%
\pgfpathlineto{\pgfqpoint{4.896871in}{0.662003in}}%
\pgfpathlineto{\pgfqpoint{4.923912in}{0.652889in}}%
\pgfpathlineto{\pgfqpoint{4.954558in}{0.645851in}}%
\pgfpathlineto{\pgfqpoint{4.990612in}{0.640853in}}%
\pgfpathlineto{\pgfqpoint{5.037482in}{0.637683in}}%
\pgfpathlineto{\pgfqpoint{5.111393in}{0.636172in}}%
\pgfpathlineto{\pgfqpoint{5.268228in}{0.635913in}}%
\pgfpathlineto{\pgfqpoint{5.268228in}{0.635913in}}%
\pgfusepath{stroke}%
\end{pgfscope}%
\begin{pgfscope}%
\pgfpathrectangle{\pgfqpoint{0.653314in}{0.551844in}}{\pgfqpoint{4.614914in}{1.849515in}} %
\pgfusepath{clip}%
\pgfsetbuttcap%
\pgfsetroundjoin%
\pgfsetlinewidth{1.003750pt}%
\definecolor{currentstroke}{rgb}{0.000000,0.000000,0.000000}%
\pgfsetstrokecolor{currentstroke}%
\pgfsetstrokeopacity{0.500000}%
\pgfsetdash{{1.000000pt}{1.650000pt}}{0.000000pt}%
\pgfpathmoveto{\pgfqpoint{2.499279in}{0.551844in}}%
\pgfpathlineto{\pgfqpoint{2.499279in}{2.401359in}}%
\pgfusepath{stroke}%
\end{pgfscope}%
\begin{pgfscope}%
\pgfpathrectangle{\pgfqpoint{0.653314in}{0.551844in}}{\pgfqpoint{4.614914in}{1.849515in}} %
\pgfusepath{clip}%
\pgfsetbuttcap%
\pgfsetroundjoin%
\pgfsetlinewidth{1.003750pt}%
\definecolor{currentstroke}{rgb}{0.000000,0.000000,0.000000}%
\pgfsetstrokecolor{currentstroke}%
\pgfsetstrokeopacity{0.500000}%
\pgfsetdash{{1.000000pt}{1.650000pt}}{0.000000pt}%
\pgfpathmoveto{\pgfqpoint{3.422262in}{0.551844in}}%
\pgfpathlineto{\pgfqpoint{3.422262in}{2.401359in}}%
\pgfusepath{stroke}%
\end{pgfscope}%
\begin{pgfscope}%
\pgfsetrectcap%
\pgfsetmiterjoin%
\pgfsetlinewidth{0.803000pt}%
\definecolor{currentstroke}{rgb}{0.000000,0.000000,0.000000}%
\pgfsetstrokecolor{currentstroke}%
\pgfsetdash{}{0pt}%
\pgfpathmoveto{\pgfqpoint{0.653314in}{0.551844in}}%
\pgfpathlineto{\pgfqpoint{0.653314in}{2.401359in}}%
\pgfusepath{stroke}%
\end{pgfscope}%
\begin{pgfscope}%
\pgfsetrectcap%
\pgfsetmiterjoin%
\pgfsetlinewidth{0.803000pt}%
\definecolor{currentstroke}{rgb}{0.000000,0.000000,0.000000}%
\pgfsetstrokecolor{currentstroke}%
\pgfsetdash{}{0pt}%
\pgfpathmoveto{\pgfqpoint{5.268228in}{0.551844in}}%
\pgfpathlineto{\pgfqpoint{5.268228in}{2.401359in}}%
\pgfusepath{stroke}%
\end{pgfscope}%
\begin{pgfscope}%
\pgfsetrectcap%
\pgfsetmiterjoin%
\pgfsetlinewidth{0.803000pt}%
\definecolor{currentstroke}{rgb}{0.000000,0.000000,0.000000}%
\pgfsetstrokecolor{currentstroke}%
\pgfsetdash{}{0pt}%
\pgfpathmoveto{\pgfqpoint{0.653314in}{0.551844in}}%
\pgfpathlineto{\pgfqpoint{5.268228in}{0.551844in}}%
\pgfusepath{stroke}%
\end{pgfscope}%
\begin{pgfscope}%
\pgfsetrectcap%
\pgfsetmiterjoin%
\pgfsetlinewidth{0.803000pt}%
\definecolor{currentstroke}{rgb}{0.000000,0.000000,0.000000}%
\pgfsetstrokecolor{currentstroke}%
\pgfsetdash{}{0pt}%
\pgfpathmoveto{\pgfqpoint{0.653314in}{2.401359in}}%
\pgfpathlineto{\pgfqpoint{5.268228in}{2.401359in}}%
\pgfusepath{stroke}%
\end{pgfscope}%
\begin{pgfscope}%
\pgftext[x=2.960771in,y=2.484692in,,base]{\rmfamily\fontsize{12.000000}{14.400000}\selectfont Asymptotic expansion of solutions for small \(\displaystyle \alpha\)}%
\end{pgfscope}%
\begin{pgfscope}%
\pgfsetrectcap%
\pgfsetroundjoin%
\pgfsetlinewidth{1.505625pt}%
\definecolor{currentstroke}{rgb}{0.000000,0.000000,0.000000}%
\pgfsetstrokecolor{currentstroke}%
\pgfsetdash{}{0pt}%
\pgfpathmoveto{\pgfqpoint{2.613173in}{1.332065in}}%
\pgfpathlineto{\pgfqpoint{2.835395in}{1.332065in}}%
\pgfusepath{stroke}%
\end{pgfscope}%
\begin{pgfscope}%
\pgftext[x=2.924284in,y=1.293176in,left,base]{\rmfamily\fontsize{8.000000}{9.600000}\selectfont \(\displaystyle \alpha = 0\)}%
\end{pgfscope}%
\begin{pgfscope}%
\pgfsetbuttcap%
\pgfsetroundjoin%
\pgfsetlinewidth{1.003750pt}%
\definecolor{currentstroke}{rgb}{0.121569,0.466667,0.705882}%
\pgfsetstrokecolor{currentstroke}%
\pgfsetdash{{3.700000pt}{1.600000pt}}{0.000000pt}%
\pgfpathmoveto{\pgfqpoint{2.613173in}{1.177132in}}%
\pgfpathlineto{\pgfqpoint{2.835395in}{1.177132in}}%
\pgfusepath{stroke}%
\end{pgfscope}%
\begin{pgfscope}%
\pgftext[x=2.924284in,y=1.138243in,left,base]{\rmfamily\fontsize{8.000000}{9.600000}\selectfont \(\displaystyle \alpha = 0.1\)}%
\end{pgfscope}%
\begin{pgfscope}%
\pgfsetbuttcap%
\pgfsetroundjoin%
\pgfsetlinewidth{1.003750pt}%
\definecolor{currentstroke}{rgb}{1.000000,0.498039,0.054902}%
\pgfsetstrokecolor{currentstroke}%
\pgfsetdash{{6.400000pt}{1.600000pt}{1.000000pt}{1.600000pt}}{0.000000pt}%
\pgfpathmoveto{\pgfqpoint{2.613173in}{1.022199in}}%
\pgfpathlineto{\pgfqpoint{2.835395in}{1.022199in}}%
\pgfusepath{stroke}%
\end{pgfscope}%
\begin{pgfscope}%
\pgftext[x=2.924284in,y=0.983310in,left,base]{\rmfamily\fontsize{8.000000}{9.600000}\selectfont \(\displaystyle \alpha = 0.2\)}%
\end{pgfscope}%
\begin{pgfscope}%
\pgfsetbuttcap%
\pgfsetroundjoin%
\pgfsetlinewidth{1.003750pt}%
\definecolor{currentstroke}{rgb}{0.172549,0.627451,0.172549}%
\pgfsetstrokecolor{currentstroke}%
\pgfsetdash{{1.000000pt}{1.650000pt}}{0.000000pt}%
\pgfpathmoveto{\pgfqpoint{2.613173in}{0.867265in}}%
\pgfpathlineto{\pgfqpoint{2.835395in}{0.867265in}}%
\pgfusepath{stroke}%
\end{pgfscope}%
\begin{pgfscope}%
\pgftext[x=2.924284in,y=0.828377in,left,base]{\rmfamily\fontsize{8.000000}{9.600000}\selectfont \(\displaystyle \alpha = 0.5\)}%
\end{pgfscope}%
\begin{pgfscope}%
\pgfsetrectcap%
\pgfsetroundjoin%
\pgfsetlinewidth{1.003750pt}%
\definecolor{currentstroke}{rgb}{0.839216,0.152941,0.156863}%
\pgfsetstrokecolor{currentstroke}%
\pgfsetdash{}{0pt}%
\pgfpathmoveto{\pgfqpoint{2.613173in}{0.712332in}}%
\pgfpathlineto{\pgfqpoint{2.835395in}{0.712332in}}%
\pgfusepath{stroke}%
\end{pgfscope}%
\begin{pgfscope}%
\pgftext[x=2.924284in,y=0.673444in,left,base]{\rmfamily\fontsize{8.000000}{9.600000}\selectfont \(\displaystyle \alpha = 1.0\)}%
\end{pgfscope}%
\end{pgfpicture}%
\makeatother%
\endgroup%

%% file: SamplingDomainNaive.tikz
\begin{tikzpicture}

\begin{axis}[
title={The Sampling Domain},
xlabel={Space},
ylabel={Heading},
xmin=-0.1, xmax=3.1,
ymin=-1.3, ymax=1.3,
width=12.6cm,
height=8cm,
xtick={0,1.5,3},
xticklabels={$0$,$L/2$,$L$},
ytick={-1,0,1},
yticklabels={$-R$,$0$,$R$},
tick align=outside,
tick pos=left,
x grid style={lightgray!92.026143790849673!black},
y grid style={lightgray!92.026143790849673!black}
]
\path [draw=black, draw opacity=0.75, line width=0.83630769230769231pt, dash pattern=on 2pt off 3pt] (axis cs:0,-1.3)
--(axis cs:0,1.3);

\path [draw=black, draw opacity=0.75, line width=0.83630769230769231pt, dash pattern=on 2pt off 3pt] (axis cs:3,-1.3)
--(axis cs:3,1.3);

\path [draw=black, draw opacity=0.75, line width=0.83630769230769231pt, dash pattern=on 2pt off 3pt] (axis cs:-0.5,-1)
--(axis cs:3.5,-1);

\path [draw=black, draw opacity=0.75, line width=0.83630769230769231pt, dash pattern=on 2pt off 3pt] (axis cs:-0.5,0)
--(axis cs:3.5,0);

\path [draw=black, draw opacity=0.75, line width=0.83630769230769231pt, dash pattern=on 2pt off 3pt] (axis cs:-0.5,1)
--(axis cs:3.5,1);

\path [draw=lightgray, fill=lightgray] (axis cs:0,1)
--(axis cs:0,0)
--(axis cs:0.015625,-0.015625)
--(axis cs:0.03125,-0.03125)
--(axis cs:0.046875,-0.046875)
--(axis cs:0.0625,-0.0625)
--(axis cs:0.078125,-0.078125)
--(axis cs:0.09375,-0.09375)
--(axis cs:0.109375,-0.109375)
--(axis cs:0.125,-0.125)
--(axis cs:0.140625,-0.140625)
--(axis cs:0.15625,-0.15625)
--(axis cs:0.171875,-0.171875)
--(axis cs:0.1875,-0.1875)
--(axis cs:0.203125,-0.203125)
--(axis cs:0.21875,-0.21875)
--(axis cs:0.234375,-0.234375)
--(axis cs:0.25,-0.25)
--(axis cs:0.265625,-0.265625)
--(axis cs:0.28125,-0.28125)
--(axis cs:0.296875,-0.296875)
--(axis cs:0.3125,-0.3125)
--(axis cs:0.328125,-0.328125)
--(axis cs:0.34375,-0.34375)
--(axis cs:0.359375,-0.359375)
--(axis cs:0.375,-0.375)
--(axis cs:0.390625,-0.390625)
--(axis cs:0.40625,-0.40625)
--(axis cs:0.421875,-0.421875)
--(axis cs:0.4375,-0.4375)
--(axis cs:0.453125,-0.453125)
--(axis cs:0.46875,-0.46875)
--(axis cs:0.484375,-0.484375)
--(axis cs:0.5,-0.5)
--(axis cs:0.515625,-0.515625)
--(axis cs:0.53125,-0.53125)
--(axis cs:0.546875,-0.546875)
--(axis cs:0.5625,-0.5625)
--(axis cs:0.578125,-0.578125)
--(axis cs:0.59375,-0.59375)
--(axis cs:0.609375,-0.609375)
--(axis cs:0.625,-0.625)
--(axis cs:0.640625,-0.640625)
--(axis cs:0.65625,-0.65625)
--(axis cs:0.671875,-0.671875)
--(axis cs:0.6875,-0.6875)
--(axis cs:0.703125,-0.703125)
--(axis cs:0.71875,-0.71875)
--(axis cs:0.734375,-0.734375)
--(axis cs:0.75,-0.75)
--(axis cs:0.765625,-0.765625)
--(axis cs:0.78125,-0.78125)
--(axis cs:0.796875,-0.796875)
--(axis cs:0.8125,-0.8125)
--(axis cs:0.828125,-0.828125)
--(axis cs:0.84375,-0.84375)
--(axis cs:0.859375,-0.859375)
--(axis cs:0.875,-0.875)
--(axis cs:0.890625,-0.890625)
--(axis cs:0.90625,-0.90625)
--(axis cs:0.921875,-0.921875)
--(axis cs:0.9375,-0.9375)
--(axis cs:0.953125,-0.953125)
--(axis cs:0.96875,-0.96875)
--(axis cs:0.984375,-0.984375)
--(axis cs:1,-1)
--(axis cs:1.015625,-1)
--(axis cs:1.03125,-1)
--(axis cs:1.046875,-1)
--(axis cs:1.0625,-1)
--(axis cs:1.078125,-1)
--(axis cs:1.09375,-1)
--(axis cs:1.109375,-1)
--(axis cs:1.125,-1)
--(axis cs:1.140625,-1)
--(axis cs:1.15625,-1)
--(axis cs:1.171875,-1)
--(axis cs:1.1875,-1)
--(axis cs:1.203125,-1)
--(axis cs:1.21875,-1)
--(axis cs:1.234375,-1)
--(axis cs:1.25,-1)
--(axis cs:1.265625,-1)
--(axis cs:1.28125,-1)
--(axis cs:1.296875,-1)
--(axis cs:1.3125,-1)
--(axis cs:1.328125,-1)
--(axis cs:1.34375,-1)
--(axis cs:1.359375,-1)
--(axis cs:1.375,-1)
--(axis cs:1.390625,-1)
--(axis cs:1.40625,-1)
--(axis cs:1.421875,-1)
--(axis cs:1.4375,-1)
--(axis cs:1.453125,-1)
--(axis cs:1.46875,-1)
--(axis cs:1.484375,-1)
--(axis cs:1.5,-1)
--(axis cs:1.515625,-1)
--(axis cs:1.53125,-1)
--(axis cs:1.546875,-1)
--(axis cs:1.5625,-1)
--(axis cs:1.578125,-1)
--(axis cs:1.59375,-1)
--(axis cs:1.609375,-1)
--(axis cs:1.625,-1)
--(axis cs:1.640625,-1)
--(axis cs:1.65625,-1)
--(axis cs:1.671875,-1)
--(axis cs:1.6875,-1)
--(axis cs:1.703125,-1)
--(axis cs:1.71875,-1)
--(axis cs:1.734375,-1)
--(axis cs:1.75,-1)
--(axis cs:1.765625,-1)
--(axis cs:1.78125,-1)
--(axis cs:1.796875,-1)
--(axis cs:1.8125,-1)
--(axis cs:1.828125,-1)
--(axis cs:1.84375,-1)
--(axis cs:1.859375,-1)
--(axis cs:1.875,-1)
--(axis cs:1.890625,-1)
--(axis cs:1.90625,-1)
--(axis cs:1.921875,-1)
--(axis cs:1.9375,-1)
--(axis cs:1.953125,-1)
--(axis cs:1.96875,-1)
--(axis cs:1.984375,-1)
--(axis cs:2,-1)
--(axis cs:2.015625,-1)
--(axis cs:2.03125,-1)
--(axis cs:2.046875,-1)
--(axis cs:2.0625,-1)
--(axis cs:2.078125,-1)
--(axis cs:2.09375,-1)
--(axis cs:2.109375,-1)
--(axis cs:2.125,-1)
--(axis cs:2.140625,-1)
--(axis cs:2.15625,-1)
--(axis cs:2.171875,-1)
--(axis cs:2.1875,-1)
--(axis cs:2.203125,-1)
--(axis cs:2.21875,-1)
--(axis cs:2.234375,-1)
--(axis cs:2.25,-1)
--(axis cs:2.265625,-1)
--(axis cs:2.28125,-1)
--(axis cs:2.296875,-1)
--(axis cs:2.3125,-1)
--(axis cs:2.328125,-1)
--(axis cs:2.34375,-1)
--(axis cs:2.359375,-1)
--(axis cs:2.375,-1)
--(axis cs:2.390625,-1)
--(axis cs:2.40625,-1)
--(axis cs:2.421875,-1)
--(axis cs:2.4375,-1)
--(axis cs:2.453125,-1)
--(axis cs:2.46875,-1)
--(axis cs:2.484375,-1)
--(axis cs:2.5,-1)
--(axis cs:2.515625,-1)
--(axis cs:2.53125,-1)
--(axis cs:2.546875,-1)
--(axis cs:2.5625,-1)
--(axis cs:2.578125,-1)
--(axis cs:2.59375,-1)
--(axis cs:2.609375,-1)
--(axis cs:2.625,-1)
--(axis cs:2.640625,-1)
--(axis cs:2.65625,-1)
--(axis cs:2.671875,-1)
--(axis cs:2.6875,-1)
--(axis cs:2.703125,-1)
--(axis cs:2.71875,-1)
--(axis cs:2.734375,-1)
--(axis cs:2.75,-1)
--(axis cs:2.765625,-1)
--(axis cs:2.78125,-1)
--(axis cs:2.796875,-1)
--(axis cs:2.8125,-1)
--(axis cs:2.828125,-1)
--(axis cs:2.84375,-1)
--(axis cs:2.859375,-1)
--(axis cs:2.875,-1)
--(axis cs:2.890625,-1)
--(axis cs:2.90625,-1)
--(axis cs:2.921875,-1)
--(axis cs:2.9375,-1)
--(axis cs:2.953125,-1)
--(axis cs:2.96875,-1)
--(axis cs:2.984375,-1)
--(axis cs:3,-1)
--(axis cs:3,0)
--(axis cs:3,0)
--(axis cs:2.984375,0.015625)
--(axis cs:2.96875,0.03125)
--(axis cs:2.953125,0.046875)
--(axis cs:2.9375,0.0625)
--(axis cs:2.921875,0.078125)
--(axis cs:2.90625,0.09375)
--(axis cs:2.890625,0.109375)
--(axis cs:2.875,0.125)
--(axis cs:2.859375,0.140625)
--(axis cs:2.84375,0.15625)
--(axis cs:2.828125,0.171875)
--(axis cs:2.8125,0.1875)
--(axis cs:2.796875,0.203125)
--(axis cs:2.78125,0.21875)
--(axis cs:2.765625,0.234375)
--(axis cs:2.75,0.25)
--(axis cs:2.734375,0.265625)
--(axis cs:2.71875,0.28125)
--(axis cs:2.703125,0.296875)
--(axis cs:2.6875,0.3125)
--(axis cs:2.671875,0.328125)
--(axis cs:2.65625,0.34375)
--(axis cs:2.640625,0.359375)
--(axis cs:2.625,0.375)
--(axis cs:2.609375,0.390625)
--(axis cs:2.59375,0.40625)
--(axis cs:2.578125,0.421875)
--(axis cs:2.5625,0.4375)
--(axis cs:2.546875,0.453125)
--(axis cs:2.53125,0.46875)
--(axis cs:2.515625,0.484375)
--(axis cs:2.5,0.5)
--(axis cs:2.484375,0.515625)
--(axis cs:2.46875,0.53125)
--(axis cs:2.453125,0.546875)
--(axis cs:2.4375,0.5625)
--(axis cs:2.421875,0.578125)
--(axis cs:2.40625,0.59375)
--(axis cs:2.390625,0.609375)
--(axis cs:2.375,0.625)
--(axis cs:2.359375,0.640625)
--(axis cs:2.34375,0.65625)
--(axis cs:2.328125,0.671875)
--(axis cs:2.3125,0.6875)
--(axis cs:2.296875,0.703125)
--(axis cs:2.28125,0.71875)
--(axis cs:2.265625,0.734375)
--(axis cs:2.25,0.75)
--(axis cs:2.234375,0.765625)
--(axis cs:2.21875,0.78125)
--(axis cs:2.203125,0.796875)
--(axis cs:2.1875,0.8125)
--(axis cs:2.171875,0.828125)
--(axis cs:2.15625,0.84375)
--(axis cs:2.140625,0.859375)
--(axis cs:2.125,0.875)
--(axis cs:2.109375,0.890625)
--(axis cs:2.09375,0.90625)
--(axis cs:2.078125,0.921875)
--(axis cs:2.0625,0.9375)
--(axis cs:2.046875,0.953125)
--(axis cs:2.03125,0.96875)
--(axis cs:2.015625,0.984375)
--(axis cs:2,1)
--(axis cs:1.984375,1)
--(axis cs:1.96875,1)
--(axis cs:1.953125,1)
--(axis cs:1.9375,1)
--(axis cs:1.921875,1)
--(axis cs:1.90625,1)
--(axis cs:1.890625,1)
--(axis cs:1.875,1)
--(axis cs:1.859375,1)
--(axis cs:1.84375,1)
--(axis cs:1.828125,1)
--(axis cs:1.8125,1)
--(axis cs:1.796875,1)
--(axis cs:1.78125,1)
--(axis cs:1.765625,1)
--(axis cs:1.75,1)
--(axis cs:1.734375,1)
--(axis cs:1.71875,1)
--(axis cs:1.703125,1)
--(axis cs:1.6875,1)
--(axis cs:1.671875,1)
--(axis cs:1.65625,1)
--(axis cs:1.640625,1)
--(axis cs:1.625,1)
--(axis cs:1.609375,1)
--(axis cs:1.59375,1)
--(axis cs:1.578125,1)
--(axis cs:1.5625,1)
--(axis cs:1.546875,1)
--(axis cs:1.53125,1)
--(axis cs:1.515625,1)
--(axis cs:1.5,1)
--(axis cs:1.484375,1)
--(axis cs:1.46875,1)
--(axis cs:1.453125,1)
--(axis cs:1.4375,1)
--(axis cs:1.421875,1)
--(axis cs:1.40625,1)
--(axis cs:1.390625,1)
--(axis cs:1.375,1)
--(axis cs:1.359375,1)
--(axis cs:1.34375,1)
--(axis cs:1.328125,1)
--(axis cs:1.3125,1)
--(axis cs:1.296875,1)
--(axis cs:1.28125,1)
--(axis cs:1.265625,1)
--(axis cs:1.25,1)
--(axis cs:1.234375,1)
--(axis cs:1.21875,1)
--(axis cs:1.203125,1)
--(axis cs:1.1875,1)
--(axis cs:1.171875,1)
--(axis cs:1.15625,1)
--(axis cs:1.140625,1)
--(axis cs:1.125,1)
--(axis cs:1.109375,1)
--(axis cs:1.09375,1)
--(axis cs:1.078125,1)
--(axis cs:1.0625,1)
--(axis cs:1.046875,1)
--(axis cs:1.03125,1)
--(axis cs:1.015625,1)
--(axis cs:1,1)
--(axis cs:0.984375,1)
--(axis cs:0.96875,1)
--(axis cs:0.953125,1)
--(axis cs:0.9375,1)
--(axis cs:0.921875,1)
--(axis cs:0.90625,1)
--(axis cs:0.890625,1)
--(axis cs:0.875,1)
--(axis cs:0.859375,1)
--(axis cs:0.84375,1)
--(axis cs:0.828125,1)
--(axis cs:0.8125,1)
--(axis cs:0.796875,1)
--(axis cs:0.78125,1)
--(axis cs:0.765625,1)
--(axis cs:0.75,1)
--(axis cs:0.734375,1)
--(axis cs:0.71875,1)
--(axis cs:0.703125,1)
--(axis cs:0.6875,1)
--(axis cs:0.671875,1)
--(axis cs:0.65625,1)
--(axis cs:0.640625,1)
--(axis cs:0.625,1)
--(axis cs:0.609375,1)
--(axis cs:0.59375,1)
--(axis cs:0.578125,1)
--(axis cs:0.5625,1)
--(axis cs:0.546875,1)
--(axis cs:0.53125,1)
--(axis cs:0.515625,1)
--(axis cs:0.5,1)
--(axis cs:0.484375,1)
--(axis cs:0.46875,1)
--(axis cs:0.453125,1)
--(axis cs:0.4375,1)
--(axis cs:0.421875,1)
--(axis cs:0.40625,1)
--(axis cs:0.390625,1)
--(axis cs:0.375,1)
--(axis cs:0.359375,1)
--(axis cs:0.34375,1)
--(axis cs:0.328125,1)
--(axis cs:0.3125,1)
--(axis cs:0.296875,1)
--(axis cs:0.28125,1)
--(axis cs:0.265625,1)
--(axis cs:0.25,1)
--(axis cs:0.234375,1)
--(axis cs:0.21875,1)
--(axis cs:0.203125,1)
--(axis cs:0.1875,1)
--(axis cs:0.171875,1)
--(axis cs:0.15625,1)
--(axis cs:0.140625,1)
--(axis cs:0.125,1)
--(axis cs:0.109375,1)
--(axis cs:0.09375,1)
--(axis cs:0.078125,1)
--(axis cs:0.0625,1)
--(axis cs:0.046875,1)
--(axis cs:0.03125,1)
--(axis cs:0.015625,1)
--(axis cs:0,1)
--cycle;

\path [draw=darkgray, fill=darkgray] (axis cs:2.203125,0.796875)
--(axis cs:2.203125,-1)
--(axis cs:2.21875,-1)
--(axis cs:2.234375,-1)
--(axis cs:2.234375,0.765625)
--(axis cs:2.234375,0.765625)
--(axis cs:2.21875,0.78125)
--(axis cs:2.203125,0.796875)
--cycle;

\path [draw=black, line width=1.6726153846153846pt] (axis cs:0,0)
--(axis cs:0,1);

\path [draw=black, line width=1.6726153846153846pt] (axis cs:3,-1)
--(axis cs:3,0);

\path [draw=black, line width=0.83630769230769231pt] (axis cs:2.203125,-1)
--(axis cs:2.203125,0.796875);

\path [draw=black, line width=0.83630769230769231pt] (axis cs:2.234375,-1)
--(axis cs:2.234375,0.765625);

\addplot [line width=1.6726153846153846pt, black, forget plot]
table {%
0 1
0.015625 1
0.03125 1
0.046875 1
0.0625 1
0.078125 1
0.09375 1
0.109375 1
0.125 1
0.140625 1
0.15625 1
0.171875 1
0.1875 1
0.203125 1
0.21875 1
0.234375 1
0.25 1
0.265625 1
0.28125 1
0.296875 1
0.3125 1
0.328125 1
0.34375 1
0.359375 1
0.375 1
0.390625 1
0.40625 1
0.421875 1
0.4375 1
0.453125 1
0.46875 1
0.484375 1
0.5 1
0.515625 1
0.53125 1
0.546875 1
0.5625 1
0.578125 1
0.59375 1
0.609375 1
0.625 1
0.640625 1
0.65625 1
0.671875 1
0.6875 1
0.703125 1
0.71875 1
0.734375 1
0.75 1
0.765625 1
0.78125 1
0.796875 1
0.8125 1
0.828125 1
0.84375 1
0.859375 1
0.875 1
0.890625 1
0.90625 1
0.921875 1
0.9375 1
0.953125 1
0.96875 1
0.984375 1
1 1
1.015625 1
1.03125 1
1.046875 1
1.0625 1
1.078125 1
1.09375 1
1.109375 1
1.125 1
1.140625 1
1.15625 1
1.171875 1
1.1875 1
1.203125 1
1.21875 1
1.234375 1
1.25 1
1.265625 1
1.28125 1
1.296875 1
1.3125 1
1.328125 1
1.34375 1
1.359375 1
1.375 1
1.390625 1
1.40625 1
1.421875 1
1.4375 1
1.453125 1
1.46875 1
1.484375 1
1.5 1
1.515625 1
1.53125 1
1.546875 1
1.5625 1
1.578125 1
1.59375 1
1.609375 1
1.625 1
1.640625 1
1.65625 1
1.671875 1
1.6875 1
1.703125 1
1.71875 1
1.734375 1
1.75 1
1.765625 1
1.78125 1
1.796875 1
1.8125 1
1.828125 1
1.84375 1
1.859375 1
1.875 1
1.890625 1
1.90625 1
1.921875 1
1.9375 1
1.953125 1
1.96875 1
1.984375 1
2 1
2.015625 0.984375
2.03125 0.96875
2.046875 0.953125
2.0625 0.9375
2.078125 0.921875
2.09375 0.90625
2.109375 0.890625
2.125 0.875
2.140625 0.859375
2.15625 0.84375
2.171875 0.828125
2.1875 0.8125
2.203125 0.796875
2.21875 0.78125
2.234375 0.765625
2.25 0.75
2.265625 0.734375
2.28125 0.71875
2.296875 0.703125
2.3125 0.6875
2.328125 0.671875
2.34375 0.65625
2.359375 0.640625
2.375 0.625
2.390625 0.609375
2.40625 0.59375
2.421875 0.578125
2.4375 0.5625
2.453125 0.546875
2.46875 0.53125
2.484375 0.515625
2.5 0.5
2.515625 0.484375
2.53125 0.46875
2.546875 0.453125
2.5625 0.4375
2.578125 0.421875
2.59375 0.40625
2.609375 0.390625
2.625 0.375
2.640625 0.359375
2.65625 0.34375
2.671875 0.328125
2.6875 0.3125
2.703125 0.296875
2.71875 0.28125
2.734375 0.265625
2.75 0.25
2.765625 0.234375
2.78125 0.21875
2.796875 0.203125
2.8125 0.1875
2.828125 0.171875
2.84375 0.15625
2.859375 0.140625
2.875 0.125
2.890625 0.109375
2.90625 0.09375
2.921875 0.078125
2.9375 0.0625
2.953125 0.046875
2.96875 0.03125
2.984375 0.015625
3 0
};
\addplot [line width=1.6726153846153846pt, black, forget plot]
table {%
0 -0
0.015625 -0.015625
0.03125 -0.03125
0.046875 -0.046875
0.0625 -0.0625
0.078125 -0.078125
0.09375 -0.09375
0.109375 -0.109375
0.125 -0.125
0.140625 -0.140625
0.15625 -0.15625
0.171875 -0.171875
0.1875 -0.1875
0.203125 -0.203125
0.21875 -0.21875
0.234375 -0.234375
0.25 -0.25
0.265625 -0.265625
0.28125 -0.28125
0.296875 -0.296875
0.3125 -0.3125
0.328125 -0.328125
0.34375 -0.34375
0.359375 -0.359375
0.375 -0.375
0.390625 -0.390625
0.40625 -0.40625
0.421875 -0.421875
0.4375 -0.4375
0.453125 -0.453125
0.46875 -0.46875
0.484375 -0.484375
0.5 -0.5
0.515625 -0.515625
0.53125 -0.53125
0.546875 -0.546875
0.5625 -0.5625
0.578125 -0.578125
0.59375 -0.59375
0.609375 -0.609375
0.625 -0.625
0.640625 -0.640625
0.65625 -0.65625
0.671875 -0.671875
0.6875 -0.6875
0.703125 -0.703125
0.71875 -0.71875
0.734375 -0.734375
0.75 -0.75
0.765625 -0.765625
0.78125 -0.78125
0.796875 -0.796875
0.8125 -0.8125
0.828125 -0.828125
0.84375 -0.84375
0.859375 -0.859375
0.875 -0.875
0.890625 -0.890625
0.90625 -0.90625
0.921875 -0.921875
0.9375 -0.9375
0.953125 -0.953125
0.96875 -0.96875
0.984375 -0.984375
1 -1
1.015625 -1
1.03125 -1
1.046875 -1
1.0625 -1
1.078125 -1
1.09375 -1
1.109375 -1
1.125 -1
1.140625 -1
1.15625 -1
1.171875 -1
1.1875 -1
1.203125 -1
1.21875 -1
1.234375 -1
1.25 -1
1.265625 -1
1.28125 -1
1.296875 -1
1.3125 -1
1.328125 -1
1.34375 -1
1.359375 -1
1.375 -1
1.390625 -1
1.40625 -1
1.421875 -1
1.4375 -1
1.453125 -1
1.46875 -1
1.484375 -1
1.5 -1
1.515625 -1
1.53125 -1
1.546875 -1
1.5625 -1
1.578125 -1
1.59375 -1
1.609375 -1
1.625 -1
1.640625 -1
1.65625 -1
1.671875 -1
1.6875 -1
1.703125 -1
1.71875 -1
1.734375 -1
1.75 -1
1.765625 -1
1.78125 -1
1.796875 -1
1.8125 -1
1.828125 -1
1.84375 -1
1.859375 -1
1.875 -1
1.890625 -1
1.90625 -1
1.921875 -1
1.9375 -1
1.953125 -1
1.96875 -1
1.984375 -1
2 -1
2.015625 -1
2.03125 -1
2.046875 -1
2.0625 -1
2.078125 -1
2.09375 -1
2.109375 -1
2.125 -1
2.140625 -1
2.15625 -1
2.171875 -1
2.1875 -1
2.203125 -1
2.21875 -1
2.234375 -1
2.25 -1
2.265625 -1
2.28125 -1
2.296875 -1
2.3125 -1
2.328125 -1
2.34375 -1
2.359375 -1
2.375 -1
2.390625 -1
2.40625 -1
2.421875 -1
2.4375 -1
2.453125 -1
2.46875 -1
2.484375 -1
2.5 -1
2.515625 -1
2.53125 -1
2.546875 -1
2.5625 -1
2.578125 -1
2.59375 -1
2.609375 -1
2.625 -1
2.640625 -1
2.65625 -1
2.671875 -1
2.6875 -1
2.703125 -1
2.71875 -1
2.734375 -1
2.75 -1
2.765625 -1
2.78125 -1
2.796875 -1
2.8125 -1
2.828125 -1
2.84375 -1
2.859375 -1
2.875 -1
2.890625 -1
2.90625 -1
2.921875 -1
2.9375 -1
2.953125 -1
2.96875 -1
2.984375 -1
3 -1
};
\node at (axis cs:2.115,-1.2)[
  scale=0.5,
  anchor=base west,
  text=black,
  rotate=0.0
]{ $E(x)$};
\end{axis}

\end{tikzpicture}

%% file: SamplingDomainNoFlux.tikz
\begin{tikzpicture}

\begin{axis}[
title={The Sampling Domain},
xlabel={Space},
ylabel={Heading},
xmin=-0.1, xmax=3.1,
ymin=-1.3, ymax=1.3,
width=12.6cm,
height=8cm,
xtick={0,1.5,3},
xticklabels={$0$,$L/2$,$L$},
ytick={-1,0,1},
yticklabels={$-R$,$0$,$R$},
tick align=outside,
tick pos=left,
x grid style={lightgray!92.026143790849673!black},
y grid style={lightgray!92.026143790849673!black}
]
\path [draw=black, draw opacity=0.75, line width=0.83630769230769231pt, dash pattern=on 2pt off 3pt] (axis cs:0,-1.3)
--(axis cs:0,1.3);

\path [draw=black, draw opacity=0.75, line width=0.83630769230769231pt, dash pattern=on 2pt off 3pt] (axis cs:3,-1.3)
--(axis cs:3,1.3);

\path [draw=black, draw opacity=0.75, line width=0.83630769230769231pt, dash pattern=on 2pt off 3pt] (axis cs:-0.5,-1)
--(axis cs:3.5,-1);

\path [draw=black, draw opacity=0.75, line width=0.83630769230769231pt, dash pattern=on 2pt off 3pt] (axis cs:-0.5,0)
--(axis cs:3.5,0);

\path [draw=black, draw opacity=0.75, line width=0.83630769230769231pt, dash pattern=on 2pt off 3pt] (axis cs:-0.5,1)
--(axis cs:3.5,1);

\path [draw=lightgray, fill=lightgray] (axis cs:0,1)
--(axis cs:0,1)
--(axis cs:0.015625,0.96875)
--(axis cs:0.03125,0.9375)
--(axis cs:0.046875,0.90625)
--(axis cs:0.0625,0.875)
--(axis cs:0.078125,0.84375)
--(axis cs:0.09375,0.8125)
--(axis cs:0.109375,0.78125)
--(axis cs:0.125,0.75)
--(axis cs:0.140625,0.71875)
--(axis cs:0.15625,0.6875)
--(axis cs:0.171875,0.65625)
--(axis cs:0.1875,0.625)
--(axis cs:0.203125,0.59375)
--(axis cs:0.21875,0.5625)
--(axis cs:0.234375,0.53125)
--(axis cs:0.25,0.5)
--(axis cs:0.265625,0.46875)
--(axis cs:0.28125,0.4375)
--(axis cs:0.296875,0.40625)
--(axis cs:0.3125,0.375)
--(axis cs:0.328125,0.34375)
--(axis cs:0.34375,0.3125)
--(axis cs:0.359375,0.28125)
--(axis cs:0.375,0.25)
--(axis cs:0.390625,0.21875)
--(axis cs:0.40625,0.1875)
--(axis cs:0.421875,0.15625)
--(axis cs:0.4375,0.125)
--(axis cs:0.453125,0.09375)
--(axis cs:0.46875,0.0625)
--(axis cs:0.484375,0.03125)
--(axis cs:0.5,0)
--(axis cs:0.515625,-0.03125)
--(axis cs:0.53125,-0.0625)
--(axis cs:0.546875,-0.09375)
--(axis cs:0.5625,-0.125)
--(axis cs:0.578125,-0.15625)
--(axis cs:0.59375,-0.1875)
--(axis cs:0.609375,-0.21875)
--(axis cs:0.625,-0.25)
--(axis cs:0.640625,-0.28125)
--(axis cs:0.65625,-0.3125)
--(axis cs:0.671875,-0.34375)
--(axis cs:0.6875,-0.375)
--(axis cs:0.703125,-0.40625)
--(axis cs:0.71875,-0.4375)
--(axis cs:0.734375,-0.46875)
--(axis cs:0.75,-0.5)
--(axis cs:0.765625,-0.53125)
--(axis cs:0.78125,-0.5625)
--(axis cs:0.796875,-0.59375)
--(axis cs:0.8125,-0.625)
--(axis cs:0.828125,-0.65625)
--(axis cs:0.84375,-0.6875)
--(axis cs:0.859375,-0.71875)
--(axis cs:0.875,-0.75)
--(axis cs:0.890625,-0.78125)
--(axis cs:0.90625,-0.8125)
--(axis cs:0.921875,-0.84375)
--(axis cs:0.9375,-0.875)
--(axis cs:0.953125,-0.90625)
--(axis cs:0.96875,-0.9375)
--(axis cs:0.984375,-0.96875)
--(axis cs:1,-1)
--(axis cs:1.015625,-1)
--(axis cs:1.03125,-1)
--(axis cs:1.046875,-1)
--(axis cs:1.0625,-1)
--(axis cs:1.078125,-1)
--(axis cs:1.09375,-1)
--(axis cs:1.109375,-1)
--(axis cs:1.125,-1)
--(axis cs:1.140625,-1)
--(axis cs:1.15625,-1)
--(axis cs:1.171875,-1)
--(axis cs:1.1875,-1)
--(axis cs:1.203125,-1)
--(axis cs:1.21875,-1)
--(axis cs:1.234375,-1)
--(axis cs:1.25,-1)
--(axis cs:1.265625,-1)
--(axis cs:1.28125,-1)
--(axis cs:1.296875,-1)
--(axis cs:1.3125,-1)
--(axis cs:1.328125,-1)
--(axis cs:1.34375,-1)
--(axis cs:1.359375,-1)
--(axis cs:1.375,-1)
--(axis cs:1.390625,-1)
--(axis cs:1.40625,-1)
--(axis cs:1.421875,-1)
--(axis cs:1.4375,-1)
--(axis cs:1.453125,-1)
--(axis cs:1.46875,-1)
--(axis cs:1.484375,-1)
--(axis cs:1.5,-1)
--(axis cs:1.515625,-1)
--(axis cs:1.53125,-1)
--(axis cs:1.546875,-1)
--(axis cs:1.5625,-1)
--(axis cs:1.578125,-1)
--(axis cs:1.59375,-1)
--(axis cs:1.609375,-1)
--(axis cs:1.625,-1)
--(axis cs:1.640625,-1)
--(axis cs:1.65625,-1)
--(axis cs:1.671875,-1)
--(axis cs:1.6875,-1)
--(axis cs:1.703125,-1)
--(axis cs:1.71875,-1)
--(axis cs:1.734375,-1)
--(axis cs:1.75,-1)
--(axis cs:1.765625,-1)
--(axis cs:1.78125,-1)
--(axis cs:1.796875,-1)
--(axis cs:1.8125,-1)
--(axis cs:1.828125,-1)
--(axis cs:1.84375,-1)
--(axis cs:1.859375,-1)
--(axis cs:1.875,-1)
--(axis cs:1.890625,-1)
--(axis cs:1.90625,-1)
--(axis cs:1.921875,-1)
--(axis cs:1.9375,-1)
--(axis cs:1.953125,-1)
--(axis cs:1.96875,-1)
--(axis cs:1.984375,-1)
--(axis cs:2,-1)
--(axis cs:2.015625,-1)
--(axis cs:2.03125,-1)
--(axis cs:2.046875,-1)
--(axis cs:2.0625,-1)
--(axis cs:2.078125,-1)
--(axis cs:2.09375,-1)
--(axis cs:2.109375,-1)
--(axis cs:2.125,-1)
--(axis cs:2.140625,-1)
--(axis cs:2.15625,-1)
--(axis cs:2.171875,-1)
--(axis cs:2.1875,-1)
--(axis cs:2.203125,-1)
--(axis cs:2.21875,-1)
--(axis cs:2.234375,-1)
--(axis cs:2.25,-1)
--(axis cs:2.265625,-1)
--(axis cs:2.28125,-1)
--(axis cs:2.296875,-1)
--(axis cs:2.3125,-1)
--(axis cs:2.328125,-1)
--(axis cs:2.34375,-1)
--(axis cs:2.359375,-1)
--(axis cs:2.375,-1)
--(axis cs:2.390625,-1)
--(axis cs:2.40625,-1)
--(axis cs:2.421875,-1)
--(axis cs:2.4375,-1)
--(axis cs:2.453125,-1)
--(axis cs:2.46875,-1)
--(axis cs:2.484375,-1)
--(axis cs:2.5,-1)
--(axis cs:2.515625,-1)
--(axis cs:2.53125,-1)
--(axis cs:2.546875,-1)
--(axis cs:2.5625,-1)
--(axis cs:2.578125,-1)
--(axis cs:2.59375,-1)
--(axis cs:2.609375,-1)
--(axis cs:2.625,-1)
--(axis cs:2.640625,-1)
--(axis cs:2.65625,-1)
--(axis cs:2.671875,-1)
--(axis cs:2.6875,-1)
--(axis cs:2.703125,-1)
--(axis cs:2.71875,-1)
--(axis cs:2.734375,-1)
--(axis cs:2.75,-1)
--(axis cs:2.765625,-1)
--(axis cs:2.78125,-1)
--(axis cs:2.796875,-1)
--(axis cs:2.8125,-1)
--(axis cs:2.828125,-1)
--(axis cs:2.84375,-1)
--(axis cs:2.859375,-1)
--(axis cs:2.875,-1)
--(axis cs:2.890625,-1)
--(axis cs:2.90625,-1)
--(axis cs:2.921875,-1)
--(axis cs:2.9375,-1)
--(axis cs:2.953125,-1)
--(axis cs:2.96875,-1)
--(axis cs:2.984375,-1)
--(axis cs:3,-1)
--(axis cs:3,-1)
--(axis cs:3,-1)
--(axis cs:2.984375,-0.96875)
--(axis cs:2.96875,-0.9375)
--(axis cs:2.953125,-0.90625)
--(axis cs:2.9375,-0.875)
--(axis cs:2.921875,-0.84375)
--(axis cs:2.90625,-0.8125)
--(axis cs:2.890625,-0.78125)
--(axis cs:2.875,-0.75)
--(axis cs:2.859375,-0.71875)
--(axis cs:2.84375,-0.6875)
--(axis cs:2.828125,-0.65625)
--(axis cs:2.8125,-0.625)
--(axis cs:2.796875,-0.59375)
--(axis cs:2.78125,-0.5625)
--(axis cs:2.765625,-0.53125)
--(axis cs:2.75,-0.5)
--(axis cs:2.734375,-0.46875)
--(axis cs:2.71875,-0.4375)
--(axis cs:2.703125,-0.40625)
--(axis cs:2.6875,-0.375)
--(axis cs:2.671875,-0.34375)
--(axis cs:2.65625,-0.3125)
--(axis cs:2.640625,-0.28125)
--(axis cs:2.625,-0.25)
--(axis cs:2.609375,-0.21875)
--(axis cs:2.59375,-0.1875)
--(axis cs:2.578125,-0.15625)
--(axis cs:2.5625,-0.125)
--(axis cs:2.546875,-0.09375)
--(axis cs:2.53125,-0.0625)
--(axis cs:2.515625,-0.03125)
--(axis cs:2.5,0)
--(axis cs:2.484375,0.03125)
--(axis cs:2.46875,0.0625)
--(axis cs:2.453125,0.09375)
--(axis cs:2.4375,0.125)
--(axis cs:2.421875,0.15625)
--(axis cs:2.40625,0.1875)
--(axis cs:2.390625,0.21875)
--(axis cs:2.375,0.25)
--(axis cs:2.359375,0.28125)
--(axis cs:2.34375,0.3125)
--(axis cs:2.328125,0.34375)
--(axis cs:2.3125,0.375)
--(axis cs:2.296875,0.40625)
--(axis cs:2.28125,0.4375)
--(axis cs:2.265625,0.46875)
--(axis cs:2.25,0.5)
--(axis cs:2.234375,0.53125)
--(axis cs:2.21875,0.5625)
--(axis cs:2.203125,0.59375)
--(axis cs:2.1875,0.625)
--(axis cs:2.171875,0.65625)
--(axis cs:2.15625,0.6875)
--(axis cs:2.140625,0.71875)
--(axis cs:2.125,0.75)
--(axis cs:2.109375,0.78125)
--(axis cs:2.09375,0.8125)
--(axis cs:2.078125,0.84375)
--(axis cs:2.0625,0.875)
--(axis cs:2.046875,0.90625)
--(axis cs:2.03125,0.9375)
--(axis cs:2.015625,0.96875)
--(axis cs:2,1)
--(axis cs:1.984375,1)
--(axis cs:1.96875,1)
--(axis cs:1.953125,1)
--(axis cs:1.9375,1)
--(axis cs:1.921875,1)
--(axis cs:1.90625,1)
--(axis cs:1.890625,1)
--(axis cs:1.875,1)
--(axis cs:1.859375,1)
--(axis cs:1.84375,1)
--(axis cs:1.828125,1)
--(axis cs:1.8125,1)
--(axis cs:1.796875,1)
--(axis cs:1.78125,1)
--(axis cs:1.765625,1)
--(axis cs:1.75,1)
--(axis cs:1.734375,1)
--(axis cs:1.71875,1)
--(axis cs:1.703125,1)
--(axis cs:1.6875,1)
--(axis cs:1.671875,1)
--(axis cs:1.65625,1)
--(axis cs:1.640625,1)
--(axis cs:1.625,1)
--(axis cs:1.609375,1)
--(axis cs:1.59375,1)
--(axis cs:1.578125,1)
--(axis cs:1.5625,1)
--(axis cs:1.546875,1)
--(axis cs:1.53125,1)
--(axis cs:1.515625,1)
--(axis cs:1.5,1)
--(axis cs:1.484375,1)
--(axis cs:1.46875,1)
--(axis cs:1.453125,1)
--(axis cs:1.4375,1)
--(axis cs:1.421875,1)
--(axis cs:1.40625,1)
--(axis cs:1.390625,1)
--(axis cs:1.375,1)
--(axis cs:1.359375,1)
--(axis cs:1.34375,1)
--(axis cs:1.328125,1)
--(axis cs:1.3125,1)
--(axis cs:1.296875,1)
--(axis cs:1.28125,1)
--(axis cs:1.265625,1)
--(axis cs:1.25,1)
--(axis cs:1.234375,1)
--(axis cs:1.21875,1)
--(axis cs:1.203125,1)
--(axis cs:1.1875,1)
--(axis cs:1.171875,1)
--(axis cs:1.15625,1)
--(axis cs:1.140625,1)
--(axis cs:1.125,1)
--(axis cs:1.109375,1)
--(axis cs:1.09375,1)
--(axis cs:1.078125,1)
--(axis cs:1.0625,1)
--(axis cs:1.046875,1)
--(axis cs:1.03125,1)
--(axis cs:1.015625,1)
--(axis cs:1,1)
--(axis cs:0.984375,1)
--(axis cs:0.96875,1)
--(axis cs:0.953125,1)
--(axis cs:0.9375,1)
--(axis cs:0.921875,1)
--(axis cs:0.90625,1)
--(axis cs:0.890625,1)
--(axis cs:0.875,1)
--(axis cs:0.859375,1)
--(axis cs:0.84375,1)
--(axis cs:0.828125,1)
--(axis cs:0.8125,1)
--(axis cs:0.796875,1)
--(axis cs:0.78125,1)
--(axis cs:0.765625,1)
--(axis cs:0.75,1)
--(axis cs:0.734375,1)
--(axis cs:0.71875,1)
--(axis cs:0.703125,1)
--(axis cs:0.6875,1)
--(axis cs:0.671875,1)
--(axis cs:0.65625,1)
--(axis cs:0.640625,1)
--(axis cs:0.625,1)
--(axis cs:0.609375,1)
--(axis cs:0.59375,1)
--(axis cs:0.578125,1)
--(axis cs:0.5625,1)
--(axis cs:0.546875,1)
--(axis cs:0.53125,1)
--(axis cs:0.515625,1)
--(axis cs:0.5,1)
--(axis cs:0.484375,1)
--(axis cs:0.46875,1)
--(axis cs:0.453125,1)
--(axis cs:0.4375,1)
--(axis cs:0.421875,1)
--(axis cs:0.40625,1)
--(axis cs:0.390625,1)
--(axis cs:0.375,1)
--(axis cs:0.359375,1)
--(axis cs:0.34375,1)
--(axis cs:0.328125,1)
--(axis cs:0.3125,1)
--(axis cs:0.296875,1)
--(axis cs:0.28125,1)
--(axis cs:0.265625,1)
--(axis cs:0.25,1)
--(axis cs:0.234375,1)
--(axis cs:0.21875,1)
--(axis cs:0.203125,1)
--(axis cs:0.1875,1)
--(axis cs:0.171875,1)
--(axis cs:0.15625,1)
--(axis cs:0.140625,1)
--(axis cs:0.125,1)
--(axis cs:0.109375,1)
--(axis cs:0.09375,1)
--(axis cs:0.078125,1)
--(axis cs:0.0625,1)
--(axis cs:0.046875,1)
--(axis cs:0.03125,1)
--(axis cs:0.015625,1)
--(axis cs:0,1)
--cycle;

\path [draw=darkgray, fill=darkgray] (axis cs:2.203125,0.59375)
--(axis cs:2.203125,-1)
--(axis cs:2.21875,-1)
--(axis cs:2.234375,-1)
--(axis cs:2.234375,0.53125)
--(axis cs:2.234375,0.53125)
--(axis cs:2.21875,0.5625)
--(axis cs:2.203125,0.59375)
--cycle;

\path [draw=black, line width=0.83630769230769231pt] (axis cs:2.203125,-1)
--(axis cs:2.203125,0.59375);

\path [draw=black, line width=0.83630769230769231pt] (axis cs:2.234375,-1)
--(axis cs:2.234375,0.53125);

\addplot [line width=1.6726153846153846pt, black, forget plot]
table {%
0 1
0.015625 1
0.03125 1
0.046875 1
0.0625 1
0.078125 1
0.09375 1
0.109375 1
0.125 1
0.140625 1
0.15625 1
0.171875 1
0.1875 1
0.203125 1
0.21875 1
0.234375 1
0.25 1
0.265625 1
0.28125 1
0.296875 1
0.3125 1
0.328125 1
0.34375 1
0.359375 1
0.375 1
0.390625 1
0.40625 1
0.421875 1
0.4375 1
0.453125 1
0.46875 1
0.484375 1
0.5 1
0.515625 1
0.53125 1
0.546875 1
0.5625 1
0.578125 1
0.59375 1
0.609375 1
0.625 1
0.640625 1
0.65625 1
0.671875 1
0.6875 1
0.703125 1
0.71875 1
0.734375 1
0.75 1
0.765625 1
0.78125 1
0.796875 1
0.8125 1
0.828125 1
0.84375 1
0.859375 1
0.875 1
0.890625 1
0.90625 1
0.921875 1
0.9375 1
0.953125 1
0.96875 1
0.984375 1
1 1
1.015625 1
1.03125 1
1.046875 1
1.0625 1
1.078125 1
1.09375 1
1.109375 1
1.125 1
1.140625 1
1.15625 1
1.171875 1
1.1875 1
1.203125 1
1.21875 1
1.234375 1
1.25 1
1.265625 1
1.28125 1
1.296875 1
1.3125 1
1.328125 1
1.34375 1
1.359375 1
1.375 1
1.390625 1
1.40625 1
1.421875 1
1.4375 1
1.453125 1
1.46875 1
1.484375 1
1.5 1
1.515625 1
1.53125 1
1.546875 1
1.5625 1
1.578125 1
1.59375 1
1.609375 1
1.625 1
1.640625 1
1.65625 1
1.671875 1
1.6875 1
1.703125 1
1.71875 1
1.734375 1
1.75 1
1.765625 1
1.78125 1
1.796875 1
1.8125 1
1.828125 1
1.84375 1
1.859375 1
1.875 1
1.890625 1
1.90625 1
1.921875 1
1.9375 1
1.953125 1
1.96875 1
1.984375 1
2 1
2.015625 0.96875
2.03125 0.9375
2.046875 0.90625
2.0625 0.875
2.078125 0.84375
2.09375 0.8125
2.109375 0.78125
2.125 0.75
2.140625 0.71875
2.15625 0.6875
2.171875 0.65625
2.1875 0.625
2.203125 0.59375
2.21875 0.5625
2.234375 0.53125
2.25 0.5
2.265625 0.46875
2.28125 0.4375
2.296875 0.40625
2.3125 0.375
2.328125 0.34375
2.34375 0.3125
2.359375 0.28125
2.375 0.25
2.390625 0.21875
2.40625 0.1875
2.421875 0.15625
2.4375 0.125
2.453125 0.09375
2.46875 0.0625
2.484375 0.03125
2.5 0
2.515625 -0.03125
2.53125 -0.0625
2.546875 -0.09375
2.5625 -0.125
2.578125 -0.15625
2.59375 -0.1875
2.609375 -0.21875
2.625 -0.25
2.640625 -0.28125
2.65625 -0.3125
2.671875 -0.34375
2.6875 -0.375
2.703125 -0.40625
2.71875 -0.4375
2.734375 -0.46875
2.75 -0.5
2.765625 -0.53125
2.78125 -0.5625
2.796875 -0.59375
2.8125 -0.625
2.828125 -0.65625
2.84375 -0.6875
2.859375 -0.71875
2.875 -0.75
2.890625 -0.78125
2.90625 -0.8125
2.921875 -0.84375
2.9375 -0.875
2.953125 -0.90625
2.96875 -0.9375
2.984375 -0.96875
3 -1
};
\addplot [line width=1.6726153846153846pt, black, forget plot]
table {%
0 1
0.015625 0.96875
0.03125 0.9375
0.046875 0.90625
0.0625 0.875
0.078125 0.84375
0.09375 0.8125
0.109375 0.78125
0.125 0.75
0.140625 0.71875
0.15625 0.6875
0.171875 0.65625
0.1875 0.625
0.203125 0.59375
0.21875 0.5625
0.234375 0.53125
0.25 0.5
0.265625 0.46875
0.28125 0.4375
0.296875 0.40625
0.3125 0.375
0.328125 0.34375
0.34375 0.3125
0.359375 0.28125
0.375 0.25
0.390625 0.21875
0.40625 0.1875
0.421875 0.15625
0.4375 0.125
0.453125 0.09375
0.46875 0.0625
0.484375 0.03125
0.5 0
0.515625 -0.03125
0.53125 -0.0625
0.546875 -0.09375
0.5625 -0.125
0.578125 -0.15625
0.59375 -0.1875
0.609375 -0.21875
0.625 -0.25
0.640625 -0.28125
0.65625 -0.3125
0.671875 -0.34375
0.6875 -0.375
0.703125 -0.40625
0.71875 -0.4375
0.734375 -0.46875
0.75 -0.5
0.765625 -0.53125
0.78125 -0.5625
0.796875 -0.59375
0.8125 -0.625
0.828125 -0.65625
0.84375 -0.6875
0.859375 -0.71875
0.875 -0.75
0.890625 -0.78125
0.90625 -0.8125
0.921875 -0.84375
0.9375 -0.875
0.953125 -0.90625
0.96875 -0.9375
0.984375 -0.96875
1 -1
1.015625 -1
1.03125 -1
1.046875 -1
1.0625 -1
1.078125 -1
1.09375 -1
1.109375 -1
1.125 -1
1.140625 -1
1.15625 -1
1.171875 -1
1.1875 -1
1.203125 -1
1.21875 -1
1.234375 -1
1.25 -1
1.265625 -1
1.28125 -1
1.296875 -1
1.3125 -1
1.328125 -1
1.34375 -1
1.359375 -1
1.375 -1
1.390625 -1
1.40625 -1
1.421875 -1
1.4375 -1
1.453125 -1
1.46875 -1
1.484375 -1
1.5 -1
1.515625 -1
1.53125 -1
1.546875 -1
1.5625 -1
1.578125 -1
1.59375 -1
1.609375 -1
1.625 -1
1.640625 -1
1.65625 -1
1.671875 -1
1.6875 -1
1.703125 -1
1.71875 -1
1.734375 -1
1.75 -1
1.765625 -1
1.78125 -1
1.796875 -1
1.8125 -1
1.828125 -1
1.84375 -1
1.859375 -1
1.875 -1
1.890625 -1
1.90625 -1
1.921875 -1
1.9375 -1
1.953125 -1
1.96875 -1
1.984375 -1
2 -1
2.015625 -1
2.03125 -1
2.046875 -1
2.0625 -1
2.078125 -1
2.09375 -1
2.109375 -1
2.125 -1
2.140625 -1
2.15625 -1
2.171875 -1
2.1875 -1
2.203125 -1
2.21875 -1
2.234375 -1
2.25 -1
2.265625 -1
2.28125 -1
2.296875 -1
2.3125 -1
2.328125 -1
2.34375 -1
2.359375 -1
2.375 -1
2.390625 -1
2.40625 -1
2.421875 -1
2.4375 -1
2.453125 -1
2.46875 -1
2.484375 -1
2.5 -1
2.515625 -1
2.53125 -1
2.546875 -1
2.5625 -1
2.578125 -1
2.59375 -1
2.609375 -1
2.625 -1
2.640625 -1
2.65625 -1
2.671875 -1
2.6875 -1
2.703125 -1
2.71875 -1
2.734375 -1
2.75 -1
2.765625 -1
2.78125 -1
2.796875 -1
2.8125 -1
2.828125 -1
2.84375 -1
2.859375 -1
2.875 -1
2.890625 -1
2.90625 -1
2.921875 -1
2.9375 -1
2.953125 -1
2.96875 -1
2.984375 -1
3 -1
};
\node at (axis cs:2.115,-1.2)[
  scale=0.5,
  anchor=base west,
  text=black,
  rotate=0.0
]{ $E(x)$};
\end{axis}

\end{tikzpicture}

%% file: chapterDiscussion.tex
\chapter{Discussion and future directions}\label{Chapter:Discussion}
The central building block to include adhesive interactions between cells in
reaction-advection-diffusion model of tissues is to use a non-local term. The
simplest scalar non-local partial differential equation including such a term is
given by:
\begin{equation}\label{eq:adh}
    u_t(x,t) = u_{xx}(x, t) - \alpha \lb u(x,t) \int_{-R}^{R} h(u(x+r,t)) \Omega(r) \dd r \rb,
\end{equation}
here $u(x,t)$ denotes the density of a cell population which adheres to itself,
$h(u)$ measures the adhesive strength\index{adhesive strength} with the
background population, and $\Omega(r)$ describes the distribution of adhesive
sites on the cell membrane\index{cell membrane} and on cell
protrusions\index{cell protrusions}. This served as our prototypical example to
develop the mathematical tools to investigate the steady-state structure of such
an equation, and their stability.  In addition, having biological applications
we posed the question how to correctly formulate no-flux boundary conditions for
such non-local models, and how these boundary conditions affect the equation's
steady states.

In \cref{chapter:BasicProperties} we study the properties of the non-local
operator $\K$. We observe that the non-local operator is a generalization of the
classical first-order derivative, and shares many of its properties. In fact, as
the sensing radius tends to zero, the non-local operator converges to the
classical derivative. It also has eigenvalues that are strongly reminiscent of
the classical derivative. Unsurprisingly, however, since it is a non-local
operator, it is indeed compact, and its eigenvalues must accumulate at zero.

The steady-states of the periodic problem are discussed in \cref{part:periodic}.
The advantage of the periodic problem is that it allows us to discuss the steady
states in absence of any boundary effects. Since the trivial solution exists for
all values of $\alpha$, we use the abstract global bifurcation\index{abstract
global bifurcation} results pioneered
by Rabinowitz~\cite{Rabinowitz1971}, to identify bifurcation points of
non-trivial solutions. Since the non-local operator inherits the same symmetry
properties\index{symmetry properties} as the classical first order derivative,
we can classify the globally existing solutions branches, using the symmetries
of the adhesion equation.  The
symmetries of each bifurcation branch allow us to show that each branch is
contained in a separate space of spiky functions (\ie\ piecewise monotone
functions having a fixed number of zeros). Due to the symmetries the location's
of the peaks and valleys of these solutions are fixed, and uniformly spaced. Due
to the symmetries all bifurcations are of pitch-fork type, and their criticality
is solely determined by the properties of the kernel $\omega(r)$ in the non-local
term.

The limitations of our main global bifurcation result, namely the restriction to
linear functions $h(u) = u$ and constant kernel $\omega(r)$, is invitation for
much future work. Indeed, it is not trivial to improve on either limitation. The
main challenge one encounters trying to improve upon both limitations is
generalizing the ``non-local'' maximum principle\index{non-local maximum
principle}. The current proof relies on the relationship between the area
function of a solution $u(x)$ and the non-local operator $\K$.

The next steps in the analysis of the dynamics of this scalar non-local equation
is to study its time-dependent solutions. This is both of mathematical and
biological interest. Questions from a modelling perspective include:
Can cell aggregates merge? Do there exist meta-stable solutions? What
happens to aggregates as $t\to\infty$. In mathematical terms we are seeking a
full characterization of the attractor\index{attractor} of the scalar non-local adhesion
equation. We believe that some of the mathematical methods we developed in the
analysis here will be invaluable in answering such questions in the future.

\section{Further thoughts}

A different kind of result, has been shown for semilinear scalar
parabolic\index{semilinear scalar parabolic equation}
equation for which it is known that the so called
\textit{lap-number}\index{lap-number} (i.e.\ the
number of zeros of function) is a non-increasing function with respect to
time~\cite{Matano1982}. This means that solutions can only become simpler as
time progresses, and all complexity must be present in the initial conditions.
In~\cite{matano1988} this result was extended to scalar equations on the unit
circle $S^1$. There however we have two possibilities either the solutions
approach a set of fixed points, or the solutions resemble a rotating wave.
Further, the rotating wave solution becomes impossible when the equation's
nonlinearity posses a reflection symmetry in its derivative dependence. In the
case of the scalar non-local equation~\eqref{eq:adh} extensive numerical
exploration suggest that the number of spikes (i.e.\ number of zeros of the
solution's derivatives) is non-increasing as a function of time.


The Chafee-Infante equation\index{Chafee-Infante equation} is one of the most
studied nonlinear scalar reaction diffusion equations. Non-local extensions to
the Chafee-Infante problem have
been studied in~\cite{Raquepas1999}, and similar non-local problems 
in~\cite{Davidson2006c,Freitas2000a,Freitas1994,Freitas1994a}.
The non-local term interacts with the regular spectrum of the Laplacian by
moving some eigenvalues~\cite{Freitas1994,Freitas1994a}. Depending on the
severity of this ``scrambling'' this leads to the stabilization or
destabilization of spatial patterns, and to changes in flow direction on the
global attractor~\cite{Raquepas1999}.
These spectral results are similar to the result we observed for our canonical
example of a non-local equation~\eqref{eq:adh}. Whether similar results hold for
the global attractor remains to be seen.

For scalar nonlinear reaction diffusion equations it is proven
in~\cite{Henry1985,Hale1988} that if $\phi, \psi$ are two hyperbolic equilibria then
$W^u(\phi)$ is transversal to $W^s(\psi)$. This means that dynamical changes on
the global attractor due to variations in the nonlinearity can only occur due to
a bifurcation. Does a similar result apply to our canonical non-local
equation~\eqref{eq:adh}?

As shown in \cref{appendix:adh_pot} the adhesion model~\eqref{eq:adh} can be
formulated as an aggregation equation, and as a gradient flow on the appropriate
space of probability measures. A rich theory is available to study the
properties of such equations, and it is of great interest to observe how the
properties uncovered in this monograph fit with that theory.

\section{Systems and higher dimensions}

The scalar non-local adhesion equation studied in this monograph served as a
kind of canonical example, on which we developed the mathematical tools
necessary to tackle more complicated and more biologically realistic models.
Originally the development of the non-local adhesion model was motivated by the
well studied phenomena of cell-sorting\index{cell-sorting}, in which cells of two distinct
populations can form different tissue structures depending solely on their
adhesive properties. These effects are crucial during the development of
organisms. The simplest model of cell-sorting considered in~\cite{Armstrong2006} is
\begin{equation}\label{adhSystem}
\begin{split}
    u_t &= u_{xx} - \alpha_{uu} \lb u(x,t) \K[u](x,t) \rb_x - \alpha_{uv} \lb u(x,t) \K[v](x, t) \rb  \\
    v_t &= v_{xx} - \alpha_{vv} \lb v(x,t) \K[v](x,t) \rb_x - \alpha_{vu} \lb v(x,t) \K[u](x, t) \rb,
\end{split}
\end{equation}
where $u(x,t)$ and $v(x,t)$ denote the densities of the two cell populations,
and $\alpha_{ij}$ denotes the strength of the 
adhesive strength\index{adhesive strength} of cell
population $i$ to $j$. Similarly, to the case of scalar reaction-diffusion
equations we expect the set of possible steady states to be much more
complicated in this case (see for instance~\cite{Nishiura1982,Fujii1983}, in
both symmetries were crucial to understand the equation's steady-states and time
dependent solutions). Numerical simulations in these cases hint that
symmetries remain important.  As discussed in the Introduction, systems of the
form~\eqref{adhSystem} have been used extensively in cancer
modelling\index{cancer modelling}.

Other extensions to higher spatial dimensions would enable us to explore
situations that are more closely connected with experimental situations, such as
models of wound-healing and tissue formation.

Yet another different avenue would be to rigorously study the dynamics of this
equation with spatio-temporally varying adhesion coefficients. Such variations
either due to environmental changes or genetic or epigenetic variations are
highly important in biological settings. For instance, it is clear that upon
wounding cells change their behaviour to close the wound, and then return to
their initial state. Another situation in which this is relevant is cancer
growth and cancer metastasis, this has been demonstrated by numerical studies
in~\cite{Domschke2014}. Other extensions include considering
volume-filling\index{volume-filling} effects
and nonlinear diffusion\index{nonlinear diffusion}, which both have been shown
to be important in biological applications~\cite{Murakawa2015,Carrillo2019}.

In conclusion, in this monograph we made significant progress in understanding
the mathematical properties of non-local models of cell adhesion. We proposed a
framework to derive such models from an underlying stochastic random walk,
studied the steady states spawning through bifurcations from the constant
solution, and finally considered the non-local cell adhesion model on a bounded
domain with no-flux boundary conditions. Finally, perhaps one of the most
important outcomes of this monograph are the new possible directions of research
that have been identified.  The main questions are: identification of secondary
bifurcations, analysis of the time-dependent solutions of
equation~\eqref{Eqn:ArmstrongModelIntro}, several cell populations, higher space
dimensions, and multi-scale modelling including intra-cellular signalling
networks.

%% file: AppendixAmorim.tex
%
\chapter[Adhesion Potential]{Adhesion Potential}\label{appendix:adh_pot}
{\em We thank Paulo Amorim (Federal University of Rio de Janeiro) for a fruitful discussion on this topic}\\

It has been discussed in the community if the adhesion model can be written as
an aggregation equation\index{aggregation equation}, as used in the swarming
literature (see for example \cite{Muntean,Carrillo2019}). Here we show that this
is indeed the case in a
very general $n$-dimensional framework. This relation opens the doors to the
rich theory of measure-valued gradient flows\index{measure-valued gradient flows} 
and many of their results will apply here as well.

The adhesion model of Armstrong, Painter, Sherratt \cite{Armstrong2006} in $n$
dimensions has the following form.
\begin{equation}\label{adhesion}
u_t = D\Delta u - \alpha \nabla \cdot \left(u \int_{B_R(0)} h(u(x+y,t)) \Omega(y) \, \dd y \right)
\end{equation}
where the adhesion force density $h(u)$ is an increasing smooth function in $u$
and the weight function $\Omega(y)$ satisfies:
\begin{eqnarray*}
\Omega(y) &=& \omega(\|y\|) \frac{y}{\|y\|}\\
\omega(\|y\|) \geq 0, &&  \omega\in L^1(B_R(0))\cap L^\infty(B_R(0)), \quad \int_{B_R(0)^+} \omega(\|y\|) \dd y = \frac{1}{2},
\end{eqnarray*}
where $B_R(0)^+$ denotes the upper-half of the ball of radius $R$, where we can
define the upper half by using any coordinate axis and consider non-negative
components on this axis.

We assume that  $\omega(r)$ has an antiderivative $V$. For physical reasons, we
work with the negative of the antiderivative. I.e.\ we assume
\[ V'(r) = - \omega(r).\]
Then
\[ \nabla_y  V (\|y\|) =- \omega(\|y\|) \frac{y}{\|y\|}= -\Omega(y).\]
Hence $\Omega(y)$ has a potential, which we call the {\it adhesion potential}.
As the adhesion potential\index{adhesion potential} is unique up to addition 
of a constant, we can always chose $V$ such that
\[V(R)=0.\]
Then we can write the integral term in \eqref{adhesion} as
\begin{eqnarray*}
\int_{B_R(0)} h(u(x+y,t)) \Omega(y) \dd y &=& -\int_{B_R(0)} h(u(x+y, t)) \nabla_y V(\|y\|) \dd y \\
 &=& \int_{B_R(0)} \nabla_y h(u(x+y,t)) V(\|y\|) \dd y \\
 && - \int_{\partial B_R(0)} h(u(x+y,t)) V(\|y\|) \dd y \\
 &=& \int_{B_R(0)} \nabla_x h(u(x+y,t)) V(\|y\|) \dd y\\
 &=& \nabla_x \int_{B_R(0)}  h(u(x+y,t)) V(\|y\|) \dd y
\end{eqnarray*}
With this transformation, the adhesion equation \eqref{adhesion} becomes
\[
u_t = D\Delta u - \alpha \nabla \cdot \left( u \nabla \int_{B_R(0)} h(u(x+y),t) V(\|y\|) \dd y\right). \]
Finally, to write the integral term as a convolution, we use the fact that
$V(R)=0$ and we cut the potential at radius $R$:
\[ W(r) := V(r) \chi_{[0,R]}(r). \]
Then we obtain a classical aggregation equation
\begin{equation}\label{aggregation}
u_t = D\Delta u - \alpha \nabla \cdot \bigl[u \,\nabla  (W \ast h(u))\bigr].
\end{equation}
Using the known results for the aggregation equation, we then know that
(\ref{aggregation}) is the gradient flow of the energy
\[ J(u) = \int D\frac{u^2}{2} \dd x - \frac{\alpha}{2} \int u \, (W\ast h(u)) \dd x \]
using the $d_2$ Wasserstein metric \cite{Ambrosio,Muntean}.

The aggregation equation\index{aggregation equation} on bounded
domains\index{bounded domains} has recently been studied in
\cite{fetecau2017swarm,wu2015nonlocal} and boundary conditions were defined,
based on the underlying energy principle. In some cases, these boundary
conditions are similar to the adhesive or repulsive boundary
conditions\index{repulsive boundary conditions} that we
discussed in \cref{Chapter:AdhNoFlux}.

Here we simply state the fact that the adhesion model can be seen as an
aggregation model, where the potential is given as the (negative)
anti-derivative of the sensing function $\Omega$. The implications of this
analogy have not been worked out yet, and future research will make the rich
theory of aggregation equations available to study cell-cell adhesion.